\documentclass[11pt]{article}

\usepackage{amssymb}
\usepackage{amstext}
\usepackage{amsmath}
\usepackage{latexsym}
\usepackage[all]{xy}
\usepackage{tikz}
\usetikzlibrary{decorations.markings}
\usetikzlibrary{calc,intersections,through,backgrounds}
\usetikzlibrary{arrows.meta,positioning,fit,petri}
\usetikzlibrary{decorations.pathmorphing}
\usetikzlibrary{plothandlers}
\usetikzlibrary{plotmarks}
\usetikzlibrary{shadings}
\usetikzlibrary{shapes.geometric}
\usetikzlibrary{datavisualization.formats.functions}

\newfont{\gothique}{eufm10 scaled 1100}  

\newcommand{\PP}{{\mathbb{P}}}
\newcommand{\CC}{{\bf{C}}}
\newcommand{\ZZ}{{\mathbb{Z}}}
\newcommand{\RR}{{\mathbb{R}}}
\newcommand{\NN}{{\mathbb{N}}}  
\newcommand{\HH}{H^1(\Theta_{C})}
\newcommand{\TET}{\Theta_{X}}
\newcommand{\KX}{K_{X}}
\newcommand{\HKC}{H^0(\OO_C (K_{C}))}
\newcommand{\TT}{{\cal{T}}}
\newcommand{\GG}{{\cal{G}}}
\newcommand{\EE}{{\cal{E}}}
\newcommand{\OM}{\Omega_{C}}

\newcommand{\OO}{{\cal{O}}}

\newcommand{\HA}{\frac{1}{2}}

\newcommand{\FF}{{\cal{F}}}
\newcommand{\TA}{\tau_{{}_{\Pi,l}}}

\newtheorem{thm}{Theorem}[section]
\newtheorem{lem}[thm]{Lemma}
\newtheorem{pro}[thm]{Proposition}
\newtheorem{cor}[thm]{Corollary}
\newtheorem{rem}[thm]{Remark}
\newtheorem{cl}[thm]{Claim}
\newtheorem{example}[thm]{Example}
\newtheorem{defi}[thm]{Definition}  

\newlength{\myskip}
\newenvironment{pf}{
     \addvspace{\myskip}  

     \noindent {\it Proof.$\, $}}
     {$\Box$

     \addvspace{\myskip}
     }

\makeatletter
\renewcommand{\@seccntformat}[1]{\S \/ {\csname the#1\endcsname}\hspace{0.5em}}
\makeatother

\title{\bf REFINEMENT OF THE INFINITESIMAL VARIATION OF HODGE STRUCTURE: THE CASE OF CANONICAL CURVES
           }
\author{Igor Reider}

\begin{document}

\bibliographystyle{amsplain}
\maketitle

\setcounter{section}{-1}
\numberwithin{equation}{section}
\begin{abstract}
Let $C$ be a smooth complex projective curve with canonical divisor $K_C$ very ample.
The general objective of the paper is to study the relation between the cup-product
$$
H^1 (\Theta_C ) \longrightarrow \HKC^{\ast} \otimes H^1 (\OO_C )
$$
where $\Theta_C =\OO_C (-K_C)$ is the holomorphic tangent bundle of $C$, and the geometry of the canonical embedding of $C$. The cup-product, following Griffiths, stratifies $\PP(H^1 (\Theta_C ))$ by the subvarieties $\Sigma_r$ according to the rank $r$ of $\xi \in H^1 (\Theta_C )$, viewed as the linear map
$$
\xi:\HKC \longrightarrow  H^1 (\OO_C),
$$
or, equivalently, by the dimension of the kernel of $\xi$
$$
W_{\xi}=ker(\xi).
$$
The pair $(\xi, W_{\xi})$ is one of the invariants of the Infinitesimal Variation of Hodge Structure (IVHS) of Griffiths and his collaborators. We construct a refinement of $(\xi, W_{\xi})$ by attaching an intrinsic filtration  $W^{\bullet}_{\xi} ([\phi])$  of $W_{\xi}$, varying holomorphically with $[\phi] \in \PP(W_{\xi})$.
The filtration has geometric meaning:

1) it is related to special divisors on $C$,

2) it `counts' certain rational normal curves in the canonical embedding
of $C$.

The filtration arises from the skew-symmetric pairing
$$
\xymatrix{
	\alpha^{(2)}_{\xi}: \bigwedge^2 W_{\xi} \ar[r]&\HKC
}
$$
naturally attached to the IVHS pair $(\xi,W_{\xi})$; this could be considered as a `symplectic' structure on $W_{\xi}$. 

We illustrate how our refinement gives the well known results about the strata $\Sigma_0$ and $\Sigma_1$:

- $\Sigma_0$ is empty, if $C$ has no $g^1_2$,

- the stratum $\Sigma_1$ is the image of the bicanonical embedding of $C$ as long as $C$ has no $g^1_2$, $g^1_3$, $g^2_5$.

That is, one recovers as corollaries the classical theorems of Max Noether on projective normality of the canonical embedding and Babbage-Enriques-Petri about the canonical curve being cut out by quadrics.

Conceptually, our refinement concerns the incidence correspondence
$$
{\bf P}=\{ ([\xi], [\phi]) \in \PP(H^1 (\Theta_{C}))\times \PP(\HKC) | \text{$\xi \phi=0$ in $H^1(\OO_C)$}\}
$$
which comes with the projections on each factor
$$
\xymatrix{
	&{\bf P}\ar_{p_1}[dl]\ar^{p_2}[dr]&\\
	\PP(H^1 (\Theta_{C}))&&\PP(\HKC)
}
$$
The left side of the diagram controls the (infinitesimal) variations of complex structure on $C$. The refinement exhibits additional structures on the fibres of $p_1$. In particular, one obtains the stratification of the inverse images $p^{-1}(\Sigma_r)$ of the Griffiths' strata $\Sigma_r$ into substrata where the length of the filtrations $W^{\bullet}_{\xi} (\phi)$, $([\xi],[\phi])\in p^{-1}(\Sigma_r)$, remains constant. On each such substratum new aspects emerge:

- quiver representations,

- Fano toric variety with a distinguished anti-canonical divisor,

- dimer models.

The quiver emerges from the construction and properties of the refinement; the Fano toric variety arises formally from the graph underling the quiver, but it also has a meaning of moduli of a sort of Higgs structures of the linear algebra of the refinement. The graph underlying the refinement becomes an important part of the theory: it connects to topics such as the Topological Quantum field theory, moduli spaces of elliptic curves with marked points, modular curves, higher categorical structures of the refinements.
%
\end{abstract}
\tableofcontents
\section{Introduction}

One of the basic and remarkable manifestation of complex projective structure on a compact 
$C^{\infty}$ manifold is the Hodge decomposition of the de Rham complex cohomology groups. Namely, for a complex projective manifold $X$ of dimension $n$, the cohomology groups of $X$ with complex coefficients $H^k(X,\CC)$ admit the decomposition into the direct sum
$$
H^k(X,\CC)=\bigoplus_{p+q=k} H^{p,q}(X),
$$
for $k\in[0,2n]$. This is the celebrated Hodge decomposition with the summands $H^{p,q}(X)$, the Hodge groups of $X$. Those groups are in turn identified as the cohomology groups of the exterior powers of the basic sheaf $\Omega_X$, the holomorphic cotangent bundle of $X$: for $\Omega^p_X:=\bigwedge^p \Omega_X$, the sheaf of germs of holomorphic $p$-forms, one has the Dolbeault isomorphism
$$
H^{p,q}(X)\cong H^q(\Omega^p_X).
$$
One of the fundamental ideas of algebraic geometry which goes back to Riemann is to study not just individual varieties but also consider all possible ways to make the underlying $C^{\infty}$ manifold into a complex projective variety. This means that one varies complex structure of $X$. This leads to the problem of the moduli space of isomorphism classes of complex structures. This is important because moduli spaces often provide
examples of varieties with interesting properties and also because properties of individual varieties themselves are encoded in a beautiful way in moduli spaces. The variation of complex structure on $X$
raises another important problem: how does the Hodge structure on the de Rham cohomology groups vary as one varies the complex structure? This leads to the theory of Variation of Hodge structure (VHS), period maps, period domains initiated by Griffiths in 1960's, see \cite{G1} for an overview. 

On the infinitesimal level the variation of complex structure on $X$ is given by the first cohomology $H^1(\Theta_X)$ of the holomorphic tangent bundle $\Theta_X$ of $X$, the dual of the cotangent bundle $\Omega_X$. This goes back to the works of Kodaira and Spencer, see \cite{KS}. The linearization of the period map is its differential: it ties together the Kodaira-Spencer space $H^1(\Theta_X)$ with the Hodge decomposition of $H^k(X,\CC)$ via the cohomology cup product
\begin{equation}\label{diff-period-intro}
H^1 (\Theta_X) \longrightarrow \bigoplus_{p+q=k} Hom (H^q (\Omega^p_X), H^{q+1} (\Omega^{p-1}_X))
\end{equation}
coming from the obvious contraction morphism
$$
\Theta_X \longrightarrow {\cal H}om (\Omega^p_X, \Omega^{p-1}_X).
$$     
The map in \eqref{diff-period-intro} and its properties are referred to as the Infinitesimal Variation of Hodge structure, the IVHS for short.  Heuristically, the differential of the period map turns the Kodaira-Spencer classes, the elements of $H^1 (\Theta_X)$, into dynamical entities by viewing them as endomorphisms of the graded vector space
\begin{equation}\label{Hk-graded-intro}
H^k(X,\CC)=\bigoplus_{p+q=k}H^q(\Omega^p_X):
\end{equation}
the dynamics is expressed by how a given cohomology class sitting in a graded piece $H^q(\Omega^p_X)$ propagates through the space under the (iterated) action of the Kodaira-Spencer classes.

In a series of papers Griffiths and his collaborators proposed several constructions  naturally attached to the Infinitesimal Variation of Hodge structure. Those are known collectively as the IVHS constructions or invariants, see \cite{CGrGH}, \cite{GH}, \cite{G} and also \cite{G1}, for an overview and additional references on the subject.
  As an illustration and one of the successful applications of the IVHS constructions one often cites the case $n=1$, that is $X$ is a smooth complex projective curve and we switch to using $C$ instead of $X$. The  cup product \eqref{diff-period-intro} for $k=1$ becomes
$$
H^1 (\Theta_C) \longrightarrow  Hom (H^0 (\Omega_C), H^{1} (\OO_C)).
$$
Here the cotangent bundle $\Omega_C$ is the canonical line bundle of $C$ which we write
$\OO_C (K_C)$, where $K_C$ is the canonical divisor class of $C$, the 
tangent bundle $\Theta_C=\OO_X (-K_C)$, the dual of $\OO_C (K_C)$, and the above cup product becomes
\begin{equation}\label{cp-C}
	\begin{gathered}
H^1 (\OO_C (-K_C)) \longrightarrow  Hom (H^0 (\OO_C (K_C)), H^{1} (\OO_C)) =H^0 (\OO_C (K_C))^{\ast} \otimes H^1 (\OO_C) \\
\cong H^0 (\OO_C (K_C))^{\ast} \otimes H^0 (\OO_C (K_C))^{\ast},
\end{gathered}
\end{equation}
where the last equality comes from Serre duality
$H^1 (\OO_C) \cong H^0 (\OO_C (K_C))^{\ast}$.

One of the insights of Griffiths about the differential of the period map could be summarized as follows:

  the study of the {\it degeneracies} of the cup product involved in the differential of the period map leads to important facts about geometry of the projective varieties as well as their moduli spaces.
  
In fact, Griffiths in \cite{G} attaches to the cup product \eqref{cp-C} the determinantal varieties
\begin{equation}\label{r-stratum}
	\Sigma_r:=\{[\xi]\in \PP(H^1(\Theta_C)) |\,\, \xi: \HKC \longrightarrow H^1 (\OO_C),\,\, rk(\xi) \leq r \},
\end{equation}
for all $r \in [0,g]$, which are the degeneracy loci of the morphism of sheaves
$$
\HKC \otimes \OO_{\PP(H^1(\Theta_C))} (-1)  \longrightarrow H^1(\OO_C) \otimes \OO_{\PP(H^1(\Theta_C))}.
$$
This gives the stratification
\begin{equation}\label{rkstrat}
	\PP(H^1(\Theta_C)) =\Sigma_g \supset \Sigma_{g-1} \supset \cdots \supset \Sigma_1 \supset \Sigma_0.
\end{equation}
 The stratum $\Sigma_0$ controls the kernel of the differential of the period  map. Thus the statement
 \begin{equation}\label{Sigma0-intro}
 \text{\it $\Sigma_0$ is empty}
\end{equation}
 	is equivalent to the Infinitesimal Torelli for curves. The dual version of 
\eqref{cp-C} together with \eqref{Sigma0-intro} becomes the classical theorem of Max Noether: the canonical embedding of a curve is projectively normal, see \cite{G-H}. Thus the knowledge about the degeneracy stratum $\Sigma_0$ is equivalent to an important geometric fact - the projective normality of the canonical embedding, {\it and} the Infinitesimal Torelli theorem - the period map is an immersion away from the locus of hyperelliptic curves. 

Griffiths goes further and shows that the next stratum, 
$\Sigma_1 $, is precisely the image of $C$ under the bicanonical morphism, the one defined by $\OO_C (2K_C)$, provided one knows another fundamental fact in the theory of curves, the Babbage-Enriques-Petri
theorem. That theorem says that for smooth, complex projective curves with no special linear systems $g^1_2$, $g^1_3$, $g^2_5$, the canonical curve is determined as the intersection of all quadrics passing through it, see \cite{St-D}, also \cite{M} for an inspiring overview. Thus Griffiths theorem
$$
\text{\it \{ $C$ has no $g^1_2$, $g^1_3$, $g^2_5$\} $\Leftrightarrow$ \{$\Sigma_1 =$ the bicanonical image of $C$\}}
$$
 tells us that the degeneracy locus $\Sigma_1$ is tied to the important fact about the geometry of the canonical embedding, the ideal is generated by quadrics, and provides a proof of the generic global Torelli theorem for curves. 

Heuristically, Griffiths stratification of the (projectivized) space of Kodaira-Spencer classes measures the `energy' of those via the size of the {\it kernel} of the corresponding endomorphism: a Kodaira-Spencer class of `higher' energy has bigger kernel - it annihilates a larger subspace of cohomology classes. There is a clear change of point of view on the differential of the period map: instead of looking at the `energy' of the cohomology classes - how far they propagate through the graded de Rham cohomology space - one
is concerned with the cohomology classes of {\it zero} energy with respect to the action of Kodaira-Spencer classes and it is the Kodaira-Spencer classes of high energy that are supposed to contain an important data about

- the moduli space of curves,

- the geometry of the curve, via the subspace of zero energy cohomology classes. 

\vspace{0.2cm}  
One would expect that the higher strata in \eqref{rkstrat}, in a way similar to the strata $\Sigma_0$ and $\Sigma_1$, contain interesting facts about the geometry of curves as well as their moduli spaces. It is with this thinking that the author embarked on the investigations of the strata $\Sigma_r$.
In this paper we propose a refinement of the Griffiths stratification \eqref{rkstrat}. This, in particular, will recover known results about the strata $\Sigma_0$ and $\Sigma_1$, {\it independently} of the old classical theorems quoted above. The refinement concerns the Kodaira-Spencer classes of rank $r$ parametrized by the locally closed strata
$$
\Sigma^0_r:=\Sigma_r \setminus \Sigma_{r-1}.
$$
By definition each stratum $\Sigma^0_r$ comes equipped with a distinguished vector bundle of rank $(g-r)$ whose fibre at a point $[\xi] \in \Sigma^0_r$ is the subspace
$$
W_{\xi}:=ker\Big(\HKC \stackrel{\xi}{\longrightarrow} H^1(\OO_C)\Big),
$$
of the canonical space $\HKC$. For every nonzero $\phi \in W_{\xi}$ our refinement attaches an intrinsic filtration of $W_{\xi}$ 
\begin{equation}\label{Wxiphi-filt-intro}
W_{\xi}=W^0_{\xi} ([\phi]) \supset W^1_{\xi} ([\phi]) \supset \cdots \supset  W^l_{\xi} ([\phi]) \supset W^{l+1}_{\xi} ([\phi])=0,
\end{equation}
where $[\phi]$ stands for the class of $\phi$ in the projective space $\PP(W_{\xi})$. Some of the properties of this filtration are summarized in the following statement.
\begin{thm}\label{pro:refine}
	The filtration \eqref{Wxiphi-filt-intro} is subject to the following properties.
	
	$\bullet$ $W^l_{\xi} ([\phi]) \supset \CC\phi$ and if the inclusion is strict, then there is a {\rm special} line bundle
	${\cal L}$ on $C$ such that 
	$
	H^0 ({\cal L}) 
	$
	is at least two dimensional: there is a plane $P_{\phi}$ in $W^l_{\xi} ([\phi])$ containing $\CC\phi$ and injecting into $H^0 ({\cal L}) $.
	
	$\bullet$ For every $i\in [2,l-1]$ (one assumes here that $l\geq 3$), the graded pieces $W^i_{\xi} ([\phi])/W^{i+1}_{\xi} ([\phi])$ of the filtration parametrize cones over rational normal curves of degree $i$ containing the hyperplane section
	$$
	Z_{\phi}=(\phi=0) \subset C \subset \PP(\HKC^{\ast})
	$$
	of $C$ in its canonical embedding.
\end{thm}

This could be viewed as a generalization of the features encountered in the discussion of the strata $\Sigma_i$, for $i=0,1$. Namely, the line bundle in the first item of the theorem should be viewed as a special linear system on $C$ (in the case of $\Sigma_0$ it is $g^1_2$, and for $\Sigma_1$ these are $g^1_2$, $g^1_3$, $g^2_5$ of the classical theorem of Babbage-Enriques-Petri), and the second item is related to quadrics through the canonical image of $C$.

\vspace{0.2cm}
More conceptually, our refinement takes place on the cohomological incidence correspondence
$$
{\bf P}=\{ ([\xi], [\phi]) \in \PP(H^1 (\Theta_{C}))\times \PP(\HKC) | \text{$\xi \phi=0$ in $H^1(\OO_C)$}\}
$$
which comes with the projections on each factor
\begin{equation}\label{two-proj}
\xymatrix{
	&{\bf P}\ar_{p_1}[dl]\ar^{p_2}[dr]&\\
	\PP(H^1 (\Theta_{C}))&&\PP(\HKC)
}
\end{equation}
The left side of the diagram controls the (infinitesimal) variations of complex structure on $C$ with the fibre over a point $[\xi]$ being the projectivization of the space $W_{\xi}$. The right hand side of the diagram is related to the canonical embedding of $C$ via the geometric incidence correspondence
$$
{\cal Z}_C = \{ ( [\phi],x) \in \PP(\HKC)\times C  | \phi (x)=0\} 
$$ 
and its two projections
$$
\xymatrix{
	&{\cal Z}_C\ar_{q_1}[dl]\ar^{q_2}[dr]&\\
	\PP(\HKC)&&C \subset \PP(\HKC^{\ast})
}
$$
Putting the two incidences together gives the diagram
\begin{equation}\label{two-inc-diag}
	\xymatrix{
		&{\bf P}\ar_{p_1}[dl]\ar^{p_2}[dr]&&{\cal Z}_C\ar_{q_1}[dl]\ar^{q_2}[dr]&\\
		\PP(H^1 (\Theta_{C}))&&\PP(\HKC)&&C \subset \PP(\HKC^{\ast})
	}
\end{equation}
On the left side we have the cohomological incidence correspondence and on the right side the geometric correspondence. One of the goals of IVHS constructions could be summarized as gaining information on the geometric side from relevant data on the cohomological side.

 Theorem \ref{pro:refine} provides an extra structure along the the fibres of the projection $p_1$ on the cohomological side and gives a partial translation of these data in terms of the geometry of $C$ and its canonical embedding. 
One can make a case, and this will be done in the main body of the paper, that our refinement is a sort of algebraic `K\"ahler structure' on $W_{\xi}$ depending on $[\phi]$ in $\PP(W_{\xi})$, a sort of hard Lefschetz decomposition of $W_{\xi}$ varying with $[\phi]$ in $\PP(W_{\xi})$. In a nutshell, once the `zero energy' space $W_{\xi}$ is sufficiently big, a new structure emerges in the form of the filtrations \eqref{Wxiphi-filt-intro}. Those filtrations could be viewed as a sort of secondary VHS. It comes with its own secondary differential which turn our zero energy vectors of $W_{\xi}$ into dynamical entities of their own: each nonzero vector of $W_{\xi}$ gives rise to a linear map which acts on $W_{\xi}$. Our refinement could be envisaged as a new form of `energy': a nonzero vector in $ W_{\xi}$ which has zero `energy' with respect to $\xi$ in the classical cup product \eqref{diff-period-intro} now becomes a linear operator whose `energy' is measured by

- the length of the secondary filtration \eqref{Wxiphi-filt-intro},

- the dimension of the last piece of that filtration.

All this leads to uncovering `hidden' structures of the filtrations \eqref{Wxiphi-filt-intro}  such as quivers, Fano toric varieties, dimer models. This in turn provides intriguing connections with such topics as quantum-type invariants, elliptic curves, number theory, mirror symmetry. Perhaps the unifying theme of all those connections is the {\it higher categorical structures} of the IVHS. The paper does not attempt to address this possibility formally, however, various aspects of the refinement presented here strongly point in this direction. The main objective of this writing is to explain

- the constructions of the refinement of IVHS,

- its connections with classical topics in the curve theory,

- how all of the above mentioned `hidden' structures arise,

- indicate how some of the topics listed above connect with the refinement of IVHS. 

\vspace{0.2cm}
In the rest of the introduction we give a brief overview of the results and constructions of the monograph and sketch the logical connections between them. Hopefully, this will provide a guide for
the reader to browse through the sections of the paper.

\vspace{0.2cm}
{\bf 1. Algebraic K\"ahler structure of $W_{\xi}$ and the quiver attached to it.}
The filtrations of $W_{\xi}$ in \eqref{Wxiphi-filt-intro} attached to points
$([\xi],[\phi])$ in the cohomological incidence correspondence ${\bf P}$ in \eqref{two-inc-diag} arise from a certain linear map
\begin{equation}\label{alphaxi-intro}
\mbox{$\alpha^{(2)}_{\xi}: \bigwedge^2 W_{\xi} \longrightarrow \HKC$}
\end{equation}
naturally attached to the pair $([\xi],W_{\xi})$, a two-form with values in $\HKC$.
This turns out to be only part of the structure. Using the Hodge metric on $\HKC$ we can break $W_{\xi}$ into orthogonal pieces, that is, there is an intrinsically defined orthogonal decomposition
\begin{equation}\label{Wxi-orthdec-intro}
	W_{\xi}=\bigoplus^{l_{\xi}([\phi])}_{s=0} P^s ([\xi],[\phi]),
\end{equation}
where $l_{\xi}([\phi])$ is the length of the filtration in \eqref{Wxiphi-filt-intro}. This decomposition is related to the filtration by the formula
$$
W^s_{\xi} ([\phi])=\bigoplus^{l_{\xi}([\phi])}_{t=s} P^t ([\xi],[\phi])
$$
which is analogous to the relation between Hodge filtration and Hodge decomposition in the classical VHS. In addition to the decomposition, the linear map \eqref{alphaxi-intro} gives rise to the linear map
$$
\alpha^{(2)}_{\xi} (\phi, \bullet): W_{\xi} \longrightarrow \HKC.
$$
The last summand $P^{l_{\xi}([\phi])} ([\xi],[\phi])$ of the decomposition turns out to be $\alpha^{(2)}_{\xi} (\phi, \bullet)$-invariant. More precisely, we have the endomorphism
\begin{equation}\label{invpart-intro}
\widehat{\alpha^{(2)}}_{\xi} (\phi, \bullet): P^{l_{\xi}([\phi])} ([\xi],[\phi])/\CC\phi \longrightarrow P^{l_{\xi}([\phi])} ([\xi],[\phi])/\CC\phi.
\end{equation} 
 The study of this endomorphism underlies the assertion about the linear system in Theorem \ref{pro:refine}: we establish a precise dictionary between the eigen spaces of the endomorphism and certain special linear systems on $C$
 \begin{equation}\label{spectrum-intro}
 \text{\{\it the spectrum of $\widehat{\alpha^{(2)}}_{\xi} (\phi, \bullet)$\}} \longrightarrow \text{\{\it special line bundles on $C$\}}.
\end{equation}
 For details and precise statement see Proposition \ref{pro:xi-Wlatleast2-eigen}. Thus the last piece of our filtration establishes links, on the one hand, to classical problems in curve theory such as special divisors, and on the other hand, to Brill-Noether loci in the moduli space of curves.
 
 We now turn to other summands of the decomposition \eqref{Wxi-orthdec-intro}. Every summand $P^s([\xi],[\phi])$, for $s\in [0,l_{\xi}([\phi])-1]$, is sent by $\alpha^{(2)}_{\xi} (\phi, \bullet)$ to the direct sum
$$
\bigoplus^{l_{\xi}([\phi])}_{t=s-1} P^t ([\xi],[\phi]),
$$
and where the summand $P^{-1}([\xi],[\phi])$ is understood as the orthogonal complement $W^{\perp}_{\xi}$ of $W_{\xi}$ in $\HKC$ with respect to the Hodge metric. Observe that this means that the step $W^s_{\xi} ([\phi])$ of the filtration is sent to the step $W^{s-1}_{\xi} ([\phi])$, for every $s\geq 1$, a sort of Griffiths transversality property in the classical variation of Hodge structure.

For the integers $t,s$ in $[0,l_{\xi}([\phi])-1]$, denote by $\alpha^{t,s}$ the block of $\alpha^{(2)}_{\xi} (\phi, \bullet)$ mapping the summand $P^s ([\xi],[\phi])$ to $P^t ([\xi],[\phi])$. We encode part of these data in the following graph: 

$\bullet$ for every summand $P^s([\xi],[\phi])$, $s$ in the range of $[0,l_{\xi}([\phi])-1]$ of the orthogonal decomposition of $W_{\xi}$, we draw two vertices on the same vertical level; color the top one white $\circ$, and the bottom one black $\bullet$; place the labels $s$ at the white (top) vertex and the label $s'$ for the black (bottom) in the decreasing order from left to right  for $s\in [0,l_{\xi}([\phi])-1]$; 

$\bullet$ the edges of the graph are drawn between the vertices of different colors only; this is done according to the rule: the top vertex $s$ is connected to the bottom $s'$ and its immediate neighbors $(s-1)'$ and $(s+1)'$.

The resulting graph is denoted $PG_{l_{\xi}([\phi])}$. The following drawing illustrates the graph $PG_{5}$, for $l_{\xi}([\phi])=5$:

$$
\begin{tikzpicture}
	[place/.style={circle,draw=black,thick},
	transition/.style={circle,draw=black,fill=black}]
	\node (white4) at (0,2) [place] [label={above:$4$}] {};
	\node (black4) at (0,0) [transition] [label={below:$4'$}]{};
	\node (white3) at (2,2) [place] [label={above:$3$}]{};
	\node (black3) at (2,0) [transition] [label={below:$3'$}]{};
	\node (white2) at (4,2) [place] [label={above:$2$}]{};
	\node (black2) at (4,0) [transition] [label={below:$2'$}]{};
	\node (white1) at (6,2) [place] [label={above:$1$}] {};
	\node (black1) at (6,0) [transition] [label={below:$1'$}]{};
	\node (white0) at (8,2) [place] [label={above:$0$}] {};
	\node (black0) at (8,0) [transition] [label={below:$0'$}]{};
	\begin{scope}[ultra thick]
		\draw[red] (white1) to (black1);
		\draw[red] (white0) to (black0);
		\draw[red] (white2) to (black2);
		\draw[red] (white3) to (black3);
		\draw[red] (white4) to (black4);
	\end{scope}
	\begin{scope}[thick]
		\draw[blue][-] (white4) to (black3);
		\draw[blue][-] (white3) to (black2);	
		\draw[blue][-] (white2) to (black1);
		\draw[blue][-] (white1) to (black0);
	\end{scope}
	\draw[red, dotted,thick] (white3) to (black4);
	\draw[red, dotted,thick] (white2) to (black3);
	\draw[red,dotted,thick] (white1) to (black2);
	\draw[red,dotted,thick] (white0) to (black1);
\end{tikzpicture}
$$

By orienting the edges from `white' to `black' we think of $PG_{l_{\xi}([\phi])}$ as a quiver. It should be clear, that the graph is constructed so that the orthogonal decomposition of $W_{\xi}$ in \eqref{Wxi-orthdec-intro} together with the maps $\alpha^{t,s}$ provide a representation of the quiver $PG_{l_{\xi}([\phi])}$. Namely, we place the summand $P^s ([\xi],[\phi])$ at the vertices $s$ and $s'$, for every $s$ in the range $[0,l_{\xi}([\phi])-1]$; the vertical edges of the graph are decorated with the endomorphisms
$$
\alpha^{s,s}: P^s ([\xi],[\phi]) \longrightarrow P^s ([\xi],[\phi]),
$$
the edges $(s)\to (s-1)'$ are decorated with the maps
$$
\alpha^{s-1,s}: P^s ([\xi],[\phi]) \longrightarrow P^{s-1} ([\xi],[\phi]), 
$$
for every $s\in [1, l_{\xi}([\phi])-1]$, and the edges $s \to (s+1)'$
with the maps
$$
\alpha^{s+1,s}: P^s ([\xi],[\phi]) \longrightarrow P^{s+1} ([\xi],[\phi]), 
$$
for every $s\in [0, l_{\xi}([\phi])-2]$.

The following picture emerges:

- the filtrations attached to points of the cohomological
incidence correspondence ${\bf P}$ in \eqref{two-inc-diag} equip it with the length function
$$
{\mathfrak{l}}: {\bf P} \longrightarrow \NN
$$
sending a point $([\xi],[\phi])$ to the length $l_{\xi}([\phi])$ of the filtration of $W_{\xi}$ in \eqref{Wxiphi-filt-intro}; this is a constructible function and it stratifies the incidence correspondence $p^{-1}_1 (\Sigma^0_r)$ overlying each stratum $\Sigma^0_r$ of Griffiths stratification into the finer strata according to the values of the length function ${\mathfrak{l}}$;

- 
 fix a value $l$ in the range of the length function and consider the stratum
$$
{\mathfrak{L}}_l:={\mathfrak{l}}^{-1}(l);
$$
for the rest of the discussion we assume that $l$ is at least $3$; from the above it follows:

$\bullet$ the stratum ${\mathfrak{L}}_l$ comes along with the
graph $PG_l$, that graph is bipartite - its vertices are colored white/black and edges connect vertices of different colors only; with the orientation of edges `white'-to-`black' $PG_l$ is a quiver; 

$\bullet$ the points of ${\mathfrak{L}}_l$ provide quiver representations of  $PG_l$ according to the rules discussed above; more precisely, for every complex number $c$ there is a distinguished family of representations
\begin{equation}\label{rhoc-intro}
\rho_c: \OO^{\times}_{{\mathfrak{L}}_l}(-1) \longrightarrow {\mathfrak{Reps}}(PG_l)
\end{equation}
parametrized by $\OO^{\times}_{{\mathfrak{L}}_l}(-1)$, the total space of the tautological line bundle $\OO_{\bf P}(-1)$ restricted to ${\mathfrak{L}}_l$ and the zero section removed;

$\bullet$ the stratum ${\mathfrak{L}}_l$ can be decomposed further into the substrata where the dimensions of the graded pieces of the filtration $W^{\bullet}_{\xi}([\phi])$ remain constant; each such substratum is denoted in the main body of the paper by ${\mathfrak{L}}_l(h^l,\lambda)$: 
$h^l$ is the dimension of the last piece of the filtration and $\lambda$ is
the partition associated to the graded module
$$
W_{\xi}/W^l_{\xi}[\phi]=\bigoplus^{l-1}_{s=1} P^s([\xi],[\phi]).
$$

The length strata ${\mathfrak{L}}_l$ and their substrata ${\mathfrak{L}}_l(h^l,\lambda)$ provide a geometric part of our refinement of IVHS; the quiver $PG_l$ and $\OO^{\times}_{{\mathfrak{L}}_l}(-1)$, viewed as a parameter space for representations of $PG_l$, constitute an algebraic part of the refinement.

\vspace{0.2cm}
{\bf 2. The graph underlying the quiver and the associated Fano toric variety.}
The graph $PG_l$ is used to introduce another geometric object of our refinement: a Fano toric variety. It can be done  formally starting from $PG_l$ as follows: assign the indeterminate $T$ to all vertical edges of the graph, the indeterminates $X_r$ to the edges $(r)\to (r-1)'$,  and $Y_r$ to the edges $(r-1) \to (r)'$, for $r\in [1,l-1]$; we think of
$$
\{T,X_r, Y_r \,|\, r\in [1,l-1]\}
$$
as linear forms on the vector space $V_{2l-1}$ freely generated by the vertices $\{(r), (r)' |\,r\in [1,l-1] \}$ and the symbol $\{e^0\}$ which stands for the vertical edge $(0)\to (0)'$; the form $T$ vanishes on all vertices and $T(e^0)=1$, the forms $X_r$'s (resp., $Y_r$'s) form the dual basis for the white (resp., black) vertices and vanishes on the black (resp., white) ones.
Next we write down the system of quadratic equations
$$
X_r Y_r=T^2, \,r=1,\ldots,l-1;
$$
this gives a system of quadrics in $\PP(V_{2l-1})\cong \PP^{2(l-1)}$ and we define our variety as the complete intersection of those quadrics; this variety is denoted ${\bf H^{1,0}}(PG_l)$ in the main body of the article. The following statement summarizes its properties.
\begin{pro}\label{pro:NADolb-intro}
	The variety ${\bf H^{1,0}}(PG_{l})$ is a singular Fano toric variety in the projective space $\PP(V_{2l-1})=\PP^{2(l-1)}.$ It has dimension $(l-1)$ and degree $2^{l-1}$.
	
	The hyperplane $T=0$ intersects ${\bf H^{1,0}}(PG_{l})$ along the divisor denoted $H_0$. This divisor is the union
	of $2^{l-1}$ projective subspaces. More precisely, for every subset
	$A \subset [1,l-1]$ denote by $A^c$ its complement in $[1,l-1]$ and let $\Pi_A$ be the projective subspace of the hyperplane
	$T=0$ spanned by the points of $\PP(V_{2l-1})$ underlying the vectors
	$$
	\{(s),(t)' | s\in A,\, t\in A^c\}  \subset V_{2l-1};
	$$
	then we have
	$$
	H_0 =\bigcup_{A\subset [1,l-1]} \Pi_A,
	$$
	where the union is taken over all subsets $A$ of $[1,l-1]$.
\end{pro}

The hyperplane sections of ${\bf H^{1,0}}(PG_{l})$ are (singular) Calabi-Yau varieties. The divisor $H_0$ in the proposition is a particular degenerate case. For the reasons explained in the article, it is called Lagrangian cycle and its irreducible components, the projective spaces $\Pi_A$ in the proposition, could be considered as Lagrangians. We suggest that these should be objects of a Fukaya-type category. This is a subject for a future research. For now we wish to explain the relation of the variety
${\bf H^{1,0}}(PG_{l})$ to the quiver representation side of the graph $PG_l$. 

So far the variety ${\bf H^{1,0}}(PG_{l})$ and its Lagrangian cycle $H_0$ are defined formally and their relevance to algebraic K\"ahler structures is not clear. So the appearance of this Fano variety may strike the reader as completely artificial. This is not so, since the quadratic equations defining ${\bf H^{1,0}}(PG_{l})$ arise from natural deformations of the quiver representations of $PG_l$. Let us briefly explain this. 

For this we return to the middle of the incidence diagram \eqref{two-inc-diag}, the projective space $\PP(\HKC)$. If we start with a nonzero subspace $W$  in $\HKC$, then the cohomological incidence correspondence $P$ produces the subspace of
Kodaira-Spencer classes annihilating $W$:
$$
\Xi_W:=\{\xi \in H^1 (\Theta_{C}) | \xi \phi =0, \forall \phi \in W\}.
$$
We assume that the dimension of $W$ is at least two and $\Xi_W$ is nonzero. Our algebraic K\"ahler structure gives linear maps
$$
\alpha^{(2)}_{\xi}(\phi,\bullet): W_{\xi} \longrightarrow \HKC.
$$
By construction $W$ is contained in every $W_{\xi}$ as $\xi$ varies in $\Xi_W$. Restricting the maps  $\alpha^{(2)}_{\xi}(\phi,\bullet)$ to $W$, we obtain the map
$$
W \longrightarrow (\Xi_W)^{\ast} \otimes Hom(W,\HKC)
$$
The main point is  that this is a sort of Higgs field. More concretely, the operators $\alpha^{(2)}_{\xi}(\phi,\bullet)$ are commuting in the following sense: for $\xi \neq \xi'$ in $\Xi_W$ we have the equality
\begin{equation}\label{Higgs-intro}
	\alpha^{(2)}_{\xi}(\phi,\alpha^{(2)}_{\xi'}(\phi,\psi)) \equiv \alpha^{(2)}_{\xi'}(\phi,\alpha^{(2)}_{\xi}(\phi,\psi)) \mod \CC\phi,
\end{equation}
whenever the compositions are defined. This in turn produces commutation relations between the blocks $\{\alpha^{m,k}_{\xi} (\phi) \}$ and $\{\alpha^{m,k}_{\xi'} (\phi) \}$ of $\alpha^{(2)}_{\xi}(\phi,\bullet)$ and $\alpha^{(2)}_{\xi'}(\phi,\bullet)$. Now we deform the operators
$\alpha^{(2)}_{\xi} (\phi, \bullet)$ and $\alpha^{(2)}_{\xi'} (\phi, \bullet)$
$$
D_{\xi} (t,{\bf x}, {\bf y}) =t\sum^{l-1}_{k=0}  \alpha^{k,k}_{\xi} (\phi) + \sum^{l-1}_{k=1} x_k \alpha^{k-1,k}_{\xi}(\phi) + \sum^{l-2}_{k=0} y_{k+1}\alpha^{k+1,k}_{\xi}(\phi) + \sum_{m-k\geq 2}\alpha^{m,k}_{\xi} (\phi),
$$
by introducing deformation parameters ${\bf x}=(x_k)$ and ${\bf y}=(y_k)$  in $\CC^{l-1}$ and $t\in \CC$, for blocks decorating the edges of the graph $PG_l$. We wish to deform so that the commutativity condition \eqref{Higgs-intro} would be conserved. The direct computation gives the 
system of quadratic equations in the formal definition of
${\bf H^{1,0}}(PG_{l})$. We also gain a conceptual understanding: ${\bf H^{1,0}}(PG_{l})$ is a parameter space of Higgs fields. This is the reason for the notation alluding to $(1,0)$-forms as well as the name used in the main body of the paper - nonabelian Dolbeault variety.

Though the defining equations of ${\bf H^{1,0}}(PG_{l})$ come from the deformations of quiver representations arising from our refinement, the variety itself depends on the graph only. It completely forgets our curve $C$. The situation is reminiscent of the classical Jacobian of a curve which takes into account only cohomological data of a curve. The connection with the curve is restored through the Abel-Jacobi maps. Following the analogy with the classical situation, we construct several naive versions
of Abel-Jacobi maps:

- the maps from $\OO^{\times}_{{\mathfrak L}_l}(-1)$ to ${\bf H^{1,0}}(PG_{l})\setminus H_0$, the complement of the Lagrangian cycle;
$\OO^{\times}_{{\mathfrak L}_l}(-1)$ is as in \eqref{rhoc-intro} - it stands for the total space of the tautological bundle $\OO_{{\mathfrak L}_l}(-1)$ with the zero section removed;

- the dual of the above assigning to the points of $\OO^{\times}_{{\mathfrak L}_l}(-1)$ hyperplane sections of ${\bf H^{1,0}}(PG_{l})$,

- for a nonzero section $\phi \in \HKC$ we distinguish the stratum
${\mathfrak L}_l (\phi)$ of ${\mathfrak L}_l$ together with the map
of $\OO^{\times}_{{\mathfrak L}_l(\phi)}(-1)$ into the ring $\CC(C)[u,u^{-1}]$ of Laurent polynomials with coefficients
in the function field $\CC(C)$ of $C$.

The constructions are given in the subsections \S6.3 and \S6.4. The first map is an analogue of the classical Abel-Jacobi map. Here we wish to point out the meaning of the second and the third maps. The reader may recall that the hyperplane sections of ${\bf H^{1,0}}(PG_{l})$ are singular Calabi-Yau varieties, so the second map connects the points of ${\mathfrak L}_l$ with Calabi-Yau varieties, that is, the tautological fibration $\OO^{\times}_{{\mathfrak L}_l}(-1)$ over ${\mathfrak L}_l$ becomes a parameter space for families of Calabi-Yau varieties. This could be a source of many interesting invariants supported on ${\mathfrak L}_l$.

For the third map, the reader will need to recall the diagram in \eqref{two-inc-diag} connecting the cohomological and geometrical incidence correspondences. Conceptually, the third map establishes a connection between the two parts of the diagram: a point $[\phi]$ in the middle space $\PP(\HKC)$ of the diagram defines the stratum ${\mathfrak L}_l(\phi)$ in the cohomological incidence $P$ and relates it to the function field of $C$ on the geometrical part of the diagram. The `strange' appearance of Laurent polynomials will become more natural and interesting once we discuss the next aspect of the refinement.

\vspace{0.2cm}
{\bf 3. From algebraic K\"ahler structure of $W_{\xi}$ to dimer models.}
We already said that the graphs $PG_l$ associated to the strata ${\mathfrak{L}}_l$ are bipartite, that is, the vertices come in two different colors (white and black) and each edge of the graph connects vertices of different colors. The graphs of this type are called dimer models and have been intensively studied in recent years in connection with mirror symmetry. There is a vast physics and mathematics literature on the subject. We cite \cite{Go-K} and references therein for more details on the subject.

The graph $PG_l$ is a dimer model with the special property that it is `almost' trivalent: all vertices have valance $3$ except the first and the last pair, that is the vertices with labels $(0),(0)'$ and $(l-1), (l-1)'$.
We turn $PG_l$ into a {\it trivalent dimer model} by connecting the vertex $(0)$ to $(l-1)'$ and the vertex $(l-1)$ to $(0)'$, in accordance to the pattern of edges at the intermediate vertices. The graph obtained this way is denoted $\widehat{PG}_l$. Below is the illustration of the graph $\widehat{PG}_4$.

$$
\begin{tikzpicture}
	[place/.style={circle,draw=black,thick},
	transition/.style={circle,draw=black,fill=black}]
	\node (white3) at (0,2) [place] [label={above:$3$}] {};
	\node (black3) at (0,0) [transition] [label={below:$3'$}]{};
	\node (white2) at (2,2) [place] [label={above:$2$}]{};
	\node (black2) at (2,0) [transition] [label={below:$2'$}]{};
	\node (white1) at (4,2) [place] [label={above:$1$}]{};
	\node (black1) at (4,0) [transition] [label={below:$1'$}]{};
	\node (white0) at (6,2) [place] [label={above:$0$}] {};
	\node (black0) at (6,0) [transition] [label={below:$0'$}]{};
	\draw[red,ultra thick] (white0) to (black0);
	\draw[red,ultra thick] (white1) to (black1);
	\draw[red,ultra thick] (white2) to (black2);
	\draw[red,ultra thick] (white3) to (black3);
	\begin{scope}[thick]
		\draw[blue][-] (white3) to (black2);	
		\draw[blue][-] (white2) to (black1);
		\draw[blue][-] (white1) to (black0);
	\end{scope}
	\draw[red, dotted,thick][-] (white2) to (black3);
	\draw[red,dotted,thick][-] (white1) to (black2);
	\draw[red,dotted,thick][-] (white0) to (black1);
	\draw[blue,thick][-] (white0) .. controls (10,-0.5) and (4,-1)  .. (3,-1) .. controls (-1,-1) and (-1,1).. (black3);
	\draw[red,dotted,thick][-] (white3) .. controls (-4,-0.5) and (2.5,-2.1)  .. (3,-2) .. controls (8,-1.5) and (7,2).. (black0);
\end{tikzpicture}
$$

The graph $\widehat{PG}_l$ is endowed with a structure of a ribbon graph, that is, we choose a cyclic order of edges at every vertex: our choice is counterclockwise, starting with the vertical edge. This ribbon graph structure turns $\widehat{PG}_l$ into a retract of a topological torus $\mathbb{T}$ with $l$ disks removed; the boundary of the disks are so called boundary cycles $B_0, B_1, \ldots, B_{l-1}$ of the ribbon graph $\widehat{PG}_l$; they are labeled by the white vertices. We obtain what in physics literature is called `brane tilings', the `tiles' in reference to the polygons formed by the boundary cycles, see \cite{Ken} for physics oriented overview of the subject. In the situation of $\widehat{PG}_l$ all  tiles happen to be hexagons; all this is discussed and worked out in \S7. 

While the graph $PG_l$ underlies the quiver representations of our IVHS refinement, the completed graph seems to lack this connection. In fact,
one of the main features of $\widehat{PG}_l$ is the following extension property.
\begin{thm}\label{pro:repsext-intro}
	Given a bipartite finite dimensional representation $\rho=\{\alpha^{t,s}:P^s \longrightarrow P^t\}$\footnote{{\it bipartite} means that the {\it same} vector space $P^s$ is assigned to  pairs $(s) -(s')$ of white-black vertices.}
	 of the quiver $PG_l$, it can be extended in a canonical way to the representation of the quiver $\widehat{PG}_l$. Those canonical extensions are parametrized by pairs of linear functionals on the space of Laurent polynomials $\CC[q,q^{-1}]$.
\end{thm}
The details and consequences of this statement are given in \S11 of the paper. Let us discuss the result and its proof informally. 

What has been done so far is a finer stratification of the incidence correspondence $\PP({\cal W}_{\Sigma^0_r})$ over $\Sigma^0_r$: it is divided into the strata ${\mathfrak{L}}_l$ where the length of filtrations is fixed and equal to $l$. The graph $PG_l$ appears as a bookkeeping device
for arranging those filtrations as quiver representations of $PG_l$. The graph could be viewed as `open' in a sense that its quiver representation 
$$
\rho=\{\alpha^{t,s}:P^s \longrightarrow P^t\}
$$
 has `ends', the spaces $P^0$ and $P^{l-1}$; those could be viewed as `output' {\it objects} of the quiver representations associated to our filtrations; the obvious {\it morphisms} are the vector spaces
 $$
 \begin{gathered}
 P^-:=Hom_{\CC} (P^0,P^{l-1}), \\ P^+:=Hom_{\CC} (P^{l-1},P^{0})
 \end{gathered}
 $$
 of linear maps between the end objects.
 
 The graph $\widehat{PG}_l$ connects the ends of the graph $PG_l$ combinatorially/topologically by introducing the additional edges allowing to go {\it directly} between the vertices $(0)$ (resp. $(0)'$) and $(l-1)'$ (resp. $(l-1)$). The proposition tells us that we can consistently label those edges using the data of the representations of the `open' graph and some `quantum' parameters, the space
$Hom_{\CC} (\CC[q,q^{-1}],\CC)$. Concretely, the proof of the proposition tells us that there are distinguished linear maps
\begin{equation}\label{tau-maps-intro}
\tau^{\pm}_{\rho}:	Hom_{\CC} (\CC[q,q^{-1}],\CC) \longrightarrow P^{\pm}
\end{equation}
which allow to connect the `end' objects of the quiver representation $\rho$ of $PG_l$. Observe the appearance of Laurent polynomials: here they appear naturally as traces of loop versions of certain Lie algebras attached to the quiver representations. The technical part of the construction of $\tau^{\pm}_{\rho}$ involves certain closed paths on $\widehat{PG}_l$ called zig-zag paths, see \S11 for details.

Thus the above maps is an extra structure on the {\it morphisms} between the end objects of the representations of the quiver $PG_l$. One could wonder if we have morphisms between morphisms and if those have additional structures. In other words, we are asking for {\it higher categorical structures of the IVHS}. In the main body of the paper we provide indications that this might be the case. There are two different directions: one connects to symplectic geometry and the other to the ideas of the Topological Quantum field theory (TQFT). The first stems from the fact that the space
$$
P=P^- \oplus P^+
$$
combining the morphisms between the end objects is naturally a symplectic vector space: this is due to the fact that the two summands are dual to each other. In the main body of the paper we construct Lagrangian subspaces in
$P$ which could be viewed as $1$-morphism; symplectic morphisms between those Lagrangians could be $2$-morphisms and a possible connection to Fukaya categories arises.

The second direction uses the fact that the graph $\widehat{PG}_l$ is a dimer model and as such it comes with the set of perfect matchings - configurations of edges which cover all vertices of the graph and precisely once, see \S10 for details. We suggest that this set can be categorified and this could provide higher order
`multiplications' in the tensor algebra of the graded space $W_{\xi}/W^l_{\xi}([\phi])$, see \S11.3 for details.

 The two directions just outlined will need further elaboration, but a concrete output of Theorem \ref{pro:repsext-intro} seems to be worth  stressing:
 \begin{equation}\label{repsext-intro}
 	\begin{gathered}
 	\text{\it the refinement mechanism  provides each point of the strata ${\mathfrak{L}}_l$ lying over}
 	\\
 	\text{\it the stratum $\Sigma^{\circ}_r$ with {\rm the moduli} of representations of the extended quiver $\widehat{PG}_l$;}
 	\\
 	\text{\it the moduli is}
 	\\
 	Hom(\CC[q,q^{-1}],\CC)\times Hom(\CC[q,q^{-1}],\CC), 
 	\\
 \text{\it the space of pairs of linear functionals on the algebra of polynomials $\CC[q,q^{-1}]$.}
\end{gathered}
\end{equation}

The way those extensions arise is also meaningful: the graph $\widehat{PG}_l$ comes equipped with a distinguished closed paths, known as zig-zag paths, see Lemma \ref{lem:zzpaths}, and the maps \eqref{tau-maps-intro} are obtained by a sort of integration over those paths. The zig-zag paths are also closely linked with perfect matchings on $\widehat{PG}_l$ whose role in the theory has been mentioned earlier. There is a geometric realization of zig-zag paths as a certain toric surface denoted $X(\Delta_l)$ and described in \S8. Its fan $\Delta_l$ is determined by the cohomology classes of zig-zag paths, as closed paths on the topological torus $\mathbb{T}$. The surface comes along with a distinguished family of curves determined by the Kasteleyn determinant of the graph $\widehat{PG}_l$ and known to be the {\it generating function of matchings of the graph}. Those curves, call them {\it Kasteleyn curves of $X(\Delta_l)$}, were studied in \cite{KeOkSh} in the general setting of bipartite graphs. In our situation we calculate explicitly the equations of those curves in $\CC^{\times}\times\CC^{\times}$, see Proposition \ref{pro:detKastelyan} and Remark \ref{rem:Kastelyanl}. This algebraic torus is viewed as the open orbit of the torus action on $X(\Delta_l)$ and provides a link with
the extensions of quiver representations in \eqref{repsext-intro}: the space of linear functionals $Hom(\CC[q,q^{-1}],\CC)$ naturally contains
$\CC^{\times}$ as functionals of evaluation at a point of $\CC^{\times}$.  So we have the natural inclusion
$$
\CC^{\times}\times\CC^{\times} \subset  	Hom(\CC[q,q^{-1}],\CC)\times Hom(\CC[q,q^{-1}],\CC)
$$
providing a {\it geometric moduli} for the extensions of representations of
$PG_l$ to the ones of $\widehat{PG}_l$. With this the Kasteleyn curves   
come along with representations of the quiver $\widehat{PG}_l$.

\vspace{0.2cm}
The extension of the representations of $PG_l$ to the ones of the quiver
$\widehat{PG}_l$ in \eqref{repsext-intro} has a geometric counterpart: the strata ${\mathfrak L}_l$
become related to the moduli spaces ${\mathfrak M}^l_1$ of elliptic curves with
$l$ marked points. This relation comes from metrics on the ribbon graph $\widehat{PG}_l$, that is positive valued functions on the set of edges $E_{\widehat{PG}_l}$ of $\widehat{PG}_l$. The general theory relates metrized ribbon graphs with the moduli spaces of curves with marked points, the theory which goes back to the works of Harer, Mumford, Penner, Thurston, see \cite{Har}, \cite{MuP} for overviews. In our case we obtain the map
$$
\RR^{E_{\widehat{PG}_l}}_{+} \longrightarrow {\mathfrak{M}}^l_1
$$
where ${\mathfrak{M}}^l_1$ is the moduli space of curves of genus one with $l$ marked points. The main point is that  a finite dimensional {\it Hermitian} quiver representation\footnote{this means that the vector spaces of a representation are equipped with Hermitian metrics.} of $\widehat{PG}_l$ {\it determines} a metric on
$\widehat{PG}_l$ and this metric turns the topological torus ${\mathbb{T}}$ into a genus one curve with $l$ marked points, the barycenters of the boundary cycles of $\widehat{PG}_l$:
$$
{\mathfrak{Reps}}^{Herm}(\widehat{PG}_l) \longrightarrow \RR^{E_{\widehat{PG}_l}}_{+} \longrightarrow {\mathfrak{M}}^l_1.
$$
 Combining this with the extension property of $\widehat{PG}_l$ in Theorem \ref{pro:repsext-intro} gives the following.  

\begin{thm}\label{pro:repsext-intro1}
	Let $\rho=\{\alpha^{t,s}:P^s \longrightarrow P^t\}$ be a bipartite, finite dimensional representation of the quiver $PG_l$, it can be extended in a canonical way to the representation of the quiver $\widehat{PG}_l$. Those canonical extensions are parametrized by pairs of linear functionals on the algebra of Laurent polynomials $\CC[q,q^{-1}]$. If, in addition, a representation $\rho$ is assumed to be Hermitian, that is the vector spaces $P^s$ of $\rho$ are equipped with Hermitian metrics, then each extension of $\rho$ to a representation of $\widehat{PG}_l$ determines a metric on $\widehat{PG}_l$ and hence a point in the moduli space ${\mathfrak M}^l_1$, that is, for every pair $(F,G)$ of linear functionals on $\CC[q,q^{-1}]$ one obtains a map
	$$
	\wp_{l} (F,G): {\mathfrak{Reps}}^{Herm}(PG_l) \longrightarrow {\mathfrak M}^l_1
	$$
	from the space ${\mathfrak{Reps}}^{Herm}(PG_l)$ of bipartite, finite dimensional Hermitian representations of $PG_l$ to the moduli space ${\mathfrak M}^l_1$ of curves of genus one with $l$ marked points.
	
	In particular, to the stratum ${\mathfrak L}_l$ are attached the families of representations of $PG_l$
	$$
	\rho_c: \OO^{\times}_{{\mathfrak{L}}_l}(-1) \longrightarrow {\mathfrak{Reps}}(PG_l)
	$$
	see \eqref{rhoc-intro}; those representations are Hermitian and give rise to continuous maps
	$$
	\wp_{l,\rho_c}: \OO^{\times}_{{\mathfrak{L}}_l}(-1)\times\left(Hom_{\CC} (\CC[q,q^{-1}],\CC)\right)^2 \longrightarrow {\mathfrak M}^l_1 .
	$$
	which send a point $(\xi\otimes\phi,(F,G))$ in the domain to the point	$\wp_{l} (F,G)(\rho_c(\xi\otimes\phi))$ in the moduli space ${\mathfrak M}^l_1$ of curves of genus one with $l$ marked points.
\end{thm}

 One can also phrase the above result as attaching to each point of ${\mathfrak L}_l$ {\it families of elliptic curves with $l$ marked points,} where the families are parametrized by the the fibre of the bundle $\OO^{\times}_{{\mathfrak{L}}_l}(-1)$ over the point and the space $\left(Hom_{\CC} (\CC[q,q^{-1}],\CC)\right)^2$. In other words, every point $([\xi],[\phi])$ of ${{\mathfrak L}_l}$ acquires {\it geometric moduli}
 $$
 \wp_{l,\rho_c}|_{([\xi],[\phi])}:\OO^{\times}_{{\mathfrak{L}}_l,([\xi],[\phi])}(-1)\times\left(Hom_{\CC} (\CC[q,q^{-1}],\CC)\right)^2 \longrightarrow {\mathfrak M}^l_1,
 $$ 
 where $\OO^{\times}_{{\mathfrak{L}}_l,([\xi],[\phi])}(-1)$ is the fibre of $\OO^{\times}_{{\mathfrak{L}}_l}(-1)$ over $([\xi],[\phi])$. Remembering the inclusion 
 $$
 \CC^{\times} \times\CC^{\times} \subset \left(Hom_{\CC} (\CC[q,q^{-1}],\CC)\right)^2
 $$
 we obtain maps
 $$
  {\wp^{X(\Delta_l)}_{l,\rho_c}}|_{([\xi],[\phi])}:\CC^{\times}\times\CC^{\times}\times \CC^{\times}  \longrightarrow {\mathfrak M}^l_1,
  $$ 
  where the first factor in the domain is an identification $\OO^{\times}_{{\mathfrak{L}}_l,([\xi],[\phi])}(-1)\cong\CC^{\times}$; the superscript refers to the toric surface $X(\Delta_l)$ defined by zig-zag paths and the remaining product $\CC^{\times}\times \CC^{\times}$ is thought as the open orbit of that surface under the torus action. We can restate \eqref{repsext-intro} in geometric terms
  \begin{equation}\label{repsext-geom-intro}
  	\begin{gathered}
  		\text{\it the refinement mechanism  provides each point $([\xi],[\phi])$ of the strata ${\mathfrak{L}}_l$ lying over}
  		\\
  		\text{\it the stratum $\Sigma^{\circ}_r$ with  the {\it geometric moduli}, the continuous maps}
  		\\
  		 \wp^{X(\Delta_{l})}_{l,\rho_c}|_{([\xi],[\phi])}:\CC^{\times}\times\CC^{\times}\times \CC^{\times}  \longrightarrow {\mathfrak M}^l_1. 
  	\end{gathered}
  \end{equation}

 The graph $\widehat{PG}_l$ has another particular feature: it is endowed with a distinguished automorphism of order $l$. This automorphism is denoted by $\sigma_l$. By choosing a metric on $\widehat{PG}_l$ invariant with respect to $\sigma_l$ allows us to be more precise about the marked points: if $l$ is at least four, the marked points form a cyclic subgroup of order $l$ in the corresponding elliptic curve and the subgroup comes with a particular choice of a generator. In other words, the $\sigma_l$-invariant metrics on $\widehat{PG}_l$ take us to the modular curve $Y_1(l)$, the moduli of pairs $(E,p)$ where $E$ is an elliptic curve  and $p$ is a point of order $l$ in $E$, see Proposition \ref{pro:torsion}. Actually, any metric of $\widehat{PG}_l$ can be made $\sigma_l$-invariant by choosing a triple of
 standard symmetric polynomials in $l$ indeterminates. It is well-known that the standard symmetric functions, such as monomial, elementary, power functions are parametrized by partitions, see \cite{Mac}. Thus once we specify which symmetric functions we associate to partitions, the $\sigma_l$-invariant metrics of $\widehat{PG}_l$ can be obtained by starting from {\it any} metric on
 $\widehat{PG}_l$ and then making it $\sigma_l$-invariant according to a choice of a triple of partitions having at least $l$ parts.

 \begin{pro}\label{pro:torsion-intro}
 	Assume $l\geq 4$. For any triple of partitions 
 	$$
 	{\bf \mu}=(\mu_0,\mu_{+}, \mu_{-})
 	$$
 	with at least $l$ parts, the map  $\wp_{l,\rho_c}$ in Theorem \ref{pro:repsext-intro1} gives rise to the continuous map
 	$$
 	\wp_{l,\rho_c,{\bf \mu}}:\OO^{\times}_{{\mathfrak{L}}_l}(-1)\times\left(Hom_{\CC} (\CC[q,q^{-1}],\CC)\right)^2 \longrightarrow Y_1(l) .
 	$$
 \end{pro}
  
  The reader may recall that the strata ${{\mathfrak L}_l}$ are further partitioned into the substrata ${{\mathfrak L}_l}(h^l, \lambda)$ where the labels $(h^l, \lambda)$ are respectively the dimension of the $l$-th piece of filtration and $\lambda$ is the partition associated to the dimensions of the graded pieces of the quotient
  $W_{\xi}/W_{\xi}([\phi])$. In addition, the partition $\lambda'$ conjugate to $\lambda$ has exactly $l$ parts. Thus the substrata ${{\mathfrak L}_l}(h^l, \lambda)$ become {\it canonical labels} for the
  maps of $\OO^{\times}_{{\mathfrak{L}}_l}(-1)$ into the modular curve $Y_1(l)$.
 
 \begin{thm}\label{pro:torsion1-intro}
 	Assume $l\geq 4$. Then each substratum ${{\mathfrak L}_l}(h^l, \lambda)$ of ${{\mathfrak L}_l}$ gives rise to distinguished continuous map
 	$$
 	\wp_{l,\rho_c,{ \lambda}}:\OO^{\times}_{{\mathfrak{L}}_l}(-1)\times\left(Hom_{\CC} (\CC[q,q^{-1}],\CC)\right)^2 \longrightarrow Y_1(l).
 	$$
 \end{thm} 
  
  The connection of the IVHS refinement with the moduli spaces of elliptic curves with marked points and the modular curves $Y_1(l)$ suggests other possibilities of categorifications as well as links to number theory, automorphic forms.

\vspace{0.2cm}
{\bf 4. Further research direction and applications.} There are several obvious directions to explore:

$\bullet$ to apply the refinement constructions to study families of canonically polarized curves; in particular, to apply them to study surfaces with nonhyperelliptic fibrations and, more generally, smooth projective varieties fibred by nonhyperelliptic curves of genus $g$,

$\bullet$ to develop IVHS refinement for higher dimensional varieties, in particular, compact complex manifolds with the canonical bundle very ample,

$\bullet$ to study the categorical aspects of the IVHS refinement.

More specific questions arise:

 can the refined invariants `see' that a curve 
 is 
 
 $\bullet$ an ample divisor 
on a smooth projective surface,

$\bullet$ a complete intersection in a projective space and, more generally, in a homogeneous space.

More generally, 

$\bullet$ what aspects of {\it extrinsic} geometry of a curve are captured by refined invariants?

\vspace{0.2cm}
As we alluded to at various points of the introduction, the connections with number theory emerge: this comes from the possibility to label the extensions  
of quiver representations in Theorem \ref{pro:repsext-intro} by arithmetic objects coming from number fields or from the fields of functions of curves over finite fields, see \S11.4 for more details.
The connection of the refinement with the moduli spaces ${\mathfrak{M}}^l_1$ and modular curves $Y_1(l)$ is another manifestation of
number theoretic ties.

\vspace{0.2cm}
Our refinement of IVHS naturally connects to quivers, Fano toric varieties, dimers on trivalent graphs. Those topics figure prominently in actual research in connection with mirror symmetry. The following is a speculative perspective on the IVHS refinement constructions from the point of view of mirror symmetry
\footnote{The reader who is not interested in this perspective can skip the next paragraph since no further mention of this aspect will appear in the paper.}. 

\vspace{0.2cm}
 {\bf 5. The mirror symmetry perspective.} Let us go back to the diagram \eqref{two-inc-diag}. As we explained earlier in the introduction, our refinement consists of attaching to the IVHS pair $([\xi], W_{\xi}) $ on the cohomological side of the diagram the filtrations of $W_{\xi}$ varying holomorphically with $[\phi]\in \PP(W_{\xi})$. It was also alluded to that this could be envisaged as a sort of K\"ahler structure varying with $[\phi]$.
Part of this structure, according to the second item of Theorem \ref{pro:refine}, `counts'
certain rational normal curves in the projective space $\PP(\HKC^{\ast})$,
the bottom right of the diagram \eqref{two-inc-diag}. 
 This invokes  some of the ingredients of Mirror symmetry: the objects, say Calabi-Yau varieties,
 are equipped with two structures - holomorphic and symplectic (B-side and A-side in physics parlor); for two mirror objects $M$ and $M'$, the mirror symmetry predicts `isomorphisms' between sides of different nature, that is, the B-side (resp. A-side) of $M$ is `isomorphic' to the A-side (resp. B-side) of $M'$. From this perspective  our suggestion is
 \begin{equation}\label{guide}
 	\begin{gathered}
 	\text{$\bullet$\it every part of the diagram \eqref{two-inc-diag} should come with its A-  and B- sides;}
 	\\
 	\text{$\bullet$\it the adjacent parts of the diagram are related by functorial equivalences}
 	\\
 	\text{\it between the sides of different nature.}
 	\end{gathered}
 \end{equation} 
For example, on the bottom left, $\PP(H^1(\Theta_{C}))$ is the projectivization of the space of the infinitesimal variation of complex structures on $C$; this qualifies for the B-side; our refinement could be viewed as algebraic K\"ahler structures along the arrow $p_1$ of the diagram and this is (a part of) A-side. This A-side should be functorially related to the B-side in the middle of the diagram and that one to the A-side on the bottom right. Rational normal curves in $\PP(\HKC^{\ast})$ related to the canonical embedding of $C$ should qualify for (a part of) the A-side of the canonical embedding
$C\subset \PP(\HKC^{\ast})$. The second item of Theorem \ref{pro:refine}
relating algebraic K\"ahler structures on $\PP(H^1(\Theta_{C}))$ to the rational normal curves appearing in the canonical embedding of $C$
could be viewed as an example of the hypothetical functorial relations in \eqref{guide}. The Hilbert polynomial
$$
H_{[\xi],[\phi]} (q)= \sum^l_{i=0} dim \left(W^i_{\xi}/W^{i+1}_{\xi} \right) q^i
$$
of the filtration \eqref{Wxiphi-filt-intro} could be viewed as a numerical output of this functorial relation. According to Theorem \ref{pro:refine}, the function counts the dimensions of families of certain rational curves related to the canonical embedding of $C$, a sort of Gromov-Witten type invariant.

The homological Mirror symmetry conjecture of Kontsevich recasts Mirror symmetry as equivalences between pairs of categories attached to each of the mirror partners, \cite{K}. In the case of Calabi-Yau varieties the two categories are the derived category of coherent sheaves, the holomorphic or B-side, and Fukaya category, the symplectic or A-side; for mirror partners
$M$ and $M'$ the conjecture predicts the equivalences of 

- the derived category of $M$ with the derived Fukaya category  of $M'$,

- the derived Fukaya category of $M$ with the derived category of $M'$.

In view of the homological mirror symmetry conjecture of Kontsevich the suggestion \eqref{guide} can be restated as follows:

\begin{equation}\label{guide-Kont}
	\begin{gathered}
		\text{$\bullet$\it every part of the diagram \eqref{two-inc-diag} should come with two types of categories:}\\
		\text{\it one related to the derived category of coherent sheaves,}
		\\
		 \text{\it the other to the Fukaya category of some (auxiliary) varieties;}
		\\
		\text{$\bullet$\it the adjacent parts of the diagram are related by equivalences}
		\\
		\text{\it between the categories  of different type.}
	\end{gathered}
\end{equation} 
Let us see how our situation fits into this pattern. The starting point for our constructions is to view the Kodaira-Spencer classes $\xi$ in $H^1(\Theta_{C})$ as the extension classes via the identification
$$
H^1 (\Theta_{C})=H^1(\OO_C (-K_C)) \cong Ext^1 (\OO_C (K_C), \OO_C).
$$
In other words, a nonzero Kodaira-Spencer class $\xi$ is viewed as the corresponding extension sequence
\begin{equation}\label{ext-seq-intro}
	\xymatrix{
		0\ar[r]& \OO_C \ar[r]& \EE_{\xi} \ar[r]& \OO_C (K_C)\ar[r]&0.
	}
\end{equation}
This way the B-side of $\PP(H^1 (\Theta_{C}))$ becomes related to the category of (short) exact complexes (up to the $\CC^{\times}$-action by the multiplication on the arrows of the complex) on $C$. Furthermore,
on the cohomology level we have the exact complex
$$
\xymatrix{
	0\ar[r]& H^0(\OO_C) \ar[r]& H^0(\EE_{\xi}) \ar[r]& H^0(\OO_C (K_C))\ar[r]^(.55){\xi}&H^1(\OO_C),
}
$$
where the rightmost arrow is the cup product with $\xi$. In particular, if the rank of the map
$$
\xymatrix{
	H^0(\OO_C (K_C))\ar[r]^(.55){\xi}&H^1(\OO_C)
}
$$
is $r$, that is the point $[\xi]$ lies in the open stratum 
$$
\Sigma^0_r:=\Sigma_r \setminus \Sigma_{r-1},
$$
then the rank $2$ bundle $\EE_{\xi}$ in the middle of the extension sequence has the property
$$
h^0(\EE_{\xi})=g-r+1.
$$
More precisely, our space
$$
W_{\xi}=ker\Big(\HKC \stackrel{\xi}{\longrightarrow} H^1(\OO_C)\Big),
$$
is a part of the following exact sequence
\begin{equation}\label{glsec-seq}
	\xymatrix{
		0\ar[r]& H^0(\OO_C) \ar[r]& H^0(\EE_{\xi}) \ar[r]&W_{\xi}\ar[r]&0. 
	}
\end{equation}
The pair of two exact sequences
\begin{equation}\label{cat-data}
	\xymatrix{
	0\ar[r]& \OO_C \ar[r]& \EE_{\xi} \ar[r]& \OO_C (K_C)\ar[r]&0,
	\\
		0\ar[r]& H^0(\OO_C) \ar[r]& H^0(\EE_{\xi}) \ar[r]&W_{\xi}\ar[r]&0, 
	}
\end{equation}
 could be viewed as a lifting of the IVHS data
 $$
 ([\xi], W_{\xi})
 $$
 to the category of perfect complexes on $C$ with an extra data on the level of global sections. Thus we have a natural way to relate $\PP(H^1(\Theta_{C}))$ equipped with the Griffiths stratification and the arrow $p_1$ of the diagram \eqref{two-inc-diag} to the derived category of coherent sheaves on $C$.
 
  The categorical data \eqref{cat-data} has an additional structure: the space of the global sections $H^0(\EE_{\xi})$ comes along with the exterior product
 \begin{equation}\label{ex-prod}
 \xymatrix{
 	\bigwedge^2 H^0(\EE_{\xi})\ar[r]& H^0(\bigwedge^2 \EE_{\xi}) =\HKC, 
 }
\end{equation}
 since the extension sequence \eqref{ext-seq-intro} implies $\bigwedge^2 \EE_{\xi}=\OO_C(K_C)$. This could be viewed as an $\HKC$-valued two-form
 attached to $([\xi], W_{\xi})$. The above exterior product is related to the one we have seen in \eqref{alphaxi-intro}. It is an essential ingredient in constructing the filtration \eqref{Wxiphi-filt-intro}. We suggest to view it as a sort of `symplectic' structure attached to $([\xi], W_{\xi})$. This should eventually lead to Fukaya-type categories related to the left side of the diagram \eqref{two-inc-diag}.
 
  \vspace{0.2cm}
  From the discussion on the nonabelian Dolbeault space ${\bf H^{1,0}}(PG_{l})$ as the parameter space of Higgs structures we also have a glimpse at the categories in the middle of the diagram \eqref{two-inc-diag}: 
   the objects are projective subspaces $\PP(W)$ in $\PP(\HKC)$  with the subspace $\Xi_W$ of $H^1(\Theta_{C})$ annihilating $W$ - the data delivered by the cohomological correspondence $P$. Over $\PP(\HKC)$  we also have the geometric correspondence 
 $$
 q_1: {\cal Z}_C \longrightarrow \PP(\HKC).
 $$
 On this side the basic sheaf object is the direct image
 $$
 (q_1)_{\ast} (\OO_{{\cal Z}_C})=(q_1)_{\ast} (q^{\ast}_2\OO_{C})
 	$$
of the structure sheaf $\OO_{{\cal Z}_C}$ of ${\cal Z}_C$. This is locally free with the fibre over a point $[\phi] 	\in \PP(\HKC)$ being the space
$$
H^0 (\OO_{Z_{\phi}})
$$
of functions on the divisor $Z_{\phi}=(\phi=0)$, the zero locus of $\phi$.
The fibre of the projection
$$
p_2: {\bf P} \longrightarrow \PP(\HKC)
$$ 
in \eqref{two-inc-diag} over $[\phi]$ is the projectivization of the kernel $\Xi_{\phi}$ of the map
$$
H^1 (\Theta_{C}) =H^1 (\OO_C (-K_C))\stackrel{\phi}{\longrightarrow} H^1(\OO_C)
$$
and we can see its relation to the space $H^0 (\OO_{Z_{\phi}})$. Namely, the above multiplication by $\phi$ on the cohomology level comes from
the exact sequence of sheaves on $C$
$$
\xymatrix{
0\ar[r]& \OO_C (-K_C) \ar[r]^(.63){\phi} & \OO_C \ar[r]& \OO_{Z_{\phi}} \ar[r]&0.
}
$$
This gives rise to the exact sequence
$$
\xymatrix{
	0\ar[r]  &H^0 (\OO_C) \ar[r]& H^0(\OO_{Z_{\phi}}) \ar[r]&\Xi_{\phi} \ar[r]&0.
} 
$$
The subspace $\Xi_W$ is contained in $\Xi_{\phi}$ for all nonzero $\phi \in W$. Lifting this subspace to $H^0(\OO_{Z_{\phi}})$ gives a distinguished subspace $\widetilde{\Xi}_W (\phi)$ of functions on $Z_{\phi} $ containing the subspace of constant functions. As $[\phi]$ varies in $\PP(W)$ the subspaces $\widetilde{\Xi}_W (\phi)$ fit together to form a subsheaf, call it ${\cal X}_W$, of the restriction of the sheaf $ (q_1)_{\ast} (\OO_{{\cal Z}_C})$ to $\PP(W)$.  Thus the cohomological pair
$(\PP(W), \Xi_W)$ lifts to the level of coherent sheaves supported on $\PP(W)$, the pair $(\PP(W), {\cal X}_W)$:
$$
\xymatrix{
(\PP(W), \Xi_W) \ar@{~>}[r]&(\PP(W), {\cal X}_W).
}
$$
We can do more: the subspaces $\widetilde{\Xi}_W (\phi)$ produce the filtrations of $H^0(\OO_{Z_{\phi}})$ which are cousins of the filtrations
\eqref{Wxiphi-filt-intro} on the side of $\PP(H^1(\Theta_C))$. This in turn leads to the graphs of the {\it same} type as $PG_l$ as well as to its own nonabelian Dolbeault varieties with their Lagrangian cycles! Furthermore, according to Theorem \ref{pro:repsext-intro1} and Theorem \ref{pro:torsion1-intro}, both sides are related to the moduli spaces of elliptic curves with marked points and the modular curves $Y_1(l)$, possibly for different values of $l$. This should lead to interesting correspondences between two sides. The suggestion \eqref{guide-Kont} gains somewhat in substance.

At this stage the suggestion formulated in \eqref{guide-Kont} serves as a working hypothesis. However, by following it, we believe that one can discover interesting new structures in the theory of curves and beyond.  

\vspace{0.2cm}
\noindent
 {\bf Organization of the paper.}
 The paper starts with a somewhat dubious premise of treating special cases before exposing general theory:  the first two sections treat the cases of Kodaira-Spencer classes of rank zero (\S1) and one (\S2). The reader will find here the essential technical features appearing in general situation:
 
 - the relation of invariant subspaces of the filtration with proportional global sections of $\EE_{\xi}$, the rank two bundle in \eqref{cat-data}, see the proof of Proposition \ref{pro:MaxN-gen};
 
 - the filtration on $W_{\xi}$ and its relation to rational normal curves
 in the canonical embedding of a curve, see Lemma \ref{lem:filt-phi}, the proof of Theorem \ref{thm:rk1KS}.
 
 In particular, the two special cases (re)prove old classical theorems of Max Noether on projective normality of the canonical embedding of curves and the Babbage-Enriques-Petri theorem about quadrics passing through canonical curves. The reader who wants to see general constructions and not interested in reproving the old results should go directly to \S3.
 
 In \S3 we give the construction and main properties of what we call
 $([\xi],[\phi])$-filtration of $W_{\xi}$, the refinement of the IVHS invariant $W_{\xi}$.
 
 The section \S4 recasts the $([\xi],[\phi])$-filtration as an orthogonal decomposition of the space $W_{\xi}$ as well as the decomposition into irreducible ${\mathfrak sl}(2,\CC)$-modules, a sort of Hard Lefschetz decomposition of $W_{\xi}$; the resulting structure is the one we alluded to as an algebraic K\"ahler structure on $W_{\xi}$.
 
 In \S5 we  consider the cohomological incidence correspondence ${\bf P}$ from the point of view of the projection onto $\PP(\HKC)$. The IVHS pair
 $(\xi, W_{\xi})$ is reinterpreted in terms of the sheaf data on the projective space $\PP(W_{\xi})$. On the technical level, the pair $([\xi],[\phi])$ in the incidence correspondence ${\bf P}$ is related to a function $F(\phi,\xi)$ on the zero divisor $Z_{\phi}=(\phi=0)$. On the conceptual level the results here establish a precise link between the left and right sides of the two part diagram \eqref{two-inc-diag}: the cohomological data on the left side is translated into spectral data on the right side.
 
 In \S6 we recast the algebraic K\"ahler structure on $W_{\xi}$ as quiver representations; the graph $PG_l$ underlying the quiver is defined as well as the Fano toric variety ${\bf H^{1,0}}(PG_l)$, the one we call the nonabelian Dolbeault space.
 
 \S7 introduces the trivalent completion $\widehat{PG}_l$ of the graph $PG_l$. For the natural ribbon structure of $\widehat{PG}_l$ we determine its boundary cycles and deduce that ribbon graph structure gives us a topological torus $\mathbb{T}$. Then we determine so called zig-zag paths of $\widehat{PG}_l$. Those play a key role in defining later on quantum-type invariants, the technical tool which allows to extend quiver representations of $PG_l$ to the ones of $\widehat{PG}_l$.
 
 In \S8 we use the zig-zag paths to construct a toric surface $X(\Delta_l)$.
 
 In \S9 is discussed the dual graph $Q_l$ of $\widehat{PG}_l$ which is naturally a quiver; its path algebra is equipped with a potential $W_l$ which is also determined. Following Ginzburg, \cite{Gi}, one defines the Jacobi algebra $A_{Q_l}$ associated to the pair $(Q_l, W_l)$. The general results \cite{Bo}, \cite{Br}, on the Jacobi algebras of dimer models imply
 that $A_{Q_l}$ is a Calabi-Yau algebra of dimension three. 

 \S10 discusses the set of perfect matchings of $\widehat{PG}_l$; those have interesting combinatorial structure related to closed paths on $\widehat{PG}_l$. We suggest the set of perfect matchings of $\widehat{PG}_l$ should admit categorifications. As a manifestation of this we show that the perfect matchings of $\widehat{PG}_l$ can be organized into a complex.
 
 \S11 devoted to the construction of what we call quantum-type invariants: this is a technical tool to prove the extension property stated in Theorem \ref{pro:repsext-intro}; it also connects to symplectic vector spaces and Lagrangians as well as suggests higher categorical structures of the refinement of IVHS. 
 
 In \S12 we add metrics to the ribbon graph $\widehat{PG}_l$ and relate our refinement of IVHS to the moduli spaces of elliptic curves with $l$ marked points as well the modular curve $Y_1(l)$.
 
 In \S13 we show how to construct local systems on the elliptic curves attached to the refinement of IVHS.
 
 The concluding section, \S14, is of technical nature: we give a sheaf version of the $([\xi],[\phi])$-filtrations.
 
 For convenience of the reader we summarize the main topics discussed in the paper in the following chart; the question marks below stand for connections to be explored.
 $$
 \xymatrix@R=16pt{
 *+[F]\txt{Refinement of IVHS}\ar^{\txt{Spectral\\data}}[rr]\ar[dd]_{\txt{Hodge-type\\decomposition}}&&*+[F]\txt{Special divisors,\\rational normal curves,\\Vector bundles}&\\
 &&&\\
 *+[F]\txt{Quivers $PG_l$,\\their representations}\ar[dd]_{\txt{Quantum-type\\invariants}}\ar[rr]^(.45){\txt{Higgs\\ data}}&&*+[F]\txt{Fano Toric ${\bf H^{1,0}}(PG_l)$\\and its Lagrangian cycle $H_0$}\ar[r]^(.7){\txt{?}}&*+[F]\txt{Fukaya\\ category}\\
 &&&\\
*+[F]\txt{Quivers $\widehat{PG}_l$\\and their representations}\ar[rr]\ar[d]_{\txt{metrics on\\ribbon graphs}}&&*+[F]\txt{Dimers,\\ perfect matchings,\\hexagon tilings,\\$3d$-partitions,\\Topological vertex}\ar[r]^(.6){\txt{?}}&*+[F]\txt{Higher\\ categories} \\
*+[F]\txt{Moduli spaces ${\mathfrak{M}}^l_1$,\\ modular curves $Y_1(l)$}\ar[rr]^{\txt{?}}&&*+[F]\txt{Automorphic forms,\\Number theory}
	&
}
$$

\section{Theorem of Max Noether}
Let $C$ be a smooth complex projective curve of genus $g \geq 2$. The differential of the period map, as we reminded in the introduction, is the cohomology
cup product
\begin{equation}\label{cp}
H^1(\Theta_C) \longrightarrow Hom (H^0 (\OO_C (K_C)), H^1 (\OO_C)).
\end{equation}
For each nonzero $\xi$ in $H^1 (\Theta_C)$, the above attaches to $\xi$ the
homomorphism
$$
\xi : H^0 (\OO_C (K_C)) \longrightarrow  H^1 (\OO_C).
$$
Its kernel $ker(\xi)$ will be denoted $W_{\xi}$:
$$
W_{\xi}:=ker(\xi).
$$
We want to understand the geometry of $W_{\xi}$. The starting point is the following.

\begin{thm}\label{th-Max}
	Assume $C$ is not hyperelliptic, that is, the canonical map
	$$
	|K_C|: C \longrightarrow \PP(H^0 (\OO_C (K_C))^{\ast})
	$$
	is an embedding. Then $dim(W_{\xi}) \leq g-1$, for every nonzero $\xi \in H^1(\Theta_C) $. In other words, the cup product in \eqref{cp} is injective.
\end{thm}

We will prove a somewhat more general statement.
\begin{pro}\label{pro:MaxN-gen}
	Let $\OO_C (L)$ be a very ample line bundle on $C$. Assume
	$$
	d=deg(L) \geq 2h^0 (\OO_C (L)) -3.
	$$
	Then the cup product
	$$
	H^1 (\OO_C (-L)) \longrightarrow Hom(	H^0 (\OO_C (L)), H^1 (\OO_C))
	$$
	is injective.
\end{pro}

Thus Theorem \ref{th-Max} is obtained from the above proposition by taking
$L=K_C$, for $C$ nonhyperelliptic.
In the rest of the section we prove Proposition \ref{pro:MaxN-gen}. We will need some preliminary considerations. 

We begin with the identification
\begin{equation}\label{H1=Ext1}
H^1 (\OO_C (-L))\cong Ext^1 (\OO_C (L),\OO_{C}).
\end{equation}
So a cohomology class $\xi \in H^1 (\OO_C (-L))$ is thought of as the extension sequence
\begin{equation}\label{ext-pf}
	\xymatrix{
		0\ar[r] & \OO_C \ar[r]&{\cal E}_{\xi} \ar[r]& \OO_C (L)\ar[r]&0
	}
\end{equation}
corresponding to $\xi$ under the identification in \eqref{H1=Ext1}. In particular, the extension sequence splits if and only if $\xi=0$.

From the extension sequence we obtain the exact complex of cohomology groups
$$
\xymatrix{
	0\ar[r] & H^0 (\OO_C )\ar[r]&H^0 ({\cal E}_{\xi} )\ar[r]& H^0(\OO_C (L))\ar[r]^{\xi}&H^1 (\OO_C),
}
$$
where the rightmost arrow is the cup-product with $\xi$. This gives the exact sequence
$$
	\xymatrix{
	0\ar[r] & H^0 (\OO_C) \ar[r]&H^0 ({\cal E}_{\xi}) \ar[r]& ker(\xi)\ar[r]&0.
}
$$
Set $W:=ker(\xi)$ and assume it to be nonzero. Choose a splitting
$$
\alpha: W \longrightarrow H^0 ({\cal E}_{\xi}).
$$
Its image is a codimension $1$ subspace of $H^0 ({\cal E}_{\xi})$ and we denote
$$
\alpha(\phi)= \text{\it the global section of ${\cal E}_{\xi}$ corresponding
	to $\phi \in  W$ under $\alpha$.}
$$
Observe the following.
\begin{lem}\label{lem:cphiphi'}
	For every pair $\phi$ and $\phi'$ in $W$ we have
	$$
	\phi' \alpha(\phi)-\phi \alpha(\phi')=\alpha^{(2)}(\phi,\phi') e_{\xi}
	$$
	for a unique $\alpha^{(2)}(\phi,\phi') \in H^0 (\OO_C (L))$ and where $e_{\xi}$ is the global section of ${\cal E}_{\xi}$, the image of the constant $1$ under the monomorphism in the extension sequence; the equality above is
	in the space of global sections of $\EE_{\xi} (L)$.
\end{lem}
\begin{pf}
	The extension sequence \eqref{ext-pf} can be viewed as the Koszul sequence of the global section $e_{\xi}$ as in the statement of the lemma
	$$
	\xymatrix{
		0\ar[r] & \OO_C \ar[r]^{e_{\xi}}&{\cal E}_{\xi} \ar[r]^(.4){\wedge e_{\xi}}& \OO_C (L)\ar[r]&0.
	}
	$$
	Tensor with $\OO_C (L)$ and observe that the global section
	$$
	\phi' \alpha(\phi)-\phi \alpha(\phi')
	$$
	of ${\cal E}_{\xi} (L) $ goes to the section
	$$
	\phi' \phi-\phi \phi' =0
	$$
	in $H^0(\OO_C (2L))$. Hence there is a unique global section of $\OO_C (L)$ denoted $\alpha^{(2)}(\phi,\phi')$ such that
	$$
	\phi' \alpha(\phi)-\phi \alpha(\phi')=\alpha^{(2)}(\phi,\phi')e_{\xi}.
	$$
\end{pf}

The sections $\alpha^{(2)}(\phi,\phi')$, as $\phi$ and $\phi'$ ranges through $W$, are organized into the linear map
$$
\alpha^{(2)}: \mbox{ $\bigwedge^2 W \longrightarrow H^0 (\OO_C (L))$}
$$
where $\phi \wedge \phi'$ is mapped to $\alpha^{(2)}(\phi,\phi').$

\begin{lem}\label{lem:c-Kosz}
	The map	$\alpha^{(2)}$ constructed above is a Koszul cocycle, that is, for every triple $\phi, \phi', \phi''$ in $W$ we have the cocycle relation:
	$$
	\phi'' \alpha^{(2)}(\phi,\phi')-\phi'\alpha^{(2)}(\phi,\phi'')+\phi \alpha^{(2)}(\phi',\phi'')=0.
	$$	 
\end{lem}
\begin{pf}
	Follows immediately from the formula in Lemma \ref{lem:cphiphi'}.
\end{pf}

We now turn to the proof of Proposition \ref{pro:MaxN-gen}.

\vspace{0.2cm}
\noindent
{\it Proof of Proposition \ref{pro:MaxN-gen}.} Let $\xi$ be a cohomology class in the kernel of the cup-product
$$
	H^1 (\OO_C (-L)) \longrightarrow Hom(	H^0 (\OO_C (L)), H^1 (\OO_C)).
$$
This means that $W$ in the above considerations is $H^0 (\OO_C (L))$ and we have constructed the map
$$
\alpha^{(2)}: \mbox{ $\bigwedge^2 H^0 (\OO_C (L)) \longrightarrow H^0 (\OO_C (L))$.}
$$

Let us fix a nonzero $\phi$ in $H^0 (\OO_C (L))$ and consider the linear map
$$
\alpha^{(2)}(\phi, \bullet): H^0 (\OO_C (L)) \longrightarrow H^0 (\OO_C (L)).
$$
Since $\alpha^{(2)}(\phi,\phi)=0$ this induces an endomorphism
$$
\overline{\alpha^{(2)}}(\phi, \bullet): H^0 (\OO_C (L))/\CC\phi \longrightarrow H^0 (\OO_C (L))/\CC \phi.
$$
Let $\{\phi_0\} \in  H^0 (\OO_C (L))/\CC\phi$ be an eigenvector of
$ \overline{\alpha^{(2)}}(\phi, \bullet)$ with the eigenvalue $\lambda_0$. This gives the formula
\begin{equation}\label{phi0-form}
	\alpha^{(2)}(\phi, \phi_0)=\lambda_0 \phi_0 +\mu_0 \phi,
\end{equation}
for some constant $\mu_0 \in \CC$.

For any $\psi \in H^0 (\OO_C (L))$ we have the cocycle relation for the triple $\phi,\phi_0, \psi$, see Lemma \ref{lem:c-Kosz} :
$$
\begin{gathered}
	0=	\psi \alpha^{(2)}(\phi,\phi_0)-\phi_0 \alpha^{(2)}(\phi,\psi)+\phi \alpha^{(2)}(\phi_0,\psi)=\\
	\psi (\lambda_0 \phi +\mu_0 \phi_0)-\phi_0 \alpha^{(2)}(\phi,\psi)+\phi \alpha^{(2)}(\phi_0,\psi)=
	\phi (\alpha^{(2)}(\phi_0,\psi)+\lambda_0 \psi) - \phi_0 (\alpha^{(2)}(\phi,\psi) -\mu_0 \psi),
\end{gathered}
$$
where the second equality uses the formula \eqref{phi0-form}. Thus we obtain
\begin{equation}\label{phi0psi-form}
	\phi (\alpha^{(2)}(\phi_0,\psi)+\lambda_0 \psi) - \phi_0 (\alpha^{(2)}(\phi,\psi) -\mu_0 \psi) =0, \,\, \forall \psi \in H^0 (\OO_C (L)).
\end{equation}

So far we did not impose any conditions on the nonzero global section $\phi$ of $\OO_C (L)$. We now use the assumption that $\OO_C (L)$ is very ample. So we think of $C\subset \PP((H^0 (\OO_C (L))^{\ast}))$ and 
 chose $\phi$ in $H^0 (\OO_C (L))$  so that its divisor of zeros 
$$
Z_{\phi} =\{\phi=0\}
$$
is a set of $d=deg(L)$ distinct points in general position in the hyperplane 
$$
H_{\phi} \subset \PP(H^0 (\OO_C (L))^{\ast})
$$
 corresponding to $\phi$, that is, for any $k \in [1,n]$, any subset of $k$ points in $Z_{\phi}$ spans the projective subspace of dimension $(k-1)$; the value
$n=h^0(\OO_C (L))-1$. 
With this choice of $\phi$ we return to the equation \eqref{phi0psi-form}.

Let $B$ be the base locus of the pencil $\{\phi,\phi_0\}$. From the general position of $Z_{\phi}$, we deduce 
$$
deg(B)\leq (n -1).
$$
From this it follows that the complement $Z_{\phi} \setminus B$ consists of at least $d-(n-1)\geq (2n-1)-(n-1)=n$ points. Hence those points 

\vspace{0.2cm}
- lie outside of the hyperplane
$H_{\phi_0}$ corresponding to $\phi_0$;

- they span the hyperplane $H_{\phi}$.

\vspace{0.2cm}
The first item and the equation \eqref{phi0psi-form} tell us that the global section $(\alpha^{(2)}(\phi,\psi) -\mu_0 \psi)$ must vanish on
$Z_{\phi} \setminus B$ and the second item says that the hyperplane
corresponding to $(\alpha^{(2)}(\phi,\psi) -\mu_0 \psi)$ coincides with $H_{\phi}$.
The latter means that
$$
\alpha^{(2)}(\phi,\psi) -\mu_0 \psi =\mu(\psi)\phi,
$$
for some scalar $\mu(\psi)$ depending on $\psi$. Since the left hand side is linear with respect to $\psi$, we deduce that 
$
\mu(\psi) 
$
is a linear function of $\psi$. The above equation gives the formula
$$
\alpha^{(2)}(\phi,\psi) =\mu_0 \psi +\mu(\psi)\phi, \,\forall \psi \in H^0 (\OO_C (L)).
$$
Substitute it into the formula in Lemma \ref{lem:cphiphi'} to obtain
$$
\psi \alpha(\phi)-\phi \alpha(\psi)=(\mu_0 \psi +\mu(\psi)\phi) e_{\xi}.
$$
This can be rewritten as follows
$$
\psi (\alpha(\phi)-\mu_0 e_{\xi}) -\phi (\alpha(\psi)+\mu(\psi)e_{\xi})=0, \,\forall \psi \in H^0 (\OO_C (L)).
$$
This means that all global sections $(\alpha(\psi)+\mu(\psi)e_{\xi})$ of ${\cal E}_{\xi}$ in the extension sequence are proportional to $(\alpha(\phi)-\mu_0 e_{\xi}) $. Hence they generate a line subsheaf of ${\cal E}_{\xi}$ which we denote by ${\cal L}$,
that is, we have a morphism
$$
{\cal L} \longrightarrow {\cal E}_{\xi}.
$$
Combining this with the extension sequence gives the diagram
$$
\xymatrix{
	&&{\cal L} \ar[d] \ar[rd]&&\\
	0\ar[r] & \OO_C \ar[r]&{\cal E}_{\xi} \ar[r]& \OO_C (L)\ar[r]&0
}
$$
where the slanted arrow of the diagram is the composition of the vertical and the horizontal arrows. From the definition of $\alpha$ and $e_{\xi}$ we deduce that the slanted arrow is an isomorphism on the level global sections.
Since $\OO_C (L)$ is globally generated, it follows that the slanted arrow is an isomorphism and hence gives the splitting of the extension sequence.  This means that the class $\xi$ is zero. $\Box$

\vspace{0.2cm}
As we already said Theorem \ref{th-Max} is obtained from Proposition \ref{pro:MaxN-gen} by setting $L=K_C$. For the convenience of the reader we recall how Theorem \ref{th-Max} implies the classical theorem of Max Noether
about projective normality of the canonical embedding.

\vspace{0.2cm}
We rewrite the cup-product in \eqref{cp} as follows
$$
\begin{gathered}
H^1(\Theta_{C}) \longrightarrow Hom (\HKC, H^1 (\OO_C))=(\HKC)^{\ast}\otimes H^1 (\OO_C)\\
\cong (\HKC)^{\ast}\otimes (\HKC)^{\ast},
\end{gathered}
$$
where the last identification uses the Serre duality
$H^1 (\OO_C) \cong \HKC^{\ast}.$ Furthermore, the image lies in the space of symmetric tensors due to the following obvious identity
$$
0=\xi(\phi \wedge \psi')=(\xi\phi) \wedge \phi' - \phi \wedge (\xi\phi'),
$$
for all $\phi,\phi' \in \HKC$. Hence the above cup-product takes the form
$$
H^1(\OO_C (-K_C))=H^1(\Theta_{C}) \longrightarrow S^2 (\HKC)^{\ast}
$$
and Theorem \ref{th-Max} says that this is injective. Equivalently, the dual map
$$
S^2 \HKC \longrightarrow H^1(\OO_C (-K_C))^{\ast} \cong H^0 (\OO_C (2K_C))
$$
is surjective. This is the quadratic normality of the classical theorem of Max Noether, see \cite{G-H}. The surjectivity
$$
S^l \HKC \longrightarrow  H^0 (\OO_C (lK_C)),
$$
for $l\geq 3$, follows easily by induction on $l$, since for a hyperplane section $Z_{\phi} =(\phi=0)$
of $C$ with  $Z_{\phi} $ in general position, the map
$$
S^l \HKC \longrightarrow   H^0 (\OO_{Z_{\phi}} (lK_C))
$$
is surjective for all $l\geq 3$ by a well known argument of Castelnuovo, see \cite{G-H}. Thus Theorem \ref{th-Max} implies theorem of Max Noether on the projective normality of the canonical embedding of a smooth projective curve.

\section{Kodaira-Spencer classes of rank 1}

We turn to the study of Kodaira-Spencer classes of rank $1$. That is we assume that $\xi \in H^1 (\Theta_{C})$ gives rise to the linear map
$$
\xi: H^0 (\OO_C (K_C)) \longrightarrow H^1 (\OO_C)
$$
whose rank is $1$. We let
$$
W_{\xi}:=ker(\xi).
$$
This is a subspace of codimension $1$ in $H^0 (\OO_C (K_C))$ or, dually, a point in $\PP(H^0 (\OO_C (K_C))^{\ast})$. Our aim is to show that the point lies on the canonical image of $C$, unless $C$ has special linear systems.

\begin{thm}\label{thm:rk1KS}
	Let $C$ be a smooth projective curve of genus $g\geq 3$ and let  $\xi \in H^1 (\Theta_{C})$ be a Kodaira-Spencer class of rank $1$. Then
	the kernel $W_{\xi}=ker(\xi)$ corresponds to the linear subsystem
	$|W_{\xi}|$ with a base point, unless $C$ carries the linear systems
	$g^1_2$, $g^1_3$, $g^2_5$.
\end{thm}

The rest of the section is devoted to the proof of this theorem and its corollaries.
Our approach as with Noether's theorem is the extension sequence which is synonym of the Koszul sequence
	$$
\xymatrix{
	0\ar[r] & \OO_C \ar[r]^{e_{\xi}}&{\cal E}_{\xi} \ar[r]^(.4){e_{\xi}\wedge }& \OO_C (K_C)\ar[r]&0
}
$$
where $e_{\xi}$ is as in Lemma \ref{lem:cphiphi'}. Passing to the global sections gives the exact sequence
\begin{equation}\label{gl-sec-xi}
\xymatrix{
	0\ar[r] & H^0 (\OO_C) \ar[r]&H^0 ({\cal E}_{\xi}) \ar[r]& W_{\xi}\ar[r]&0
}
\end{equation}
As before we take a splitting
$$
\alpha: W_{\xi} \longrightarrow H^0 ({\cal E}_{\xi})
$$
which gives the direct sum decomposition
$$
H^0({\cal E}_{\xi}) =\CC e_{\xi} \oplus \alpha(W_{\xi}).
$$

The proof of Lemma \ref{lem:cphiphi'} applies to give  the relation
\begin{equation}\label{alpah-c}
	\phi' \alpha(\phi)-\phi \alpha(\phi') =\alpha^{(2)}(\phi,\phi')e_{\xi}, \, \forall \phi,\phi' \in W_{\xi},
\end{equation}
 This gives the linear map
$$
\alpha^{(2)}: \mbox{$\bigwedge^2 W_{\xi} \longrightarrow H^0(\OO_C (K_C))$}
$$
which sends $\phi \wedge \phi'$ to $\alpha^{(2)}(\phi,\phi')$ and it is subject to the cocycle relation
\begin{equation}\label{c-cocycle}
	\phi'' \alpha^{(2)}(\phi,\phi')-\phi'\alpha^{(2)}(\phi,\phi'')+\phi \alpha^{(2)}(\phi',\phi'')=0,
\end{equation}
for any triple $\phi,\phi',\phi''$ in $W_{\xi}$.

\vspace{0.2cm}
From now on we assume:

$\bullet$ the linear system $|W_{\xi}|$ has no base points,

$\bullet$ $C$ has no $g^1_2$, $g^1_3$, $g^2_5$.

The first item tells us that the linear system $|W_{\xi}|$ defines a morphism
\begin{equation}\label{morph-proj}
|W_{\xi}|: C \longrightarrow \PP(W^{\ast}_{\xi})
\end{equation}
which is the projection of the canonical image of $C$ in $\PP(\HKC^{\ast})$ from the point $p_{\xi}=\PP(W^{\perp}_{\xi})$ in
$\PP(\HKC^{\ast})\setminus C$. And the second item implies the following.

\begin{lem}\label{lem:proj-birational}
	Assume $C$ has no $g^1_2$, $g^1_3$, $g^2_5$. Then the projection morphism
	$$
	|W_{\xi}|: C \longrightarrow \PP(W^{\ast}_{\xi})
	$$
	is birational onto its image.
\end{lem}
\begin{pf}
	Let $C'$ be the image of $C$ under the morphism defined by $|W_{\xi}|$.
	The condition that $C$ has no $g^1_2$ tells us that $C$ can be identified with its canonical image, while the condition no $g^1_3$ means that $C\subset \PP(\HKC^{\ast})$ has no trisecant lines. Hence the morphism
	$$
	C \longrightarrow C'
	$$ 
	defined by $|W_{\xi}|$ is of degree at most $2$. Assume it is of degree $2$. Taking the normalization $C'_{norm}$ of $C'$ gives
	a double covering
	$$
	f:C \longrightarrow  C'_{norm}
	$$
	such that $K_C =f^{\ast} (D')$, where $D'$ is a divisor of degree
	$(g_C -1)$ on $C'_{norm}$; furthermore, we have
	 $$
	 h^0 (\OO_{C'_{norm}} (D')) \geq dim(W_{\xi})=g_C-1.
	 $$
	  From the Riemann-Roch for $\OO_{C'_{norm}} (D')$
	$$
	g_C-1 -h^1 (\OO_{C'_{norm}} (D'))\leq \chi(\OO_{C'_{norm}} (D'))=g_C -g_{C'_{norm}},
	$$
	we deduce
	$$
	h^1 (\OO_{C'_{norm}} (D')) \geq g_{C'_{norm}} -1.
	$$
	If $g_{C'_{norm}} \geq 2$ we obtain
	$$
	\begin{gathered}
	g_{C'_{norm}} =h^0 (\OO_{C'_{norm}} (K_{C'_{norm}}))) \geq 	h^0 (\OO_{C'_{norm}} (D')) + 	h^0 (\OO_{C'_{norm}} (K_{C'_{norm}}-D')) -1 \\
	\geq g_C -1 +g_{C'_{norm}} -2=g_C +g_{C'_{norm}} -3,
	\end{gathered}
	$$
	and this tells us that $g_C \leq 3$ which is clearly impossible under the assumptions on $C$. Thus $g_{C'_{norm}} \leq 1$. Since $C$ can not be a double cover of a rational curve, the assumption no $g^1_2$, we deduce that $C'$ is a smooth curve of genus $1$ and $C$ is a double cover of $C'$ branched along a reduced effective divisor $B' \in |2D'|$. This means that the curve $C$ admits an involution, call it $\sigma$, and $C'=C/\sigma$. In particular, $\sigma$ acts on $H^0 (\OO_C (K_C))$ and we have the direct sum decomposition
	$$
	H^0 (\OO_C (K_C)) =H^+ \oplus H^-,
	$$
	where  $H^{\pm}$ is invariant/anti-invariant subspace of $\HKC$ under the action of $\sigma$. From the formula of the canonical divisor

	$$
	K_C =f^{\ast}(D')
	$$
	we deduce 
	$$
	H^+= H^0(\OO_{C'}), \,\, H^-=H^0 (\OO_{C'} (D'))
	$$
	Similarly, we have the action of $\sigma$ on all the spaces involved in the cup-product map
	$$
	H^1 (\Theta_C) \longrightarrow \HKC^{\ast} \otimes H^1 (\OO_C)
	$$
	and the map is $\sigma$-equivariant:
	$$
	\begin{gathered}
		H^1 (\Theta_C) =	H^1 (\Theta_C)^+ \oplus 	H^1 (\Theta_C)^- \cong H^1 (\OO_{C'} (-B')) \oplus H^1 (\OO_{C'} (-\HA B'))\\
		=H^1 (\OO_{C'} (-2D')) \oplus H^1 (\OO_{C'} (-D')),
	\\
	H^1 (\OO_C) =	H^1 (\OO_C)^+ \oplus 	H^1 (\OO_C)^- \cong H^1 (\OO_{C'}) \oplus H^1 (\OO_{C'} (-D'));
	\\
	\\
H^1 (\Theta_{C})=	H^1 (\OO_{C'} (-2D')) \oplus H^1 (\OO_{C'} (-D')) \longrightarrow
\\
	\Big(H^0 (\OO_{C'}) \oplus H^0 (\OO_{C'} (D'))\Big)^{\ast} \otimes 
	\Big(H^1 (\OO_{C'}) \oplus H^1 (\OO_{C'} (-D'))\Big).
	\end{gathered}
	$$
	We write our Kodaira-Spencer class $\xi$ according the above decomposition
	$$
	\xi=\xi^+ +\xi^-.
	$$	
By the assumption $\xi$ is annihilated by the anti-invariant part
$
W_{\xi}\cong H^0 (\OO_{C'} (D')) \cong H^-
$
of $\HKC$. Hence $\xi^+$ and $\xi^-$ are respectively in the kernels of the following maps
$$
	H^1 (\OO_{C'} (-2D')) \longrightarrow  H^0 (\OO_{C'} (D'))^{\ast} \otimes H^1 (\OO_{C'} (-D')),
	$$
	$$
		H^1 (\OO_{C'} (-D')) \longrightarrow  H^0 (\OO_{C'} (D'))^{\ast} \otimes H^1 (\OO_{C'} ).
	$$
But both of those maps are injective: the second one is Serre duality isomorphism; the first one, after Serre duality, becomes
$$
H^0 (\OO_{C'} (2D'))^{\ast} \longrightarrow \Big(S^2(H^0 (\OO_{C'}(D')))\Big)^{\ast};
$$
its dual
$$
S^2(H^0 (\OO_{C'}(D'))) \longrightarrow H^0 (\OO_{C'} (2D'))
$$
is the multiplication map and this is surjective since $\OO_{C'} (D')$  defines a projectively normal embedding of $C'$. Thus $\xi^{\pm}=0$. But this means that the Kodaira-Spencer class $\xi$ is zero as well and this is contrary to our assumption that $\xi$ is nonzero.
\end{pf}

It is plausible to assume that the properties of the 
 vector bundle ${\cal E}_{\xi}$ are related to the geometry of $C$. This will be confirmed shortly. To begin with, we have the obvious homomorphism
 \begin{equation}\label{wedge2}
 	\mbox{$\bigwedge^2 H^0 ({\cal E}_{\xi}) \longrightarrow H^0 (\bigwedge^2 {\cal E}_{\xi})=\HKC.$}
 \end{equation}
Observe that this is related to the map
$$
\mbox{$\alpha^{(2)}: \bigwedge^2 W_{\xi} \longrightarrow \HKC$}
$$
 defined previously: we have the formula
 \begin{equation}\label{a2-wedge}
 \text{\it $\alpha^{(2)}(\phi,\phi')=$ the image of $\alpha(\phi)\wedge\alpha(\phi')$ under the map \eqref{wedge2}, $\forall \phi,\phi' \in W_{\xi}$.}
 \end{equation} 
 In other words the following diagram commutes
 \begin{equation}\label{a2-wedge-diag}
 	\xymatrix{
 	&\bigwedge^2 W_{\xi}\ar[ld]_{\alpha \wedge \alpha} \ar[rd]^{\alpha^{(2)}}&\\
 	\bigwedge^2H^0 (\EE_{\xi})\ar[rr]^{\eqref{wedge2}}&&\HKC
 }
 \end{equation}
A link of the map \eqref{wedge2} to geometry is the following.
\begin{lem}\label{lem:Exi-glgen}
	Assume $g_C \geq 4$. Then the kernel of the homomorphism in \eqref{wedge2} contains a nonzero decomposable tensor, that is, there exist two linear independent global sections $e$ and $e'$ of ${\cal E}_{\xi}$ whose exterior product $e\wedge e'$ maps to zero under the  
	homomorphism \eqref{wedge2}.
\end{lem}
\begin{pf}
	This is a simple dimension count: the projectivization of the decomposable tensors in $\bigwedge^2 H^0 ({\cal E}_{\xi})$ is the Grassmannian ${\bf Gr}(1,\PP(H^0 ({\cal E}_{\xi})))$ of lines in
	$\PP(H^0 ({\cal E}_{\xi}))$. This tells us that the necessary condition for the kernel of 
	\eqref{wedge2} to avoid the set of decomposable tensors is the inequality
	$$
	2(g_C -2)=dim({\bf Gr}(1,\PP(H^0 ({\cal E}_{\xi})))) \leq dim(\PP(\HKC))=g_C -1.
	$$
Hence the inequality $g_C \leq 3$.
\end{pf}

Our assumptions on $C$ imply that $g_C \geq 5$. So the above lemma tells us that there exists a pair of linearly independent global sections $e$ and $e'$
of ${\cal E}_{\xi}$ which are proportional and hence span a subsheaf of rank $1$ of ${\cal E}_{\xi}$. Call the saturation of that subsheaf ${\cal L}$. This gives us the following exact sequence
$$
\xymatrix{
	0\ar[r] & {\cal L} \ar[r]&{\cal E}_{\xi} \ar[r]& \OO_C (B)\ar[r]&0,
}
$$
where $h^0 ({\cal L}) \geq 2$.  Combining with our extension sequence
$$
\xymatrix{
	0\ar[r] & \OO_C \ar[r]&{\cal E}_{\xi} \ar[r]& \OO_C (K_C)\ar[r]&0
}
$$
gives rise to the diagram
\begin{equation}\label{diag}
\xymatrix{
	&&0\ar[d]&&\\
	&&{\cal L}\ar[d]\ar[dr]&&\\
	0\ar[r] & \OO_C \ar[r]\ar[dr]&{\cal E}_{\xi} \ar[r]\ar[d]& \OO_C (K_C)\ar[r]&0\\
	&&\OO_C (B)\ar[d]&&\\
	&&0&&
}
\end{equation}
where the composition arrows, the slanted arrows of the diagram, are nonzero. In particular, the divisor $B$ in the vertical sequence is effective and nonzero. Next lemma shows that both line bundles in the column of the diagram are quite special.

\begin{lem}\label{lem:L-B-g}
	The line bundles ${\cal L}$ and $\OO_C (B)$ of the vertical sequence in \eqref{diag} are subject to the following properties.
	
	\vspace{0.2cm}
	1). $\displaystyle{h^0 ({\cal L})=\HA (deg({\cal L})+1)}$.
	
	\vspace{0.2cm}
	2). $\displaystyle{h^0 ({\cal O}_C (B))=\HA (deg(B)+1)}$.
	
	\vspace{0.2cm}
	3). $g_C=h^0({\cal E}_{\xi})=h^0 ({\cal L}) + h^0 ({\cal O}_C (B))$.
\end{lem} 
\begin{pf}
	From the vertical sequence of the diagram \eqref{diag} we obtain
	\begin{equation}\label{B-L-g}
	g_C=h^0 ({\cal E}_{\xi}) \leq h^0 ({\cal L}) +h^0 ({\cal O}_C (B)).
\end{equation}
	 The Riemann-Roch for ${\cal L}$ tells us
	 $$
	h^0 ({\cal L}) -h^0 ({\cal O}_C (B))=deg({\cal L})-g_C +1.
	$$
	Solving for $h^0 ({\cal L})$ and substituting into the inequality above we obtain
	$$
g_C \leq 2h^0 ({\cal O}_C (B)) + deg({\cal L})-g_C +1.
$$
Hence the inequality
$$
2h^0 ({\cal O}_C (B)) \geq 2g_C -1 -deg({\cal L})=deg(B)+1.
$$
From Clifford inequality we have
$$
2h^0 ({\cal O}_C (B)) < deg(B)+2,
$$
where the inequality is strict because $B$ is neither $0$ or $K_C$ and $C$ is not hyperelliptic. Combining this with the previous inequality we deduce
the equality
$$
2h^0 ({\cal O}_C (B)) =deg(B)+1.
$$
This in turn implies that the inequality we started with in \eqref{B-L-g} must be equality. This proves assertions 2) and 3) of the lemma. The equality 1) follows immediately from those two. 	 
\end{pf}

\begin{lem}\label{lem:9}
	If $C$ has no $g^1_2$, $g^1_3$ and $g^2_5$, the line bundles 
	${\cal O}_C (B)$ and ${\cal L}$ in Lemma \ref{lem:L-B-g} have the degree at least $9$. In particular, $g_C \geq 10$.	
\end{lem}
\begin{pf}
	The assumption that $C$ has no $g^1_2$, $g^1_3$ and $g^2_5$ and the formula in Lemma \ref{lem:L-B-g} imply that the smallest possible value
	for $deg({\cal L})$ is $7$. Assume such a situation, then the line bundle
	${\cal L}$ must be very ample and hence embeds $C$ as a curve of degree $7$ in $\PP^3$.  The Castelnuovo bound on the genus of the space curves, see \cite{G-H}, tells us
	$$
	g_C \leq \frac{7^2 -1}{4} -7 +1=6.
	$$
	This and the lower bound $g_C \geq 5$ imply 
	$$
	g_C=5 \,\text{or} \, 6.
	$$
	From this and Lemma \ref{lem:L-B-g}, 3), we deduce
	
	$$
	h^0(\OO_C (B)) =1 \,\text{or} \, 2.
	$$
	The first possibility implies that $B$ is a point and $H^0 ({\cal L}) \cong W_{\xi}$ and then $B$ is a fixed point of the linear system $|W_{\xi}|$ and this contrary to our assumption. The second possibility and the formula in Lemma \ref{lem:L-B-g}, 2), tell us that the degree of $B$ is $3$ and hence the line bundle $\OO_C (B)$ gives $g^1_3$ on $C$. Again a contradiction. Thus we conclude that the smallest possible value of
	$deg({\cal L})$ is $9$.
	
	The same argument applied to $\OO_C (B)$ gives
	$$
	deg({\cal L}) \leq 3,
	$$
contrary to the lower bound $9$ obtained in the first part of the proof. This completes the proof of the lemma.
\end{pf}

Let $\phi \in W_{\xi}$ be chosen so that the zero divisor 
$Z_{\phi}=(\phi=0)$ is a reduced scheme of $2(g_C -1)$ points and it is mapped by $|W_{\xi}|$ bijectively into the set of points in $\PP(W^{\ast}_{\xi})$ in general position. We call such a $\phi$ {\it general}. From what we have learned about $|W_{\xi}|$, the points $[\phi]$ in $|W_{\xi}|$ corresponding to general $\phi \in W_{\xi}$ form a Zariski dense open subset.

\begin{lem}\label{lem:cphi-gen}
	Let $\phi \in W_{\xi}$ be general. Then the linear map
	$$
	\overline{\alpha^{(2)}}(\phi,\bullet): W_{\xi}/\CC \phi \longrightarrow \HKC
	$$
	is injective; the map 
	$\overline{\alpha^{(2)}}(\phi,\bullet)$ denotes the factorization of
	$$
	\alpha^{(2)}(\phi,\bullet): W_{\xi} \longrightarrow \HKC
	$$
	through the quotient space $W_{\xi}/\CC \phi$.
\end{lem}
\begin{pf}
	Assume $\overline{\alpha^{(2)}}(\phi,\bullet)$ fails to be injective. Then $\alpha^{(2)}(\phi,\phi')=0$ for some $\phi'$ linearly independent of $\phi$.
	 The cocycle relation \eqref{c-cocycle} becomes
	$$
	\phi' \alpha^{(2)}(\phi, \psi) -\phi \alpha^{(2)}(\phi', \psi)=0
	$$
	for all $\psi \in W_{\xi}$. Arguing as in the proof of Proposition \ref{pro:MaxN-gen}, we deduce
	$$
	\alpha^{(2)}(\phi, \psi)=\mu(\psi)\phi , \forall \psi\in W_{\xi},
	$$
	where $\mu$ is a linear function of $W_{\xi}$.
	This tells us that $ker(\mu)$ parametrizes the subspace of global sections of ${\cal E}_{\xi}$ which are proportional to $\alpha(\phi)$.
	Hence the diagram
	$$
	\xymatrix{
		&&0\ar[d]&&\\
		&&{\cal L} \ar[d] \ar[rd]&&\\
		0\ar[r] & \OO_C \ar[r]&{\cal E}_{\xi} \ar[r]\ar[d]& \OO_C (K_C)\ar[r]&0\\
		&&\OO_C (B) \ar[d]&&\\
		&&0&&
	}
	$$
	with the knowledge $h^0 ({\cal L}) = g_C-2$. This and Lemma \ref{lem:L-B-g} imply $deg(B) = 3$. From 3) of the same lemma we also deduce
	$$
	g_C =h^0 ({\cal E}_{\xi})  = h^0 ({\cal L}) + h^0 (\OO_C (B))=g_C -2 + h^0 (\OO_C (B)).
	$$
	Thus $ h^0 (\OO_C (B)) = 2$. But then $C$ has $g^1_3$ which is contrary to our assumptions.	
\end{pf}

Let $\overline{\alpha^{(2)}}(\phi,\bullet) :  W_{\xi}/\CC \phi \longrightarrow \HKC$ be as in Lemma \ref{lem:cphi-gen}. Composing with the projection
$$
\HKC \longrightarrow \HKC /W_{\xi} 
$$
gives the linear form
$$ 
l_{\phi} :W_{\xi}/\CC \phi \longrightarrow \HKC /W_{\xi} \cong \CC.
$$

\begin{lem}\label{lem:l-phi}
	For a general $\phi$ in $W_{\xi}$ the linear form $l_{\phi}$ is non zero.
\end{lem}
\begin{pf}
	Assume $l_{\phi}$ is zero. Then $\overline{\alpha^{(2)}}(\phi, \bullet)$ factors through $W_{\xi}$:
	$$
	\overline{\alpha^{(2)}}(\phi, \bullet) :  W_{\xi}/\CC \phi \longrightarrow W_{\xi}.
	$$
	Composing with the projection
	$$
	W_{\xi} \longrightarrow W_{\xi}/\CC \phi
	$$
	gives the endomorphism
	$$
	\{\overline{\alpha^{(2)}}\}(\phi, \bullet) :  W_{\xi}/\CC \phi \longrightarrow W_{\xi}/\CC \phi.
	$$
	Let $\{\phi_0\} \in W_{\xi}/\CC \phi$ be an eigenvector of $\{\overline{\alpha^{(2)}}\}(\phi, \bullet)$ with the eigenvalue $\lambda_0$. This gives the formula
	$$
	\alpha^{(2)}(\phi,\phi_0)=\lambda_0 \phi_0 + \mu_0 \phi
	$$
	which we already encountered in \eqref{phi0-form}; from this point on we
	repeat the argument there.
\end{pf}

We use the above construction to obtain a maximal ladder of $W_{\xi}/\CC \phi$.
\begin{lem}\label{lem:filt-phi}
	Let $\phi$ in $W_{\xi}$ be general. Then the vector space 
	$W_{\xi}/\CC \phi$ admits a filtration
	$$
	W_{\xi}/\CC \phi=F^{g_C -2}(\phi) \supset F^{g_C -3}(\phi) \supset \cdots \supset F^2(\phi) \supset  F^1 (\phi) \supset F^0(\phi) =0
	$$
	such that each quotient $F^p (\phi) /F^{p-1} (\phi)$, for $p\in [1,(g_C-2)]$, is one dimensional and the map
	$\overline{\alpha^{(2)}}(\phi, \bullet)$ takes $F^p (\phi)$ to $F^{p+1} (\phi)$ and induces isomorphisms
	$$
	 \overline{\alpha^{(2)}}_p(\phi, \bullet) : F^p (\phi) /F^{p-1} (\phi) \longrightarrow F^{p+1} (\phi) /F^{p} (\phi)
	 $$
	for each $p\in [1, (g_C -3)]$. 
\end{lem}
\begin{pf}
	We set $F^{g_C -2}(\phi):=W_{\xi}/\CC \phi$ and define the next step
	of the filtration by the formula
	$$
	F^{g_C -3}(\phi):=ker(l_{\phi}),
	$$
	where $l_{\phi}$ is the linear form in Lemma \ref{lem:l-phi}; since 
	$l_{\phi}$ is nonzero, this is a codimension $1$ subspace of $F^{g_C -2}(\phi)$. In addition, by definition of $l_{\phi}$ we have that
	$
	 \overline{\alpha^{(2)}}(\phi, \bullet)
	 $
	 maps 
	$F^{g_C -3}(\phi)$ to $F^{g_C -2}(\phi)$. This gives the linear form
	$$
	l^{g_C-3}_{\phi}: F^{g_C -3}(\phi) \longrightarrow F^{g_C -2}(\phi) / F^{g_C -3}(\phi).
	$$
	This is nonzero: otherwise we have an endomorphism
	$$
	\overline{\alpha^{(2)}}(\phi, \bullet): F^{g_C -3}(\phi) \longrightarrow F^{g_C -3}(\phi)
	$$
	and the argument with an eigenvector can be repeated. We now define
	$$
	F^{g_C -4}(\phi):=ker(l^{g_C-3}_{\phi})
	$$
	and obtain the map induced by $\overline{\alpha^{(2)}}(\phi, \bullet)$
	$$
	F^{g_C -3}(\phi)/	F^{g_C -4}(\phi)  \longrightarrow F^{g_C -2}(\phi)/	F^{g_C -3}(\phi).
	$$
	  Since this is nonzero and each quotient space is one dimensional, the construction gives an isomorphism. This completes the inductive step of the construction.
\end{pf}

{\it Proof of {\bf Theorem \ref{thm:rk1KS}}.}
We are going to use the filtration in Lemma \ref{lem:filt-phi} to write
down the quadratic equations for $C \subset \PP(\HKC^{\ast})$. For this
we choose a generator of the subspace $F^1(\phi)$ of that filtration and denote its lifting to $W_{\xi}$ by $\phi_1$. From the proof of Lemma \ref{lem:filt-phi} it follows that the string of vectors 
$$
\{{\phi}_1,\,{\phi}_2,\ldots,{\phi}_{g_C-1}\}
$$
defined inductively starting with $\phi_1$ by the condition
\begin{equation}\label{phip-string}
{\phi}_p={\alpha^{(2)}}(\phi,\phi_{p-1}),
\end{equation}
for every $p\in [2,{g_C-1}]$, is linear independent in $\HKC$. Furthermore,
all vectors with exception of the last one are in $W_{\xi}$:
$$
\{{\phi}_1,\,{\phi}_2,\ldots,{\phi}_{g_C-2}\} \subset W_{\xi} .
$$
For a triple
$\phi,\phi_i, \phi_j$, where $i,j \in [1, g_C-2]$, we have the cocycle relation 
$$
\phi_j \alpha^{(2)}(\phi,\phi_i)-\phi_i \alpha^{(2)}(\phi,\phi_j)+\phi \alpha^{(2)}(\phi_i,\phi_j)=0.
$$
Using the equation \eqref{phip-string} the above reads
$$
\phi_j \phi_{i+1}-\phi_i \phi_{j+1}+\phi \alpha^{(2)}(\phi_i,\phi_j)=0,
$$ 
for all $i, j\in [1,g_C-2]$. These are quadratic equations in $\PP(\HKC^{\ast})$ for $C$. The pairs $i<j$ in $[1,g_C-2]$ give the collection
\begin{equation}\label{Qij}
Q_{ij}=\phi_i \phi_{j+1}-\phi_j \phi_{i+1}-\phi \alpha^{(2)}(\phi_i,\phi_j)
\end{equation} 
of linear independent quadratic polynomials in $S^2 \HKC$ vanishing on $C$.
The number of those polynomials is
$$
\binom{g_C -2}{2}=dim(S^2 \HKC) -h^0 (\OO_C (2K_C)) =dim (I_C (2)).
$$
Hence 
$$
Q_{ij}=0, \,\, 1\leq i<j\leq g_C -2.
$$
is a complete set of quadratic equations for $C\subset \PP(\HKC^{\ast})$.
 To see what the subscheme defined by those quadrics is, we consider the restriction to the hyperplane
$H_{\phi}$ in $\PP(W^{\ast}_{\xi})$ corresponding to $\phi$. From \eqref{Qij} we deduce 
$$
Q_{ij}\Big|_{H_{\phi}}=\overline{\phi}_i \overline{\phi}_{j+1} - \overline{\phi}_j \overline{\phi}_{i+1} ,
$$
where $\overline{\phi}_k$ stands for the restriction of $\phi_k$ to $H_{\phi}$. Thus the hyperplane section
$Z_{\phi}=H_{\phi} \cdot C$ of $C$ lies on the subscheme defined by the vanishing of the $2\times 2$-minors of the matrix
$$
\begin{pmatrix}
	\overline{\phi}_1&\overline{\phi}_2&\cdots&\overline{\phi}_{g_C -2}\\
	\overline{\phi}_2&\overline{\phi}_3&\cdots&\overline{\phi}_{g_C -1}
\end{pmatrix}.
$$
It is well known that the variety determined by the vanishing of those minors is a rational normal
curve $R_{g_C-2}$ of degree $(g_C-2)$ in $H_{\phi} \cong \PP^{g_C-2}$, see \cite{FH}. As we said before, that rational normal curve is a hyperplane section of the subscheme $S$ in 
$\PP(\HKC^{\ast})$ defined by the quadratic polynomials in \eqref{Qij}. Hence $S$ is a surface of degree $(g_C-2)$ in $\PP(\HKC^{\ast}) \cong \PP^{g_C -1}$. This is a nondegenerate surface of the minimal degree and those are known to be either scrolls, or Veronese embedding of $\PP^2$ by the linear system of quadrics, or the projection from a point on the Veronese surface. The latter two cases imply that $(g_C -1)= 5$ or $4$. Hence $g_C=6$ or $5$. But in Lemma \ref{lem:9} we have learned that $g_C$ is at least $10$. So $S$ must be a rational normal scroll,
that is, $S$ is the image of a birational morphism of a ruled surface
\begin{equation}\label{scroll-pr}
\Sigma:=\PP(\FF) \longrightarrow \PP^1
\end{equation}
over $\PP^1$, where $\FF=\OO_{\PP^1} \oplus \OO_{\PP^1} (-s) $ is a rank two bundle on $\PP^1$ for some nonnegative integer $s$. More precisely, we have a morphism
$$
\Sigma \longrightarrow S \subset \PP(\HKC^{\ast})
$$
which is either an embedding or is a contraction of a section of $\Sigma$ to a point; in the latter case $S$ is a cone over a rational normal curve of degree $(g_C -2)$.

We denote by $h_0$ the section of $\Sigma$ of minimal self-intersection $(-s)$ and $l$ the class of a fibre of the structure projection in \eqref{scroll-pr}. Those form a basis of the N\'eron-Severi group of $\Sigma$. In particular, $C$ is the image
of a reduced irreducible divisor $\Gamma=dh_0 +ml$, for positive integers
$d$ and $m$. The integer $d$ is the degree of the intersection of a ruling of the scroll with $C$. It gives rise to a $g^1_d$ on our curve.
Thus 
\begin{equation}\label{datleast4}
	d\geq 4,
\end{equation}
 because we are assuming that $C$ has no $g^1_2$ or $g^1_3$.
 Let $D$ be the divisor on $C$ cut out by a ruling of the scroll. Since $D$ is not a point we have
 $$
 h^0 (\OO_C (K_C -D))=g_C -2.
 $$
 This and the Riemann-Roch for $\OO_C (D)$ tell us
 \begin{equation}\label{h0D}
 h^0 (\OO_C (D)) =h^0(\OO_C (K_C -D))+d-g_C+1=d-1.
\end{equation} 
 Combining with the Clifford inequality gives
 $$
 \HA d +1 \geq h^0 (\OO_C (D)) =d-1.
 $$
 Hence we obtain
 $$
 d\leq 4.
 $$
 This and the lower bound \eqref{datleast4} give the equality
 $$
 d=4.
 $$
  Then the Clifford inequality also must be equality. But this occurs only if $D=0, K_C$ or a multiple of
 $g^1_2$. Neither possibility can take place under our assumptions. 
  The proof of Theorem \ref{thm:rk1KS} is completed.
$\Box$

\vspace{0.2cm}
The Kodaira-Spencer classes of rank $1$, from the point of view of extensions, can be characterized as follows.
\begin{cor}\label{cor:rk1-ext}
	Let $\xi \in \HH$ be of rank $1$. Then the linear $|W_{\xi}|$ system has a unique base point, call it $p_{\xi}$, and in the extension sequence
	$$
	\xymatrix{
		0\ar[r] & \OO_C \ar[r]&{\cal E}_{\xi} \ar[r]& \OO_C (K_C)\ar[r]&0
	}
	$$
	the vector bundle ${\cal E}_{\xi}$ fits into the sequence
	$$
	\xymatrix{
		0\ar[r] & \OO_C (K_C -p_{\xi}) \ar[r]&{\cal E}_{\xi} \ar[r]& \OO_C (p_{\xi}) \ar[r]&0
	}
	$$
	and that sequence splits, that is we have
	$$
	{\cal E}_{\xi} \cong \OO_C (K_C -p_{\xi}) \oplus \OO_C (p_{\xi}).
	$$
	\end{cor}
\begin{pf}
	From Theorem \ref{thm:rk1KS} we know that the linear subsystem $|W_{\xi}|$ has a unique base point, call it $p_{\xi}$. Denote by
	$\delta$ a unique, up to a nonzero scalar, global section of $\OO_C (p_{\xi})$. Then we have
	$$
	W_{\xi}=\delta H^0(\OO_C (K_C -p_{\xi}))
	$$
	and we know 
	\begin{equation}\label{xidelattau=0}
	\xi (\delta \tau)=0, \,\forall \tau\in H^0 (\OO_C (K_C -p_{\xi})).
\end{equation}
	We claim that
	\begin{equation}\label{xideltap=0}
		\text{$\xi \delta =0$ in $H^1(\OO_C (-K_C +p_{\xi})$.}
	\end{equation}
Indeed, set $\eta:=\xi \delta$, then the equation \eqref{xidelattau=0} says
$$
0=\xi (\delta \tau)=(\xi \delta) \tau =\eta \tau, \, \forall \tau \in H^0 (\OO_C (K_C -p_{\xi})).
$$
This means that the cohomology class $\eta$ lies in the kernel of the cup-product
\begin{equation}\label{cp-pxi}
H^1(\OO_C (-K_C +p_{\xi}) \longrightarrow Hom (H^0 (\OO_C (K_C -p_{\xi})), H^1 (\OO_C)).
\end{equation}
The line bundle $\OO_C (K_C -p_{\xi})$ is subject to the assumptions of Proposition \ref{pro:MaxN-gen}: it is very ample because $C$ has no $g^1_3$ and we have
$$
deg(K_C -p_{\xi})=2g-3=2(g-1)-1=2h^0(\OO_C (K_C -p_{\xi}))-1.
$$
Hence according to Proposition \ref{pro:MaxN-gen} the cup-product \eqref{cp-pxi} is injective and the vanishing
$$
\xi \delta=0
$$ 
follows. This vanishing, from the point of view of the extension
$$
	\xymatrix{
	0\ar[r] & \OO_C \ar[r]&{\cal E}_{\xi} \ar[r]& \OO_C (K_C)\ar[r]&0,
}
$$
means the following: tensor the exact sequence with $\OO_C (p_{\xi}-K_C)$
to obtain the exact sequence of cohomology groups
$$
	\xymatrix{
	0\ar[r] &H^0({\cal E}_{\xi} (-(K_C -p_{\xi}))\ar[r]& H^0(\OO_C (p_{\xi}) )\ar[r]^(.42){\xi}&H^1(\OO_C (p_{\xi}-K_C));
}
$$
the global section $\delta$ of $\OO_C (p_{\xi})$ lies in the kernel of the coboundary map and hence comes from a nonzero global section of ${\cal E}_{\xi} (-(K_C -p_{\xi})$. That section is the same as a monomorphism
$$
\xymatrix{
	0\ar[r] & \OO_C (K_C -p_{\xi}) \ar[r]&{\cal E}_{\xi} 
}
$$
which we complete to an exact sequence
$$
\xymatrix{
	0\ar[r] & \OO_C (K_C -p_{\xi}) \ar[r]&{\cal E}_{\xi} \ar[r]& {\cal Q}\ar[r]&0.
}
$$
This together with the defining extension sequence give
 the diagram
	$$
		\xymatrix{
			&&0\ar[d]&&\\
			&&{\cal O} (K_C-p_{\xi})\ar[d]\ar[dr]&&\\
			0\ar[r] & \OO_C \ar[r]\ar[dr]&{\cal E}_{\xi} \ar[r]\ar[d]& \OO_C (K_C)\ar[r]&0\\
			&&{\cal Q}\ar[d]&&\\
			&&0&&
		}
$$
where the slanted arrows, the composition arrows, are nonzero. Since the horizontal sequence does not split, the quotient in the vertical sequence is torsion free and hence 
$$
{\cal Q}=\OO_C (p_{\xi}).
$$
The above diagram now reads
$$
	\xymatrix{
	&&0\ar[d]&&\\
	&&{\cal O}_C (K_C-p_{\xi})\ar[d]\ar[dr]^{\delta}&&\\
	0\ar[r] & \OO_C \ar[r]\ar[dr]_{\delta}&{\cal E}_{\xi} \ar[r]\ar[d]& \OO_C (K_C)\ar[r]&0\\
	&&{\cal O}_C(p_{\xi})\ar[d]&&\\
	&&0&&
}
$$
The vertical sequence of the diagram is the one stated in the corollary.
Furthermore, the epimorphism of the vertical sequence induces the surjective map on the level of global sections. This implies that the sequence splits:
from the identification
$$
Ext^1 (\OO_C (p_{\xi}), \OO_C (K_C -p_{\xi})) \cong H^1(\OO_C (K_C -2p_{\xi}))
$$
the vertical exact sequence of the diagram corresponds to a cohomology class in $H^1(\OO_C (K_C -2p_{\xi}))$, call that class $\lambda$; the surjectivity
$$
H^0(\EE_{\xi})\longrightarrow H^0 (\OO_C (p_{\xi}))
$$
from the vertical sequence is equivalent to saying that the coboundary map
$$
H^0 (\OO_C (p_{\xi})) \stackrel{\lambda}{\longrightarrow} H^1 (\OO_C (K_C -p_{\xi}))
$$
is zero; the coboundary map is the cup product with the class $\lambda$ so the above says that $\lambda$ lies in the kernel of the map
\begin{equation}\label{deltap}
H^1 (\OO_C (K_C -2p_{\xi})) \stackrel{\delta}{\longrightarrow} H^1 (\OO_C (K_C -p_{\xi}));
\end{equation}
 this map is dual to  
$$
H^0(\OO_C (p_{\xi})) \stackrel{\delta}{\longrightarrow} H^0(\OO_C (2p_{\xi})),
$$
the multiplication by $\delta$, which is an isomorphism since both spaces are one dimensional. 
\end{pf}
\begin{rem}\label{rem:Shiff}
	The claim in \eqref{xideltap=0} that 
	\begin{equation}\label{Schiff}
	\xi \delta=0
\end{equation}
	is interpreted classically as a Shiffer variation: write the exact sequence
	$$
	\xymatrix{
		0\ar[r] & \Theta_C \ar[r]^(.4){\delta}&\Theta_{C} (p_{\xi}) \ar[r]& \Theta_C \otimes \OO_{p_{\xi}} (p_{\xi})\ar[r]&0;
	}
	$$
	the associated exact sequence of the cohomology groups gives
	$$
		\xymatrix{
	0\ar[r] &H^0 (\Theta_C \otimes \OO_{p_{\xi}} (p_{\xi})) \ar[r]& H^1 (\Theta_{C}) \ar[r]^(.45){\delta}&H^1 (\Theta_{C} (p_{\xi}));
}
$$
the equation \eqref{Schiff} says that $\xi$ is in the image of the coboundary map of the above sequence: one of the definitions of Shiffer variations, see \cite{G}. Thus Corollary \ref{cor:rk1-ext} (re)proves Griffiths result stating that all Kodaira-Spencer classes of rank $1$ are Shiffer variations of $C$ provided the curve has no $g^1_2$, $g^1_3$, $g^2_5$.

Equivalently, for a point $p\in C$ the line
$$
\xymatrix{
	0\ar[r] &H^0 (\Theta_C \otimes \OO_{p} (p)) \ar[r]& H^1 (\Theta_{C}) =H^1 (\OO_C (-K_C)) \cong H^0 (\OO_C (2K_C)) ^{\ast}
}
$$
is identified with the image of $p$ under the bicanonical embedding of $C$.
So Corollary \ref{cor:rk1-ext} also means that the stratum
$$
\Sigma_1=\{[\xi]\in \PP(H^1 (\OO_C (-K_C)))| \, rk(\xi)=1\}
$$
in \eqref{rkstrat} is precisely the bicanonical image of $C$ under the same
provision of no $g^1_2$, $g^1_3$, $g^2_5$.
\end{rem}

The Enriques-Babbage-Petri theorem about $C \subset \PP(\HKC^{\ast})$ being cut out by quadrics now follows from the above results. This is proved in \cite{G}. For convenience of the reader we give a proof by essentially reproducing Griffiths' argument.
\begin{cor}
	Let $C$ be a smooth projective curve of genus $g$ with no $g^1_2$, $g^1_3$, $g^2_5$. Then the image of $C$ in $\PP(\HKC^{\ast})$ is cut out by quadrics
	in $\PP(\HKC^{\ast})$ passing through $C$.
\end{cor}
 \begin{pf}
 Set $W:=\HKC$ and	consider the Veronese embedding
 	$$
 	\nu: \PP(W^{\ast})\longrightarrow \PP(S^2 W^{\ast}).
 	$$
 Viewing $S^2 V^{\ast}=Hom^{s}(V,V^{\ast})$ as the space of symmetric linear maps from $V$ to $V^{\ast}$, the image of $\nu$, call it $V_2$, parametrizes maps of rank $1$:
 $$
 V_2=im(\nu)=\{ [v]\in \PP(Hom^{s}(V,V^{\ast}))\Big| rk(v)=1\}.
 $$
 	According to Theorem \ref{th-Max} the map
 	$$
 H^1 (\OO_C (-K_C))\cong	H^0 (\OO_C (2K_C))^{\ast} \longrightarrow S^2 V^{\ast}
 $$
 is injective. Hence the inclusion
 $$
 \PP(H^1 (\OO_C (-K_C)))=\PP(H^0 (\OO_C (2K_C))^{\ast}) \hookrightarrow \PP(S^2 V^{\ast}).
 $$
 The intersection of this subspace with $V_2$
 $$
 \PP(H^1 (\OO_C (-K_C))) \bigcap V_2 
 $$
 is precisely the stratum $\Sigma_1$ of Kodaira-Spencer classes of rank $1$:
 $$
 \Sigma_1=\PP(H^1 (\OO_C (-K_C))) \bigcap V_2. 
 $$
 From Remark \ref{rem:Shiff} we know that the stratum $\Sigma_1$ is the bicanonical image of $C$. Thus we have the commutative diagram
 $$
 \xymatrix{
 C \ar[r]^(.3){|2K_C|} \ar[d]_{|K_C|}& \PP(H^1 (\OO_C (-K_C)))\ar@{^{(}->}[d]\\
\PP(W^{\ast}) \ar^{\nu}[r]&\PP(S^2 W^{\ast})
}
$$
The inclusion on the right side of the diagram can be restated as follows
$$
\PP(H^1 (\OO_C (-K_C)))=\PP(I_C (2)^{\perp}),
$$
where $I_C (2)$ is the space of quadrics in $\PP(W^{\ast})$ through the canonical image of $C$. From the commutativity of the above diagram it follows
$$
\nu^{\ast} (\PP(I_C (2)^{\perp}))=\nu^{\ast}\Big(\PP(I_C (2)^{\perp})\bigcap V_2\Big)=\nu^{\ast}(\Sigma_1) =C,
$$
the canonical image of $C$ is the vanishing locus of quadrics  through $C \subset \PP(W^{\ast})$. 
 \end{pf}

Once we know the stratum $\Sigma_1$ of the stratification \eqref{rkstrat},
we can produce Kodaira-Spencer classes in the higher strata $\Sigma_k$, for $k \geq 2$. Indeed, take any $d$ points on $C$, for $d\in [1,g_C]$, and consider their linear span $L$ in the {\it bicanonical} embedding of $C$, then $L$ is contained in the stratum $\Sigma_d$. This is well known, see \cite{B}, \cite{G}; we give a proof for the convenience of the reader.
\begin{lem}\label{lem:Ld}
	Let $D$ be an effective divisor of degree $d \in [1,g_C]$ and let $L_D$
	be the linear span of $D$ in the bicanonical embedding of $C$, that is
	$$
	L_D:=\PP(H^0(\OO_C (2K_C-D))^{\perp}),
	$$ 
the intersection of all hyperplanes in $\PP(H^0(\OO_C (2K_C))^{\ast})$ passing through $D$. Then $L_D$	is contained in the stratum $\Sigma_d$.
\end{lem}
\begin{pf}
	Let $D$ be any effective divisor $C$ without any restriction on its degree. Take $\delta_D$, a global section of $\OO_C (D)$ defining $D$:
	$$
	D=(\delta_D=0).
	$$
	Consider the exact sequence
	$$
\xymatrix{
	0\ar[r] & \OO_C (-K_C)\ar[r]^{\delta_D}&\OO_C (D-K_C)\ar[r]& \OO_D (D-K_C)\ar[r]&0.
}
	$$
	On the level of global sections this gives
	$$
\xymatrix@C=10pt{
	0\ar[r] &H^0 (\OO_C (D-K_C))\ar[r]& H^0(\OO_D (D-K_C)) \ar[r]&H^1(\OO_C (-K_C))\ar[r]^{\delta_D}& H^1(\OO_C (D-K_C))\ar[r]&0.
}
$$
This implies that the Kodaira-Spencer classes $\xi$ subject to the condition
$$
\delta_D \xi =0
$$
are precisely the ones lying in the image of the coboundary map
$$
\xymatrix@C=10pt{
	 H^0(\OO_D (D-K_C)) \ar[r]&H^1(\OO_C (-K_C)).
	}
$$
By duality, that image consists of linear forms on $H^0 (\OO_C (2K_C))$ vanishing on the subspace $\delta_D H^0(\OO_C (K_C-D))$. Thus the identification
$$
\begin{gathered}
ker\big(H^1(\OO_C (-K_C))\stackrel{\delta_D}{\longrightarrow}H^1(\OO_C (D-K_C))\big)=H^0(\OO_C (2K_C-D))^{\perp}\\
=\big(H^0 (\OO_C (2K_C))/ \delta_D H^0(\OO_C (2K_C-D))\big)^{\ast}.
\end{gathered}
$$
Geometrically, the projectivization of the vector space $H^0(\OO_C (2K_C-D))^{\perp}$ is
the linear span $L_D$ of the subscheme $D\subset C$ in the bicanonical embedding of $C$.

From the point of view of the extension
$$
\xymatrix{
	0\ar[r]&\OO_C \ar[r]& \EE_{\xi}\ar[r]&\OO_C (K_C)\ar[r]&0
}
$$
corresponding to 
 a nonzero Kodaira-Spencer class $\xi$, the cohomological condition
$$
\delta_D\xi =0
$$
means that $\EE_{\xi} (-(K_C -D))$ has a nonzero global section: tensor the extension sequence with $\OO_C (-(K_C -D))$ to obtain the exact sequence
$$
\xymatrix{
	 H^0(\EE_{\xi}(-(K_C-D)))\ar[r]&H^0(\OO_C (D))\ar[r]^(.43){\xi}&H^1(\OO_C (D-K_C));
}
$$
	the cohomological condition above says that the global section $\delta_D$ comes from a nonzero global section of $\EE_{\xi} (-(K_C -D))$.
	
	A nonzero global section of $\EE_{\xi} (-(K_C -D))$ is viewed as a monomorphism
	$$
	\xymatrix@R=12pt{
		0\ar[r]&\OO_C (K_C -D) \ar[r] &\EE_{\xi}.
	}
		$$
Putting it together with the extension sequence gives the commutative diagram

\begin{equation}\label{deltaD-diag}
\xymatrix{
	&&0\ar[d]&&\\
	&&\OO_C (K_C -D) \ar[d] \ar[dr]^{\delta_D}&&\\
	0\ar[r]&\OO_C \ar[r]& \EE_{\xi}\ar[r]\ar[d]&\OO_C (K_C)\ar[r]&0\\
	&&{\cal Q}\ar[d]&&\\
	&&0&&
}
\end{equation}
From this follows the inclusion
$$
\delta_D H^0 (\OO_C (K_C -D)) \subset W_{\xi}
$$
and an estimate on the rank of $\xi$
\begin{equation}\label{rkxi-estim}
	\begin{gathered}
rk(\xi)=g-dim(W_{\xi}) \leq g-h^0(\OO_C (K_C -D))\\
=deg(D)-h^0(\OO_C (D))+1=d-h^0(\OO_C (D))+1,
\end{gathered}
\end{equation}
where the first equality is the Riemann-Roch for $\OO_C(D)$. Hence
the inequality
$$
rk(\xi) \leq d,
$$
for all $\xi$ lying in the subspace $H^0(\OO_C (2K_C-D))^{\perp}$. In particular, for $d\in [1,g_C]$, the inclusion
$$
L_D=\PP(H^0(\OO_C (2K_C-D))^{\perp}) \subset \Sigma_d
$$
asserted in the lemma follows.
\end{pf} 
\begin{rem}\label{rem:D-D'}
1)	The diagram \eqref{deltaD-diag} encapsulates the relation between the secant spaces of the bicanonical embedding of $C$ and the extensions parametrized by secant spaces. For example, one easily sees that the quotient sheaf ${\cal Q}$ in that diagram is locally free if and only if the point $[\xi]$ corresponding to the Kodaira-Spencer class $\xi$ lies in the linear span $L_D$ but not in the subspace $L_{D'}$, the linear span of any proper subscheme $D'$ of $D$.
	
2)	Let $L_D$ be the linear span of an effective divisor $D$ of degree $d \in [1,g_C]$ in the bicanonical embedding of $C$. From Lemma \ref{lem:Ld} we know
$$
L_D \subset \Sigma_d.
$$
If, in addition, $L_D$ is not contained in the smaller rank stratum $\Sigma_{d-1}$, then the proof of the lemma implies that the intersection
	$$
	L_D \bigcap \Big(\Sigma_d \setminus \Sigma_{d-1}\Big)
	$$
	parametrizes projectivized Kodaira-Spencer classes $[\xi]$ with
	$$
	\text{$W_{\xi}=\delta_D H^0(\OO_C (K_C -D))$ and $H^0 (\OO_C (D))=\CC \delta_D.$}
	$$
	In other words, the linear span $\{D\}_{K_C}$ of $D$ in the {\rm canonical} embedding of $C$ is a $\PP^{d-1}$-plane intersecting $C\subset \PP(\HKC^{\ast})$ along $D$.
	
	3) Let $D_{\phi}=(\phi=0)$ be the zero divisor of a nonzero global section $\phi$ of $\OO_C(K_C)$. The linear span $L_{D_{\phi}}$ of $D$ in the bicanonical embedding is a subspace of dimension $(2g-4)$ and that linear span parametrizes the Kodaira-Spencer classes $\xi$ subject to
	$$
	\xi \phi=0.
	$$
	From the diagram \eqref{deltaD-diag} we see that the intersection of
	$L_{D_{\phi}}$ with the bicanonical curve is precisely the subscheme $D_{\phi}$:
	$$
	D_{\phi}=C_{2K_C} \bigcap L_{D_{\phi}}.
	$$
	From the estimate in \eqref{rkxi-estim} it follows that $L_{D_{\phi}}$ lies in the stratum $\Sigma_{g-1}$. As $D_{\phi}$ moves in the canonical linear system, the linear subspaces $L_{D_{\phi}}$ sweep out
	$(3g-5)$ dimensional subvariety contained in $\Sigma_{g-1}$. We will see that this is precisely $\Sigma_{g-1}$, see Proposition \ref{pro:Sigma(g-1)}.   
\end{rem}

From Lemma \ref{lem:Ld} we deduce that the higher strata of \eqref{rkstrat} are related to the secant varieties of $C\subset \PP(H^0(\OO_C (2K_C))^{\ast})$. Namely, for an integer $k\geq 1$ denote by $Sec_k(C_{2K_C})$ the subvariety of $\PP(H^0(\OO_C (2K_C))^{\ast})$ swept out by $(k-1)$-planes in $\PP(H^0(\OO_C (2K_C))^{\ast})$ intersecting $C $ along a subscheme of degree at least $k$. Then the above lemma implies the following.

\begin{pro}\label{pro:SkinSigmak}
	Let $C$ be a smooth, complex projective curve of genus $g_C$ with the canonical bundle
	$\OO_C (K_C)$ very ample. Then for every $k \in [1,g_C]$ we have the inclusion
	$$
	Sec_k (C_{2K_C})\subset \Sigma_k.
	$$
\end{pro}

For $k=1$, the variety $Sec_1(C_{2K_C})$ is the bicanonical image of $C$ and Theorem \ref{thm:rk1KS} says that we have the equality
$$
Sec_1(C_{2K_C})=\Sigma_1,
$$
provided our curve is neither trigonal, nor a plane quintic. Such equalities can be extended to other values of $k$, if the curve $C$ is general in the moduli space $\mathfrak{M}_g$ of curves of genus $g$. This
will be the subject treated in the sequel to this paper. For now we give an illustration of the fact that even for general curves in the moduli the higher rank strata become more involved and intimately related to the geometry of rank $2$ bundles on the curve.
\begin{example}\label{ex:g=5}
The smallest value of the genus $g$ for a curve to satisfy the conditions of Theorem \ref{thm:rk1KS} is five. Our example concerns a smooth projective curve $C$ of genus $g=5$ with no $g^1_2$, $g^1_3$, $g^2_5$.
Thus according to Theorem \ref{thm:rk1KS} the stratum
$\Sigma_1$ is the bicanonical image of $C$. We wish to understand the next stratum $\Sigma_2$. From Proposition \ref{pro:SkinSigmak} we know the inclusion
$$
\Sigma_2 \supset Sec_2=Sec_2(C_{2K_C}),
$$
where $Sec_2$ is the secant variety of $\Sigma_1$, the bicanonical image of $C$. We will see in a moment that $\Sigma_{2}$ is much larger than the secant variety $Sec_2$ and geometrically it is related to quadrics passing through $C$. We begin by recalling that the space of quadrics $I_C (2)$ passing
through $C$ in the canonical embedding
$$
I_C (2):=ker \left(S^2 \HKC \longrightarrow H^0(\OO_C (2K_C))\right)
$$
is three dimensional, and $C \subset \PP(\HKC^{\ast})=\PP^4$ is a complete intersection of any three quadrics spanning that space. Conversely, a general plane in $\PP(S^2 \HKC)$ gives rise to a smooth complete intersection in $\PP^4$ which is a canonical curve of genus $5$. 

A general plane  in $\PP(S^2 \HKC)$ intersects the subvariety of singular quadrics along a quintic curve. Furthermore, since the subvariety parametrizing quadrics of rank at most $3$ is of dimension $11$, a general plane in $\PP(S^2 \HKC)$ intersects the locus of singular quadrics along the quintic curve which parametrizes quadrics of rank precisely $4$. From now on we assume that the plane $\PP(I_C (2)) $ satisfies this genericity assumption: 
\begin{equation}\label{g=5genassump}
	\text{\rm $\PP(I_C (2)) $ contains no quadrics of rank $3$,}
\end{equation}
 that is $\PP(I_C (2)) $ intersects the locus of singular quadrics along a curve of degree $5$, call that curve $\Gamma_C$, and that curve parametrizes singular quadrics of rank precisely $4$ passing through $C$. This means that
$\Gamma_C$ is a smooth irreducible plane curve of degree $5$.

The quadrics $Q_{\gamma}$ parametrized by closed points $\gamma \in \Gamma_C$ are cones over smooth quadrics in $\PP^3$. So they are ruled by two distinct pencils of $\PP^2$'s passing through the vertex $p_{\gamma}$ of $Q_{\gamma}$. Each of those $\PP^2$'s intersect $C$ along a divisor of degree $4$. Hence singular quadrics $Q_{\gamma}$ give rise to two {\rm distinct} $g^1_4$ on $C$. Thus we have the unramified double cover $\widetilde{\Gamma}_C$ of $\Gamma_C$ parametrizing distinct $g^1_4$ on $C$. Hence the inclusion
$$
\widetilde{\Gamma}_C \subset W^1_4 (C),
$$
where $W^1_4 (C)$ is the scheme parametrizing $g^1_4$ on $C$. For $C$ general, the dimension of $W^r_d (C)$, the scheme parametrizing $g^r_d$'s on a curve of genus $g$ is given by the Brill-Noether number, see \cite{ACGH},
$$
dim(W^r_d)=g_C -(r+1)(r+g_C -d);
$$
hence in the case at hand
$$
dim(W^1_4(C))=1
$$
and $W^1_4(C)$ is smooth and connected, see \cite{FL}. Hence {\rm all} $g^1_4$ on $C$ are parametrized by $\widetilde{\Gamma}_C$:
$$
\widetilde{\Gamma}_C=W^1_4 (C).
$$

We now proceed with connecting the above facts with the rank two bundles of the extension construction.
Let $[\xi]$ be a closed point in $\Sigma_2 \setminus \Sigma_1$. In particular,  the rank of $\xi$ is precisely $2$ or, equivalently, the kernel $W_{\xi}$ of $\xi$ is a subspace of $\HKC$ of dimension $3$. The extension
$$
\xymatrix{
	0\ar[r]&\OO_C \ar[r]&\EE_{\xi}\ar[r]& \OO_C (K_C)\ar[r]&0
}
$$
gives us the rank $2$ bundle $\EE_{\xi}$ with the space of global sections
fitting into the following exact sequence
$$
\xymatrix{
	0\ar[r]&H^0(\OO_C) \ar[r]&H^0(\EE_{\xi})\ar[r]& W_{\xi}\ar[r]&0
}
$$
In particular, $H^0(\EE_{\xi})$ is of dimension $4$. From the map
$$
\xymatrix{
	\bigwedge^2 H^0(\EE_{\xi}) \ar[r]& \HKC
}
	$$
	we learn that two situations may arise.
	
	{\rm Situation 1.} There are no decomposable tensors in the kernel of the above map and then the map is surjective; in particular, $\EE_{\xi}$ is generated by its global sections and fits into the following exact sequence
	$$
	\xymatrix{
		0\ar[r]&\EE^{\ast}_{\xi} \ar[r]&H^0(\EE_{\xi})\otimes \OO_C\ar[r]& \EE_{\xi}\ar[r]&0.
}
$$
Geometrically, the vector bundle $\EE_{\xi}$ can be interpreted through its Gauss map
$$
\xymatrix{
C\ar[r]& {\bf Gr}(1,\PP(H^0(\EE_{\xi})^{\ast})) \ar@{^{(}->}[r] &\PP(\bigwedge^2 H^0(\EE_{\xi})^{\ast}),
}
$$
where the inclusion arrow is the Pl\"ucker embedding of the Grassmannian;
this is related to the canonical embedding of $C$ by the following diagram
\begin{equation}\label{Pl-K-ex}
\xymatrix{
	C\ar[r] \ar[rd]& {\bf Gr}(1,\PP(H^0(\EE_{\xi})^{\ast})) \ar@{^{(}->}[r] &\PP(\bigwedge^2 H^0(\EE_{\xi})^{\ast})\\
	&\PP(\HKC^{\ast}) \ar[ru]&
}
\end{equation}
where the slanted arrow on the right is the dual of the surjective map
\begin{equation}\label{wedge2HKC-ex}
\xymatrix{
	\bigwedge^2 H^0(\EE_{\xi}) \ar[r]& \HKC.
}
\end{equation}
Thus  $\PP(\HKC^{\ast}) \subset \PP(\bigwedge^2 H^0(\EE_{\xi})^{\ast})$ gives a distinguished hyperplane in the Pl\"ucker embedding of the Grassmannian ${\bf Gr}(1,\PP(H^0(\EE_{\xi})^{\ast}))$. Since the latter is a quadric we obtain
$$
C\subset \PP(\HKC^{\ast}) \bigcap {\bf Gr}(1,\PP(H^0(\EE_{\xi})^{\ast})) \subset \PP(\HKC^{\ast})
$$
a distinguished quadric, call it $Q_{\xi}$ containing $C$. This quadric is nonsingular, since it corresponds to a generator of the kernel of the map \eqref{wedge2HKC-ex} which is a tensor of rank $4$. Thus geometrically, the rank $2$ bundle $\EE_{\xi}$ with no decomposable tensors in the kernel of the map \eqref{wedge2HKC-ex} corresponds to a smooth quadric through $C$ in the canonical embedding.
 Hence the Zariski open set 
$$
U_C:=\PP(I_C (2)) \setminus \Gamma_C
$$
 parametrizes a family of rank $2$ bundles as described above. 
Over $C\times U_C$ we have the universal bundle ${\bf E}$. Its direct image $\pi_{U_C \ast} ({\bf E})$ under the projection
$$
\pi_{U_C} :C\times U_C \longrightarrow U_C
$$
is a rank $4$ bundle on $U_C$. The projectivization
$$
\pi:\PP(\pi_{U_C \ast} ({\bf E})) \longrightarrow U_C 
$$
parametrizes pairs $(\{\EE\}, [e])$, where $\{\EE\}$ is a point of $U_C$ corresponding the vector bundle $\EE$ and $[e]\in \PP(H^0(\EE))$.

We notice that the $\PP^3$-bundle $\PP(\pi_{U_C \ast} ({\bf E}))$ comes with a distinguished divisor denoted ${\bf D}_C$. Set-theoretically, that divisor can be described as follows:
$$
{\bf D}_C=\{(\{\EE\},[e])\in \PP(\pi_{U_C \ast} ({\bf E})) | \, \text{\rm the zero locus $(e=0)$ is nonempty}\}.
$$
This is another geometric interpretation of rank $2$ bundles $\EE$ parametrized by $U_C$. Namely, the projectivization $\PP(\EE^{\ast})$ is mapped into
$\PP(H^0(\EE)^{\ast})$ by the tautological line bundle $\OO_{\PP(\EE^{\ast})} (1)$:
$$
\PP(\EE^{\ast}) \longrightarrow \PP(H^0(\EE)^{\ast})\cong \PP^3.
$$
The morphism is birational onto its image which is a surface $S_{\EE}$ of degree $8$ in $\PP^3$. The birationality follows from the genericity assumption - no quadrics of rank $3$ in $\PP(I_C(2))$, see \eqref{g=5genassump}. The dual surface $\Check{S}_{\EE}$ in $\PP(H^0(\EE))$ is the fibre of the divisor ${\bf D}_C$ over the point $\{\EE\}$ of $U_C$, that is, we have the formula
$$
\Check{S}_{\EE}={\bf D}_C \bigcap \pi^{-1} (\{\EE\})={\bf D}_C \bigcap \PP(H^0(\EE)):
$$
it parametrizes nonzero global sections $e$ of $\EE$ whose scheme of zeros
$$
Z_e=(e=0)
$$
is nonempty. Thus the complement
$$
Y_C:=\PP(\pi_{U_C \ast} ({\bf E})) \setminus {\bf D}_C 
$$
is a Zariski dense open subset parametrizing pairs $(\{\EE\}, [e]) \in \PP(\pi_{U_C \ast} ({\bf E}))$ with
$e$ nowhere vanishing.
This is fibred over $U_C$
$$
\pi:Y_C \longrightarrow U_C
$$
with the fibre $\pi^{-1}(\{\EE\})$ over $\{\EE\} \in U_C$ being the Zariski dense open subset of
$\PP(H^0(\EE))$ of nowhere vanishing global sections of $\EE$, the complement of the dual surface $\Check{S}_{\EE}$.
Thus it can be identified with the Koszul sequence of $(\EE,e)$
$$
\xymatrix{
0\ar[r]&\OO_C \ar^{e}[r]& \EE \ar^(.38){\wedge e}[r]&\OO_C (K_C)\ar[r]&0
}
$$
up to the $\CC^{\ast}$-action (of multiplying the arrows by a scalar).
This gives the inclusion
$$
Y_C \hookrightarrow \PP(H^1(\OO_C (-K_C)))
$$
which sends $(\{\EE\}, [e])$ to the projectivized Kodaira-Spencer class corresponding to the Koszul sequence above. Furthermore, by construction the image of the above morphism is in $\Sigma_2$. Hence we deduce
$$
\text{\rm $Y_C$ is  a substratum of $\Sigma_{2}$ of dimension $5$}.
$$
Since the secant variety $Sec_2$ has dimension $3$ we deduce that
$$
\Sigma_2 \setminus Sec_2
$$
is nonempty. In fact, we will see in a moment that $Y_C$ is disjoint from $Sec_2$, that is it is contained in $\Sigma_2 \setminus Sec_2$.

\vspace{0.2cm}
 {\rm Situation 2.}	We assume that the kernel of the map
 \begin{equation}\label{wedgemap-xi}
 	\xymatrix{
 		\bigwedge^2 H^0(\EE_{\xi}) \ar[r]& \HKC
 	}
\end{equation}
contains a decomposable tensor.	 It turns out that the geometric interpretation of 
$\EE_{\xi}$ in this case is related to $g^1_4$ and hence, as we explained at the beginning of the example, to singular quadrics through the canonical image of $C$, the ones parametrized by the curve $\Gamma_C \subset \PP(I_C(2))$. 

Since the kernel of the map in \eqref{wedgemap-xi} contains a decomposable tensor, say $e\wedge e'$, we have two linearly independent global sections $e$ and $e'$ of $\EE_{\xi}$ which are proportional, that is, they generate a subsheaf of rank one and give rise to the diagram

\begin{equation}\label{Exi-D-diag-ex}
\xymatrix{
	&&0\ar[d]&&\\
	&&{\cal O}_C (K_C -D)\ar[d]\ar[dr]^{\delta}&&\\
	0\ar[r] & \OO_C \ar[r]\ar[dr]_{\delta}&{\cal E}_{\xi} \ar[r]\ar[d]& \OO_C (K_C)\ar[r]&0\\
	&&\OO_C (D)\ar[d]&&\\
	&&0&&
}
\end{equation}
where $h^0 ({\cal O}_C (K_C -D))\geq 2$ and $\delta \in H^0 (\OO_C (D))$, see the argument in the paragraph following Lemma \ref{lem:Exi-glgen}.
If $h^0 ({\cal O}_C (K_C -D))> 2$, then
$$
h^0 ({\cal O}_C (K_C -D)) =3,
$$
since otherwise $[\xi] \in \Sigma_1$, and then the divisor $D_0=(\delta=0)$ is on the secant line of the {\rm canonical} curve. Since we are assuming that $C$ has no $g^1_3$, we deduce that
$deg(D)=2$ and $[\xi]$ lies in the secant variety $Sec_2$ of the bicanonical curve. 
Thus if $[\xi] \in \Sigma_2 \setminus Sec_2$, then
\begin{equation}\label{2-ex}
h^0 ({\cal O}_C (K_C -D))= 2,
\end{equation}
This in turn implies that 
$h^0 ({\cal O}_C (D))= 2$ as well and the degree of $D$ (resp. $(K_C-D)$) is four.
Thus the rank $2$ bundles $\EE_{\xi}$, for $[\xi] \in \Sigma_2 \setminus Sec_2$ and having a decomposable tensor in the kernel of the map \eqref{wedge2HKC-ex} gives rise to $g^1_4$ and, as discussed in the beginning of the example, this is related to singular quadrics through the canonical image of $C$; later on we will see how to write down equations of those quadrics from an appropriate parametrization of global sections of $\EE_{\xi}$. 

Observe also that the vertical sequence in \eqref{Exi-D-diag-ex} produces two $g^1_4$, one corresponding to $\OO_C(K_C-D) $ and the other to $\OO_C (D)$. From the genericity assumption
of no quadrics of rank $3$ in $\PP(I_C(2))$ it follows that two line bundles are distinct, or equivalently, 
\begin{equation}\label{nontrivial}
\OO_C(K_C -2D) \neq \OO_C.
\end{equation}

From the vertical sequence in the diagram \eqref{Exi-D-diag-ex} the vector bundle $\EE_{\xi}$ is defined by a cohomology classes in
$
H^1 (\OO_C (K_C -2D))
$
lying in the kernel of the cup-product map
\begin{equation}\label{Dcupprod-ex}
H^1 (\OO_C (K_C -2D)) \longrightarrow H^0(\OO_C (D))^{\ast} \otimes H^1(\OO_C (K_C-D)).
\end{equation}
The kernel is one dimensional: follows immediately from the dual of the above map and \eqref{nontrivial}. This gives us two possibilities: 

$\bullet$ either the vertical sequence in \eqref{Exi-D-diag-ex} is nonsplit or, equivalently, the decomposable tensor we fixed to construct that sequence is unique up to a nonzero scalar, 

$\bullet$ or $\EE_{\xi}$ is a direct sum
$$
\EE_{\xi} \cong \OO_C (D)\oplus \OO_C (K_C-D),
$$
equivalently, the kernel of the map 
\begin{equation}\label{wedge2-dec-ex}
\xymatrix{
	\bigwedge^2 H^0(\EE_{\xi}) \ar[r]& \HKC
}
\end{equation}
is two dimensional and spanned by {\rm two} decomposable tensors.

The first possibility gives us 
 the variety of pairs
$$
Y^{nsplit}_C=\{(\{\EE\}_{\widetilde{\gamma}},[e])\}_{\widetilde{\gamma} \in \widetilde{\Gamma}_C},
$$
 where $\{\EE\}_{\widetilde{\gamma}}$ is determined by a nonsplit vertical sequence in \eqref{Exi-D-diag-ex} corresponding to a generator of the kernel of the cup-product in \eqref{Dcupprod-ex} and $e$ is a nowhere vanishing global section of $\EE$. It is fibred over the curve $\widetilde{\Gamma}_C$
 $$
 \pi':Y^{nsplit}_C \longrightarrow \widetilde{\Gamma}_C
 $$
 with fibres Zariski dense open subsets of $\PP^3$, and we have the inclusion
 $$
 Y^{nsplit}_C \hookrightarrow \Sigma_{2}\setminus Y_C.
 $$

The second possibility gives the variety of pairs
 $$
 Y^{split}_C=\{(\{\EE\}_{\gamma},[e])\}_{\gamma \in \Gamma_C}
 $$
with $\EE_{\gamma}\cong \OO_C(D)\oplus \OO_C (D')$, where $D$ and $D'$ are two $g^1_4$ over a point $\gamma \in \Gamma_C$.

From the point of view of vector bundle theory, we can say that $\EE_{\xi}$ in $Y^{nsplit}_C$ are strictly semistable, indecomposable, while in $Y^{split}_C$ we have strictly semistable, decomposable bundles; in {\rm Situation 1}, the vector bundles $\EE_{\xi}$ parametrized by $Y_C$ are all stable. To complete the picture, the vector bundles $\EE_{\xi}$ parametrized by the points in $Sec_2\setminus \Sigma_1$ are all unstable with the destabilizing sequence
$$
\xymatrix{
		0\ar[r] & \OO_C (K_C-p-p')\ar[r]&{\cal E}_{\xi} \ar[r]& \OO_C (p+p')\ar[r]&0,
} 
$$
where $(p+p')$ is the divisor cut out by a secant line of $C\subset \PP(\HKC^{\ast})$, see Lemma \ref{lem:Ld}. Thus we obtain the decomposition of $\Sigma_2$ into the  union of mutually disjoint strata
$$
\Sigma_2 =Y_C \coprod Y^{nsplit}_C \coprod Y^{split}_C \coprod Sec_2.
$$
\end{example}

\vspace{0.2cm}
The proof of Theorem \ref{thm:rk1KS} and the example above indicate that the study of Griffiths stratification \eqref{rkstrat} ties together several strands in the theory of curves: geometry of canonical and bicanonical embeddings of a curve, the Brill-Noether theory and rank $2$ vector bundles and their moduli. 
We will now put on the general footing several aspects of the proof of Theorem \ref{thm:rk1KS}. This will lead to our refinement invariants
of the stratification \eqref{rkstrat} as discussed in the introduction. 

\section{A Refinement of IVHS invariants}
Throughout this section we assume that $C$ is a smooth projective curve of genus $g$ an it has no $g^1_2$, $g^1_3$ and $g^2_5$. In other words, our curve is as in Theorem \ref{thm:rk1KS}. Then according to that theorem in the stratification
	\begin{equation}\label{rkstrat-again}
\PP(H^1(\Theta_C)) =\Sigma_g \supset \Sigma_{g-1} \supset \cdots \supset \Sigma_1 \supset \Sigma_0
\end{equation}
the stratum $\Sigma_0$ is empty and $\Sigma_1$ is the bicanonical image of 
$C$ in $\PP(H^1(\Theta_C))$. We now turn to the study of the strata of higher rank: fix $r\geq 2$
 and consider the open stratum
 $$
\Sigma^0_r:=\Sigma_r \setminus \Sigma_{r-1}.
$$ 
By definition this parametrizes projectivized Kodaira-Spencer classes $[\xi]$ for which the cup product with $\xi$
$$
\xi :\HKC \longrightarrow H^1(\OO_{C})
$$
has rank precisely $r$. Equivalently, the dimension of the kernel 
$$
W_{\xi}=ker(\xi)
$$
is constant and equals $(g-r)$, for all $[\xi]$ in $\Sigma^0_r$. This means that the stratum $\Sigma^0_r$ is a parameter space of $(g-r)$-planes in $\HKC$. Hence it comes with the classifying morphism
\begin{equation}\label{kappa-r}
\kappa_{g-r}:\Sigma^0_r \longrightarrow {\bf Gr}((g-r), \HKC)
\end{equation}
into the Grassmannian ${\bf Gr}((g-r), \HKC)$ of $(g-r)$-planes in $\HKC$: that map sends a closed point $[\xi]$ in $\Sigma^0_r$ to the point $[W_{\xi}]$ of the Grassmannian corresponding to the subspace $W_{\xi} \subset \HKC$. To make the discussion geometrically meaningful we assume
 $$
 g-r \geq 3.
 $$
 Thus we will be working under the assumption
 \begin{equation}\label{r-assump}
 	2\leq r \leq g-3.
 \end{equation}

The pull back of the universal bundle of the Grassmannian ${\bf Gr}((g-r), \HKC)$ equips $\Sigma^0_r$ with a locally free sheaf denoted ${\cal W}_{\Sigma^0_r }$. This can also be defined as follows. 
 
 On $\PP(H^1(\Theta_C))$ we have the morphism of sheaves
 $$
 {\mathfrak{c}}:\HKC \otimes \OO_{\PP(H^1(\Theta_C))} (-1) \longrightarrow H^1(\OO_{C})\otimes \OO_{\PP(H^1(\Theta_C))},
 		$$
 		the sheaf version of the cup product
 		$$
 		H^1(\Theta_C) \longrightarrow Hom(\HKC,H^1(\OO_{C})).
 		$$
 Observe: the stratification \eqref{rkstrat-again} corresponds to the degeneracy loci of the morphism ${\mathfrak{c}}$. In particular,
 the stratum $\Sigma^0_r$ is the subscheme of $\PP(H^1(\Theta_C))$, where the cokernel of ${\mathfrak{c}}$ has rank precisely $(g-r)$. Namely, denote
 by ${\mathfrak{c}}_{\Sigma^0_r}$ the restriction of ${\mathfrak{c}}$ to $\Sigma^0_r$, then we obtain the exact complex of locally free sheaves on $\Sigma^0_r$:  
 $$
\xymatrix@C=12pt{
0\ar[r]& ker({\mathfrak{c}}_{\Sigma^0_r}) \ar[r]&\HKC \otimes \OO_{\Sigma^0_r} (-1) \ar[r]^(.6){{\mathfrak{c}}_{\Sigma^0_r}} &H^1(\OO_{C})\otimes \OO_{\Sigma^0_r} \ar[r]&coker({\mathfrak{c}}_{\Sigma^0_r})\ar[r]&0. 		
}
$$
The sheaf ${\cal W}_{\Sigma^0_r }$ defined previously is identified with
the twist of $ker({\mathfrak{c}}_{\Sigma^0_r})$:
$$
{\cal W}_{\Sigma^0_r } =ker({\mathfrak{c}}_{\Sigma^0_r})(1):=ker({\mathfrak{c}}_{\Sigma^0_r}) \otimes \OO_{\Sigma^0_r} (1),
$$
while the sequence
\begin{equation}\label{pullback-univseq}
\xymatrix{
	0\ar[r]& ker({\mathfrak{c}}_{\Sigma^0_r})(1) \ar[r]&\HKC \otimes \OO_{\Sigma^0_r} \ar[r]^(.6){{\mathfrak{c}}_{\Sigma^0_r}}  \ar[r]&im({\mathfrak{c}}_{\Sigma^0_r}) (1)\ar[r]&0 		
}
\end{equation}
is the pull back of the universal sequence on the Grassmannian ${\bf Gr}((g-r), \HKC)$ under the classifying map $\kappa_{g-r}$ in \eqref{kappa-r}. All of the above is the IVHS construction of Griffiths in \cite{G} attached to the differential of the period map for curves. We will now produce a refinement of this construction. It will take place on the projectivization $\PP({\cal W}_{\Sigma^0_r })$ of the sheaf ${\cal W}_{\Sigma^0_r }$. In down to earth terms this projectivization is the incidence correspondence
$$
\PP({\cal W}_{\Sigma^0_r })=\{([\xi],[\phi]) \in \Sigma^0_r \times \PP(\HKC) | \, \text{$\xi \phi=0$ in $H^1 (\OO_C)$}\}.
$$
Thus it comes with two projections
$$
\xymatrix{
&\PP({\cal W}_{\Sigma^0_r })\ar[ld]_{p_1} \ar[rd]^{p_2}&\\
\Sigma^0_r&&\PP(\HKC)
}
$$
 
The projection on the first factor
$$
p_1: \PP({\cal W}_{\Sigma^0_r }) \longrightarrow \Sigma^0_r 
$$
is the structure projection and the projection on the second factor
$$
p_2:\PP({\cal W}_{\Sigma^0_r }) \longrightarrow \PP(\HKC)
$$
corresponds to the inclusion of sheaves
$$
\xymatrix@R=12pt{
	0\ar[r]& {\cal W}_{\Sigma^0_r }\ar@{=}[d] \ar[r]&\HKC \otimes \OO_{\Sigma^0_r}\\
&	ker({\mathfrak{c}}_{\Sigma^0_r})(1)&
}
$$ 
in the exact sequence \eqref{pullback-univseq}.

The aforementioned refinement consists of attaching to every point $([\xi],[\phi]) \in \PP({\cal W}_{\Sigma^0_r })$ a canonical weakly decreasing filtration of $W_{\xi}$ 
\begin{equation}\label{phi-filt}
W_{\xi}=W^0_{\xi}([\phi]) \supset W^1_{\xi}([\phi]) \supset \cdots \supset W^{i}_{\xi}([\phi]) \supset  W^{i+1}_{\xi}([\phi]) \supset \cdots.
\end{equation}
To relate this with the contents of previous sections: this is a generalization of the maximal ladder in Lemma \ref{lem:filt-phi}. We have seen that this was a key step in the proof of Theorem \ref{thm:rk1KS}, so it should be expected that various properties of the above filtration contain a certain amount of geometry of the curve $C$. 

The plan of this section is to construct the above filtration, discuss some of its basic properties and its relevance to the geometry of the canonical embedding of $C$.

\vspace{0.2cm}
The starting point  of our construction is the canonical morphism
\begin{equation}\label{wedge2W}
	\mbox{$\bigwedge^2 {\cal W}_{\Sigma^0_r } \otimes \OO_{\Sigma^0_r } (-1) \longrightarrow \Big(\HKC \otimes \OO_{\Sigma^0_r} \Big) / {\cal W}_{\Sigma^0_r },$}
\end{equation}
or, more concretely, for every closed point $[\xi] \in {\Sigma^0_r }$ we produce a canonical, up to a scalar multiple, linear map
\begin{equation}\label{wedge2atxi}
\mbox{$\bigwedge^2 { W}_{\xi }  \longrightarrow \HKC / { W}_{\xi },$}
\end{equation}
where ${ W}_{\xi }$ is the fibre of $ker({\mathfrak{c}}_{\Sigma^0_r})={\cal W}_{\Sigma^0_r } (-1)$ at $[\xi]$.

Let us begin with the concrete version, the map \eqref{wedge2atxi}. We fix $\xi \in H^1(\Theta_{C})$ lying over the point $[\xi]\in {\Sigma^0_r }$ and consider an extension sequence
$$
\xymatrix{
	0\ar[r]&\OO_C \ar[r]^s& \EE_{\xi} \ar[r]^(.3)t& \OO_{C} (K_C)\ar[r]&0
}
$$
representing the class $\xi$. Our construction starts by choosing a Koszul representative of $\xi$. Observe, the monomorphism $s$ is determined by a choice of a nowhere vanishing global section of $\EE_{\xi}$, call it $e$, and then the epimorphism $t$ is determined by the exterior product with a nonzero scalar multiple of $e$, that is, the exact sequence has the form
$$
\xymatrix{
	0\ar[r]&\OO_C \ar[r]^e& \EE_{\xi} \ar[r]^(.35){\lambda e \wedge}& \OO_{C} (K_C)\ar[r]&0,
}
$$ 
for some $\lambda \in \CC^{\times}$. It is well known that two exact sequences 
 $$
 \xymatrix{
 	0\ar[r]&\OO_C \ar[r]^e& \EE_{\xi} \ar[r]^(.35){\lambda e \wedge}& \OO_{C} (K_C)\ar[r]&0\\
 		0\ar[r]&\OO_C \ar[r]^{e'}& \EE'_{\xi} \ar[r]^(.35){\lambda'e'\wedge}& \OO_{C} (K_C)\ar[r]&0
 }
 $$
 represent the same extension class $\xi$, if the they are related by the commutative diagram
 $$
 \xymatrix{
 	0\ar[r]&\OO_C \ar[r]^e \ar@{=}[d]& \EE_{\xi} \ar[r]^(.35){\lambda e\wedge}\ar[d]^f& \OO_{C} (K_C)\ar[r]\ar@{=}[d]&0\\
 	0\ar[r]&\OO_C \ar[r]^{e'}& \EE'_{\xi} \ar[r]^(.35){\lambda'e'\wedge}& \OO_{C} (K_C)\ar[r]&0
 }
 $$
 In particular, we can always choose an exact sequence representing $\xi$ of
 the form
 $$
 \xymatrix{
 	0\ar[r]&\OO_C \ar[r]^{e_{\xi}}& \EE_{\xi} \ar[r]^(.4){e_{\xi}\wedge  }& \OO_{C} (K_C)\ar[r]&0,
 }
 $$
 where $e_{\xi}$ is a nowhere vanishing global section of $\EE_{\xi}$, that is $\xi$ can be represented by the Koszul sequence of some nowhere vanishing global section $e_{\xi}$ of $\EE_{\xi}$ corresponding to the image of the constant section $1 \in H^0(\OO_C)$ under the monomorphism in the extension sequence. The two Koszul sequences
$$
\xymatrix{
	0\ar[r]&\OO_C \ar[r]^{e_{\xi}}& \EE_{\xi} \ar[r]^(.4){e_{\xi}\wedge  }& \OO_{C} (K_C)\ar[r]&0\\
		0\ar[r]&\OO_C \ar[r]^{e'_{\xi}}& \EE'_{\xi} \ar[r]^(.4){e'_{\xi}\wedge  }& \OO_{C} (K_C)\ar[r]&0,
	}
	$$
represent the same $\xi$ if they are related by the commutative diagram
$$
\xymatrix{
	0\ar[r]&\OO_C \ar@{=}[d]\ar[r]^{e_{\xi}}& \EE_{\xi}\ar[d]^f \ar[r]^(.4){e_{\xi}\wedge  }& \OO_{C} (K_C)\ar[r]\ar@{=}[d]&0\\
	0\ar[r]&\OO_C \ar[r]^{e'_{\xi}}& \EE'_{\xi} \ar[r]^(.4){e'_{\xi}\wedge  }& \OO_{C} (K_C)\ar[r]&0.
}
$$
Thus we have the following.
\begin{lem}\label{lem:Koszrep-xi}
	Every nonzero Kodaira-Spencer class $\xi \in H^1(\OO_C (-K_C))$ can be represented by the Koszul sequence
	$$
	\xymatrix{
		0\ar[r]&\OO_C \ar[r]^{e_{\xi}}& \EE_{\xi} \ar[r]^(.4){e_{\xi}\wedge  }& \OO_{C} (K_C)\ar[r]&0
	}
		$$
		of a nowhere vanishing global section $e_{\xi}$ of $\EE_{\xi}$. Furthermore, such a representation is unique up to an isomorphism:
		two pairs $(\EE_{\xi},e_{\xi})$ and $(\EE'_{\xi},e'_{\xi})$ represent the same cohomology class if and only if there is a bundle isomorphism
		$$
		f:\EE_{\xi}\longrightarrow \EE'_{\xi}
		$$
		such that
		
		1) $f(e_{\xi})=e'_{\xi}$,
		
		2) $f(s)\wedge f(t)=s\wedge t$, for all global sections $s$ and $t$ of $\EE_{\xi}$.
\end{lem}
From now on we always fix a Koszul representative for the extension sequence associated to a nonzero cohomology class $\xi$ in $H^1(\OO_C (-K_C))$:
\begin{equation}\label{Koszulseq-xi}
\xymatrix{
	0\ar[r]&\OO_C \ar[r]^{e_{\xi}}& \EE_{\xi} \ar[r]^(.4){e_{\xi}\wedge  }& \OO_{C} (K_C)\ar[r]&0.
}
\end{equation}
 
 Since we work on the projectivization $\PP(H^1(\OO_C (-K_C)))$ of $H^1(\OO_C (-K_C))$ we need to address the issue of the Koszul representatives of multiples of $\xi$ by a nonzero scalar. The following is easy to work out.
 \begin{lem}\label{lem:xi-scale}
 	Let $\xi$ be a nonzero Kodaira-Spencer class of $C$ and let $t\in \CC^{\times}$. Then a Koszul representative of $t\xi$ is obtain by scaling the arrows of the Koszul sequence of $(\EE_{\xi},e_{\xi})$ by the reciprocal of a square root
 	$\sqrt{t}$ of $t$, that is, setting $e'_{\xi}=\frac{1}{\sqrt{t}} e_{\xi}$, the exact sequence
 	$$
 	\xymatrix{
 		0\ar[r]&\OO_C \ar[r]^{e'_{\xi}}& \EE_{\xi} \ar[r]^(.4){e'_{\xi}\wedge  }& \OO_{C} (K_C)\ar[r]&0
 	}
 	$$
 	represents the cohomology class $t\xi$.
\end{lem}	 
This will imply that most of the constructions below depend only on the projective class $[\xi]$ of $\xi$.

We return to the Koszul sequence \eqref{Koszulseq-xi} to obtain the exact sequence of the cohomology groups 
$$
\xymatrix{
	0\ar[r]&H^0 (\OO_C )\ar[r]^{e_{\xi}}& H^0(\EE_{\xi}) \ar[r]^(.37){e_{\xi}\wedge }&\HKC \ar[r]^(.55){\xi}& H^1(\OO_C),
}
$$
where the last arrow is the cup product with the class $\xi$. In particular,
the kernel of this map is the subspace $W_{\xi}$ and we obtain 
\begin{equation}\label{xi-glsec}
\xymatrix{
	0\ar[r]&H^0 (\OO_C )\ar[r]^{e_{\xi}}& H^0(\EE_{\xi}) \ar[r]^(.57){e_{\xi}\wedge }& W_{\xi}\ar[r]&0.
}
\end{equation}
Choose a splitting 
$$
\alpha_{\xi}: W_{\xi} \longrightarrow  H^0(\EE_{\xi})
$$
in the above sequence. Thus we have
\begin{equation}\label{split-eq}
	 e_{\xi}\wedge\alpha_{\xi}(\phi) =\phi, \, \forall \phi \in W_{\xi}.
\end{equation}
The `wedge' in the above equation stands for the image under the map
\begin{equation}\label{wedge2H0EExi}
\mbox{$\bigwedge^2 H^0 (\EE_{\xi})\longrightarrow H^0(det(\EE_{\xi}))=\HKC.$}
\end{equation}
The resulting composition
$$
\mbox{$\bigwedge^2 W_{\xi} \hookrightarrow \bigwedge^2 H^0 (\EE_{\xi})\longrightarrow H^0(det(\EE_{\xi}))=\HKC.$}
$$
will be denoted
$$
\mbox{$\alpha_{\xi}^{(2)}: \bigwedge^2 W_{\xi} \longrightarrow \HKC.$}
$$
In particular, we write
$$
\alpha_{\xi}^{(2)}(\phi,\phi')
	$$
	for the value of $\alpha_{\xi}^{(2)}$ on the decomposable tensor $\phi \wedge \phi'$ in $\bigwedge^2 W_{\xi}$.

The map $\alpha_{\xi}^{(2)}$ just constructed is a representative of the map announced in \eqref{wedge2atxi}: it is a lift of \eqref{wedge2atxi} to
$\HKC$ which depends on the choice of the splitting map $\alpha_{\xi}$. We claim
the following.
\begin{lem}\label{lem:cxi-map}
	The composition 
	$$
	\xymatrix{
	 \bigwedge^2 W_{\xi} \ar[r]^(.4){\alpha_{\xi}^{(2)}} \ar@/_1.5pc/[rr]&\HKC \ar[r]& \HKC/ W_{\xi}
}
$$
of $\alpha_{\xi}^{(2)}$ with the projection
	$
	\HKC \longrightarrow \HKC/ W_{\xi}
	$
	is the map in \eqref{wedge2atxi}. It is independent on the splitting $\alpha_{\xi}$ and it will be denoted by
	$$
	\mbox{$c_{\xi}:\bigwedge^2 W_{\xi} \longrightarrow \HKC/ W_{\xi}.$}
	$$
	We write $c_{\xi}(\phi,\phi')$ for its value on $\phi\wedge \phi'$ in $\bigwedge^2 W_{\xi}$. 
	
\end{lem}
 \begin{pf}
 	Let $\alpha': W_{\xi} \longrightarrow H^0(\EE_{\xi})$ be another splitting of the exact sequence \eqref{gl-sec-xi}. Then we have
 	$$
 	\alpha'-\alpha_{\xi}=fe_{\xi},
 	$$
 	where $f: W_{\xi} \longrightarrow \CC$ is a linear function. Thus we have
 	$$
 	\alpha'(\phi)=\alpha_{\xi} (\phi)+f(\phi)e_{\xi}, \,\forall \phi \in W_{\xi}.
 	$$
 	From this it follows
 	$$
 	\alpha^{'(2)}(\phi,\psi)=\alpha^{(2)}_{\xi} (\phi,\psi)-f(\psi)e_{\xi} \wedge\alpha_{\xi}(\phi) +f(\phi)e_{\xi}\wedge\alpha_{\xi}(\psi) =\alpha^{(2)}_{\xi} (\phi,\psi) +f(\phi)\psi-f(\psi)\phi ,
 	$$
 	for all $\phi,\psi \in W_{\xi}$; the second equality uses the splitting equation \eqref{split-eq}. Thus the formula
 	\begin{equation}\label{change-split-pf}
 		\alpha^{'(2)}(\phi,\psi)-\alpha^{(2)}_{\xi} (\phi,\psi)=f(\phi)\psi-f(\psi)\phi, \forall \phi,\psi \in W_{\xi},
 		\end{equation}
 		for a change of splitting. From this it follows
 		$$
 		\alpha^{'(2)}(\phi,\psi)\equiv \alpha^{(2)}_{\xi} (\phi,\psi)\!\!\mod(W_{\xi}).
 		$$
 		This proves that the composition
 			$$
 		\xymatrix{
 			\bigwedge^2 W_{\xi} \ar[r]^(.4){\alpha_{\xi}^{(2)}} \ar@/_1.5pc/[rr]&\HKC \ar[r]& \HKC/ W_{\xi}
 		}
 		$$
 		
 		\vspace{0.2cm}
 		\noindent
 		is independent on the choice of the splitting $\alpha_{\xi}$.
 \end{pf}
	\begin{rem}\label{rem:alpha-scalar}
		If we change $\xi$ by a nonzero scalar multiple $t\xi$, Lemma \ref{lem:xi-scale} tells us that the Koszul sequence representing
		$t\xi$ becomes
		$$
		\xymatrix{
			0\ar[r]&\OO_C \ar[r]^{\frac{1}{\sqrt{t}}e_{\xi}}& \EE_{\xi} \ar[r]^(.4){\frac{1}{\sqrt{t}}e_{\xi}\wedge  }& \OO_{C} (K_C)\ar[r]&0
		}
	$$
	Hence the splitting $\alpha_{t\xi}$ is $\sqrt{t} \alpha_{\xi}$, while
	$\alpha^{(2)}_{t\xi}$ becomes $t\alpha^{(2)}_{\xi}$.
	\end{rem}
	We now turn to the sheaf version \eqref{wedge2W} of the map $c_{\xi}$ in
	Lemma \ref{lem:cxi-map}.
To define it we use the family version of the extension sequence attached to a Kodaira-Spencer class in the proof of the above lemma. Namely, we construct the universal extension on the Cartesian product
$Y_C:=C\times \PP(H^1(\Theta_{C}))$: denote by $\pi_i$ the projection on the $i$-th factor for $i=1,2$ and consider the group of extensions
$$
Ext^1( \pi^{\ast}_1 \OO_C (K_C), \pi^{\ast}_2 \OO_{\PP(H^1(\Theta_{C}))}(1));
$$
this is easily identified as follows
$$ 
Ext^1( \pi^{\ast}_1 \OO_C (K_C), \pi^{\ast}_2 \OO_{\PP(H^1(\Theta_{C}))}(1)) \cong H^1 (\pi^{\ast}_1 \OO_C (-K_C)\otimes \pi^{\ast}_2 \OO_{\PP(H^1(\Theta_{C}))}(1)) \cong End (H^1(\Theta_C)).
$$
The class $\xi_{univ}$ in $Ext^1( \pi^{\ast}_1 \OO_C (K_C), \pi^{\ast}_2 \OO_{\PP(H^1(\Theta_{C}))}(1))$ corresponding to the identity endomorphism
of $H^1(\Theta_C)$ is called the universal extension. It gives the exact sequence of locally free sheaves on $Y_C$
\begin{equation}\label{ext-univ}
\xymatrix{
0\ar[r]& \pi^{\ast}_2 \OO_{\PP(H^1(\Theta_{C}))}(1) \ar[r]&{\EE_{\xi_{univ}}} \ar[r]& \pi^{\ast}_1 \OO_C (K_C) \ar[r]&0.  
}
\end{equation}
By construction, the restriction of the above sequence to the fibre
$C\times\{[\xi]\}$ of $\pi_2$ over $[\xi]$ is identified via $\pi_1$ with the extension sequence on $C$
\begin{equation}\label{univ-at-xi}
\xymatrix{
	0\ar[r]& (\CC\xi)^{\ast}\otimes \OO_C \ar[r]&{\EE_{\xi}} \ar[r]&\OO_C (K_C) \ar[r]&0.  
}
\end{equation}
corresponding to the natural inclusion 
$$
 \CC\xi \hookrightarrow H^1(\Theta_{C})=H^1(\OO_C (-K_C))
 $$
  via the identifications
  $$
  Ext^1 (\OO_C (K_C),(\CC\xi)^{\ast}\otimes \OO_C)\cong (\CC\xi)^{\ast} \otimes H^1(\OO_C (-K_C)) \cong Hom(\CC\xi,H^1(\OO_C (-K_C)) ).
  $$
  In particular, the coboundary map in \eqref{univ-at-xi} is
  $$
  \HKC \longrightarrow (\CC\xi)^{\ast}\otimes H^1(\OO_C)
  $$
  is the cup product with $\xi$. Hence on the level of the global sections we obtain
  $$
 \xymatrix{
 	0\ar[r]& (\CC\xi)^{\ast}\otimes H^0(\OO_C) \ar[r]&{H^0(\EE_{\xi})} \ar[r]&W_{\xi} \ar[r]&0,  
 } 
$$
where as before
$$
W_{\xi}=ker(\HKC \stackrel{\xi}{\longrightarrow} H^1(\OO_{C})).
$$
Taking the direct image under $\pi_2$ of the universal extension in \eqref{ext-univ} we obtain
$$
\xymatrix@C=10pt{
	0\ar[r]&  \OO_{\PP(H^1(\Theta_{C}))}(1) \ar[r]&\pi_{2 \ast} {\EE_{\xi_{univ}}} \ar[r]& H^0 (\OO_C (K_C)) \otimes \OO_{\PP(H^1(\Theta_{C}))} \ar[r]&H^1(\OO_C) \otimes \OO_{\PP(H^1(\Theta_{C}))} (1).  
}
$$
Restricting to the stratum $\Sigma^0_r$ we deduce the exact sequence
\begin{equation}\label{WSig-seq}
 \xymatrix{
 	0\ar[r]&  \OO_{\Sigma^0_r}(1) \ar[r]&\pi_{2 \ast} {\EE_{\xi_{univ}}} \otimes \OO_{\Sigma^0_r} \ar[r]& {\cal W}_{\Sigma^0_r} \ar[r]&0.  
 }
\end{equation}
Its second exterior power gives
\begin{equation}\label{extpower2}
\xymatrix{
	0\ar[r]&  {\cal W}_{\Sigma^0_r}\otimes\OO_{\Sigma^0_r}(1) \ar[r]&\bigwedge^2 \Big(\pi_{2 \ast} {\EE_{\xi_{univ}}} \otimes \OO_{\Sigma^0_r}\Big) \ar[r]& \bigwedge^2 {\cal W}_{\Sigma^0_r} \ar[r]&0.  
}
\end{equation}
We have the natural morphism
$$
\pi^{\ast}_2 \Big(\pi_{2 \ast} {\EE_{\xi_{univ}}} \otimes \OO_{\Sigma^0_r}\Big)=\pi^{\ast}_2 \pi_{2 \ast} \Big(\EE_{\xi_{univ}} \otimes \OO_{\pi^{-1}_2 (\Sigma^0_r)}\Big) \longrightarrow \EE_{\xi_{univ}} \otimes \OO_{\pi^{-1}_2 (\Sigma^0_r)}
$$
coming from the fact that the functors $\pi^{\ast}_2$ and $\pi_{2 \ast}$ are adjoint. Taking its second exterior power gives the morphism
$$
\xymatrix@C=10pt@R=10pt{
	\pi^{\ast}_2 \Big( \bigwedge^2 \pi_{2 \ast} {\EE_{\xi_{univ}}} \otimes \OO_{\Sigma^0_r}\Big) \ar[r]&\bigwedge^2 \EE_{\xi_{univ}} \otimes \OO_{\pi^{-1}_2 (\Sigma^0_r)} \ar@{=}[d]\\ 
	& \pi^{\ast}_1 (\OO_C (K_C)) \otimes \pi^{\ast}_2 ( \OO_{\PP(H^1(\Theta_{C}))}(1)) \otimes \OO_{\pi^{-1}_2 (\Sigma^0_r)} 
}
$$
where the identification on the right comes from the universal extension sequence \eqref{ext-univ}. The above morphism corresponds, under adjoint property of functors $\pi^{\ast}_2$ and $\pi_{2 \ast}$, to the morphism
$$
\xymatrix@R=10pt{
\bigwedge^2 \Big(\pi_{2 \ast} {\EE_{\xi_{univ}}} \otimes \OO_{\Sigma^0_r}\Big) \ar[r]& \pi_{2\ast} \Big(\pi^{\ast}_1 (\OO_C (K_C)) \otimes \pi^{\ast}_2 ( \OO_{\PP(H^1(\Theta_{C}))}(1)) \otimes \OO_{\pi^{-1}_2 (\Sigma^0_r)} \Big)\ar@{=}[d]\\
&H^0(\OO_C (K_C))\otimes \OO_{\Sigma^0_r} (1)
}
$$
Thus we obtain the morphism
$$
\xymatrix{
	\bigwedge^2 \Big(\pi_{2 \ast} {\EE_{\xi_{univ}}} \otimes \OO_{\Sigma^0_r}\Big) \ar[r]&H^0(\OO_C (K_C))\otimes \OO_{\Sigma^0_r} (1)
}
$$
which is the sheaf version of the map in \eqref{wedge2H0EExi}. Combining this morphism with the exact sequence \eqref{extpower2} we obtain
$$
\xymatrix@C=12pt{
	0\ar[d]&0\ar[d]\\
	{\cal W}_{\Sigma^0_r}\otimes\OO_{\Sigma^0_r}(1) \ar[d] \ar@{=}[r]&	{\cal W}_{\Sigma^0_r}\otimes\OO_{\Sigma^0_r}(1) \ar[d]\\
\bigwedge^2 \Big(\pi_{2 \ast} {\EE_{\xi_{univ}}} \otimes \OO_{\Sigma^0_r}\Big) \ar[r] \ar[d]&H^0(\OO_C (K_C))\otimes \OO_{\Sigma^0_r} (1) \ar[d]\\
\bigwedge^2 {\cal W}_{\Sigma^0_r} \ar[d]&\Big( H^0(\OO_C (K_C))\otimes \OO_{\Sigma^0_r} /{\cal W}_{\Sigma^0_r} \Big) \otimes \OO_{\Sigma^0_r} (1)\ar[d]\\
0&0
}
$$
where the square on the top is commutative: this is the sheaf version of the
equation \eqref{split-eq}. Thus we can complete the diagram by the morphism
$$
\xymatrix{
\bigwedge^2 {\cal W}_{\Sigma^0_r} \ar[r]&\Big( H^0(\OO_C (K_C))\otimes \OO_{\Sigma^0_r} /{\cal W}_{\Sigma^0_r} \Big) \otimes \OO_{\Sigma^0_r} (1)
}
$$
at the bottom making the whole diagram commutative. Tensoring with $\OO_{\Sigma^0_r} (-1)$ gives the morphism
$$
\xymatrix{
	\bigwedge^2 {\cal W}_{\Sigma^0_r}\otimes \OO_{\Sigma^0_r} (-1) \ar[r]&\Big( H^0(\OO_C (K_C))\otimes \OO_{\Sigma^0_r} \Big)/{\cal W}_{\Sigma^0_r}
}
$$
asserted in \eqref{wedge2W}.
\begin{rem}\label{rem:alpha}
	In Lemma \ref{lem:cxi-map} is given an explicit lift of the map
	$$
	\xymatrix{
	c_{\xi}:	\bigwedge^2 { W}_{\xi} \ar[r]& H^0(\OO_C (K_C)) /{ W}_{\xi}
	}
$$
to the map
$$
 \xymatrix{
 		\bigwedge^2 { W}_{\xi} \ar[r]& H^0(\OO_C (K_C)) 
 }
 $$
 in terms of the splitting map
 $$
 \alpha_{\xi} : W_{\xi} \longrightarrow H^0(\EE_{\xi}).
 $$
 This is the map
 $$
\xymatrix{
	\alpha^{(2)}_{\xi}:	\bigwedge^2 { W}_{\xi} \ar[r]& H^0(\OO_C (K_C)).
}
$$
This is a useful computational mechanism and also gives a connection between the properties of the global sections of the vector bundle $\EE_{\xi}$ and the linear subsystem $|W_{\xi}|$. 

The map $\alpha^{(2)}_{\xi}$ can be viewed as a higher order product on the subspace $W_{\xi}$ or a vector valued, possibly degenerate, symplectic structure.

 Let us give two more formulas we will need in the sequel:
\begin{equation}\label{alpha-aplha2-eq}
\phi '\alpha_{\xi}(\phi)-\phi \alpha_{\xi}(\phi') =\alpha^{(2)}_{\xi} (\phi,\phi') e_{\xi},\,\forall \phi,\phi'\in W_{\xi},
\end{equation}
where both sides are viewed as global sections of $\EE_{\xi} (K_C)$;
\begin{equation}\label{aplha2-eq}
\phi''\alpha^{(2)}_{\xi} (\phi,\phi') - \phi'\alpha^{(2)}_{\xi} (\phi,\phi'') +\phi \alpha^{(2)}_{\xi} (\phi',\phi'')=0, \, \forall \phi,\phi',\phi'' \in W_{\xi}.
\end{equation}
The last formula follows immediately from the equation \eqref{alpha-aplha2-eq} and it takes place in the space $H^0(\OO_C (2K_C))$. We will call it the Koszul cocycle relation: view  the map
$$
\mbox{$\alpha^{(2)}_{\xi}:\bigwedge^2 { W}_{\xi} \longrightarrow H^0(\OO_C (K_C))$}
$$
 as an element of $ \bigwedge^2 { W}^{\ast}_{\xi} \otimes \HKC$; it sits in the Koszul complex associated to the pair $(\OO_C (K_C), W_{\xi}) $:
$$
\xymatrix{
	W^{\ast}_{\xi}\otimes H^0(\OO_{C})\ar[r]^(.45){d_{Kosz}}& \bigwedge^2 { W}^{\ast}_{\xi} \otimes \HKC \ar[r]^{d_{Kosz}}& \bigwedge^3 { W}^{\ast}_{\xi} \otimes H^0(\OO_C(2K_C)) \ar[r]&\cdots
}
$$
and the equation \eqref{aplha2-eq} means that
$$
d_{Kosz}(\alpha^{(2)}_{\xi})=0,
$$
that is, $\alpha^{(2)}_{\xi}$ is a Koszul cocycle, see \cite{Gr} for details about Koszul formalism. It should be noticed that the cocycle equation is independent of the splitting map: the proof of
Lemma \ref{lem:cxi-map} shows that for another splitting $\alpha'$ the difference
$$
 \alpha^{'(2)} -\alpha^{(2)}_{\xi}=d_{Kosz}(f),
 $$
 for $f\in 	W^{\ast}_{\xi}\otimes H^0(\OO_{C})$, that is, it is a Koszul coboundary.
\end{rem}

\vspace{0.2cm}
We are now ready to define the filtration of $W_{\xi}$ as stated in \eqref{phi-filt}.

\begin{lem}\label{lem:phi-filt}
	For every pair $([\xi],[\phi]) \in \PP({\cal W}_{\Sigma^0_r})$ the subspace 
	$$
	W_{\xi}=ker\Big(\HKC \stackrel{\xi}{\longrightarrow}H^1 (\OO_C)\Big)
	$$
	admits a canonical decreasing filtration
	$$
	W_{\xi}=W^0_{\xi}([\phi]) \supset W^1_{\xi}([\phi]) \supset \cdots \supset W^{i}_{\xi}([\phi]) \supset  W^{i+1}_{\xi}([\phi]) \supset \cdots \supset W^{l_{\xi}([\phi])}([\phi]),
	$$
	where $l_{\xi}([\phi])$ is the smallest index at which the filtration stabilizes, that is,
	$$
	l_{\xi}([\phi])=min\{j|\, W^{j}_{\xi}([\phi])=W^{j+1}_{\xi}([\phi])\}.
	$$
	
	Every subspace $W^{i}_{\xi}([\phi])$ contains the line $\CC \phi$ and comes equipped with a distinguished linear map
	$$
	c^{(i)}_{\xi}(\phi, \bullet): W^{i}_{\xi}([\phi]) \longrightarrow  W^{i-1}_{\xi}([\phi]) /  W^{i}_{\xi}([\phi]),
	$$
	for every $i\in [1, l([\phi])]$, such that
	$$
	 W^{i+1}_{\xi}([\phi])=ker(c^{(i)}_{\xi}([\phi], \bullet)).
	 $$
\end{lem}
\begin{pf}
	In Lemma \ref{lem:cxi-map} we have constructed the linear map
	$$
	\xymatrix{
	c_{\xi}: \bigwedge^2 W_{\xi} \ar[r]& \HKC/W_{\xi}.
}
	$$
	We rewrite it in the form
	$$
	\xymatrix{
	c_{\xi}:W_{\xi} \ar[r]& Hom (W_{\xi},\HKC/W_{\xi} ):
}
	$$
	it takes $\phi\in W_{\xi} $ to the linear map
	$$
	c_{\xi}(\phi,\bullet):W_{\xi}\longrightarrow \HKC/W_{\xi}
	$$
	defined by the rule
	$$
	\psi\mapsto c_{\xi}(\phi,\psi),  \,\forall \psi\in W_{\xi}.
	$$
	
	For a point $([\xi],[\phi]) \in \PP({\cal W}_{\Sigma^0_r})$ we set
	\begin{equation}
		\text{$W^0_{\xi}([\phi]):=W_{\xi}$ and $W^1_{\xi}([\phi]):= ker (c_{\xi}(\phi,\bullet))$,}
	\end{equation}
where $\phi$ is any representative of $[\phi]$. This gives
$$
W^0_{\xi}([\phi]) \supset W^1_{\xi}([\phi]) \supset \CC \phi,
$$
where the last inclusion uses the skew-symmetry of $c_{\xi}$. This is the first step of our construction. If the first inclusion from the left is the equality, the process stops. If not, we continue. For this we need to recall that $c_{\xi}$ admits a lifting
$$
 	\xymatrix{
	\alpha^{(2)}_{\xi}: \bigwedge^2 W_{\xi} \ar[r]& \HKC
}
$$
coming from a splitting map
$$
\alpha_{\xi}: W_{\xi} \longrightarrow H^0(\EE_{\xi}),
$$
see the proof of Lemma \ref{lem:cxi-map}. In particular, we have
$$
\alpha^{(2)}_{\xi}(\phi,\bullet): W^0_{\xi}([\phi]) \longrightarrow \HKC,
$$
a lifting of $c_{\xi}(\phi,\bullet)$. By definition, this map, restricted to $W^1_{\xi}([\phi])$, takes values in  $W^0_{\xi}([\phi])$:
$$
\alpha^{(2)}_{\xi}(\phi,\bullet): W^1_{\xi}([\phi]) \longrightarrow  W^0_{\xi}([\phi]).
$$
Composing with the projection
$$
W^0_{\xi}([\phi]) \longrightarrow W^0_{\xi}([\phi])/W^1_{\xi}([\phi])
$$
gives the linear map
$$
W^1_{\xi}([\phi]) \longrightarrow W^0_{\xi}([\phi])/W^1_{\xi}([\phi])
$$
which will be denoted by $c^{(1)}_{\xi}(\phi,\bullet)$.
We claim that the map 
\begin{equation}\label{c1xi-phi}
c^{(1)}_{\xi}(\phi,\bullet):W^1_{\xi}([\phi]) \longrightarrow W^0_{\xi}([\phi])/W^1_{\xi}([\phi])
\end{equation}
is independent on the choice of the splitting map $\alpha_{\xi}$. Indeed, another splitting 
$$
\alpha':W_{\xi}\longrightarrow H^0(\EE_{\xi})
$$
gives the formula
$$
\alpha^{'{(2)}}(\psi,\psi')-\alpha^{(2)}_{\xi}(\psi,\psi')=f(\psi)\psi'- f(\psi') \psi, \,\forall \psi,\psi' \in W_{\xi},
$$
for some linear map $f \in W^{\ast}_{\xi}$, see the proof of Lemma \ref{lem:cxi-map}. This implies that for $\psi,\psi'$ in $W^1_{\xi}([\phi])$,
the difference lies in $W^1_{\xi}([\phi])$. In particular, we have
$$
\alpha^{'{(2)}}(\phi,\psi)-\alpha^{(2)}_{\xi}(\phi,\psi)=\Big( f(\phi)\psi- f(\psi) \phi \Big) \in W^1_{\xi}([\phi]), \,\forall \psi \in W^1_{\xi} ([\phi]),
$$
where the last inclusion because by construction $\CC \phi \subset W^1_{\xi}([\phi])$. Thus the map
$c^{(1)}_{\xi}(\phi,\bullet)$ in \eqref{c1xi-phi} is well defined and given by the formula
$$
c^{(1)}_{\xi}(\phi,\psi)\equiv \alpha^{(2)}_{\xi}(\phi,\psi)\!\!\!\mod (W^1_{\xi}([\phi])),\,\,\forall \psi \in W^1_{\xi}([\phi]),
$$
where $\alpha_{\xi}: W_{\xi}\longrightarrow H^0(\EE_{\xi})$ any splitting map. 
We go on to define the next step of the filtration as the kernel of $c^{(1)}_{\xi}(\phi,\bullet)$:
$$
 W^2_{\xi}([\phi]):=ker(c^{(1)}_{\xi}(\phi,\bullet)).
 $$
 Again in terms of the splitting $\alpha_{\xi}$ we have
 $$
\alpha^{(2)}_{\xi}(\phi,\bullet): W^2_{\xi}([\phi]) \longrightarrow  W^1_{\xi}([\phi])
$$
and we define
$$
c^{(2)}_{\xi}(\phi,\bullet):W^2_{\xi}([\phi]) \longrightarrow W^1_{\xi}([\phi])/W^2_{\xi}([\phi])
$$
by the formula
$$
c^{(2)}_{\xi}(\phi,\psi)\equiv \alpha^{(2)}(\phi,\psi)\!\!\!\mod (W^2_{\xi}([\phi])),\,\,\forall \psi \in W^2_{\xi}([\phi]).
$$
If $W^2_{\xi}([\phi])=W^1_{\xi}([\phi])$, we stop. Otherwise we continue.
Since the dimension of the space drops every time we continue, the process terminates after finite number of steps. The final point is that the filtration of $W_{\xi}$ constructed above depends only on the point $[\xi]$ in $\Sigma^0_r$ since the multiplication of $\xi$ by a scalar $t \in \CC^{\times}$ changes the map $\alpha_{\xi}$ to $\sqrt{t}\alpha_{\xi}$, while the map $\alpha^{(2)}_{\xi}$ changes to $t\alpha^{(2)}_{\xi}$, Remark \ref{rem:alpha-scalar}.
\end{pf}

The above result tells us that at every point $([\xi],[\phi])$ of the incidence correspondence $\PP({\cal W}_{\Sigma^0_r })$ the space
$W_{\xi}$, naturally attached to the point $([\xi],[\phi])$, admits a distinguished filtration
$$
W^{\bullet}_{\xi}([\phi])=\{W^0_{\xi}([\phi]) \supset W^1_{\xi}([\phi]) \supset \cdots \supset W^{i}_{\xi}([\phi]) \supset  W^{i+1}_{\xi}([\phi]) \supset \cdots \supset W^{l_{\xi}([\phi])} ([\phi])\}
$$
as described in Lemma \ref{lem:phi-filt}. We call it  {\it $([\xi],[\phi])$-filtration of $W_{\xi}$} and its length $l_{\xi}([\phi])$ will be called {\it the length of $W_{\xi}$ at $([\xi],[\phi])$} or simply
{\it $([\xi],[\phi])$-length of $W_{\xi}$}. From now on we let
$$
W^{l_{\xi}([\phi])+1} ([\phi]):=0
$$
and write  $([\xi],[\phi])$-filtration in the augmented form
\begin{equation}\label{xi-phi-filt}
	W^{\bullet}_{\xi}([\phi])=\{W^0_{\xi}([\phi]) \supset W^1_{\xi}([\phi]) \supset \cdots \supset \cdots \supset W^{l_{\xi}([\phi])} ([\phi])\supset W^{l_{\xi}([\phi])+1} ([\phi])=0\}.
\end{equation}

 Set 
\begin{equation}\label{xi-phi-gr}
Gr^{\bullet}_{W^{\bullet}_{\xi}([\phi])}=\bigoplus^{l_{\xi}([\phi])}_{i=0}Gr^{i}_{W^{\bullet}_{\xi}([\phi])}:= \bigoplus^{l_{\xi}([\phi])}_{i=0} W^i_{\xi}([\phi])/ W^{i+1}_{\xi}([\phi]), 
\end{equation}
the graded space attached to $([\xi],[\phi])$-filtration; the dimension of the graded pieces are denoted
$$
h^i ([\xi],[\phi]):=dim \Big(W^i_{\xi}([\phi])/ W^{i+1}_{\xi}([\phi]) \Big), \,\forall i\in [0,l_{\xi}([\phi])].
$$
Those values form
$$
{\bf h} ([\xi],[\phi])=(h^0 ([\xi],[\phi]),h^1 ([\xi],[\phi]), \ldots, h^{l_{\xi}([\phi])} ([\xi],[\phi]))
$$
a partition of $dim(W_{\xi})$, that is, we have the equality
$$
dim(W_{\xi})=\sum^{l_{\xi}([\phi])}_{i=0} h^i ([\xi],[\phi]).
$$
In addition, the linear maps
$$
c^{(i)}_{\xi}(\phi,\bullet):W^i_{\xi}([\phi]) \longrightarrow W^{i-1}_{\xi}([\phi])/W^i_{\xi}([\phi])
$$
in Lemma \ref{lem:phi-filt} give rise to the {\it injective} maps
$$
gr(c^{(i)}_{\xi}(\phi,\bullet)): W^i_{\xi}([\phi])/W^{i+1}_{\xi}([\phi]) \longrightarrow W^{i-1}_{\xi}([\phi])/W^i_{\xi}([\phi]), \,\forall i\in [1,l_{\xi}([\phi])-1]:
$$
for the injectivity we recall that $W^{i+1}_{\xi}([\phi])$ is defined as the kernel of the map $c^{(i)}_{\xi}(\phi,\bullet)$ for $i$ in the range $[1,l_{\xi}([\phi])-1]$; the last map
$$
gr(c^{(l_{\xi}([\phi]))}_{\xi}(\phi,\bullet)):W^{l_{\xi}([\phi])}_{\xi}([\phi]) =W^{l_{\xi}([\phi])}_{\xi}([\phi])/W^{l_{\xi}([\phi])+1}_{\xi}([\phi]) \longrightarrow W^{l_{\xi}([\phi])-1}([\phi])/W^{l_{\xi}([\phi])}_{\xi}([\phi])
$$
is identically zero.
We summarize the above discussion in the following.
\begin{pro}\label{pro:xi-phi-filt}
1)	To every point $([\xi],[\phi])$ of the incidence correspondence 
	$$
	\PP({\cal W}_{\Sigma^0_r })=\{([\xi],[\phi]) \in \Sigma^0_r \times \PP(\HKC)|\,\text{$\xi \phi=0$ in $H^1(\OO_X)$}\}
	$$
	is canonically attached $([\xi],[\phi])$-filtration
	$$
	\begin{gathered}
		W^{\bullet}_{\xi}([\phi])=\{W^0_{\xi}([\phi]) \supset W^1_{\xi}([\phi]) \supset \cdots \supset W^{i}_{\xi}([\phi]) \supset  W^{i+1}_{\xi}([\phi])\\
		 \supset \cdots \supset W^{l_{\xi}([\phi])} ([\phi])\supset W^{l_{\xi}([\phi])+1} ([\phi])=0\}
		\end{gathered}
		$$
		of $W_{\xi}$. This comes with the natural inclusion
		$$
		\CC\phi \subset W^{l_{\xi}([\phi])} ([\phi]).
		$$
		
		2) The above filtration is equipped with linear maps
		$$
		c^{(i)}_{\xi}(\phi,\bullet):W^i_{\xi}([\phi]) \longrightarrow W^{i-1}_{\xi}([\phi])/W^i_{\xi}([\phi])
		$$
		 for every $i\in [1,l_{\xi}([\phi])]$ subject to the property
		 $$
		 ker(c^{(i)}_{\xi}(\phi,\bullet))=W^{i+1}_{\xi}([\phi]), \, \forall i\in [1,l_{\xi}([\phi])-1],
		 $$
		 and the last map $c^{(l_{\xi}([\phi]))}_{\xi}(\phi,\bullet)$ is identically zero.

		 3)
		  The associated graded space
		  $$
		 Gr^{\bullet}_{W^{\bullet}_{\xi}([\phi])}=\bigoplus^{l_{\xi}([\phi])}_{i=0}Gr^{i}_{W^{\bullet}_{\xi}([\phi])}= \bigoplus^{l_{\xi}([\phi])}_{i=0} W^i_{\xi}([\phi])/ W^{i+1}_{\xi}([\phi])
		 $$
		 is equipped with the graded map of degree $(-1)$
		 $$
		gr(c^{(\bullet)}_{\xi}(\phi,\bullet)): Gr^{\bullet}_{W^{\bullet}_{\xi}([\phi])} \longrightarrow Gr^{\bullet}_{W^{\bullet}_{\xi}([\phi])} [-1],
		 $$
		 where 
		 $$
		 gr(c^{(i)}_{\xi}(\phi,\bullet)): Gr^{i}_{W^{\bullet}_{\xi}([\phi])}= W^i_{\xi}([\phi])/ W^{i+1}_{\xi}([\phi])\longrightarrow Gr^{i}_{W^{\bullet}_{\xi}([\phi])} [-1] =W^{i-1}_{\xi}([\phi])/ W^{i}_{\xi}([\phi])
		 $$
		 for every $i\in [1,l_{\xi}([\phi])-1]$ is injective and identically zero for $i=0$ and $l_{\xi}([\phi])$.
		 
		 4) The dimensions
		 $$
		 h^i ([\xi],[\phi]):=dim \Big(W^i_{\xi}([\phi])/ W^{i+1}_{\xi}([\phi]) \Big), \,\forall i\in [0,l_{\xi}([\phi])].
		 $$
		 are subject to
		 $$
		 \begin{gathered}
		 a)\,\sum^{l_{\xi}([\phi])}_{i=0}h^i ([\xi],[\phi])=dim(W_{\xi}),\\
		 b)\,\, h^i ([\xi],[\phi]) \geq h^{i+1} ([\xi],[\phi]),\,\,\forall i\in [0,l_{\xi}([\phi])-2].
		 \end{gathered}
		 $$		 
\end{pro}

\vspace{0.2cm}
On the conceptual level the above result can be viewed as follows: we take the pull back 
$
p^{\ast}_1 ({\cal W}_{\Sigma^0_r })
$
of the vector bundle ${\cal W}_{\Sigma^0_r }$ under the structure projection 
$$
p_1: 	\PP({\cal W}_{\Sigma^0_r }) \longrightarrow \Sigma^0_r; 
$$
Proposition \ref{pro:xi-phi-filt} says that $p^{\ast}_1 ({\cal W}_{\Sigma^0_r })$ comes equipped with a period filtration; the map of degree $(-1)$ could be viewed as Griffiths transversality property of the classical variation of Hodge structure, see \cite{G1}. We will  give a sheaf version of the above construction in \S14. For now let
us see what the $([\xi],[\phi])$-filtrations tell us about the geometry.

\begin{lem}\label{lem:Wlatleast2}
	 Assume the dimension $h^{l_{\xi}([\phi])}=dim(W^{l_{\xi}([\phi])}_{\xi}([\phi]))$ at least two.
	Then the rank $2$ bundle $\EE_{\xi}$ in the middle of the extension sequence
	$$
	\xymatrix{
	0\ar[r]& \OO_{C} \ar[r]&\EE_{\xi}\ar[r]& \OO_{C}(K_C)\ar[r]&0
}
$$
has two linearly independent global sections $e$ and $e'$ which are proportional, that is, the image of $e\wedge e'$ under the map
$$
\xymatrix{
	\bigwedge^2 H^0(\EE_{\xi}) \ar[r]& H^0(det(\EE_{\xi}))=\HKC
}
$$
is zero. In addition, $e$ and $e'$ can be chosen so that under the 
map
$$
H^0(\EE_{\xi}) \longrightarrow W_{\xi}
$$
induced by the epimorphism of the extension sequence the section $e$ goes to $\phi$ and $e'$ goes to an element $\phi'$ of  $W^{l_{\xi}([\phi])}_{\xi}([\phi]) \setminus \CC\phi$.
\end{lem}

\begin{pf}
	Fix a splitting map
	$$
	\alpha: W_{\xi} \longrightarrow H^0(\EE_{\xi}).
	$$
	and recall the map
	$$
	\xymatrix{
	\alpha^{(2)}: \bigwedge^2W_{\xi} \ar[r]& H^0(\bigwedge^2\EE_{\xi})=\HKC.
}
$$	
	By definition of the $([\xi],[\phi])$-filtration we have
	$$
	\alpha^{(2)}(\phi, \bullet): W^i_{\xi} ([\phi]) \longrightarrow W^{i-1}_{\xi}([\phi])
	$$
	for all $i\in [1,l_{\xi}([\phi])]$. Furthermore, for $l_{\xi}([\phi])$ 
	the composition with the projection
	$$
	 W^{l_{\xi}([\phi])}_{\xi}([\phi]) \longrightarrow  W^{l_{\xi}([\phi])-1}_{\xi}([\phi])/ W^{l_{\xi}([\phi])}_{\xi}([\phi])
	 $$
	 is identically zero, see Proposition \ref{pro:xi-phi-filt}, 3). This means that $\alpha^{(2)}(\phi,\bullet)$ induces an endomorphism of $W^{l_{\xi}([\phi])}_{\xi}([\phi])$
	 $$
	 	\alpha^{(2)}(\phi,\bullet): W^{l_{\xi}([\phi])}_{\xi} ([\phi]) \longrightarrow W^{l_{\xi}([\phi])}_{\xi}([\phi]).
	 	$$
	 	By construction the line $\CC\phi$ is contained in $W^{l_{\xi}([\phi])}_{\xi} ([\phi])$, Proposition \ref{pro:xi-phi-filt}, 1). Since $\alpha^{(2)}(\phi,\phi)=0$, the above endomorphism descends to the endomorphism of the quotient space
	 	$$
	 	W^{l_{\xi}([\phi])}_{\xi} ([\phi])/\CC\phi \longrightarrow W^{l_{\xi}([\phi])}_{\xi}([\phi])/\CC \phi.
	 	$$
	 	By our assumption the dimension of $W^{l_{\xi}([\phi])}_{\xi} ([\phi])$ is at least two. So the the quotient is nonzero. For an eigenvalue $\lambda_0$ of the above endomorphism we choose an eigenvector in $W^{l_{\xi}([\phi])}_{\xi} ([\phi])/\CC\phi$ and call $\phi' \in W^{l_{\xi}([\phi])}_{\xi} ([\phi])$ its lifting. Then $\phi'$ lies in $W^{l_{\xi}([\phi])}_{\xi} ([\phi]) \setminus \CC \phi$ and the value $\alpha^{(2)}(\phi,\phi')$ has the form
	 	$$
	 	\alpha^{(2)}(\phi,\phi')=\lambda_0 \phi' -\mu_{0}\phi
	 	$$
	 	for some constant $\mu_0 \in \CC$. We now use it in the equation
	 	$$
	 	\phi'\alpha(\phi)-\phi \alpha(\phi')=\alpha^{(2)}(\phi,\phi') e_{\xi},
	 	$$
	 	see \eqref{alpha-aplha2-eq}, to obtain
	 	$$
	 	\phi'\alpha(\phi)-\phi \alpha(\phi')=\alpha^{(2)}(\phi,\phi') e_{\xi}=\Big(\lambda_0 \phi' -\mu_{0}\phi \Big)e_{\xi}.
	 	$$
	 	From this it follows
	 	\begin{equation}\label{eq-twosec}
	 	\phi' \Big(\alpha(\phi)-\lambda_0 e_{\xi}\Big)-\phi \Big(\alpha(\phi')-\mu_0 e_{\xi}\Big)=0.
	 \end{equation}
	 	The two global sections
	 	$$
	 	\text{$e:=\alpha(\phi)-\lambda_0 e_{\xi}$ and $e':=\alpha(\phi')-\mu_0 e_{\xi}$}
	 	$$
	 	in the above equation have all the properties asserted in the lemma:
	 	they are linearly independent since under the map
	 	$$
	 	H^0(\EE_{\xi})\longrightarrow W_{\xi}
	 	$$
	 	they go to $\phi$ and $\phi'$ respectively and those are chosen to be linearly independent in the subspace $W^{l_{\xi}([\phi])}_{\xi} ([\phi]) \subset W_{\xi}$; the equation \eqref{eq-twosec} gives the proportionality. 
\end{pf}

The two global sections in Lemma \ref{lem:Wlatleast2} generate a subsheaf of rank $1$ in $\EE_{\xi}$. We denote by ${\cal L}$ its saturation and obtain another extension sequence for $\EE_{\xi}$:
$$
\xymatrix{
		0\ar[r]& {\cal L} \ar[r]&\EE_{\xi}\ar[r]& \OO_{C}(D)\ar[r]&0.
}
$$
Putting it together with the defining extension gives the diagram
$$
\xymatrix{
	&&0\ar[d]&&\\
	&& {\cal L} \ar[d] \ar[rd]&&\\
	0\ar[r]& \OO_{C} \ar[r]\ar[rd]&\EE_{\xi}\ar[r] \ar[d]& \OO_{C}(K_C)\ar[r]&0\\
	&&\OO_{C} (D)\ar[d]&&\\
	&&0&&
}
$$
where the slanted arrows of the diagram are the obvious composition arrows.
Those are nonzero since by construction
$$
\CC\{e,e'\} \hookrightarrow H^0({\cal L})
$$
and the slanted arrow in the upper right corner of the diagram takes takes $e$ and $e'$ to the sections $\phi$ and $\phi'$ of $\OO_{C}(K_C)$, see the proof of Lemma \ref{lem:Wlatleast2}. In particular, the divisor $D$ is effective and the slanted arrows of the diagram are defined by the multiplication with a distinguished global section $\delta$ of $\OO_{C} (D)$. From now on we write the above diagram as follows

\begin{equation}\label{xi-D-diag}
	\xymatrix{
		&&0\ar[d]&&\\
		&& {\cal O}_C (K_C-D) \ar[d] \ar[rd]^{\delta}&&\\
		0\ar[r]& \OO_{C} \ar[r]\ar[rd]_{\delta}&\EE_{\xi}\ar[r] \ar[d]& \OO_{C}(K_C)\ar[r]&0\\
		&&\OO_{C} (D)\ar[d]&&\\
		&&0&&
	}
\end{equation}

We summarize all of the above in the following.
\begin{pro}\label{pro:xi-Wlatleast2}
	Let $[\xi]\in\Sigma^0_r $ and assume there is a closed point
	$([\xi],[\phi])\in \PP({\cal W}_{\Sigma^0_r }) $ with
	the $([\xi],[\phi])$-filtration of $W_{\xi}$ subject to
	$$
	h^{l_{\xi}([\phi])}=dim (W^{l_{\xi}([\phi])}_{\xi} ([\phi]))\geq 2.
	$$
	Then the extension sequence
	$$
	\xymatrix{
		0\ar[r]& \OO_{C} \ar[r]&\EE_{\xi}\ar[r] & \OO_{C}(K_C)\ar[r]&0
	}
$$
defined by $\xi$ fits into the diagram \eqref{xi-D-diag}. This means that
the global section $\delta \in H^0(\OO_C (D))$ annihilates $\xi$, that is,
$$
\xi \in ker \Big(H^1 (\Theta_{C})=H^1(\OO_C (-K_C)) \stackrel{\delta}{\longrightarrow} H^1(\OO_C (D-K_C)) \Big).
$$
Furthermore, $h^0 (\OO_C (K_C-D)) \geq 2$ and we have the inclusions
$$
\CC\phi \subset \delta H^0 (\OO_C (K_C-D)) \subset W^{l_{\xi}([\phi])}_{\xi} ([\phi]).
$$
\end{pro}

\begin{rem}\label{rem:secant}
1).	The statement about $\xi$ being annihilated by the section $\delta$ in Proposition \ref{pro:xi-Wlatleast2} can be rephrased in terms of secant linear subspaces of the bicanonical embedding of $C$:
	set $D_{\delta}:= (\delta=0)$ to be the divisor of zeros of $\delta$ and denote
	$$
	\{D_{\delta}\}_{2K_C} 
	$$
 the linear span of $D_{\delta}$ in the bicanonical embedding of $C$;
 this is the secant subspace of dimension $(deg(D)-1)$ of $C\subset \PP(H^1(\OO_C (-K_C)))$ intersecting
 $C$  along a subscheme containing $D_{\delta}$, and $[\xi]$ is contained in the secant space $\{D_{\delta}\}_{2K_C}$ and no other secant subspace of  $\{D_{\delta}\}_{2K_C}$ spanned by a proper subscheme of $D_{\delta}$, see Remark \ref{rem:D-D'}, 1). 
 
 2). The converse to  Proposition \ref{pro:xi-Wlatleast2} is valid as well: given the diagram \eqref{xi-D-diag} with the line subbundle $\OO_C (K_C-D)$ subject to
 $$
 h^0(\OO_C (K_C-D))\geq 2,
 $$
 the subspace $\delta H^0 (\OO_C (K_C-D))$ of $W_{\xi}$ lies in the
 last step $W^{l_{\xi}([\phi])}_{\xi} ([\phi])$ of $([\xi],[\phi])$-filtration for any nonzero $\phi \in \delta H^0 (\OO_C (K_C-D))$. 
\end{rem}

In case the dimension of the space $W^{l_{\xi}([\phi])}_{\xi}$ is at least two, the proof of Lemma \ref{lem:Wlatleast2} establishes a precise relation between
the properties of the endomorphism
$$
\alpha^{(2)} (\phi, \bullet): W^{l_{\xi}([\phi])}_{\xi} \longrightarrow W^{l_{\xi}([\phi])}_{\xi}
$$
and the line subsheaves of $\EE_{\xi}$ appearing in the vertical exact sequence of the diagram \eqref{xi-D-diag}. This is spelled out in the following.

\begin{pro}\label{pro:xi-Wlatleast2-eigen}
	Let $[\xi]\in\Sigma^0_r $ and let 	$([\xi],[\phi])$ be a closed point
in	$\PP({\cal W}_{\Sigma^0_r })$.
	Then for any splitting
	$$
	\alpha: W_{\xi} \longrightarrow H^0(\EE_{\xi})
	$$
	the map
	$$
	\alpha^{(2)}(\phi,\bullet):W_{\xi}\longrightarrow \HKC
	$$
	preserves the subspace  $W^{l_{\xi}([\phi])}_{\xi} ([\phi])$ and hence induces the endomorphism
	$$
	\alpha^{(2)}(\phi,\bullet): W^{l_{\xi}([\phi])}_{\xi} ([\phi]) \longrightarrow W^{l_{\xi}([\phi])}_{\xi} ([\phi]).
	$$
	 with the line $\CC\phi$ contained in $ker (\alpha^{(2)}(\phi,\bullet))$. Assume 
	 	$$
	 h^{l_{\xi}([\phi])}=dim (W^{l_{\xi}([\phi])}_{\xi} ([\phi]))\geq 2
	 $$
	 and consider the induced endomorphism
	 $$
	 \widehat{\alpha^{(2)}}(\phi,\bullet):W^{l_{\xi}([\phi])}_{\xi} ([\phi])/\CC\phi \longrightarrow W^{l_{\xi}([\phi])}_{\xi} ([\phi])/\CC\phi
	 	$$
	 	The following properties hold.
	 	
1) For each eigen value $\lambda$ of $\widehat{\alpha^{(2)}}(\phi,\bullet)$ let
$V_{\lambda}$ be the corresponding eigen subspace of $\widehat{\alpha^{(2)}}(\phi,\bullet)$ in $W^{l_{\xi}([\phi])}_{\xi} ([\phi])/\CC\phi$; denote by $\widetilde{V}_{\lambda}$ the inverse image of $V_{\lambda}$ under the quotient homomorphism
$$
W^{l_{\xi}([\phi])}_{\xi} ([\phi]) \longrightarrow W^{l_{\xi}([\phi])}_{\xi} ([\phi])/\CC\phi.
$$
Then there is the subspace
$E_{\lambda}$ of $H^0(\EE_{\xi})$ subject to the following properties.

(i) The global section $(\alpha(\phi) -\lambda e_{\xi})$ belongs to $E_{\lambda}$.

(ii) The homomorphism
$$
 e_{\xi} \wedge: H^0 (\EE_{\xi}) \longrightarrow W_{\xi}
$$
maps $E_{\lambda}$ isomorphically onto $\widetilde{V}_{\lambda}$; in particular, $dim(E_{\lambda})=dim(V_{\lambda}) +1$.

(iii) All global sections in $E_{\lambda}$ are proportional to each other; in particular, they give rise to the line subbundle ${\cal L}_{\lambda}$
of $\EE_{\xi}$ fitting into the following diagram	
	\begin{equation}\label{xi-lambd-diag}
	\xymatrix{
		&&0\ar[d]&&\\
		&&{\cal L}_{\lambda} \ar[d] \ar[dr]&&\\
		0\ar[r]& \OO_{C} \ar[r]^{e_{\xi}} \ar[dr]&\EE_{\xi}\ar[r]^(.42){ e_{\xi}\wedge} \ar[d] & \OO_{C}(K_C)\ar[r]&0\\
		&&\OO_{C}(K_C) \otimes {\cal L}^{-1}_{\lambda} \ar[d]&&\\
		&&0&& 
	}
\end{equation}
	Furthermore, the inclusion $H^0({\cal L}_{\lambda}) \hookrightarrow H^0(\EE_{\xi})$ induces the identification
	$$
	H^0({\cal L}_{\lambda})\cong E_{\lambda}.
	$$
	
	Conversely, every line subbundle ${\cal L}$ of $\EE_{\xi}$ with
	$h^0({\cal L})$ at least two gives rise to a nonzero $\phi \in W_{\xi}$ such that ${\cal L}$ corresponds to an eigen value of
	$$
	\widehat{\alpha^{(2)}}(\phi,\bullet):W^{l_{\xi}([\phi])}_{\xi} ([\phi])/\CC\phi \longrightarrow W^{l_{\xi}([\phi])}_{\xi} ([\phi])/\CC\phi.
	$$
	  
	\vspace{0.2cm}
	2) Each eigen value $\lambda$ of $\widehat{\alpha^{(2)}}(\phi,\bullet)$ determines a distinguished global section $\delta_{\lambda}$ of the line bundle $\OO_{C}(K_C) \otimes {\cal L}^{-1}_{\lambda}$ in the diagram \eqref{xi-lambd-diag}. That section has the following properties.
	
	(i) The zero divisor $D_{\lambda}=(\delta_{\lambda}=0)$ of $\delta_{\lambda}$ is nonzero.
	
	(ii) $\delta_{\lambda} H^0({\cal L}_{\lambda}) =\widetilde{V}_{\lambda}$. In particular, the base locus 
	$B_{|\widetilde{V}_{\lambda}|}$ of the linear subsystem $|\widetilde{V}_{\lambda}|$ of $|K_C|$ has the form
	$$
	B_{|\widetilde{V}_{\lambda}|}=D_{\lambda}+B_{\lambda},
	$$
	where $B_{\lambda}$ is the base locus of ${\cal L}_{\lambda}$. In addition, we have the inclusion
	$$
	Z_{\phi} =(\phi=0) \supset B_{|\widetilde{V}_{\lambda}|}=D_{\lambda}+B_{\lambda}.
	$$
	
	(iii) The Kodaira-Spencer class $\xi$ is annihilated by $\delta_{\lambda}$, that is, $\xi$ lies in the kernel of the
	map
	$$
	\xymatrix{
	H^1 (\OO_C(-K_C)) \ar[r]^(.45){\delta_{\lambda}}& H^1 (\OO_C(D_{\lambda}-K_C)).
}
	$$
	Geometrically, this means that the point $[\xi] \in \PP(H^1(\OO_C(-K_C)))=\PP(H^0(\OO_C (2K_C))^{\ast})$ lies in the projective subspace
	$$
	\{D_{\lambda}\}_{2K_C}\cong \PP^{deg(D_{\lambda})-1}=\text{the linear span of $D_{\lambda}$ in the bicanonical embedding of $C$} .
	$$
\end{pro}
\begin{pf}
	The construction of the linear map
	$$
	\alpha^{(2)} (\phi,\bullet): W_{\xi} \longrightarrow \HKC
	$$
	from a splitting
	$$
	\alpha: W_{\xi} \longrightarrow H^0(\EE_{\xi})
	$$
	is a part of the construction of $([\xi],[\phi])$-filtration of $W_{\xi}$. By the construction the last nonzero term of the filtration
	$W^{l_{\xi}([\phi])}_{\xi} ([\phi])$ is preserved by $\alpha^{(2)} (\phi,\bullet)$ and hence induces the endomorphism
	$$
	\alpha^{(2)} (\phi,\bullet):W^{l_{\xi}([\phi])}_{\xi} ([\phi]) \longrightarrow W^{l_{\xi}([\phi])}_{\xi} ([\phi]). 
	$$
	The inclusion 
	$$
	\CC\phi \subset ker(\alpha^{(2)} (\phi,\bullet)) \subset W^{l_{\xi}([\phi])}_{\xi} ([\phi])
	$$
	is again a part of the construction of the filtration. Hence the induced
	endomorphism
	$$
	 \widehat{\alpha^{(2)}}(\phi,\bullet):W^{l_{\xi}([\phi])}_{\xi} ([\phi])/\CC\phi \longrightarrow W^{l_{\xi}([\phi])}_{\xi} ([\phi])/\CC\phi.
	$$
	We now turn to the proof of the properties of this endomorphism listed in the proposition.
	
	Let $V_{\lambda}$ be the eigen space of $\widehat{\alpha^{(2)}}(\phi,\bullet)$ associated to an eigen value $\lambda$, that is, we have
	$$
	\widehat{\alpha^{(2)}}(\phi,v)=\lambda v,\,\forall v\in V_{\lambda}.
	$$
	We now consider the pull back of that equation under the projection
	$$
	W^{l_{\xi}([\phi])}_{\xi} ([\phi]) \longrightarrow W^{l_{\xi}([\phi])}_{\xi} ([\phi])/\CC\phi.
	$$
	Let $\widetilde{V}_{\lambda}$ be the the inverse image of $V_{\lambda}$ under that map. Then the above equation becomes
	\begin{equation}\label{eq:lift}
	\alpha^{(2)}(\phi,\widetilde{v})=\lambda \widetilde{v} +\mu(\widetilde{v}) \phi,\,\forall \widetilde{v}\in \widetilde{V} _{\lambda},
\end{equation}
where $\mu: \widetilde{V}_{\lambda} \longrightarrow \CC$ is some linear function. Observe for $\widetilde{v}=\phi$ the above equation tells us
$$
\mu(\phi)=-\lambda.
$$
We define the subspace $E_{\lambda}$ of $H^0(\EE_{\xi})$
$$
E_{\lambda}:=\{\alpha(\widetilde{v})+\mu(\widetilde{v})e_{\xi} | \widetilde{v}\in \widetilde{V}_{\lambda}\, \}.
$$
This and $\mu(\phi)=-\lambda$ imply that
$$
(\alpha(\phi)-\lambda e_{\xi}) \in E_{\lambda},
$$
the property $(i)$ in 1). From the formula
$$
 e_{\xi}\wedge (\alpha(\widetilde{v})+\mu(\widetilde{v})e_{\xi})=\widetilde{v}, \forall\widetilde{v}\in \widetilde{V}_{\lambda},
 $$ 
the property $(ii)$ follows.

We now consider
$$
(\alpha(\phi)-\lambda e_{\xi}) \wedge (\alpha(\widetilde{v})+\mu(\widetilde{v})e_{\xi})=
\alpha^{(2)} (\phi,\widetilde{v}) -\lambda \tilde{v} -\mu(\tilde{v})\phi=0,
$$
where the last equality follows from the equation \eqref{eq:lift}. Thus all global sections in $E_{\lambda} $ are proportional to $(\alpha(\phi)-\lambda e_{\xi})$ and hence they are proportional to each other. This also shows that the subsheaf of $\EE_{\xi}$ generated by the subspace $E_{\lambda}$ has rank one. Its saturation is the subbundle
${\cal L}_{\lambda}$ of $\EE_{\xi}$ asserted in $(iii)$ of the proposition.
It remains to show the last assertion - the identification $H^0({\cal L}_{\lambda}) \cong E_{\lambda}$.

By definition of ${\cal L}_{\lambda}$ we have the inclusion
$$
E_{\lambda} \subset im \left( H^0({\cal L}_{\lambda}) \hookrightarrow H^0(\EE_{\xi}) \right).
$$
Let $s$ be a global section of ${\cal L}_{\lambda}$.  Denote by $e_s$ the corresponding global section of $\EE_{\xi}$. Its image under the homomorphism
$$
e_{\xi}\wedge: H^0(\EE_{\xi}) \longrightarrow W_{\xi}
 	$$
	will be denoted $\psi$ and we can write
	$$
	e_s=\alpha(\psi) +c e_{\xi},
	$$
	for some constant $c\in \CC$. Since this is proportional to $(\alpha(\phi)-\lambda e_{\xi})$ we obtain
	$$
	0=(\alpha(\phi)-\lambda e_{\xi})\wedge e_s=(\alpha(\phi)-\lambda e_{\xi})\wedge (\alpha(\psi) +c e_{\xi})=\alpha^{(2)}(\phi,\psi) -\lambda \psi -c \phi.
	$$
	Hence
	$$
	\alpha^{(2)}(\phi,\psi) =\lambda \psi +c \phi,
	$$
	which means that $\psi$ belongs to the subspace $\widetilde{V}_{\lambda}$. Thus we have
	$$
	e_s=\alpha(\psi) +c e_{\xi}
	$$
	with $\psi \in \widetilde{V}_{\lambda}$. We compare it with the section
	$$
	\alpha(\psi)+\mu(\psi)e_{\xi} \in E_{\lambda}.
	$$
	Since both global sections are proportional, they come from global section of the line subbundle ${\cal L}_{\lambda}$, we have
	$$
	0=(\alpha(\psi) +c e_{\xi})\wedge (\alpha(\psi)+\mu(\psi)e_{\xi} )=\alpha^{(2)}(\psi,\psi)+c\psi- \mu(\psi)\psi=c\psi- \mu(\psi)\psi,
	$$
	where the last equality uses the skew-symmetry of $\alpha^{(2)}$. Hence
	$c=\mu(\psi)$ and the global section
	$$
	e_s=\alpha(\psi)+\mu(\psi)e_{\xi}
	$$
	lies in $E_{\lambda}$. This proves that under the inclusion
	$$
	H^0({\cal L}_{\lambda}) \hookrightarrow H^0(\EE_{\xi})
	$$
	all global sections of ${\cal L}_{\lambda}$ land in the subspace $E_{\lambda}$. This completes the proof of the isomorphism
	$$
	H^0({\cal L}_{\lambda}) \cong E_{\lambda}.
	$$
	 
	 Conversely, let ${\cal L}$ be a line subbundle of $\EE_{\xi}$ with
	 the space of global sections $H^0({\cal L})$ at least two dimensional.
	 Then the inclusion
	 $$
	 {\cal L} \hookrightarrow \EE_{\xi}
	 $$
	 gives rise to the commutative diagram
	 \begin{equation}\label{xi-L-diag-pf}
	 	\xymatrix{
	 		&&0\ar[d]&&\\
	 		&&{\cal L} \ar[d] \ar[dr]&&\\
	 		0\ar[r]& \OO_{C} \ar[r]^{e_{\xi}} \ar[dr]&\EE_{\xi}\ar[r]^(.42){ e_{\xi}\wedge} \ar[d] & \OO_{C}(K_C)\ar[r]&0\\
	 		&&\OO_{C}(K_C) \otimes {\cal L}^{-1} \ar[d]&&\\
	 		&&0&& 
	 	}
	 	\end{equation}
	with the slanted arrows nonzero. In particular, the one in the top right corner gives the inclusion
	\begin{equation}\label{Lincl-pf}
	H^0({\cal L}) \hookrightarrow \HKC.
\end{equation}
	Denote its image by $W_{\cal L}$. For any $\phi$ in that subspace denote by $s(\phi)$ the corresponding global section of ${\cal L}$. Under the monomorphism of the vertical sequence of the diagram \eqref{xi-L-diag-pf} that global section goes to the one of $\EE_{\xi}$ which has the form
	$$
	\alpha(\phi)-\lambda(\phi) e_{\xi}
	$$
	for some linear function $\lambda$ on  $W_{\cal L}$. Any two of those global sections are proportional in $\EE_{\xi}$. Hence the exterior product
	$$
	\left(\alpha(\phi')-\lambda(\phi') e_{\xi}\right)\wedge \left(\alpha(\phi)-\lambda(\phi) e_{\xi})\right),
	$$
	for all $\phi,\phi'$ in $W_{\cal L}$, goes to zero under the homomorphism
	$$
	\xymatrix{
		\bigwedge^2H^0(\EE_{\xi}) \ar[r]&H^0(det(\EE_{\xi}))=\HKC.
	}
$$
This gives the equation in $\HKC$
$$
\alpha^{(2)}(\phi',\phi)-\lambda(\phi')\phi + \lambda(\phi)\phi' =0, \,\forall \phi, \phi' \in W_{\cal L}.
$$
Fix any nonzero $\phi \in W_{\cal L}$ and rewrite the above equation as follows
$$
\alpha^{(2)}(\phi,\phi')= \lambda(\phi)\phi' -\lambda(\phi')\phi, \, \forall \phi' \in W_{\cal L}.
$$
From this it follows:

- the subspace $W_{\cal L}$ is $\alpha^{(2)}(\phi,\bullet)$-invariant,

- $W_{\cal L} /\CC \phi$ is the $\lambda(\phi)$-eigen space of the induced endomorphism
$$
\widehat{\alpha}^{(2)}(\phi,\bullet): W^{l_{\xi}}_{\xi}([\phi])/\CC \phi \longrightarrow W^{l_{\xi}}_{\xi}([\phi])/\CC \phi.
$$

	Next we turn to part 2) of the proposition. The diagram \eqref{xi-lambd-diag} and the fact that $h^0({\cal L}_{\lambda})$ is  at least two imply that the slanted arrows of the diagram are nonzero. The slanted arrow in the low left corner delivers the distinguished global section $\delta_{\lambda}$ of $\OO_C (K_C) \otimes {\cal L}^{-1}_{\lambda}$ - the image of the constant section $1$ of $\OO_C$.
	This in turn implies that the slanted arrow in the upper right corner
	is the multiplication by $\delta_{\lambda}$
	$$
	\xymatrix{
	{\cal L}_{\lambda}\ar[r]^(.4){\delta_{\lambda}}& \OO_C(K_C).
}
$$
The divisor of zeros $D_{\lambda}=(\delta_{\lambda}=0)$ is nonzero because the defining extension sequence is nonsplit. This proves 2), $(i)$.

On the level of global sections we obtain
$$
\xymatrix{
	H^0({\cal L}_{\lambda})\ar[r]^(.35){\delta_{\lambda}}& W_{\xi} \subset \HKC.
}
$$
From the part 1), $(ii)$ and the identification $H^0({\cal L}_{\lambda}) \cong E_{\lambda}$ in 1), $(iii)$, we deduce the equality
$$
\delta_{\lambda} H^0({\cal L}_{\lambda}) =\widetilde{V}_{\lambda}.
$$
The statements about the base locus is obvious, while the inclusion
$$
Z_{\phi}=(\phi=0) \supset B_{|\widetilde{V}_{\lambda}|}
$$
follows from the fact that $\phi \in \widetilde{V}_{\lambda}$.

The statement in 2), $(iii)$, is the same as in Proposition \ref{pro:xi-Wlatleast2}.
\end{pf}

We can be more precise if we assume that the subspace $W_{\xi}$ contains a nonzero global section 
$\phi$ whose divisor of zeros
$$
Z_{\phi} =(\phi=0)
$$
is in general position in $\PP(\HKC^{\ast})$, that is, $Z_{\phi}$ is a hyperplane section of $C\subset \PP(\HKC^{\ast})$ consisting of $2(g-1)$ distinct points and any subset of $(g-1)$ distinct points of $Z_{\phi}$ spans the hyperplane $H_{\phi}=\{Z_{\phi} \}_{|K_C|}$.

\begin{pro}\label{pro:phi-gp}
	Let $[\xi]\in \Sigma^0_r $ and assume there is $\phi \in W_{\xi} $ whose divisor of zeros 
	$$
	\text{\rm $Z_{\phi}=(\phi=0)$ is in general position in $\PP(\HKC^{\ast})$.}
	$$
	Assume furthermore that the $([\xi],[\phi])$-filtration of $W_{\xi}$
	has the last step $W^{l_{\xi}([\phi])}_{\xi} ([\phi])$ of dimension at least two. Then we have the equality
	$$
	W_{\xi}=\delta H^0 (\OO_C (K_C-D)),
	$$
	where $\delta$ and $\OO_C (D)$ are as in the diagram \eqref{xi-D-diag} and we have
	$$
	h^0 (\OO_C (D))=1.
	$$
	Furthermore, the divisor $D_{\delta}=(\delta=0)$ is a set of $r$ distinct points and its span $\{D_{\delta}\}_{|K_C|}$ in $\PP(\HKC^{\ast})$ is the projective subspace $\PP(W^{\perp}_{\xi})\cong \PP^{r-1}$. In particular, $D_{\delta}$ is the base locus of the linear subsystem
	$|W_{\xi}|$ and $[\xi]$ lies in the variety $Sec_{r}(C_{2K_C})$ of $r$-secants of the bicanonical image of $C$. 
\end{pro}
\begin{pf}
	Fix $\phi \in W_{\xi}$ with the divisor $Z_{\phi}=(\phi=0)$ in general position in $\PP(\HKC^{\ast})$. From here we continue to argue as in the proof of Lemma \ref{lem:Wlatleast2}: we take a splitting
	$$
	\alpha:W_{\xi}\longrightarrow H^0(\EE_{\xi})
	$$
	construct the map
	$$
	\xymatrix{
	\alpha^{(2)}:\bigwedge^2 W_{\xi} \ar[r]& \HKC
}
$$
and evaluate it on $\phi$ to obtain the linear map
$$
\alpha^{(2)}(\phi, \bullet): W_{\xi} \longrightarrow \HKC
$$
which restricts to the endomorphism
$$
\alpha^{(2)}(\phi, \bullet): W^{l_{\xi}([\phi])}_{\xi} ([\phi]) \longrightarrow W^{l_{\xi}([\phi])}_{\xi} ([\phi]);
$$
This enables us to find $\phi'$ in $W^{l_{\xi}([\phi])}_{\xi} ([\phi])$ linearly independent of $\phi$ such that
\begin{equation}\label{phi-phi'-eq}
\alpha^{(2)}(\phi, \phi')=\lambda_0 \phi' -\mu_0 \phi,
\end{equation}
for some constants $\lambda_0, \mu_0$ in $\CC$. We now recall the Koszul cocycle relation \eqref{aplha2-eq}:
$$
\psi \alpha^{(2)}(\phi, \phi') -\phi'\alpha^{(2)}(\phi, \psi) +\phi \alpha^{(2)}(\phi', \psi)=0
$$
which is valid for all $\psi\in W_{\xi}$. This and the equation \eqref{phi-phi'-eq} imply
\begin{equation}\label{phi-phi'-eq2}
\phi \Big( \alpha^{(2)}(\phi', \psi) -\mu_0 \psi\Big) -\phi'\Big( \alpha^{(2)}(\phi, \psi) -\lambda_0 \psi\Big)=0,\, \forall \psi \in W_{\xi}.
\end{equation}
Let $B$ denote the base locus of the pencil $\{\phi,\phi'\}$. This is contained in $Z_{\phi}$ and spans the subspace of dimension at most $(g-3)$.
The general position assumption for $Z_{\phi} $ implies that $deg(B)\leq (g-2)$. From \eqref{phi-phi'-eq2} we deduce that
$(\alpha^{(2)}(\phi, \psi) -\lambda_0 \psi)$ must vanish on the complement
$B^c=Z_{\phi} \setminus B$. This has degree at least $g$ and hence, by general position, $B^c$ spans the hyperplane $\{Z_{\phi}\}$. From this it follows that $(\alpha^{(2)}(\phi, \psi) -\lambda_0 \psi)$ is a scalar multiple of $\phi$:
$$
\alpha^{(2)}(\phi, \psi) -\lambda_0 \psi=\mu(\psi)\phi.
$$
Thus we obtain 
$$
\alpha^{(2)}(\phi, \psi)=\mu(\psi)\phi +\lambda_0 \psi,\,\forall \psi\in W_{\xi}.
$$
Substituting this into the equation
$$
\psi \alpha(\phi)-\phi \alpha(\psi)=\alpha^{(2)}(\phi, \psi) e_{\xi},
$$
see \eqref{alpha-aplha2-eq}, we deduce
$$
\psi \Big(\alpha(\phi)-\lambda_0 e_{\xi}\Big) -\phi \Big(\alpha(\psi)+\mu(\psi) e_{\xi}\Big)=0, \forall \psi\in W_{\xi}.
$$
This tells us that the space of sections
$$
\{\alpha(\psi)+\mu(\psi) e_{\xi} |\,\psi \in W_{\xi}\}
$$
are all proportional and hence generate a line subsheaf of $\EE_{\xi}$. In other words we obtain the diagram 
$$
	\xymatrix{
	&&0\ar[d]&&\\
	&& {\cal O}_C (K_C-D) \ar[d] \ar[rd]^{\delta}&&\\
	0\ar[r]& \OO_{C} \ar[r]\ar[rd]_{\delta}&\EE_{\xi}\ar[r] \ar[d]& \OO_{C}(K_C)\ar[r]&0\\
	&&\OO_{C} (D)\ar[d]&&\\
	&&0&&
}
$$
as in \eqref{xi-D-diag} with the bonus:
$$
h^0(\OO_C (K_C-D)) \geq dim \{\alpha(\psi)+\mu(\psi) e_{\xi} |\,\psi \in W_{\xi}\} =dim (W_{\xi}).
$$
On the other hand we know the inclusion
$$
\delta H^0(\OO_C (K_C-D)) \subset W_{\xi}.
$$
Hence the equality
\begin{equation}\label{Wxi=delta(K-C)}
\delta H^0(\OO_C (K_C-D)) = W_{\xi}.
\end{equation}
From this it follows
$$
\PP^{r-1}\cong \PP(W^{\perp}_{\xi}) = \{D_{\delta}\}
$$
and hence $deg(D_{\delta})\geq r$. 
From the equality \eqref{Wxi=delta(K-C)} we also have
$$
\phi=\delta \tau,
$$
for a unique $\tau \in H^0(\OO_C (K_C-D))$. Hence
$$
Z_{\phi} =(\phi=0)=D_{\delta} + (\tau=0).
$$
Since $Z_{\phi}$ is in general position in $\PP(\HKC^{\ast})$ we deduce that
the subscheme $D_{\delta}$ is reduced and has precisely $r$ points. Thus 
$$
deg(D)=deg(D_{\delta})=r
$$
and the secant space $\PP(W^{\perp}_{\xi})\cong \PP^{r-1}$ intersects $C$ precisely along $D_{\delta}$. 
From the Riemann-Roch for $\OO_C (D)$ and \eqref{Wxi=delta(K-C)} we also deduce
$$
h^0(\OO_C (D))=1.
$$ 
\end{pf}

It will be useful to introduce the following notion.
\begin{defi}\label{def:gpp}
	1) A Kodaira-Spencer class $\xi \in H^1(\Theta_{C})$ is called base point free if it is nonzero and the linear subsystem $|W_{\xi}|$ is base point free.
	
	2)  A Kodaira-Spencer class $\xi \in H^1(\Theta_{C})$ has the  {\rm general position property (gpp)} if it is nonzero and there is $\phi \in W_{\xi}$ whose
	zero locus $Z_{\phi}$ is in general position with respect to the canonical embedding of $C$.
\end{defi}

\begin{rem}\label{rem:gpp-geom}
	Geometrically, the condition that $\xi$ is base point free means that
	the projective subspace $\PP(W^{\perp}_{\xi})$ is disjoint from the canonical image of $C$. 
	
	Assume $C$ is not hyperelliptic. The canonical map is an embedding and we view $C$ as a subvariety of $\PP(\HKC^{\ast})$. A general hyperplane $H$ in $\PP(\HKC^{\ast})$ intersect $C$ along $(2g-2)$ distinct points $Z_H=H \cdot C$ in general position. Choose any projective subspace $\Pi$ of dimension $(r-1)$ in $H$ and disjoint from $Z_H$. Let $W_{\Pi}$
	be the subspace of $\HKC$ annihilating $\Pi$, that is,
	$$
	W_{\Pi}= \text{the space of liner forms on $\HKC^{\ast}$ vanishing on $\Pi$.}
	$$ 
	Any nonzero $\xi \in  H^1(\Theta_{C})$ annihilating $W_{\Pi}$ is a Kodaira-Spencer class of rank at most $r$, it is base point free and has general position property. 	That is we have:
	either the map
	$$
	H^1(\Theta_C) \longrightarrow W^{\ast}_{\Pi} \otimes H^1(\OO_C)
	$$
	is injective, or every nonzero $\xi$ in the kernel of the above map is
	is base point free and has gpp.
\end{rem}

From Proposition \ref{pro:xi-Wlatleast2} we deduce the following.
\begin{pro}\label{pro:bpf-gp}
	Assume $\xi \in  H^1(\Theta_{C})$ is base point free and has general position property. Then $([\xi],[\phi])$-filtration of $W_{\xi}$ is subject to
	$$
	W^{l_{\xi}([\phi])}_{\xi} ([\phi])=\CC\phi
	$$
	for all $[\phi]$ in some Zariski dense open subset of $|W_{\xi}|$.
	
	In particular, the values $\{h^i([\xi],[\phi])\}$ form the partition
	$$
	{\bf h}([\xi],[\phi])=(h^0([\xi],[\phi]), h^1([\xi],[\phi]), \ldots,  h^{l_{\xi}([\phi])-1}([\xi],[\phi]),1)
	$$
	of weight  
	$$
	|{\bf h}([\xi],[\phi])|=
	\sum^{l_{\xi}([\phi])}_{i=0}h^i([\xi],[\phi])=dim(W_{\xi}),
	$$
	 where the parts are weakly decreasing 
	$$
	h^0([\xi],[\phi])\geq h^1([\xi],[\phi])\geq \cdots \geq  h^{l_{\xi}([\phi])-1}([\xi],[\phi]) \geq h^{l_{\xi}([\phi])}([\xi],[\phi])=1.
	$$
\end{pro}

\vspace{0.2cm}
The previous results establish geometric meaning of the last step of $([\xi],[\phi])$-filtrations. Its remaining part, provided that the length  $l_{\xi}([\phi])$ is at least three, can be used to produce quadrics passing through $C \subset \PP(\HKC^{\ast})$. This is a generalization of the argument used in the proof of Theorem \ref{thm:rk1KS}.

\begin{lem}\label{lem:rnc-basis}
	Let $([\xi],[\phi])\in \PP({\cal W}_{\Sigma^0_r })$ and let the length
$l_{\xi}([\phi])$ of $([\xi],[\phi])$-filtration of $W_{\xi}$ be at least three. Then for every nonzero equivalence class $\{\psi\}$ in the quotient space 
$$
W^{l_{\xi}([\phi])-1}_{\xi}([\phi])/W^{l_{\xi}([\phi])}_{\xi}([\phi])
$$
 can be found  $l_{\xi}([\phi])$ vectors
$$
\psi=\psi^{(l_{\xi}([\phi])-1)}, \psi^{(l_{\xi}([\phi])-2)}, \ldots, \psi^{(0)}
$$
in $W_{\xi}$ subject to the following properties

\vspace{0.2cm}
1) $\psi^{(i)} \in W^{i}_{\xi}([\phi]) \setminus  W^{i+1}_{\xi}([\phi])$, for each $i\in [0, l_{\xi}([\phi])-1]$,

\vspace{0.2cm}
2) for every $0\leq i < j \leq l_{\xi}([\phi])-2$, there exists $h_{ij} \in \HKC$ such that the quadrics
$$
Q_{ij}:=\{\psi^{(i)} \psi^{(j+1)} -\psi^{(j)}\psi^{(i+1)} + \phi h_{ij} =0 \}
$$
pass through $C\subset \PP(\HKC^{\ast})$. 

\vspace{0.2cm}
3) All of the above entities are uniquely determined once a splitting map
$$
\alpha: W_{\xi} \longrightarrow H^0(\EE_{\xi})
$$
is chosen.  
\end{lem}
\begin{pf}
	Fix a splitting map
	$$
	\alpha: W_{\xi} \longrightarrow H^0(\EE_{\xi}).
	$$
	This gives us the map
	$$
	\xymatrix{
		\alpha^{(2)}: \bigwedge^2 W_{\xi} \ar[r]& \HKC
	}
$$
which we evaluate at $\phi$ to obtain
	$$
\alpha^{(2)}(\phi,\bullet): W_{\xi} \longrightarrow \HKC.
$$
 Recall that the $([\xi],[\phi])$-filtration comes together with the  linear maps
	$$
	c^{(i)}_{\xi} (\phi,\bullet): W^{i}_{\xi} ([\phi]) \longrightarrow W^{i-1}_{\xi} ([\phi])/ W^{i}_{\xi} ([\phi]), \, \forall i\in [1,l_{\xi}([\phi])-1]
	$$
	with the associated graded
	$$
	gr(	c^{(i)}_{\xi} (\phi,\bullet)): W^{i}_{\xi} ([\phi])/ W^{i}_{\xi} ([\phi]) \longrightarrow W^{i-1}_{\xi} ([\phi]) / W^{i}_{\xi} ([\phi])
	$$
	injective. Furthermore, those maps are determined by the map $	\alpha^{(2)}(\phi,\bullet)$:
	$$
	c^{(i)}_{\xi} (\phi,\bullet) \equiv \alpha^{(2)}(\phi,\bullet) \mod (W^{i}_{\xi} ([\phi])).
	$$
	We now define the collection of vectors as follows. Choose a representative $\psi \in W^{l_{\xi}([\phi])-1}_{\xi} ([\phi])$ of the class $\{\psi\}$ and set
	$$
	\psi^{(l_{\xi}([\phi])-1)} :=\psi.
	$$
	Now propagate through $W_{\xi}$ by applying $\alpha^{(2)}(\phi,\bullet)$. Namely, we define by descending induction
	\begin{equation}\label{induct-psi}
	\psi^{(i-1)}:=\alpha^{(2)}(\phi,\psi^{(i)}),
\end{equation}
	starting with $i=l_{\xi}([\phi])-1$ and ending at $i=1$. This ensures that the string of vectors thus obtained has the property 1) of the lemma.
	
	The quadratic equations in 2) come from the Koszul cocycle relation
	\eqref{aplha2-eq} applied to the triple $\phi, \psi^{(i+1)}, \psi^{(j+1)}$, for $0\leq i<j \leq l_{\xi}([\phi])-2$:
		$$
	\psi^{(j+1)} 	\alpha^{(2)}(\phi,\psi^{(i+1)}) -\psi^{(i+1)} 	\alpha^{(2)}(\phi,\psi^{(j+1)}) + \phi \alpha^{(2)}(\psi^{i+1},\psi^{(j+1)})=0
		$$
		This together with \eqref{induct-psi} give
		$$
		\psi^{(j+1)} \psi^{(i)} -\psi^{(i+1)}\psi^{(j)} + \phi \alpha^{(2)}(\psi^{(i+1)},\psi^{(j+1)})=0
		$$
		as asserted in 2) of the lemma. The assertion in 3) is obvious from the construction.
\end{pf}
\begin{rem}\label{rem-rnc}
	Restricting the quadratic polynomials in Lemma \ref{lem:rnc-basis}, 2), to the hyperplane $H_{\phi}$ in $\PP(\HKC^{\ast})$ corresponding to $\phi$ gives the quadrics
	\begin{equation}\label{quad-in-Hphi}
	Q_{ij}|_{H_{\phi}} =\{\psi^{(j+1)}_{\phi} \psi^{(i)}_{\phi} -\psi^{(i+1)}_{\phi}\psi^{(j)}_{\phi} =0\},
\end{equation}
	passing through the hyperplane section $Z_{\phi} =\{\phi=0\}$; here $\psi^{(\bullet)}_{\phi}$ denotes the image of $\psi^{(\bullet)}$ under the projection
	$$
	W_{\xi} \longrightarrow W_{\xi} /\CC\phi.
	$$
	If we organize the collection
	$$
	\{\psi^{(l_{\xi}([\phi])-1)}_{\phi},\psi^{(l_{\xi}([\phi])-2}_{\phi}, \ldots, \psi^{(1)}_{\phi}, \psi^{(0)}_{\phi} \}
	$$
	in the form of $2\times(l_{\xi}([\phi])-1)$ matrix 
	$$
	M_{\phi}(\{\psi\})=\begin{pmatrix}
		\psi^{(l_{\xi}([\phi])-2)}_{\phi}&\psi^{(l_{\xi}([\phi])-3)}_{\phi}&\cdots&\psi^{(i)}_{\phi}&\cdots&\psi^{(1)}_{\phi}&\psi^{(0)}_{\phi}\\
		\psi^{(l_{\xi}([\phi])-1)}_{\phi}&\psi^{(l_{\xi}([\phi])-2)}_{\phi}&\cdots&\psi^{(i+1)}_{\phi}&\cdots&\psi^{(2)}_{\phi}&\psi^{(1)}_{\phi}
	\end{pmatrix}
$$
we recognize the quadratic equations in \eqref{quad-in-Hphi} as the equations of vanishing of $2\times 2$ minors of that matrix. Here the label $\{\psi\}$ stands for an element in $W^{l_{\xi}([\phi])-1}_{\xi} ([\phi])/W^{l_{\xi}([\phi])}_{\xi}([\phi])$ which gives rise to
the initial element of the collection. It is well known that the vanishing of those minors determines a rational normal curve $R_{l_{\xi}([\phi])-1}$ of degree $(l_{\xi}([\phi])-1)$ in $\PP^{l_{\xi}([\phi])-1}$. Thus we find the following geometric meaning of $W^{l_{\xi}([\phi])-1}_{\xi}/W^{l_{\xi}([\phi])}_{\xi}$:
\begin{equation}
	\begin{gathered}
	\text{\rm the space $\PP(W^{l_{\xi}([\phi])-1}_{\xi}([\phi])/W^{l_{\xi}([\phi])}_{\xi}([\phi]))$ parametrizes a family}\\
	\text{\rm of cones over rational normal curves}\\
		\text{\rm  of degree $(l_{\xi}([\phi])-1)$ containing the divisor $Z_{\phi}=(\phi=0)$.}
		\end{gathered}
\end{equation}
This obviously generalizes to give the similar meaning of $W^{i}_{\xi}([\phi])/W^{i+1}_{\xi}([\phi])$ (always under the assumption $l_{\xi}([\phi])\geq 3$):
\begin{equation}
	\begin{gathered}
		\text{\rm for every $i\in [2, l_{\xi}([\phi])-1]$, the space $\PP(W^{i}_{\xi}([\phi])/W^{i+1}_{\xi} ([\phi]))$}\\
		\text{\rm parametrizes a family of cones over rational normal curves}\\
		\text{\rm  of degree $i$ containing the divisor $Z_{\phi}=(\phi=0)$.}
	\end{gathered}
\end{equation}
\end{rem}

In down to earth terms the constructions above can be summarized as follows:
for a point $([\xi],[\phi])$ of $\PP({\cal W}_{\Sigma^{\circ}_r})$ we have the linear map
$$
 \alpha^{(2)}(\phi, \bullet): W_{\xi} \longrightarrow \HKC
 $$
 together with the maximal invariant subspace, the subspace $W^{l_{\xi}([\phi])}_{\xi}([\phi])$ of the $([\xi],[\phi])$-filtration; to the quotient space
 $W_{\xi}/W^{l_{\xi}([\phi])}_{\xi}([\phi])$ can be attached a nilpotent operator
 $$
 N_{\phi}: W_{\xi}/W^{l_{\xi}([\phi])}_{\xi}([\phi]) \longrightarrow W_{\xi}/W^{l_{\xi}([\phi])}_{\xi}([\phi])
 $$
 with the index of nilpotency $l_{\xi}([\phi])$, that is,
 $$
 l_{\xi}([\phi])=\min \{k \in \NN | N^k_{\phi}=0\};
 $$
 the $([\xi],[\phi])$-filtration is essentially the standard filtration associated to a nilpotent operator:
 $$
 W^{i}_{\xi}([\phi])/W^{l_{\xi}([\phi])}_{\xi}([\phi]) =im (N^i_{\phi}),
 $$
 for $i\in [0, l_{\xi}([\phi])-1]$; the graded pieces
 $$
 W^{i}_{\xi}([\phi])/W^{i+1}_{\xi} ([\phi])
 $$
 for $i\in [0, l_{\xi}([\phi])-1]$, count the number of Jordan blocks of size
 $(i+1)\times (i+1)$ for the Jordan canonical form of $N_{\phi}$. 
 In the next section we show that there is essentially a canonical choice for the nilpotent operator. This involves replacing the filtration
  $W^{\bullet}_{\xi}([\phi])$ by an associated orthogonal decomposition of
  $W_{\xi}$, a sort of Hodge decomposition of $W_{\xi}$. 
 
 \section{Hodge-type decomposition on $W_{\xi}$}
 Let us recall that $\HKC$ carries the Hodge metric:
 $$
 (\psi,\psi'):= \sqrt{-1} \int_C \psi \wedge \overline{\psi'}, \, \forall \psi,\psi' \in \HKC,
 $$
 where $\overline{\psi'}$ denotes the complex conjugate of $\psi'$. This is a Hermitian metric on $\HKC$, see \cite{G-H}. With respect to the Hodge metric the $([\xi],[\phi])$-filtrations of $W_{\xi}$ can be turned into the orthogonal decompositions of $W_{\xi}$. Namely, we fix $([\xi],[\phi])$ in 
 $\PP({\cal W}_{\Sigma^0_r})$ and consider the $([\xi],[\phi])$-filtration
 $$
 W_{\xi}=W^0_{\xi}([\phi]) \supset W^1_{\xi}([\phi]) \supset \cdots \supset W^{l_{\xi}([\phi])}_{\xi} \supset W^{l_{\xi}([\phi])+1}_{\xi}=0.
 $$
 For every subspace $W^i_{\xi}([\phi])$ of the filtration we set
 $$
 \Big(W^i_{\xi}([\phi])\Big)^{\perp},
 $$
  the orthogonal complement of $W^i_{\xi}([\phi])$ in $\HKC$ with respect to the Hodge metric. Define
 $$
 P^i ([\xi],[\phi]) :=W^i_{\xi}([\phi]) \bigcap \Big(W^{i+1}_{\xi}([\phi])\Big)^{\perp}
 $$
 for every $i\in [0, l_{\xi}([\phi])]$. This gives the identification
 \begin{equation}\label{Pi=Gri}
 	P^i ([\xi],[\phi]) \cong W^i_{\xi}([\phi]) / W^{i+1}_{\xi}([\phi]), \, \forall i \in [0, l_{\xi}([\phi])].
 \end{equation} 
Hence the orthogonal decomposition
$$
W^i_{\xi}([\phi])=W^{i+1}_{\xi}([\phi]) \oplus P^i ([\xi],[\phi]).
$$
Continuing in this fashion gives the orthogonal decomposition
 \begin{equation}\label{Wxi-orth}
  W^i_{\xi}([\phi])= \bigoplus^{l_{\xi}([\phi])}_{r=i} P^r ([\xi],[\phi]),
  \end{equation}
  for every $i \in [0,l_{\xi}([\phi])]$.
 \begin{rem}\label{rem:Cinf}
 	The orthogonal decomposition above, as $([\xi],[\phi])$ varies, is only
 	$C^{\infty}$ decomposition. Only the filtration 
 	$ W^{\bullet}_{\xi}([\phi])$ depends holomorphically on the parameters
 	$[\phi]$ and $[\xi]$. This is analogous to the situation of the classical Hodge decomposition of the cohomology of the complex projective variety. 
 \end{rem}
  
Next we recall the map
$$
\alpha^{(2)}_{\xi} (\phi,\bullet) : W_{\xi} \longrightarrow \HKC
$$
which we used to define the $([\xi],[\phi])$-filtration.
 It has the property of having the degree $(-1)$ with respect to the filtration, that is, we have
$$
\alpha^{(2)}_{\xi} (\phi,\bullet) : W^i_{\xi}([\phi]) \longrightarrow W^{i-1}_{\xi}([\phi]),
$$
for every $i\geq 1$. In addition, the maps 
$$
gr (\alpha^{(2)}_{\xi} (\phi,\bullet)):  W^i_{\xi}([\phi])/W^{i+1}_{\xi}([\phi]) \longrightarrow W^{i-1}_{\xi}([\phi]) /W^{i}_{\xi}([\phi])
$$
on the successive quotient spaces are injective, for all $i \in [0,l_{\xi}([\phi])-1]$.
Taking into account the orthogonal decompositions in \eqref{Wxi-orth} we obtain
 $$
 \alpha^{(2)}_{\xi} (\phi,\bullet) :  W^{i}_{\xi}([\phi])= W^{i+1}_{\xi}([\phi]) \oplus P^i ([\xi],[\phi]) \longrightarrow W^{i-1}_{\xi}([\phi]) = \bigoplus^{l_{\xi}([\phi])}_{s=i-1} P^s ([\xi],[\phi]).
 $$
 This gives the decomposition of $\alpha^{(2)}_{\xi} (\phi,\bullet)$ into the direct sum
 $$
 \alpha^{(2)}_{\xi} (\phi,\bullet)=\bigoplus_{s,i} \alpha^{s,i} (\xi,\phi),
 $$
 where the summands $\alpha^{s,i} (\xi,\phi)$ are the `block' maps
 $$
 \alpha^{s,i} (\xi,\phi): P^i ([\xi],[\phi]) \longrightarrow P^s ([\xi],[\phi])
 $$
 obtained by restricting $\alpha^{(2)}_{\xi} (\phi,\bullet)$ to $P^i ([\xi],[\phi])$ and then composing with the projection onto $P^s ([\xi],[\phi])$. From the properties of $\alpha^{(2)}_{\xi} (\phi,\bullet)$
 the blocks are subject to the following.
 \begin{lem}\label{lem:alpha-sum}
 	The components  $\alpha^{s,i} (\xi,\phi)$ are subject to the following conditions.
 	
 	\vspace{0.2cm}
 	1) $\alpha^{s,i} (\xi,\phi)=0$, for all $s \leq i-2$.
 	
 	\vspace{0.2cm}
 	2) $\alpha^{i-1,i} (\xi,\phi): P^i ([\xi],[\phi]) \longrightarrow P^{i-1} ([\xi],[\phi])$ is injective for every $i\in [1,l_{\xi}([\phi])-1]$.
 	
 	\vspace{0.2cm}
 	3) $\alpha^{s,l_{\xi}([\phi])} (\xi,\phi)=0$, for all $s<l_{\xi}([\phi])$.
 \end{lem}
 \begin{rem}\label{rem:choicesplit}
 	It should be recalled that the map $\alpha^{(2)}_{\xi}$ depends on the choice of a splitting of the exact sequence
 	$$
 	\xymatrix{
 	0\ar[r]&H^0(\OO_C)\ar[r]^{e_{\xi}}& H^0(\EE_{\xi}) \ar[r]^(.58){ e_{\xi}\wedge}& W_{\xi} \ar[r]&0.
 }
$$
Namely, if we choose
$$
\alpha_{\xi}: W_{\xi} \longrightarrow H^0(\EE_{\xi}) 
$$
as our reference splitting, then any other splitting $\alpha$ differs from
$\alpha_{\xi}$ by the multiple of $e_{\xi}$ by a linear function
$$
f:W_{\xi} \longrightarrow \CC,
$$
that is we have
$$
\alpha(\psi)-\alpha_{\xi} (\psi)=f(\psi)e_{\xi},\, \forall \psi\in W_{\xi}.
$$
The maps $\alpha^{(2)}$ and $\alpha^{(2)}_{\xi}$ are related as follows
$$
\alpha^{(2)}(\psi,\psi')- \alpha^{(2)}_{\xi}(\psi,\psi')=f(\psi)\psi'- f(\psi')\psi,\,\forall \psi,\psi' \in W_{\xi}.
$$
For $\psi=\phi$ fixed we have
\begin{equation}\label{eq:changesplit}
\alpha^{(2)}(\phi,\psi)- \alpha^{(2)}_{\xi}(\phi,\psi)=f(\phi)\psi- f(\psi)\phi,\,\forall \psi \in W_{\xi}.
\end{equation}
We know that the filtration ${W^{\bullet}_{\xi}([\phi])}$ of $W_{\xi}$ is independent of the splitting. Thus the orthogonal decomposition
$$
W_{\xi}=\bigoplus^{l_{\xi}([\phi])}_{s=0} P^s ([\xi],[\phi])
$$
 is independent of the splitting as well. The restriction of the equation \eqref{eq:changesplit} to the summand $ P^s ([\xi],[\phi])$ gives
 $$
 (\alpha^{(2)}(\phi,\psi^s), \psi^t)- (\alpha^{(2)}_{\xi}(\phi,\psi^s),\psi^t)=(f(\phi)\psi^s- f(\psi)\phi,\psi^t), \forall \psi^s \in P^s, \psi^t \in P^t.
 $$
 In particular, for all $s\neq t$ in the range $[0,l_{\xi}([\phi])-1]$ the above equation becomes
 $$
 (\alpha^{(2)}(\phi,\psi^s), \psi^t)- (\alpha^{(2)}_{\xi}(\phi,\psi^s),\psi^t)=0.
 $$
  This tells us that the blocks
 $$
 \alpha^{t,s}: P^s ([\xi],[\phi]) \longrightarrow  P^t ([\xi],[\phi])
  $$
  for $t\neq s$ in $[0,l_{\xi}([\phi])-1]$ are independent of the splitting.
  
  \vspace{0.2cm}
  For $t=s$ in $[0,l_{\xi}([\phi])-1]$ the equation \eqref{eq:changesplit} becomes
  $$
 (\alpha^{(2)}(\phi,\psi), \psi')- (\alpha^{(2)}_{\xi}(\phi,\psi),\psi')=(f(\phi)\psi- f(\psi)\phi,\psi')=f(\phi)(\psi,\psi'), \forall \psi, \psi' \in P^s,
 $$ 
 where the last equality uses the fact that $\phi$ lies in the component
 $P^{l_{\xi}([\phi])}$ of the orthogonal decomposition. That equation tells us that the two blocks $(\alpha^{(2)})^{s,s}$ and $(\alpha^{(2)}_{\xi})^{s,s}$, for all $s\in [0,l_{\xi}([\phi])-1]$ differ by  the multiple of identity
 $$
 (\alpha^{(2)})^{s,s} - (\alpha^{(2)}_{\xi})^{s,s} =f(\phi) id_{P^s}.
 $$ 
 \end{rem}
 We define
 \begin{equation}\label{N}
 	N(\xi,\phi):  W_{\xi}/W^{l_{\xi}([\phi])}_{\xi} ([\phi]) \longrightarrow W_{\xi}/W^{l_{\xi}([\phi])}_{\xi} ([\phi])
 \end{equation}
using the orthogonal decomposition
\begin{equation}\label{WmodWl-orthdecomp}
W_{\xi}/W^{l_{\xi}([\phi])}_{\xi} ([\phi]) =\bigoplus^{l_{\xi}([\phi])-1}_{s=0} P^s ([\xi],[\phi])
\end{equation}
by setting
\begin{equation}\label{Ndef}
N([\xi],[\phi])=\begin{cases}
	\displaystyle \bigoplus^{l_{\xi}([\phi])-1}_{s=1} \alpha^{s-1,s} (\xi,\phi),& \text{on $\displaystyle \bigoplus^{l_{\xi}([\phi])-1}_{s=1} P^s ([\xi],[\phi])$},\\
	0,& \text{on $ P^0 ([\xi],[\phi])$.}
\end{cases}
\end{equation}
Thus we can describe the action of $N(\xi,\phi)$ on $W_{\xi}/W^{l_{\xi}([\phi])}_{\xi} ([\phi])$ via the following diagram
$$
\xymatrix{
P^{l_{\xi}([\phi])-1}  \ar@/^1pc/[r]^N& P^{l_{\xi}([\phi])-2}  \ar@/^1pc/[r]^N& \cdots \ar@/^1pc/[r]^N&P^{1}\ar@/^1pc/[r]^N&P^{0}\ar@/^1pc/[r]^N&0,
}
$$
where to simplify the notation we omitted to specify $([\xi],[\phi])$.
From the construction of $N$ the following properties are immediate.
\begin{lem}\label{lem:Nprop}
	The linear map
	$$
	N(\xi,\phi): \bigoplus^{l_{\xi}([\phi])-1}_{s=0} P^s ([\xi],[\phi]) \longrightarrow \bigoplus^{l_{\xi}([\phi])-1}_{s=0} P^s ([\xi],[\phi])
	$$
defined in \eqref{Ndef}	has degree $(-1)$ with respect to the grading. It is independent on the choice of a splitting map
$$
\alpha_{\xi}: W_{\xi} \longrightarrow H^0(\EE_{\xi})
$$
which is used to define $\alpha^{(2)}_{\xi}$. Furthermore, the maps
	$$
	N(\xi,\phi): P^s ([\xi],[\phi]) \longrightarrow  P^{s-1} ([\xi],[\phi])
	$$
	are injective for all $s \in [1,l_{\xi}([\phi])-1]$. In particular,
	$	N(\xi,\phi)$ is nilpotent and its index of nilpotency is $l_{\xi}([\phi])$.
\end{lem}

  We use the nilpotent operator $N(\xi,\phi)$ to define ${\bf \mathfrak sl}(2,\CC)$-module structure on the space $W_{\xi}/W^{l_{\xi}([\phi])}_{\xi}([\phi])$. Namely, by Jacobson-Morozov theorem, see \cite{Ch-Gi}, the nilpotent endomorphism $N(\xi,\phi)$ can be completed to an ${\bf \mathfrak sl}(2,\CC)$-triple
 $$
 \{N(\xi,\phi), h(\xi,\phi), M(\xi,\phi)\} \subset End(W_{\xi}/W^{l_{\xi}([\phi])}_{\xi}([\phi])),
 $$
 where $h(\xi,\phi)$ is a semisimple endomorphism of the triple. To simplify the notation, we will omit $(\xi,\phi)$ whenever no ambiguity is likely. Thus in this simplified form the triple is $\{N, h, M\}$ and it is subject to the relations
 $$
 [h, N]=2N,\,\, [h, M]=-2M,\,\,
 h=[N,M].
 $$
 Furthermore, a triple $\{N, h, M\}$ can be chosen to respect the grading of the decomposition \eqref{WmodWl-orthdecomp} of $W_{\xi}/W^{l_{\xi}([\phi])}_{\xi}([\phi])$, that is, $h$ (resp., $M$) is a homomorphism of degree $0$ (resp., $1$) of the graded module on the right side of \eqref{WmodWl-orthdecomp}. 
 
 Let $V_d$ denote an irreducible ${\bf \mathfrak sl}(2,\CC)$-representation of the highest weight $d$, then
 $$
 W^{i}_{\xi}([\phi])/W^{i+1}_{\xi} ([\phi]) \cong Hom_{{\bf \mathfrak sl}(2,\CC)} (V_i, W_{\xi}/W^{l_{\xi}([\phi])}_{\xi}([\phi]))
 $$
 is the multiplicity module for the irreducible $V_i$ to occur in $W_{\xi}/W^{l_{\xi}([\phi])}_{\xi}([\phi])$. The decomposition of 
 $W_{\xi}/W^{l_{\xi}([\phi])}_{\xi}([\phi])$ into irreducible ${\bf \mathfrak sl}(2,\CC)$-representations has the form
 $$
 W_{\xi}/W^{l_{\xi}([\phi])}_{\xi}([\phi]) \cong \bigoplus^{l_{\xi}([\phi])-1}_{i=0} W^{i}_{\xi}([\phi])/W^{i+1}_{\xi} ([\phi]) \otimes V_i.
 $$
 Thus the weights of the decomposition are as follows
$$
 W_{\xi}/W^{l_{\xi}([\phi])}_{\xi}([\phi]) \cong \bigoplus_{\substack{0\leq i \leq l_{\xi}([\phi])-1, \\ 0\leq k \leq i}} (-i+2k)^{h^i ([\xi],[\phi])}, 
 $$
 where $(w)^m$ stands for the weight $w$ occurring with multiplicity $m$.
 We can now define the subspace spanned by vectors of weight $(n-(l_{\xi}([\phi])-1))$ in $W_{\xi}/W^{l_{\xi}([\phi])}_{\xi}([\phi])$: 
 $$
  H^n([\xi],[\phi]):=\{x\in  W_{\xi}/W^{l_{\xi}([\phi])}_{\xi}([\phi]) \Big|\, h (v)=(n-l_{\xi}([\phi])+1)v\} =\!\!\!\!\!\bigoplus_{\substack{-i+2k=n-l_{\xi}([\phi])+1,\\ 0\leq i \leq l_{\xi}([\phi])-1,\\ 0\leq k \leq i}
   } \!\!\!\!\!\!\!\!\!\!\!\!\!\!\!\!(-i+2k)^{h^i ([\xi],[\phi])}
  $$
  
  Hence the weight decomposition 
  \begin{equation}\label{wdecomp}
  W_{\xi}/W^{l_{\xi}([\phi])}_{\xi}([\phi])=H^{\ast}([\xi],[\phi]):=\bigoplus^{2(l_{\xi}([\phi])-1)}_{n=0}  H^n([\xi],[\phi]),
\end{equation}
  the `cohomology' module attached to  $([\xi],[\phi])$. The nilpotent operator
  $$
  N: H^{\ast}([\xi],[\phi]) \longrightarrow H^{\ast}([\xi],[\phi])
  $$
  with respect to this grading appears as an endomorphism of degree $2$: it takes
 $H^n([\xi],[\phi])$ to $H^{n+2}([\xi],[\phi])$, for all $n$ and satisfies  the symmetry:
 \begin{equation}\label{Lefschetz}
 	\begin{gathered}
 \text{\it $N^n: H^{l_{\xi}([\phi])-1-n}([\xi],[\phi]) \longrightarrow H^{l_{\xi}([\phi])-1 +n}([\xi],[\phi])$}\\
 \text{ \it is an isomorphism, for all $n\in [0, l_{\xi}([\phi])-1]$.}
\end{gathered}
\end{equation}
 
 The weight decomposition \eqref{wdecomp} together with the orthogonal decomposition 
  $$
   W_{\xi}/W^{l_{\xi}([\phi])}_{\xi}([\phi])=\bigoplus^{l_{\xi}([\phi])-1}_{s=0}  P^s([\xi],[\phi])
   $$
   in \eqref{WmodWl-orthdecomp} gives the bigrading
   \begin{equation}\label{Hk,n-k}
   	H^{k,n-k}([\xi],[\phi])=H^n \bigcap P^{l-1 -k},
   \end{equation} 
where we omitted the reference to  $([\xi],[\phi])$ to simplify the notation. This gives the orthogonal decomposition
\begin{equation}\label{Hdecomp}
	H^n=\bigoplus_k H^{k,n-k}([\xi],[\phi]), 
\end{equation}
  a sort of Hodge decomposition of the weight spaces. 
 \begin{lem}
 1)	The graded module
 $$
 H^{\ast}([\xi],[\phi]) =\bigoplus^{2(l_{\xi}([\phi])-1)}_{n=0}  H^{n}([\xi],[\phi])
 $$
 admits a bigrading
 $$
 H^{\ast}([\xi],[\phi]) =\bigoplus_{k,n} H^{k,n-k}([\xi],[\phi])
 $$
 
 2) The nilpotent operator
 	$$
 	N: H^{\ast}([\xi],[\phi]) \longrightarrow H^{\ast}([\xi],[\phi])
 	$$
 	has bidegree $(1,1)$, that is, $N$ takes $H^{k,n-k}([\xi],[\phi])$ to
 	$H^{k+1,n-k+1}([\xi],[\phi])$.
 	
 	\vspace{0.2cm}
 	3) For every $n\in [0,l_{\xi}([\phi])-1]$ the map $N^n$ induces an isomorphism
 	$$
 	N^n: H^{k,l-1-n-k}([\xi],[\phi]) \longrightarrow H^{k+n,l-1-k}([\xi],[\phi]).
 	$$
 	
 	4) The Hodge spaces $H^{k,n-k}$ vanish unless $\displaystyle{\frac{n}{2} \leq k\leq max\{n,l-1\}}$.
 \end{lem} 
  \begin{pf}
  	By definition of triple $\{N,h,M\}$ the semisimple endomorphism $h$ preserves the spaces $P^r$ of the orthogonal decomposition \eqref{WmodWl-orthdecomp}. Hence the summand $P^{l-1-k}$ admits the following weight decomposition
  	\begin{equation}\label{wP}
  	P^{l-1-k}=\bigoplus^k_{i=0} (-(l-1)+2k -i)^{m_i}.
  \end{equation}
  	Hence the Hodge space
  	$$
  	H^{k,n-k}=H^n \bigcap P^{l-1-k}=(-(l-1)+2k-(2k-n))^{m_{2k-n}}=(-(l-1)+n)^{m_{2k-n}}
  	$$
  	corresponds to the value $i=2k-n$ in the direct sum \eqref{wP}. Since $i\in [0,k]$, we deduce
  	$$
  	0\leq 2k-n \leq k.
  	$$
  	Hence the range of values for $k$ in 4) of the lemma.
  	
  	The statement 3) follows from 2) and the isomorphism in \eqref{Lefschetz}.
  	
  	For 2) we recall that $N$ has degree $2$ with respect to the weight decomposition and the degree $(-1)$ with respect to the orthogonal decomposition of  $W_{\xi}/W^{l_{\xi}([\phi])}_{\xi}([\phi])$. Hence
  	$$
  	N:H^{k,n-k}=H^n \bigcap P^{l-1-k} \longrightarrow H^{n+2} \bigcap P^{l-1-k-1} =H^{k+1,n-k+1},
  	$$
  	that is, $N$ has the bidegree $(1,1)$ as asserted. 
  \end{pf}
  

 The above is reminiscent of the Hard Lefschetz theorem for the cohomology of a compact K\"ahler manifold, see \cite{G-H}: the operator $N(\xi,\phi)$ plays the role of the exterior multiplication by a K\"ahler form. So it is suggestive to think of it as a `K\"ahler' structure on $W_{\xi}/W^{l_{\xi}([\phi])}_{\xi}([\phi])$. The following statement gives the summary of the above discussion.
 \begin{pro}
 	Let $([\xi], [\phi]) \in \PP({\cal W}_{\Sigma^{\circ}_r})$ and let
 	$$
 	W_{\xi}=W^0_{\xi} ([\phi]) \supset W^1_{\xi} ([\phi]) \supset \cdots \supset W^{l_{\xi}([\phi])}_{\xi}([\phi]) \supset W^{l_{\xi}([\phi])+1}_{\xi}([\phi]) =0,
 	$$
 	be the $([\xi], [\phi])$-filtration of $W_{\xi}$. Then the quotient space
 	$	W_{\xi}/W^{l_{\xi}([\phi])}_{\xi}([\phi])$ has a structure of
 	${\bf \mathfrak sl}(2,\CC)$-module such that
 	$$
 	W_{\xi}/W^{l_{\xi}([\phi])}_{\xi}([\phi]) \cong \bigoplus^{l_{\xi}([\phi])-1}_{i=0} W^{i}_{\xi}([\phi])/W^{i+1}_{\xi} ([\phi]) \otimes V_i
 	$$
 	is the decomposition into the irreducible ${\bf \mathfrak sl}(2,\CC)$-modules, where $V_i$ is the irreducible ${\bf \mathfrak sl}(2,\CC)$ representation of the highest weight $i$. 
 	
 	Geometrically, every copy of $V_i$ in the above decomposition, for $i \geq 1$, gives rise to a cone over a rational normal curve of degree $i$ in $\PP^i$ passing through the hyperplane section $Z_{\phi}=(\phi=0)$ of $C \subset \PP(\HKC^{\ast})$; the projectivization of the multiplicity module
 	$$
 	\PP(W^{i}_{\xi}([\phi])/W^{i+1}_{\xi} ([\phi]))
 	$$
 	is a parameter space for such cones.   
 	
 	The decomposition
 	of $W_{\xi}/W^{l_{\xi}([\phi])}_{\xi}([\phi])$ into the weight spaces
 	$$
 		W_{\xi}/W^{l_{\xi}([\phi])}_{\xi}([\phi]) =H^{\ast}([\xi],[\phi]):=\bigoplus^{2(l_{\xi}([\phi])-1)}_{n=0}  H^n([\xi],[\phi])
 		$$
 	comes along with an endomorphism
 	$$
 	N(\xi,\phi): H^{\ast}([\xi],[\phi]) \longrightarrow H^{\ast}([\xi],[\phi])
 	$$
 	of degree $2$, that is, it sends $H^n([\xi],[\phi])$ to $H^{n+2}([\xi],[\phi])$ and induces an isomorphism
 	$$
 	N^{n}:  H^{l_{\xi}([\phi])-1-n}([\xi],[\phi]) \longrightarrow H^{l_{\xi}([\phi])-1+n}([\xi],[\phi]),
 	$$
 	for every $n\in [0,l_{\xi}([\phi])-1]$.
 	
 	Furthermore, the space $W_{\xi}/W^{l_{\xi}([\phi])}_{\xi}([\phi])$
 	admits the orthogonal decomposition
 	$$
 	\bigoplus^{l_{\xi}([\phi])-1}_{s=0}  P^s([\xi],[\phi])
 	$$
 	on which $N(\xi,\phi)$ acts as an endomorphism of degree $(-1)$.
 	
 	The orthogonal and the weight decompositions define a structure of bigraded module  on $W_{\xi}/W^{l_{\xi}([\phi])}_{\xi}([\phi])$ 
 	$$
 	H^n([\xi],[\phi])=\bigoplus_{\frac{n}{2}\leq k \leq n} H^{k,n-k},
 	$$
 	where $H^{k,n-k}=H^n \bigcap P^{l-1-k}$.
 	
 	The endomorphism $N(\xi,\phi)$ has bidegree $(1,1)$ with the respect to the bigrading.
 
 	We refer to the above properties of $W_{\xi}$ as ${\bf \mathfrak sl}(2,\CC)$-structure of $W_{\xi}$ or, more informally, as {\rm algebraic K\"ahler structure of $W_{\xi}$}.
 \end{pro}
 
 Schematically all of the above can be depicted in an inverted Young diagram of boxes stacked in the columns numbered in the decreasing order from $(l-1):=l_{\xi}([\phi])-1$ to $0$, from left to right - the column  number $s$ represents the summand $P^s([\xi],[\phi])$ of the orthogonal decomposition; the rows (of boxes) of the decreasing length, disposed from top to bottom - they represent the irreducible 
 ${\bf \mathfrak sl}(2,\CC)$-modules of the decomposition. Thus the diagram starts with
 the row of length $l$ repeated $h^{l-1}=dim(P^{l-1})$ number of times,
 follows by the rows of length $(l-1)$ repeated $(h^{l-2}-h^{l-1}) $ and running from the column $(l-2)$ to $0$; continuing to move from left to right and to descend from top to bottom we arrive to the column number $0$ which we complete
 with $(h^{0}-h^{1}) $ lines filled with one box.
 \begin{example}\label{ex:Ydiagram}
 	Consider the following diagram associated to a filtration of $W_{\xi}$
 	$$
 	\begin{tikzpicture}
 		\draw[black,very thick] (0,0) rectangle (0.5,0.5)
 		[xshift=14pt] (0,0) rectangle (0.5,0.5)
 		[xshift=14pt] (0,0) rectangle (0.5,0.5)
 		[xshift=14pt] (0,0) rectangle (0.5,0.5)
 		[xshift=14pt] (0,0) rectangle (0.5,0.5)
 		[xshift=14pt] (0,0) rectangle (0.5,0.5)
 		[yshift=-14pt] (0,0) rectangle (0.5,0.5)
 		[xshift=-14pt] (0,0) rectangle (0.5,0.5)
 		[xshift=-14pt] (0,0) rectangle (0.5,0.5)
 		[xshift=-14pt] (0,0) rectangle (0.5,0.5)
 		[xshift=-14pt] (0,0) rectangle (0.5,0.5)
 		[xshift=-14pt] (0,0) rectangle (0.5,0.5)
 		[xshift=14pt,yshift=-14] (0,0) rectangle (0.5,0.5)
 		[xshift=14pt] (0,0) rectangle (0.5,0.5)
 		[xshift=14pt] (0,0) rectangle (0.5,0.5)
 		[xshift=14pt] (0,0) rectangle (0.5,0.5)
 		[xshift=14pt] (0,0) rectangle (0.5,0.5)
 		[yshift=-14pt] (0,0) rectangle (0.5,0.5)
 		[xshift=-14pt] (0,0) rectangle (0.5,0.5)
 		[xshift=-14pt] (0,0) rectangle (0.5,0.5)
 		[xshift=-14pt] (0,0) rectangle (0.5,0.5)
 		[xshift=-14pt] (0,0) rectangle (0.5,0.5)
 		[yshift=-14pt] (0,0) rectangle (0.5,0.5)
 		[xshift=14pt] (0,0) rectangle (0.5,0.5)
 		[xshift=14pt] (0,0) rectangle (0.5,0.5)
 		[xshift=14pt] (0,0) rectangle (0.5,0.5)
 		[xshift=14pt] (0,0) rectangle (0.5,0.5)
 		[yshift=-14pt] (0,0) rectangle (0.5,0.5)
 		[xshift=-14pt] (0,0) rectangle (0.5,0.5)
 		[xshift=-14pt] (0,0) rectangle (0.5,0.5)
 		[xshift=-14pt] (0,0) rectangle (0.5,0.5)
 		[xshift=14pt,yshift=-14pt] (0,0) rectangle (0.5,0.5)
 		[xshift=14pt] (0,0) rectangle (0.5,0.5)
 		[xshift=14pt] (0,0) rectangle (0.5,0.5)
 		[yshift=-14pt] (0,0) rectangle (0.5,0.5)
 		[xshift=-14pt] (0,0) rectangle (0.5,0.5)
 		[xshift=-14pt] (0,0) rectangle (0.5,0.5)
 		[xshift=14pt,yshift=-14pt] (0,0) rectangle (0.5,0.5)
 		[xshift=14pt] (0,0) rectangle (0.5,0.5)
 		[yshift=-14pt] (0,0) rectangle (0.5,0.5)
 		[yshift=-14pt] (0,0) rectangle (0.5,0.5)
 		[yshift=-14pt] (0,0) rectangle (0.5,0.5);
 	\end{tikzpicture}
 $$
 The diagram has $6$ columns; this tells us that the length of filtration $l=6$ and the columns correspond, from left to right, to the summands
 $P^5,P^4,P^3,P^2,P^1,P^0$ of the orthogonal decomposition of
 $W_{\xi}/W^6_{\xi} ([\phi])$. 
 
 There are two rows of length $6$, meaning that
 $h^5-h^6=dim(P^5)=2$;
 there are $3$ rows of length $5$, hence
 $$
 h^4-h^5=3;
 $$
 one row of length $4$:
 $$
 h^3-h^4=1,
 $$
 two rows of length $3$:
  $$
  h^2-h^3=2,
  $$
  one row of length $2$:
  $$
 h^1-h^2=1,
 $$
 three rows of length $1$:
 $$
 h^0-h^1=3.
 $$   
 \end{example}
  
  We will refer to the diagram associated to the orthogonal decomposition of
the quotient $W_{\xi}/W^l_{\xi} ([\phi])$ according to the rules described above as the {\it decomposition diagram of $W_{\xi}/W^l_{\xi} ([\phi])$}. Reflecting that diagram with respect to the rightmost column, we obtain the usual Young diagram of the partition
$$
(l^{h^{l-1}}, (l-1)^{h^{l-2}-h^{l-1}},\ldots,2^{h^1-h^2}, 1^{h^0-h^1})
$$
where $i^k$ denotes the part $i$ of a partition repeated $k$ times. 
This partition (resp., its Young diagram) will be called the partition (resp. Young diagram) of the $([\xi],[\phi])$-filtration of $W_{\xi}$ and it will be denoted ${\lambda}_{([\xi],[\phi])}$.
\begin{example}
	The partition $\lambda$ of the decomposition diagram in Example \ref{ex:Ydiagram} has the form
	$$
	\lambda=(6^2,5^3,4,3^2,2,1^3);
	$$
	Reflecting the diagram depicted in that example with respect to the rightmost column gives
		$$
	\begin{tikzpicture}
		\draw[black,very thick] (0,0) rectangle (0.5,0.5)
		[xshift=14pt] (0,0) rectangle (0.5,0.5)
		[xshift=14pt] (0,0) rectangle (0.5,0.5)
		[xshift=14pt] (0,0) rectangle (0.5,0.5)
		[xshift=14pt] (0,0) rectangle (0.5,0.5)
		[xshift=14pt] (0,0) rectangle (0.5,0.5)
		[yshift=-14pt] (0,0) rectangle (0.5,0.5)
		[xshift=-14pt] (0,0) rectangle (0.5,0.5)
		[xshift=-14pt] (0,0) rectangle (0.5,0.5)
		[xshift=-14pt] (0,0) rectangle (0.5,0.5)
		[xshift=-14pt] (0,0) rectangle (0.5,0.5)
		[xshift=-14pt] (0,0) rectangle (0.5,0.5)
		[yshift=-14] (0,0) rectangle (0.5,0.5)
		[xshift=14pt] (0,0) rectangle (0.5,0.5)
		[xshift=14pt] (0,0) rectangle (0.5,0.5)
		[xshift=14pt] (0,0) rectangle (0.5,0.5)
		[xshift=14pt] (0,0) rectangle (0.5,0.5)
		[yshift=-14pt] (0,0) rectangle (0.5,0.5)
		[xshift=-14pt] (0,0) rectangle (0.5,0.5)
		[xshift=-14pt] (0,0) rectangle (0.5,0.5)
		[xshift=-14pt] (0,0) rectangle (0.5,0.5)
		[xshift=-14pt] (0,0) rectangle (0.5,0.5)
		[yshift=-14pt] (0,0) rectangle (0.5,0.5)
		[xshift=14pt] (0,0) rectangle (0.5,0.5)
		[xshift=14pt] (0,0) rectangle (0.5,0.5)
		[xshift=14pt] (0,0) rectangle (0.5,0.5)
		[xshift=14pt] (0,0) rectangle (0.5,0.5)
		[xshift=-14pt,yshift=-14pt] (0,0) rectangle (0.5,0.5)
		[xshift=-14pt] (0,0) rectangle (0.5,0.5)
		[xshift=-14pt] (0,0) rectangle (0.5,0.5)
		[xshift=-14pt] (0,0) rectangle (0.5,0.5)
		[yshift=-14pt] (0,0) rectangle (0.5,0.5)
		[xshift=14pt] (0,0) rectangle (0.5,0.5)
		[xshift=14pt] (0,0) rectangle (0.5,0.5)
		[yshift=-14pt] (0,0) rectangle (0.5,0.5)
		[xshift=-14pt] (0,0) rectangle (0.5,0.5)
		[xshift=-14pt] (0,0) rectangle (0.5,0.5)
		[yshift=-14pt] (0,0) rectangle (0.5,0.5)
		[xshift=14pt] (0,0) rectangle (0.5,0.5)
		[xshift=-14pt,yshift=-14pt] (0,0) rectangle (0.5,0.5)
		[yshift=-14pt] (0,0) rectangle (0.5,0.5)
		[yshift=-14pt] (0,0) rectangle (0.5,0.5);
	\end{tikzpicture}
	$$
	the Young diagram of the partition $\lambda$.
\end{example}

We summarize the above discussion in the following.
\begin{pro}
	Let 
	$$
	W^{\bullet}_{\xi}([\phi])=\{W^0_{\xi} ([\phi])\supset W^1_{\xi} ([\phi]) \supset \cdots \supset W^{l-1}_{\xi} ([\phi]) \supset W^l_{\xi} ([\phi]) \supset W^{l+1}_{\xi} ([\phi])=0 \}
	$$
	be an $([\xi],[\phi])$-filtration of $W_{\xi}$. Then the values 
	$$
	h^i=dim(W^i_{\xi}([\phi])/W^{i+1}_{\xi}([\phi]))
	$$
	are positive and form a weakly decreasing function of $i \in [0,l]$. The partition ${\lambda}_{([\xi],[\phi])}$ associated with the filtration has the form
	$$
{\lambda}_{([\xi],[\phi])}=(l^{h^{l-1}}, (l-1)^{h^{l-2}-h^{l-1}},\ldots,2^{h^1-h^2}, 1^{h^0-h^1}).
$$
Its length $l({\lambda}_{([\xi],[\phi])}) $, that is, the number of rows in its Young diagram is the dimension of the quotient space $W^0_{\xi}([\phi])/W^{1}_{\xi}([\phi])$:
$$
l({\lambda}_{([\xi],[\phi])}) =h^0=dim(W^0_{\xi}([\phi])/W^{1}_{\xi}([\phi]))	.
$$
The number of columns of ${\lambda}_{([\xi],[\phi])}$ is $l$, the length of
the filtration 	$W^{\bullet}_{\xi}([\phi])$; the length of the $i$-th column of ${\lambda}_{([\xi],[\phi])}$ is the dimension
$h^{i-1}=dim(W^{i-1}_{\xi}([\phi])/W^{i}_{\xi}([\phi]))$.

The partition ${\lambda}'_{([\xi],[\phi])}$, the conjugate of ${\lambda}_{([\xi],[\phi])}$, is a partition of length $l$, the length of 
the filtration 	$W^{\bullet}_{\xi}([\phi])$:
$$
l({\lambda}'_{([\xi],[\phi])})=l.
$$   
\end{pro}

The simplest situation  arises when $([\xi],[\phi])$-filtrations are the maximal ladders in $W_{\xi}$.
\begin{pro}\label{pro:maxladder}
	Assume $\xi$ is a base point free Kodaira-Spencer class of rank $r$ subject to
	$$
	g-r=dim(W_{\xi})\geq 4.
	$$
	  Assume that for a general $[\phi] \in \PP(W_{\xi})$ the $([\xi],[\phi])$-filtration of $W_{\xi}$ is a maximal ladder
	$$
	W^{\bullet}_{\xi}([\phi])=\{W^0_{\xi} ([\phi])\supset W^1_{\xi} ([\phi]) \supset \cdots \supset W^{l-1}_{\xi} ([\phi]) \supset W^{l}_{\xi} ([\phi]) \supset W^{l+1}_{\xi} ([\phi])=0 \},
	$$
	that is 
	$$
	h^i=dim(W^i_{\xi}([\phi])/W^{i+1}_{\xi}([\phi]))=1, \forall i\in [0,l].
	$$
	Then the following hold.
	
	1) $l=dim(W_{\xi})-1=g-r-1$.
	
	2) The morphism $f_{|W_{\xi}|}: C \longrightarrow \PP(W^{\ast}_{\xi})=\PP^l$
	defined by the linear subsystem $|W_{\xi}| \subset |K_C|$ maps $C$ onto a reduced irreducible curve denoted $C'$. For every $[\phi]$ in $\PP(W_{\xi})$ with $W^{\bullet}_{\xi}([\phi])$ a maximal ladder, the hyperplane section $Z'_{\phi}$ of $C'$ corresponding to $[\phi]$ lies on a rational normal curve of degree $(l-1)$ in the hyperplane $H_{[\phi]}\cong \PP^{l-1}$ corresponding to $[\phi]$. 
	
 3) Let $d'=deg(C')$ and let $M_l ([\xi])$ be the subset of points $[\phi]$ in $\PP(W_{\xi})$ where the $([\xi],[\phi])$-filtration $W^{\bullet}_{\xi}([\phi])$ is a maximal ladder in $W_{\xi}$. Then $M_l ([\xi])$ is a parameter space of
 smooth rational curves with $d'$ marked unordered points. More precisely, there is
 a $M_l ([\xi])$-scheme ${\cal R}_l $ with the structure morphism
 $$
 \pi_{ {\cal R}_l} :{\cal R}_l \longrightarrow M_l([\xi])
 $$
 which is a $\PP^1$-fibration; it comes along with the morphism
 $$
 	\xymatrix{
 	{\cal R}_l \ar[r]^(.4){\sigma} & \PP(W^{\ast}_{\xi})  
 }
 	$$
 	which is part of the diagram
 $$
 	\xymatrix{
 		{\cal R}_l \ar[r]^(.45){\sigma} \,\, \ar[d]_{\pi_{ {\cal R}_l} }  &\PP(W^{\ast}_{\xi})\\
 		M_l([\xi])\ar@{^{(}->}[r]&\PP(W_{\xi}).
 	}
 	$$
 The morphism $\sigma$ embeds  the fibres of $\pi_{ {\cal R}_l}$ as the rational normal curves as in the item 2) of the proposition. The pullback
 $\sigma^{\ast}(C')$ is a divisor in ${\cal R}_l $ intersecting the fibres
 of $\pi_{ {\cal R}_l}$ in subschemes of degree $d'$, the markings on the fibres of the $\PP^1$-fibration $\pi_{ {\cal R}_l}$. 
  
\end{pro}
\begin{pf}
	The first assertion is the summation formula
	$$
	g-r=dim(W_{\xi})=\sum^l_{i=0} dim(W^i_{\xi}([\phi])/W^{i+1}_{\xi}([\phi]))=\sum^l_{i=0} 1=l+1.
	$$
	
	We now turn to the second part of the proposition. For $[\phi] \in \PP(W_{\xi})$ the filtration
	$	W^{\bullet}_{\xi}([\phi])$ is a maximal ladder. In particular,
	$W^l_{\xi}([\phi])=\CC{\phi}$. We now follow the proof of Lemma \ref{lem:rnc-basis}: choose a splitting
	$$
	\alpha:W_{\xi} \longrightarrow \HKC;
	$$
	construct the map
	$$
	\alpha^{(2)}(\phi,\bullet): W_{\xi} \longrightarrow \HKC
	$$
	and use it to obtain a string of vectors
	$$
	{\psi}^{(l-1)}, \psi^{(l-2)}, \ldots, \psi^{(1)}, \psi^{(0)}
	$$
	in $W_{\xi}$ subject to
	$$
	\psi^{(i)}=\alpha^{(2)}(\phi,\psi^{i+1}),
	$$
	for all $i\in [0,l-1]$. From Lemma \ref{lem:rnc-basis}, 2), we obtain quadrics
	\begin{equation}\label{quadrics-phi}
	Q_{ij}=\{ \psi^{(i)}\psi^{(j+1)}-\psi^{(j)}\psi^{(i+1)}+\phi h_{ij} =0\},
	\end{equation}
	passing through the canonical curve $C\subset \PP(\HKC^{\ast})$, for all $0\leq i<j\leq l-1$. 
	Consider the projection of $C$ from the linear subspace $\PP(W^{\perp})$. This gives the map
	$$
	f_{|W_{\xi}|}: C \longrightarrow C' \subset \PP(W^{\ast}_{\xi})=\PP^l
	$$
	defined by the linear subsystem $|W_{\xi}|$. Since $\xi$ is assumed to be base point free, this is a morphism and $C'$ denotes the image of this morphism. Hence $C'$ is reduced and irreducible curve spanning $\PP^l$. Taking the hyperplane $H_{\phi}$ in $\PP^l$ corresponding to $\phi$ the equations \eqref{quadrics-phi} tell us that the hyperplane section
	$$
	Z'_{\phi}=C' \cdot H_{\phi}
	$$
	lies on the quadrics
	$$
	Q'_{ij}=\{ \psi^{(i)}_{\phi}\psi^{(j+1)}_{\phi}-\psi^{(j)}_{\phi}\psi^{(i+1)}_{\phi} =0\},
$$
for all $0\leq i<j \leq l-1$, where $\psi^{(p)}_{\phi}$ denotes the image of $\psi^{(p)}$ under the projection 
$$
W_{\xi}\longrightarrow W_{\xi}/\CC\phi.
$$  
This in turn tells us $Z'_{\phi}$ lies on the subscheme of $H_{\phi}$
determined by the minors of the $2\times (l-2)$ matrix
$$
\begin{pmatrix}
	\psi^{(0)}_{\phi}&\psi^{(1)}_{\phi}&\cdots&\psi^{(l-2)}_{\phi}\\
	\psi^{(1)}_{\phi}&\psi^{(1)}_{\phi}&\cdots&\psi^{(l-1)}_{\phi}
\end{pmatrix}.
$$
That subscheme is known to be a rational normal curve of degree $(l-1)$ in the hyperplane $H_{\phi} \cong \PP^{l-1}$ in $\PP^l$, see \cite{FH}. We denote this curve by $R_{[\phi]}$. 

Next we consider what happens when $[\phi]$ varies in $\PP(W_{\xi})$. Let $M_l([\xi])$ be the subset of $[\phi] \in \PP(W_{\xi})$ for which the $([\xi], [\phi])$-filtration is a maximal ladder in $W_{\xi}$. By the assumption this contains a nonempty Zariski open subset of $\PP(W_{\xi})$. From what has just been proved
follows that $M_l([\xi])$ parametrizes rational normal curves $R_{[\phi]}$ as $[\phi]$ moves in $M_l([\xi])$. Namely, consider the incidence correspondence
$$
{\cal R}_l:=\{([\phi], p) \in M_l([\xi]) \times \PP(W^{\ast}_{\xi})| \,p\in R_{[\phi]} \subset \PP(W^{\ast}_{\xi}) \} \subset M_l([\xi]) \times \PP(W^{\ast}_{\xi}).
$$
It comes equipped with two morphisms, the restrictions to ${\cal R}_l$ of the projections of $M_l([\xi]) \times \PP(W^{\ast}_{\xi})$. Those morphisms are denoted by $\pi_{ {\cal R}_l} $ and $\sigma$ respectively
$$
\xymatrix{
&{\cal R}_l	\ar[dl]_{\pi_{ {\cal R}_l}} \ar[dr]^{\sigma} &\\
M_l([\xi])&&\PP(W^{\ast}_{\xi})
}
$$
The morphism $\pi_{ {\cal R}_l}$ gives to ${\cal R}_l$ a structure of a smooth $\PP^1$-fibration, while $\sigma$ embeds those fibres as rational normal curves of degree $(l-1)$ as asserted in the proposition.

By construction, $\sigma $ maps the fibre of $\pi_{ {\cal R}_l}$ over each point $[\phi]\in M_l([\xi])$ into the hyperplane $H_{[\phi]}$ in $\PP(W^{\ast}_{\xi})$ corresponding to $[\phi]$. This implies that   ${\cal R}_l$ comes along with an embedding into the $\PP^{l-1}$-bundle 
$$
\PP(\Omega_{\PP(W_{\xi})})=\{([\phi],[\tau]) \in \PP(W_{\xi}) \times \PP(W^{\ast}_{\xi}) | \phi(\tau) =0\}:
$$

we have the commutative diagram
$$
\xymatrix{
{\cal R}_l \ar@{^{(}->}[r] \ar[d]_{\pi_{ {\cal R}_l} }&\PP(\Omega_{\PP(W_{\xi})}) \ar[d]^(.58){\pi_{\PP(\Omega_{\PP(W_{\xi})})}}\\
M_l([\xi]) \ar@{^{(}->}[r]&\PP(W_{\xi})
}
$$
where for every $[\phi]\in M_l([\xi])$ the fibre $\pi^{-1}_{ {\cal R}_l} ([\phi])=\PP^1$ is embedded as the rational normal curve $R_{[\phi]}$ into $H_{[\phi]}$, the hyperplane in $\PP(W^{\ast}_{\xi})$ corresponding to $[\phi]$. In addition, composing with the projection
$$
p_{\PP(W^{\ast}_{\xi})}: \PP(\Omega_{\PP(W_{\xi})}) \longrightarrow \PP(W^{\ast}_{\xi})
$$
each hyperplane $H_{[\phi]}$ cuts the curve $C'$ along the hyperplane section
$$
Z'_{\phi}:=C' \cdot H_{[\phi]}
$$
and the rational normal curve $R_{[\phi]}$ passes through $Z'_{\phi}$. In other words we obtain the factorization of $\sigma$ as presented  in the diagram below
$$
\xymatrix{
	{\cal R}_l \ar@/^2.1pc/[rr]^{\sigma} \,\, \ar@{^{(}->}[r] \ar[d]_{\pi_{ {\cal R}_l} } &\PP(\Omega_{\PP(W_{\xi})})\ar[r] \ar[d]^(.58){\pi_{\PP(\Omega_{\PP(W_{\xi})})}} &\PP(W^{\ast}_{\xi})\\
	M_l([\xi])\ar@{^{(}->}[r]&\PP(W_{\xi})&
}
$$
It also follows that ${\cal R}_l$ comes with a distinguished divisor
$$
D_l:=\{([\phi], p) \in {\cal R}_l| \,p\in Z'_{\phi}=C' \cdot H_{[\phi]} \subset R_{[\phi]}  \}
$$   
which is the pullback of $C'$ under the morphism $\sigma$:
$$
D_l=\sigma^{\ast}(C').
$$
In particular, it has the relative degree 
$d'=deg(C') $, the degree of $C'$ in $\PP(W^{\ast}_{\xi})$, that is, the projection $\pi_{ {\cal R}_l}$ restricted to $D_l$ 
$$
\pi_{ {\cal R}_l}: D_l \longrightarrow 	M_l([\xi])
$$
is a finite morphism of degree $d'$. On the other hand the map $\sigma$ restricted to $D_l$ maps it onto the curve $C'$ with at most finite number of points removed:
$$
\sigma: D_l \longrightarrow \sigma(D_l) \subset C'.
$$
with the complement $(C'\setminus \sigma(D_l))$ at most finite set. This is a morphism of relative dimension $(l-1)$: over every point $p'\in \sigma(D_l) $
the fibre $\sigma^{-1}(p')=D_l (p')$ is identified via $\pi_{ {\cal R}_l}$ with the set $(p')^{\perp} \cap M_l([\xi])$, where $(p')^{\perp}$ is the
hyperplane in $\PP(W_{\xi})$ dual to $p'$:
$$
(p')^{\perp}=\{[\psi]\in \PP(W^{\ast}_{\xi}) | p'\in H_{[\psi]} \}.
$$
All the assertions of the proposition are now proved.
\end{pf}

We can draw some geometric consequences from the above considerations.
\begin{cor}\label{cor:Ml-complement}
	Let $M_l([\xi])$ be as in Proposition \ref{pro:maxladder} and let $Z'$ be a subscheme of $(l-1)$ distinct closed points on $C'$ in general position, that is the linear span $\Pi_{Z'}$ of $Z'$ is a linear subspace of codimension two in $\PP^l=\PP(W^{\ast}_{\xi})$. Then the dual line
	$$
	\Pi^{\perp}_{Z'}=\{[\psi] \in \PP(W_{\xi}) | H_{\psi}\supset \Pi_{Z'}\}
	$$
parametrizing the hyperplanes in $\PP(W^{\ast}_{\xi})$ passing through $\Pi_{Z'}$, is not contained in $M_l([\xi])$. This in turn implies that the complement
	$$
	M^c_l([\xi]):=\PP(W_{\xi}) \setminus M_l([\xi])
	$$
	is a subscheme of codimension one in $\PP(W_{\xi})$.
	In addition, every line $\Pi^{\perp}_{Z'}$ which is not entirely contained in 	$M^c_l([\xi])$ defines an irreducible rational surface in $\PP(W^{\ast}_{\xi})$ containing $C'$. 
\end{cor}
\begin{pf}
	Let the line $\Pi^{\perp}_{Z'}$ be as in the statement of the corollary and assume it is contained in $M_l ([\xi])$. Then from Proposition \ref{pro:maxladder} it follows that over $\Pi^{\perp}_{Z'}$ we have the $\PP^1$-fibration
	$$
	\pi_{ {\cal R}_l}: \pi^{-1}_{ {\cal R}_l} (\Pi^{\perp}_{Z'})
	 \longrightarrow \PP^1=	\Pi^{\perp}_{Z'}.
	 $$
	 Denote $\pi^{-1}_{ {\cal R}_l} (\Pi^{\perp}_{Z'})$ by $X$. This is a ruled surface over $\PP^1$ and according to the same proposition
	 it comes with the morphism
	 $$
	 \sigma_X: X \longrightarrow \PP(W^{\ast}_{\xi})=\PP^l,
	 $$
	 the restriction to $X=\pi^{-1}_{ {\cal R}_l} (\Pi^{\perp}_{Z'})$ of the morphism $\sigma$ in Proposition \ref{pro:maxladder}. We also know
	 that the $\PP^1$'s, the fibres of $X$, are embedded as rational normal curves of degree $(l-1)$ into the hyperplanes of the pencil
	 $\Pi^{\perp}_{Z'}$. This means that every rational curve
	 $R_{[\psi]}$, for $[\psi]\in \Pi^{\perp}_{Z'}$, intersects the codimension two subspace $\Pi_{Z'}$ precisely along the subscheme $Z'$. From this it follows that the surface $\sigma_X(X)$ is an irreducible surface of degree $(l-1)$ in $\PP^l$ and the pencil
	 $\{R_{[\psi]}\}_{[\psi]\in \Pi^{\perp}_{Z'}} $ is a pencil of hyperplane sections of  $\sigma_X(X)$. Furthermore, for every point $p'$ in $Z'$, the inverse image $\sigma^{-1}_X (p')$ is a section of the fibration 
	 $$
	 \pi_{ {\cal R}_l}: X
	 \longrightarrow \PP^1=	\Pi^{\perp}_{Z'}.
	 $$
	 Thus it has $(l-1)$ disjoint sections. From
	 $$
	 l-1=g-r -2\geq 4-2=2,
	 $$
	 it follows that $X$ has at least two disjoint sections contracted by
	 $\sigma_X$. Let $E_0$ be a section of $X$ of minimal self-intersection $n_0=E^2_0 \leq 0$ and $f$ is the equivalence class of a fibre of $X$ in the N\'eron-Severi group of $X$. It is well known that those divisor classes form a $\ZZ$-basis of the N\'eron-Severi group of $X$.  Then the morphism $\sigma_X$ is defined
	 by the linear system $|M|$, where $M$ is a divisor on $X$ of the form
	 $$
	 M=(l-1)E_0 +bf,
	 $$
	 for an integer $b$, subject to the equation
	 $$
	 l-1=M^2=(l-1)^2n_0 +2b(l-1).
	 $$
	 Hence
	 \begin{equation}\label{b}
	 2b=1-(l-1)n_0.
\end{equation} 
	 Let $\{E_1,\ldots, E_{l-1}\}$ be sections of $X$ contracted by $\sigma_X$ onto the set $Z'$. Each has the form
	 $$
	 E_i=E_0 +s_if
	 $$
	 for the integer $s_i$ subject to
	 $$
	 0=E_i \cdot M=(E_0 +s_if)\cdot ((l-1)E_0 +bf)=(l-1)n_0 +(l-1)s_i +b.
	 $$
	We deduce
	$$
	s_i=-n_0 - \frac{b}{l-1}, \forall i\in[1,\ldots,l-1].
	$$
	Thus they all have the form
	$$
E=:	E_0 +(-n_0 - \frac{b}{l-1})f.
$$
	Since the sections are mutually disjoint and there at least two of them we also have
	$$
	0=E^2=\left(E_0 +(-n_0 - \frac{b}{l-1})f \right)^2=n_0 +2\left(-n_0 - \frac{b}{l-1}\right)=-n_0 -\frac{2b}{l-1}.
	$$
	Substituting the value of $2b$ from the equation \eqref{b} gives
	$$
	0=-\frac{1}{l-1}
	$$
	which is clearly impossible. This proves the first assertion of the corollary. The argument also supplies the proof for the last statement that every line $\Pi^{\perp}_{Z'}$ not entirely contained in the complement of $M_l([\xi])$ gives rise to an irreducible rational surface containing $C'$. Indeed, let $\Pi^{\perp}_{Z'}$ be such a line.
	Then it intersects $M^c_l([\xi])$ along nonempty finite subscheme, call it $A_{Z'}$. Set 
	$$
	U_{Z'}:=\Pi^{\perp}_{Z'} \setminus A_{Z'}=\PP^1 \setminus A_{Z'} \subset M_l([\xi]).
$$	
	Then going through the first part of the above argument we obtain the
	$\PP^1$-fibration over $U_{Z'}$
	$$
	\pi_{ {\cal R}_l}: \pi^{-1}_{ {\cal R}_l} (	U_{Z'})
	\longrightarrow 	U_{Z'}.
	$$
	Denote $\pi^{-1}_{ {\cal R}_l} (	U_{Z'})$ by $X_{U_{Z'}}$. This a rational smooth connected quasi-projective surface. By Proposition \ref{pro:maxladder} the restriction of the morphism $\sigma$ to $X_{U_{Z'}}$ gives the map
	$$
	\sigma: X_{U_{Z'}} \longrightarrow \PP(W^{\ast}_{\xi})
	$$
	birational onto its image. Taking its closure in $\PP(W^{\ast}_{\xi})$
	$$
	\Sigma_{Z'}:=\overline{\sigma( X_{U_{Z'}})}
	$$
	gives an irreducible rational surface in $\PP(W^{\ast}_{\xi})$. Since
	the hyperplanes $H_{[\psi]}$, for $[\psi]\in U_{Z'}$, intersect $C'$
	along the hyperplane sections $Z'_{[\psi]} $, those move with $[\psi]$ varying in $U_{Z'}$. Hence those hyperplane sections sweep out a nonempty Zariski open subset of $C'$. This implies that $C'$ is contained in the surface $\Sigma_{Z'}$.
	
	We now turn to the statement of the corollary asserting that the complement of $M_l([\xi])$ is a subvariety of codimension one in $\PP(W^{\ast}_{\xi})$. For this consider the normalization $\widetilde{C'}$ of $C'$. This is a smooth projective curve. Take its
	symmetric product $\widetilde{C'}^{(l-1)}$. There is the largest Zariski open subset parametrizing the subschemes of degree $(l-1)$ in $\widetilde{C'}$ which are mapped bijectively by the normalization map
	$$
	\eta:\widetilde{C'} \longrightarrow C'
	$$
	onto $(l-1)$ points of $C'$ which span a projective subspace of codimension two in $\PP(W^{\ast}_{\xi}) $. Denote that subset by $T$ 
	and consider the incidence correspondence
	$$
	L:=\{([Z'], [\psi])\in T\times \PP(W_{\xi}) | [\psi]\in \Pi^{\perp}_{Z'} \}.
	$$
	It comes with two morphisms
	$$
	\xymatrix{
	&L\ar[dl]_{p_T} \ar[dr]^{p_{ \PP(W_{\xi})}}&\\
T& &\PP(W_{\xi}) 
}
$$
the restriction to $L$ of the projections of $T\times \PP(W_{\xi})$ onto each of its factors. The morphism $p_T$ makes $L$ into a $\PP^1$ fibration with the fibre $p^{-1}_T ([Z'])=\Pi^{\perp}_{Z'}$, for every $[Z']\in T$.
In particular, $L$ is a smooth irreducible variety of dimension $l$.
The morphism $p_{\PP(W_{\xi})}$ is the natural inclusion
	$\Pi^{\perp}_{Z'} \subset \PP(W_{\xi})$ on each of the fibres of $p_T$ and those lines sweep out a Zariski dense open subset of $\PP(W_{\xi})$. This tells us that $p_{\PP(W_{\xi})}$ is generically finite over its image. It is actually finite, since for any $[\psi]$ in $\PP(W_{\xi})$ the hyperplane $H_{\psi}$ in $\PP(W^{\ast}_{\xi})$ corresponding to $[\psi]$ intersects $C'$ along $Z'_{[\psi]}$, a zero-dimensional subscheme of $C'$, hence there there can be only finite number of subschemes $Z'$ of degree $(l-1)$ in $\eta^{\ast}(Z'_{[\psi]}) $. 
	
	The first part of the proof tells us that a general line $\Pi^{\perp}_{Z'}$ must intersect the complement of $M_l$ along a finite non empty subscheme. Thus there is a nonempty Zariski open subset
	$T'$ in $T$ parametrizing such lines and over $T'$ the fibration
	$$
	p_T: p^{-1}_T (T') \longrightarrow T'
	$$
	has a multisection. This determines a divisor, call it $D'$, in $p^{-1}_T (T')$. The image of $D'$ under $p_{\PP(W_{\xi})}$  is contained in the complement of $M_l([\xi])$:
	$$
	p_{\PP(W_{\xi})} (D') \subset  M^c_l ([\xi]).
	$$
	Since the dimension of $D'$ is $(l-1)$ and $p_{\PP(W_{\xi})}$ is finite
	we deduce that the complement $M^c_l ([\xi])$ has a component of dimension $(l-1)$.
\end{pf}

The $([\xi],[\phi])$-filtrations for $[\phi] \in M^c_l ([\xi])$ are no longer maximal ladders. By semicontinuity it is easy to understand how the ladders degenerate.
\begin{lem}\label{lem:filt-degenerate}
	Assume the set up of Proposition \ref{pro:maxladder}, that is, a general
	$([\xi],[\phi])$-filtration of $W_{\xi}$ is a maximal ladder in $W_{\xi}$ and $dim(W_{\xi})=g-r \geq 4$. Then for $[\psi]\in  M^c_l ([\xi])$ the filtration $W^{\bullet}_{\xi}([\psi])$ has the form

	$$
	W^{\bullet}_{\xi}([\psi])=\{W_{\xi}=W^0_{\xi} ([\psi])\supset W^1_{\xi} ([\psi]) \supset \cdots \supset W^{l'-1}_{\xi} ([\psi]) \supset W^{l'}_{\xi} ([\psi]) \supset W^{l'+1}_{\xi} ([\psi])=0 \},
	$$
	where the length of the filtration
	$$
	l'=dim(W_{\xi})-dim(W^{l'}_{\xi} ([\psi]))=g-r-dim(W^{l'}_{\xi} ([\psi]))
	$$ 
	and the dimension of $W^{l'}_{\xi} ([\psi])$ is at least two.
\end{lem}

 We can stratify $M^c_l ([\xi])$ according to the length of filtration:
 define the map
 $$
 {\mathfrak{l}} : M^c_l ([\xi]) \longrightarrow \ZZ_{\geq 0}
 	$$ 
	which assigns to a point $[\psi]\in M^c_l ([\xi]) $ the length
	of the filtration $	W^{\bullet}_{\xi}([\psi])$; it takes values in the set
	$$
	[l-1]:=\{0,1,\ldots, l-1\}.
	$$
	If $l'$ is a value of ${\mathfrak{l}}$ we denote by $M_{l'}([\xi])$ the preimage of $l'$ under the map ${\mathfrak{l}}$:
	$$
	 M_{l'}([\xi])={\mathfrak{l}}^{-1}(l').
	 $$
From Lemma \ref{lem:filt-degenerate} it follows that the stratum $M_{l'}([\xi])$ comes along with the vector bundle of rank $(g-r-l')$
formed by the subspaces $W^{l'}_{\xi}([\psi])$ as $[\psi]$ varies in $M_{l'}([\xi])$. More precisely, we have the classifying map
$$
p_{M_{l'}([\xi])}: M_{l'}([\xi]) \longrightarrow {\bf Gr}((g-r-l'), W_{\xi}),
$$
where ${\bf Gr}((g-r-l'), W_{\xi})$ is the Grassmannian of $(g-r-l')$-planes in $W_{\xi}$ and the map is defined by sending the point
$[\psi]\in M_{l'}([\xi])$ into $[W^{l'}_{\xi}([\psi])]$ the point of the Grassmannian corresponding to the subspace $W^{l'}_{\xi}([\psi])$. The aforementioned vector bundle on $M_{l'}([\xi])$ is the pullback of the universal bundle ${\cal U}$ of the Grassmannian under the map $p_{M_{l'}([\xi])}$:
$$
{\cal W}^{l'} ([\xi]):=p^{\ast}_{M_{l'}([\xi])} ({\cal U}).
$$
We recall the inclusion
$$
\CC{\psi} \subset W^{l'}_{\xi}([\psi]), \,\,\forall [\psi]\in  M_{l'}([\xi])
$$
and the endomorphism
$$
\alpha(\psi, \bullet): W^{l'}_{\xi}([\psi]) \longrightarrow W^{l'}_{\xi}([\psi]).
$$
This translates into the inclusion
$$
\OO_{M_{l'}([\xi])} (-1):=\OO_{\PP(W_{\xi})} (-1)\otimes \OO_{M_{l'}([\xi])} \hookrightarrow {\cal W}^{l'} ([\xi]).
$$
and the morphism
$$
\alpha_{M_{l'}([\xi])}:\OO_{M_{l'}([\xi])} (-1) \longrightarrow {\cal{E}}nd({\cal W}^{l'} ([\xi])).
$$
Thus on every substratum we can use the dictionary between the spectral data of $\alpha_{M_{l'}([\xi])}$ and the line subbundles of $\EE_{\xi}$ as described in Proposition \ref{pro:xi-Wlatleast2-eigen}. The next example will be an illustration of the above considerations.
\begin{example}\label{ex:g=6-ladder}
	Let $C$ be a smooth projective curve of genus $g=6$. We assume that it has no $g^1_2$, $g^1_3$, $g^2_5$. By Theorem \ref{thm:rk1KS} we are in control of the strata $\Sigma_0$ and $\Sigma_1$: they are respectively empty and the bicanonical curve. We consider the stratum $\Sigma_2 \setminus Sec_2 (C_{2K_C})$ consisting of points $[\xi]$ of rank precisely $2$ and which are base point free, that is, the linear subsystem $|W_{\xi}|$ is base point free: the morphism
	\begin{equation}\label{ex6}
	f_{|W_{\xi}|} : C\longrightarrow C' \subset \PP(W^{\ast}_{\xi})=\PP^3
\end{equation}
	is the projection from the line $L_{\xi}:=\PP(W^{\perp}_{\xi})$ in $\PP^5=\PP(\HKC^{\ast})$. In addition, we assume the set-up of Proposition \ref{pro:maxladder}, that is for a general
	$[\phi] \in \PP(W_{\xi})$ the filtration 
	$$
		W^{\bullet}_{\xi}([\phi])=\{W_{\xi}=W^0_{\xi} ([\phi])\supset W^1_{\xi} ([\phi]) \supset W^{2}_{\xi} ([\phi]) \supset W^{3}_{\xi} ([\phi]) \supset W^{4}_{\xi} ([\phi])=0 \}
	$$
	is a maximal ladder in $W_{\xi}\cong \CC^4$, that is $l$ in Proposition \ref{pro:maxladder} equals three. We denote by $M_3 ([\xi])$
	the Zariski open subset of $\PP(W_{\xi})$ parametrizing such points $[\phi]$. For every $[\phi]\in M_3 ([\xi])$, the proposition tells us that the plane
	$H_{[\phi]}$ contains the smooth conic $R_{[\phi]}$ passing through the
	plane section $Z_{[\phi]}=C' \cdot H$. Hence the rational
	surface in Corollary \ref{cor:Ml-complement} arising from a general secant line $\Pi_{Z'}$ of $C'$ is either a quadric or a cubic. The latter can only occur if the degree of $C'$ at most $6$, since in this case a general plane section of $C'$ is contained in the complete intersection of an irreducible conic and an irreducible cubic. This in turn can only occur if the projection morphism $f_{|W_{\xi}|} $ in \eqref{ex6} is of degree $2$ onto its image and then $C'$ is a space curve of degree $5$. 
	That curve can be either a smooth curve of genus $2$, the Castelnuovo bound for the genus of a space curve of degree $5$, or a projection from a point of a normally embedded curve $\widetilde{C'}$ of genus $1$ of degree $5$ in $\PP^4$. The two possibilities
	
	- $C'$ is a smooth curve of genus two,
	
	- $C'$ is a smooth curve of genus one,
	
	\noindent
	can not occur: in both cases $C'$ is the quotient of $C$ by an involution; we already encountered such a situation, see the proof of Lemma \ref{lem:proj-birational}; arguing in the same way the cases are discarded. Thus we are left with the only possibility:
	$$
	\begin{gathered}
	\text{\rm $C'$ is singular quintic curve: it is obtained from $\widetilde{C'}$ by projecting from a point}
	\\
	\text{\rm lying on a (unique) secant line of  $\widetilde{C'}$ but not on $\widetilde{C'}$.}
\end{gathered}
	$$
It is known that such $C'$ {\rm must} lie on a unique quadric surface, see \cite{Hu}. Conversely, a quintic curve in $\PP^3$ lying on quadric has an irreducible cubic passing through it: this is seen by estimating the number of cubics passing through such a curve. Thus we obtain: 
	 \begin{equation}\label{ex:deg2}
	 	\begin{gathered}
	 	\text{\rm in the case 
	 $f_{|W_{\xi}|} $ is of degree $2$, its image, the curve $C'$, is the component}\\
 \text{\rm  residual to a line of a complete intersection of quadric and cubic surfaces in $\PP^3$;}
 \\
 \text{\rm in addition, $C'$ is a projection of a normal quintic curve of genus $1$ in $\PP^4$;}
 \\
 \text{\rm the center of projection is a point lying on a unique secant line}
 \\
 \text{\rm of a normal quintic curve but not on that curve.} 
\end{gathered}
 \end{equation}

	  In case $f_{|W_{\xi}|} $ is of degree one onto its image, the curve $C'$ is a space curve of degree $10$ and from the analysis above the rational surface in Corollary \ref{cor:Ml-complement} is a quadric. Since it is unique it is independent of a secant line $\Pi_{Z'}$ in that corollary, in other words the quadric surface is the closure of the image of the morphism
	  $$
	  \sigma: {\cal R}_3 \longrightarrow \PP(W^{\ast}_{\xi})=\PP^3
	  $$
	  in Proposition \ref{pro:maxladder}.
	  
	  Let us now examine what can be gained from the substrata $M_{l'}([\xi])$, for $l'\in [0,1,2]$, contained in the complement of $M_{3}([\xi])$ in $\PP(W_{\xi})$. Those strata give rise to eigen spaces of the linear maps 
	  $$
	  \widehat{\alpha^{(2)}}(\psi,\bullet): W^{l'}_{\xi} ([\psi]) /\CC\psi \longrightarrow W^{l'}_{\xi} ([\psi]) /\CC\psi,
	  $$
	   for every $[\psi]\in M_{l'}([\xi])$. An eigen space $V_{\lambda}$ of dimension $m_{\lambda} $ produces the subbundle ${\cal L}_{\lambda}$ of $\EE_{\xi}$ with $h^0({\cal L}_{\lambda})=m_{\lambda}+1$,  that is
	   $\EE_{\xi}$ fits into the exact sequence
	   $$
	   \xymatrix{
	   	0\ar[r]&{\cal L}_{\lambda}\ar[r]& \EE_{\xi} \ar[r]&\OO_C (K_C)\otimes {\cal L}^{-1}_{\lambda} \ar[r]&0.
	   }
   $$
   This together with the defining extension sequence
   $$
   \xymatrix{
   	0\ar[r]&\OO_C\ar[r]& \EE_{\xi} \ar[r]&\OO_C (K_C)\ar[r]&0
   }
   $$
   gives the diagram
   \begin{equation}\label{ex:diagram}
    \xymatrix{
    	&&0\ar[d]&&\\
    	&&{\cal L}_{\lambda}\ar[d]\ar[dr]& &\\
   	0\ar[r]&\OO_C\ar[r]\ar[dr]& \EE_{\xi} \ar[r]\ar[d]&\OO_C (K_C)\ar[r]&0\\
   	&&\OO_C (K_C)\otimes{\cal L}^{-1}_{\lambda}\ar[d]&&\\
   	&&0&& 
   }
\end{equation}
	    Since we assuming that $\xi$ is base point free, this imposes the upper bound
	   \begin{equation}\label{ex:mult-eigen}
	   m_{\lambda} \leq 2,
	\end{equation}
	   for all eigen spaces. 
	   
	   The stratum  $M_2([\xi])\setminus(M_0([\xi])\cup M_1([\xi])) $ consists of filtrations of the form
	   $$
	   W_{\xi}=W^0_{\xi}([\psi])\supset W^1_{\xi}([\psi]) \supset W^2_{\xi}([\psi]) \supset W^3_{\xi}([\psi])=0, \forall [\psi]\in M_2([\xi])\setminus(M_0([\xi])\cup M_1([\xi])),
	   $$
	    where $ W^2_{\xi}([\psi])$ is two dimensional; so it produces only one dimensional eigen spaces: 
	     $$
	    V_{\lambda}=W^2_{\xi}([\psi])/\CC\psi,
	    $$
	    for every point $[\psi]\in M_2([\xi])\setminus(M_0([\xi])\cup M_1([\xi]))$.
	     The subspace $W^2_{\xi}([\psi])$ is identified with  the dimension two
	   space $H^0({\cal L}_{\lambda})$. This implies that the quotient line bundle has at least three dimensional space of global sections. The assumption no $g^1_2$, $g^1_3$, $g^2_5$ implies
	   $$
	   deg({\cal L}_{\lambda})\geq 4,\, deg(\OO_C (K_C)\otimes{\cal L}^{-1}_{\lambda} ) \geq 6.
	   $$ 
	   Since those degrees add up to $10$, the above inequalities must be equalities. We write $\OO_C (D_{\lambda})$ for $\OO_C (K_C)\otimes{\cal L}^{-1}_{\lambda}$ and $\OO_C (K_C-D_{\lambda})$ for ${\cal L}_{\lambda}$. The slanted arrows of the diagram above give a distinguished global section $\delta_{\lambda}$ of $\OO_C (D_{\lambda})$ and identifies the subspace $W^2_{\xi}([\psi])$ as
	   $$
	   W^2_{\xi}([\psi])=\delta_{\lambda} H^0(\OO_C (K_C-D_{\lambda})) \subset W_{\xi} \subset \HKC,
	   $$
	   that is the divisor $D_{\delta_{\lambda}}=(\delta_{\lambda}=0)$ of zeros of $\delta_{\lambda}$ spans the three-dimensional subspace
	   $\{D_{\lambda}\}_{K_C}$ containing the line $L_{\xi}$. Under the projection from $L_{\xi}$ the span $\{D_{\lambda}\}_{K_C}$ maps onto a line and the image $D'_{\lambda}$ of $D_{\lambda}$ is the intersection of that line with $C'$. Since the degree of $D'_{\lambda}$ is at least three, the line must be a ruling of the quadric containing $C'$.
	   
	   Observe that we can replace $\psi$ above by {\rm any} nonzero vector $\psi'$ in $ W^2_{\xi}([\psi])$ and the subspace will be invariant with respect to $\alpha_{\xi}(\psi',\bullet)$, that is the  projective line $ \PP(W^2_{\xi}([\psi]))$ is contained in the stratum $M_2([\xi])$.  In addition, such a line corresponds to the line subbundle of $\EE_{\xi}$ fitting into the exact sequence 
	   \begin{equation}\label{ex:anotherext}
	   	\xymatrix{
	   		0\ar[r]&\OO_C (K_C-{D}_{\lambda})\ar[r]& \EE_{\xi} \ar[r]&\OO_C ({D}_{\lambda})  \ar[r]&0,
	   	}
	   \end{equation} 
	  the vertical sequence of the diagram \eqref{ex:diagram}. The subbundle and the quotient bundle of the sequence define respectively $g^1_4$ and $g^2_6$ on $C$.
	    Since there can be only finite number of such subbundles in $\EE_{\xi}$ we deduce:
	   \begin{equation}\label{M2onlylines}
	   \begin{gathered}
	   \text{\rm the closure of $M_2([\xi])\setminus(M_0([\xi])\cup M_1([\xi]))$ in $\PP(W_{\xi})$ }
	   \\
	   \text{\rm consists
	   	of at most finite number of $\PP^1$'s;}
   	\\
   	\text{\rm in particular, it is at most of dimension $1$.} 
   \end{gathered}
	\end{equation}
	   In Corollary \ref{cor:Ml-complement} we have learned that the complement of $M_3([\xi])$ in $\PP(W_{\xi})=\PP^3$ must be of dimension two. This and \eqref{M2onlylines} imply that the stratum $M_1([\xi])$ must have a two dimensional component and furthermore, there must be $[\psi]\in M_1([\xi]$ with the endomorphism
	   $$
	   \widehat{\alpha^{(2)}}(\psi,\bullet): W^{l'}_{\xi} ([\psi]) /\CC\psi \longrightarrow W^{l'}_{\xi} ([\psi]) /\CC\psi,
	   $$	   
	   having an eigen space $V_{\mu}$ of dimension two. Then
	   $\widetilde{V}_{\mu}$, the preimage of $V_{\mu}$ under the projection morphism
	   $$
	   W^{l'}_{\xi} ([\psi])  \longrightarrow
	   W^{l'}_{\xi} ([\psi]) /\CC\psi,
	   $$
	   gives the plane $\PP(\widetilde{V}_{\mu}) \subset M_0([\xi]) \cup M_1([\xi])$. 
	    We unravel the geometry of the eigen space $\widetilde{V}_{\mu}$: it gives us the line subbundle ${\cal L}_{\mu}$ fitting into the diagram \eqref{ex:diagram} with ${\cal L}_{\lambda}$ replaced by ${\cal L}_{\mu}$.   That line subbundle
	    gives rise to $g^2_d$ on $C$, where $d=deg({\cal L}_{\mu})$. By assumption
	   $C$ has no $g^2_5$, so $d\geq 6$. Hence the quotient bundle has degree
	   at most $4$. From the vertical sequence of the diagram we also know that the space of sections of the quotient bundle is at least two dimensional. Since $C$ has no $g^1_m$ for $m\leq 3$, we deduce
	   $$
	   \text{$deg(\OO_C (K_C)\otimes{\cal L}^{-1}_{\mu})=4$ and $h^0(\OO_C (K_C)\otimes{\cal L}^{-1}_{\mu})=2$.}
	   $$
	   As before we write $\OO_C(D_{\mu})$ for $\OO_C (K_C)\otimes{\cal L}^{-1}_{\mu}$. Then ${\cal L}_{\mu}=\OO_C(K_C-D_{\mu})$ and obtain the extension
	   \begin{equation}\label{ex:mu-ext}
	   	\xymatrix{
	   		0\ar[r]&\OO_C (K_C-{D}_{\mu})\ar[r]& \EE_{\xi} \ar[r]&\OO_C ({D}_{\mu})  \ar[r]&0.
	   	} 	   
	   \end{equation}	  
	   Only now the roles of sub- and quotient bundles are switched: the first one defines $g^2_6$ and the second $g^1_4$. The slanted arrow in the lower left corner of the diagram \eqref{ex:diagram} singles out
	   a nonzero global section $\delta_{\mu} \in H^0(\OO_C(D_{\mu}))$, then the slanted arrow in the upper right corner of the diagram identifies $\widetilde{V}_{\mu}$ as the subspace of hyperplanes in $\PP(\HKC^{\ast})$ passing through the line 
	   $L_{\xi}=\PP(W^{\perp}_{\xi})$ and cutting out the divisor
	   $D_{\delta_{\mu}}=(\delta_{\mu}=0)$. Thus we deduce that
	   the secant plane $\{D_{\delta_{\mu}}\}_{K_C}$ of the canonical curve $C$ intersects $C$ along the divisor $D_{\delta_{\mu}}$ of degree $4$ and under the morphism $f_{|W_{\xi}|} $ that divisor is mapped onto a singular point of $C'$. 
	   Also observe, that the subbundle is the maximal destabilizing for $\EE_{\xi}$. Hence it is unique. This implies that the plane $\PP(\widetilde{V}_{\mu})$ is unique as well. Furthermore,
	   if we have two extensions \eqref{ex:anotherext} and \eqref{ex:mu-ext} simultaneously, then both split and we have the equality
	   $$
	   D_{\lambda}=K_C -D_{\mu}.
	   $$
	    Thus we deduce 
	    \begin{equation}\label{ex:complementM3}
	    	\begin{gathered}
	    	\text{\rm either the complement $\PP(W_{\xi})\setminus M_{3}([\xi])$ is the plane $\PP(\widetilde{V}_{\mu}) $,}
	    	\\
	    	\text{\rm or the sequence \eqref{ex:mu-ext} splits and 
	    		$
	    		\EE_{\xi} =\OO_C (D_{\mu})\oplus \OO_C (K_C-D_{\mu}).
	    		$
	    	}
	    		\end{gathered} 
	    \end{equation}
    
     Assume the destabilizing sequence \eqref{ex:mu-ext} is nonsplit. However, on the level of global sections 
$$
 \xymatrix{
	0\ar[r]&H^0 (\OO_C (K_C-{D}_{\mu}))\ar[r]& H^0(\EE_{\xi}) \ar[r]&H^0 (\OO_C ({D}_{\mu}))  \ar[r]&0
} 	   
$$	  
it continues to be exact. This means the extension class $\eta \in H^1(\OO_C (K_C -2D_{\mu}))$ corresponding to the extension sequence \eqref{ex:mu-ext} lies in the kernel of the linear map
\begin{equation}\label{ex:ker}
H^1(\OO_C (K_C -2D_{\lambda})) \longrightarrow (H^0 (\OO_C ({D}_{\mu})))^{\ast} \otimes H^1(\OO_C (K_C -D_{\mu}))
\end{equation}	   
induced by the multiplication morphism
$$
H^1(\OO_C (K_C -2D_{\mu})) \otimes H^0 (\OO_C ({D}_{\mu})) \longrightarrow  H^1(\OO_C (K_C -D_{\mu})).
$$
Dualizing the map \eqref{ex:ker} we see that the kernel of that map is dual to the cokernel of the multiplication map
$$
Sym^2 H^0(\OO_C (D_{\mu})) \longrightarrow H^0(\OO_C (2D_{\mu})).
$$
This is injective and the dimension of the space on the right is at most $4$. Hence we obtain
\begin{equation}\label{ex:nonsplit}
	\text{\rm the extension \eqref{ex:mu-ext} is nonsplit if and only if $h^0(\OO_C (K_C -2D_{\mu}))=1$.}
\end{equation}
In this situation one obtains more geometry: the quadric containing the curve $C'$, the image of $C$ under the projection $f_{|W_{\xi}|}$ from the line $L_{\xi}$ is a quadric cone with the vertex being the image of the divisor $D_{\delta_{\mu}}=(\delta_{\mu}=0)$.

	  In conclusion we deduce:
	  \begin{pro}\label{ex:6-ladder}
	  	 Let $C$ be a smooth projective curve of genus $g=6$ and let $\xi$ be a base point free Kodaira-Spencer class of rank $2$. Assume that
	  	$M_3 ([\xi])$ contains a Zariski dense open subset of $\PP(W_{\xi})=\PP^3$. Then the vector bundle $\EE_{\xi}$ fits into the destabilizing extension sequence
	  	$$
	  	\xymatrix{
	  		0\ar[r]&\OO_C (K_C-{D}_{\mu})\ar[r]& \EE_{\xi} \ar[r]&\OO_C ({D}_{\mu})  \ar[r]&0.
	  	} 	   
  	$$
  	where the subbundle and quotient bundle define respectively $g^2_6$ and $g^1_4$ on $C$. Furthermore, there is a distinguished global section $\delta_{\mu} \in H^0(\OO_C ({D}_{\mu}))$ such that
  	$$
  	 \widetilde{V}_{\mu}:=\delta_{\mu} H^0 (\OO_C (K_C-{D}_{\mu})) \subset W_{\xi}
  	 $$
  	 and the plane 
  	$\PP(\widetilde{V}_{\mu})$ is contained in the substrata $M_0 ([\xi])\cup M_1 ([\xi])$ of the length stratification of the complement of $M_3 ([\xi])$ in $\PP(W_{\xi})$. In addition, the above extension splits unless the complement
  	$$
  	 \PP(W_{\xi}) \setminus M_3 ([\xi]) =\PP(\widetilde{V}_{\mu}).
  	 $$
  	 Under this condition the extension might be nonsplit and then it corresponds to a unique, up to a nonzero scalar multiple, cohomology class lying in the kernel of the map
  	 $$
  	 H^1(\OO_C (K_C -2D_{\mu})) \longrightarrow (H^0 (\OO_C ({D}_{\mu})))^{\ast} \otimes H^1(\OO_C (K_C -D_{\mu})).
  	 $$
  	Geometrically, the morphism 
  	$$
  	f_{|W_{\xi}|}: C\longrightarrow C' \subset \PP(W^{\ast}_{\xi})\cong \PP^3
  	$$
  	is the projection from the line $L_{\xi}=\PP(W^{\perp}_{\xi})$. The image curve $C'$ is contained in a quadric surface. Furthermore, if the
  	destabilizing extension is nonsplit that quadric is singular with the vertex at the point $p_{\xi}$, the image of the divisor $(\delta_{\mu}=0)$.  
	  \end{pro}  
\end{example}

\vspace{0.2cm}
We return to the general considerations. To each point $([\xi], [\phi])$ of $\PP({\cal W}_{\Sigma^0_r})$ we have attached the $([\xi], [\phi])$-filtration. We have accumulated a certain evidence that the properties of those filtrations are related to the geometry of the curve.
The numerical data of the filtrations, such as their length and shape, obviously play a role in understanding of this relation. In addition, these data induce further stratifications on $\PP({\cal W}_{\Sigma^0_r})$: the stratification according to the length of filtrations and then further into substrata of constant partition diagrams. More precisely, define the {\it filtration length function}
$$
{\mathfrak{l}}:\PP({\cal W}_{\Sigma^0_r}) \longrightarrow \ZZ_{\geq 0}
$$
which attaches to each point $([\xi], [\phi])$ the length
$$
 {\mathfrak{l}}([\xi], [\phi]):=l_{\xi}([\phi])
 $$
 of the filtration $W^{\bullet}_{\xi}([\phi])$. The values of ${\mathfrak{l}}$ belong to the set 
 $$
 [g-r-1]=\{0,1,\ldots, g-r-1\}.
 $$
 If $l$ is a value of ${\mathfrak{l}}$ its preimage
 ${\mathfrak{l}}^{-1}(l)$ will be denoted by ${\mathfrak L}_l$. This is the stratum parametrizing the points $([\xi], [\phi])$ having the filtration
 $W^{\bullet}_{\xi}([\phi])$ precisely of length $l$.

 Each stratum ${\mathfrak L}_l$ can be further partitioned into substrata where
 the dimension vector
 $$
 {\bf h}([\xi], [\phi])=(h^0([\xi], [\phi]),h^1([\xi], [\phi]),\ldots, h^l([\xi], [\phi]))
 $$
 remains constant; we recall that $h^i([\xi], [\phi])$ denotes the dimension
 of the $i$-th graded piece of the filtration $W^{\bullet}_{\xi}([\phi])$:
 $$
h^i([\xi], [\phi])=dim(W^{i}_{\xi}([\phi])/W^{i+1}_{\xi}([\phi])) 
$$
for $i\in [0,l] $. To reflect the ties with geometry we will label those strata by ${\mathfrak L}_l(h^l,  \lambda)$ - the subset of ${\mathfrak L}_l$
parametrizing the points $([\xi], [\phi])$ where the filtration $W^{\bullet}_{\xi}([\phi])$ has length $l$, the $l$-th piece of the filtration has dimension $h^l$ and $\lambda$ is the partition/Young diagram
attached to the quotient spaces $W_{\xi}/W^{l}_{\xi}([\phi])$, that is, $\lambda$ is as follows
$$
\lambda=(l^{h^{l-1}}, (l-1)^{h^{l-2}-h^{l-1}},\ldots, 2^{h^1-h^2}, 1^{h^0-h^1}).
$$
Thus in the pair $(h^l,  \lambda)$ the first entry takes values in the set 
$$
[g-r]^{\times}=\{1,\ldots,g-r\} 
$$
and $\lambda$ belongs to the set of partitions of $(g-r-h^l)$ whose Young diagram can be inscribed into the $h^0 \times l$-rectangle. 
\begin{example}\label{ex:stratum}
	The stratum ${\mathfrak L}_{g-r-1} ={\mathfrak L}_{g-r-1} (1,(l))$ parametrizes the points $([\xi], [\phi])$ in $\PP({\cal W}_{\Sigma^0_r})$ where the filtration $W^{\bullet}_{\xi}([\phi])$ is a maximal ladder.
\end{example}

Until now we considered the left part of the cohomological incidence correspondence
$$
\xymatrix@C=18pt@R=15pt{
&{\bf P}\ar[rd]^{p_2}\ar[ld]_{p_1}&\\
\PP(H^1(\Theta_{C}))& &\PP(\HKC)
}
$$
The outcome is the extra structure along the fibres of the morphism $p_1$ which resulted in stratifying ${\bf P}$ into the length strata and partition shape strata.
We will now turn to the study the right part of the above diagram.
\section{ The morphism $p_2:{\bf P}\longrightarrow \PP(\HKC)$} 
We begin by determining the fibres of the morphism 
$$
p_2:{\bf P}\longrightarrow \PP(\HKC).
$$
\begin{lem}\label{lem:Zs-Xis}
	Let $\phi$ be a nonzero global section of $\OO_C (K_C)$ and let $Z_{\phi}=(\phi=0)$ be its scheme of zeros. Then the kernel $\Xi_{\phi}$ of the
	map
	$$
	\xymatrix{
		H^1(\OO_C (-K_C))  \ar[r]^(.6)\phi & 	H^1(\OO_C)
	}
	$$
	is canonically identified with the kernel of the trace morphism $H^0(\OO_{Z_{\phi}}) \stackrel{tr}{\longrightarrow} \CC$.
\end{lem}
\begin{pf}
	The morphism
	$$
	\xymatrix{
		\OO_C (-K_C)  \ar[r]^(.6)\phi & 	\OO_C
	}
	$$
	of the multiplication by $\phi$ can be completed to the following exact sequence
	$$	\xymatrix{
		0\ar[r]&	\OO_C (-K_C)  \ar[r]^(.6)\phi & 	\OO_C \ar[r] & \OO_{Z_{\phi}} \ar[r]&0.
	}
	$$
	This gives rise to the sequence of the cohomology groups
	$$
	\xymatrix{
		0 \ar[r]&H^0(\OO_C) \ar[r]& H^0 (\OO_{Z_\phi}) \ar[r]&  	H^1(\OO_C (-K_C))  \ar[r]^(.6)\phi & 	H^1(\OO_C) \ar[r]&0.
	}
	$$
	From this it follows that the kernel $\Xi_{\phi}$ of the arrow induced by  $\phi$ on the cohomology level fits into the exact sequence
	$$
	\xymatrix{
		0 \ar[r]&H^0(\OO_C) \ar[r]& H^0 (\OO_{Z_{\phi}}) \ar[r]&  \Xi_{\phi} \ar[r]&0.
	}
	$$
	This gives the identification 
	$$
	H^0 (\OO_{Z_{\phi}}) /H^0(\OO_C) \cong \Xi_{\phi}.
	$$
	On $ H^0 (\OO_{Z_{\phi}})$ we also have the trace map
	$$
	tr: H^0 (\OO_{Z_{\phi}}) \longrightarrow \CC
	$$
	defined by the composition
	$$
	H^0 (\OO_{Z_{\phi}}) \hookrightarrow End_{\CC}(H^0 (\OO_{Z_{\phi}})) \stackrel{tr}{\longrightarrow} \CC,
	$$
	where the arrow on the right is the trace map of endomorphisms and the inclusion is defined by sending a function $f \in H^0 (\OO_{Z_{\phi}})$ onto the endomorphism
	$$
	m_f: H^0 (\OO_{Z_{\phi}}) \longrightarrow H^0 (\OO_{Z_{\phi}})
	$$
	of the multiplication by $f$:
	$m_f(a)=fa$, for all $a\in  H^0 (\OO_{Z_{\phi}})$. 
	
	The subspace $H^0 (\OO_{C}) \subset H^0 (\OO_{Z_{\phi}})$ of constant functions is mapped onto $\CC id_{H^0 (\OO_{Z_{\phi}})}$. Hence the direct sum decomposition
	$$
	H^0 (\OO_{Z_{\phi}}) =H^0 (\OO_{C}) \oplus ker(tr)
	$$
	and the identification
	$$
	ker(tr)\cong H^0 (\OO_{Z_{\phi}})/H^0 (\OO_{C}) \cong \Xi_{\phi}.
	$$	
\end{pf}

The lemma above implies the following.
\begin{pro}\label{pro:Pbundle}
	The morphism
	$$
	p_2:{\bf P}\longrightarrow \PP(\HKC)
	$$
	is a $\PP^{2(g-2)}$-bundle over $\PP(\HKC)$. In particular, 
	${\bf P}$ is a smooth connected projective variety of dimension $(3g-5)$.	
\end{pro}

 In Lemma \ref{lem:Zs-Xis} we established an isomorphism between the vector space overlying the fibre of $p_2$ over a point $[\phi]\in \PP(\HKC)$ and a certain codimension one subspace of $H^0(\OO_{Z_{\phi}})$; that map will be denoted 
\begin{equation}\label{Fs-map}
	F(\phi,\bullet):\Xi_{\phi} \longrightarrow H^0(\OO_{Z_{\phi}});
\end{equation}
for $\xi \in \Xi_{\phi} $ the value $F(\phi,\xi)$ is the function on the divisor 
$Z_{\phi}$ canonically associated to $\xi$. This map is important in understanding
the two part diagram
\begin{equation}\label{twopart-diag}
\xymatrix@C=18pt{
	&{\bf P}\ar_{p_1}[dl]\ar^{p_2}[dr]&&{\cal Z}_C\ar_{q_1}[dl]\ar^{q_2}[dr]&&\\
	\PP(H^1 (\Theta_{C}))&&\PP(\HKC)&&C \ar@{^{(}->}[r]& \PP(\HKC^{\ast})
}
\end{equation}
relating the cohomological and geometric incidence correspondences: the projectivization of the space $\Xi_{\phi}$ is the fibre of $p_2$ over the point $[\phi]\in \PP(\HKC)$ while the space $H^0(\OO_{Z_{\phi}})$ is the fibre at $[\phi]$ of the direct image sheaf
$q_{1\ast} \OO_{{\cal Z}_C}$. We give the sheaf version of the map \eqref{Fs-map}. This will also
 identify ${\bf P}$ as the projectivization of a certain vector bundle 
on $\PP(\HKC)$. 

To simplify the notation set $W:=\HKC$. On the Cartesian product $\PP(W)\times C$ we have the morphism
$$
q^{\ast}_1(\OO_{\PP(W)}(-1)) \otimes q^{\ast}_2 (\OO_C(-K_C)) \longrightarrow \OO_{\PP(W)\times C}
$$
defining the geometric incidence ${\cal Z}_C$, that is we have the exact sequence
$$
\xymatrix{
	0\ar[r]& q^{\ast}_1(\OO_{\PP(W)}(-1)) \otimes q^{\ast}_2 (\OO_C(-K)) \ar[r]& \OO_{\PP(W)\times C} \ar[r]&\OO_{{\cal Z}_C} \ar[r]&0.
}
$$
Taking its direct image with respect to $q_1$ gives
$$
\xymatrix@C=12pt{
	0\ar[r]&\OO_{\PP(W)} \ar[r]&q_{1\ast}\OO_{{\cal Z}_C} \ar[r]& H^1(\OO(-K_C))\otimes \OO_{\PP(W)}(-1)\ar[r]& H^1(\OO_C) \otimes \OO_{\PP(W)} \ar[r]&0 .
}
$$
This is a sheaf version of the long exact cohomology sequence in the proof
of Lemma \ref{lem:Zs-Xis}. This exact complex is broken into two short exact sequences
$$
\xymatrix@C=12pt{
	0\ar[r]&\OO_{\PP(W)} \ar[r]&q_{1\ast}\OO_{{\cal Z}_C} \ar[r]&{\cal P} \ar[r]&0,
}
$$
$$
\xymatrix@C=12pt{
	0\ar[r]&{\cal P}\ar[r]& H^1(\OO(-K_C))\otimes \OO_{\PP(W)}(-1)\ar[r]& H^1(\OO_C) \otimes \OO_{\PP(W)} \ar[r]&0.
}
$$
In addition, there is the trace map
$$
q_{1\ast}\OO_{{\cal Z}_C} \stackrel{tr}{\longrightarrow} \OO_{\PP(W)}
$$
splitting the first exact sequence above. This gives the identification
\begin{equation}\label{P=kertr}
	{\cal P} \cong ker\left(q_{1\ast}\OO_{{\cal Z}_C} \stackrel{tr}{\longrightarrow} \OO_{\PP(W)}\right),
\end{equation}
which is the sheaf version of the isomorphism \eqref{Fs-map}. Thus we obtain the following.
\begin{pro}\label{pro:Pproj}
	The cohomological incidence correspondence ${\bf P}$ in \eqref{twopart-diag}
	is the projectivization of the vector bundle ${\cal P}$ which admits
	the following resolutions:
	$$
	\xymatrix@C=12pt{
		0\ar[r]&\OO_{\PP(W)} \ar[r]&q_{1\ast}\OO_{{\cal Z}_C} \ar[r]&{\cal P} \ar[r]&0,
	}
	$$
	
	$$
	\xymatrix@C=12pt{
		0\ar[r]&{\cal P}\ar[r]& H^1(\OO(-K_C))\otimes \OO_{\PP(W)}(-1)\ar[r]& H^1(\OO_C) \otimes \OO_{\PP(W)} \ar[r]&0.
	}
	$$
Furthermore, the first, a projective resolution of ${\cal P}$, splits
$$
\OO_{\PP(W)} \oplus 	{\cal P} \cong q_{1\ast}\OO_{{\cal Z}_C}
$$
and identifies ${\cal P}$ as the kernel of the trace map
$$
{\cal P} \cong ker\left(q_{1\ast}\OO_{{\cal Z}_C} \stackrel{tr}{\longrightarrow} \OO_{\PP(W)}\right).
$$
\end{pro}

With the above facts we return to the map
$$
p_1: {\bf P} \longrightarrow \PP(H^1(\Theta_C))
$$
to obtain the following.
\begin{pro}\label{pro:Sigma(g-1)}
	The image of $p_1$ is the stratum $\Sigma_{g-1}$ of the Griffiths stratification \eqref{rkstrat}. It is a determinantal hypersurface of degree $g$ in 
	$\PP(H^1(\Theta_C))$. The map
	$$
	p_1: {\bf P} \longrightarrow \Sigma_{g-1} \subset\PP(H^1(\Theta_C))
	$$
	onto its image is birational, that is, ${\bf P}$ is a resolution of singularities of the hypersurface $\Sigma_{g-1}$.
\end{pro}
 \begin{pf}
 	The image of ${\bf P}$ under $p_1$ is precisely the set of points of $\PP(H^1(\Theta_C))$ where the morphism
 	\begin{equation}\label{p-cp-morphism}
 	\HKC \otimes \OO_{\PP(H^1(\Theta_C))} (-1)\longrightarrow H^1(\OO_C)\otimes \OO_{\PP(H^1(\Theta_C))}
 \end{equation}
 	drops its rank. Thus
 set-theoretically $\Sigma_{g-1}$ is equal to the image of ${\bf P}$ under $p_1$.
 The latter, by Proposition \ref{pro:Pbundle}, is an irreducible subvariety
 of $\PP(H^1(\Theta_C))=\PP^{3g-4}$ of codimension at least one. This implies that the morphism of sheaves \eqref{p-cp-morphism} is injective and
 $\Sigma_{g-1}$ is the scheme of zeros of the determinant of that morphism
 $$
 \mbox{
 	$\bigwedge^g \HKC \otimes\OO_{\PP(H^1(\Theta_C))} (-g)\longrightarrow \bigwedge^g H^1(\OO_C)\otimes \OO_{\PP(H^1(\Theta_C))}$.
 }
$$
Hence $\Sigma_{g-1}$ is the subscheme of zeros of a global section of
$\OO_{\PP(H^1(\Theta_C))}(g)$, that is a hypersurface in $\PP(H^1(\Theta_C))$ of degree $g$. This in turn implies that the image of
${\bf P}$ under $p_1$ is a hypersurface. Hence $p_1$ is generically finite onto its image. Let $[\xi]$ be a closed point in the image of $p_1$. Recall that the fibre of $p_1$ over $[\xi]$ is the projectivization of the vector space
$$
W_{\xi}=ker\left(\HKC \stackrel{\xi}{\longrightarrow}H^1(\OO_C)\right).
$$
This means that $W_{\xi}$ is one dimensional at a general point $[\xi]$ of the image of $p_1$. Equivalently, the morphism \eqref{p-cp-morphism} at a general point of the image of $p_1$ drops its rank precisely by one and this means that the determinantal hypersurface $\Sigma_{g-1}$ is reduced and
hence equals the image $p_1({\bf P})$ scheme theoretically.
 \end{pf}

The picture we get from the consideration of the morphism
$$
p_1: {\bf P}\longrightarrow \PP(H^1(\Theta_C))
$$
is reminiscent of the classical Abel-Jacobi map
$$
AJ:C^{(g-1)}  \longrightarrow J(C)
$$
of $C^{(g-1)}$, the $(g-1)$-symmetric power of $C$, into the Jacobian $J(C)$ of $C$;
 the image $W_{g-1}$, up to translation in $J(C)$, is the theta divisor of
the Jacobian $J(C)$ of $C$, while the special linear systems $g^r_{g-1}$ on $C$ appear as fibres of $AJ$: for a line bundle $\OO_C (D)$ of degree $(g-1)$ the fibre of $AJ$ over $[\OO_C (D)] \in J(C)$ is the linear system $|H^0(\OO_C(D))|$. From this perspective the hypersurface $\Sigma_{g-1}$ in Proposition \ref{pro:Sigma(g-1)} is an analogue of the theta divisor of $J(C)$ and the higher degeneracy loci $\Sigma_r$'s are related to the singularities of the map $p_1$ and the `linear systems' of rank two bundles.  Namely, we have the following.
\begin{pro}\label{pro:diffp1}
	Let $([\xi],[\phi])$ be a point of ${\bf P}$. Then the differential
	$$
	d{p_1}_{([\xi],[\phi])}:T_{{\bf P},([\xi],[\phi])} \longrightarrow T_{\PP(H^1(\Theta_C)),[\xi]}
	$$
	of $p_1$ at $([\xi],[\phi])$ has the kernel isomorphic to
	$W_{\xi}/\CC\phi$, the tangent space at $[\phi]$ of the fibre
	$$
	p^{-1}([\xi])=\PP(W_{\xi})
	$$
	of $p_1$ over $[\xi]$. Furthermore, identify $[\xi]$ with the extension sequence 
	$$
	\xymatrix{
	0\ar[r]&\OO_C \ar[r]& \EE_{\xi}\ar[r]& \OO_C (K_C)\ar[r]&0,
}
	$$
	up to the action of $\CC^{\times}$ on the arrows of the sequence; thus the vector bundle $\EE_{\xi}$ is intrinsically associated to the point $[\xi]$; its space of global sections $H^0(\EE_{\xi})$ is related to $W_{\xi}$ by the exact sequence
	$$
		\xymatrix{
		0\ar[r]&H^0(\OO_C) \ar[r]& H^0(\EE_{\xi})\ar[r]& W_{\xi}\ar[r]&0.
	}
	$$
	In particular, the following equality holds
	$$
dim(ker(d{p_1}_{([\xi],[\phi])}))=	corank(d{p_1}_{([\xi],[\phi])})=h^0(\EE_{\xi})-2.
	$$
	
	The projectivized tangent cone of the determinantal hypersurface
	$\Sigma_{g-1}$ at $[\xi]$ is the image of the morphism
	$$
	\PP({\cal N}_{\PP(W_{\xi})/{\bf P}}) \longrightarrow \PP(T_{\PP(H^1(\Theta_C)),[\xi]}),
	$$
	where ${\cal N}_{\PP(W_{\xi})/{\bf P}}$ is the normal bundle of the fibre
	$p^{-1}_1([\xi])=\PP(W_{\xi})$ in ${\bf P}$; that normal bundle fits into the following exact sequence
	$$
	\xymatrix{
		0\ar[r]& {\cal T}_{p_2} \otimes \OO_{\PP(W_{\xi})} \ar[r] &{\cal N}_{\PP(W_{\xi})/{\bf P}} \ar[r]& W^{\perp}_{\xi} \otimes \OO_{ \PP(W_{\xi})} (1) \ar[r]&0.
	}
$$
\end{pro}
\begin{pf}
	From the projection $p_2:{\bf P}\longrightarrow \PP(\HKC)$ the tangent bundle
	${\cal T}_{\bf P}$ has the form
	$$
	\xymatrix{
	0\ar[r]&{\cal T}_{p_2} \ar[r]&{\cal T}_{\bf P} \ar[r]& p^{\ast}_2\left({\cal T}_{\PP(W)}\right) \ar[r]&0,
}
$$
where ${\cal T}_{p_2}$ is the relative tangent bundle of $p_2$ and $W$ denotes $\HKC$. From Proposition \ref{pro:Pproj} the relative tangent bundle
fits into the relative Euler sequence
$$
	\xymatrix{
	0\ar[r]&\OO_P \ar[r]& {\cal P} \otimes \OO_P (1) \ar[r]&{\cal T}_{p_2} \ar[r]&0.
}
$$
Thus at a point $([\xi],[\phi])$ of ${\bf P}$ the tangent space
$T_{{\bf P},([\xi],[\phi])}$ fits into the the following sequence
$$
\xymatrix{
	0 \ar[r]& \Xi_{\phi}/\CC\xi \ar[r]&T_{{\bf P},([\xi],[\phi])} \ar[r]&W/\CC\phi\ar[r]&0,
}
$$
where $\Xi_{\phi}$ is the fibre of ${\cal P}$ at $[\phi]$ and we use the identification of the fibre of ${\cal T}_{p_2}$ at $([\xi],[\phi])$ with the quotient space $\Xi_{\phi}/\CC\xi$. The differential $dp_1$ at $([\xi],[\phi])$ takes us to the tangent space of $\PP(H^1(\Theta_C))$ at $[\xi]$:
$$
\xymatrix{
	0\ar[d]&\\
	\Xi_{\phi}/\CC\xi \ar[d]&\\
 T_{{\bf P},([\xi],[\phi])} \ar[r]^{dp_{1,([\xi],[\phi])}}\ar[d]& H^1(\Theta_C)/\CC\xi\\
 \HKC/\CC\phi \ar[d]&\\
 0&
}
$$
From the definition of the incidence correspondence ${\bf P}$ the above is compatible with the cup products maps
$$
\xymatrix{
	&H^1(\Theta_{C})\ar[d]^{\phi}\\
	\HKC \ar[r]^(.55){\xi} &H^1(\OO_C)
}
$$
that is, we have the commutative diagram
$$
\xymatrix{
	0\ar[d]&\\
	\Xi_{\phi}/\CC\xi \ar[d]&\\
	T_{{\bf P},([\xi],[\phi])} \ar[r]^{dp_{1,([\xi],[\phi])}}\ar[d]& H^1(\Theta_C)/\CC\xi \ar[d]^{\phi}\\
	\HKC/\CC\phi \ar[d] \ar[r]^(.6){\xi}&H^1(\OO_C)\\
	0&
}
$$ 
This can be completed in the following way
$$
\xymatrix{
&&	0\ar[d]&0\ar[d]\\
&&	\Xi_{\phi}/\CC\xi \ar[d] \ar@{=}[r]&ker(\phi)/\CC\xi \ar[d]\\
0\ar[r]&ker(dp_{1,([\xi],[\phi])})\ar[r]&	T_{{\bf P},([\xi],[\phi])} \ar[r]^{dp_{1,([\xi],[\phi])}}\ar[d]& H^1(\Theta_C)/\CC\xi \ar[d]^{\phi}\\
0\ar[r]&W_{\xi}/\CC\phi\ar[r]&\HKC/\CC\phi \ar[d] \ar[r]^(.6){\xi}&H^1(\OO_C)\ar[d]\\
	&&0&0
}
$$
The commutativity of this diagram implies that the kernel $ker(dp_{1,([\xi],[\phi])})$ factors through $W_{\xi}/\CC\phi$ and the resulting linear map
$$
ker(dp_{1,([\xi],[\phi])}) \longrightarrow W_{\xi}/\CC\phi
$$
is an isomorphism. This proves the first assertion of the proposition. The second follows from the definition of the extension sequence
$$
\xymatrix{
	0\ar[r]&\OO_C \ar[r]& \EE_{\xi}\ar[r]& \OO_C (K_C)\ar[r]&0.
}
$$

We now turn to the statement about the projectivized tangent cone of $\Sigma_{g-1}$.  The differential $dp_1$ restricted to the fibre
$p^{-1}_1([\xi])=\PP(W_{\xi})$ gives rise to the map
$$
\PP({\cal N}_{\PP(W_{\xi})/{\bf P}}) \longrightarrow \PP(T_{\PP(H^1(\Theta_C)),[\xi]} \otimes \OO_{\PP(W_{\xi})}) =\PP(T_{\PP(H^1(\Theta_C)),[\xi]}) \times\PP(W_{\xi}),
	$$
where ${\cal N}_{\PP(W_{\xi})/{\bf P}}$ is the normal bundle of the fibre $p^{-1}_1([\xi])=\PP(W_{\xi})$ in ${\bf P}$. Composing with the projection onto the first factor gives the morphism
	$$
	\PP({\cal N}_{\PP(W_{\xi})/{\bf P}}) \longrightarrow \PP(T_{\PP(H^1(\Theta_C)),[\xi]})
	$$
	and it is well known that its image is the projectivized tangent cone
	of $\Sigma_{g-1}$ at the point $[\xi]$, see \cite{ACGH}. Thus only the description of the normal bundle ${\cal N}_{\PP(W_{\xi})/{\bf P}}$ requires a proof. For this write the normal sequence of $\PP(W_{\xi}) \subset {\bf P}$
	$$
	\xymatrix{
		0\ar[r]&{\cal T}_{\PP(W_{\xi})} \ar[r]& {\cal T}_{\bf P} \otimes \OO_{ \PP(W_{\xi})}\ar[r]& {\cal N}_{\PP(W_{\xi})/{ \bf P}}\ar[r]&0.
	}
	$$
	Combining with the exact sequence for the tangent bundle ${\cal T}_{\bf P}$ of ${\bf P}$ we obtain
	$$
	\xymatrix{
		&&0\ar[d]&&\\
		&&{\cal T}_{p_2}\otimes \OO_{ \PP(W_{\xi})}\ar[d]&&\\
		0\ar[r]&{\cal T}_{\PP(W_{\xi})} \ar[r]& {\cal T}_{\bf P} \otimes \OO_{ \PP(W_{\xi})}\ar[r]\ar[d]& {\cal N}_{\PP(W_{\xi})/{\bf P}}\ar[r]&0\\
		&&p^{\ast}_2 {\cal T}_{\PP(W)}\otimes \OO_{ \PP(W_{\xi})}\ar[d]&&\\
		&&0&&
	}
	$$
	where $W=\HKC$. From the point of view of the projection $p_2$ the fibre
	$p^{-1}_1 ([\xi])=\PP(W_{\xi})$ is a section of $p_2$ over $\PP(W_{\xi}) \subset \PP(W)$. This allows to complete the above diagram by the normal sequence of $\PP(W_{\xi})$ in $\PP(W)$:
		$$
	\xymatrix{
		&&0\ar[d]&&\\
		&&{\cal T}_{p_2}\otimes \OO_{ \PP(W_{\xi})}\ar[d]&&\\
		0\ar[r]&{\cal T}_{\PP(W_{\xi})} \ar@{=}[d] \ar[r]& {\cal T}_{\bf P} \otimes \OO_{ \PP(W_{\xi})}\ar[r]\ar[d]& {\cal N}_{\PP(W_{\xi})/{\bf P}}\ar[r]&0\\
		0\ar[r]&{\cal T}_{\PP(W_{\xi})} \ar[r]& {\cal T}_{\PP(W)}\otimes \OO_{ \PP(W_{\xi})}\ar[d] \ar[r]&{\cal N}_{\PP(W_{\xi})/\PP(W)}\ar[r]&0\\
		&&0&&
	}
	$$
	From this follows the exact sequence
	$$
	\xymatrix{
		0\ar[r]&{\cal T}_{p_2}\otimes \OO_{ \PP(W_{\xi})} \ar[r]& {\cal N}_{\PP(W_{\xi})/{\bf P}}\ar[r]& {\cal N}_{\PP(W_{\xi})/\PP(W)}\ar[r]&0.
	}
$$
Since $\PP(W_{\xi})$ in $\PP(W)$ is a complete intersection of hyperplanes parametrized by $W^{\perp}_{\xi}$ we deduce the formula for the normal bundle
${\cal N}_{\PP(W_{\xi})/\PP(W)}$:
$$
{\cal N}_{\PP(W_{\xi})/\PP(W)}=W^{\perp}_{\xi} \otimes\OO_{ \PP(W_{\xi})} (1).
$$
Hence the exact sequence
$$
\xymatrix{
	0\ar[r]&{\cal T}_{p_2}\otimes \OO_{ \PP(W_{\xi})} \ar[r]& {\cal N}_{\PP(W_{\xi})/{\bf P}}\ar[r]&W^{\perp}_{\xi} \otimes\OO_{ \PP(W_{\xi})} (1) \ar[r]&0
}
$$
asserted in the statement of the proposition. 
\end{pf}

The above results illustrate the idea of Griffiths that the IVHS formalism, the map $p_1$ in our situation, should be a substitute for the classical theta divisor and the higher degeneracy strata $\Sigma_r$'s should be related to geometry of $C$.

Let us now reinterpret a part of our refinement of IVHS using the above considerations. We assume that $\xi$ is a nonzero Kodaira-Spencer class of rank $r \leq g-2$. Then the space $W_{\xi}$ has dimension 
$g-r \geq 2$. Fix a nonzero $\phi \in W_{\xi}$. Using the inclusion
\begin{equation}\label{Fphixi-endo}
F(\phi,\bullet): \Xi_{\phi} \hookrightarrow H^0(\OO_{Z_{\phi}})
\end{equation}
	in \eqref{Fs-map} we identify $\xi$ with the function $F(\phi,\xi)$ on $Z_{\phi}$, the zero divisor of $\phi$. The multiplication by $F(\phi,\xi)$ gives the endomorphism
	$$
	m_{F(\phi,\xi)}: H^0(\OO_{Z_{\phi}}(K_C)) \longrightarrow  H^0(\OO_{Z_{\phi}}(K_C)),
	$$
	that is, $ m_{F(\phi,\xi)}(s)=F(\phi,\xi) s$, $\forall s\in H^0(\OO_{Z_{\phi}}(K_C))$.
The filtration $W^{\bullet}_{\xi}([\phi])$ can be realized via this endomorphism. For this we tensor the exact sequence
$$
\xymatrix{
0\ar[r]&\OO_C(-K_C)\ar[r]^(.6){\phi}& \OO_C \ar[r] &\OO_{Z_{\phi}}\ar[r]&0
}
$$
with $\OO_C(K_C)$ to obtain
$$
\xymatrix{
	0\ar[r]&\OO_C\ar[r]^(.35){\phi}& \OO_C (K_C) \ar[r] &\OO_{Z_{\phi}}(K_C)\ar[r]&0.
}
$$
Writing out the sequence of the cohomology groups gives the exact sequence
$$
\xymatrix{
	0\ar[r]&\HKC/\CC\phi \ar[r] &H^0(\OO_{Z_{\phi}}(K_C))\ar[r]&H^1(\OO_C) \ar[r]&H^1(\OO_C (K_C) ) \ar[r]&0.
}
$$
The essential point of the construction is to study the action of the endomorphism $m_{F(\phi,\xi)}$ on the inclusion
$$
\xymatrix{
	0\ar[r]&\HKC/\CC\phi \ar[r] &H^0(\OO_{Z_{\phi}}(K_C)).
}
$$
Namely, we want to understand which vectors of  $\HKC/\CC\phi$ remain in that subspace under the action of  $m_{F(\phi,\xi)}$. The answer is provided in the following.
\begin{lem}\label{lem:Fphixi-Wxi}
	Let $\psi\in \HKC$ and denote by $\{\psi\}$ its equivalence class in 
	the quotient $\HKC/\CC\phi$. Then the image $m_{F(\phi,\xi)}(\{\psi\})$ lies in
		$\HKC/\CC\phi$ if and only if $\psi \in W_{\xi}$.
\end{lem} 

\begin{pf}
	The multiplication by $\psi$ connects two exact sequences
	$$
	\xymatrix@C=12pt{
		0\ar[r]&\OO_C(-K_C)\ar[r]^(.6){\phi}& \OO_C \ar[r] &\OO_{Z_{\phi}}\ar[r]&0\\
		0\ar[r]&\OO_C\ar[r]^(.35){\phi}& \OO_C (K_C) \ar[r] &\OO_{Z_{\phi}}(K_C)\ar[r]&0
	}
	$$
	by a commutative diagram
	$$
	\xymatrix@C=12pt{
		0\ar[r]&\OO_C(-K_C)\ar[r]^(.6){\phi} \ar[d]^{\psi}& \OO_C \ar[r]\ar[d]^{\psi} &\OO_{Z_{\phi}}\ar[r] \ar[d]^{\psi}&0\\
		0\ar[r]&\OO_C\ar[r]^(.35){\phi}& \OO_C (K_C) \ar[r] &\OO_{Z_{\phi}}(K_C)\ar[r]&0
	}
	$$
	Passing to cohomology we obtain
	$$
	\xymatrix@C=12pt{
		&0\ar[r] \ar[d]^{\psi}& H^0(\OO_C) \ar[r]\ar[d]^{\psi} &H^0(\OO_{Z_{\phi}})\ar[r]^(.4){\delta} \ar[d]^{\psi}&H^1(\OO_C(-K_C)) \ar[d]^{\psi}\\
		0\ar[r]&H^0(\OO_C)\ar[r]^(.4){\phi}& H^0(\OO_C (K_C) )\ar[r] &H^0(\OO_{Z_{\phi}}(K_C))\ar[r]^(.6){\delta'}&H^1(\OO_C)
	}
	$$
	We examine the rightmost commutative square of the above diagram:
	by definition the function $F(\phi,\xi)$ under the top horizontal arrow goes to $\xi$; hence the formula
	$$
	\psi \xi =\delta'(F(\phi,\xi)\psi).
	$$
	Observe: $F(\phi,\xi)\psi$ lies in the subspace $\HKC/ \CC\phi$ if and only if the right hand side vanishes or, equivalently, $\psi \xi =0$ and this means that $\psi \in W_{\xi}$.	
\end{pf}

The above lemma produces the commutative diagram
\begin{equation}\label{Wxi-Fphixi}
	\xymatrix{
0\ar[r]&	W_{\xi}/ \CC\phi \ar[r]\ar[d]& \HKC/ \CC\phi \ar[d]^{m_{F(\phi,\xi)}}\\
0\ar[r]&	\HKC/ \CC\phi \ar[r]&H^0(\OO_{Z_{\phi}}(K_C))
}
\end{equation}
The vertical arrow on the left will be denoted 
\begin{equation}\label{beta-phi-map}
\beta(\phi,\bullet): W_{\xi}/ \CC\phi \longrightarrow \HKC/ \CC\phi.
\end{equation}
The above diagram tells us that the restriction to $Z_{\phi}$ is subject to
the identity
\begin{equation}\label{beta-formula}
	\beta(\phi,\psi)\big|_{Z_{\phi}} =F(\phi,\xi) \psi \big|_{Z_{\phi}}, \,\,\forall \psi \in W_{\xi}.
\end{equation}
The map can be lifted to a homomorphism
$$
\widetilde{\beta}(\phi, \bullet): W_{\xi} \longrightarrow \HKC
$$
subject to

- $\widetilde{\beta}(\phi, \phi)=0$,

- $\widetilde{\beta}(\phi, \bullet)\big|_{Z_{\phi}} =\beta(\phi,\bullet)\big|_{Z_{\phi}}$;

\noindent
furthermore, two different liftings $\widetilde{\beta}(\phi, \bullet)$ and $\widetilde{\beta'}(\phi, \bullet)$ differ as follows
$$
\widetilde{\beta'}(\phi, \bullet)-\widetilde{\beta}(\phi, \bullet)=f(\bullet) \phi,
$$
for some linear function $f$ on $W_{\xi}$. In addition, we have 
$$
\psi' \widetilde{\beta}(\phi,\psi) - \psi \widetilde{\beta}(\phi,\psi') \in \phi \HKC.
$$
Hence we define $\widetilde{\beta}(\phi',\psi) \in \HKC$ so that the equality
\begin{equation}\label{beta-def}
\psi' \widetilde{\beta}(\phi,\psi) - \psi \widetilde{\beta}(\phi,\psi')=\phi \widetilde{\beta}(\phi',\psi), \forall \phi, \psi' \in W_{\xi}.
\end{equation}
This gives a bilinear skew-symmetric map
$$
\widetilde{\beta}(\bullet,\bullet): W_{\xi}\times W_{\xi} \longrightarrow \HKC
$$
or, equivalently, the linear map
$$
\mbox{$\widetilde{\beta}(\bullet,\bullet): \bigwedge^2 W_{\xi} \longrightarrow \HKC$}
$$
Furthermore, the identity \eqref{beta-def} becomes the Koszul cycle relation
$$
\psi' \widetilde{\beta}(\phi,\psi) - \psi \widetilde{\beta}(\phi,\psi') +\phi \widetilde{\beta}(\psi,\psi')=0, \,\,\forall \phi,\psi,\psi'\in W_{\xi}.
$$
The map $\widetilde{\beta}(\bullet,\bullet)$ is formally similar to $\alpha^{(2)}_{\xi}$ we attached to $\xi$ via the extension construction, see Remark \ref{rem:alpha}. So we can use it in the same way as in the proof of Lemma \ref{lem:phi-filt} to produce a filtration on $W_{\xi}$. The two constructions become identical under the assumption $\xi$ is base point free; we remind the reader that this means that the linear subsystem $|W_{\xi}|$ is base point free, see Definition \ref{def:gpp}.
\begin{pro}\label{pro:beta=alpha}
	Assume a nonzero Kodaira-Spencer class $\xi$ is base point free. Then we have equality
	$$
	\widetilde{\beta}(\bullet,\bullet) \equiv\alpha^{(2)}_{\xi}  mod(\text{Koszul coboundary}).
	$$
\end{pro}
\begin{pf}
	Consider the Koszul complex associated to $(\OO_C (-K_C), W_{\xi})$:
	$$
	\xymatrix@C=12pt{
		0\ar[r]& \OO_C(-K_C) \ar[r]&  W^{\ast}_{\xi} \otimes \OO_C \ar[r] & \bigwedge^2 W^{\ast}_{\xi} \otimes \OO_C (K_C) \ar[r]& \bigwedge^3 W^{\ast}_{\xi} \otimes \OO_C (2K_C) \ar[r]&\cdots
	}
$$
Since the linear subsystem $|W_{\xi}|$ is base point free the complex is exact and the cohomology of the complex
$$
\xymatrix{
0 \ar[r]&  W^{\ast}_{\xi} \ar[r] & \bigwedge^2 W^{\ast}_{\xi} \otimes H^0(\OO_C (K_C)) \ar[r]& \bigwedge^3 W^{\ast}_{\xi} \otimes H^0(\OO_C (2K_C)) \ar[r]&\cdots
}
$$
can be identified with the cohomology of the complex
$$
\xymatrix{
	H^1 (\OO_C(-K_C)) \ar[r]&  W^{\ast}_{\xi} \otimes H^1(\OO_C )\ar[r] & \bigwedge^2 W^{\ast}_{\xi} \otimes H^1(\OO_C (K_C)) \ar[r]& 0\ar[r]&\cdots
}
$$
via the differentials of the spectral sequence. In particular, the group of 
Koszul cocycles in $\bigwedge^2 W^{\ast}_{\xi} \otimes H^0(\OO_C (K_C))$
modulo the Koszul coboundaries, the group we are interested in, is identified as the term $E^{2,0}_2=E^{2,0}_{\infty}$ of the spectral sequence and the 
 differential
$$
d_2: E^{0,1}_{\infty}=E^{0,1}_2 \longrightarrow E^{2,0}_2=E^{2,0}_{\infty}
$$
is an isomorphism. The term $E^{0,1}_2$ is the kernel of the map
$$
H^1(\OO_C(-K_C)) \longrightarrow W^{\ast}_{\xi} \otimes H^1(\OO_C ).
$$
Hence the identification
$$
\begin{gathered}
E^{0,1}_{\infty}=ker\left(H^1(\OO_C(-K_C))\right) \longrightarrow W^{\ast}_{\xi} \otimes H^1(\OO_C )) \cong E^{2,0}_{\infty} 
\\
= \frac{ker\left(\bigwedge^2 W^{\ast}_{\xi} \otimes H^0(\OO_C (K_C)) \longrightarrow \bigwedge^3 W^{\ast}_{\xi} \otimes H^0(\OO_C (2K_C))\right)}{im\left(W^{\ast}_{\xi} \longrightarrow\bigwedge^2 W^{\ast}_{\xi} \otimes H^0(\OO_C (K_C))\right)}
\end{gathered}
$$
Our Kodaira-Spencer class $\xi$ is an element of $E^{0,1}_{\infty}$ and we attached to it two Koszul cocycles $\widetilde{\beta}$ and $\alpha^{(2)}_{\xi}$. Since $d_2$ is an isomorphism, the two cocycle must represent the same cohomology class in $E^{2,0}_{\infty}$ and hence they differ by a Koszul coboundary.
\end{pf}

Let $[\xi]\in \Sigma^{\circ}_r$ for $r\in[2,g-2]$ and assume it is base point free. We have seen that the last piece $W^{l_{\xi}([\phi])}_{\xi}([\phi])$ of the filtration $W^{\bullet}_{\xi}([\phi])$, if its dimension is at least two, is related to special divisors on $C$. Recall that the link is furnished by the fact that $W^{l_{\xi}([\phi])}_{\xi}([\phi])$ is invariant with respect to the map
$\alpha^{(2)}_{\xi}(\phi,\bullet)$ and the study of the eigen spaces of the
endomorphism
$$
\widehat{\alpha^{(2)}_{\xi}}(\phi,\bullet):
W^{l_{\xi}([\phi])}_{\xi}([\phi])/\CC\phi \longrightarrow W^{l_{\xi}([\phi])}_{\xi}([\phi])/\CC\phi,
$$
see Proposition \ref{pro:xi-Wlatleast2-eigen}. By Proposition \ref{pro:beta=alpha}, up to adding a multiple of identity, this is given by the endomorphism
$$
{\beta}(\phi, \bullet):W^{l_{\xi}([\phi])}_{\xi}([\phi])/\CC\phi \longrightarrow W^{l_{\xi}([\phi])}_{\xi}([\phi])/\CC\phi.
$$
This has the property
\begin{equation}\label{beta=multF}
{\beta}(\phi, \psi)\big|_{Z_{\phi}}=F(\phi,\xi)\psi\big|_{Z_{\phi}}, \,\forall \psi\in W^{l_{\xi}([\phi])}_{\xi}([\phi]) \,\, \text{and all $\phi \neq 0$ in $W_{\xi}$}.
\end{equation}
This equation gives a simple way to connect the spectral properties of
$\widetilde{\beta}(\phi, \bullet)$ (resp. $\alpha^{(2)}_{\xi}(\phi,\bullet)$) and geometry of $Z_{\phi}$.
\begin{pro}\label{pro:beta-eigen}
	Let $[\xi]\in \Sigma^{\circ}_r$ for $r\in[2,g-2]$ and assume it is base point free. For every nonzero $\phi \in W_{\xi}$ the endomorphism
	$$
	{\beta}(\phi, \bullet):W^{l_{\xi}([\phi])}_{\xi}([\phi])/\CC \phi \longrightarrow W^{l_{\xi}([\phi])}_{\xi}([\phi])/ \CC \phi.
	$$
in \eqref{beta-phi-map}	is subject to
	$$
	{\beta}(\phi, \{\psi\})\big|_{Z_{\phi}}=F(\phi,\xi)\psi\big|_{Z_{\phi}}, \,\forall \psi\in W^{l_{\xi}([\phi])}_{\xi}([\phi]) ,
	$$
where $\{\psi\}$ is the equivalence class in $W^{l_{\xi}([\phi])}_{\xi}([\phi])/ \CC \phi$ corresponding to $\psi \in W^{l_{\xi}([\phi])}_{\xi}([\phi])$. 

Let $\lambda$ be an eigen value of ${\beta}(\phi, \bullet)$ and 
$V_{\lambda}$ the corresponding eigen subspace. Denote by $\widetilde{V_{\lambda}}$ the preimage of $V_{\lambda}$ under the projection
$$
W_{\xi}\longrightarrow W_{\xi}/\CC \phi.
$$
  Then the zero divisor $Z_{\phi}$ admits the following decomposition
$$
Z_{\phi}=Z_{\phi}(\lambda)+ Z^{V_{\lambda}}_{\phi}
$$
into proper subschemes, where $Z_{\phi}(\lambda)$ is the subscheme of zeros
of the function $(F(\phi,\xi)-\lambda)$ on $Z_{\phi}$ and $Z^{V_{\lambda}}_{\phi}$ is the subscheme of $Z_{\phi}$ complementary to 
$Z_{\phi}(\lambda)$. Furthermore, the space $H^0(\OO_C (K_C -Z^{V_{\lambda}}_{\phi}))=H^0(\OO_C (Z_{\phi}(\lambda)))$ is isomorphic to $\widetilde{V_{\lambda}}$. 

In addition, the line bundles $\OO_C (K_C -Z^{V_{\lambda}}_{\phi})$ and $\OO_C (Z^{V_{\lambda}}_{\phi})$ are the same as in Proposition \ref{pro:xi-Wlatleast2-eigen}: they give rise to the following commutative diagram
$$
\xymatrix{
	&&0\ar[d]&&\\
	&& \OO_C (K_C -Z^{V_{\lambda}}_{\phi})\ar[d] \ar[dr]&&\\
	0\ar[r]& \OO_{C} \ar[r]^{e_{\xi}} \ar[dr]&\EE_{\xi}\ar[r]^(.42){ e_{\xi}\wedge} \ar[d] & \OO_{C}(K_C)\ar[r]&0\\
	&&\OO_{C}(Z^{V_{\lambda}}_{\phi})  \ar[d]&&\\
	&&0&& 
}
$$
where the slanted arrows are the multiplication by a global section of
$\OO_{C}(Z^{V_{\lambda}}_{\phi})$ corresponding to the divisor $Z^{V_{\lambda}}_{\phi}$.
 If, furthermore, the divisor $Z_{\phi}=(\phi=0)$ is reduced, then the endomorphism
	$$
{\beta}(\phi, \bullet):W^{l_{\xi}([\phi])}_{\xi}([\phi])/\CC \phi \longrightarrow W^{l_{\xi}([\phi])}_{\xi}([\phi])/ \CC \phi.
$$
is semisimple and we have the eigen space decomposition
$$
 W^{l_{\xi}([\phi])}_{\xi}([\phi])/\CC \phi =\bigoplus_{\lambda} V_{\lambda},
 $$
 where the sum is over the distinct eigenvalues of ${\beta}(\phi, \bullet)$ on $W^{l_{\xi}([\phi])}_{\xi}([\phi])/\CC \phi$. 
\end{pro} 
\begin{pf}
	We assume that $W^{l_{\xi}([\phi])}_{\xi}([\phi])$ is at least two dimensional since otherwise ${\beta}(\phi, \bullet)$ is zero and there nothing to prove. Let $\lambda$ be an eigen value of ${\beta}(\phi, \bullet)$. Denote $V_{\lambda}$ the $\lambda$-eigen space of ${\beta}(\phi, \bullet)$ and $\widetilde{V}_{\lambda}$ its preimage in
	$W_{\xi}$. Then we have
	$$
	{\beta}(\phi, \{\psi\})=\lambda \{\psi\}, \forall \psi \in \widetilde{V}_{\lambda};
	$$
	as before $\{s\} $ stands for the class in $W_{\xi}/\CC \phi$ for an $s$ in $W_{\xi}$. This and the equation
	$$
	{\beta}(\phi, \{\psi\})\big|_{Z_{\phi}} =F(\phi,\xi) \psi \big|_{Z_{\phi}}
	$$
	imply the equality
	$$
	\lambda  \psi \big|_{Z_{\phi}}=F(\phi,\xi) \psi \big|_{Z_{\phi}}.
	$$
	Equivalently, this can be rewritten as the equation
	\begin{equation}\label{F-lambda=0-pf}
	(F(\phi,\xi)-\lambda)\psi \big|_{Z_{\phi}}=0,\,\forall \psi \in \widetilde{V_{\lambda}}.
\end{equation}
	Set $Z_{\phi}(\lambda)$ the subscheme of $Z_{\phi}$ defined by the vanishing 
	$$
	F(\phi,\xi)-\lambda=0.
	$$
	Since $F(\phi,\xi)$ is not a constant function, the subscheme is a proper subscheme of $Z_{\phi}$. Let $Z^{V_{\lambda}}_{\phi}$ denote the
	subscheme of $Z_{\phi}$ complementary to $Z_{\phi}(\lambda)$. Then the equation \eqref{F-lambda=0-pf} implies that all $\psi$ in $\widetilde{V}_{\lambda}$ vanish on $Z^{V_{\lambda}}_{\phi}$. Conversely, if $\psi \in \HKC$ vanishes on $Z^{V_{\lambda}}_{\phi}$, then the equation
	$$
	(F(\phi,\xi)-\lambda)\psi \big|_{Z_{\phi}}=0
	$$
	holds. Hence
	$$
	F(\phi,\xi)\psi \big|_{Z_{\phi}}=\lambda \psi \big|_{Z_{\phi}}
	$$
	and this in turn means that $\{\psi\}$ lies in the eigen space
	$V_{\lambda}$ of ${\beta}(\phi,\bullet)$. Hence the identification
	$$
	H^0(\OO_C (K_C-Z^{V_{\lambda}}_{\phi}))=\widetilde{V}_{\lambda}.
	$$
	
	From the equivalence of the Koszul cocycles $\widetilde{\beta}$ and $\alpha^{(2)}_{\xi}$, see Proposition \ref{pro:beta=alpha}, it follows that $\widehat{\alpha^{(2)}_{\xi}}(\phi,\bullet)$ has the same eigen spaces as 
	${\beta}(\phi, \bullet)$. Hence all the results of Proposition \ref{pro:xi-Wlatleast2-eigen} apply.
	
	\vspace{0.2cm}
	We turn to the last assertion of the proposition stating that ${\beta}(\phi,\bullet)$ is semisimple provided the zero divisor $Z_{\phi}$ of $\phi$ is reduced.
	For this we recall that the homomorphism ${\beta}(\phi,\bullet)$ is obtained from the endomorphism
	$$
	m_{F(\phi,\xi)}: H^0(\OO_{Z_{\phi}} (K_C)) \longrightarrow H^0(\OO_{Z_{\phi}}(K_C))
	$$
	of the multiplication by $F(\phi,\xi)$. For $Z_{\phi}$ reduced, that is a collection of $2(g-1)=deg(Z_{\phi})$ distinct points we have the basis
	formed by sections $\psi_p$ of $\OO_{Z_{\phi}} (K_C)$ supported precisely at $p$, as $p$ runs through $Z_{\phi}$. From the equation
	$$
	m_{F(\phi,\xi)}(\psi_p)=F(\phi,\xi)\psi_p=F(\phi,\xi)(p)\psi_p 
	$$
	it follows that each $\psi_p$ is an eigen vector of $m_{F(\phi,\xi)}$ with the eigen value $F(\phi,\xi)(p)$. Hence $m_{F(\phi,\xi)}$ is semisimple on $H^0(\OO_{Z_{\phi}} (K_C))$ with the eigen values
	the set of values of $F(\phi,\xi)$; for a value $\lambda$ of $F(\phi,\xi)(p)$ the corresponding eigen space
	$$
	W_{\lambda}=\bigoplus_{p \in Z_{\phi}(\lambda)} \CC\psi_p
	$$
	where $Z_{\phi}(\lambda)=\{p\in Z_{\phi}| F(\phi,\xi)(p)=\lambda\}$; hence and eigen space decomposition
	$$
	H^0(\OO_{Z_{\phi}} (K_C))=\bigoplus_{\lambda} W_{\lambda},
	$$
	the direct sum is taken over the set of distinct values of $F(\phi,\xi)$. We show that for any $W$, an $m_{F(\phi,\xi)}$-invariant subspace of $H^0(\OO_{Z_{\phi}} (K_C))$ the induced endomorphism
	$$ 
	m_{F(\phi,\xi)}|_W :W \longrightarrow W
	$$
	remains semisimple. For this it will be enough to show that for any eigen value $\lambda$ of $m_{F(\phi,\xi)}|_W$ the equality
	$$
	ker\left((m_{F(\phi,\xi)}|_W -\lambda id_W)^n\right)=ker(m_{F(\phi,\xi)}|_W -\lambda id_W)
	$$
	holds for any $n\geq 2$. Indeed, for $x$ in the kernel of 
	$(m_{F(\phi,\xi)}|_W -\lambda id_W)^n$ for some $n\geq 2$ we have
	$$
	\text{$\left(F(\phi,\xi)-\lambda\right)^n x=0$ in  $H^0(\OO_{Z_{\phi}} (K_C))$},
	$$
	that is, the equality
	$$
	\left(F(\phi,\xi)(p)-\lambda\right)^n x_p=0, \,\forall p\in Z_{\phi},
	$$
	holds; here $x_p$'s are the coordinates of $x$ with respect to the basis $\{\psi_p\}$. If $x\neq0$, then for all $p \in Z_{\phi}$ with $x_p \neq 0$, the above equality implies
	$$
	\left(F(\phi,\xi)(p)-\lambda\right)^n =0 \Leftrightarrow F(\phi,\xi)(p)-\lambda=0.
	$$
	Hence 
	$$
	(F(\phi,\xi)(p)-\lambda)x_p=0, \,\forall p\in Z_{\phi},
	$$
	and this means that $x$ is in the kernel of $(m_{F(\phi,\xi)}|_W -\lambda id_W)$.
\end{pf}

The above result combines together what have been learned about the last step of
$([\xi],[\phi])$-filtrations defined along the fibres of $p_1: {\bf P} \longrightarrow \PP(H^1(\Theta_{C}))$, together with the properties of the second map, $p_2: {\bf P} \longrightarrow \PP(\HKC)$, of the cohomological incidence correspondence.
The main technical feature added is the assignment 
\begin{equation}\label{Fphixi-assignment}
(\xi,\phi) \rightarrow F(\phi,\xi)\in H^0(\OO_{Z_{\phi}})
\end{equation}
	which attaches to a vector $(\xi,\phi)$ overlying a point $([\xi],[\phi])$ in ${\bf P}$ a specific function $F(\phi,\xi)$ on the subscheme
	of zeros $Z_{\phi} =(\phi=0)$. This assignment has also a conceptual meaning of relating the cohomological and geometrical correspondences.
	We wish to discuss this point since it seems to be an important part of the story.
	
	Recall the two part diagram
$$ 
 \xymatrix@C=18pt{
 	&{\bf P}\ar_{p_1}[dl]\ar^{p_2}[dr]&&{\cal Z}_C\ar_{q_1}[dl]\ar^{q_2}[dr]&&\\
 	\PP(H^1 (\Theta_{C}))&&\PP(\HKC)&&C \ar@{^{(}->}[r]& \PP(\HKC^{\ast})
 }
$$
which we have encountered on several occasions. The whole philosophy of the
IVHS can be summarized as gaining insight on the geometrical side, the right part of the diagram, from some functorial data on the cohomological side, the left part of the diagram. We started with a point $[\xi]$ in the rank $r$ stratum
$\Sigma^{\circ}_r$ on the cohomological side. The fibre $\PP(W_{\xi})$ of $p_1$ over $[\xi]$ is a basic IVHS invariant of $[\xi]$.  The projection $p_2$ embeds $\PP(W_{\xi})$ into $\PP(\HKC)$. The basic objects on that projective space which come from the geometric correspondence ${\cal Z}_C$ are the direct image sheaves 
$$
\text{$q_{1 \ast} \OO_{{\cal Z}_C}=q_{1 \ast} \left(q^{\ast}_2\OO_{C} \right)$ and $q_{1 \ast} \left(\OO_{{\cal Z}_C} \otimes q^{\ast}_2\OO_{C} (K_C)\right)= q_{1 \ast} \left(q^{\ast}_2\OO_{C} (K_C) \right)$}.
$$
The assignment \eqref{Fphixi-assignment} can be rephrased as a morphism
\begin{equation}\label{xi-sheafdata}
F(\bullet,\xi) : \OO_{\PP(W_{\xi})} (-1) \longrightarrow q_{1 \ast} \OO_{{\cal Z}_C} \otimes \OO_{\PP(W_{\xi})}.
\end{equation}
which at a point $[\phi] \in \PP(W_{\xi})$ is the map
$$
F(\bullet,\xi)_{[\phi]} : \OO_{\PP(W_{\xi}),[\phi]} (-1)=\CC\phi \longrightarrow \left(q_{1 \ast} \OO_{{\cal Z}_C} \otimes \OO_{\PP(W_{\xi})}\right)_{[\phi]} =H^0(\OO_{Z_{\phi}})
$$
defined by the rule
$$
F(\bullet,\xi)_{[\phi]} (u\phi)=uF(\phi, \xi), \,\forall u\in \CC.
$$
Furthermore, we recast the function $F(\phi, \xi)$ as the endomorphism
$$
m_{F(\phi, \xi)}: H^0(\OO_{Z_{\phi}}(K_C)) \longrightarrow  H^0(\OO_{Z_{\phi}}(K_C))
$$
of multiplication by $F(\phi, \xi)$. This can now be promoted to the sheaf morphism
\begin{equation}\label{mFxi-morphism}
m_{F(\bullet, \xi)}:\OO_{\PP(W_{\xi})} (-1) \longrightarrow {\cal{E}}nd(q_{1 \ast} \left(q^{\ast}_2\OO_{C} (K_C) \right))
\end{equation}
which at a point $[\phi] \in \PP(W_{\xi})$ is the map
$$
(m_{F(\bullet, \xi)})_{[\phi]} : \OO_{\PP(W_{\xi}),[\phi]} (-1)=\CC\phi \longrightarrow \left({\cal{E}}nd(q_{1 \ast} \left(q^{\ast}_2\OO_{C} (K_C) \right))\right)_{[\phi]} =End(H^0(\OO_{Z_{\phi}}(K_C)).
$$
defined by the rule
$$
(m_{F(\bullet, \xi)})_{[\phi]} (u\phi)=um_{F(\phi, \xi)}, \,\forall u\in \CC.
$$
Thus the IVHS data $(\xi, W_{\xi})$ attached to a point $[\xi]\in \Sigma^{\circ}_r$ translates into the data of the category of coherent sheaves on
$\PP(W_{\xi})$:
\begin{equation}\label{xi-catdata}
	\xymatrix{
		\Sigma^{\circ}_r \ni [\xi] \ar[rr]&&
		*+[F]{\txt{Categorical data:\\
				$\OO_{\PP(W_{\xi})} (-1)
				\stackrel{F(\bullet,\xi) }{\longrightarrow} q_{1 \ast} \OO_{{\cal Z}_C} \otimes \OO_{\PP(W_{\xi})}$
			\\
			$\OO_{\PP(W_{\xi})} (-1)
			\stackrel{m_{F(\bullet,\xi)} }{\longrightarrow} {\cal{E}}nd(q_{1 \ast} \left(q^{\ast}_2\OO_{C} (K_C) \right))$}}			
	}
\end{equation}
In addition, if $[\xi]$ is base point free, Proposition \ref{pro:beta=alpha} tells us that at every point $[\phi] \in \PP(W_{\xi})$ the $([\xi],[\phi])$- filtration $W^{\bullet}_{\xi}([\phi])$ can be recovered from the action of the morphism $m_{F(\phi,\xi)}$ on the 
	subspace $\HKC/\CC\phi$ of $H^0(\OO_{Z_{\phi}}(K_C))$. The link to the geometric correspondence comes from the invariant part of this action, the last piece
	$W^{l_{\xi}([\phi])}_{\xi}([\phi])$ of $([\xi],[\phi])$-filtration of $W_{\xi}$: Proposition \ref{pro:beta-eigen} delivers a precise dictionary between the spectral properties of the restriction of $m_{F(\phi,\xi)}$ to $W^{l_{\xi}([\phi])}_{\xi}([\phi])$ and geometry of the finite morphism
	$$
	q_1: {\cal Z}_C \longrightarrow \PP(\HKC) 
	$$ 
	over $\PP(W_{\xi})$; this tells that for $W^{l_{\xi}([\phi])}_{\xi}([\phi])$ at least two dimensional, every eigen value of $m_{F(\phi,\xi)}$ on $W^{l_{\xi}([\phi])}_{\xi}([\phi])/\CC\phi$ determines a proper eigen subcycle of $Z_{\phi}$, the fibre of $q_1$ over $[\phi]\in \PP(W_{\xi})$. In other words  the morphism
	$m_{F(\bullet, \xi)}$ in \eqref{mFxi-morphism} determines the spectral
	scheme ${\cal S}_{m_{F(\bullet, \xi)}}$ over $\PP(W_{\xi})$ and
	$q_1$ factors through that scheme.  This is recorded in the  following diagram 
	\begin{equation}\label{xi-catdata+geom}
		\xymatrix{
			\Sigma^{\circ}_r \ni [\xi] \ar@2{->}[rr]&&
			*+[F]{\txt{Categorical data:\\
					$\OO_{\PP(W_{\xi})} (-1)
					\stackrel{F(\bullet,\xi) }{\longrightarrow} q_{1 \ast} \OO_{{\cal Z}_C} \otimes \OO_{\PP(W_{\xi})}$
					\\
					$\OO_{\PP(W_{\xi})} (-1)
					\stackrel{m_{F(\bullet,\xi)} }{\longrightarrow} {\cal{E}}nd(q_{1 \ast} \left(q^{\ast}_2\OO_{C} (K_C) \right))$}} \ar@2{->}[d]
				\\
				&& *+[F]{\txt{ Spectral scheme ${\cal S}_{m_{F(\bullet, \xi)}} \longrightarrow \PP(W_{\xi})$ of $m_{F(\bullet,\xi)}$,\\
				the factorization of $q_1$ over $\PP(W_{\xi})$:
			\\
		$q_1: q^{-1}_1(\PP(W_{\xi})) \longrightarrow {\cal S}_{m_{F(\bullet, \xi)}} \longrightarrow \PP(W_{\xi})$ }}			
		}
	\end{equation}
Furthermore, the multiplicities of eigenvalues are geometrically meaningful: they are related to the dimensions
of the linear systems of eigen subcycles; namely, for an eigen value $\lambda$ of
$m_{F(\phi,\xi)}$ on the quotient $W^{l_{\xi}([\phi])}_{\xi}([\phi])/\CC\phi$ we have the equality:
$$
dim(V_{\lambda})=dim(|Z_{\phi}(\lambda)|),
$$
here $V_{\lambda}$ is the $\lambda$-eigen space of $\beta(\phi,\bullet)$
in Proposition \ref{pro:beta-eigen}, $Z_{\phi}(\lambda)$ the $\lambda$-eigen subcycle of $Z_{\phi}$ and $|Z_{\phi}(\lambda)|$, the linear system defined by $Z_{\phi}(\lambda)$. In addition, outside of the branch locus of $q_1$, the intersection of  
$\PP(W_{\xi})$ with the dual variety $C^{\vee}$, the endomorphism $\beta(\phi,\bullet)$ is semisimple and we can encode its spectral data
in the associated zeta function
$$\label{Zetabeta}
{\mathfrak z}(\beta(\phi,\bullet),t)=exp\left(\sum^{\infty}_{n=1} \frac{tr(\beta^n(\phi,\bullet))t^n}{n} \right) =\prod_{\lambda} (1-\lambda t)^{-dim(|Z_{\phi}(\lambda)|)}=\frac{1}{det({\bf 1} -t\beta(\phi,\bullet))},
$$
where the product is taken over the distinct eigen values of $\beta(\phi,\bullet)$ and ${\bf 1}$ in the determinant stands for the identity endomorphism of 
$W^{l_{\xi}([\phi])}_{\xi}([\phi])/\CC\phi$. Of course the above encodes only the last piece of $([\xi],[\phi])$-filtrations $W^{\bullet}_{\xi}([\phi])$. One wonders if the remaining part of the filtrations leads to other structures.
  We will now turn to that remaining part. This will uncover other facets of the theory which are linked to topics ranging from toric geometry to quantum type invariants. The main underlying object providing those links is a certain bipartite (almost)trivalent graph which will be described next.
  
\section{Bipartite graph of algebraic K\"ahler structure of $W_{\xi}$}
In the previous sections we learned that the space $W_{\xi}$  comes with the direct sum decomposition
\begin{equation}\label{eq:orthdecWxi}
W_{\xi}=\bigoplus^{l_{\xi}([\phi])}_{s=0}  P^s([\xi],[\phi]),
\end{equation}
for every point $([\xi],[\phi])$ in ${\bf P}$ lying over the stratum $\Sigma^{\circ}_r$. Furthermore, the decomposition is equipped with  the collection $A_{[\xi]} ([\phi])$ of linear maps 
\begin{equation}\label{alpha2xiphi-assign}
\alpha^{(2)}_{\xi}(\phi,\bullet):W_{\xi} \longrightarrow \HKC; 
\end{equation}
for any nonzero $\xi$ and $\phi$ lying over $[\xi]$ and $[\phi]$ respectively. We wish to describe $A_{[\xi]} ([\phi])$. Let $\OO_{\bf P} (-1)$ be the tautological line bundle over ${\bf P}$, that is,
$$
\OO_{\bf P} (-1):=p^{\ast}_1 (\OO_{\PP(H^1(\Theta_C))} (-1))\otimes p^{\ast}_{1}(\OO_{\PP(\HKC)} (-1)) \otimes \OO_{\bf P}.
$$
Restrict it to the stratum $\widetilde{\Sigma^0_r}=p^{-1}_1 (\Sigma^0_r)$; this restriction will be denoted $\OO_{\widetilde{\Sigma^0_r}} (-1)$. The assignment \eqref{alpha2xiphi-assign} can be rephrased as a map
$$
\OO_{\PP(H^1(\Theta_C)),[\xi]} (-1) \otimes \OO_{\PP(W_{\xi}),[\phi]} (-1)  =\OO_{\widetilde{\Sigma^0_r}, ([\xi], [\phi]))} (-1) \longrightarrow Hom(W_{\xi}, \HKC).
 $$
 and the set $A_{[\xi]} ([\phi])$ is the image of the above map.
 The first thing to notice is that the above map intertwines the natural $\CC^{\times}$-action on the both sides:
 \begin{equation}
 	\alpha^{(2)}_{t\xi} (u\phi, \bullet)=ut	\alpha^{(2)}_{\xi} (\phi, \bullet), \forall t, u \in \CC^{\times};
 \end{equation}
 with respect to the factor $\OO_{\PP(W_{\xi}),[\phi]} (-1)$ follows from the fact that $\alpha^{(2)}_{\xi}$ is a linear map
 $$
 \mbox{$\bigwedge^2 W_{\xi} \longrightarrow \HKC$};
 	$$
 with respect to the $\CC^{\times}$-action on the factor $\OO_{\PP(H^1(\Theta_C)),[\xi]} (-1)$ comes from Remark \ref{rem:alpha-scalar}. Furthermore, 	    
two members $\alpha^{(2)}_{\xi}(\phi,\bullet)$ and $\alpha^{'(2)}_{\xi}(\phi,\bullet)$ of the collection lying over a fixed nonzero $\xi\otimes \phi \in \OO_{\PP(H^1(\Theta_C)),[\xi]} (-1)\otimes \OO_{\PP(W_{\xi}),[\phi]} (-1)$  are related by the equation
\begin{equation}\label{split-change-6}
\alpha^{'(2)}_{\xi}(\phi,\psi) -\alpha^{(2)}_{\xi}(\phi,\psi)=f(\phi)\psi -f(\psi)\phi, \,\,\forall \psi \in W_{\xi},
\end{equation}
where $f: W_{\xi} \longrightarrow \CC$ is a linear map, see \eqref{change-split-pf} in the proof of Lemma \ref{lem:cxi-map}. 

 The orthogonal decomposition \eqref{eq:orthdecWxi} gives the decomposition of the maps $\alpha^{(2)}_{\xi'}(\phi',\bullet)$ in $A_{[\xi]} ([\phi])$ into the direct sum
$$
\alpha^{(2)}_{\xi'}(\phi',\bullet)=\bigoplus \alpha^{t,s}
$$
and where the block $\alpha^{t,s}$ acts on the direct sum decomposition \eqref{eq:orthdecWxi} by mapping $ P^s([\xi],[\phi])$ to $ P^t([\xi],[\phi])$:
$$
\alpha^{t,s}: P^s([\xi],[\phi]) \longrightarrow P^t([\xi],[\phi]),
$$
the $(t,s)$-block of $\alpha^{(2)}_{\xi'}(\phi',\bullet)$ obtained by restricting 
$\alpha^{(2)}_{\xi'}(\phi',\bullet)$ to $P^s([\xi],[\phi])$ and then projecting to 
$P^t([\xi],[\phi])$. 

We also have learned that the summand $P^{l_{\xi}([\phi])}([\xi],[\phi])=W^{l_{\xi}([\phi])}_{\xi}$ is
$\alpha^{(2)}_{\xi'}(\phi',\bullet)$-invariant. So the quotient
$$
W_{\xi}/W^{l_{\xi}([\phi])}_{\xi}=\bigoplus^{l_{\xi}([\phi])-1}_{s=0}  P^s([\xi],[\phi])
$$
comes with the endomorphisms
\begin{equation}\label{alpha2-dirsum}
\alpha^{(2)}_{\xi'}(\phi',\bullet)=\bigoplus \alpha^{t,s} : 
\bigoplus^{l_{\xi}([\phi])-1}_{s=0}  P^s([\xi],[\phi]) \longrightarrow \bigoplus^{l_{\xi}([\phi])-1}_{s=0}  P^s([\xi],[\phi]).
\end{equation}
where we continue to use the same notation for the restriction of $\alpha^{(2)}_{\xi'}(\phi',\bullet)$ to $W_{\xi}/W^{l_{\xi}([\phi])}_{\xi}$.
The only components $\alpha^{t,s}$ appearing in the direct sum are those
allowed by Lemma \ref{lem:alpha-sum}. Furthermore, the blocks $\alpha^{t,s}$
above, for $t\neq s$, are independent of the choice of $\alpha^{(2)}_{\xi'}(\phi',\bullet)$ in the collection $A_{[\xi]}([\phi])$ for every fixed $\xi'$ and $\phi'$ lying over $[\xi]$ and $[\phi]$ respectively, while the blocks $\alpha^{s,s}$ differ by a scalar multiple of identity: for 
$\alpha^{(2)}_{\xi'}(\phi',\bullet)$ and $\alpha^{'(2)}_{\xi'}(\phi',\bullet)$ in $A_{[\xi]}([\phi])$ with
$$
\alpha^{'(2)}_{\xi'}(\phi',\psi) -\alpha^{(2)}_{\xi'}(\phi',\psi)=f(\phi')\psi,
$$
for $f\in W^{\ast}_{\xi}$, we have
$$
(\alpha^{'(2)})^{s,s} -(\alpha^{(2)})^{s,s} =f(\phi')id_{P^s}
$$
for every $s\in [0,l_{\xi}([\phi])-1]$, see Remark \ref{rem:choicesplit}.
Thus the map
$$
A_{[\xi]}([\phi]) \longrightarrow Hom (W_{\xi}/W^{l}_{\xi}([\phi]), \HKC/W^{l}_{\xi}([\phi]))
$$
identifies the subset
$$
A_{[\xi]}([\phi])(\xi',\phi' ) =\{\alpha^{'(2)}_{\xi'} (\phi', \bullet) | \phi' \in [\phi], \xi'\in [\xi] \} \cong \alpha^{(2)}_{\xi'} (\phi',\bullet) + f(\phi')id_{W_{\xi}/W^{l}_{\xi}([\phi])}
$$
where $f: W_{\xi} \longrightarrow \CC$ is a linear map. 
 Furthermore, taking the trace
\begin{equation}\label{tr-map}
	tr([\xi],[\phi]): A_{[\xi]}([\phi])(\xi',\phi') \longrightarrow \CC
\end{equation}
defined by the sum of traces of the degree preserving blocks of $\alpha^{'(2)}_{\xi'}(\phi',\bullet)$
$$
tr([\xi],[\phi]) (\alpha^{'(2)}_{\xi'}(\phi',\bullet)):=\sum^{l_{\xi}([\phi])-1}_{s=0} tr((\alpha^{'(2)}_{\xi'}(\phi',\bullet))^{s,s}),\,\, \forall \phi' \in [\phi], \xi'\in [\xi] ,
$$
gives a unique element of  $A_{[\xi]}([\phi])(\xi',\phi')$ of a given trace, that is, the trace map \eqref{tr-map} is an isomorphism.

Let us fix a value $c$ of the trace and denote by $A^c_{[\xi]}([\phi])$ the subset of maps in $A_{[\xi]}([\phi])$ having trace equals $c$.
The above tells us that $A^c_{[\xi]}([\phi])(\xi',\phi')$ contains a unique element which we denote $\alpha^{(2)}_{\xi',c} (\phi',\bullet)$. Furthermore, the $\CC^{\times}$-actions on vectors lying over $[\xi]$ and $[\phi]$ give the following relation
\begin{equation}\label{uv-action-on-c}
	\alpha^{(2)}_{u\xi,c} (v\phi,\bullet)=uv \alpha^{(2)}_{\xi,c} (\phi,\bullet)+\frac{c(1-uv)}{dim(W_{\xi}/W^{l}_{\xi})} id_{W_{\xi}/W^{l}_{\xi}}, \,\forall u,v \in \CC^{\times},
\end{equation}
where $\xi$ and $\phi$ are nonzero vectors lying over $[\xi]$ and $[\phi]$ respectively. We summarize the above discussion in the following.
\begin{lem}\label{lem:traceparam}
	Let $([\xi],[\phi])$ be a point of the cohomological correspondence ${\bf P}$
	with $[\xi]$ in the stratum $\Sigma^0_r$ and $[\phi] \in \PP(W_{\xi})$.
	Let $W^{\bullet}_{\xi}([\phi])$ be the $([\xi],[\phi])$-filtration of $W_{\xi}$ attached to $[\phi]$ and let
	$$
	W_{\xi}/W^{l_{\xi}([\phi])}_{\xi}=\bigoplus^{l_{\xi}([\phi])-1}_{s=0}  P^s([\xi],[\phi]).
	$$
	be the corresponding orthogonal decomposition of 	$W_{\xi}/W^{l_{\xi}([\phi])}_{\xi}$. Then the set of maps
	 $$
	A_{[\xi]}([\phi])(\xi,\phi)
	$$ 
	attached to nonzero vectors $\xi$ and $\phi$ lying over the point $([\xi],[\phi]) \in p^{-1}_1(\Sigma^0_r)$, is parametrized by the values of the trace map
	$$
	tr([\xi],[\phi]): A_{[\xi]}([\phi])(\xi,\phi) \longrightarrow \CC.
	$$
	More precisely, any two members $\alpha^{(2)}$ and $\alpha^{'(2)}$ of $A_{[\xi]}([\phi])(\xi,\phi)$ differ by a multiple of identity
	$$
	\alpha^{'(2)} - \alpha^{(2)}=f(\phi)id_{W_{\xi}/W^{l_{\xi}([\phi])}_{\xi}},
	$$
	where $f:W_{\xi} \longrightarrow \CC$ is a linear map. 
In particular, the off diagonal blocks of elements in $A_{[\xi]}([\phi])(\xi,\phi)$ mapping the summand $P^s([\xi],[\phi])$ to $P^t([\xi],[\phi])$, $s\neq t$,
	$$
	\alpha^{t,s}:P^s([\xi],[\phi]) \longrightarrow P^t([\xi],[\phi])
	$$
	are constant, while the diagonal ones are subject to the equation
	$$
	(\alpha^{'(2)})^{s,s} - (\alpha^{(2)})^{s,s}=f(\phi)id_{P^s([\xi],[\phi])},
	$$
	for every $s\in[0,l_{[\xi]}([\phi])-1]$.
	
	In addition, the natural actions of $\CC^{\times}$ on the lines $[\xi]$ and $[\phi]$ act on the elements of $A_{[\xi]}([\phi])$ as follows:
	$$
	\alpha^{(2)}_{t\xi}(u\phi,\bullet)=tu\alpha^{(2)}_{\xi}(\phi,\bullet),\,\,\forall t,u\in \CC^{\times}, \forall  \alpha^{(2)}_{\xi}(\phi,\bullet) \in A_{[\xi]}([\phi]).
	$$  
\end{lem}

\vspace{0.5cm}
 We will now recast the collection of maps $A_{[\xi]}([\phi])$ as a subset of representations of a quiver whose underlying graph denoted
$G_{l_{\xi}([\phi])}$ is defined as follows. Set $l:=l_{\xi}([\phi])$; our graph will have $l$ pairs of vertices colored black and white; we place white vertices along the horizontal line numbering them in the decreasing order, from left to right, by integers $i\in [0,l-1]$; on a horizontal level below the white vertices we place the black vertices labeled by $i'$'s so that each pair $i$ and $i'$ is vertically aligned. The edges of $G_{l}$ are drawn according to the blocks
$\alpha^{t,s}$: we draw the edge connecting the white vertex labeled $s$ with the black vertex labeled $t'$ whenever the block $\alpha^{t,s}(\xi,\phi)$ appears in \eqref{alpha2-dirsum}. Below is an illustration of the graph $G_4$.  
$$
\begin{tikzpicture}
	[place/.style={circle,draw=black,thick},
	transition/.style={circle,draw=black,fill=black}]
	\node (white3) at (0,2) [place] [label={above:$3$}] {};
	\node (black3) at (0,0) [transition] [label={below:$3'$}]{};
		\node (white2) at (2,2) [place] [label={above:$2$}]{};
		\node (black2) at (2,0) [transition] [label={below:$2'$}]{};
		\node (white1) at (4,2) [place] [label={above:$1$}]{};
		\node (black1) at (4,0) [transition] [label={below:$1'$}]{};
		\node (white0) at (6,2) [place] [label={above:$0$}] {};
		\node (black0) at (6,0) [transition] [label={below:$0'$}]{};
		\draw[red, ultra thick][-] (white1) to (black1);
		\draw[red,ultra thick][-] (white2) to (black2);
		\draw[red,ultra thick][-] (white3) to (black3);
		\draw[red, ultra thick][-] (white0) to (black0);
				\draw[blue, thick][-] (white3) to (black2);	
				\draw[blue, thick][-] (white2) to (black1);
				\draw[blue, thick][-] (white1) to (black0);
				\begin{scope}[thick]
				\draw[red, dotted][-] (white2) to (black3);
				\draw[red,dotted][-] (white1) to (black2);
				\draw[red,dotted][-] (white1) to (black3);
				\draw[red,dotted][-] (white0) to (black1);
				\draw[red,dotted][-] (white0) to (black2);
				\draw[red,dotted][-] (white0) to (black3);
				\end{scope}
\end{tikzpicture}
$$
The edges `solid red' are the ones corresponding to the grading preserving blocks $\alpha^{s,s}$, the edges colored in blue correspond to
blocks $\alpha^{s-1,s}$ lowering the degree by $1$. All other edges, dotted red, correspond to the degree increasing blocks. Thus the graph $G_{l}$ is a bipartite directed graph: the vertices are colored black and white and the edges of the graph connect the vertices of different colors only and we agree that they are oriented from `white' to `black'. 

By construction $G_{l}$ comes with the set of `defining' representations $\rho_c (\xi,\phi)$, attached to every pair nonzero vectors $\xi$ and $\phi$ lying over a point $([\xi],[\phi]) \in p^{-1}_1 (\Sigma^0_r)$ and where $c$ is a parameter taking values in $\CC$: we put the direct summand $P^s([\xi],[\phi])$ at the vertices $s$ and $s'$; given a value $c\in \CC$, according to Lemma \ref{lem:traceparam}, there is a unique element in $A_{[\xi]}([\phi])(\xi,\phi)$ having the trace equals $c$. That map will is denoted $\alpha^{(2)}_{\xi,c} (\phi,\bullet)$;
the orthogonal decomposition 
\begin{equation}\label{Wxi-Wlxi-orthdecopm}
W_{\xi}/W^{l_{\xi}([\phi])}_{\xi}=\bigoplus^{l_{\xi}([\phi])-1}_{s=0}  P^s([\xi],[\phi]);
\end{equation}
produces the decomposition of $\alpha^{(2)}_{\xi,c} (\phi,\bullet)$ into blocks
 denoted $\alpha^{t,s}_c (\xi,\phi)$ and we color the edge 
$(s) \to (t')$ of the graph $G_l$
 with the map $\alpha^{t,s}_c(\xi,\phi):P^s([\xi],[\phi])\longrightarrow P^t([\xi],[\phi]) $. This determines the map
 \begin{equation}\label{rhoxiphi}
 	\rho ([\xi],[\phi]):\OO^{\times}_{\widetilde{\Sigma^0_r},([\xi],[\phi])} (-1) \times \CC \longrightarrow {\mathfrak Reps}(G_l)
 \end{equation}
which assigns to a triple $(\xi \otimes \phi,c)$ in $\OO^{\times}_{\widetilde{\Sigma^0_r},([\xi],[\phi])} (-1) \times \CC $ the representation 
$$
\rho_c(\xi,\phi):=\{\alpha^{t,s}_c (\xi,\phi): P^s([\xi],[\phi])\longrightarrow P^t([\xi],[\phi])\}
$$
of the quiver $G_l$. 
We also remind the reader that once $\xi$ and $\phi$ are fixed, the blocks  $\alpha^{t,s}_c(\xi,\phi)$, for $s\neq t$, are {\it independent} of $c$ and the diagonal blocks differ by a multiple of identity when $c$ varies, that is if $c'$ is another trace value, then we have
 \begin{equation}\label{tr-v-change}
 \alpha^{s,s}_{c'}(\xi,\phi) -\alpha^{s,s}_c(\xi,\phi) =\frac{c'-c}{dim(W_{\xi}/W^{l}_{\xi})} id_{P^s},
\end{equation}
for all $s\in [0,l-1]$. We omit the value of trace parameter whenever a block is independent of that value. 

  \vspace{0.2cm}
  Next we give an `almost' trivalent approximation of the quiver $G_l$. To motivate this step, we fix
 a trace value $c$ and consider an element $\psi \in P^s([\xi],[\phi])$, for $s\in [0,l-1]$. For  $\alpha^{(2)}_{\xi,c}(\phi,\bullet)$ in $A_{[\xi]}([\phi])(\xi,\phi)$ of trace $c$, the  
 image $\alpha^{(2)}_{\xi,c} (\phi,\psi)$ of  $\alpha^{(2)}_{\xi,c}(\phi,\psi)$ in $W_{\xi}/W^{l}_{\xi}$
 admits the expansion
 $$
\alpha^{(2)}_{\xi,c} (\phi,\psi) =\sum_{t\geq s-1}\alpha^{t,s}_c (\psi)
 $$ 
  with respect to the orthogonal decomposition \eqref{Wxi-Wlxi-orthdecopm}; the parameters $(\xi,\phi)$ on the right are omitted to simplify the notation. We distinguish the first three terms of the expansion
 $$
 \alpha^{(2)}_{\xi,c} (\phi,\psi) =\alpha^{s-1,s} (\psi) + \alpha^{s,s}_c (\psi) + \alpha^{s+1,s} (\psi) +\sum_{t\geq s+2}\alpha^{t,s} (\psi)
 $$
 and view them as the principal part of $\alpha^{(2)}_{\xi,c} (\phi,\psi)$ and the remaining terms as `higher order' terms. By analogy, suppressing all edges
 of the graph $G_{l}$ to the left from the edge 
 $(s) \mapsto (s+1)'$ at every white vertex we obtain the {\it principal part} of the graph $G_{l}$. It will be denoted $PG_{l}$. For example, the principal part $PG_{4}$ of $G_4$ is depicted below.
$$
\begin{tikzpicture}
	[place/.style={circle,draw=black,thick},
	transition/.style={circle,draw=black,fill=black}]
	\node (white3) at (0,2) [place] [label={above:$3$}] {};
	\node (black3) at (0,0) [transition] [label={below:$3'$}]{};
	\node (white2) at (2,2) [place] [label={above:$2$}]{};
	\node (black2) at (2,0) [transition] [label={below:$2'$}]{};
	\node (white1) at (4,2) [place] [label={above:$1$}]{};
	\node (black1) at (4,0) [transition] [label={below:$1'$}]{};
	\node (white0) at (6,2) [place] [label={above:$0$}] {};
	\node (black0) at (6,0) [transition] [label={below:$0'$}]{};
	\begin{scope}[ultra thick]
	\draw[red] (white1) to (black1);
	\draw[red] (white0) to (black0);
	\draw[red] (white2) to (black2);
	\draw[red] (white3) to (black3);
	\end{scope}
\begin{scope}[thick]
	\draw[blue][-] (white3) to (black2);	
	\draw[blue][-] (white2) to (black1);
	\draw[blue][-] (white1) to (black0);
	\end{scope}
	\draw[red, dotted,thick] (white2) to (black3);
	\draw[red,dotted,thick] (white1) to (black2);
	\draw[red,dotted,thick] (white0) to (black1);
\end{tikzpicture}
$$
The graph $PG_l$ is almost trivalent in the sens that all vertices, with the exception of $0$, $0'$, $l-1$, $(l-1)'$, are trivalent.
\subsection{Nonabelian Dolbeault variety ${\bf H^{1,0}}(PG_{l})$}

Denote the length $l_{[\xi]}([\phi])$ of the $([\xi],[\phi])$-filtration of $W_{\xi}$ by $l$ to simplify the notation. To the graph $PG_{l}$ is naturally attached a toric projective variety. This is done as follows: we label the edges 
$(s) \mapsto (s-1)'$ by the indeterminate $X_s$ and
 the edges $(s-1) \mapsto (s)'$ by the indeterminate $Y_s$, for every
 $s=1,\ldots, l-1$; all vertical edges $(s)\mapsto (s)'$ are labeled by the indeterminate $T$.
We think of the indeterminates ${T,X_1,Y_1,\ldots, X_{l-1},Y_{l-1}}$ as linear forms
on the complex vector space
$V_{2l-1}$ spanned by $2(l-1)$ vertices $\{(s),(s)' | s=1,\ldots,l-1 \}$ and the vertical edge $e_0:=(0)\mapsto (0)'$, that is, $\{e_0, (1),(1)', \ldots, (l-1),(l-1)'\}$ is a basis of $V_{2l-1}$ and the linear forms
$$
\{T,X_1,Y_1,\ldots, X_{l-1},Y_{l-1}\}
$$
is the corresponding dual basis of $V^{\ast}_{2l-1}$.
  
For every $s \in [1,l-1]$ we write the quadratic relation
$$
X_sY_s=T^2. 
$$
 These are viewed as quadrics in the projective space $\PP(V_{2l-1})$.

Define the {\it nonabelian Dolbeault variety} ${\bf H^{1,0}}(PG_{l})$ as the complete intersection in the projective space $\PP(V_{2l-1})$ of $(l-1)$ quadrics above. A similar variety appears in the previous work of the author in a context of vector bundles on projective surfaces, see \cite{R3} and \cite{R4} for an overview. All the properties of ${\bf H^{1,0}}(PG_{l})$ discussed below could be found in those works.

\begin{pro}\label{pro:NADolb}
	The variety ${\bf H^{1,0}}(PG_{l})$ is a singular Fano toric variety in the projective space $\PP(V_{2l-1})=\PP^{2(l-1)}.$ It has dimension $(l-1)$ and degree $2^{l-1}$.
	
	The hyperplane $T=0$ intersects ${\bf H^{1,0}}(PG_{l})$ along the divisor denoted $H_0$. This divisor is the union
	of $2^{l-1}$ projective subspaces. More precisely, for every subset
	$A \subset [1,l-1]$ denote by $A^c$ its complement in $[1,l-1]$ and let $\Pi_A$ be the projective subspace of the hyperplane
	$T=0$ spanned by the points of $\PP(V_{2l-1})$ underlying the vectors
	$$
	\{(s),(t)' | s\in A,\, t\in A^c\}  \subset V_{2l-1};
	$$
	then we have
	$$
	H_0 =\bigcup_{A\subset [1,l-1]} \Pi_A,
	$$
	where the union is taken over all subsets $A$ of $[1,l-1]$.
	
	The complement $U_0$ of $H_0$ in ${\bf H^{1,0}}(PG_{l})$ has natural identification with the algebraic torus $(\CC^{\times})^{l-1}$
	$$
	U_0 ={\bf H^{1,0}}(PG_{l}) \setminus H_0 \cong (\CC^{\times})^{l-1}.
	$$
\end{pro}
\begin{pf}
	The assertions of the proposition about the dimension and the degree of the complete intersection ${\bf H^{1,0}}(PG_{l})$ are immediate from the definition.
	By adjunction formula, the dualizing sheaf $\omega_{{\bf H^{1,0}}(PG_{l})}$ is given by the formula
	$$
	\omega_{{\bf H^{1,0}}(PG_{l})}=\OO_{\PP^{2(l-1)}}(2(l-1)-2(l-1)-1) \otimes \OO_{{\bf H^{1,0}}(PG_{l})}=\OO_{{\bf H^{1,0}}(PG_{l})} (-1).
		$$
		Hence ${{\bf H^{1,0}}(PG_{l})}$ is a Fano variety.
		
		The torus $(\CC^{\times})^{l-1}$ acts on ${\bf H^{1,0}}(PG_{l})$ as follows. Let ${\bf z}=(t,x_1, y_1,\ldots,x_{l-1}, y_{l-1})$ be a closed point of ${\bf H^{1,0}}(PG_{l})$ and let
		${\lambda}=(\lambda_1,\ldots,\lambda_{l-1}) \in (\CC^{\times})^{l-1}$, then ${\bf \lambda}$ acts on ${\bf z}$ by the formula
		$$
		{\bf \lambda} {\bf z}:=(t,\lambda_1x_1, \lambda^{-1}_1y_1,\ldots,\lambda_{l-1}x_{l-1}, \lambda^{-1}_{l-1}y_{l-1}).
		$$
	We turn to the divisor $H_0$.
	The hyperplane $T=0$ in the projective space $\PP(V_{2l-1})$ corresponds to the projectivization of the subspace of $V_{2l-1}$ spanned by the vectors
	$$
	\{ (1),(1)', \ldots, (l-1),(l-1)'\}. 
	$$
	The divisor $H_0$ is the subscheme in the hyperplane determined by the equations
	$$
	X_sY_s=0, \, s=1,\ldots,l-1.
	$$
	From this it follows that for each subset $A$ in $[1,l-1]$ the projective subspace $\Pi_A$ spanned by the vectors
	$$
	\{(s),(t)'| s\in A, t\in A^c\}
	$$
	is contained in $H_0$. The union of $\Pi_A$, as $A$ runs through all subsets of $[1,l-1]$, is a subscheme of $H_0$. The degree of that subscheme
	is $2^{l-1}$. From the first assertion of the proposition this is the degree of $H_0$. Hence the asserted equality
	$$
	H_0=\bigcup_{A\subset[1,l-1]} \Pi_A.
	$$
	
	Next we consider the complement of the divisor $H_0$ in 
	${\bf H^{1,0}}(PG_{l})$. It contains the point
	$$
	{\bf z_1}=(1,1,1,\ldots,1,1)
	$$
	all of whose coordinates equal $1$. The orbit of that point under the action of the torus $(\CC^{\times})^{l-1}$ is
	$$
	{Orbit({\bf z_1})}=(\CC^{\times})^{l-1} {\bf z_1}=\{(1,\lambda_1,\lambda^{-1}_1,\ldots,\lambda_{l-1},\lambda^{-1}_{l-1} ) | \, \forall \lambda=(\lambda_1,\ldots,\lambda_{l-1}) \in (\CC^{\times})^{l-1}\}.
	$$
	The isotropy group of ${\bf z_1}$ is obviously trivial. Hence
	the isomorphism
	$$
	 {Orbit({\bf z_1})} \cong (\CC^{\times})^{l-1} .
	 $$
	 The points of the orbit $Orbit({\bf z_1})$ have the first coordinate nonzero, that is, we have an inclusion
	 $$
	  {Orbit({\bf z_1})} \subset {\bf H^{1,0}}(PG_{l})\setminus H_0.
	  $$
	  On the other hand a point 
	  $$
	  [{\bf z}]=[(t,x_1,y_1,\ldots,x_{l-1},y_{l-1})] \in {\bf H^{1,0}}(PG_{l})\setminus H_0
	  $$
	  has the coordinate $t\neq 0$. Hence
	  $$
	  [{\bf z}]=[(t,x_1,y_1,\ldots,x_{l-1},y_{l-1})]=[(1,\frac{x_1}{t},\frac{y_1}{t},\ldots, \frac{x_{l-1}}{t},\frac{y_{l-1}}{t} )].
	  $$
	  Since $x_i y_i=t^2$, for all $i\in[1,l-1]$, we deduce
	  $$
	  \frac{y_i}{t}=\left( \frac{x_i}{t} \right)^{-1}, \forall i\in [1,l-1].
	  $$
	 Setting $\displaystyle \lambda_i=\frac{x_i}{t}$ gives
	 $$
 [{\bf z}]=[(1,\frac{x_1}{t},\frac{y_1}{t},\ldots, \frac{x_{l-1}}{t},\frac{y_{l-1}}{t} )]=[(1,\lambda_1,\lambda^{-1}_1,\ldots,\lambda_{l-1},\lambda^{-1}_{l-1} )] \in Orbit({\bf z_1}).
 $$
 Hence the inclusion in the opposite direction
 $$
 {\bf H^{1,0}}(PG_{l})\setminus H_0 \subset Orbit({\bf z_1})).
 $$	  
\end{pf}

The vector subspace
$$
V_{2l-2} =\CC\{(1),(1)', \ldots, (l-1),(l-1)'\}
$$
has a natural symplectic structure: define the pairing
$$
\Omega:V_{2l-2} \times V_{2l-2} \longrightarrow \CC
$$
by the rules
$$
\Omega((r), (r)')=-\Omega((r)', (r))=1, \forall r \in[1,l-1]
$$
and zero for any other pair of the basis
$$
\{(1),(2),\ldots,(l-1), (1)',(2)',\ldots,(l-1)'\}
$$
of $V_{2l-2}$, that is, the matrix of the pairing with the respect to the above basis
$$
\Omega=\begin{pmatrix}
	0_{l-1}&I_{l-1}\\
	-I_{l-1}&	0_{l-1}
\end{pmatrix}
$$
gives the standard symplectic form; $I_{l-1}$ (resp. $0_{l-1})$) is the  identity  (resp. zero) matrix of order ${(l-1)\times(l-1)}$.

For each subset $A \subset[1,l-1]$ the vector subspace
$$
V_A =\CC\{(s),(t)' |s\in A, t\in A^c\}
$$
overlying the irreducible component $\Pi_A$ of the divisor $H_0$ is an isotropic subspace of $V_{2l-2}$ with respect to $\Omega$. The dimension of 
$V_{A}$ is $(l-1) =\HA dim(V_{2l-2})$. Hence it is a Lagrangian subspace of the symplectic space $(V_{2l-2},\Omega)$. Thus we obtain
\begin{lem}\label{lem:Lagrangian}
	For every subset $A\subset[1,l-1]$, the vector subspace
	$$
	V_A =\CC\{(s),(t)' |s\in A, t\in A^c\}
	$$
	overlying the irreducible component $\Pi_A$ of $H_0$ is a Lagrangian subspace of the symplectic space $(V_{2l-2},\Omega)$. For this reason the divisor $H_0$ in Proposition \ref{pro:NADolb} will be called {\rm Lagrangian cycle} of ${\bf H^{1,0}}(PG_{l})$.
\end{lem}

The Fano variety ${\bf H^{1,0}}(PG_{l})$ is singular. The following describes its singular locus.
\begin{pro}\label{pro:singNADolb}
	The singular locus $Sing({\bf H^{1,0}}(PG_{l}))$ of 
	${\bf H^{1,0}}(PG_{l})$ is the union
	$$
	Sing({\bf H^{1,0}}(PG_{l}))= \bigcup_{A\neq B} \Pi_A \cap \Pi_B,
	$$
	where the union is over all pairs of distinct subsets $A$ and $B$ of $[1,l-1]$. Each intersection $\Pi_A \cap \Pi_B$ is the projectivization of
	the linear subspace of $V_{2l-1}$ spanned by the vectors
	$$
	\{(s), (t)' | s\in A\cap B,\, t\in A^c \cap B^c\}.
	$$
\end{pro} 
\begin{pf}
	The proof is straightforward and we omit it. The interested reader could consult \cite{R3} for details.
\end{pf}
\begin{example}\label{exDolbl=3}
	We look at the case $l=3$, the first admissible value in our considerations. This gives $V_{2l-1}=V_5$ the vector space of dimension
	$5$ with a distinguished basis spanned by the symbols
	$$
	e^0, (1), (2), (1)', (2)'.
	$$
	The dual basis 
	$$
	T,X_1,X_2,Y_1,Y_2
	$$
	of $V^{\ast}_5$ provides the homogeneous linear forms on $\PP(V_5)\cong \PP^4$. The Dolbeault variety ${\bf H^{1,0}}(PG_{3})$ is the complete intersection of $2$ quadrics in $\PP^4$:
	$$
	X_1 Y_1 =T^2, \,\,\, X_2 Y_2 =T^2.
	$$
	This a toric singular del Pezzo surface of degree $4$ in $\PP^4$. Its Lagrangian
	cycle
	$$
	H_0 =\{T=0\}
	$$
	is the `square' formed by four lines
	$$
	\begin{tikzpicture}
		\draw[very thick, black]
		(-3,0) -- (3,0) 
		(-3.5,0) node {$\Pi_{(1),(2)}$}; 
		\draw[thick, black]
		(-3,3) -- (3,3)
		(-3.5,3) node {$\Pi_{\emptyset}$};
		\draw[thick, black]
		(-2,-1) -- (-2,4)
		(-1.8,-1.3) node {$\Pi_{(1)}$};
		\draw[very thick, black]
		(2,-1) to (2,4)
		(2.1,-1.3) node {$\Pi_{(2)}$};
		\filldraw[red]  (-2,0) circle [radius=3pt] [label={below left:$3$}] {};
		\filldraw[red]  (-2,3) circle [radius=3pt]; 
		\filldraw[red]  (2,0) circle [radius=3pt];
		\filldraw[red]  (2,3) circle [radius=3pt];
			\end{tikzpicture}
	$$		
	The red dots are the points of intersections of the lines:
	$$
	\begin{gathered}
	\Pi_{(1),(2)} \cap \Pi_{(1)}=[0:1:0:0:0],\,\,\, \Pi_{(1),(2)} \cap \Pi_{(2)}=[0:0:1:0:0],
	\\
	 \Pi_{(2)} \cap \Pi_{\emptyset}=[0:0:0:1:0],\,\,\,
	\Pi_{(1)} \cap \Pi_{\emptyset}=[0:0:0:0:1].
	\end{gathered}
	$$
These four points form the singular locus of ${\bf H^{1,0}}(PG_{3})$.
The Lagrangian cycle $H_0$ is a degenerate curve of genus one. A general hyperplane section of ${\bf H^{1,0}}(PG_{3})$ is a smooth curve of genus $1$ with four marked points, the points of intersection with $H_0$.
\end{example}

\begin{example}\label{exDolbl=4}
	For $l=4$, we have $V_{2l-1}=V_7$ the vector space of dimension
	$7$ with a distinguished basis spanned by the symbols
	$$
	e^0, (1), (2), (3), (1)', (2)', (3)'.
	$$
	The dual basis 
	$$
	T,X_1,X_2,X_3, Y_1,Y_2,Y_3
	$$
	of $V^{\ast}_5$ provides the homogeneous linear forms on $\PP(V_7)\cong \PP^6$. The Dolbeault variety ${\bf H^{1,0}}(PG_{4})$ is the complete intersection of $3$ quadrics in $\PP^6$:
	$$
	X_1 Y_1 =T^2, \,\,\, X_2 Y_2 =T^2,\,\, X_3 Y_3 =T^2.
	$$
	It is a singular toric Fano threefold of degree $8$. The Lagrangian cycle
	$H_0$ is a union of $8=2^3$ planes in $\PP^5=\{T=0\}$; the planes $\Pi_A$, the irreducible components of $H_0$, are labeled by the subsets $A$ of $\{1,2,3\}$. 
	
	The singular locus $Sing({\bf{ H^{1,0}}}(PG_{4}))$ consists of the union
of $12$ lines: there are three lines in the plane $\Pi_{(1),(2),(3)}$ spanned by the pairs of distinct points $(i)$, $(j)$, for $i,j \in \{1,2,3\}$; another three lines in the `opposite' plane $\Pi_{\emptyset}$, spanned by the pairs of distinct points $(i)'$, $(j)'$, for $i,j \in \{1,2,3\}$; the remaining lines join each point $(i)$ in $\Pi_{(1),(2),(3)}$ with a point $(j)'$ in $\Pi_{\emptyset}$, for $j \in \{1,2,3\}\setminus \{i\}$.

 Each point $(i)$ (resp, $(i)'$) lies precisely on four distinct irreducible components of $H_0$. Hence the six points
 $$
 \{(1),(2),(3), (1)',(2)',(3)',\}
 $$
 are the singular points of multiplicity four. All other points of $Sing({\bf{ H^{1,0}}}(PG_{4}))$ are ordinary double points. In particular, a general hyperplane section of ${\bf H^{1,0}}(PG_{4})$ is a singular $K3$ surface with $12$ distinct ordinary double points.  
 
 Denote by $H$ such a general hyperplane section. In addition to 12 ordinary double points, $H$ contains $8$ lines, the intersections of the hyperplane defining $H$ with the planes $\Pi_A$'s forming the Lagrangian cycle $H_0$. These $8$ lines and $12$ ordinary double points form the configuration $(8_3, 12_2)$, that is,
 every line passes through exactly $3$ double points and every double point
 lies in precisely 2 lines. Thus the configuration can be visualized as a three dimensional cube : the vertices of the cube correspond to $8$ lines and the edges of the cube correspond to the points of intersections of the lines, the twelve ordinary double points.
 
 Let $\sigma: \widetilde{H} \longrightarrow H$ be the minimal desingularization of $H$. So $\widetilde{H}$ is a $K3$ surface with twenty
 $(-2)$ curves: there are $12$ curves which are the inverse images of $\sigma$ of the $12$ double points of $H$, and $8$ strict transforms of lines.  The dual graph of this configuration is formed by $12$ edges of a cube and $20$ vertices: $8$ vertices of the cube and $12$ midpoints of the edges. In the following drawing
 \begin{equation}\label{cubeK3}
 \begin{tikzpicture}
 	\draw[very thick, black]
 	(-2,0) -- (2,0) --(2,4) -- (-2,4) -- (-2,0);
 	\draw[very thick, black]
 	(2,0) -- (4,1) -- (4,5) -- (2,4); 
 	\draw[very thick, black]
 	(4,5) -- (0,5) -- (-2,4);
 	\draw[thin, black]
 	(-2,0) -- (0,1) -- (4,1);
 	\draw[thin, black]
 	(0,1) -- (0,5);
 	\filldraw[red]  (0,0) circle [radius=2pt] {};
 	\filldraw[red]  (-2,2) circle [radius=2pt]; 
 	\filldraw[red]  (2,2) circle [radius=2pt];
 	\filldraw[red]  (0.2,4) circle [radius=2pt];
 	\filldraw[red]  (-1,4.5) circle [radius=2pt];
 	\filldraw[red]  (2,5) circle [radius=2pt];
 		\filldraw[red]  (3,4.5) circle [radius=2pt];
 			\filldraw[red]  (3,0.5) circle [radius=2pt];
 			\filldraw[red]  (4,3) circle [radius=2pt];
 			\filldraw[red]  (0,3) circle [radius=2pt];
 			\filldraw[red]  (-1,0.5) circle [radius=2pt];
 			\filldraw[red]  (1.6,1) circle [radius=2pt] {};
 			\filldraw[black]  (-2,0) circle [radius=2pt] {};
 			\filldraw[black]  (-2,4) circle [radius=2pt] {};
 				\filldraw[black]  (2,0) circle [radius=2pt] {};
 					\filldraw[black]  (2,4) circle [radius=2pt] {};
 						\filldraw[black]  (4,5) circle [radius=2pt] {};
 							\filldraw[black]  (0,5) circle [radius=2pt] {};
 				\filldraw[black]  (0,1) circle [radius=2pt] {};
 					\filldraw[black]  (4,1) circle [radius=2pt] {};
 \end{tikzpicture}
\end{equation}
the red dots are the midpoints;  two dots are connected by the half-edge of the cube if and only if the corresponding $(-2)$ curves intersect; in this case the intersection is transversal and equals one. 

The dual graph of the configuration of twenty $(-2)$ can be used to uncover other special features of $\widetilde{H}$:

- the $(-2)$ curves of each face $f$ of the cube form the extended Dynkin diagram of type $\widetilde{A}_8$; this gives the nef effective divisor
$L_f$; the corresponding line bundle $\OO_{\widetilde{H}} (L_f)$ is base point free and defines the genus one fibration
$$
|L_f|:\widetilde{H} \longrightarrow \PP^1;
$$
the pencil $|L_f|$ is spanned by the pair of divisors $L_f$ and $L_{f^{opp}}$, where $f^{opp}$ is the face of the cube opposite to $f$;
those two divisors are the singular fibres of type $I_8$ on the Kodaira list of classification of singular fibres of an elliptic surface; in addition, the midpoints of four edges of the cube connecting the opposite faces $f$ and $f^{opp}$ correspond to four sections of the fibration $|L_f|$; thus $|L_f|$ is an elliptic fibration with precisely $2$ singular fibres of type $I_8$; the surface $\widetilde{H}$ comes
with three distinct elliptic fibrations of this type - each for the pair of opposite faces of the cube;

\vspace{0.2cm}
- each edge $e$ of the cube gives an elliptic fibration
$$
|L_e|:\widetilde{H} \longrightarrow \PP^1
$$
with one singular fibre $L_e$ of type $I_{12}$; that divisor is seen on the cube as follows: take the two faces of the cube {\rm adjacent along the edge $e$} and consider the boundary with the edge $e$ {\rm deleted}; the divisor $L_e$ is formed  by summing up once the $(-2)$ curves found on that boundary;
look at the edge $e^{opp}$, the common edge of two faces {\rm opposite} to the ones adjacent to $e$; the $(-2)$ curves of that edge are disjoint from $L_e$, hence they are part of another reducible fibre of the fibration $|L_e|$; call that fibre $L'_e$, from the Kodaira's list of singular fibres we deduce that  $L'_e$ must contain $(-2)$ curves {\rm different} from the ones of the cube; thus each edge fibration $|L_e|$ contributes at least one more $(-2)$ curve different from the $20$ accounted by the dots on the cube; this way we obtain at least $12$ more $(-2)$ curves on $\widetilde{H}$.

 At this point it should be clear that $K3$ surface $\widetilde{H}$ is very special. To summarize the discussion, we list some of its features:

$\bullet$ $ \widetilde{H}$ admits at least $32$ distinct $(-2)$ curves; $20$ of them form the dual graph of the cube, see \eqref{cubeK3};

\vspace{0.2cm}
$\bullet$ three distinct elliptic fibrations   
$$
\pi_i : \widetilde{H} \longrightarrow \PP^1
$$
surjective morphisms, $i=1,2,3$; the fibres of each morphism are curves of genus one; each morphism comes with two distinguished fibres of type $I_8$ and four disjoint sections;

\vspace{0.2cm}
$\bullet$ the N\'eron-Severi group $NS(\widetilde{H})$ of $\widetilde{H}$ has rank at least $17$.

\vspace{0.2cm}
 The details and further considerations of these $K3$ surfaces will be given elsewhere; the $K3$-surfaces with twenty $(-2)$-curves in the configuration of a cube have been studied in \cite{Ka-Ko}\footnote{the author thanks Xavier Roulleau for pointing out this reference.} and \cite{Pe-S}  
in different contexts.
 		
\end{example}
	
Though the pair $({\bf H^{1,0}}(PG_{l}),H_0)$ is interesting from geometrical and combinatorial points of view, its definition so far has been formal. We will now discuss more conceptual way to understand this construction.

\subsection{${\bf H^{1,0}}(PG_{l})$ and Higgs structures}
Let us fix a nonzero subspace $W$ in $\HKC$ and let
$\Xi_W$ be the subspace of Kodaira-Spencer classes annihilating $W$, that is
$$
\Xi_W:=\{\xi \in H^1(\Theta_{C})|\, \xi \phi =0,\,\forall \phi \in W\}.
$$
We assume that the dimensions of $W$ and $\Xi_W$ are at least two. 

By definition we have the inclusions 
\begin{equation}\label{WinWxi}
W\subset W_{\xi}, \forall \xi \in \Xi_W.
\end{equation}
For every nonzero $\xi \in \Xi_W$ and every nonzero $\phi\in W_{\xi}$ we have the collection $A_{[\xi]} ([\phi])$ of linear maps
$$
\alpha^{(2)}_{\xi} (\phi,\bullet): W_{\xi} \longrightarrow \HKC 
$$
which come with the decomposition
$$
\alpha^{(2)}_{\xi} (\phi,\bullet)= \bigoplus \alpha^{t,s}(\xi,\phi): W_{\xi}/W^l_{\xi} ([\phi]) =\bigoplus^{l-1}_{s=0} P^s ([\xi],[\phi]) \longrightarrow \bigoplus^{l-1}_{s=0} P^s ([\xi],[\phi]) = W_{\xi}/W^l_{\xi},
$$
where $\alpha^{t,s}(\xi,\phi) : P^s ([\xi],[\phi]) \longrightarrow P^t ([\xi],[\phi])$ is obtained by restricting $\alpha^{(2)}_{\xi} (\phi,\bullet)$ to
$P^s ([\xi],[\phi])$ and then projecting the result to $P^t ([\xi],[\phi])$.
 
We have the trace map
$$
tr([\xi],[\phi]): A_{[\xi]} ([\phi]) \longrightarrow \CC
$$
defined by the formula
$$
tr([\xi],[\phi])(\alpha^{(2)}_{\xi} (\phi,\bullet))=\sum^{l-1}_{s=0}tr(\alpha^{s,s}(\xi,\phi)).
$$
Let us recall that once a nonzero lift $\xi$ and $\phi$ over $[\xi]$ and $[\phi]$ are chosen and the value of the trace is fixed all the blocks
$\alpha^{t,s}$ in the above decomposition are uniquely defined. We assume that a value of the trace $c$ is fixed and we denote a unique map in $A_{[\xi]}([\phi])(\xi,\phi) $ with this value of the trace by  $\alpha_{\xi,c} (\phi,\bullet)$. In addition, the off diagonal blocks
$\alpha^{t,s}(\xi,\phi)$ are independent of the trace value, while the diagonal ones change by multiples of identity, see \eqref{tr-v-change} for the precise formula.

Using the inclusions $W\subset W_{\xi}$ in \eqref{WinWxi} gives the linear maps
$$
\alpha_{\xi} (\phi,\bullet): W \longrightarrow \HKC 
$$
depending linearly on $\xi \in \Xi_W$ and $\phi \in W$. We organize all this in  the linear map
\begin{equation}\label{D_W}
D_W: W \longrightarrow (\Xi_W)^{\ast} \otimes Hom(W,\HKC),
\end{equation}
defined by the rule 
$$
D_W (\phi) (\xi)=\alpha_{\xi}(\phi,\bullet)|_W: W\longrightarrow \HKC.
$$
A more conceptual way to understand the quadratic equations defining ${\bf H^{1,0}}(PG_{l})$ is to observe that the map $D_W$ is a sort of a {\it Higgs} field. Namely, for a nonzero $\phi \in W$ fixed, we have
\begin{equation}\label{D_W(phi)}
D_W (\phi): \Xi_W \longrightarrow Hom(W,\HKC).
\end{equation}
By `Higgs' we mean that the operators $D_W (\phi)(\xi)$ and $D_W (\phi)(\xi')$, for any pair of distinct $\xi$ and $\xi'$ in $\Xi_W$, are subject to certain commutation relations. We will now explain this.

We begin by constructing the filtration of $W$ in a way similar to the one defined for $W_{\xi}$. For this we rewrite $D_W (\phi)$ in \eqref{D_W(phi)}
in the form
$$
W \longrightarrow (\Xi_W)^{\ast}\otimes \HKC.
$$
We continue to denote it by $D_W (\phi)$. Composing with the projection
$$
\HKC \longrightarrow \HKC/W
$$
gives the map
$$
D'_W(\phi):W \longrightarrow (\Xi_W)^{\ast}\otimes \HKC/W.
$$
We start defining our filtration by
$$
W^0(\phi):=W, \,\,\, W^1(\phi):=ker(D'_W(\phi)).
$$
From this it follows
$$
D_W(\phi): W^1(\phi) \longrightarrow (\Xi_W)^{\ast}\otimes W^0(\phi).
$$
Composing with the quotient
$$
W^0(\phi) \longrightarrow W^0(\phi)/W^1(\phi)
$$
gives the map 
$$
D''_W(\phi): W^1(\phi) \longrightarrow (\Xi_W)^{\ast}\otimes W^0(\phi)/W^1(\phi)
$$
and we define
$$
W^2(\phi):=ker(D''_W(\phi)).
$$
Again we obtain
$$
D_W(\phi): W^2(\phi) \longrightarrow (\Xi_W)^{\ast}\otimes W^1(\phi).
$$
Continuing in this fashion gives a decreasing filtration
\begin{equation}\label{W-filt}
\HKC= W^{-1}(\phi)\supset W^0(\phi)\supset W^1(\phi) \supset W^2(\phi)\supset \cdots \supset W^i(\phi) \supset W^{i+1}(\phi) \supset \cdots
\end{equation}

This is subject to the property
$$
D_W(\phi):W^i(\phi) \longrightarrow (\Xi_W)^{\ast}\otimes W^{i-1}(\phi),
$$
for all $i\geq 1$.
\begin{rem}\label{rem:Wi(phi)=inter}
	From the definitions it is easy to see that $W^i(\phi)$ is the intersection of $W^i_{\xi}([\phi])$ as $[\xi]$ runs through $\PP(\Xi_W)$
	$$
	W^i(\phi)= \bigcap_{[\xi]\in \PP(\Xi_W)}W^i_{\xi}([\phi]),
	$$
	for every $i$. In particular, we have the inclusion
	$$
	\CC\phi \subset W^i(\phi)
	$$
	for all $i$.
\end{rem}

 We can now state the Higgs property of the map
 $$
 D_W (\phi): W \longrightarrow (\Xi_W)^{\ast}\otimes \HKC.
 $$
\begin{pro}\label{pro:Higgs}
Assume that the linear subsystem $|W|$ corresponding to the subspace $W$ is base point free. Then for every $i\geq 1$, the composition
	$$
\mbox{$	W^i(\phi) \longrightarrow (\Xi_W)^{\ast}\otimes W^{i-1}(\phi) \longrightarrow \bigwedge^2 (\Xi_W)^{\ast}\otimes W^{i-2}(\phi)/\CC\phi $}
	$$
	is identically zero. Explicitly, for two distinct elements $\xi$ and $\xi'$ in $\Xi_W$ we have the following commutation relation
	\begin{equation}\label{Higgs-comutrel}
	\alpha_{\xi} (\phi, \alpha_{\xi'} (\phi, \psi))\equiv \alpha_{\xi'} (\phi, \alpha_{\xi} (\phi, \psi)) \mod \CC\phi,
\end{equation}
	for all $\psi \in W^i(\phi)$ and all $i\geq 1$.
\end{pro}
\begin{pf}
	The commutation relation \eqref{Higgs-comutrel} of the proposition can be deduced from the Koszul cocycle relation \eqref{aplha2-eq}
	$$
	\phi''\alpha_{\xi} (\phi,\phi')-\phi'\alpha_{\xi} (\phi,\phi'')+\phi\alpha_{\xi} (\phi',\phi'')=0, \,\,\forall \phi,\phi',\phi'' \in W_{\xi}.
	$$
	We fix $\phi$ as in the proposition and let $\phi'=\alpha_{\xi'} (\phi,\psi)$, for $\psi \in W^i(\phi)$. Substituting into the equation above and restricting to $Z_{\phi}=(\phi=0)$ gives the equation
	$$
	\left(	\phi''\alpha_{\xi} (\phi,\alpha_{\xi'} (\phi,\psi))-\alpha_{\xi'} (\phi,\psi)\alpha_{\xi} (\phi,\phi'') \right) \big|_{Z_{\phi}}=0,\,\,\forall \phi''\in W_{\xi}.
	$$
	Exchanging the roles of $\xi$ and $\xi'$ we obtain
	$$
	\left(	\phi''\alpha_{\xi'} (\phi,\alpha_{\xi} (\phi,\psi))-\alpha_{\xi} (\phi,\psi)\alpha_{\xi'} (\phi,\phi'') \right) \big|_{Z_{\phi}}=0,\,\,\forall \phi''\in W_{\xi'}.
	$$
	Both equations are valid on the intersection of $W_{\xi}$ and $W_{\xi'}$. Subtracting the two equations gives
	$$
	\phi''\left(	\alpha_{\xi} (\phi,\alpha_{\xi'} (\phi,\psi))-\alpha_{\xi'} (\phi,\alpha_{\xi} (\phi,\psi)) \right) \big|_{Z_{\phi}}=0,\,\,\forall \phi''\in W_{\xi} \bigcap W_{\xi'}.
	$$
	From the inclusion
	$$
	W\subset  W_{\xi} \bigcap W_{\xi'}
	$$
	and the assumption that the linear subsystem $|W|$ corresponding to $W$ is base point free, we can choose $\phi'' \in W$ whose divisor of zeros is disjoint from $Z_{\phi}$. Hence the last equation implies
	$$
	\left(	\alpha_{\xi} (\phi,\alpha_{\xi'} (\phi,\psi))-\alpha_{\xi'} (\phi,\alpha_{\xi} (\phi,\psi)) \right) \big|_{Z_{\phi}}=0.
	$$
	Equivalently, this means
	 $$
	 \left(	\alpha_{\xi} (\phi,\alpha_{\xi'} (\phi,\psi))-\alpha_{\xi'} (\phi,\alpha_{\xi} (\phi,\psi)) \right) \in \CC\phi, 
	 $$
	 and this holds for all $\psi \in W^i$ and $i\geq 1$.
	\end{pf}

Let $l_W$ be the length of the filtration in \eqref{W-filt}, that is we have
$$
\HKC=W^{-1}(\phi) \supset W=W^{0}(\phi) \supset W^{1}(\phi) \supset \cdots W^{l_W}(\phi) \supset W^{l_W+1}(\phi)=0.
$$    
As in the case of the filtration of $W_{\xi}$ we use the Hodge metric on $\HKC$ to obtain the orthogonal decomposition 
\begin{equation}\label{Wphi-ortdecomp}
W=\bigoplus^{l_W}_{s=0} P^s (\phi)
\end{equation}
such that 
$$
W^i(\phi)=\bigoplus^{l_W}_{s=i} P^s (\phi).
$$
By definition of the filtration the maps 
$$
\alpha_{\xi}(\phi,\bullet): W \longrightarrow \HKC,
$$
for every $\xi\in \Xi_W$, have the property of sending $W^i(\phi)$ to $W^{i-1}(\phi)$:
\begin{equation}\label{Grifftrans}
 \alpha_{\xi}(\phi,\bullet):W^i (\phi)\longrightarrow W^{i-1}(\phi), \forall i\geq 0.
\end{equation}
 In addition, the last nonzero term $W^{l_W}$ of the filtration is $\alpha_{\xi}(\phi,\bullet)$-invariant for all $\xi \in \Xi_W$. So  we restrict to the complementary component
 $$
 W'=\bigoplus^{l_W-1}_{s=0} P^s (\phi)
 $$
 and from now on we consider the restriction of $\alpha_{\xi}(\phi,\bullet)$ to this part only.
 
 We decompose $\alpha_{\xi}(\phi,\bullet)$ into blocks according to the orthogonal decomposition of $W$ in \eqref{Wphi-ortdecomp}
 $$
 \alpha_{\xi}(\phi,\bullet)=\sum \alpha^{t,s}_{\xi},
 $$
where the block $\alpha^{t,s}_{\xi}$ maps the summand $P^s(\phi)$ to $P^t(\phi)$: 
$$
\alpha^{t,s}_{\xi}: P^s(\phi) \longrightarrow P^t(\phi).
$$

We define the degree of a block $\alpha^{t,s}_{\xi}$ as the difference $(t-s)$. From \eqref{Grifftrans} it follows that all blocks have degree at least $(-1)$. Grouping together the block of the same degree we obtain the decomposition of $\alpha_{\xi}(\phi,\bullet)$ into the sum of maps
$$
\alpha_{\xi}(\phi,\bullet)=D^{-1}_{\xi} + D^0_{\xi} + D^{1}_{\xi}+\cdots +D^{l_W}_{\xi}
$$
where $D^i_{\xi}$ is the sum of blocks of degree $i$:
$$
D^i_{\xi}=\sum_{t-s=i} \alpha^{t,s}_{\xi}.
$$
We call a summands $D^i_{\xi}$ the homogeneous component of degree $i$. Let us write down the first tree homogeneous components:
\begin{equation}\label{D-101}
	D^{-1}_{\xi}=\sum^{l_W-1}_{s=1} \alpha^{s-1,s}_{\xi},\hspace{0.2cm}	D^{0}_{\xi}=\sum^{l_W-1}_{s=0} \alpha^{s,s}_{\xi},\hspace{0.2cm}
		D^{1}_{\xi}=\sum^{l_W-1}_{s=0} \alpha^{s+1,s}_{\xi}.
\end{equation}

The Higgs condition of Proposition \ref{pro:Higgs} translates into commutativity equations between  the homogeneous components.
\begin{pro}\label{pro:Higgs-2-10}
	Assume the length of the filtration $l_W$ is at least $3$. Then
	for every $\xi$ and $\xi'$ in $\Xi_W$ we have the following relations.
	$$
	D^{-1}_{\xi} D^{-1}_{\xi'}=D^{-1}_{\xi'}D^{-1}_{\xi},
	$$
	$$
	D^{-1}_{\xi} D^{0}_{\xi'}+ D^{0}_{\xi} D^{-1}_{\xi'}=	D^{-1}_{\xi'} D^{0}_{\xi}+ D^{0}_{\xi'} D^{-1}_{\xi},
	$$
	$$
		D^{-1}_{\xi} D^{1}_{\xi'}+ D^{0}_{\xi} D^{0}_{\xi'} +D^{1}_{\xi} D^{-1}_{\xi'}=	D^{-1}_{\xi'} D^{1}_{\xi}+ D^{0}_{\xi'} D^{0}_{\xi'}+D^{1}_{\xi} D^{-1}_{\xi} .
	$$
\end{pro}
\begin{pf}
	To simplify the notation we write $\alpha_{\eta}$ instead of
	$\alpha_{\eta}(\phi,\bullet)$.
		The Higgs condition in Proposition \ref{pro:Higgs} says that
	$
	\alpha_{\xi}
	$
	and $ \alpha_{\xi'}$ are commuting as long we take the values modulo $\CC\phi$. That line lies in $W^{l_W}$. So as long as the values are considered in $W'$, the orthogonal complement of $W^{l_W}$, we have the commutativity
	$$
	\alpha_{\xi} \alpha_{\xi'}=\alpha_{\xi'}\alpha_{\xi}.
	$$
	Taking the homogeneous decomposition of factors on each side we have
	$$
	 \left(\sum^{l_W}_{i=-1} D^i_{\xi} \right)\left(\sum^{l_W-1}_{j=-1} D^j_{\xi'} \right) =\left(\sum^{l_W}_{j=-1} D^j_{\xi'} \right)\left(\sum^{l_W-1}_{i=-1} D^i_{\xi} \right) \mod W^{l_W} 
	 $$
Expanding each side as the sum of homogeneous components we obtain the relations. The ones asserted in the proposition correspond to the components of degrees $-2$, $-1$, $0$.   	 
\end{pf}
\begin{rem}\label{rem:Higgsblocks}
	Substituting the expressions from \eqref{D-101}, the equations in Proposition \ref{pro:Higgs-2-10} are translated into equations between the blocks. In degree $(-2)$ we obtain:
	$$
	\alpha^{s-2,s-1}_{\xi}	\alpha^{s-1,s}_{\xi'}=\alpha^{s-2,s-1}_{\xi'}	\alpha^{s-1,s}_{\xi},\,\forall s\geq 2.
	$$
	In degree $(-1)$:
	$$
		\alpha^{s-1,s}_{\xi}	\alpha^{s,s}_{\xi'}+\alpha^{s-1,s-1}_{\xi}	\alpha^{s-1,s}_{\xi'} =\alpha^{s-1,s}_{\xi'}	\alpha^{s,s}_{\xi}+\alpha^{s-1,s-1}_{\xi'}	\alpha^{s-1,s}_{\xi},\,\forall s\geq 1.
	$$
	In degree $0$:
	$$
	\alpha^{s,s+1}_{\xi}\alpha^{s+1,s}_{\xi'}+	\alpha^{s,s}_{\xi}\alpha^{s,s}_{\xi'}+\alpha^{s,s-1}_{\xi}	\alpha^{s-1,s}_{\xi'} =\alpha^{s,s+1}_{\xi'}\alpha^{s+1,s}_{\xi}+	\alpha^{s,s}_{\xi'}\alpha^{s,s}_{\xi}+\alpha^{s,s-1}_{\xi'}	\alpha^{s-1,s}_{\xi} ,\,\forall s\geq 1.
	$$
\end{rem}
We will now deform the operator $D_W (\phi)$ as follows. Set 
$$
t,\,\, {\bf x}=(x_1,\ldots, x_{l_W-1}),\,\,{\bf y}=(y_1,\ldots, y_{l_W-1}),
$$
be deformation parameters, where $t$, $x_i$, $y_i$ take values in $\CC$, and consider the following perturbation of $D_W (\phi)$:
\begin{equation}\label{Ddeform}
	D_W(\phi; {\bf x}, {\bf y}, t) (\xi)=\alpha_{\xi} ({\bf x}, {\bf y}, t)=
	\sum^{l_W-1}_{s=1}x_s \alpha^{s-1,s}_{\xi} +t\sum^{l_W-1}_{s=0}\alpha^{s,s}_{\xi} + \sum^{l_W-2}_{s=0}y_{s+1} \alpha^{s+1,s}_{\xi} +\sum_{j\geq 2} D^j_{\xi},
\end{equation}
that is only the homogeneous components of degree $(-1)$, $0$ and $1$ are perturbed:
\begin{equation}\label{-101}
	\begin{gathered}
	D^{-1}(\phi; {\bf x}, {\bf y}, t) (\xi)=\sum^{l_W-1}_{s=1}x_s \alpha^{s-1,s}_{\xi},
	\\
		D^0(\phi; {\bf x}, {\bf y}, t) (\xi)=t\sum^{l_W-1}_{s=0}\alpha^{s,s}_{\xi},
		\\
	D^1(\phi; {\bf x}, {\bf y}, t) (\xi)=	\sum^{l_W-2}_{s=0}y_{s+1} \alpha^{s+1,s}_{\xi}.
	\end{gathered}
\end{equation}

 We want to find conditions for our deformation to remain Higgs up to the homogeneous degree zero.
\begin{pro}\label{pro:deformHiggs}
	The map
	$$
 	D_W(\phi; {\bf x}, {\bf y}, t): W \longrightarrow (\Xi_W)^{\ast}\otimes \HKC
 	$$
 	is Higgs in degrees up to zero, if the deformation parameters ${\bf x}, {\bf y}, t$ satisfy the equations
 	$$
 	x_s y_s=t^2,
 	$$
 	for all $s\in [1,l_W-1]$.
\end{pro}
\begin{pf}
	Write the decomposition of 	$D_W(\phi; {\bf x}, {\bf y}, t) (\xi)$ into the sum of homogeneous components
	$$
	 D_W(\phi; {\bf x}, {\bf y}, t) (\xi)= \sum^{l_W-1}_{i=-1} D^i(\phi; {\bf x}, {\bf y}, t) (\xi).
	 $$
	 We seek conditions on the deformation parameters ${\bf x}, {\bf y}, t$ for the Higgs condition
	 $$
	 D_W(\phi; {\bf x}, {\bf y}, t) (\xi) D_W(\phi; {\bf x}, {\bf y}, t) (\xi')=D_W(\phi; {\bf x}, {\bf y}, t) (\xi')D_W(\phi; {\bf x}, {\bf y}, t) (\xi)
	 $$
	 to hold in degrees $(-2)$, $(-1)$ and $0$, for all $\xi,\xi'\in \Xi_W$. Using the decomposition into homogeneous components we expand the product on each side into the homogeneous pieces. In particular, the first three degrees, -2, -1, 0, give the equations analogous to the ones in Proposition \ref{pro:Higgs-2-10}:
	 	$$
	 	\begin{gathered}
\text{degree $-2$:}\\
	D^{-1}(\phi; {\bf x}, {\bf y}, t) (\xi) D^{-1}(\phi; {\bf x}, {\bf y}, t) (\xi')=	D^{-1}(\phi; {\bf x}, {\bf y}, t) (\xi') D^{-1}(\phi; {\bf x}, {\bf y}, t)(\xi),
\end{gathered}
	 $$
	 $$
	 \begin{gathered}
\text{degree $-1$:}\\
		D^{-1}(\phi; {\bf x}, {\bf y}, t) (\xi) 	D^{0}(\phi; {\bf x}, {\bf y}, t) (\xi') + D^{0}(\phi; {\bf x}, {\bf y}, t) (\xi) 	D^{-1}(\phi; {\bf x}, {\bf y}, t) (\xi')
		\\
		=	D^{-1}(\phi; {\bf x}, {\bf y}, t) (\xi') 	D^{0}(\phi; {\bf x}, {\bf y}, t) (\xi) + D^{0}(\phi; {\bf x}, {\bf y}, t) (\xi') 	D^{-1}(\phi; {\bf x}, {\bf y}, t) (\xi),
			\end{gathered}
	 $$
	 $$
	 \begin{gathered}
	 \text{degree $0$:}\\
	 	D^{-1}(\phi; {\bf x}, {\bf y}, t) (\xi) 	D^{1}(\phi; {\bf x}, {\bf y}, t) (\xi') + D^{0}(\phi; {\bf x}, {\bf y}, t) (\xi) 	D^{0}(\phi; {\bf x}, {\bf y}, t) (\xi')
	 	\\
	 	+D^{1}(\phi; {\bf x}, {\bf y}, t) (\xi) 	D^{-1}(\phi; {\bf x}, {\bf y}, t) (\xi')
	 \\
	= D^{-1}(\phi; {\bf x}, {\bf y}, t) (\xi') 	D^{1}(\phi; {\bf x}, {\bf y}, t) (\xi) + D^{0}(\phi; {\bf x}, {\bf y}, t) (\xi') 	D^{0}(\phi; {\bf x}, {\bf y}, t) (\xi) 
	\\
	+ D^{1}(\phi; {\bf x}, {\bf y}, t) (\xi') 	D^{-1}(\phi; {\bf x}, {\bf y}, t) (\xi).
 \end{gathered}
	 $$
Substituting the expressions from \eqref{-101} will impose the conditions on the parameters. In degree $(-2)$ we obtain
$$
\left(\sum^{l_W-1}_{s=1}x_s \alpha^{s-1,s}_{\xi}\right) \left( \sum^{l_W-1}_{s=1}x_s \alpha^{s-1,s}_{\xi'}\right)=\left(\sum^{l_W-1}_{s=1}x_s \alpha^{s-1,s}_{\xi'}\right) \left( \sum^{l_W-1}_{s=1}x_s \alpha^{s-1,s}_{\xi}\right).
$$
Multiplying out on each side we obtain
$$
\sum_{s\geq 2}x_{s-1}x_s \alpha^{s-2,s-1}_{\xi}\alpha^{s-1,s}_{\xi'}=
\sum_{s\geq 2}x_{s-1}x_s \alpha^{s-2,s-1}_{\xi'}\alpha^{s-1,s}_{\xi}
$$
In view of the first equation in Remark \ref{rem:Higgsblocks} this imposes no condition on ${\bf x}$.

In degree $(-1)$:
$$
\begin{gathered}
\left(\sum^{l_W-1}_{s=1}x_s \alpha^{s-1,s}_{\xi}\right) \left( \sum^{l_W-1}_{s=1}t \alpha^{s,s}_{\xi'}\right) +\left( \sum^{l_W-1}_{s=1}t \alpha^{s,s}_{\xi}\right)\left(\sum^{l_W-1}_{s=1}x_s \alpha^{s-1,s}_{\xi'}\right)
\\
 =\left(\sum^{l_W-1}_{s=1}x_s \alpha^{s-1,s}_{\xi'}\right) \left( \sum^{l_W-1}_{s=1}t \alpha^{s,s}_{\xi}\right) +\left( \sum^{l_W-1}_{s=1}t \alpha^{s,s}_{\xi'}\right)\left(\sum^{l_W-1}_{s=1}x_s \alpha^{s-1,s}_{\xi}\right)
\end{gathered}
$$
Multiplying out on each side gives
$$
\sum_s tx_s \left( \alpha^{s-1,s}_{\xi}\alpha^{s,s}_{\xi'} +  \alpha^{s-1,s-1}_{\xi}\alpha^{s,s}_{\xi'} \right)
 =\sum_s tx_s \left( \alpha^{s-1,s}_{\xi'}\alpha^{s,s}_{\xi} +  \alpha^{s-1,s-1}_{\xi'}\alpha^{s,s}_{\xi} \right).
$$
In view of the second equation in Remark \ref{rem:Higgsblocks}, the equation is valid for all values of $t$ and $\bf x$.

We now turn to the degree $0$. The left hand side of the equation gives
$$
\begin{gathered}
\left(\sum_{s}x_s \alpha^{s-1,s}_{\xi'}\right)\left(\sum_{s}y_{s+1} \alpha^{s+1,s}_{\xi'}\right) +t^2\sum_{s}  \alpha^{s,s}_{\xi}\alpha^{s,s}_{\xi'} +\left(\sum_{s}y_{s+1} \alpha^{s+1,s}_{\xi}\right)\left(\sum_{s}x_s \alpha^{s-1,s}_{\xi'}\right)
\\
=\sum_s x_s y_{s}\alpha^{s-1,s}_{\xi} \alpha^{s,s-1}_{\xi'}+ t^2\sum_{s}  \alpha^{s,s}_{\xi}\alpha^{s,s}_{\xi'} +
\sum_s x_s y_{s}\alpha^{s,s-1}_{\xi} \alpha^{s-1,s}_{\xi'}
\end{gathered}
$$
The right hand side has the same form with $\xi$ and $\xi'$ interchanged:
$$
\sum_s x_s y_{s}\alpha^{s-1,s}_{\xi'} \alpha^{s,s-1}_{\xi}+ t^2\sum_{s}  \alpha^{s,s}_{\xi'}\alpha^{s,s}_{\xi} +
\sum_s x_s y_{s}\alpha^{s,s-1}_{\xi'} \alpha^{s-1,s}_{\xi}.
$$
Taking the difference of the last two expressions we obtain
$$
\begin{gathered}
\sum_s x_s y_{s} \left(\alpha^{s-1,s}_{\xi} \alpha^{s,s-1}_{\xi'} -\alpha^{s-1,s}_{\xi'} \alpha^{s,s-1}_{\xi}\right)
\\
+ t^2\sum_{s}  \left(\alpha^{s,s}_{\xi}\alpha^{s,s}_{\xi'} -\alpha^{s,s}_{\xi'}\alpha^{s,s}_{\xi} \right) +
\sum_s x_s y_{s}\left( \alpha^{s,s-1}_{\xi} \alpha^{s-1,s}_{\xi'}-\alpha^{s,s-1}_{\xi'} \alpha^{s-1,s}_{\xi}\right).
\end{gathered}
$$
 
The third equation in Remark \ref{rem:Higgsblocks} tells us that
$$
\alpha^{s,s}_{\xi}\alpha^{s,s}_{\xi'} -\alpha^{s,s}_{\xi'}\alpha^{s,s}_{\xi} =\alpha^{s,s+1}_{\xi'}\alpha^{s+1,s}_{\xi} -\alpha^{s,s+1}_{\xi}\alpha^{s+1,s}_{\xi'}	+ \alpha^{s,s-1}_{\xi'}\alpha^{s+1,s}_{\xi} -\alpha^{s,s-1}_{\xi}\alpha^{s-1,s}_{\xi'}.
$$
Substituting this for the terms of the sum in the middle of the difference computed above gives
$$
\sum_s (x_s y_{s}-t^2) \left(\alpha^{s-1,s}_{\xi} \alpha^{s,s-1}_{\xi'} -\alpha^{s-1,s}_{\xi'} \alpha^{s,s-1}_{\xi}\right) +
\sum_s (x_s y_{s}-t^2) \left( \alpha^{s,s-1}_{\xi} \alpha^{s-1,s}_{\xi'}-\alpha^{s,s-1}_{\xi'} \alpha^{s-1,s}_{\xi}\right).
$$
Hence the condition
$$
x_sy_s =t^2, \forall s,
$$
implies that the Higgs condition holds in degree $0$. This completes the proof of the proposition.   
\end{pf}

The map 
$$
D_W (\phi): W \longrightarrow \Xi^{\ast}_W \otimes \HKC
$$
is reminiscent of a holomorphic one-form with values in the vector space
of linear maps $Hom(W,\HKC)$, a sort of nonabelian analogue of a one-form $\phi$ on $C$. This is the reason for notation ${\bf H^{1,0}}(PG_{l})$ and the terminology `nonabelian Dolbeault variety'.

The pair $({\bf H^{1,0}}(PG_{l}),H_0)$ could be thought of as an analogue of the Jacobian of a curve and its theta divisor. Observe that $({\bf H^{1,0}}(PG_{l}),H_0)$ depends on the graph $PG_{l}$ only and completely forgets our curve $C$. To remedy this we relate $({\bf H^{1,0}}(PG_{l}),H_0)$ to the incidence variety $\PP({\cal W}_{\Sigma^0_r})$. 

\subsection{Relating  $({\bf H^{1,0}}(PG_{l}),H_0)$ to $\PP({\cal W}_{\Sigma^0_r})$}

At every point of $([\xi],[\phi])$ in $\PP({\cal W}_{\Sigma^0_r})$ we have the filtrations
on $W_{\xi}={\cal W}_{\Sigma^0_r,[\xi]}$, the fibre of ${\cal W}_{\Sigma^0_r}$ at $[\xi] \in \Sigma^0_r$. Those vary as $[\phi]$ varies in the projective space
$\PP(W_{\xi})$. The length of the filtration denoted $l_{\xi}([\phi])$ gives the integer valued function
$$
{\mathfrak{l}}:\PP({\cal W}_{\Sigma^0_r}) \longrightarrow \ZZ_{\geq 0}.
$$
It partitions $\PP({\cal W}_{\Sigma^0_r})$ into strata where ${\mathfrak{l}}$ is constant. Fix the value $l$ of that function which we assume to be at least three and consider the corresponding stratum
\begin{equation}\label{l-stratum}
{\mathfrak{L}}_l:=	{\mathfrak{l}}^{-1}(l)
\end{equation}
parametrizing the points $([\xi],[\phi]) \in p^{-1}_1 (\Sigma^0_r)$ with $([\xi],[\phi])$-filtration having length $l$. We have learned that it comes equipped with the graph $PG_l$; this in turn defines the pair
$
 ({\bf H^{1,0}}(PG_{l}),H_0),
 $
  the nonabelian Dolbeault variety ${\bf H^{1,0}}(PG_{l})$ and its toric divisor $H_0$. We relate  ${\mathfrak{L}}_l$ to ${\bf H^{1,0}}(PG_{l})$ via the diagrams of maps
  \begin{equation}\label{AJl-1}
  	\xymatrix{
  		& \OO^{\times}_{{\mathfrak{L}}_l} (-1)\times \CC \ar[ld]_{\pi_{\mathfrak{L}_l}\times id_{\CC}} \ar[dr]^{\tau^0_l (c)}&\\
  		{\mathfrak{L}}_l \times \CC && {\bf H^{1,0}}(PG_{l})\setminus H_0
  	}
  \end{equation}
In this diagram $\OO^{\times}_{{\mathfrak{L}}_l} (-1)$ is the total space of the tautological line bundle $\OO_{P} (-1)$ minus its zero section and restricted to the stratum ${\mathfrak{L}}_l$ and the map $\pi_{\mathfrak{L}_l}$ is the structure projection
$$
\pi_{\mathfrak{L}_l}:\OO^{\times}_{{\mathfrak{L}}_l} (-1) \longrightarrow \mathfrak{L}_l ;
	$$
	 the map $\tau^0_l (c)$ on the right side of the diagram is the one we construct below; in the notation $c$ is the trace parameter and the extra factor $\CC$ in the diagram refers to spectral parameter which will be a part of the construction. 
\begin{rem}\label{rem:loops}
	Conceptually, the maps $\tau^0_l (c)$ relate ${\mathfrak{L}}_l \times \CC$ to the space $L({\bf H^{1,0}}(PG_{l})\setminus H_0)$ of free loops in ${\bf H^{1,0}}(PG_{l})\setminus H_0$. Indeed,  
	we can replace $\OO^{\times}_{{\mathfrak{L}}_l} (-1)$ by the associated circle bundle
	$$
	\pi^{S^1}_{\mathfrak{L}_l}:S^1\left(\OO^{\times}_{{\mathfrak{L}}_l} (-1)\right) \longrightarrow \mathfrak{L}_l, 
	$$
	where $S^1\left(\OO^{\times}_{{\mathfrak{L}}_l} (-1)\right)$ is the fibre subspace of $\OO^{\times}_{{\mathfrak{L}}_l} (-1)$ formed by the unit length vectors with respect to the natural (product of) Fubini-Study metrics on $\OO_P (-1)$; $\pi^{S^1}_{\mathfrak{L}_l}$ is the
	corresponding fibre projection, the restriction of $\pi_{\mathfrak{L}_l}$ to $S^1\left(\OO^{\times}_{{\mathfrak{L}}_l} (-1)\right)$.
	Then the map $\tau^0_l (c)$ associates a free loop in ${\bf H^{1,0}}(PG_{l})\setminus H_0$ with every point of ${\mathfrak{L}}_l \times \CC$:
	 $$
	{\mathfrak{L}}_l \times \CC \ni ([\xi],[\phi];u) \mapsto \left[\tau^0_l (c):
	(\pi^{S^1}_{\mathfrak{L}_l})^{-1}([\xi],[\phi])\times \{u\} \longrightarrow {\bf H^{1,0}}(PG_{l})\setminus H_0 \right].
	$$
	This connects to the ideas of string topology of Chas and Sullivan, \cite{Ch-Su}, \cite{CoVo} for overview. Exploring those connections should lead to new features of IVHS refinement. This will be investigated elsewhere.   
\end{rem}
 
The construction proceeds in several stages.
Recall that for a fixed value $c\in \CC$ we have assigned representations of the quiver $PG_l$:
$$
\rho_c:\OO^{\times}_{{\mathfrak{L}}_l} (-1) \longrightarrow {\mathfrak{Reps}}(PG_l).
$$
Namely, for $\xi \otimes \phi$ in $\OO^{\times}_{{\mathfrak{L}}_l} (-1)$ lying over $([\xi], [\phi])$ in ${\mathfrak{L}}_l$ we have the representation of $PG_l$
$$
\rho_c(\xi,\phi)=\{\{P^s([\xi],[\phi])\}_s; \, \alpha^{t,s}_c(\xi,\phi): P^s([\xi],[\phi]) \rightarrow P^t([\xi],[\phi])\}, 
$$
which labels every pair of white-black vertices $s$ and $s'$ of $PG_l$ with the vector spaces
$P^s([\xi],[\phi])$ of the decomposition
$$
W_{\xi}/W^l_{\xi}([\phi])=\bigoplus^{l-1}_{s=0} P^s([\xi],[\phi])
$$ 
and each edge $(s) \to (t)'$ of the graph is decorated with the linear map $\alpha^{t,s}_c(\xi,\phi)$, the $(t,s)$-blocks of the map  $\alpha^{(2)}_{\xi} (\phi,\bullet)$ in the collection $A_{[\xi]}([\phi])(\xi,\phi)$ having the trace $c$, see Lemma \ref{lem:traceparam}. With this representation at hand we pass to the second stage of the construction.

For each diagonal block $\alpha^{s,s}_c(\xi,\phi)$ of the representation $\rho_c(\xi,\phi)$ let
\begin{equation}\label{c_r}
	ch_s (c,\xi,\phi;u):=det(\alpha^{s,s}_c (\xi,\phi) -u{\bf 1}_{P^s([\xi],[\phi])})
\end{equation}  
be the value of the characteristic polynomial of $\alpha^{s,s}_c (\xi,\phi)$ at $u\in \CC$. Thus we obtain functions
$$
ch_s(c) : \OO^{\times}_{{\mathfrak{L}}_l} (-1) \times \CC \longrightarrow \CC
$$
defined by the formula
$$
(\xi\otimes \phi, u) \mapsto ch_s (c,\xi,\phi;u)=det(\alpha^{s,s}_c (\xi,\phi) -u{\bf 1}_{P^s([\xi],[\phi])}).
$$
For obvious reasons we refer to $u$ as {\it spectral parameter}.
From the construction it follows that the functions are $C^{\infty}$ along ${\mathfrak{L}}_l$ and holomorphic in ${\CC}$ direction as well as in the fibre direction of the structure projection $\pi_{\mathfrak{L}_l}$. 

We use these functions to define coordinates of points in ${\bf H^{1,0}}(PG_{l})$.  
Set
$$
x_s (c,\xi,\phi;u):=\frac{exp(ch_s (c,\xi,\phi;u))}{exp( ch_{s-1} (c,\xi,\phi;u))},\,\,\, y_s(c,\xi,\phi;u):=\frac{exp(ch_{s-1} (c,\xi,\phi;u))}{exp( ch_{s} (c,\xi,\phi;u))},
$$
for every $s\in [1,l-1]$. Remember that the space $V_{2l-1}$, whose projectivization contains ${\bf H^{1,0}}(PG_{l})$, comes with a distinguished basis
$$
\{e_0, (s),(s)' |\, s\in [1,l-1] \}.
$$
With this in mind, define the map $\tau^0_l(c)$ in \eqref{AJl-1} by the formula
\begin{equation}\label{tau0l-formula}
	\tau^0_l (c) (\xi,\phi;u)=[e_0 +\sum^{l-1}_{s=1} \left(x_s(c,\xi,\phi;u) (s) +y_s(c,\xi,\phi;u)(s)'\right)] \in \PP(V_{2l-1}). 
\end{equation}
In homogeneous coordinates we have
$$
\tau^0_l (c)(\xi,\phi;u)=[1:x_1(c;u):y_1(c;u):\cdots : x_{l-1}(c;u):y_{l-1}(c;u)] \in \PP(V_{2l-1}),
$$ 
where we omitted $(\xi,\phi)$ in the notation of coordinates. From the definition of the coordinates $x_r(c;u)$ and $y_r(c;u)$ it follows
$$
x_s (c;u) y_s(c;u)=1, \forall s\in[1,l-1].
$$
Hence $\tau^0_l (c)([\xi],[\phi];u)$ lies in ${\bf H^{1,0}}(PG_{l})$. Furthermore, since the first coordinate 
$$
T(\tau^0_l(c) ([\xi],[\phi]))=1 \neq 0,
$$
we deduce that $\tau^0_l(c) ([\xi],[\phi];u)$ lies in the complement of the divisor $H_0$. Thus we obtain the map
$$
\tau^0_l(c): \OO^{\times}_{{\mathfrak{L}}_l} (-1) \times \CC \longrightarrow
{\bf H^{1,0}}(PG_{l}) \setminus H_0,
$$
for every value $c \in \CC$ of the trace parameter.

Recall that ${\bf H^{1,0}}(PG_{l})$ is a toric variety and the complement of the toric divisor $H_0$ is an algebraic torus
$$
{\bf H^{1,0}}(PG_{l})\setminus H_0 \cong (\CC^{\times})^{l-1}.
$$
Hence the maps $\tau^0_l(c)$ can be multiplied. 
\begin{pro}\label{pro:tau0l-mult}
	Fix a positive integer $n$ and let ${\bf c}=(c_1,\ldots,c_n) \in \CC^n$.
	Then we have the maps
	$$
	(\tau^0_l)^n({\bf c}): \left(\OO^{\times}_{{\mathfrak{L}}_l} (-1) \times \CC\right)^n \longrightarrow {\bf H^{1,0}}(PG_{l})\setminus H_0 \cong (\CC^{\times})^{l-1}
	$$
	defined by the product formula
	$$
	(\tau^0_l)^n({\bf c})=\prod^n_{i=1} \tau^0_l( c_i),
	$$
	where on the right hand side the multiplication of factors is with respect to the operation in $(\CC^{\times})^{l-1}$. 
\end{pro}

 The above gives the map
 \begin{equation}\label{tau0-n}
 	(\tau^0_l)^n: \CC^n \longrightarrow Maps \left(\left(\OO^{\times}_{{\mathfrak{L}}_l} (-1) \times \CC\right)^n, (\CC^{\times})^{l-1}\right)
 \end{equation}
defined by the rule
$$
\CC^n \ni {\bf c} \mapsto (\tau^0_l)^n({\bf c}) \in  Maps \left(\left(\OO^{\times}_{{\mathfrak{L}}_l} (-1) \times \CC\right)^n, (\CC^{\times})^{l-1}\right).
$$
In addition, both sides of \eqref{tau0-n} are equipped with $(\CC^{\times})^n$-actions: on the left a vector $v=(v_1,\ldots,v_n)\in (\CC^{\times})^n$ acts on $\CC^n$ by scaling $i$-th coordinate by $v_i$; on the right the action on the space of maps is defined through the action on the domain, that is, for $f\in  Maps \left(\left(\OO^{\times}_{{\mathfrak{L}}_l} (-1) \times \CC\right)^n, (\CC^{\times})^{l-1}\right)$ we define
$$
(v\cdot f)\big((\xi_1 \otimes \phi_1,u_1),\ldots,(\xi_n \otimes \phi_n,u_n)\big)=
f\big((v_1(\xi_1 \otimes \phi_1),u_1),\ldots,(v_n(\xi_n \otimes \phi_n),u_n)\big),
$$
where $v_i (\xi_i \otimes \phi_i) =(v_i\xi_i)\otimes \phi_i = \xi_i \otimes (v_i \phi_i)$, the natural $\CC^{\times}$-action on the fibres of $\OO^{\times}_{{\mathfrak{L}}_l} (-1)$ over ${\mathfrak{L}}_l$. With those actions on both sides of \eqref{tau0-n} understood, we show that
 the map $(\tau^0_l)^n$ is $(\CC^{\times})^n$-equivariant.
\begin{pro}\label{pro:tau0-n-equivariant}
The map $(\tau^0_l)^n$ in \eqref{tau0-n} is $(\CC^{\times})^n$-equivariant:
		for every ${\bf v} \in (\CC^{\times})^n$ the diagram
		$$
		\xymatrix@C=55pt{
		\CC^n \ar^(.3){(\tau^0_l)^n}[r] \ar_{{\bf v}}[d]&
	Maps(\left(\OO^{\times}_{{\mathfrak{L}}_l} (-1) \times \CC\right)^n, (\CC^{\times})^{l-1}) \ar^{\bf v}[d]\\
	\CC^n \ar^(.3){(\tau^0_l)^n}[r]&
	Maps(\left(\OO^{\times}_{{\mathfrak{L}}_l} (-1) \times \CC\right)^n, (\CC^{\times})^{l-1})
}
$$
commutes; the vertical arrows of the diagram correspond to the actions of 
$v\in (\CC^{\times})^{n}$ described in the paragraph preceding the proposition.  
\end{pro}
 \begin{pf}
 	It is enough to see how the multiplication of the trace parameter $c$
 	by a scalar $v\in \CC^{\times}$ affects the map $\tau^0_l$. For this
 	we write the expressions which determine the coordinates of $\tau^0_l(vc)$:
 	$$
 	\begin{gathered}
 ch_s(vc,(\xi \otimes\phi,u)):=ch_s(vc,\xi,\phi;u)=	det(\alpha^{s,s}_{vc}(\xi,\phi)-u{\bf 1}_{P^{s}})=	det(v\alpha^{s,s}_{c}(\xi,\phi)-u{\bf 1}_{P^{s}})\\
 =det(\alpha^{s,s}_{c}(v\xi,\phi)-u{\bf 1}_{P^{s}}) =ch_s(c,v\xi,\phi;u)=ch_s(c,v((\xi\otimes\phi),u))=(v\cdot ch_s )(c,(\xi\otimes\phi,u)).
 \end{gathered}
 $$
 This implies the formula
 $$
 \tau^0_l(vc)(\xi,\phi;u)=(v\cdot\tau^0_l(c))(\xi,\phi;u)
 $$
 intertwining two actions.
 \end{pf}

In addition to the torus action considered above the trace parameters form the additive group $\CC$. This gives an additive action on the trace parameters:
$$
T_{t}:c \mapsto c+t,
$$
the translation by $t\in\CC$. On the space of maps
$Maps \left(\OO^{\times}_{{\mathfrak{L}}_l} (-1) \times \CC, (\CC^{\times})^{l-1}\right)$ we can also define an additive action of $\CC$
via the additive action on the spectral parameter. Namely, on the domain
$\OO^{\times}_{{\mathfrak{L}}_l} (-1) \times \CC$ of the space of maps $\CC$ acts by shifts of the spectral parameter:
$$
S_a: (\xi\otimes \phi, u) \mapsto (\xi\otimes \phi, u+a),\,\,\forall (\xi\otimes \phi, u) \in \OO^{\times}_{{\mathfrak{L}}_l} (-1) \times \CC;
$$
this induces the shift operator on the space of maps
$$
\widehat{S}_a:Maps \left(\OO^{\times}_{{\mathfrak{L}}_l} (-1) \times \CC, (\CC^{\times})^{l-1}\right) \longrightarrow Maps \left(\OO^{\times}_{{\mathfrak{L}}_l} (-1) \times \CC, (\CC^{\times})^{l-1}\right)
$$
defined by the formula
\begin{equation}\label{shift-op}
	\widehat{S}_a(f)=f\circ S_{-a}, \forall f\in Maps \left(\OO^{\times}_{{\mathfrak{L}}_l} (-1) \times \CC, (\CC^{\times})^{l-1}\right).
\end{equation}
Explicitly, for $f\in Maps \left(\OO^{\times}_{{\mathfrak{L}}_l} (-1) \times \CC, (\CC^{\times})^{l-1}\right)$ we have
$$
\widehat{S}_a(f)(\xi\otimes\phi,u)=f(\xi\otimes\phi,u-a),\,\,\forall (\xi\otimes\phi,u) \in \OO^{\times}_{{\mathfrak{L}}_l} (-1) \times \CC.
$$
We wish to understand how the map
$$
\tau^0_l: \CC \longrightarrow Maps \left(\OO^{\times}_{{\mathfrak{L}}_l} (-1) \times \CC, (\CC^{\times})^{l-1}\right)
$$
relates those two actions. This is answered in the following statement.
\begin{pro}\label{pro:trace-spectral}
	Let ${\mathfrak{L}}_l (h^l,\lambda)$ be the substratum of ${\mathfrak{L}}_l$ parametrizing the points $([\xi],[\phi])$ where the dimension of $W^l_{\xi}([\phi])$ is $h^l$ and the dimensions of the remaining graded pieces of $([\xi],[\phi])$-filtrations are determined by the partition $\lambda$, that is, we have
	$$
	dim(P^i([\xi],[\phi]))=dim W^i_{\xi}([\phi])/W^{i+1}_{\xi}([\phi])=\lambda_i,\,\,\forall i\in[0,l-1].
	$$
	Then the restriction $\tau^0_l|_{{\mathfrak{L}}_l (h^l,\lambda)}$ to
	the substratum ${\mathfrak{L}}_l (h^l,\lambda)$
	$$
	\tau^0_l|_{{\mathfrak{L}}_l (h^l,\lambda)}: \CC \longrightarrow
	 Maps \left(\OO^{\times}_{{\mathfrak{L}}_l(h^l,\lambda)} (-1) \times \CC, (\CC^{\times})^{l-1}\right)
	 $$
	 relates the two additive actions via the commutative diagram
	 $$
	 \xymatrix@C=45pt{
	 \CC \ar[r]^(.25){\tau^0_l|_{{\mathfrak{L}}_l (h^l,\lambda)}}\ar[d]_{T_a}& Maps \left(\OO^{\times}_{{\mathfrak{L}}_l(h^l,\lambda)} (-1) \times \CC, (\CC^{\times})^{l-1}\right) \ar[d]^{\widehat{S}_{\frac{a}{|\lambda|}}}\\
\CC \ar[r]^(.25){\tau^0_l|_{{\mathfrak{L}}_l (h^l,\lambda)}}& Maps \left(\OO^{\times}_{{\mathfrak{L}}_l(h^l,\lambda)} (-1) \times \CC, (\CC^{\times})^{l-1}\right)	 
 } 
$$
where $|\lambda|=dim(W_{\xi}/W^l_{\xi}([\phi]))=g-r-h^l$ is the weight of the partition $\lambda$.
\end{pro}

 The above statement follows from the lemma below. It tells us that dependence of $\tau^0_l(c)$ on the trace parameter $c$ intertwines the two actions: changing $c$ to $(c+a)$ is equivalent to an appropriate shift of the spectral parameter, the second factor in the Cartesian product $\OO^{\times}_{{\mathfrak{L}}_l} (-1) \times \CC$. 

\begin{lem}\label{lem:translation}
	A translation of the trace parameter is equivalent to the shifting operation on the side of the spectral parameter $\CC$ in the Cartesian product $\OO^{\times}_{{\mathfrak{L}}_l} (-1) \times \CC$.
	More precisely, set 
	$$
	{w([\xi],[\phi])}:=dim(W^0_{\xi}/W^l_{\xi ([\phi])})=g-r-h^l ([\xi],[\phi]),
	$$
	then the following formula holds
	$$
\tau^0_l(c+a) (\xi,\phi;u) =\tau^0_l(c) (\xi,\phi;u-\frac{a}{w([\xi],[\phi])}), \forall c,a,\in \CC, \,\, \forall (\xi\otimes\phi,u) \in  \OO^{\times}_{{\mathfrak{L}}_l} (-1) \times \CC.
$$	  
\end{lem} 
\begin{pf}
	The formula follows from the observation
	$$
	\alpha^{s,s}_{c+a}(\xi,\phi) -\alpha^{s,s}_c (\xi,\phi) =\frac{a}{w([\xi],[\phi])} {\bf 1}_{P^s}, \,\forall s\in [1,l-1], \forall c,a\in \CC,
	$$
	see \eqref{tr-v-change}; from now on $\xi$ and $\phi$ are fixed and they will be omitted to simplify the notation. The above implies
	$$
	\begin{gathered}
	ch_s (c+a;u)=det(\alpha^{s,s}_{c+a}-u{\bf 1}_{P^s})=
	det(\alpha^{s,s}_c +\frac{a}{w} {\bf 1}_{P^s} -u{\bf 1}_{P^s})
	\\
	=det(\alpha^{s,s}_{c}-(u-\frac{a}{w}){\bf 1}_{P^s})=ch_s (c;u-\frac{a}{w}).
\end{gathered}
	$$
	Hence we obtain
	$$
	x_s(c+a,u) =\frac{exp(ch_s (c+a;u))}{exp(ch_{s-1} (c+a;u))}=\frac{exp(ch_s (c;u-\frac{a}{w}))}{exp(ch_{s-1} (c;u-\frac{a}{w}))}= x_s(c,u-\frac{a}{w}),
	$$
	$$
		y_s(c+a,u) =x^{-1}_s(c+a,u)= x^{-1}_s(c,u-\frac{a}{w})=y_s(c,u-\frac{a}{w}).
	$$
	From this follows the asserted formula
	$$
		\begin{gathered}
	\tau^0_l (c+a) (\xi,\phi;u)=[e_0 +\sum^{l-1}_{s=1} \left(x_s(c+a;u) (s) +y_r(c+a;u)(s)'\right)]
	\\
	= [e_0 +\sum^{l-1}_{s=1} \left(x_s(c;u-\frac{a}{w}) (s) +y_r(c;u-\frac{a}{w})(s)'\right)] =	\tau^0_l (c) (\xi,\phi;u-\frac{a}{w}).
	\end{gathered}
$$	
\end{pf} 

The above considerations can be clearly extended to the $n$-fold version of
$\tau^0_l$, the map 
$$
(\tau^0_l)^n:\CC^n \longrightarrow   Maps \left( \left(\OO^{\times}_{{\mathfrak{L}}_l(h^l,\lambda)} (-1) \times \CC\right)^n, (\CC^{\times})^{l-1}\right)
$$
in \eqref{tau0-n}. The additive actions on both sides  extend in an obvious way to give the following $n$-fold version of Proposition \ref{pro:trace-spectral}.
\begin{cor}\label{cor:trace-spectral-n}
	Let ${\mathfrak{L}}_l (h^l,\lambda)$ be the substratum of ${\mathfrak{L}}_l$ parametrizing the points $([\xi],[\phi])$ where the dimension of $W^l_{\xi}([\phi])$ is $h^l$ and the dimensions of the remaining graded pieces of $([\xi],[\phi])$-filtrations are determined by the partition $\lambda$, that is, we have
	$$
	dim(P^i([\xi],[\phi]))=dim W^i_{\xi}([\phi])/W^{i+1}_{\xi}([\phi])=\lambda_i,\,\,\forall i\in[0,l-1].
	$$
	Then the restriction $(\tau^0_l)^n|_{{\mathfrak{L}}_l (h^l,\lambda)}$ to
	${\mathfrak{L}}_l (h^l,\lambda)$
	$$
	(\tau^0_l)^n|_{{\mathfrak{L}}_l (h^l,\lambda)}: \CC^n \longrightarrow
	Maps \left(\left(\OO^{\times}_{{\mathfrak{L}}_l(h^l,\lambda)} (-1) \times \CC\right)^n, (\CC^{\times})^{l-1}\right)
	$$
	relates the two additive actions via the commutative diagram
	$$
	\xymatrix@C=55pt{
		\CC^n \ar[r]^(.25){(\tau^0_l)^n|_{{\mathfrak{L}}_l (h^l,\lambda)}}\ar[d]_{T_{\bf a}}& Maps \left(\left(\OO^{\times}_{{\mathfrak{L}}_l(h^l,\lambda)} (-1) \times \CC\right)^n, (\CC^{\times})^{l-1}\right) \ar[d]^{\widehat{S}_{\frac{\bf a}{|\lambda|}}}\\
		\CC^n \ar[r]^(.25){(\tau^0_l)^n|_{{\mathfrak{L}}_l (h^l,\lambda)}}& Maps \left(\left(\OO^{\times}_{{\mathfrak{L}}_l(h^l,\lambda)} (-1) \times \CC\right)^n, (\CC^{\times})^{l-1}\right)	 
	} 
	$$
	where ${\bf a} \in \CC^n$ and $T_{\bf a}({\bf c})={\bf c} +{\bf a}$ is the translation by ${\bf a}$ on the side of the space of trace parameters and $\widehat{S}_{\frac{\bf a}{|\lambda|}}$ is the shift operator by $\frac{\bf a}{|\lambda|}$ with $|\lambda|=dim(W_{\xi}/W^l_{\xi}([\phi]))=g-r-h^l$ the weight of the partition $\lambda$.
\end{cor}

\vspace{0.2cm}
The dual version of the above construction relates ${\mathfrak{L}}_l$ to the hyperplane sections of the nonabelian Dolbeault variety ${\bf H^{1,0}}(PG_{l})$. Since the latter are (singular) Calabi-Yau varieties, the construction gives a family of those parametrized by ${\mathfrak{L}}_l$. More precisely, the dual space $V^{\ast}_{2l-1}$ of $V_{2l-1}$ comes equipped with the basis
$$
\{T, X_1,Y_1,\ldots,X_{l-1}, Y_{l-1}\}
$$
we define the rational map
\begin{equation}\label{tauduall-map}
	{{}^{\vee}}\tau^0_l(c):\OO^{\times}_{{\mathfrak{L}}_l}(-1) \times \CC --\rightarrow \PP(V^{\ast}_{2l-1})
\end{equation}
as follows.
  Again we let 
$$
ch_s (c,\xi,\phi;u):=det(\alpha^{s,s}_c (\xi,\phi) -u{\bf 1}_{P^s([\xi],[\phi])})
$$ 
to be the value of the characteristic polynomial of $\alpha^{s,s}_c (\xi,\phi)$ at $u \in \CC$, for $s\in [0,l-1]$, and then define
\begin{equation}\label{dualtau0l-formula}
	\begin{gathered}
	{{^{\vee}}}\tau^0_l(c) ((\xi,\phi), u):=
	\\
	[T +\sum^{l-1}_{s=1} \left(ch_s (c,\xi,\phi;u) X_s +ch_{s-1} (c,\xi,\phi;u)Y_s\right)] \in \PP(V^{\ast}_{2l-1}).
\end{gathered} 
\end{equation}
This is defined outside the locus of the simultaneous vanishing of all characteristic polynomials. Varying the trace parameter gives the map
\begin{equation}
{{^{\vee}}}\tau^0_l:	\CC \longrightarrow rMaps(\OO^{\times}_{{\mathfrak{L}}_l}(-1) \times \CC, \PP(V^{\ast}_{2l-1})),
\end{equation}
where $rMaps(\bullet,\bullet)$ denotes the space of rational maps and ${{^{\vee}}}\tau^0_l$ defined by the rule
$$
\CC\ni c\mapsto {{^{\vee}}}\tau^0_l(c) \in rMaps(\OO^{\times}_{{\mathfrak{L}}_l}(-1) \times \CC, \PP(V^{\ast}_{2l-1})).
$$
 
Similar to $\tau^0_l$ the above map is $\CC^{\times}$-equivariant and on each substratum ${\mathfrak{L}}_l (h^l,\lambda)$ it satisfies the translation formula
\begin{equation}
	{{}^{\vee}}\tau^0_l(c+t) (\cdot, u) ={{}^{\vee}}\tau^0_l(c) (\cdot, u-\frac{t}{|\lambda|}),\,\,\forall c,t \in \CC,
\end{equation}
where `$\cdot$' stands for $(\xi\otimes\phi) \in \OO^{\times}_{{\mathfrak{L}}_l (h^l,\lambda)}$ and the equality is valid whenever one of the sides is well defined.

As in the case of the map $\tau^0_l (c)$ the maps ${{}^{\vee}}\tau^0_l(c)$ can interpreted via the diagram
\begin{equation}
	\xymatrix{
		&\OO^{\times}_{{\mathfrak{L}}_l}(-1) \times \CC\ar[dl]_{\pi_{\mathfrak{L}_l}\times id_{\CC}} \ar@{-->}[dr]^{{{}^{\vee}}\tau^0_l(c)}&\\ 
	{\mathfrak{L}}_l\times \CC&	&\PP(V^{\ast}_{2l-1})
	}
\end{equation}
This in turn, in the spirit of Remark \ref{rem:loops}, can be viewed as a map from ${\mathfrak{L}}_l\times \CC$ to the space $L(\PP(V^{\ast}_{2l-1}))$
of free loops in $\PP(V^{\ast}_{2l-1})$. Since the latter space parametrizes hyperplane sections of ${\bf H^{1,0}}(PG_{l})$ and those are (singular) Calabi-Yau varieties, the map ${{}^{\vee}}\tau^0_l(c)$ attaches 
to every point of  ${\mathfrak{L}}_l\times \CC$ a family of Calabi-Yau varieties parametrized by a loop. This may lead to interesting invariants related to the periods of Calabi-Yau varieties.

\subsection{Relating  $\HKC$ to $({\bf H^{1,0}}(PG_{l}),H_0)$}
Given a nonzero global section $\phi$ of the canonical bundle $\OO_C(K_C)$ we construct a $C^{\infty}$ map
$$
\lambda_{\phi}:\OO_{{\mathfrak{L}}_l (\phi)} (-1) \longrightarrow \CC(C)[u,u^{-1}];
$$
in the above ${\mathfrak{L}}_l (\phi)$ is a certain Zariski open subset of ${\mathfrak{L}}_l$ intrinsically associated to $\phi$ and $\OO_{{\mathfrak{L}}_l (\phi)} (-1)$ is the total space of the tautological line bundle $\OO_{\bf P} (-1)$ restricted to ${\mathfrak{L}}_l (\phi)$; the codomain $\CC(C)[u,u^{-1}]$ of the map is the ring of Laurent polynomials with coefficients in the function field $\CC(C)$ of $C$. 

We begin by recalling that $\PP({\cal W}_{\Sigma^0_r})$ is contained in the incidence correspondence
$$
{\bf P}=\{([\xi],[\phi]) \in \PP(H^1(\Theta_X))\times \PP(\HKC) | \xi \phi =0\}
$$
which comes with the projections
$$
\xymatrix{
&{\bf P} \ar_{p_1}[ld] \ar^{p_2}[rd]&\\
\PP(H^1(\Theta_X))&& \PP(\HKC
}
$$
on each factor. Our constructions so far involved the left side of the correspondence: we consider the strata $\PP({{\cal W}_{\Sigma^0_r}})$ as projective bundles over $\Sigma^0_r \subset \PP(H^1(\Theta_X))$ via the structure projection 
$$
p_{\Sigma^0_r}: \PP({{\cal W}_{\Sigma^0_r}}) \longrightarrow \Sigma^0_r
$$
which is the restriction of $p_1$ to $\PP({{\cal W}_{\Sigma^0_r}})$. The basic object has been the sheaf ${\cal W}_{\Sigma^0_r}$ on the base $\Sigma^0_r$ and our refinement invariants involve an additional structure on the pull back $p^{\ast}_{\Sigma^0_r} ({\cal W}_{\Sigma^0_r})$. Here we move to the right side of the correspondence. A natural object here is the space $\HKC$, equivalently, the space of global sections of the trivial bundle 
$$
\HKC \otimes \OO_{\PP(\HKC)}.
$$
Other natural sheaves come from the universal cycle
$$
{\cal Z}_C:=\{([\phi],z)\in \PP(\HKC)\times C| \, \phi(z)=0\}
$$
equipped with two natural projections
$$
\xymatrix{
	&{\cal Z}_C \ar_{q_1}[ld] \ar^{q_2}[rd]&\\
	\PP(\HKC)&&C
}
$$
The projection $q_1$ is a finite morphism of degree $2(g_C -1)$ and its branch locus is the dual variety of the canonical curve, that is, the branch locus of $q_1$ parametrizes the hyperplanes in $\PP(\HKC^{\ast})$ which are tangent to the canonical image of $C$. This of course is a classical object and it figures prominently in the proof of the Torelli theorem for curves: it is the branch locus of the Gauss map of the theta divisor of $C$ in the Jacobian $J(C)$, see \cite{G-H}. 

The sheaves on $\PP(\HKC)$ we have in mind are the direct images under $q_1$ of $\OO_{{\cal Z}_C}$ and $q^{\ast}_2 (\OO_C (K_C))$:
	$$
	q_{1 \ast} (\OO_{{\cal Z}_C}), \,\,  q_{1 \ast} (q^{\ast}_2 (\OO_C (K_C))).
	$$
The function field $\CC(C)$ can be viewed as the space of rational sections of $\OO_{{\cal Z}_C}$ via the identification
$$
 \OO_{{\cal Z}_C} =q^{\ast}_2(\OO_C).
 $$
 So the map we are about to define is a certain correspondence
 $$
 \xymatrix{
 &\phi\ar@{~>}[d]&\\
{\bf P} &{\mathfrak L}(\phi)\ar@{_{(}->}[l] \ar[r]^(.35){\lambda_{\phi}}&\PP\big(\CC(C)[u,u^{-1}]\big)
}
$$
relating the geometric object ${\mathfrak L}(\phi)$ on the side of the IVHS incidence correspondence ${\bf P}$ to the algebraic objects on the side of the geometric correspondence ${\cal Z}_C$.

 We fix a nonzero element $\phi$ in $\HKC$ and proceed with the  construction of our map $\lambda_{\phi}$.
 
 \vspace{0.2cm}
 {\it Step 1.} By definition
 the sheaf ${\cal W}_{\Sigma^0_r} $ is the subsheaf of the trivial bundle
 $\HKC \otimes \OO_{\Sigma^0_r}$. The Hodge metric on $\HKC$ gives the orthogonal $C^{\infty}$ decomposition
 $$
 \HKC \otimes \OO_{\Sigma^0_r} ={\cal W}_{\Sigma^0_r} \oplus {\cal W}^{\perp}_{\Sigma^0_r}.
 $$
 Taking the projection onto the first factor gives the $C^{\infty}$ map of vector bundles
 $$
 \HKC \otimes \OO_{\Sigma^0_r} \longrightarrow {\cal W}_{\Sigma^0_r}. 
 $$
 We consider its pull back via the projection $p:=p_{\Sigma^0_r}:\PP({\cal W}_{\Sigma^0_r}) \longrightarrow \Sigma^0_r$:
 \begin{equation}\label{p-proj-onSig0r}
 \HKC \otimes \OO_{\PP({\cal W}_{\Sigma^0_r})} \longrightarrow  p^{\ast} \left({\cal W}_{\Sigma^0_r}\right).
\end{equation}
  At this point we think of $\phi$ as a global section of  $\HKC \otimes \OO_{\PP({\cal W}_{\Sigma^0_r})}$; to make this distinction it will be denoted by $\widetilde{\phi}$. 
  
  {\it Step 2.} Denote the image of $\widetilde{\phi}$ under the morphism in \eqref{p-proj-onSig0r} by $\widehat{\phi}$. 
  On the stratum ${\mathfrak{L}}_l$ at every point $([\xi],[\psi])$ we have the orthogonal decomposition
  $$
 p^{\ast} ({\cal W}_{\Sigma^0_r} )_{([\xi],[\psi])}=W_{\xi}=W^l ([\xi],[\psi]) \oplus \left(\bigoplus^{l-1}_{s=0}P^s([\xi],[\psi]) \right).
 $$
 The component of $\widehat{\phi}$ in $P^s([\xi],[\psi])$ is denoted by
$\widehat{\phi}^s ([\xi],[\psi])$. 

Fix an auxiliary parameter $u$. Assign the weight $u^{(-1)^s s}$ to the summand $P^s([\xi],[\psi])$ and expand  $\widehat{\phi}([\xi],[\psi])$, the value of $\widehat{\phi}$ at $([\xi],[\psi])$, according to the weighted direct sum decomposition
\begin{equation}\label{phi-wdecomp}
\widehat{\phi}([\xi],[\psi])=\sum^{l-1}_{s=0}\widehat{\phi}^s ([\xi],[\psi]) u^{(-1)^s s}.
\end{equation}

\begin{defi}\label{def:phi-admiss}
	The domain of $\phi \in \HKC$ in ${\mathfrak{L}}_l$ is the subset
	$$
	{\mathfrak{L}}_l (\phi):=\{([\xi],[\psi])\in {\mathfrak{L}}_l | \, \widehat{\phi}^s ([\xi],[\psi]) \neq 0,\forall s\in [0,l-1] \}.
	$$ 
	We say that $\phi \in \HKC$ is admissible if its domain ${\mathfrak{L}}_l (\phi)$ is nonempty.
\end{defi}
Observe that ${\mathfrak{L}}_l (\phi)$ is an algebraic subvariety of ${\mathfrak{L}}_l$ since depends only on the $([\xi],[\psi])$-filtrations,
while the components $\widehat{\phi}^s ([\xi],[\psi])$ in the decomposition \eqref{phi-wdecomp} vary in a $C^{\infty}$ manner with $([\xi],[\psi])$ in
 ${\mathfrak{L}}_l (\phi)$. 

\vspace{0.2cm}
{\it Step 3.} Let $\phi$ be admissible in the sense of Definition \ref{def:phi-admiss} and let ${\mathfrak{L}}_l (\phi)$ be its domain.  Consider the total space $\OO^{\times}_{{\mathfrak{L}}_l (\phi)} (-1)$ of the tautological line bundle $\OO_{\bf P} (-1)$ restricted to ${\mathfrak{L}}_l (\phi)$ and the zero section removed. Fix a value $c$ of the trace parameter and recall the representation
$$
\rho_c :\OO^{\times}_{{\mathfrak{L}}_l (\phi)} (-1) \longrightarrow {\mathfrak{Reps}}(PG_l):
$$
at each point $([\xi],[\psi])$ of ${\mathfrak{L}}_l (\phi)$ and at every vertex $(s)$ (resp. $(s)'$) we have the vector space $P^s([\xi],[\psi])$ as well as the linear maps
$$
\alpha^{t,s}_c(\xi,\psi): P^s([\xi],[\psi]) \longrightarrow P^t([\xi],[\psi])
$$
for $t=s-1,s,s+1$, where as before $c$ is the trace parameter and $\xi\otimes \psi$ is a nonzero vector in the fibre of $\OO_{{\mathfrak{L}}_l (\phi)} (-1)$ lying over $([\xi],[\psi])$.
 For $\xi\otimes \psi \in \OO^{\times}_{{\mathfrak{L}}_l (\phi)} (-1)_{([\xi],[\psi])}$ define the rational functions on $C$
$$
\begin{gathered}
f^0_s (\xi,\psi,c) :=\frac{\alpha^{s,s}_c(\xi,\psi)(\widehat{\phi}^s ([\xi],[\psi]))}{\widehat{\phi}^s ([\xi],[\psi])}, \,\, \forall s\in [0,l-1],
\\
f^+_s (\xi,\psi) :=\frac{\alpha^{s+1,s}(\xi,\psi)(\widehat{\phi}^{s} (\xi,\psi))}{\widehat{\phi}^{s+1} ([\xi],[\psi])}, \,\,\forall s\in [0,l-2],
\\
f^-_s (\xi,\psi) :=\frac{\alpha^{s-1,s}(\xi,\psi)(\widehat{\phi}^s ([\xi],[\psi]))}{\widehat{\phi}^{s} ([\xi],[\psi])}, \,\, \forall s\in [1,l-1].
\end{gathered}
$$
Recall: only the diagonal blocks depend on the trace parameter. Hence only  the first function above depends on $c$. Furthermore, this dependence is subject to the translation formula
\begin{equation}\label{f0r-tr-change}
	f^0_s (\xi,\psi,c+t)=f^0_s (\xi,\psi,c) +\frac{t}{w([\xi],[\psi])}, \,\forall c,t \in \CC,
\end{equation}
where the second term on the right is the constant function of the value
 $\frac{t}{w([\xi],[\psi])}$ and where $w([\xi],[\psi])=dim(W_{\xi}/W^l([\xi],[\psi]))=g-r-dim (W^l([\xi],[\psi]))$; we remind the reader that this comes from the formula 
 $$
 \alpha^{s,s}_{c+t} -\alpha^{s,s}_{c}=\frac{t}{w([\xi],[\psi])} {\bf 1}_{P^s}
 $$
of changing the trace parameter, see \eqref{tr-v-change}; the \eqref{f0r-tr-change} is an immediate consequence of this formula. 

\vspace{0.2cm}
We think of the above functions as labels of the edges of the graph $PG_{l}$.
Denote by $\CC(C)$ the function field of $C$ and define the map
\begin{equation}\label{lambdaphi}
\lambda_{\phi} (c):\OO^{\times}_{{\mathfrak{L}}_l (\phi)} (-1) \longrightarrow \CC(C)[u,u^{-1}],
\end{equation}
by the formula
\begin{equation}\label{lambdaphi-formula}
\lambda_{\phi} (c) (\xi,\psi) = f^0_0 (\xi,\psi,c) +\sum^{l-1}_{s=1} \left(f^-_{s} (\xi,\psi) + f^+_{s-1} (\xi,\psi) \right) u^{(-1)^s s},
\end{equation}
for every $([\xi],[\psi]) \in {\mathfrak{L}}_l (\phi)$ and every $\xi\otimes \psi$ in $\OO^{\times}_{{\mathfrak{L}}_l (\phi)} (-1)$ lying over $([\xi],[\psi])$. 

Observe that for all $\xi$ and $\phi$ nonzero
 the rational function
$$
f^-_{s} (\xi,\psi) + f^+_{s-1} (\xi,\psi)=\frac{\alpha^{s-1,s}(\xi,\psi)(\widehat{\phi}^{s} ([\xi],[\psi])) +\alpha^{s,s-1}(\xi,\psi)(\widehat{\phi}^{s-1} ([\xi],[\psi]))}{\widehat{\phi}^{s}([\xi],[\psi])},
$$
 for each $s\in[1,l-1]$ is nonzero and has poles along the divisor
 $Z_{\widehat{\phi}^{s} ([\xi],[\psi])}=\{\widehat{\phi}^{s} ([\xi],[\psi])=0\}$. The `constant' term $f^0_{0} ([\xi],[\psi],c)$ in \eqref{lambdaphi-formula} has poles
 along $Z_{\widehat{\phi}^0 ([\xi],[\psi])}=\{\widehat{\phi}^0 ([\xi],[\psi])=0\}$. We summarize the construction in the following statement.
 \begin{pro}\label{pro:lambdaphi}
 	Let $\phi$ be admissible global section of the canonical bundle $\OO_C(K_C)$ and let ${\mathfrak{L}}_l (\phi)$ be its domain, see Definition \ref{def:phi-admiss}. Then for every value $c\in \CC$ of the trace parameter there is map
 	$$
 	\lambda_{\phi} (c): \OO^{\times}_{{\mathfrak{L}}_l (\phi)} (-1) \longrightarrow \CC(C)[u,u^{-1}]
 	$$
 	associating to every $\xi\otimes \phi \in \OO^{\times}_{{\mathfrak{L}}_l (\phi)} (-1)$ the Laurent polynomial
 	$$
 	\lambda_{\phi} (c) (\xi,\psi) = f^0_0 (\xi,\psi,c) +\sum^{l-1}_{s=1} \left(f^-_{s} (\xi,\psi) + f^+_{s-1} (\xi,\psi) \right) u^{(-1)^s s}
 	$$
 	in a formal variable $u$ with coefficients $f^0_0 (\xi,\psi,c)$, $f^{\pm}_s (\xi,\psi)$ in the function field $\CC(C)$ of $C$.
 	
 	The coefficients of $\lambda_{\phi} (c) (\xi,\psi)$ have the following properties.
 	
 	1) They are regular outside the divisor
 	$$
 	\sum^{l-1}_{s=0}Z_{\widehat{\phi}^{s} ([\xi],[\psi])} \in |lK_C|.
 	$$
 	
 	2) The coefficients $(f^-_{s} (\xi,\psi) + f^+_{s-1} (\xi,\psi))$, for
 	$s\in [1,l-1]$, are nonzero for all $\xi\otimes\phi$ in $\OO^{\times}_{{\mathfrak{L}}_l (\phi)} (-1)$ .
 	
 	Furthermore, the natural $\CC^{\times}$-action on $\xi$ and $\psi$
 	acts on the maps $\lambda_{\phi} (c)$ as follows
 	$$
 	\lambda_{\phi} (c) (t\xi,s\psi)=\lambda_{\phi} (tsc) (\xi,\psi), \forall s,t \in \CC^{\times}. 
 	$$
 	In addition, over the stratum
 	$$
 	{\mathfrak{L}}_l (\phi)(h^l,\lambda):= {\mathfrak{L}}_l (\phi) \bigcap {\mathfrak{L}}_l (h^l,\lambda),
 	$$
 	whenever nonempty, the following holds:
 	$$
 	\lambda_{\phi} (c) (t\xi,s\psi)=ts\lambda_{\phi} (c) (\xi,\psi) +\frac{(1-st)c}{|\lambda|}, \,\,\forall s,t \in \CC^{\times},
 	$$
 	where $|\lambda|=g-r -h^l$ is the weight of the partition $\lambda$.
 \end{pro}
 \begin{pf}
 	Only the assertions about the $\CC^{\times}$-action are not explicit in the construction. The formula
 	$$
 	\lambda_{\phi} (c) (t\xi,s\psi)=\lambda_{\phi} (tsc) (\xi,\psi), \forall s,t \in \CC^{\times},
 	$$
 	is the equivariance of $\CC^{\times}$-action in Proposition \ref{pro:tau0-n-equivariant}.
 	
 	For the formula
 	\begin{equation}\label{lambdaphi-add-formuala}
 	\lambda_{\phi} (c) (t\xi,s\psi)=ts\lambda_{\phi} (c) (\xi,\psi) +\frac{(1-st)c}{|\lambda|}, \,\,\forall s,t \in \CC^{\times},
 \end{equation}
 	we need to recall the formula
 	$$
 	\alpha^{(2)}_{t\xi} (s\psi,\bullet)=ts\alpha^{(2)}_{\xi} (\psi,\bullet)+f(\psi)id_{W_{\xi}/W^l_{\xi}([\psi])},
 	$$
 	for some linear function $f:W_{\xi}\longrightarrow \CC$.
 	If we impose $\alpha^{(2)}_{\xi} (\psi,\bullet)$ and $\alpha^{(2)}_{t\xi} (s\psi,\bullet)$ both to have trace $c$, then follows the relation
 	$$
 	\alpha^{(2)}_{t\xi,c} (s\psi,\bullet)=ts\alpha^{(2)}_{\xi,c} (\psi,\bullet)+\frac{(1-ts)c}{|\lambda|}id_{W_{\xi}/W^l_{\xi}([\psi])}.
 	$$
 	From this we deduce
 	$$
 	\alpha^{i,j}_c (t\xi,s\psi)=ts\alpha^{i,j}_c (\xi,\psi)
 	$$
 	for all blocks with $i\neq j$, and for the diagonal blocks
 	$$
 	\alpha^{i,i}_c (t\xi,s\psi)=ts\alpha^{i,i}_c (\xi,\psi)+ \frac{(1-ts)c}{|\lambda|}{\bf 1}_{P^i([\xi],[\psi])}.
 	$$
 	From these relations it follows
 	$$
 	f^{\pm}_s(t\xi,s\psi)=tsf^{\pm}_s(\xi,\psi),\hspace{0.2cm} f^{0}_s(t\xi,s\psi)=tsf^{0}_s(\xi,\psi)+ \frac{(1-ts)c}{|\lambda|}.
 	$$
 	The formula \eqref{lambdaphi-add-formuala} is now immediate from the definition of $\lambda_{\phi} (c)$.
 \end{pf}
 \vspace{0.2cm}
 We will now use the above constructions to relate certain dense open parts of the Cartesian product $C^l$ of the curve $C$ with the hyperplane sections of
 ${\bf H^{1,0}}(PG_{l})$.
  
  The starting point is the polar subvariety in $C\times {\mathfrak{L}}_l (\phi)$. Those are the following incidence correspondences 
 $$
 D_s (\phi):=\{(z,([\xi],[\psi]) )\in C\times {\mathfrak{L}}_l (\phi) | z\in Z_{\widehat{\phi}^s ([\xi],[\psi])}\},
 $$
 for every $s\in [0,l-1]$; they are zero loci of $C^{\infty}$ global sections of the relative canonical bundle of the projection morphism
 $$
 C\times {\mathfrak{L}}_l (\phi) \longrightarrow {\mathfrak{L}}_l (\phi).
 $$
  Consider the Cartesian product
 $C^{l}$ where the factors are ordered by the set $[0,l-1]$ of labels of pairs of aligned black-white vertices of $PG_l$. Denote the projection of $C^l$ onto the $s$-th factor by $pr_s$ and lift it to the morphism
 $$
  \xymatrix{
 	C^{l}\times{\mathfrak{L}}_l (\phi) \ar[r]^{\rho_s}& C\times {\mathfrak{L}}_l (\phi)
 }
$$
making the diagram
 $$
 \xymatrix{
 	C^{l}\times{\mathfrak{L}}_l (\phi) \ar[r]^{\rho_s}\ar[d]& C\times {\mathfrak{L}}_l (\phi)\ar[d]\\
 	C^{l} \ar[r]^{pr_s}&C
 } 
 $$ 
commutative, that is, $\rho_s=pr_s \times id_{{\mathfrak{L}}_l (\phi)}$.
 We pull back to $C^{l}\times {\mathfrak{L}}_l (\phi)$ the subvarieties $D_s (\phi)$
$$
\widetilde{D}_s (\phi)= \rho^{\ast}_s(D_s(\phi))
$$
and define the total polar subvariety of $\phi$ as the union of $\widetilde{D}_s (\phi)$:
$$
\widetilde{D} (\phi)=\bigcup^{l-1}_{s=0} \widetilde{D}_s (\phi).
$$ 
Take the complement
$$
U^l(\phi)=	\left(C^{l}\times{\mathfrak{L}}_l (\phi) \right) \setminus \widetilde{D} (\phi)
$$
  of the total polar subvariety $\widetilde{D} (\phi)$ and consider the pull back $\widetilde{U^l}(\phi)$ of $U^l(\phi)$ under the projection
  $$
  C^l\times \OO^{\times}_{{\mathfrak{L}}_l (\phi)} (-1) \longrightarrow C^{l}\times{\mathfrak{L}}_l (\phi).
$$
 Observe that we can evaluate 
  $\lambda_{\phi}(c)$ defined in \eqref{lambdaphi} over $\widetilde{U^l} (\phi)$: 
  
  for a point ${\bf u}=(z_0,z_1,\ldots,z_{l-1};\xi\otimes\psi )$ in $\widetilde{U^l}(\phi)$ set
  $$
  \begin{gathered}
  \lambda^{(l)}_{\phi} (c)({\bf u}):=
  \\
   \left(f^0_0 (\xi,\psi,c) (z_0), f^-_{s} (\xi,\psi)(z_s) , f^+_s (\xi,\psi)(z_s) | \,s=1,\ldots,l-1\right) \in \CC^{2l-1},
   \end{gathered}
   $$
  and define the map
  $$
 exp \left(\lambda^{(l)}_{\phi}(c)\right): \widetilde{U^l}(\phi) \times \CC^{\times} \longrightarrow \PP(V^{\ast}_{2l-1})
  $$
by the rule
$$
\begin{gathered}
\widetilde{U^l}(\phi) \times \CC^{\times} \ni (z_0,z_1,\ldots,z_l;\xi \otimes\psi;t ) \mapsto
\\
\left[a_0(z_0,\xi \otimes\psi,c)T+ \sum^{l-1}_{s=1} \left( a^-_s (z_s,\xi \otimes\psi) X_s + a^+(z_s,\xi\otimes\psi) Y_s\right)t^{(-1)^s s}\right] \in \PP(V^{\ast}_{2l-1}),
\end{gathered}
$$
where the coefficients $a_0(z_0,\xi\otimes\psi,c)$, $a^{\pm}_s (z_s,\xi\otimes\psi)$ are defined by exponentiating the coordinates of $\lambda^{(l)}_{\phi} (c)({\bf u})$:
$$
a_0(z_0,\xi\otimes\psi,c)=exp(f^0_0 (\xi,\psi,c) (z_0)), \,\, a^{\pm}_s =exp\left(f^{\pm}_{s} (\xi,\psi)(z_s)\right), \,\forall s\in [1,l-1].
$$
The construction can be viewed as a $C^{\infty}$ family of Calabi-Yau varieties parametrized by \linebreak   $\widetilde{U^l}(\phi) \times \CC^{\times}$.
If we fix a point $([\xi],[\psi])$ in ${\mathfrak{L}}_l (\phi)$ and restrict the above construction to the slice of $\widetilde{U^l}(\phi) \times \CC^{\times}$ over that point, we obtain the {\it holomorphic} map
$$
  exp \left(\lambda^{(l)}_{\phi} (c)\right)_{([\xi],[\psi])}: C^l(\phi;([\xi],[\psi])) \times \OO^{\times}_{{\mathfrak{L}}_l (\phi),([\xi],[\psi])} (-1) \times \CC^{\times} \longrightarrow \PP(V^{\ast}_{2l-1}),
  $$
  where $C^l(\phi;([\xi],[\psi]))$ is the complement of the divisor
  $$
  \widetilde{D} (\phi)([\xi],[\psi])=\sum^{l-1}_{s=0} pr^{\ast}_s (Z_{\widehat{\phi}^s([\xi],[\psi])})
  $$
  in the Cartesian product $C^l$, that is, we have
  $$
  C^l(\phi;([\xi],[\psi]))=\{ (z_0,z_1,\ldots,z_{l-1})\in C^l | z_i \notin \sum^{l-1}_{s=0} Z_{\widehat{\phi}^s([\xi],[\psi])},\, \forall i\in [0,l-1] \}.
  $$
  Thus the above construction attaches the family
  $$
   exp \left(\lambda^{(l)}_{\phi}(c)\right)_{([\xi],[\psi])}({\bf z}, \bullet,\bullet):\OO^{\times}_{{\mathfrak{L}}_l (\phi),([\xi],[\psi])} (-1) \times\CC^{\times} \longrightarrow \PP(V^{\ast}_{2l-1})
   $$
   of Calabi-Yau varieties to every point $\bf z$ of specific Zariski dense open subsets of $C^l$. One could expect interesting  new invariants for curves coming from the ones of these families of Calabi-Yau varieties.
  
\section{The trivalent completion of $PG_{l}$}

Observe: all vertices of $PG_{l}$ with the exception of the first
and the last pairs have precisely three edges. 
We complete  $PG_{l}$ to a {\it trivalent} graph by connecting the leftmost (resp. rightmost) white vertex labeled $(l-1)$ (resp., $0$) to the black vertex $0'$ (resp., $(l-1)'$). This graph will be denoted $\widehat{PG}_{l}$. Below is given the completion $\widehat{PG}_{4}$ of $PG_4$.
$$
 \begin{tikzpicture}
 	[place/.style={circle,draw=black,thick},
 	transition/.style={circle,draw=black,fill=black}]
 	\node (white3) at (0,2) [place] [label={above:$3$}] {};
 	\node (black3) at (0,0) [transition] [label={below:$3'$}]{};
 	\node (white2) at (2,2) [place] [label={above:$2$}]{};
 	\node (black2) at (2,0) [transition] [label={below:$2'$}]{};
 	\node (white1) at (4,2) [place] [label={above:$1$}]{};
 	\node (black1) at (4,0) [transition] [label={below:$1'$}]{};
 	\node (white0) at (6,2) [place] [label={above:$0$}] {};
 	\node (black0) at (6,0) [transition] [label={below:$0'$}]{};
 	\draw[red,ultra thick] (white0) to (black0);
 	\draw[red,ultra thick] (white1) to (black1);
 	\draw[red,ultra thick] (white2) to (black2);
 	\draw[red,ultra thick] (white3) to (black3);
 	\begin{scope}[thick]
 	\draw[blue][-] (white3) to (black2);	
 	\draw[blue][-] (white2) to (black1);
 	\draw[blue][-] (white1) to (black0);
 	\end{scope}
 	\draw[red, dotted,thick][-] (white2) to (black3);
 	\draw[red,dotted,thick][-] (white1) to (black2);
 	\draw[red,dotted,thick][-] (white0) to (black1);
 	\draw[blue,thick][-] (white0) .. controls (10,-0.5) and (4,-1)  .. (3,-1) .. controls (-1,-1) and (-1,1).. (black3);
 		\draw[red,dotted,thick][-] (white3) .. controls (-4,-0.5) and (2.5,-2.1)  .. (3,-2) .. controls (8,-1.5) and (7,2).. (black0);
 \end{tikzpicture}
 $$
 
 Next we equip the graph $\widehat{PG}_{l}$ with the structure of a {\it ribbon} graph. This means that a cyclic ordering of edges at every vertex of the graph is chosen. We choose the {\it counterclockwise} cyclic ordering of edges at every vertex and call it {\it the standard ribbon graph structure} of $\widehat{PG}_{l}$. Unless said otherwise, we always assume that our graph is equipped with this structure.
 
 Recall that a ribbon graph $G$, in addition to vertices, $V_G$, and (nonoriented) edges, $E_G$, acquires so called boundary cycles, $B_G$. This way a ribbon graph becomes a deformation retract of an oriented two manifold with boundary, where connected components of the boundary correspond to the boundary cycles. By gluing in the discs along the boundary cycles, one obtains a compact oriented two manifold which will be denoted $S_G$. The topological Euler characteristic of $S_G$ is computed by the formula
 \begin{equation}\label{SG-Euler}
 \chi_{top} (S_G)=|V_G|-|E_G| +|B_G|.
 \end{equation}
   
 We apply all this to the graph $\widehat{PG}_{l}$ with its standard ribbon structure.
 \begin{lem}\label{lem:bound-cycles}
 	Every boundary cycle of $\widehat{PG}_{l}$ is given by the
 	path in $G$
 	$$
 	i \to i' \to (i-1) \to (i-2)' \mapsto (i-2) \to (i-1)' \to i,
 	$$
 	for every $i\in [0,l-1]$, where all values are taken modulo $l$, that is, the vertex $(-1)=(l-1)$, $(-2)=l-2$, etc.
 	
 	In particular, the number of boundary cycles
 	$$
 	|B_{\widehat{PG}_{l}}| =l.
 	$$
 \end{lem}
\begin{pf}
Let $G=\widehat{PG}_{l}$.	Consider the set $E^{or}_G$ of {\it oriented} edges of $G$, that is, every edge $e$ in $E_G$ occurs twice in $E^{or}_G$ with two opposite orientations. The set  $E^{or}_G$ comes with three natural permutations:

$\bullet$ $\sigma_1$ is the involution which reverses the orientation on each edge
$e$ in $E^{or}_G$;

$\bullet$ $\sigma_0$ is defined by sending an oriented edge $e$ to the next one in the cyclic ordering around the tale vertex of $e$;

$\bullet$ $\sigma_{\infty}=\sigma^{-1}_0 \sigma_1$, that is, for an oriented edge $e$, first reverse the orientation on an edge to obtain $e^{op}$ and then take the {\it preceding} edge in the cyclic order of edges around the tale vertex of $e^{op}$. 

The orbits of the action of these permutations on the set $E^{or}_G$ give respectively the set of edges $E_G$, the set of vertices $V_G$, and the set of $B_G$ of boundary cycles of $G$. Thus to see a boundary cycle we need to determine the orbit of a edge under the action of $\sigma_{\infty}$.

Start with the oriented edge
$$
e_i=i \to i'
$$
and compute
$$
\sigma_{\infty} (e_i)= \sigma^{-1}_0 \sigma_1 (e_i) =\sigma^{-1}_0 (i'\to i)= i'\to (i-1) .
$$
Continuing in this fashion gives the path stated in the lemma.

It remains to see that the orbits
$$
B_i=\text{the $\sigma_{\infty}$-orbit of $e_i$}, 
$$
for $i\in [0,l-1]$, form the complete set of orbits. Indeed, those orbits are distinct and their union gives us the subset of $6l$ edges of $E^{or}_G$ and this is precisely the number of oriented edges of $G$.   
\end{pf}

The immediate corollary of the above is that the surface 
$S_{\widehat{PG}_{l}}$ associated to the ribbon graph $\widehat{PG}_{l}$ is the topological torus $\mathbb{T}$ and the ribbon  graph $\widehat{PG}_{l}$ can be embedded into $\mathbb{T}$ so that the complement consists of $l$ two-cells.
\begin{cor}\label{cor:torus}
	The compact oriented surface $S_{\widehat{PG}_{l}}$ associated to the ribbon graph $\widehat{PG}_{l}$ is a surface of genus $1$. More precisely, the graph $\widehat{PG}_{l}$ can be embedded into a topological torus $\mathbb{T}$ so that its complement is the disjoint union of
	$l$ open two-cells
	$$
	D_0,D_1,\ldots, D_{l-1};
	$$
	the boundaries of the cells are the boundary cycles
	$$
	B_0,B_1,\ldots, B_{l-1}
	$$
	described in Lemma \ref{lem:bound-cycles}
	\end{cor}

\begin{pf}
	The formula \eqref{SG-Euler} for the topological Euler characteristic of $S_{\widehat{PG}_{l}}$ reads
	$$
	\chi_{top} (S_{\widehat{PG}_{l}})=2l-3l +l=0.
	$$
	Hence the genus of $S_{\widehat{PG}_{l}}$ is $1$.
	The rest of the statements of the corollary follows from the description of the boundary cycles in Lemma \ref{lem:bound-cycles}.
\end{pf}

\begin{rem}\label{rem:dimer}
	We have learned that the bipartite graph $\widehat{PG}_{l}$ can be realized on a compact two dimensional oriented torus $\mathbb{T}$
	in such a way that the complement
	$$
	\mathbb{T} \setminus \widehat{PG}_{l}
	$$
	is the disjoint union of  $2$-cells. This is known as a dimer model. Lifting to $\RR^2$, the universal cover of the torus, one obtains a tiling of the plane. Hence another name for the model - brane tiling.
	There is a huge interest of physicists in those models in connection with various dualities in string theory; the literature on the subject is enormous, as a sample we refer to an excellent review \cite{Ken} and the references therein.
\end{rem}
To reflect the cyclic nature of the graph  $\widehat{PG}_{l}$
it is more convenient to draw two concentric circles, one for black vertices (the outer circle) and the other (inner circle) for white ones and place the vertices of each color on its circle so that the vertices $i$ and $i'$ are on the same radial segment. Below is the picture for $l=6$.
$$
\begin{tikzpicture}
	[place/.style={circle,draw=black,thick},
	transition/.style={circle,draw=black,fill=black}]
	\node (white5) at (90:1.5cm) [place] [label={below:$5$}] {};
	\node (black5) at (90:2.5cm) [transition] [label={above:$5'$}]{};
	\node (white4) at (150:1.5cm) [place] [label={below right:$4$}]{};
	\node (black4) at (150:2.5cm) [transition] [label={above left:$4'$}]{};
	\node (white3) at (210:1.5cm) [place] [label={above right:$3$}]{};
	\node (black3) at (210:2.5cm) [transition] [label={below left:$3'$}]{};
	\node (white2) at (270:1.5cm) [place] [label={above:$2$}] {};
	\node (black2) at (270:2.5cm) [transition] [label={below:$2'$}]{};
	\node (white1) at (330:1.5cm) [place] [label={above left:$1$}] {};
	\node (black1) at (330:2.5cm) [transition] [label={below right:$1'$}]{};
	\node (white0) at (30:1.5cm) [place] [label={below left:$0$}] {};
	\node (black0) at (30:2.5cm) [transition] [label={above right:$0'$}]{};
	\draw[ gray,ultra thin] (white5) to [bend right=30] (white4);
	\draw[gray,ultra thin] (white4) to [bend right=30](white3);
	\draw[gray,ultra thin] (white3) to [bend right=30](white2);
	\draw[gray,ultra thin] (white2) to [bend right=30](white1);
	\draw[gray,ultra thin] (white1) to [bend right=30](white0);
	\draw[gray,ultra thin] (white0) to [bend right=30](white5);
	\draw[gray,ultra thin] (black5) to [bend right=30](black4);
	\draw[gray,ultra thin] (black4) to [bend right=30](black3);
	\draw[gray,ultra thin] (black3) to [bend right=30](black2);
	\draw[gray,ultra thin] (black2) to [bend right=30](black1);
	\draw[gray,ultra thin] (black1) to [bend right=30](black0);
	\draw[gray,ultra thin] (black0) to [bend right=30](black5);
	\draw[red,ultra thick] (white0) to (black0);
	\draw[red,ultra thick] (white1) to (black1);
	\draw[red,ultra thick] (white2) to (black2);
	\draw[red,ultra thick] (white3) to (black3);
	\draw[red,ultra thick] (white4) to (black4);
	\draw[red,ultra thick] (white5) to (black5);
\end{tikzpicture} 
$$
Further, for each white vertex $i$ we draw the remaining two edges: one joining it to $(i-1)'$ and another to $(i+1)'$; as before we use the convention that the values are taken in $[0,l-1]$ modulo $l$. Below is the picture of the graph  $\widehat{PG}_{6}$.

$$
\begin{tikzpicture}
	[place/.style={circle,draw=black,thick},
	transition/.style={circle,draw=black,fill=black}]
	\node (white5) at (90:1.5cm) [place] [label={below:$5$}] {};
	\node (black5) at (90:3cm) [transition] [label={above:$5'$}]{};
	\node (white4) at (150:1.5cm) [place] [label={below right:$4$}]{};
	\node (black4) at (150:3cm) [transition] [label={above left:$4'$}]{};
	\node (white3) at (210:1.5cm) [place] [label={above right:$3$}]{};
	\node (black3) at (210:3cm) [transition] [label={below left:$3'$}]{};
	\node (white2) at (270:1.5cm) [place] [label={above:$2$}] {};
	\node (black2) at (270:3cm) [transition] [label={below:$2'$}]{};
	\node (white1) at (330:1.5cm) [place] [label={above left:$1$}] {};
	\node (black1) at (330:3cm) [transition] [label={below right:$1'$}]{};
	\node (white0) at (30:1.5cm) [place] [label={below left:$0$}] {};
	\node (black0) at (30:3cm) [transition] [label={above right:$0'$}]{};
	\begin{scope}[gray, ultra thin]
	\draw (white5) to [bend right=30] (white4);
	\draw (white4) to [bend right=30](white3);
	\draw (white3) to [bend right=30](white2);
	\draw (white2) to [bend right=30](white1);
	\draw (white1) to [bend right=30](white0);
	\draw (white0) to [bend right=30](white5);
	\draw (black5) to [bend right=30](black4);
	\draw (black4) to [bend right=30](black3);
	\draw (black3) to [bend right=30](black2);
	\draw (black2) to [bend right=30](black1);
	\draw (black1) to [bend right=30](black0);
	\draw (black0) to [bend right=30](black5);
	\end{scope}
	\draw[red,ultra thick] (white0) to (black0);
	\draw[red,ultra thick] (white1) to (black1);
	\draw[red,ultra thick] (white2) to (black2);
	\draw[red,ultra thick] (white3) to (black3);
	\draw[red,ultra thick] (white4) to (black4);
	\draw[red,ultra thick] (white5) to (black5);
		\begin{scope}[thick]
				\draw[blue][-] (white3) to [bend right=45](black2);	
				\draw[blue][-] (white2) to [bend right=45](black1);
				\draw[blue][-] (white1) to [bend right=45](black0);
					\draw[blue][-] (white5) to [bend right=45](black4);	
						\draw[blue][-] (white4) to [bend right=45](black3);	
						\draw[blue][-] (white0) to [bend right=45](black5);	
			\end{scope}
		\draw[red, dotted,thick][-] (white2) to [bend left=45](black3);
		\draw[red,dotted,thick][-] (white1) to [bend left=45](black2);
		\draw[red,dotted,thick][-] (white0) to [bend left=45](black1);
		\draw[red,dotted,thick][-] (white3) to [bend left=45](black4);
		\draw[red,dotted,thick][-] (white4) to [bend left=45](black5);
		\draw[red,dotted,thick][-] (white5) to [bend left=45](black0);
\end{tikzpicture} 
$$  

\noindent
The gray thin lines are not the part of the graph - they are drawn to see  the torus emerging in the picture: the annulus formed by those two concentric circles is
a `half' of the torus; it is bisected by the plane of the page with the other half of the torus being on the other side of the page. The edges in  `solid' colors should be thought as drawn on the visible side of the torus, while the dotted edges are on the other half of the torus.

\vspace{0.2cm}
By concatenation the oriented edges of a graph form paths, that is, a path on a graph is a word $e_m e_{m-1} \cdots e_1$ composed by oriented edges $e_i$'s of a graph and read from {\it right to left} :
$$
t(e_{i+1})=h(e_i), \,\forall i\in [1,m-1],
$$
where $t(e)$ (resp. $h(e)$) is the tail (resp. head) of $e$. We will be interested in closed paths of the graph $\widehat{PG}_{l}$. The abelian group generated by the closed paths of a graph $G$ is thought as the first homology group
$$
H_1 (G,\ZZ)
$$
of $G$. For a ribbon graph $G$, viewed as a graph embedded into the surface $S_G$, we have the identification
$$
H_1 (G,\ZZ) \cong H_1(\stackrel{\circ}{S}_G),
$$
where $\stackrel{\circ}{S}_G$ is the open surface obtained from $S_G$ by removing open two-cells bordered by the boundary cycles of $G$. This gives the surjective homomorphism
$$
\xymatrix{
H_1 (G,\ZZ) \ar[r]^{j}&H_1(S_G)\ar[r]&0.
}
$$
The kernel $ker(j)$ is the subgroup of $H_1 (G,\ZZ)$ generated by the boundary cycles of $G$. In the case of the graph $\widehat{PG}_{l}$ we obtain
\begin{equation}\label{hom-seq}
\xymatrix{
0\ar[r]&\ZZ\{B_i\}\ar[r]&H_1 (\widehat{PG}_{l},\ZZ) \ar[r]^{j}&H_1(\mathbb{T}, \ZZ)\ar[r]&0,
}
\end{equation}
where $\{B_i\}=\{B_0,B_1,\ldots,B_{l-1}\}$ is the set of boundary cycles of $\widehat{PG}_{l}$, see Lemma \ref{lem:bound-cycles}.

On any bipartite ribbon graph one distinguishes  zig-zag paths. We briefly recall this and refer the reader for more details to \cite{Go-K} and references their in. On a ribbon graph at every vertex we have a cyclic ordering of edges incident to a given vertex. We assume that our graph has valency at least three at every vertex. So when a path arrives along an edge $e$ to the vertex $h(e)$, the head of $e$, we have two neighboring edges: one preceding $e$ and the other one following $e$ in the cyclic order; one says that a path turns {\it left} (resp. {\it right}) at $h(e)$ if it continues along the edge {\it preceding} (resp. {\it following}) $e$.
If $G$ is a bipartite ribbon graph with all of its vertices of valency at least three, a {\it zig-zag path} is a path on $G$ which turns {\it right} at each black vertex and turns {\it left} at the white ones. 
\begin{example}
	On the graph $\widehat{PG}_{6}$ by starting with the oriented edge 
	$$
	(2\to 2')
	$$
	one obtains the following zig-zag path
	$$
	2\to 2' \to 3\to3'\to 4\to4'\to5\to5'\to0\to0' \to 1\to 1' \to2.
	$$
	This is  visualized on the drawing below:
	
	\begin{equation}\label{zz1}
	\begin{tikzpicture}
		[place/.style={circle,draw=black,thick},
		transition/.style={circle,draw=black,fill=black}]
		\node (white5) at (90:1.5cm) [place] [label={below:$5$}] {};
		\node (black5) at (90:3cm) [transition] [label={above:$5'$}]{};
		\node (white4) at (150:1.5cm) [place] [label={below right:$4$}]{};
		\node (black4) at (150:3cm) [transition] [label={above left:$4'$}]{};
		\node (white3) at (210:1.5cm) [place] [label={above right:$3$}]{};
		\node (black3) at (210:3cm) [transition] [label={below left:$3'$}]{};
		\node (white2) at (270:1.5cm) [place] [label={above:$2$}] {};
		\node (black2) at (270:3cm) [transition] [label={below:$2'$}]{};
		\node (white1) at (330:1.5cm) [place] [label={above left:$1$}] {};
		\node (black1) at (330:3cm) [transition] [label={below right:$1'$}]{};
		\node (white0) at (30:1.5cm) [place] [label={below left:$0$}] {};
		\node (black0) at (30:3cm) [transition] [label={above right:$0'$}]{};
		\begin{scope}[gray,ultra thin]
		\draw (white5) to [bend right=30] (white4);
		\draw (white4) to [bend right=30](white3);
		\draw (white3) to [bend right=30](white2);
		\draw (white2) to [bend right=30](white1);
		\draw (white1) to [bend right=30](white0);
		\draw (white0) to [bend right=30](white5);
		\draw (black5) to [bend right=30](black4);
		\draw (black4) to [bend right=30](black3);
		\draw (black3) to [bend right=30](black2);
		\draw (black2) to [bend right=30](black1);
		\draw (black1) to [bend right=30](black0);
		\draw (black0) to [bend right=30](black5);
		\end{scope}
		\draw[red,ultra thick][->,>= {Stealth}] (white0) to (black0);
		\draw[red,ultra thick][->,>= {Stealth}] (white1) to (black1);
		\draw[red,ultra thick][->,>= {Stealth}] (white2) to (black2);
		\draw[red,ultra thick][->,>= {Stealth}] (white3) to (black3);
		\draw[red,ultra thick][->,>= {Stealth}] (white4) to (black4);
		\draw[red,ultra thick][->,>= {Stealth}] (white5) to (black5);
		\begin{scope}[thick]
			\draw[blue][<-, >= {Stealth}] (white3) to [bend right=45](black2);	
			\draw[blue][<-,>= {Stealth}] (white2) to [bend right=45](black1);
			\draw[blue][<-,>= {Stealth}] (white1) to [bend right=45](black0);
			\draw[blue][<-,>= {Stealth}] (white5) to [bend right=45](black4);	
			\draw[blue][<-,>= {Stealth}] (white4) to [bend right=45](black3);	
			\draw[blue][<-,>= {Stealth}] (white0) to [bend right=45](black5);	
		\end{scope}
	\end{tikzpicture} 
\end{equation}
\end{example}
 
\begin{rem}\label{rem:twoz-z}
	Given a nonoriented edge of $G$ we have precisely two zig-zag paths containing that edge: one is obtained by flowing along the edge from `white' to `black' and the other by reversing the flow. For example, reversing the orientation in the above example, that is, starting with the edge
	$$
	(2'\to 2)
	$$
	the corresponding zig-zag path is as follows.
	
	\begin{equation}\label{zz2}
	\begin{tikzpicture}
		[place/.style={circle,draw=black,thick},
		transition/.style={circle,draw=black,fill=black}]
		\node (white5) at (90:1.5cm) [place] [label={below:$5$}] {};
		\node (black5) at (90:3cm) [transition] [label={above:$5'$}]{};
		\node (white4) at (150:1.5cm) [place] [label={below right:$4$}]{};
		\node (black4) at (150:3cm) [transition] [label={above left:$4'$}]{};
		\node (white3) at (210:1.5cm) [place] [label={above right:$3$}]{};
		\node (black3) at (210:3cm) [transition] [label={below left:$3'$}]{};
		\node (white2) at (270:1.5cm) [place] [label={above:$2$}] {};
		\node (black2) at (270:3cm) [transition] [label={below:$2'$}]{};
		\node (white1) at (330:1.5cm) [place] [label={above left:$1$}] {};
		\node (black1) at (330:3cm) [transition] [label={below right:$1'$}]{};
		\node (white0) at (30:1.5cm) [place] [label={below left:$0$}] {};
		\node (black0) at (30:3cm) [transition] [label={above right:$0'$}]{};
		\begin{scope}[gray,ultra thin]
		\draw (white5) to [bend right=30] (white4);
		\draw (white4) to [bend right=30](white3);
		\draw (white3) to [bend right=30](white2);
		\draw (white2) to [bend right=30](white1);
		\draw (white1) to [bend right=30](white0);
		\draw (white0) to [bend right=30](white5);
		\draw (black5) to [bend right=30](black4);
		\draw (black4) to [bend right=30](black3);
		\draw (black3) to [bend right=30](black2);
		\draw (black2) to [bend right=30](black1);
		\draw (black1) to [bend right=30](black0);
		\draw (black0) to [bend right=30](black5);
		\end{scope}
		\draw[red,ultra thick][<-,>={Stealth}] (white0) to (black0);
		\draw[red,ultra thick][<-,>={Stealth}] (white1) to (black1);
		\draw[red,ultra thick][<-,>={Stealth}] (white2) to (black2);
		\draw[red,ultra thick] [<-,>={Stealth}](white3) to (black3);
		\draw[red,ultra thick][<-,>={Stealth}] (white4) to (black4);
		\draw[red,ultra thick][<-,>={Stealth}] (white5) to (black5);
		\begin{scope}[decoration={markings,
			mark =at position 1cm with {\arrow[red,line width=1mm]{stealth}}}]
		\draw[red, dotted,thick][postaction={decorate}] (white2) to [bend left=45](black3);
		\draw[red,dotted,thick][postaction={decorate}] (white1) to [bend left=45](black2);
		\draw[red,dotted,thick][postaction={decorate}] (white0) to [bend left=45](black1);
		\draw[red,dotted,thick][postaction={decorate}] (white3) to [bend left=45](black4);
		\draw[red,dotted,thick][postaction={decorate}] (white4) to [bend left=45](black5);
		\draw[red,dotted,thick][postaction={decorate}] (white5) to [bend left=45](black0);
		\end{scope}
	\end{tikzpicture}
\end{equation}

Pursuing this example, let us reverse the direction of one of the dotted edges, say take the edge
$$
5' \to 4.
$$
The corresponding zig-zag path is as follows

\begin{equation}\label{zz-3}
\begin{tikzpicture}
	[place/.style={circle,draw=black,thick},
	transition/.style={circle,draw=black,fill=black}]
	\node (white5) at (90:1.5cm) [place] [label={below:$5$}] {};
	\node (black5) at (90:3cm) [transition] [label={above:$5'$}]{};
	\node (white4) at (150:1.5cm) [place] [label={below right:$4$}]{};
	\node (black4) at (150:3cm) [transition] [label={above left:$4'$}]{};
	\node (white3) at (210:1.5cm) [place] [label={above right:$3$}]{};
	\node (black3) at (210:3cm) [transition] [label={below left:$3'$}]{};
	\node (white2) at (270:1.5cm) [place] [label={above:$2$}] {};
	\node (black2) at (270:3cm) [transition] [label={below:$2'$}]{};
	\node (white1) at (330:1.5cm) [place] [label={above left:$1$}] {};
	\node (black1) at (330:3cm) [transition] [label={below right:$1'$}]{};
	\node (white0) at (30:1.5cm) [place] [label={below left:$0$}] {};
	\node (black0) at (30:3cm) [transition] [label={above right:$0'$}]{};
	\begin{scope}[gray,ultra thin]
	\draw (white5) to [bend right=30] (white4);
	\draw (white4) to [bend right=30](white3);
	\draw (white3) to [bend right=30](white2);
	\draw (white2) to [bend right=30](white1);
	\draw (white1) to [bend right=30](white0);
	\draw (white0) to [bend right=30](white5);
	\draw (black5) to [bend right=30](black4);
	\draw (black4) to [bend right=30](black3);
	\draw (black3) to [bend right=30](black2);
	\draw (black2) to [bend right=30](black1);
	\draw (black1) to [bend right=30](black0);
	\draw (black0) to [bend right=30](black5);
	\end{scope}
	\begin{scope}
		[decoration={markings,
			mark =at position 2cm with {\arrow[blue,line width=1mm]{stealth}}}]
		\draw[blue,thick][postaction={decorate}] (white3) to [bend right=45](black2);	
		\draw[blue,thick][postaction={decorate}] (white2) to [bend right=45](black1);
		\draw[blue,thick][postaction={decorate}] (white1) to [bend right=45](black0);
		\draw[blue,thick][postaction={decorate}](white5) to [bend right=45](black4);	
		\draw[blue,thick][postaction={decorate}] (white4) to [bend right=45](black3);	
		\draw[blue,thick][postaction={decorate}] (white0) to [bend right=45](black5);	
	\end{scope}
	\begin{scope}[decoration={markings,
			mark =at position 2cm with {\arrow[red,line width=1mm]{stealth}}}]
		\draw[red, dotted,thick][postaction={decorate}] (black3) to [bend right=45](white2);
		\draw[red,dotted,thick][postaction={decorate}] (black2) to [bend right=45](white1);
		\draw[red,dotted,thick][postaction={decorate}] (black1) to [bend right=45](white0);
		\draw[red,dotted,thick][postaction={decorate}] (black4) to [bend right=45](white3);
		\draw[red,dotted,thick][postaction={decorate}] (black5) to [bend right=45](white4);
		\draw[red,dotted,thick][postaction={decorate}] (black0) to [bend right=45](white5);
	\end{scope}
\end{tikzpicture}
\end{equation}

We can see that it has two connected components: one passing through the even white vertices $0$, $2$, $4$ and the other through the odd ones $5$, $3$, $1$.
Thus there are $4$ distinct connected zig-zag paths in the graph $\widehat{PG}_{6}$, see \eqref{zz1}, \eqref{zz2}, \eqref{zz-3}. This is the {\rm complete list}, since every edge of the graph occurs precisely two times with the opposite orientations. We will show in a moment that this is a general pattern of the graphs $\widehat{PG}_{l}$ by exhibiting their zig-zag paths. 	 
\end{rem}
  
  To determine the zig-zag paths for $\widehat{PG}_{l}$ we label the edges
  with the tail at the white vertex $i$ as follows
  $$
   e^0_i =(i\to i'), \,\, e^-_i =(i\to (i-1)'), \,\, e^+_i =(i\to (i+1)').
   $$
   To indicate the change of orientation we multiply the above labels by $(-1)$. Thus the edges with the tail at the black vertex $i'$ are 
   $$
   -e^0_i=(i'\to i),\,\, -e^-_{i+1}= (i' \to (i+1)), \,\, -e^+_{i-1} =(i'\to (i-1)).
   $$
   We remind the reader that the standard ribbon graph structure is the counterclockwise orientation at all vertices. This means that at the vertex $i$ the cyclic order is
   $$
   \xymatrix{
   &e^0_i \ar[dl]& \\
    e^-_i \ar[rr]& & e^+_i \ar[ul]
}
   $$
   while at the vertex $i'$ we have
   $$
   \xymatrix{
   	&-e^0_i \ar[dl]& \\
   	-e^-_{i+1} \ar[rr]& & -e^+_{i-1} \ar[ul]
   }
   $$
   We can now give the list of zig-zag paths of $\widehat{PG}_{l}$.
   \begin{lem}\label{lem:zzpaths}
   	On the graph $\widehat{PG}_{l}$ there are three or four connected zig-zag paths, depending on the parity of $l$.
   	
   	If $l$ is odd, then there are three connected zig-zag paths:
   	$$
   	Z_{0,-} = (-e^-_0)e^0_{l-1}\cdots(-e^{-}_3)e^0_2 (-e^{-}_2)e^0_1(-e^{-}_1)e^0_0,
   	$$
   		$$
   	Z_{0,+} =e^{+}_{l-1}(-e^0_{l-1})\cdots e^{+}_2(-e^0_2)e^{+}_1(-e^0_1)e^{+}_0(-e^0_0 ) ,
   	$$
   	$$
   		Z_{-,+} =  (-e^+_1)e^-_3\cdots(-e^+_{l-4})e^-_{l-2}(-e^+_{l-2})e^-_0\cdots(-e^+_{l-5})e^-_{l-3}(-e^+_{l-3})e^{-}_{l-1}(-e^+_{l-1})e^{-}_1 .
   		$$
   		
   		If $l$ is even, then there are four connected zig-zag paths:
   		the paths $Z_{0,\pm}$ as above and the paths
   		$$
   		Z^0_{-,+}= (-e^+_0)e^-_2\cdots (-e^+_{l-4})e^-_{l-2}(-e^+_{l-2})e^-_0 ,
   		$$
   		$$
   		Z^1_{-,+}= (-e^+_1)e^-_3\cdots(-e^+_{l-3})e^-_{l-1}(-e^+_{l-1})e^-_1,
   		$$
   		each composed of edges $e^-_i$ and $(-e^+_j)$ and where in the first (resp. second) path all indices $i,j$ are even (resp. odd), that is, the first (resp. second) path goes only through the even (resp. odd) white vertices.   
   \end{lem}
\begin{pf}
	Follows immediately from the labeling of edges and the rules of turning at every vertex of the graph.
\end{pf}
Recall the exact sequence of homology groups
$$
\xymatrix{
	0\ar[r]&\ZZ\{B_i\}\ar[r]&H_1 (\widehat{PG}_{l},\ZZ) \ar[r]^{j}&H_1(\mathbb{T}, \ZZ)\ar[r]&0.
}
$$
It tells us that every closed path of $\widehat{PG}_{l}$ is mapped by the epimorphism $j$ of the sequence to a homology class of the torus. For a closed path $\gamma$ of $\widehat{PG}_{l}$  denote $j(\gamma)$ by $[\gamma]$ and call it the homology class of $\gamma$. To determine the homology classes of the zig-zag paths of  $\widehat{PG}_{l}$ we choose a basis of the homology group of the torus $\mathbb{T}$ as indicated on the following drawing.
\begin{equation}\label{torus}
\begin{tikzpicture}
	\draw[black, thick] (0,0) circle [radius=2cm];
	\draw[black, thick] (0,0) circle [radius=1cm];
	\begin{scope}
	[decoration={markings,
		mark =at position 2.5cm with {\arrow[green,line width=0.5mm]{stealth}}}]
		\draw[green,thick][postaction={decorate}] (0,0) circle [radius=1.5cm];
		\end{scope}
	\begin{scope}
		[decoration={markings,
			mark =at position 0.5cm with {\arrow[red,line width=0.5mm]{stealth}}}]
	\draw[red,thick] [postaction={decorate}](-2,0)
	arc[start angle=180, end angle=360, x radius =0.5cm, y radius =0.2cm];
		\draw[red,dotted, thick][postaction={decorate}] (-1,0)
	arc[start angle=0, end angle=180, x radius =0.5cm, y radius =0.2cm];
	\end{scope}
		\draw (1,0)
	arc[start angle=180, end angle=360, x radius =0.5cm, y radius =0.2cm];
	\draw[dotted,thick] (2,0)
	arc[start angle=0, end angle=180, x radius =0.5cm, y radius =0.2cm];
\end{tikzpicture}
\end{equation}
In the identification
$$
H_1 (\mathbb{T}, \ZZ) \cong \ZZ^2,
$$
the class of the equatorial red circle corresponds to $(-1,0)$ and the class of the meridian green circle is $(0,1)$. 
In this basis the homology classes of the zig-zag paths of the graph $\widehat{PG}_{l}$ are as follows.
\begin{lem}\label{lem:zz-homology}
	The homology classes of the zig-zag paths of $\widehat{PG}_{l}$ are as follows.
	
	$$
	[Z_{0,-}]=(0,-1), \,\, 	[Z_{0,+}]=(-l,-1).
	$$
	If $l$ is odd, then
	$$
	[Z_{-,+}] =(l,2).
	$$
	If $l$ is even, then
	$$
	[Z^0_{-,+}]=[Z^1_{-,+}] =\left(\frac{l}{2},1 \right).
	$$
\end{lem}
\begin{pf}
	The path $Z_{0,-}$ forms a simple closed curve on the `upper' annulus of the torus which turns once in the clockwise direction around the annulus. Hence it has the class of the meridian circle in \eqref{torus}
	with the orientation reversed, that is it has the class $(0,-1)$.
	
	For $Z_{0,+}$ observe that by moving two neighboring $0$-colored edges
	together gives the equatorial loop $(-1,0)$; since there are $l$ of those edges in $Z_{0,+}$ it has multiplicity $l$ for $(-1,0)$. The $+$-colored edges make a complete turn in the clockwise direction around the annulus. Hence the homology class $[Z_{0,+}]$ of $Z_{0,+}$:
	$$
	[Z_{0,+}]=l(-1,0)+(0,-1)=(-l,-1).
	$$
	
	The path $Z_{-,+}$ turns twice around the annulus counterclockwise: once following $-$-colored edges and the second time along $+$-colored edges; while turning around the annuli, the `upper' and the `lower' ones, the path passes $l$ times from the upper to the lower annulus, that is it makes $l$ turns around the equator in \eqref{torus}. Hence the homology class $[Z_{-,+}]$ has the form
	$$
	[Z_{-,+}]=l(1,0)+2(0,1)=(l,2).
	$$
	If $l$ is odd, the path $Z_{-,+}$ is connected, and if $l$ is even it
	falls into two connected components: 
	
	- $Z^0_{-,+}$ passing only through the white vertices with even labels,
	
	- $Z^1_{-,+}$ passing only through the white vertices with odd labels.
	
	We can slide one component into another. Hence they are homologous and we obtain the formula
	$$
	[Z^0_{-,+}]=[Z^1_{-,+}]=\HA(l,2)=\left(\frac{l}{2},1\right).
	$$
\end{pf}
We are now in the position to attach a toric surface to the graph $\widehat{PG}_{l}$.
\section{The toric surface of the graph $\widehat{PG}_{l}$ }
 The integral points of the lattice $\ZZ^2$ associated to zig-zag paths found in Lemma \ref{lem:zz-homology} are used to define a toric surface.
 Namely, we fix the counterclockwise cyclic ordering of vectors
 $$
 [Z_{0,+}]=(-l,-1), \,\,[Z_{0,-}]=(0,-1), \,\, 	v_{\mp}=\begin{cases}
 	(l,2),&\text{if $l$ is odd},\\
 	\left(\frac{l}{2},1 \right), & \text{if $l$ is even}.
 \end{cases}
$$
 In $\mathbb{R}^2=\ZZ^2 \otimes \mathbb{R}$ consider the fan
$$
\Delta_l =\{\sigma_+,\sigma_{-}, \sigma_{\mp} \}
$$
defined by the cones
$$
\sigma_+ =\mathbb{R}_{+} \langle [Z_{0,+}],[Z_{0,-}] \rangle, \,\,
\sigma_{-}=\mathbb{R}_{+} \langle [Z_{0,-}],v_{\mp} \rangle, \,\,
\sigma_{\mp}=\mathbb{R}_{+} \langle v_{\mp},[Z_{0,+}], \rangle.
$$
$$
\begin{tikzpicture}[>=stealth]
	\draw
	(-4,0)--(4,0);
	\draw
	(0,-2)--(0,2);
	\filldraw[gray]
	(0,0) circle [radius=2pt]
	(0,-1) circle [radius=2pt] node[right,color=red]{$[Z_{0,-}]$}
	(-3,-1) circle [radius=2pt] node[left,color=red]{$[Z_{0,+}]$}
	(3,1.5)circle [radius=2pt] node[above left,color=red]{$v_{\mp}$};
	\draw[black, very thick][->]
	(0,0)--(0,-1);
	\draw[black, very thick][->]
	(0,0)--(3,1.5);
	\draw[black, very thick][->]
	(0,0)--(-3,-1);
	\begin{pgfonlayer}{background}
		\fill[green!20!white]
		(0,-2)--(0,0)--(-4,-1.33);
	\end{pgfonlayer}
	\draw
	(-1,-1.5) node{$\sigma_{+}$};
	\begin{pgfonlayer}{background}
		\fill[yellow!20!white]
		(0,-2)--(0,0)--(4,2)--(4,-2);
	\end{pgfonlayer}
	\draw
	(3,-1.5) node{$\sigma_{-}$};
		\begin{pgfonlayer}{background}
		\fill[blue!20!white]
		(-4,-1.33)--(0,0)--(4,2)--(-4,2);
	\end{pgfonlayer}
	\draw
	(-1,1) node{$\sigma_{\mp}$};
\end{tikzpicture}
$$

\vspace{0.5cm}
We associate to $\widehat{PG}_{l}$ the toric variety $X(\Delta_l)$ defined by the fan $\Delta_l$.
\begin{pro}\label{pro:toricsurf}
	$X(\Delta_l)$ is a toric projective surface. If $l$ is odd
	we have a distinguished finite morphism
	$$
	m_l : \PP^2 \longrightarrow  X(\Delta_l)
	$$
	of degree $l$.
	
	If $l=2k$ is even, we have a distinguished finite morphism
	$$
	m_l : Q_2 \longrightarrow  X(\Delta_l)
	$$
	of degree $k$, where $Q_2$ is the singular toric surface obtained
	from the Hirzebruch surface $F_2=\PP(\OO_{\PP^1} \oplus \OO_{\PP^1} (-2) ) $ by contracting the section $E_0$ of $F_2$ of self-intersection $(-2)$. 
\end{pro}
\begin{pf}
	Since the union of three cones defining the fan $\Delta_l$ is $\RR^2=\ZZ^2 \otimes \RR$, it follows that $X(\Delta_l)$ is projective.
	
	Let $l$ be odd. Consider the sublattice $N$ of $\ZZ^2$ spanned by the vectors
	$$
	[Z_{0,+}]=(-l,-1), \,\,[Z_{0,-}]=(0,-1), \,\, 	v_{\mp}=
		(l,2).
		$$
	We denote the cones $\sigma_{\mp}$, $\sigma_{+}$ and $\sigma_{-}$  spanned by these vectors in the lattice $N$ by 	$\sigma'_{\mp}$, $\sigma'_{+}$ and $\sigma'_{-}$. The fan formed by these cones is denoted by $\Delta'_l$. The inclusion of lattices
	$$
	N \hookrightarrow \ZZ^2
	$$
	induces the identity map of fans $\Delta'_l \longrightarrow \Delta_l$ and hence the map of the corresponding toric surfaces
	\begin{equation}\label{morph-toric}
	X(\Delta'_l) \longrightarrow X(\Delta_l)
\end{equation}
	which is surjective. To see that $	X(\Delta'_l)$ is a projective plane, observe
 that the lattice $N$ is generated over $\ZZ$ by
		$$
		e'_1=le_1 =l(1,0) \,\text{and} \,e'_2=e_2 =(0, -1).
		$$
		Furthermore, the vectors generating the cones
		$\sigma'_{\mp}$, $\sigma'_{+}$ and $\sigma'_{-}$ are
		$$
		v_{\mp} =e'_1 -2e'_2,\,\, [Z_{0,+}]=-e'_1 +e'_2,\,\, [Z_{0,-}]=e'_2.
		$$
		The vectors generating the cones $\sigma'_{\mp}$, $\sigma'_{+}$ and $\sigma'_{-}$ are now bases in $N$. So the toric surface
		$	X(\Delta'_l)$ is smooth. In addition, the sum of three vectors is zero. It is well-known that this is a projective plane. 
		
		The lattice $N$, by definition, is the sublattice of index $l$ in $\ZZ^2$. This means that the morphism in \eqref{morph-toric} is of degree $l$.
		
		In the case $l=2k$ is even we argue similarly, the only difference that we consider the sublattice $N$ of $\ZZ^2$
		generated by the vectors
		$$
		v_1 = (-k,-1), \, [Z_{0,-}]=(0,-1), \, v_{\mp}=(k,1).
		$$
		The vectors
		$$
		e'_1=k(1,0), \, \text{and} \, e'_2=(0,-1)
		$$
		form a $\ZZ$-basis of $N$ and the vectors generating the cones of
		$\Delta'_l$ are written in this basis as follows
		$$
		v_{\mp}=(k,1)=e'_1 -e'_2,\, [Z_{0,+}]=(-l,-1)=(-2k,-1)=-2e'_1 +e'_2, \,  [Z_{0,-}]=e'_2.
		$$
		In this case the cones $\sigma'_{\mp}$ and $\sigma'_{-}$ are still generated by bases in $N$, while the cone  $\sigma'_{+}$ is generated by $(-2e'_1 +e'_2 )$ and $e'_2$ and this is not a basis of $N$. The affine open patch
		$U_{\sigma_{+}}$ is singular. It is easy to see that $U_{\sigma_{+}}$ is the quadric surface
		$$
		UV=Y^2
		$$
		in $\CC^3$. Thus the toric surface
		$X(\Delta'_l)$ has a unique ordinary double point and we have a finite morphism
		$$
		X(\Delta'_l) \longrightarrow X(\Delta_l)
		$$
		of degree $k$, since the index of $N$ in $\ZZ^2$ is $k$.
		
		To resolve the singularity of  $X(\Delta'_l) $ we refine the fan
		$\Delta'_l$ by inserting the vector
		$v_0 =(-1, 1)$ between $v_{\mp}$ and $e'_2$. It is again a standard fact that the resulting fan gives the Hirzebruch surface $F_2$.
\end{pf}
\section{The dual quiver attached to $\widehat{PG}_{l}$}
It is well known that a dimer model gives rise in a natural way to a quiver.
We briefly remind the construction and then determine the quiver associated to $\widehat{PG}_{l}$.

Let $G=(V_G,E_G, F_G)$ be a dimer model, where $V_G$ (resp. $E_G$ and $F_G$) is the set of vertices (resp. edges and faces) of $G$. The quiver
of $G$, denoted $Q_G$ and called {\it dual quiver of $G$}, is the {\it dual} graph $\check{G}$ of $G$, that is,
the vertices $V_{Q_G}$ of $Q_G$ are in bijection with the set $F_G$ of faces of $G$, while the faces $F_{Q_G}$ and edges $E_{Q_G}$ are in bijection with the set $V_G$ and $E_G$ respectively. More precisely, we identify the vertices $V_{Q_G}$ of $Q_G$ with the barycenters of faces of $G$:
$$
F_G \ni F \mapsto \check{F} \in V_{Q_G},
$$
where $\check{F}$ is the barycenter of $F$; two barycenters $\check{F}$ and $\check{F'}$ are connected by an edge if and only if the faces $F$ and $F'$ are adjacent: if $e$ is a common edge of $F$ and $F'$, the dual edge $\check{e}$ of $e$ is drawn as a segment in $\mathbb{T}$ connecting the barycenters of $F$ and $F'$ and passing through the midpoint of $e$; the orientation of $\check{e}$ is chosen so that it crosses $e$ with black (resp., white) vertex on the left (resp. right) side; for every vertex $v$ of $G$ we have a cyclic ordering of edges of $G$ incident to $v$; following the corresponding dual edges around $v$ gives the boundary of the dual face $\check{v}$ of $Q_G$. The bipartite nature of $G$ gives rise to the coloring of faces of $Q_G$: for a black (resp. white) vertex $v_b$ (resp. $v_w$) the dual face
$\check{v}_b$ (resp. $\check{v}_w$) is colored black (resp. white).
As an example of all of the above we will work out the dual quiver 
$Q_l:=Q_{\widehat{PG}_{l}}$ of the graph $\widehat{PG}_{l}$.

From Lemma \ref{lem:bound-cycles} we know that the graph $\widehat{PG}_{l}$ has $l$ faces $\{F_0, F_1, \ldots, F_{l-1}\}$ where the face $F_i$ is as follows.  

\begin{equation}\label{Fi}
  	\begin{tikzpicture}
  	[place/.style={circle,draw=black,thick, inner sep=0pt, minimum size=2mm},
  	transition/.style={circle,draw=black,fill=black, inner sep=0pt, minimum size=2mm}]
  	\node (whitei) at (150:1.5cm) [place] [label={left:$i$}] {};
  	\node (blacki) at (210:1.5cm) [transition] [label={left:$i'$}]{};
  	\node (whitei1) at (270:1.5cm) [place] [label={below:$i-1$}]{};
  	\node (blacki2) at (330:1.5cm) [transition] [label={right:$(i-2)'$}]{};
  	\node (whitei2) at (30:1.5cm) [place] [label={ right:$i-2$}]{};
  	\node (blacki1) at (90:1.5cm) [transition] [label={above:$(i-1)'$}]{};
  	\node at (0,0) [circle, draw=red, fill=red, inner sep=0pt, minimum size=2mm] [label={below:$F_i$}]{};
  	\draw[black,ultra thick][-] (whitei) to (blacki);
  	\draw[black,ultra thick][-] (whitei1) to (blacki);
  	\draw[black,ultra thick][-] (whitei2) to (blacki2);
  	\draw[black,ultra thick][-] (whitei2) to (blacki1);
  	\draw[black,ultra thick][-] (whitei) to (blacki1);
  	\draw[black,ultra thick][-] (whitei1) to (blacki2);
  \end{tikzpicture} 
\end{equation} 
The red dot is the barycenter of $F_i$; this is the vertex $\check{F}_i$ of the quiver $Q_l$ dual to to the face $F_i$. Thus the set of vertices of  $Q_l$ is the set
$$
\{\check{F}_0,\check{F}_1, \ldots, \check{F}_{l-1} \}.
$$ 
The boundary of each face is a hexagon. Hence it is adjacent to six faces and this means that every vertex of the quiver $Q_l$ has degree $6$. From Lemma \ref{lem:bound-cycles} we obtain the following.
\begin{lem}\label{lem:adjFi}
	The face $F_i$ is adjacent to
	
	1) the face $F_{i-1}$ along two edges:
	$$
	\text{\rm $e^-_{i-1}=(i-1)\to (i-2)'$ and $e^+_{i-2}=(i-2)\to (i-1)'$;}
	$$
	
	2) the face $F_{i+1}$ along two edges:
	$$
	\text{\rm $e^-_{i}=i\to (i-1)'$ and $e^+_{i-1}=(i-1)\to i'$;}
$$	

	3) the face $F_{i-2}$ along the edge $e^0_{i-2}=(i-2)\to (i-2)'$;
	
	\vspace{0.2cm}
	4) the face $F_{i+2}$ along the edge $e^0_{i}=i\to i'$.
\end{lem}

The above give the following configuration of the edges of $Q_l$ at the vertex $\check{F}_i$
\begin{equation}\label{Fi-dualedges}
	\begin{tikzpicture}
		[place/.style={circle,draw=black,thick, inner sep=0pt, minimum size=2mm},
		transition/.style={circle,draw=black,fill=black, inner sep=0pt, minimum size=2mm}]
		\node (whitei) at (150:1.5cm) [place] [label={left:$i$}] {};
		\node (blacki) at (210:1.5cm) [transition] [label={left:$i'$}]{};
		\node (whitei1) at (270:1.5cm) [place] [label={below:$i-1$}]{};
		\node (blacki2) at (330:1.5cm) [transition] [label={right:$(i-2)'$}]{};
		\node (whitei2) at (30:1.5cm) [place] [label={ right:$i-2$}]{};
		\node (blacki1) at (90:1.5cm) [transition] [label={above:$(i-1)'$}]{};
		\node (redFi) at (0,0) [circle, draw=red, fill=red, inner sep=0pt, minimum size=2mm] {\small$\check{F}_i$};
		\node (redFi11)at (120:2.55cm) [circle, draw=red, fill=red, inner sep=0pt, minimum size=2mm] {\tiny$\check{F}_{i+1}$};
			\node (redFi12)at (240:2.55cm) [circle, draw=red, fill=red, inner sep=0pt, minimum size=2mm] {$\scriptscriptstyle{\check{F}_{i+1}}$};
			\node (redFi11-)at (300:2.55cm) [circle, draw=red, fill=red, inner sep=0pt, minimum size=2mm] {$\scriptscriptstyle{\check{F}_{i-1}}$};
			\node (redFi12-)at (60:2.55cm) [circle, draw=red, fill=red, inner sep=0pt, minimum size=2mm] {\tiny$\check{F}_{i-1}$};
			\node (redFi2+)at (180:2.55cm) [circle, draw=red, fill=red, inner sep=0pt, minimum size=2mm] {\tiny$\check{F}_{i+2}$};
			\node (redFi2-)at (0:2.55cm) [circle, draw=red, fill=red, inner sep=0pt, minimum size=2mm] {\tiny$\check{F}_{i-2}$};
			\begin{scope}
					[decoration={markings,
					mark =at position 1.25cm with {\arrow[blue,line width=0.7mm]{stealth}}}]
				
					\draw[blue,ultra thick][postaction={decorate}] (redFi11) to (redFi);
					\draw [blue,ultra thick][postaction={decorate}](redFi12) to (redFi);
					\draw  [blue,ultra thick][postaction={decorate}](redFi) to (redFi11-);
					\draw  [blue,ultra thick][postaction={decorate}](redFi) to (redFi12-);
					\draw  [blue,ultra thick][postaction={decorate}]  (redFi) to (redFi2+);
					\draw [blue,ultra thick][postaction={decorate}]  (redFi2-) to (redFi);
				\end{scope}
		\draw[gray,thin][-] (whitei) to (blacki);
		\draw[gray, thin][-] (whitei1) to (blacki);
		\draw[gray, thin][-] (whitei2) to (blacki2);
		\draw[gray, thin][-] (whitei2) to (blacki1);
		\draw[gray,thin][-] (whitei) to (blacki1);
		\draw[gray, thin][-] (whitei1) to (blacki2);
	\end{tikzpicture} 
\end{equation} 
From this we conclude the following.
\begin{pro}\label{pro:Ql}
	The quiver $Q_l$ has $l$ vertices 
	$$
	\{\check{F}_0,\check{F}_1, \ldots, \check{F}_{l-1} \}.
	$$ 
	Ordered cyclically, counterclockwise from  $\check{F}_{l-1}$ to $\check{F}_0$, they form a `double' $l$-gon, that is, each side of the $l$-gon is  doubled; in addition, for every $i\in [0,l-1]$, the diagonals connecting the vertex $\check{F}_i$ to $\check{F}_{i+2}$ are drawn; the portion of $Q_l$ around the vertex $F_i$ is as follows.
	
	$$
	\begin{tikzpicture}
		[place/.style={circle,draw=red,thick, fill=red, inner sep=0pt, minimum size=2mm}]
			\node (redFi+2) [place] at (75:2cm)  [label={$\check{F}_{i+2}$}]{};
			\node (redFi+1) [place] at (120:2cm)  [label={above left: $\check{F}_{i+1}$}]{};
				\node (redFi) [place] at (165:2cm)  [label={left:$\check{F}_{i}$}]{};
		\node (redFi-1) [place] at (210:2cm)  [label={below left: $\check{F}_{i-1}$}]{};
			\node (redFi-2) [place] at (255:2cm)  [label={below right:$\check{F}_{i-2}$}]{};
			\begin{scope}
					[decoration={markings,
					mark =at position 0.7cm with {\arrow[blue,line width=0.5mm]{stealth}}}]
				\draw[blue, thick][postaction={decorate}] 
				(redFi+2) to (redFi+1);
				\draw[blue, thick][postaction={decorate}] 
				(redFi+2) to [bend right=20](redFi+1);
				\draw[blue, thick][postaction={decorate}] (redFi+1) to (redFi);
				\draw[blue, thick][postaction={decorate}] 
				(redFi+1) to [bend right=20](redFi);
				\draw[blue, thick][postaction={decorate}] 
				(redFi) to (redFi-1);
				\draw[blue, thick][postaction={decorate}] 
				(redFi) to [bend right=20](redFi-1);
				\draw[blue, thick][postaction={decorate}] 
				(redFi-1) to (redFi-2);
				\draw[blue, thick][postaction={decorate}] 
				(redFi-1) to [bend right=20](redFi-2);
			\end{scope}	
			\begin{scope}
			[decoration={markings,
				mark =at position 1.5cm with {\arrow[blue,line width=0.5mm]{stealth}}}]
				\draw[blue, thick][postaction={decorate}] 
			(redFi) to (redFi+2);
				\draw[blue, thick][postaction={decorate}] 
			(redFi-2) to (redFi);
			\end{scope}
			\draw[blue, dotted] (255:2cm)
			arc[start angle=255, end angle=435, radius=2cm] ;
	\end{tikzpicture}
	$$
	\end{pro} 
	\begin{example}\label{ex:6quivers}
		We illustrate the quivers $Q_l$ for the first six values of $l$.
	
		$$
	\begin{tikzpicture}
		[place/.style={circle,draw=red,thick, fill=red, inner sep=0pt, minimum size=2mm}]
		\node at (-3,1.1) [label={$l=1$}]{};
		\node (redF0) [place] at (-3,0)  [label={below:$0$}]{};
		\begin{scope}
			[decoration={markings,
				mark =at position 0.7cm with {\arrow[blue,line width=0.5mm]{stealth}}}]
			\draw[blue, thick][postaction={decorate}] 
			(redF0) to [bend right=30](-3,1) ;
			\draw[blue, thick][postaction={decorate}] 
			(-3,1) to [bend right=30](redF0) ;
			\draw[blue, thick][postaction={decorate}] 
			(redF0) to [bend right=30](-3.9,-0.5) ;
			\draw[blue, thick][postaction={decorate}] 
			(-3.9,-0.5) to [bend right=30](redF0) ;
			\draw[blue, thick][postaction={decorate}] 
			(redF0) to [bend right=30](-2.1,-0.5) ;
			\draw[blue, thick][postaction={decorate}] 
			(-2.1,-0.5) to [bend right=30](redF0) ;
		\end{scope}
		\node at (1,1.1) [label={$l=2$}]{};
	\node (redF1-2) [place] at (0,0)  [label={below:$1$}]{};
	\node (redF0-2) [place] at (2,0)  [label={below:$0$}]{};
	\begin{scope}
		[decoration={markings,
			mark =at position 1cm with {\arrow[blue,line width=0.5mm]{stealth}}}]
		\draw[blue, thick][postaction={decorate}] 
		(redF1-2) to [bend right=20](redF0-2) ;
		\draw[blue, thick][postaction={decorate}] 
		(redF1-2) to [bend right=45](redF0-2);
		\draw[blue, thick][postaction={decorate}] 
		(redF0-2) to [bend right=20](redF1-2) ;
		\draw[blue, thick][postaction={decorate}] 
		(redF0-2) to [bend right=45](redF1-2) ;
	\end{scope}
\begin{scope}
		[decoration={markings,
			mark =at position 0.5cm  with {\arrow{stealth}},
		mark=at position 1cm  with {\arrow{stealth}},
mark=at position 1.5cm  with {\arrow{stealth}}}]
		\draw[blue, thick][postaction={decorate}] 
		(redF1-2) arc[start angle=0, end angle=360, radius=0.3cm];
			\draw[blue, thick][postaction={decorate}] 
		(redF0-2) arc[start angle=-180, end angle=180, radius=0.3cm];
	\end{scope}
\node at (5,1.1) [label={$l=3$}]{};
\node (redF2-3) [place] at (4,0)  [label={below:$2$}]{};
\node (redF1-3) [place] at (6,0)  [label={below:$1$}]{};
\node (redF0-3) [place] at (5,1)  [label={below:$0$}]{};
\begin{scope}
	[decoration={markings,
		mark =at position 0.7cm with {\arrow[blue,line width=0.5mm]{stealth}}}]
	\draw[blue, thick][postaction={decorate}] 
	(redF2-3) to [bend right=20](redF1-3) ;
	\draw[blue, thick][postaction={decorate}] 
	(redF2-3) to (redF1-3) ;
		\draw[blue, thick][postaction={decorate}] 
	(redF2-3) to [bend left=20](redF1-3) ;
	\draw[blue, thick][postaction={decorate}] 
	(redF1-3) to [bend right=20](redF0-3);
	\draw[blue, thick][postaction={decorate}] 
	(redF1-3) to (redF0-3) ;
	\draw[blue, thick][postaction={decorate}] 
	(redF1-3) to [bend left=20](redF0-3) ;
	\draw[blue, thick][postaction={decorate}] 
	(redF0-3) to [bend right=20](redF2-3) ;
	\draw[blue, thick][postaction={decorate}] 
	(redF0-3) to (redF2-3) ;
		\draw[blue, thick][postaction={decorate}] 
	(redF0-3) to [bend left=20](redF2-3) ;
\end{scope}
\node at (-3,-2.5) [label={$l=4$}]{};
\node (redF3-4) [place] at (-3.8,-3)  [label={left:$3$}]{};
\node (redF2-4) [place] at (-3.8,-4.6)  [label={below:$2$}]{};
\node (redF1-4) [place] at (-2.2,-4.6)  [label={below:$1$}]{};
\node (redF0-4) [place] at (-2.2,-3)  [label={right:$0$}]{};
\begin{scope}
	[decoration={markings,
		mark =at position 0.7cm with {\arrow[blue,line width=0.5mm]{stealth}}}]
	\draw[blue, thick][postaction={decorate}] 
	(redF3-4) to [bend right=20](redF2-4) ;
	\draw[blue, thick][postaction={decorate}] 
	(redF3-4) to (redF2-4) ;
	\draw[blue, thick][postaction={decorate}] 
	(redF2-4) to [bend right=20](redF1-4) ;
	\draw[blue, thick][postaction={decorate}] 
	(redF2-4) to (redF1-4) ;
	\draw[blue, thick][postaction={decorate}] 
	(redF1-4) to [bend right=20](redF0-4);
	\draw[blue, thick][postaction={decorate}] 
	(redF1-4) to (redF0-4) ;
	\draw[blue, thick][postaction={decorate}] 
	(redF3-4) to (redF1-4) ;
	\draw[blue, thick][postaction={decorate}] 
	(redF2-4) to (redF0-4) ;
	\draw[blue, thick][postaction={decorate}] 
	(redF0-4) to [bend right=20](redF3-4) ;
	\draw[blue, thick][postaction={decorate}] 
	(redF0-4) to (redF3-4) ;
	\draw[blue, thick][postaction={decorate}] 
		(redF0-4) to [bend right=20](redF2-4) ;
			\draw[blue, thick][postaction={decorate}] 
		(redF1-4) to [bend right=20](redF3-4) ;
\end{scope}

\node at (1,-2.5) [label={$l=5$}]{};
\node (redF0-5) [place] at (1,-3)  [label={above left:$0$}]{};
\node (redF4-5) [place] at (0.1,-3.7) [label={left:$4$}]{};
\node (redF3-5) [place] at (0.4,-4.8)  [label={below left:$3$}]{};
\node (redF2-5) [place] at (1.6,-4.8)  [label={below right:$2$}]{};
\node (redF1-5) [place] at (1.9,-3.7)  [label={right:$1$}]{};

\begin{scope}
	[decoration={markings,
		mark =at position 0.7cm with {\arrow[blue,line width=0.5mm]{stealth}}}]
	\draw[blue, thick][postaction={decorate}] 
	(redF4-5) to [bend right=20](redF3-5) ;
	\draw[blue, thick][postaction={decorate}] 
	(redF4-5) to (redF3-5) ;
		\draw[blue, thick][postaction={decorate}] 
	(redF2-5) to (redF4-5) ;
	\draw[blue, thick][postaction={decorate}] 
	(redF3-5) to [bend right=20](redF2-5) ;
	\draw[blue, thick][postaction={decorate}] 
	(redF1-5) to (redF3-5) ;
	\draw[blue, thick][postaction={decorate}] 
	(redF3-5) to (redF2-5) ;
	\draw[blue, thick][postaction={decorate}] 
	(redF2-5) to [bend right=20](redF1-5);
	\draw[blue, thick][postaction={decorate}] 
	(redF2-5) to (redF1-5) ;
	\draw[blue, thick][postaction={decorate}] 
	(redF0-5) to (redF2-5) ;
	\draw[blue, thick][postaction={decorate}] 
	(redF1-5) to (redF0-5) ;
	\draw[blue, thick][postaction={decorate}] 
	(redF1-5) to [bend right=20](redF0-5) ;
	\draw[blue, thick][postaction={decorate}] 
	(redF0-5) to [bend right=20](redF4-5) ;
	\draw[blue, thick][postaction={decorate}] 
	(redF0-5) to (redF4-5) ;
	\draw[blue, thick][postaction={decorate}] 
	(redF3-5) to (redF0-5) ;
	\draw[blue, thick][postaction={decorate}] 
	(redF4-5) to (redF1-5) ;
\end{scope}

\node at (5,-2.5) [label={$l=6$}]{};
\node (redF0-6) [place] at (5,-3)  [label={above left:$0$}]{};
\node (redF5-6) [place] at (4.1,-3.5) [label={left:$5$}]{};
\node (redF4-6) [place] at (4.1,-4.5)  [label={below left:$4$}]{};
\node (redF3-6) [place] at (5,-5)  [label={below right:$3$}]{};
\node (redF2-6) [place] at (5.9,-4.5)  [label={right:$2$}]{};
\node (redF1-6) [place] at (5.9,-3.5)  [label={right:$1$}]{};
\begin{scope}
	[decoration={markings,
		mark =at position 0.7cm with {\arrow[blue,line width=0.5mm]{stealth}}}]
	\draw[blue, thick][postaction={decorate}] 
	(redF5-6) to [bend right=20](redF4-6) ;
	\draw[blue, thick][postaction={decorate}] 
	(redF5-6) to (redF4-6) ;
		\draw[blue, thick][postaction={decorate}] 
	(redF3-6) to (redF5-6) ;
	\draw[blue, thick][postaction={decorate}] 
	(redF4-6) to (redF3-6) ;
	\draw[blue, thick][postaction={decorate}] 
	(redF4-6) to [bend right=20](redF3-6) ;
		\draw[blue, thick][postaction={decorate}] 
	(redF2-6) to (redF4-6) ;
	\draw[blue, thick][postaction={decorate}] 
	(redF3-6) to (redF2-6) ;
	\draw[blue, thick][postaction={decorate}] 
	(redF3-6) to [bend right=20](redF2-6) ;
		\draw[blue, thick][postaction={decorate}] 
	(redF1-6) to (redF3-6) ;
	\draw[blue, thick][postaction={decorate}] 
	(redF2-6) to [bend right=20](redF1-6);
	\draw[blue, thick][postaction={decorate}] 
	(redF2-6) to (redF1-6) ;
		\draw[blue, thick][postaction={decorate}] 
	(redF0-6) to (redF2-6) ;
	\draw[blue, thick][postaction={decorate}] 
	(redF1-6) to (redF0-6) ;
	\draw[blue, thick][postaction={decorate}] 
	(redF1-6) to [bend right=20](redF0-6) ;
	\draw[blue, thick][postaction={decorate}] 
	(redF0-6) to [bend right=20](redF5-6) ;
	\draw[blue, thick][postaction={decorate}] 
	(redF0-6) to (redF5-6) ;
	\draw[blue, thick][postaction={decorate}] 
	(redF4-6) to (redF0-6) ;
	\draw[blue, thick][postaction={decorate}] 
	(redF5-6) to (redF1-6) ;
\end{scope}
	\end{tikzpicture}
	$$
	\end{example}

\vspace{0.2cm}
We now return to general considerations. Around every vertex $\check{F}_i$ we have six arrows of the quiver $Q_l$ which we agree to range counterclockwise starting from the upper edge connecting
$\check{F}_{i+1}$ to $\check{F}_i $, see \eqref{Fi-dualedges}. We label it $(\check{F}_{i+1} \to \check{F}_i )^1$. The following edges are labeled
$$
(\check{F}_{i} \to \check{F}_{i+2} ),\,\,(\check{F}_{i+1} \to \check{F}_i ) ^2,\,\, (\check{F}_{i} \to \check{F}_{i-1} ) ^1,\,\, (\check{F}_{i-2} \to \check{F}_i ), \,\, (\check{F}_{i} \to \check{F}_{i-1} ) ^2. 
$$

We know that that the quiver is the dual graph of $\widehat{PG}_{l}$, so the arrows of $Q_l$ are the edges dual to the edges of $\widehat{PG}_{l}$. We recall that those are labeled by $e^0_i$, $e^+_i$, $e^-_i$, for $i\in \{0,1,\ldots, l-1\}$. This gives the following correspondence
 \begin{equation}\label{dualedg}
 	(\check{F}_{i+1} \to \check{F}_i )^1 =(e^-_i)^{\check{}}, \,\, (\check{F}_{i} \to \check{F}_{i+2} )=(e^0_i)^{\check{}},\,\, (\check{F}_{i+1} \to \check{F}_{i} )^2=(e^+_{i-1})^{\check{}}.
 \end{equation}
With this the edges of $Q_l$ acquire `colors':
$$
 -,\, 0,\, +.
 $$
  The double edges of the $l$-gon in Proposition \ref{pro:Ql} are colored by $\pm$ and the proper diagonals by $0$. We agree that the outer edges of $l$-gon are labeled by $(e^-_i)^{\check{}}$'s and the inner  ones are
  $(e^+_i)^{\check{}}$'s:
  	$$
  \begin{tikzpicture}
  	[place/.style={circle,draw=red,thick, fill=red, inner sep=0pt, minimum size=2mm}]
  	\node (redFi+2) [place] at (72:2cm)  [label={$\check{F}_{i+2}$}]{};
  	\node (redFi+1) [place] at (144:2cm)  [label={above left: $\check{F}_{i+1}$}]{};
  	\node (redFi) [place] at (216:2cm)  [label={left:$\check{F}_{i}$}]{};
  	\node (redFi-1) [place] at (288:2cm)  [label={below left: $\check{F}_{i-1}$}]{};
  	\node (redFi-2) [place] at (360:2cm)  [label={below right:$\check{F}_{i-2}$}]{};
  	\node (Fi+2-) at (108:1.8cm) [label={$\scriptscriptstyle(e^-_{i+1})^{\check{}}$}]{};
  	\node (Fi+2+) at (108:1.7cm) [label={ below:$\scriptscriptstyle(e^+_{i})^{\check{}}$}]{};
  	\node (Fi+1-) at (180:1.8cm) [label=left:{$\scriptscriptstyle(e^-_{i})^{\check{}}$}]{};
  	\node (Fi+1+) at (180:1.7cm) [label={ right:$\scriptscriptstyle(e^+_{i-1})^{\check{}}$}]{};
  	\node (Fi-1+) at (252:1.8cm) [label={$\scriptscriptstyle(e^+_{i-2})^{\check{}}$}]{};
  	\node (Fi-1-) at (252:1.6cm) [label={ below left:$\scriptscriptstyle(e^-_{i-1})^{\check{}}$}]{};
  	\node (Fi-2-) at (324:1.6cm) [label=below right:{$\scriptscriptstyle (e^-_{i-2})^{\check{}}$}]{};
  	\node (Fi-2+) at (324:1.9cm) [label={above left:$\scriptscriptstyle (e^+_{i-1})^{\check{}}$}]{};
  	\begin{scope}
  		[decoration={markings,
  			mark =at position 0.7cm with {\arrow[blue,line width=0.5mm]{stealth}}}]
  		\draw[blue, thick][postaction={decorate}] 
  		(redFi+2) to (redFi+1);
  		\draw[blue, thick][postaction={decorate}] 
  		(redFi+2) to [bend right=20](redFi+1) ;
  		\draw[blue, thick][postaction={decorate}] (redFi+1) to (redFi);
  		\draw[blue, thick][postaction={decorate}] 
  		(redFi+1) to [bend right=20](redFi);
  		\draw[blue, thick][postaction={decorate}] 
  		(redFi) to (redFi-1);
  		\draw[blue, thick][postaction={decorate}] 
  		(redFi) to [bend right=20](redFi-1);
  		\draw[blue, thick][postaction={decorate}] 
  		(redFi-1) to (redFi-2);
  		\draw[blue, thick][postaction={decorate}] 
  		(redFi-1) to [bend right=20](redFi-2);
  	\end{scope}	
  	\draw[blue, dotted] (360:2cm)
  	arc[start angle=360, end angle=432, radius=2cm] ;
  \end{tikzpicture}
  $$

\begin{rem}\label{rem:Qlandsubquiv}
	The quiver $Q_l$ is the union of several subquivers: the quivers $Q^+_l$ and $Q^-_l$ formed by edges 
	 $(e^+_i)^{\check{}}$'s and $(e^-_i)^{\check{}}$'s respectively; and the quiver $Q^0_l$ formed by the edges $(e^0_i)^{\check{}}$'s; the first two are Dynkin quivers of type
	 $\widetilde{A}_{l-1}$  and the last one depends on the parity of $l$:
	 
	 - if $l$ is odd, then  $Q^0_l$ is also the Dynkin quiver $\widetilde{A}_{l-1}$, 
	 
	 - if $l$ is even, then   $Q^0_l$ falls into disjoint union of two quivers of type $\displaystyle \widetilde{A}_{\frac{l}{2}-1}$: one with vertices $\check{F}_{i}$'s, with $i$ even, and the other with $i$ odd.
	 
	 With the convention that the vertices of $Q_l$ are placed on the circle from $\check{F}_{l-1}$ to $\check{F}_{0}$ counterclockwise, the edges $(e^{\pm}_i)^{\check{}}$'s move counterclockwise, while the edges $(e^0_i)^{\check{}}$'s clockwise.  
\end{rem}

Next we turn to the faces of $Q_l$.    
They are dual to the vertices of the graph $\widehat{PG}_{l}$: these faces
 are triangles formed by arrows of $Q_l$ dual to the edges of $\widehat{PG}_{l}$ incident to a given vertex. Thus the white face $\check{i}$ dual to the vertex $i$ of $\widehat{PG}_{l}$ is the triangle
\begin{equation}\label{i-face}
\check{i}=(\check{F}_{i+2} \to \check{F}_{i+1} )^2(\check{F}_{i} \to \check{F}_{i+2} )(\check{F}_{i+1} \to \check{F}_i )^1 =(e^+_{i})^{\check{}}(e^0_i)^{\check{}}(e^-_i)^{\check{}}.
\end{equation}

The black face $\check{i'}$ dual to the vertex $i'$ is the triangle
\begin{equation}\label{i'-face}
\check{i'}=(\check{F}_{i+2} \to \check{F}_{i+1} )^1(\check{F}_{i} \to \check{F}_{i+2} )(\check{F}_{i+1} \to \check{F}_i )^2  =(e^-_{i+1})^{\check{}}(e^0_i)^{\check{}}(e^+_{i-1})^{\check{}}.
\end{equation}

We can now write down the potential of $Q_l$ and define its Jacobi algebra. Recall that with a quiver $Q$ one associates the path algebra $\CC Q$: it is the vector space over $\CC$ with a basis given by paths in $Q$, where a path $p$ in $Q$ is a finite sequence of arrows in $Q$
$$
p=a_n a_{n-1} \cdots a_1,
$$
where the head $h(a_i) $ is equal to the tale $t(a_{i+1})$, for every $i\in [1,n-1] $; those are paths of positive length; in addition, there are paths of length zero which correspond to the vertices of $Q$: for a vertex $v$ of $Q$, the corresponding path is denoted $e_v$. The multiplication in $\CC Q$ is defined by concatenation of paths: for $p=a_n a_{n-1} \cdots a_1$ and
$p'=a'_m a'_{m-1} \cdots a'_1$ one defines
$$
pp'=\begin{cases}
	a_n a_{n-1} \cdots a_1 a'_n a'_{n-1} \cdots a'_1,&\text{if $t(a_1)=h(a'_m)$},\\
		0,& \text{otherwise}. 
\end{cases}
$$
The paths of length zero are idempotents of $\CC Q$, that is, $e^2_v=e_v$, for every vertex $v$ of $Q$, and subject to
$$
ae_v=a \,(\text{resp. $e_v a=a$}),
$$
for every arrow $a$ starting (resp. ending) at $v$, and zero otherwise.

For the quiver $Q$ associated to a bipartite ribbon graph we have the collection of paths which are formed by the closed paths which are the boundaries of faces of $Q$, call them boundary paths of $Q$.  Remember that the faces of $Q$ are colored, black and white. One defines the potential $W$ of $Q$ as the signed sum of the  boundary paths of $Q$:
$$
W=\sum_{b} p_b - \sum_{w} p_w,
$$
where $p_b$ (resp. $p_w$) is the boundary path of a black (resp., white) face $b$ (resp. $w$) and the first (resp. second) sum is taken over the black (resp. white) faces of $Q$. As an example, let us work out the potential of our quiver $Q_l$ which will be denoted $W_l$. 

We know that white and black faces of $Q_l$ are $\{\check{i}\}_{i\in[0,l-1]}$ and $\{\check{i'}\}_{i\in[0,l-1]}$ respectively. From the equations \eqref{i-face} and \eqref{i'-face} we deduce the potential $W_l$ of $Q_l$:
$$
W_l = \sum^{l-1}_{i=0} (e^-_{i+1})^{\check{}}(e^0_i)^{\check{}}(e^+_{i-1})^{\check{}} - \sum^{l-1}_{i=0} (e^+_{i})^{\check{}}(e^0_i)^{\check{}}(e^-_i)^{\check{}}.
$$
One defines the Jacobian ideal $\partial(W_l)$ as the two sided ideal in $\CC Q_l$ generated by the `partial derivatives' of $W_l$. The reader may recall that that for a path $p=a_n a_{n-1} \cdots a_1$ and an arrow $a$ of a quiver $Q$ the partial derivative $\partial_a p$ is defined as follows
$$
\partial_a p =\sum_{a_i=a} a_{i-1} \cdots a_1a_n \cdots a_{i+1}.
$$
For example, for the black boundary path $p_{\check{i'}}=(e^-_{i+1})^{\check{}}(e^0_i)^{\check{}}(e^+_{i-1})^{\check{}}$ we have
$$
\partial_{(e^-_{i+1})^{\check{}}}  \,p_{\check{i'}} =(e^0_i)^{\check{}}(e^+_{i-1})^{\check{}}, \,\,\,\,
\partial_{(e^0_{i})^{\check{}}} \,p_{\check{i'}} =(e^+_{i-1})^{\check{}} (e^-_{i+1})^{\check{}},\,  \, \, \, 
\partial_{(e^+_{i-1})^{\check{}}} \,p_{\check{i'}} =(e^-_{i+1})^{\check{}} (e^0_{i})^{\check{}},
$$
and all other partial derivatives of $p_{\check{i'}}$ equal zero. 

We now compute the partials of the potential $W_l$.
\begin{lem}
	For every $i\in [0,l-1]$ we have
\begin{equation}\label{path-rel}
	\begin{gathered}
	\partial_{(e^0_{i})^{\check{}}} W_l =(e^+_{i-1})^{\check{}} (e^-_{i+1})^{\check{}} -(e^-_{i})^{\check{}} (e^+_{i})^{\check{}},\\
	\partial_{(e^+_{i})^{\check{}}} W_l =(e^-_{i+2})^{\check{}} (e^0_{i+1})^{\check{}} -(e^0_{i})^{\check{}} (e^-_{i})^{\check{}},\\
	\partial_{(e^-_{i})^{\check{}}} W_l =(e^0_{i-1})^{\check{}} (e^+_{i-2})^{\check{}} -(e^+_{i})^{\check{}} (e^0_{i})^{\check{}}.
	\end{gathered}
\end{equation}
Thus the Jacobian ideal $\partial W_l$ is the two-sided ideal generated by elements
$$
\{(e^+_{i-1})^{\check{}} (e^-_{i+1})^{\check{}} -(e^-_{i})^{\check{}} (e^+_{i})^{\check{}},\,\, (e^-_{i+2})^{\check{}} (e^0_{i+1})^{\check{}} -(e^0_{i})^{\check{}} (e^-_{i})^{\check{}},\,\, (e^0_{i-1})^{\check{}} (e^+_{i-2})^{\check{}} -(e^+_{i})^{\check{}} (e^0_{i})^{\check{}}\}_{i \in [0,l-1]}.
$$
\end{lem}

To the pair $(Q_l, W_l)$ one associates the quotient algebra
\begin{equation}\label{AQl}
	A_{Q_l}:=\CC Q_l /\partial W_l, 
\end{equation}
the Jacobi algebra of $(Q_l, W_l)$, where $\partial W_l$ is the Jacobian ideal of $\CC Q_l$, the two sided ideal of $\CC Q_l$ generated by the partial derivatives of $W_l$.  We denote by $a^0_i$ (resp. $a^{\pm}_i$) the image of $(e^0_i)^{\check{}}$ (resp. $(e^{\pm}_{i})^{\check{}})$) in $A_{Q_l}$ under the canonical projection
$$
\CC Q_l \longrightarrow A_{Q_l}.
$$
The subquivers $Q^{\pm}_l$ and $Q^0_l$ give rise to the subalgebras
$\CC Q^{\pm}_l$ and $\CC Q^0_l$ of the path algebra $\CC Q_l$. The images of these subalgebras under the canonical projection will be denoted
$$
A^+_{Q_l},\, A^-_{Q_l},\, A^0_{Q_l}
$$
respectively. From the definition of the Jacobian ideal $\partial W_l$ we deduce
$$
\partial W_l \bigcap \CC Q^{\pm}_l =\partial W_l \bigcap \CC Q^{0}_l =0.
$$
This immediately implies isomorphism
\begin{equation}\label{AisoQ}
	A^+_{Q_l} \cong \CC Q^{+}_l,\, 	A^-_{Q_l} \cong \CC Q^{-}_l,\,	A^0_{Q_l} \cong \CC Q^{0}_l.
\end{equation}

From the partials of $W_l$ in \eqref{path-rel} we deduce the following 	commutation relations in $A_{Q_l}$.
\begin{pro}\label{pro:AQl}
	The algebra $A_{Q_l}$ is generated by the orthogonal idempotents
	$$
	\{x_i | i=0, 1, \ldots, l-1\},
	$$
	corresponding to the vertices of the quiver (paths of length zero), together with the set of elements
	$$
	\{a^0_i, a^+_i, a^-_i | i=0,1,\ldots, l-1\}
	$$
	defined just after \eqref{AQl}. The latter generators are subject to the relations
	\begin{equation}\label{path-rel-A}
		\begin{gathered}
			a^+_{i-1} a^-_{i+1} = a^-_{i} a^+_{i},\\
			a^-_{i+2} a^0_{i+1} =a^0_{i} a^-_{i},\\
			a^0_{i-1} a^+_{i-2} = a^+_{i} a^0_{i}.
		\end{gathered}
	\end{equation}
\end{pro}

Pictorially, thinking of the circle model of $Q_l$, the first relation in \eqref{path-rel-A} says that moving from the vertex $\check{F}_{i+2}$ to $\check{F}_{i}$ along outer edge ($a^-_{i+1}$) followed by the inner one ($a^+_{i-1}$) is the same as moving first by the inner edge ($a^+_{i}$) and then the outer one ($a^-_{i}$). This is illustrated bellow:
	$$
\begin{tikzpicture}
	[place/.style={circle,draw=red,thick, fill=red, inner sep=0pt, minimum size=2mm}]
	\node (redFi+2) [place] at (72:2cm)  [label={$\check{F}_{i+2}$}]{};
	\node (redFi+1) [place] at (144:2cm)  [label={above left: $\check{F}_{i+1}$}]{};
	\node (redFi) [place] at (216:2cm)  [label={left:$\check{F}_{i}$}]{};
	\node (Fi+2-) at (108:1.8cm) [label={$a^-_{i+1}$}]{};
	\node (Fi+1+) at (180:1.7cm) [label={ left:$a^+_{i-1}$}]{};
	\begin{scope}
		[decoration={markings,
			mark =at position 0.7cm with {\arrow[blue,line width=0.5mm]{stealth}}}]
		\draw[blue, thick][postaction={decorate}] 
		(redFi+2) to [bend right=20](redFi+1) ;
		\draw[blue, very thick][postaction={decorate}] (redFi+1) to (redFi);
	\end{scope}	
\begin{scope}[xshift=7cm]
	[place/.style={circle,draw=red,thick, fill=red, inner sep=0pt, minimum size=2mm}]
	\node (redFi+2) [place] at (72:2cm)  [label={$\check{F}_{i+2}$}]{};
	\node (redFi+1) [place] at (144:2cm)  [label={above left: $\check{F}_{i+1}$}]{};
	\node (redFi) [place] at (216:2cm)  [label={left:$\check{F}_{i}$}]{};
		\node (Fi+2+) at (108:1.7cm) [label={ above:$a^+_{i}$}]{};
		\node (Fi+1-) at (180:1.8cm) [label=left:{$a^-_{i}$}]{};
\end{scope}
	\begin{scope}
		[decoration={markings,
			mark =at position 0.7cm with {\arrow[blue,line width=0.5mm]{stealth}}}]
			\draw[blue, thick][postaction={decorate}] 
			(redFi+2) to (redFi+1);
			\draw[blue, thick][postaction={decorate}] 
			(redFi+1) to [bend right=20](redFi);
	\end{scope}	
\begin{scope}  
[xshift=1.5cm, yshift=-2cm]	\node (equal) at (72:2cm)  [label={$\text{\Huge$=$}$}]{};
\end{scope}
\end{tikzpicture}
$$
Similarly, the other two relations in \eqref{path-rel-A} are depicted in the following two diagrams 
	$$
\begin{tikzpicture}
	[place/.style={circle,draw=red,thick, fill=red, inner sep=0pt, minimum size=2mm}]
	\node (redFi+3) [place] at (72:2cm)  [label={$\check{F}_{i+3}$}]{};
	\node (redFi+2) [place] at (144:2cm)  [label={above left: $\check{F}_{i+2}$}]{};
	\node (redFi+1) [place] at (216:2cm)  [label={left:$\check{F}_{i+1}$}]{};
	\node (Fi+3-) at (108:1.8cm) [label={$a^-_{i+2}$}]{};
	\node (Fi+1+) at (160:1.7cm) [label={ right:$a^0_{i+2}$}]{};
	\begin{scope}
		[decoration={markings,
			mark =at position 1.5cm with {\arrow[blue,line width=0.5mm]{stealth}}}]
		\draw[blue, thick][postaction={decorate}] 
		(redFi+3) to [bend right=20](redFi+2) ;
			\draw[blue, thick][postaction={decorate}] 
			(redFi+1) to (redFi+3);
	\end{scope}	
	\begin{scope}[xshift=7cm]
		[place/.style={circle,draw=red,thick, fill=red, inner sep=0pt, minimum size=2mm}]
		\node (redFi) [place] at (288:2cm)  [label={$\check{F}_{i}$}]{};
		\node (redFi+2) [place] at (144:2cm)  [label={above left: $\check{F}_{i+2}$}]{};
		\node (redFi+1) [place] at (216:2cm)  [label={left:$\check{F}_{i+1}$}]{};
		\node (Fi+2+) at (220:0.3cm) [label={ above:$a^0_{i}$}]{};
		\node (Fi+1-) at (260:1.8cm) [label={below left: $a^-_{i}$}]{};
	\end{scope}
	\begin{scope}
		[decoration={markings,
			mark =at position 1.7cm with {\arrow[blue,line width=0.5mm]{stealth}}}]
		\draw[blue, thick][postaction={decorate}] 
		(redFi) to (redFi+2);
		\draw[blue, thick][postaction={decorate}] 
		(redFi+1) to [bend right=20](redFi);
	\end{scope}	
	\begin{scope}  
		[xshift=1.5cm, yshift=-2cm]	\node (equal) at (72:2cm)  [label={$\text{\Huge$=$}$}]{};
	\end{scope}
\end{tikzpicture}
$$
$$
\begin{tikzpicture}
	[place/.style={circle,draw=red,thick, fill=red, inner sep=0pt, minimum size=2mm}]
	\node (redFi) [place] at (180:2cm)  [label={below left: $\check{F}_{i}$}]{};
	\node (redFi-1) [place] at (252:1cm)  [label={below: $\check{F}_{i-1}$}]{};
	\node (redFi+1) [place] at (108:2cm)  [label={above:$\check{F}_{i+1}$}]{};
	\node (Fi-) at (197:2.1cm) [label={right:$a^+_{i-2}$}]{};
	\node (Fi-1+) at (130:0.7cm) [label={ right:$a^0_{i-1}$}]{};
	\begin{scope}
		[decoration={markings,
			mark =at position 1.5cm with {\arrow[blue,line width=0.5mm]{stealth}}}]
			\draw[blue, thick][postaction={decorate}] 
			(redFi) to (redFi-1);
		\draw[blue, thick][postaction={decorate}] 
		(redFi-1) to (redFi+1);
	\end{scope}	
	\begin{scope}[xshift=7cm]
			[place/.style={circle,draw=red,thick, fill=red, inner sep=0pt, minimum size=2mm}]
		\node (redFi) [place] at (180:2cm)  [label={below left:$\check{F}_{i}$}]{};
		\node (redFi+1) [place] at (140:3cm)  [label={below left: $\check{F}_{i+1}$}]{};
		\node (redFi+2) [place] at (92:3cm)  [label={right:$\check{F}_{i+2}$}]{};
		\node (Fi-) at (140:1cm) [label={$a^0_{i}$}]{};
		\node (Fi+2+) at (109:2.8cm) [label={ left:$a^+_{i}$}]{};
	\end{scope}
	\begin{scope}
		[decoration={markings,
			mark =at position 1.7cm with {\arrow[blue,line width=0.5mm]{stealth}}}]
		\draw[blue, thick][postaction={decorate}] 
		(redFi) to (redFi+2);
		\draw[blue, thick][postaction={decorate}] 
		(redFi+2) to (redFi+1);
	\end{scope}	
	\begin{scope}  
		[xshift=1.5cm, yshift=-2cm]	\node (equal) at (72:2cm)  [label={$\text{\Huge$=$}$}]{};
	\end{scope}
\end{tikzpicture}
$$

The commutation relations \eqref{path-rel-A} tell us that the generators of different colors commute, up to shifting indices of generators. This implies the following.
\begin{pro}\label{pro:prodsub}
	The Jacobi algebra $A_{Q_l}$ is equal to the product of subalgebras
	$A^0_{Q_l}$, $A^+_{Q_l}$, $A^-_{Q_l}$:
	$$
	A_{Q_l}=A^0_{Q_l}A^+_{Q_l}A^-_{Q_l}.
	$$
	More precisely, given a path $p$ in $Q_l$ of positive length, there are elements
	$p^0$, $p^+$, $p^-$ belonging respectively to the subalgebras
	$A^0_{Q_l}$, $A^+_{Q_l}$, $A^-_{Q_l}$ such that the image $[p]$
	of $p$ under the projection morphism
	$$
	\CC Q_l \longrightarrow A_{Q_l}
	$$
	is written as the product
	$$
	[p]=p^0 p^+ p^-.
	$$
\end{pro} 

The Jacobi algebra $ A_{Q_l}$ is a noncommutative analogue of a three dimensional Calabi-Yau variety: it is a Calabi-Yau algebra of dimension three in a sense of Ginzburg, see \cite{Gi}, that is, the derived category $D^b(A_{Q_l})$ of finitely generated projective $A_{Q_l}$-bimodules is a Calabi-Yau category of dimension $3$: one has natural isomorphisms
$$
f_{M,N}:Hom_{D^b(A_{Q_l})} (M,N) \stackrel{\cong}{\longrightarrow} Hom_{D^b(A_{Q_l})} (N,M[3])^{\ast},
$$
for any pair $(M,N)$ of bounded complexes of finitely generated projective $A_{Q_l}$-bimodules; $M[n]$ is the shifted complex. This will not be used in the rest of the paper, but we record for future reference the following key property of $ A_{Q_l}$.
\begin{pro}\label{pro:CY-3}
	The Jacobi algebra $ A_{Q_l}$ is a Calabi-Yau algebra of dimension three.
\end{pro}
This follows from the works of N. Broomhead \cite{Br}, R. Bocklandt, \cite{Bo}: the authors give various consistency conditions on a dimer model implying that the Jacobi algebra of the quiver of a dimer is a Calabi-Yau algebra of dimension $3$. In our case $Q_l$ comes from the dimer model
$\widehat{PG}_l$ whose faces are hexagons, so called hexagonal tilings, and it is known, see \cite{Br}, that those satisfy the consistency conditions in the works above. 

It is clear that it is interesting to relate the IVHS refinement to the representation theory of the Jacobi algebra $A_{Q_l}$. As a step in this direction we will see shortly how to build natural representations of the path algebra $\CC Q_l$ from the representations of the quiver $PG_l$ parametrized by the stratum ${\mathfrak L}_l$. 
   
\section{Perfect matchings of $\widehat{PG}_{l}$}
  A perfect matching on a bipartite graph $G$ is a collection of edges of $G$ such that every vertex of $G$ is covered by precisely one edge of the collection, see \cite{Go-K}  and references therein for more details. The set of perfect matchings has many fascinating aspects: it is related to combinatorics, closed paths of the graph in question. In this section we begin exploring the perfect matchings of the graph $\widehat{PG}_l$. Let us start with some examples.
 \begin{example}\label{exa:match}
 	For the graph $\widehat{PG}_{l}$ the collection of edges
 	$$
 	M^0:=\{e^0_i \,|\,i=0,1,\ldots,l-1\}
 	$$
 	is a perfect matching. The same for the collections
 	$$
 	M^{\pm}:=\{e^{\pm}_i \,|\,i=0,1,\ldots,l-1\}.
 	$$
 \end{example}
The relation between the perfect matchings and closed paths of a bipartite ribbon graph $G$ can be seen as follows. Denote $V_G$, $E_G$ and $F_G$ the set of vertices, edges and faces of $G$ and form the complex
$$
\xymatrix{
0 \ar[r]& \ZZ\{F_G\} \ar[r]^{\partial}& \ZZ \{E_G\} \ar[r]^{\partial}& \ZZ \{V_G\} \ar[r]&0,
}
$$
where $\partial$ denotes the boundary differential. So for an edge $e$ the differential 
$$
\partial(e)=b(e)-w(e),
$$
where $b(e)$ (resp. $w(e)$) denotes the black (resp. white) vertex incident to $e$. 

Let $M$ be a perfect matching of $G$. We attach to it the chain
$$
e_M:=\sum_{e\in M} e
$$
in the free abelian group $\ZZ \{E_G\}$. Applying the boundary differential gives
$$
\partial (e_M)=\sum_{e\in M}\partial(e) = \sum_{e\in M} (b(e) -w(e)) =\sum_{b\in B} b -\sum_{w\in W } w,
$$
where $B$ (resp. $W$) stands for the set of black (resp. white) vertices of $G$. 
From this it follows
\begin{equation}\label{match-cycle}
\text{\it	$(e_M -e_{M'})$ is a cycle for any two matchings $M$ and $M'$ of $G$.} 
\end{equation}
In view of this we wish to understand the matchings of the graph $\widehat{PG}_{l}$. Let us recall that the edges of our graph are labeled
by $e^0_i$ and $e^{\pm}_i$, where $i \in [0,l-1]$. The superscripts $0, -,+$ will be called the {\it colors} of edges. We begin with matchings which have no edges of color $0$.
\begin{lem}\label{lem:matchings}
	1) If $l$ is odd, then there are precisely two matchings which have no edges of color $0$:
	$$
	e_{M^+}=\sum^{l-1}_{i=0} e^+_i, \,\, e_{M^-}=\sum^{l-1}_{i=0} e^-_i,
	$$
	where $M^{\pm}$ are the matchings described in Example \ref{exa:match}.
	
	2) If $l=2k$ is even, then in addition to the matchings $M^{\pm}$ in 1), there are two more matchings $M^{-,+} $ and $M^{+,-}$ whose chains are as follows:
	$$
	e_{M^{-,+}}=\sum^k_{s=1}(e^-_{2s-1} +e^+_{2(s-1)}), \,\,e_{M^{+,-}}=\sum^k_{s=1}(e^+_{2s-1} +e^-_{2(s-1)}).
	$$
\end{lem}
\begin{pf}
	Let $M$ be a matching of $\widehat{PG}_{l}$ without edges of color $0$.
	Assume we have two consecutive edges in $M$ of color $(-)$, that is, for some $i$, the edges $e^-_i$ and $e^-_{i-1}$ are in $M$. Then the edge at $(i-2)$ must be also  $e^-_{i-2}$. By induction, all edges are colored $(-)$ and this gives the matching $M^-$. The same argument holds if two consecutive edges are colored $(+)$ giving the matching $M^{+}$. The argument also implies that a matching $M$ which is not mono-colored must have consecutive edges of alternating colors. This can only occur if the number $l$ is even and in this case we have two more matchings $M^{-,+}$ and $M^{+,-}$ as stated in 2) of the lemma.  
\end{pf}

We now turn to the matchings with some edges colored $0$. Let $M$ be such a matching and arrange $0$-colored edges in decreasing order
$$
M(0):=\{e^0_{i_1}, \ldots, e^0_{i_r} \},
$$
that is, the indices $\{i_s\}_{s=1,\ldots,r}$ form a decreasing sequence.

\begin{pro}\label{pro:matchings0}
	Let $M$ be a matching with some edges colored $0$ and let
	$$
	M(0)=\{e^0_{i_1}, \ldots, e^0_{i_r} \} 
	$$
	the subset of $0$-colored edges of $M$ arranged in the decreasing order. Then the following occur.
	
	1) $M(0)=M=M^0$, the matching in Example \ref{exa:match}.
	
	\vspace{0.2cm}
	2) If $M \neq M^0$, then
	
	a) the number $r\equiv l\, mod (2)$;
	
	b) any two consecutive edges in $M(0)$ are separated by an even number of edges of $M$; if that number is nonzero then the separating edges  alternate in colors from $(-)$ to $(+)$.
\end{pro}
\begin{pf}
	The first assertion is obvious. For the second, part b) implies a).
	For b), let $e^0_{i_k}$ and $e^0_{i_{k+1}}$ be two consecutive edges in 
	$M(0)$ and assume that there are edges of $M$ in between, otherwise we are done. Then all edges in the gap are colored $+$ or $-$. The argument in the proof of Lemma \ref{lem:matchings} tells us that no two consecutive edges in the gap can have the same color, that is, the colors in the gap alternate. Furthermore, the edge after $e^0_{i_k}$ must color $(-)$ and the following switches to $(+)$ and so on until we arrive to $e^0_{i_{k+1}}$.  
\end{pf}
\begin{rem}\label{rem:match-subgones}
	We think of white vertices placed on a circle in the counterclockwise order from $(l-1)$ to $0$ and thus forming an $l$-gon. Denote that $l$-gon by $R_l$. It can be thought of as the matching $M^0$. A matching $M$
	with $M(0) =\{e^0_{i_1}, \ldots, e^0_{i_r} \}$ could be identified with the $r$-gon $R_{M(0)}$ having vertices 
	$$
	\{i_1,\ldots, i_r \}.
	$$
	The sides of  $R_{M(0)}$ are either the sides of $R_l$ or proper diagonals of $R_l$. In the latter case Proposition \ref{pro:matchings0} tells us that the polygon formed by a proper diagonal and the vertices of $R_l$ between the end points of the diagonal is even sided.
\end{rem}

We can now understand the matching cycles of $\widehat{PG}_{l}$. For this we set
\begin{equation}\label{set-matchings}
{\cal M}_{\widehat{PG}_{l}} =\text{the set of perfect matchings of $\widehat{PG}_{l}$.}
\end{equation}
Fix $M^0$ as a reference matching and consider the collection
\begin{equation}\label{matchcycles}
MC_{\widehat{PG}_{l}}:=\{ e_M- e_{M^0} | M\in {\cal M}_{\widehat{PG}_{l}}\}
\end{equation}
of matching cycles.

\begin{pro}\label{pro:matchcycles}
	1) The matching cycles $(e_{M^{\pm}}- e_{M^0})$ correspond to the zig-zag paths $Z_{0,+}$ and $(-Z_{0,-})$ respectively.
	
	\vspace{0.2cm}
	2) For a matching $M$ with $M(0) = \{e^0_{i_1}, \ldots, e^0_{i_r} \}$, a proper subset of $M^0$, denote $R_{M(0)}$ the $r$-gone formed by the vertices $i_1,\cdots, i_r$ and let 
	$$
	\Delta_M =\text{the set of sides of $R_{M(0)}$ which are proper diagonals of $R_{l}$ }.
	$$
	Then the matching cycle
	$$
	 e_M- e_{M^0} =\sum_{d\in \Delta_M } C_d 
	 $$
	 is the sum of cycles $C_d$ labeled by elements of $\Delta_M$; for $d=(i_s,i_{s+1}) $, a side of  $R_{M(0)}$ which is a proper diagonal of $R_l$, the cycle $C_d$ has the form
	 $$
	 C_d =\sum_{i_{s} <i<i_{s+1}} (e^-_{i} +e^+_{i-1}- e^0_{i} -e^0_{i-1}),
	 	$$
	 	with the convention that the $i_{r+1} =i_1$ and the interval $(i_r,i_{1})$ is understood in the cyclic ordering of the vertices $\{0,1,\ldots, (l-1)\}$.
	 	
	 	\vspace{0.2cm}
	 	3) If $l=2k$ is even, then the matching cycles for $M^{-,+} $ and $M^{+,-}$ are as follows
	 	$$
	 	e_{M^{-,+}} -e_{M^0} =\sum^k_{s=1}(e^-_{2s-1} +e^+_{2(s-1)} -e^0_{2s-1} -e^0_{2(s-1)}),
	 	$$
	 	$$
	 	e_{M^{+,-}} -e_{M^0} =\sum^k_{s=1}(e^+_{2s-1} +e^-_{2(s-1)} -e^0_{2s-1} -e^0_{2(s-1)}).
	 	$$
\end{pro}
\begin{pf} 
	Immediate from Proposition \ref{pro:matchings0}.
\end{pf}
	
\begin{rem}\label{rem:matchhomology}
	For $i$ in $\{0,1,\ldots,l-1\}$ the cycle 
	$$
	e^-_{i} +e^+_{i-1}- e^0_{i} -e^0_{i-1}
	$$
	corresponds to the path

	$$
		\begin{tikzpicture}
			[place/.style={circle,draw=black,thick},
			transition/.style={circle,draw=black,fill=black}]
			\node (whitei) at (210:1.5cm) [place] [label={above right:$i$}]{};
			\node (blacki) at (210:3cm) [transition] [label={below left:$i'$}]{};
			\node (whitej) at (270:1.5cm) [place] [label={above right:$(i-1)$}] {};
			\node (blackj) at (270:3cm) [transition] [label={below:$(i-1)'$}]{};
			\begin{scope}[gray,ultra thin]
				\draw (whitei) arc[start angle= 206, end angle=-90, radius=1.5cm];
				\draw (whitei) to [bend right=30](whitej);
			\draw (blacki) arc[start angle =210, end angle= -90, radius=3cm];
				\draw (blacki) to [bend right=30](blackj);
			\end{scope}
				\draw[red,ultra thick][<-,>={Stealth}] (whitej) to (blackj);
				\draw[red,ultra thick] [<-,>={Stealth}](whitei) to (blacki);
			\begin{scope}
				[decoration={markings,
					mark =at position 2cm with {\arrow[blue,line width=1mm]{stealth}}}]
				\draw[blue,thick][postaction={decorate}] (whitei) to [bend right=45](blackj);	
			\end{scope}
			\begin{scope}[decoration={markings,
					mark =at position 2cm with {\arrow[red,line width=1mm]{stealth}}}]
				\draw[red, dotted,thick][postaction={decorate}] (whitej) to [bend left=45](blacki);
			\end{scope}
		\end{tikzpicture}
	$$
of the graph $\widehat{PG}_{l}$. This is homologous to the equatorial circle in \eqref{torus} and hence has the homology class $(-1,0)$. From this it follows that under the epimorphism in the exact sequence
$$
\xymatrix{
	0\ar[r]&\ZZ\{B_i\}\ar[r]&H_1 (\widehat{PG}_{l},\ZZ) \ar[r]^{j}&H_1(\mathbb{T}, \ZZ)\ar[r]&0.
}
$$
the matching cycle $(e_M -e_{M^0})$, for all $M\neq M^{\pm}$, is mapped to a multiple of the homology class $-(1,0)$. More precisely, let $R(d)$ be the polygon formed by a proper diagonal $d=(i_s,i_{s+1}) \in \Delta_M$ and the vertices of $R_l$ between the end points of $d$ and denote by $v_d$ the number of vertices of $R(d)$, then 
$$
[e_M-e_{M^0}]=j(e_M -e_{M^0}) =-\frac{1}{2}\left(\sum_{d\in \Delta_M} (v_d -2) \right) (1,0)=-\frac{1}{2} (l-r)(1,0).
$$	

Furthermore, for any $r\equiv l\, mod(2)$ in the interval $[0,l]$ the class $-\frac{1}{2} (l-r)(1,0)$ is realized by a matching: take $M$ with a single proper diagonal, that is, 
$$
M(0)=\{e^0_i,e^0_{i-1}, \ldots, e^0_{i-r+1}\}
$$
is an uninterrupted sequence of $r\equiv l\, mod(2)$ consecutive $0$-colored edges of $\widehat{PG}_{l}$. 
\end{rem}
 
 We can now establish the image of the matching cycles under the the map
 \begin{equation}\label{j-map}
 \xymatrix{
 H_1 (\widehat{PG}_{l},\ZZ) \ar[r]^{j}&H_1(\mathbb{T}, \ZZ)\ar[r]&0.
 }
 \end{equation}
\begin{pro}\label{pro:matchingtriangle}
	Let 
	$$
	MC_{\widehat{PG}_{l}} =\{e_M-e_{M_0} | M \in {\cal M}_{\widehat{PG}_{l}}\}
	$$
	be the the set of the matching cycles.
	The epimorphism $j$ in \eqref{j-map} maps this set into the integral points of the triangular region $T_l$ in $\RR^2=\RR\otimes_{\ZZ} H_1(\mathbb{T}, \ZZ)$ formed by the points
	$$
(0,0),\,\,	[Z_{0,+}]=(-l,-1),\,\, [-Z_{0,-}]=(0,1),
	$$
	the vertices of $T_l$. These are the homology classes in $H_1(\mathbb{T}, \ZZ)$ corresponding respectively to the matchings $M^+$, $M^-$ and $M^0$.
	The region $T_l$, besides the vertices listed above, contains precisely $k=\left[\frac{l}{2}\right]$ other integral points
	$$
	\{(-m,0)| m=1,\ldots,k\}.
	$$
	Each point $(-m,0)$ is the homology class of all matchings having exactly $(l-2m)$ edges colored by $0$.
	
	If $l=2k+1$ is odd, then the three vertices of $T_l$ are the only integral points lying on the boundary of $T_l$. 
	
	If $l=2k$ is even, then the point $(-k,0)$ is the additional integral point on the boundary of $T_l$. 
	All the remaining points
	$$
	\{(-m,0)| m=1,\ldots,k-1\}
	$$
	are internal.
\end{pro}
\begin{pf}
	Follows from Proposition \ref{pro:matchcycles} and Remark \ref{rem:matchhomology}.
\end{pf}

Recall the toric surface $X(\Delta_l)$, see Proposition \ref{pro:toricsurf}, defined by the fan $\Delta_l$. The triangle $T_l$ is contained in the cone $\sigma_{\mp}$ as illustrated in the following drawing.
$$
\begin{tikzpicture}[>=stealth]
	\draw
	(-4,0)--(4,0);
	\draw
	(0,-2)--(0,2);
	\filldraw[gray]
	(0,0) circle [radius=3pt] node[right, below,color=red]{$(0,0)$}
	(0,1) circle [radius=3pt] node[right,color=red]{$[-Z_{0,-}]$}
	(-4,-1) circle [radius=3pt] node[left,color=red]{$[Z_{0,+}]$};
	\filldraw[black]
		(-0.5,0) circle [radius=2pt]
	(-1,0) circle [radius=2pt]
	(-1.5,0) circle [radius=2pt];
	\draw[black, very thick][->]
	(0,0)--(0,1);
	\draw[black, very thick][->]
	(0,0)--(-4,-1);
	\draw[black, very thick]
	(0,1)--(-4,-1);
	\begin{pgfonlayer}{background}
		\fill[green!20!white]
		(0,-2)--(0,0)--(-4,-1);
	\end{pgfonlayer}
	\draw
	(-1,-1.5) node{$\sigma_{+}$};
	\begin{pgfonlayer}{background}
		\fill[yellow!20!white]
		(0,-2)--(0,0)--(4,2)--(4,-2);
	\end{pgfonlayer}
	\draw
	(3,-1.5) node{$\sigma_{-}$};
	\begin{pgfonlayer}{background}
		\fill[blue!20!white]
		(-4,-1)--(0,0)--(4,2)--(-4,2);
	\end{pgfonlayer}
	\draw
	(-1,1) node{$\sigma_{\mp}$};
\end{tikzpicture}
$$
 The side of the triangle joining $[Z_{0,+}]$ with $[-Z_{0,-}]$ intersects the horizontal axis at the point $(-\frac{l}{2},0)$; the integral points of that axis contained in the interior of $T_l$ are precisely the integral interior points of $T_l$ (the black dots on the picture); 
 if $l$ is even the point $(-\frac{l}{2},0)$ gives an additional integral point on the boundary of $T_l$, otherwise the vertices (gray dots) of the triangle are the only integral points on the boundary of $T_l$.
 
 The integral points in the diagram represent only `homological' degree of perfect matchings. That degree tells us only the number of $0$-colored edges of the matchings: there might be several matchings with the same degree. The general method for counting the number of perfect matchings is the Kasteleyn determinant associated to the bipartite graph, see \cite{Go-K}, \cite{Ken}. For convenience of the reader we recall the basic construction of the Kasteleyn matrix of the bipartite graph $G$ on the topological torus ${\mathbb{T}}$. 
 
 The matrix one constructs has rows indexed by the ordered set of black vertices $V^b_G$ of $G$ and the columns by the white ones $V^w_G$. Let $i_b$ and $j_w$ be respectively the row and column index; if the corresponding vertices of $G$ are related by $k$ edges, numbered in some order, one assigns the weight indeterminate $a^r_{i_b j_w}$ to the $r$-th edge; there is also a sign involved which depends on the number of edges of the boundary cycle of $G$ containing the edge, but we will ignore that since in our case all boundary cycles have the same length equals $6$; in the first approximation the $(i_b,j_w)$-entry of the Kasteleyn matrix is the sum of weights
 $$
 \sum^k_{r=1} a^r_{i_b j_w};
 $$  
 this is modified according to the intersection of each edge $e_r$ of the sum with generators of the homology group $H_1({\mathbb{T},\ZZ})$ of the torus: let $\gamma_1$ and $\gamma_2$ be a basis of $H_1({\mathbb{T},\ZZ})$
 and assign the variable $z$ (resp.,$w$) to $\gamma_1$ (resp. $\gamma_2$); the intersection weight of $e_r$ is the monomial
 $$
 z^{(\gamma_1,e_r)} w^{(\gamma_2,e_r)}
 $$
where $(\gamma_i,e_r)$ the signed intersection number of $\gamma_i$ with $e_r$: a point of intersection contributes $1$ if $\gamma_i$ traverses $e_r$ with the black vertex on the left and white on the right; if the order is reversed, the contribution is $(-1)$; if $e_r$ does not intersect $\gamma_i$ the contribution is $0$. The Kasteleyn matrix $M_G (z,w)$ is defined as follows
$$
M_G (z,w)=(M^{i_b,j_w}_G(z,w))_{i_b \in V^b_G,j_w\in V^w_G},
$$
where the entries $M^{i_b,j_w}_G(z,w)$ are determined by the formula
$$
M^{i_b,j_w}_G(z,w)=\begin{cases}\sum^k_{r=1} a^r_{i_b j_w}z^{(\gamma_1,e_r)} w^{(\gamma_2,e_r)},\\
	0,\,\, \text{if no edges between $i_b$ and $j_w$.}
	\end{cases}  
$$

\vspace{0.2cm}
We will determine the Kasteleyn matrix for the graphs $\widehat{PG}_l$. To simplify the notation we write $M_l(z,w)$ instead of $M_{\widehat{PG}_l} (z,w)$. This is $l\times l$-matrix; we number the rows and columns by the indices $i \in [0,l-1]$ of the black/white vertices of $\widehat{PG}_l$ in the decreasing order as they appear in the graph, that is, the rows (resp. columns) vary from $(l-1)$ at the top (resp., left) to $0$ at the bottom (resp. right). Recall, that the edges of $\widehat{PG}_l$ at the white vertex $i$ are labeled 
$$
e^0_i =(i)\rightarrow (i)',\,\, e^-_i=(i)\rightarrow (i-1)',\,\, e^+_i =(i)\rightarrow (i+1)',
$$
with the cyclic convention that $e^-_0=(0)\rightarrow (l-1)'$ and $e^+_{l-1}=(l-1)\rightarrow (0)'$. Our weight assignment is as follows
\begin{equation}\label{wlabels}
	e^0_i \mapsto x_i, \,\, e^-_i \mapsto x^-_i,\,\, e^+_i\mapsto x^+_i.
\end{equation}
To calculate the intersection weights we represent the topological torus ${\mathbb{T}}$ by a `fundamental' parallelogram

$$
\begin{tikzpicture}[>=Stealth]
	\draw[->][black, thick]
	(0,0)--(4,0) node[below right]{$\gamma_1$};
	\draw[black, densely dotted]
	(4,0)--(5,2)--(1,2);
	\draw[->][black, thick]
	(0,0)--(1,2) node[left]{$\gamma_2$};
\end{tikzpicture}
$$
with $\gamma_i$, for $i=1,2$ representing a basis for the homology group $H_1({\mathbb{T}, \ZZ})$ of the torus. We now draw the graph $\widehat{PG}_l$ in this fundamental domain: all vertices of the graph as well as the $0$-colored edges are in the interior of the domain; the edges colored `$-$', with the exception of $e^-_0$, are also in the interior, the following picture illustrates the case $l=4$:
$$
\begin{tikzpicture}[>=Stealth]
	\draw[->][black, thick]
	(0,0)--(6,0) node[below right]{$\gamma_1$};
	\draw[black, densely dotted]
	(6,0)--(7,3)--(1,3);
	\draw[->][black, thick]
	(0,0)--(1,3) node[left]{$\gamma_2$};
	\begin{scope}
			[place/.style={circle,draw=black,thick, inner sep=0pt, minimum size=2mm},
		transition/.style={circle,draw=black,fill=black, inner sep=0pt, minimum size=2mm}]
	\node (white3) at (2,1) [place] [label={below:$3$}] {};
	\node (black3) at (2,2.5) [transition] [label={}]{};
	\node (white2) at (3,1) [place] [label={below:$2$}]{};
	\node (black2) at (3,2.5) [transition] [label={}]{};
	\node (white1) at (4,1) [place] [label={ below:$1$}]{};
	\node (black1) at (4,2.5) [transition] [label={}]{};
		\node (white0) at (5,1) [place] [label={ below:$0$}]{};
	\node (black0) at (5,2.5) [transition] [label={}]{};
	\draw[black,ultra thick][-] (white3) to (black3);
	\draw[black,ultra thick][-] (white2) to (black2);
	\draw[black,ultra thick][-] (white1) to (black1);
	\draw[black,ultra thick][-] (white0) to (black0);
	\draw[blue,ultra thick][-] (white3) to (black2);
	\draw[blue,ultra thick][-] (white2) to (black1);
	\draw[blue,ultra thick][-] (white1) to (black0);
	\end{scope}
\end{tikzpicture}
$$
The remaining edge $e^-_0$ edge starts at the white vertex $0$ and joins the
black vertex $(l-1)'$ by crossing  the boundary of the parallelogram as illustrated below 
$$
\begin{tikzpicture}[>=Stealth]
	\draw[->][black, thick]
	(0,0)--(6,0) node[below right]{$\gamma_1$};
	\draw[black, densely dotted]
	(6,0)--(7,3)--(1,3);
	\draw[->][black, thick]
	(0,0)--(1,3) node[left]{$\gamma_2$};
	\begin{scope}
		[place/.style={circle,draw=black,thick, inner sep=0pt, minimum size=2mm},
		transition/.style={circle,draw=black,fill=black, inner sep=0pt, minimum size=2mm}]
		\node (white3) at (2,1) [place] [label={below:$3$}] {};
		\node (black3) at (2,2.5) [transition] [label={}]{};
		\node (white2) at (3,1) [place] [label={below:$2$}]{};
		\node (black2) at (3,2.5) [transition] [label={}]{};
		\node (white1) at (4,1) [place] [label={ below:$1$}]{};
		\node (black1) at (4,2.5) [transition] [label={}]{};
		\node (white0) at (5,1) [place] [label={ below:$0$}]{};
		\node (black0) at (5,2.5) [transition] [label={}]{};
		\draw[black,ultra thick][-] (white3) to (black3);
		\draw[black,ultra thick][-] (white2) to (black2);
		\draw[black,ultra thick][-] (white1) to (black1);
		\draw[black,ultra thick][-] (white0) to (black0);
		\draw[blue,ultra thick][-] (white3) to (black2);
		\draw[blue,ultra thick][-] (white2) to (black1);
		\draw[blue,ultra thick][-] (white1) to (black0);
		\draw[blue,ultra thick]
		(white0)--(6.5,1.5);
		\draw[blue,ultra thick]
		(black3)--(0.5,1.5);
	\end{scope}
\end{tikzpicture}
$$
All edges colored `$+$' must cross the boundary: to draw the edge $e^+_i$ we start at the white vertex $i$ and join it to the black vertex $(i+1)'$ by first, descending to the boundary of the parallelogram and then reappearing on the opposite side; this is illustrated below, where the `$+$'-colored edges are in red:
$$
\begin{tikzpicture}[>=Stealth]
	\draw[->][black, thick]
	(0,0)--(6,0) node[below right]{$\gamma_1$};
	\draw[black, densely dotted]
	(6,0)--(7,3)--(1,3);
	\draw[->][black, thick]
	(0,0)--(1,3) node[left]{$\gamma_2$};
	\begin{scope}
		[place/.style={circle,draw=black,thick, inner sep=0pt, minimum size=2mm},
		transition/.style={circle,draw=black,fill=black, inner sep=0pt, minimum size=2mm}]
		\node (white3) at (2,1) [place] [label={below:$3$}] {};
		\node (black3) at (2,2.5) [transition] [label={}]{};
		\node (white2) at (3,1) [place] [label={below:$2$}]{};
		\node (black2) at (3,2.5) [transition] [label={}]{};
		\node (white1) at (4,1) [place] [label={ below:$1$}]{};
		\node (black1) at (4,2.5) [transition] [label={}]{};
		\node (white0) at (5,1) [place] [label={ below:$0$}]{};
		\node (black0) at (5,2.5) [transition] [label={}]{};
		\draw[black,ultra thick][-] (white3) to (black3);
		\draw[black,ultra thick][-] (white2) to (black2);
		\draw[black,ultra thick][-] (white1) to (black1);
		\draw[black,ultra thick][-] (white0) to (black0);
		\draw[blue,ultra thick][-] (white3) to (black2);
		\draw[blue,ultra thick][-] (white2) to (black1);
		\draw[blue,ultra thick][-] (white1) to (black0);
		\draw[blue,ultra thick]
		(white0)--(6.5,1.5);
		\draw[blue,ultra thick]
		(black3)--(0.5,1.5);
	\end{scope}
\begin{scope}[red, very thick]
	\draw (white2)--(2.5,0);
	\draw (white1)--(3.5,0);
	\draw (white0)--(4.5,0);
	\draw (black3)--(2.5,3);
	\draw (black2)--(3.5,3);
	\draw (black1)--(4.5,3);
	\draw (white3)--(0,0);
	\draw (black0)--(7,3);
		\end{scope}
\end{tikzpicture}
$$
We can now write down the entries of the Kasteleyn matrix $M_l(z,w)$: in $i$-th row we have precisely three nonzero entries in the columns $(i+1)$, $i$ and $(i-1)$:
$$
\begin{gathered}
M^{i,i}_l(z,w)=x_i,
\\
 M^{i,i+1}_l (z,w) =\begin{cases}
	x^-_{i+1}, & \text{ if $i \neq l-1$},
		\\
		x^-_0 w^{-1}, &\text{if $i=l-1$}, 
\end{cases}
\\
M^{i,i-1}_l (z,w) =\begin{cases}
	x^+_{i-1}z^{-1}, & \text{ if $i \neq 0$},
	\\
	x^+_{l-1} z^{-1}w, &\text{if $i=0$}, 
\end{cases}
\end{gathered}
$$
Hence the Kasteleyn matrix has the form
\begin{equation}\label{Kastelyanmatrix}
	M_l(z,w)=\begin{pmatrix}
		x_{l-1} &x^+_{l-2}z^{-1}&0&0&\cdots&0&0&0&x^-_0 w^{-1}\\
		x^-_{l-1}&x_{l-2}&x^+_{l-3}z^{-1}&0&\cdots&0&0&0&0\\
		0&x^-_{l-2}&x_{l-3}&x^+_{l-4}z^{-1}&\cdots&0&0&0&0\\ 
		\cdots&\cdots&\cdots&\cdots&\cdots&\cdots&\cdots&\cdots&\cdots\\
		0&0&0&0&\cdots&x^-_3&x_2&x^+_1z^{-1}&0\\
		0&0&0&0&\cdots&0&x^-_2&x_1&x^+_0z^{-1}\\
		x^+_{l-1}z^{-1}w&0&0&0&\cdots&0&0&x^-_1&x_0
	\end{pmatrix}
\end{equation}
that is, on the principal diagonal run the weights $x_i$ attached to $0$-colored edges, on the next to it, lower diagonal, run the weights $x^-_i$, from $l-1$ on the top, to $1$ on the bottom, on the diagonal 
up from the principal one run $x^+_i z^{-1}$, for $i$ from $(l-2)$ to $0$,
and in the corners we have $x^-_0 w^{-1}$ and  $x^+_{l-1}z^{-1}w$, in the upper right and lower left respectively.

One of the main properties of the Kasteleyn matrix is that its determinant
counts the number of perfect matchings of a bipartite graph, see \cite{Go-K}, \cite{Ken}, for details and references on the subject. The following statement spells out the situation for the graph $\widehat{PG}_l$. 
\begin{pro}\label{pro:detKastelyan}
	Let $K_l(z,w)$ be the determinant of $M_l(z,w)$
	$$
	K_l(z,w)=det(M_l (z,w))=\begin{vmatrix}
			x_{l-1} &x^+_{l-2}z^{-1}&0&\cdots&0&0&0&x^-_0 w^{-1}\\
		x^-_{l-1}&x_{l-2}&x^+_{l-3}z^{-1}&\cdots&0&0&0&0\\
		0&x^-_{l-2}&x_{l-3}&\cdots&0&0&0&0\\ 
		\cdots&\cdots&\cdots&\cdots&\cdots&\cdots&\cdots&\cdots\\
		0&0&0&\cdots&x^-_3&x_2&x^+_1z^{-1}&0\\
		0&0&0&\cdots&0&x^-_2&x_1&x^+_0z^{-1}\\
		x^+_{l-1}z^{-1}w&0&0&\cdots&0&0&x^-_1&x_0
	\end{vmatrix}
$$
Its expansion is a polynomial in $z^{-1},w,w^{-1}$ with coefficients
in the algebra of polynomials $\CC[x_i,x^+_i,x^-_i \Big| i\in[0,l-1]]$ in edge weights. More precisely, we have
\begin{equation}\label{Kl-formula}
	K_l(z,w)=(-1)^{l-1} A^-_l w^{-1} +(-1)^{l-1} A^+_l wz^{-l} +P_l (z^{-1}),
	\end{equation}
	where the notation is as follows:
	
	$\bullet$
$\displaystyle{	A^{\pm}_l=\prod^{l-1}_{i=0} x^{\pm}_i}$,

$\bullet$ $\displaystyle{P_l (z^{-1})=\sum^{[\frac{l}{2}]}_{d=0} {\mathfrak m}_d z^{-d}}$ is a polynomial in $z^{-1}$ with coefficients ${\mathfrak m}_d$ homogeneous polynomials of degree $l$ in $\CC[x_i,x^+_i,x^-_i \Big| i\in[0,l-1]]$.

 Each monomial occurring in ${\mathfrak m}_d$ is a product
of precisely 

- $d$ indeterminates from the set $\{x^-_i\}$, 

- $d$ indeterminates from the set $\{x^+_i\}$, furthermore, $x^+_j$ occurs
in the set if and only if $x^-_{j+1}$ appears in the previous item;

 - $(l-2d)$ indeterminates from the set $\{x_i\}$, with the indices $i$  not occurring in the two previous items. 

The {\rm number} of monomials appearing in ${\mathfrak m}_d$ is the number of perfect matchings of $\widehat{PG}_l$ having precisely $(l-2d)$ edges colored $0$. In particular, the number of perfect matchings of $\widehat{PG}_l$ is the number of terms on the right hand side of \eqref{Kl-formula}:
$$
|{\cal M}_{\widehat{PG}_l}|=2 +\sum^{[\frac{l}{2}]}_{d=0} NM({\mathfrak m}_d),
$$
where  $NM({\mathfrak m}_d)$ stands for the number of monomials in ${\mathfrak m}_d$. Furthermore, for $l\geq 4$, the following formula holds
$$
\sum^{[\frac{l}{2}]}_{d=0} NM({\mathfrak m}_d)= f_{l-1}+2f_{l-2},
$$
where ${f_2,f_3,\ldots,f_{l-1}}$ is a part of the Fibonacci sequence with  $f_2=2$, $f_3=3$. 
\end{pro}
\begin{pf}
	Expand the determinant  of the matrix $M_l(z,w)$ with respect to the first column to obtain the expression
	$$
	K_l(z,w)=(-1)^{l-1}A^-_l w^{-1} +(-1)^{l-1}A^+_l wz^{-l} +P_l (z^{-1}).
	$$
	where the polynomial $P_l(z^{-1})$ is determined by the formula
	\begin{equation}\label{Pl-Kastelyan}
	P_l(z^{-1})=x_{l-1} X_{l-1}(z^{-1}) -x^-_{l-1}x^+_{l-2}z^{-1}X^{\ulcorner}_{l-1}(z^{-1}) -x^-_0x^+_{l-1}z^{-1}X^{{}^{\lrcorner}}_{l-1}(z^{-1});
	\end{equation}
	on the right hand side of the formula 
	
	- $X_{l-1}(z^{-1})$ is the determinant of the matrix
	$$
	M^{\ulcorner}_{l}=\begin{pmatrix}
		x_{l-2}&x^+_{l-3}z^{-1}&0&0&\cdots&0&0&0&0\\
		x^-_{l-2}&x_{l-3}&x^+_{l-4}z^{-1}&0&\cdots&0&0&0&0\\ 
		0&x^-_{l-3}&x_{l-4}&x^+_{l-5}z^{-1}&\cdots&0&0&0&0\\
		\cdots&\cdots&\cdots&\cdots&\cdots&\cdots&\cdots&\cdots&\cdots\\
		0&0&0&0&\cdots&x^-_3 &x_2&x^+_1 z^{-1}&0\\
			0&0&0&0&\cdots&0&x^-_2 &x_1&x^+_0 z^{-1}\\
				0&0&0&0&\cdots&0&0&x^-_1 &x_0
	\end{pmatrix}
$$
obtained from the Kasteleyn matrix $M_l(z,w)$ by deleting the first row and column; 

- $X^{\ulcorner}_{l-1}(z^{-1})$ and $X^{{}^{\lrcorner}}_{l-1}(z^{-1})$ are respectively the $(1,1)$ and $(l-1,l-1)$ minors of the matrix $	M^{\ulcorner}_{l}$.

From the formula \eqref{Pl-Kastelyan} and the definition of $X_{l-1}$ and $X^{\ulcorner}_{l-1}(z^{-1})$, $X^{{}^{\lrcorner}}_{l-1}(z^{-1})$ follows 
the assertion
$$
 P_l (z^{-1})=\sum^{[\frac{l}{2}]}_{d=0} {\mathfrak m}_d z^{-d}
 $$
 with the coefficients ${\mathfrak m}_d$ as described in the statement of the proposition. From \eqref{Pl-Kastelyan} also follows that the number of terms in $ P_l (z^{-1})$ is given by the formula
 $$
 \text{number of terms in $P_l (z^{-1})$}=f_{l-1} +2f_{l-2},
 $$
 where $f_k$ stands for the number of terms in the expansion of the determinant
 $$
X_k = \begin{vmatrix}
 	x_{k-1}&x^+_{k-2}z^{-1}&0&0&\cdots&0&0&0&0\\
 	x^-_{k-1}&x_{k-2}&x^+_{k-3}z^{-1}&0&\cdots&0&0&0&0\\ 
 	0&x^-_{k-2}&x_{k-3}&x^+_{k-4}z^{-1}&\cdots&0&0&0&0\\
 	\cdots&\cdots&\cdots&\cdots&\cdots&\cdots&\cdots&\cdots&\cdots\\
 	0&0&0&0&\cdots&x^-_3 &x_2&x^+_1 z^{-1}&0\\
 	0&0&0&0&\cdots&0&x^-_2 &x_1&x^+_0 z^{-1}\\
 	0&0&0&0&\cdots&0&0&x^-_1 &x_0
 \end{vmatrix}
 	$$
 	for $k\geq 2$. Expanding with respect to the first column one finds the recursion formula
 	$$
 	X_k =x_{k-1}X_{k-1} -x^-_{k-1}x^{+}_{k-2}z^{-1} X_{k-2}.
 	$$
 	Hence the recursion relation for the number of terms
 	$$
 	f_k =f_{k-1}+f_{k-2},
 	$$
 	for $k\geq 4$; this means that $\{f_k\}$ is the sequence of Fibonacci; for $k=2$ and $k=3$ we have
 	$$
 	X_2=\begin{vmatrix}
 		x_1 &x^+_0 z^{-1}\\
 		x^-_1 &x_0 
 	\end{vmatrix},
 \hspace{0.2cm}
 X_3=\begin{vmatrix}
 	 x_2&x^+_1 z^{-1}&0\\
 	x^-_2 &x_1&x^+_0 z^{-1}\\
 	0&x^-_1 &x_0
 \end{vmatrix}
$$
Hence the number of terms $f_2=2$ and $f_3=3$. 
\end{pf} 
\begin{rem}\label{rem:Kastelyanl}
	1) The Newton polygon of the Kasteleyn determinant $K_l(z,w)$ is easily seen to be the triangle with vertices
	$$
	(0,0), (0,-1), (-l,1).
	$$
	Reflecting with respect to the horizontal line, we find the triangle $T_l$ of matching cycles described in Proposition \ref{pro:matchingtriangle}; the reflection corresponds to renaming $w$ by $w^{-1}$.
	
	2) The equations 
	\begin{equation}\label{Kast=0}
	K_l(x,z,w)=0,
	\end{equation}
	where $x=(x_i,x^{\pm}_i)_{i\in[0,l-1]}$ is the weight variable, could be viewed as families of curves in the algebraic torus $(\CC^{\times})^2$. The latter is viewed as the complement of the toric divisor in the surface $X(\Delta_l)$. In \cite{KeOkSh}, the authors studied curves determined by the vanishing of Kasteleyn determinant for a general bipartite graph on a torus. Viewing those as plane curves they showed many remarkable properties of those. In particular, they are Harnack curves, whenever the values of weights assigned to edges are positive real numbers. We will refer to the curves in \eqref{Kast=0} as
	{\rm Kasteleyn} curves
	$$
	{\mathfrak C}_l (x):=\overline{(K_l(x,z,w)=0)},
	$$
	the closure of $(K_l(x,z,w)=0)$ in $X(\Delta_l)$.
	
	\vspace{0.2cm}
	3) One can also view the equations \eqref{Kast=0} as the equations in {\rm edge variables} with $(z,w) \in  (\CC^{\times})^2$ as parameters. From this point of view \eqref{Kast=0} determines a family of hypersurfaces of degree $l$ in $\PP(\CC E_{\widehat{PG}_l})=\PP^{3l-1}$, the projectivization of the vector space $\CC E_{\widehat{PG}_l}$ spanned by the edges $E_{\widehat{PG}_l}$ of the graph $\widehat{PG}_l$. 
	This is another way to connect to Calabi-Yau varieties: a complete intersection of three hypersurfaces of degree $l$ in $\PP^{3l-1}$. This will be explored elsewhere.     
\end{rem}  
  
  As was mentioned earlier, the Newton polygon of the $K_l(z,w)$ captures only the homology
  classes on the topological torus $\mathbb{T}$ of the matching cycles. The `multiplicity' of each homology cycle is read off the corresponding coefficient $K_l(z,w)$, the homogeneous polynomial of degree $l$ in the edge-weights
  $\{x_i,x^{\pm}_i \big| i\in[0,l-1]\}$: the monomials occurring in those homogeneous polynomials label the matchings having the same homology class of the torus $\mathbb{T}$; Proposition \ref{pro:detKastelyan} also tells us that the configuration of edges of a matching can be read off the monomials as well.  
  
All of the above indicates that the structure of perfect matchings is much more involved. Indeed the subject has been explored extensively in mathematics and physics literature, see \cite{Go-K}, \cite{Ken} and references therein. We will see later on that the perfect matchings of $\widehat{PG}_l$ have structure of categories. As an indication of this rich structure, we show that the matchings of $\widehat{PG}_l$ can be assembled into a homology complex. 
  
  We denote by $\ZZ[{\cal M}_{\widehat{PG}_l}] $ the $\ZZ$-module freely generated by the set of perfect matchings of $\widehat{PG}_l$. For a perfect matching $M$ define $deg(M)$, the degree of $M$, as the number of $0$-colored edges in $M$, that is, the number of elements $|M(0)|$ of the set $M(0)$:
  $$
  deg(M)=|M(0)|.
  $$
  This defines a grading on $\ZZ[{\cal M}_{\widehat{PG}_l}] $:
  $$
  \ZZ[{\cal M}_{\widehat{PG}_l}]^d:=\text{the submodule of  $\ZZ[{\cal M}_{\widehat{PG}_l}]$ generated by the matchings of degree $d$.}
  $$
  This turns $\ZZ[{\cal M}_{\widehat{PG}_l}]$ into a graded $\ZZ$-module
  $$
  \ZZ[{\cal M}_{\widehat{PG}_l}] =\bigoplus_{d} \ZZ[{\cal M}_{\widehat{PG}_l}]^d
  $$
  The values of the degree $d$ vary in $[0,l]$ and we know
  
  - $\ZZ[{\cal M}_{\widehat{PG}_l}]^0$ is spanned by two or four matchings according to the parity of $l$:
  
  if $l$ is odd, then the submodule is spanned by the matchings $M^-$ and $M^+$;
  
   if $l$ is even, the submodule has two more generators, $M^{-,+}$ and $M^{+,-}$;
   
   - $\ZZ[{\cal M}_{\widehat{PG}_l}]^{l}=\ZZ{M^0}$.

   In addition, all positive values of $d$ must be of the same parity as $l$. We wish to investigate the positive part
   $$
   \ZZ[{\cal M}_{\widehat{PG}_l}]^{+}:=\bigoplus_{d>0} \ZZ[{\cal M}_{\widehat{PG}_l}]^d.
   $$
   From Proposition \ref{pro:detKastelyan} we know the following.
   \begin{pro}
   For every $d$ the matchings forming a basis of $\ZZ[{\cal M}_{\widehat{PG}_l}]^d$, are in bijection with the monomials of the homogeneous polynomial
   	$\displaystyle{{\mathfrak m}_{\HA (l-d)}}$, the coefficient of
   	$z^{-\HA (l-d)}$ in the polynomial $P_l(z^{-1})$, see Proposition \ref{pro:detKastelyan}. In particular, the rank of $\ZZ[{\cal M}_{\widehat{PG}_l}]^d$ equals the number of monomials in $\displaystyle{{\mathfrak m}_{\HA (l-d)}}$: 
   	$$
   	rk (\ZZ[{\cal M}_{\widehat{PG}_l}]^d) = NM({\mathfrak m}_{\HA (l-d)}).
   	$$ 
   \end{pro}
   We will turn $\ZZ[{\cal M}_{\widehat{PG}_l}]^{+}$ into a differential graded module by defining the differentials
   $$
   \delta^d: \ZZ[{\cal M}_{\widehat{PG}_l}]^d \longrightarrow \ZZ[{\cal M}_{\widehat{PG}_l}]^{d+2},
   $$
   for every positive value of $d$. This is reminiscent of the (co)boundary operator for simplicies. First recall that for a perfect matching $M$ of degree $d>0$ the complement of $M\setminus M(0)$, if nonempty, consists of disjoint union of the strings of consecutive pairs of edges of color $\pm$:
   \begin{equation}\label{+-string}
   \begin{tikzpicture}
   	[place/.style={circle,draw=black,thick},
   	transition/.style={circle,draw=black,fill=black}]
   	\begin{scope}
   	\node (white3) at (0,2) [place] [label={above:$i$}] {};
   	\node (black3) at (0,0) [transition] [label={below:$i'$}]{};
   	\node (white2) at (2,2) [place] [label={above:$i-1$}]{};
   	\node (black2) at (2,0) [transition] [label={below:$(i-1)'$}]{};
   		\draw[blue][-] (white3) to (black2);	
   	\draw[red, dotted,thick] (white2) to (black3);
   	\end{scope}
   \begin{scope}[xshift=4cm]
   	\node (white3) at (0,2) [place] [label={above:$i-2$}] {};
   	\node (black3) at (0,0) [transition] [label={below:$(i-2)'$}]{};
   	\node (white2) at (2,2) [place] [label={above:$i-1$}]{};
   	\node (black2) at (2,0) [transition] [label={below:$(i-1)'$}]{};
   	\draw[blue][-] (white3) to (black2);	
   	\draw[red, dotted,thick] (white2) to (black3);
   \end{scope}
\begin{scope}[xshift=9cm]
	\node (white3) at (0,2) [place] [label={above:$i-2k$}] {};
	\node (black3) at (0,0) [transition] [label={below:$(i-2k)'$}]{};
	\node (white2) at (2,2) [place] [label={above:$i-2k-1$}]{};
	\node (black2) at (2,0) [transition] [label={below:$(i-2k-1)'$}]{};
	\draw[blue][-] (white3) to (black2);	
	\draw[red, dotted,thick] (white2) to (black3);
\end{scope}
\draw[black, loosely dotted, very thick]
(7,1)--(8,1);
   \end{tikzpicture}
\end{equation} 
above is depicted a string of $(k+1)$ consecutive pairs of $\pm$ edges. We think of a pair of $\pm$ edges
$$
 \begin{tikzpicture}
 	[place/.style={circle,draw=black,thick},
 	transition/.style={circle,draw=black,fill=black}]
 		\node (white3) at (0,2) [place] [label={above:$s$}] {};
 		\node (black3) at (0,0) [transition] [label={below:$s'$}]{};
 		\node (white2) at (2,2) [place] [label={above:$s-1$}]{};
 		\node (black2) at (2,0) [transition] [label={below:$(s-1)'$}]{};
 		\draw[blue][-] (white3) to (black2);	
 		\draw[red, dotted,thick] (white2) to (black3); 
  \end{tikzpicture}
$$
as a sort of braiding and define the {\it unbraiding} operation by replacing this pair by two consecutive $0$-colored edges of $\widehat{PG}_l$:
$$
\begin{tikzpicture}
	[place/.style={circle,draw=black,thick},
	transition/.style={circle,draw=black,fill=black}]
	\node (white3) at (0,2) [place] [label={above:$s$}] {};
	\node (black3) at (0,0) [transition] [label={below:$s'$}]{};
	\node (white2) at (2,2) [place] [label={above:$s-1$}]{};
	\node (black2) at (2,0) [transition] [label={below:$(s-1)'$}]{};
	\draw[blue][-] (white3) to (black2);	
	\draw[red, dotted,thick] (white2) to (black3);
	\begin{scope}[xshift=6cm]
		[place/.style={circle,draw=black,thick},
		transition/.style={circle,draw=black,fill=black}]
		\node (white3) at (0,2) [place] [label={above:$s$}] {};
		\node (black3) at (0,0) [transition] [label={below:$s'$}]{};
		\node (white2) at (2,2) [place] [label={above:$s-1$}]{};
		\node (black2) at (2,0) [transition] [label={below:$(s-1)'$}]{};
		\draw[red,very thick][-] (white3) to (black3);	
		\draw[red, very thick] (white2) to (black2); 
	\end{scope}
	\draw[->,decorate,decoration={snake, amplitude=.4mm, segment length=2mm, post length=1mm}, >={Stealth}]
(3,1)--(5,1)
node[above,align=center,midway]
{\textcolor{red} {Unbraiding}}; 
\end{tikzpicture}
$$
 Such an unbraiding operation performed with the pair $(e^-_s, e^+_{s-1})$ will be denoted $\delta_s$
 $$
 \delta_s(e^-_s, e^+_{s-1}) =(e^0_s,e^0_{s-1}).
 $$
 We will now define the map
 \begin{equation}
 \delta^d: \ZZ[{\cal M}_{\widehat{PG}_l}]^d \longrightarrow \ZZ[{\cal M}_{\widehat{PG}_l}]^{d+2},
\end{equation}
  as a combination of the unbraiding operations $\delta_s$ for 
  the values of $s$ occurring among the braided pairs of a given matching.
  Explicitly, let $M$ be a perfect matching of degree $d$ and let $M(0)$ be the $0$-colored edges in $M$:
  $$
  M(0)=\{e^0_{i_1}, e^0_{i_2}, \ldots, e^0_{i_d}\},
  $$
  where the labels $i_k$ are always listed in counterclockwise order starting from $i_1$ and ending at $i_d$. Consider the complement
  \begin{equation}\label{M(0)-comp}
  M\setminus M(0),
\end{equation}
 the collection of the braided pairs of $M$. If this is empty, then $M=M^0$, that is we are in $\ZZ[{\cal M}_{\widehat{PG}_l}]^l$ and we define
 $\delta^l$ to be zero:
 $$
 \delta^l: \ZZ[{\cal M}_{\widehat{PG}_l}]^l \longrightarrow 0.
 $$
 We now assume that $M\setminus M(0)$ is nonempty. Then it falls into disjoint union of strings of consecutive braided pairs. Let
 $$
 M_i=[(e^-_i,e^+_{i-1}),(e^-_{i-2},e^+_{i-3}),\ldots, (e^-_{i-2k},e^+_{i-2k-1})]
 $$
 be such a string, see \eqref{+-string}. Apply the unbraiding operators
 $\delta_s$ to every $s$ occurring in the string, weigh the result with the sign $\pm1$ and replace the string $M_i$ by the weighted sum
 {\small
 $$
 \begin{gathered}
 \delta(M_i):=\sum_{s}(-1)^{\frac{i-s}{2}} [(e^-_i,e^+_{i-1}),\ldots,(e^-_{s+2},e^+_{s+1}),\delta_s(e^-_{s},e^+_{s-1}), (e^-_{i-2},e^+_{s-3})\ldots, (e^-_{i-2k},e^+_{i-2k-1})]
 \\
 =\sum_{s}(-1)^{\frac{i-s}{2}} [(e^-_i,e^+_{i-1}),\ldots,(e^-_{s+2},e^+_{s+1}),(e^0_{s},e^0_{s-1}), (e^-_{s-2},e^+_{s-3})\ldots, (e^-_{i-2k},e^+_{i-2k-1})],
 \end{gathered}
 $$
}

\noindent
 where the sum is taken over $s=i-2j$, for $j\in [0,k]$. Define $\delta^d(M)_i$ as a signed sum of matchings by replacing in $M$ the string
 $M_i$ by $\delta(M_i)$; finally set $\delta^d(M)$ to be the sum of $(-1)^i\delta^d(M)_i$ as $M_i$ ranges through the strings occurring in the complement \eqref{M(0)-comp}:
 \begin{equation}\label{diffM-formula}
 \delta^d(M):=\sum_{M_i}(-1)^i\delta^d(M)_i,
 \end{equation}  
 again the numbering of the strings $M_i$ is taken in a counterclockwise order as we walk from $i_1$ around the circle.
 \begin{example}\label{ex:diff-matchings}
 	We take $l=10$; it is even and this means that $\ZZ[{\cal M}_{\widehat{PG}_{10}}]^{+}$ is evenly graded
 	$$
 	\ZZ[{\cal M}_{\widehat{PG}_{10}}]^{+}=\bigoplus^5_{k=1} \ZZ[{\cal M}_{\widehat{PG}_{10}}]^{2k}.
 	$$
 	Take a matching $M$ in $\ZZ[{\cal M}_{\widehat{PG}_{10}}]^{2}$ defined 
 	by the following configuration of edges
 	$$
 	M=[e^0_9, (e^-_8,e^+_7), e^0_6, (e^-_5,e^+_4),(e^-_3,e^+_2),(e^-_1,e^+_0)];
 	$$
 	diagrammatically this can be pictured as follows:
 	$$
 	\xymatrix@C=8pt@R=20pt{
 	9 \ar@{-}[d] &8\ar[dr]&7\ar[dl]&6\ar@{-}[d]&5\ar[dr]&4\ar[dl]&3\ar[dr]&2\ar[dl]&1\ar[dr]&0\ar[dl]\\
 	9' &8'&7'&6'&5'&4'&3'&2'&1'&0'
 }
$$
 Thus the set 
 $$
 M(0)=\{e^0_9,e^0_6\}
 $$ 
 and the complement $M\setminus M(0)$ falls into  two disjoint strings 
 $$
 M_1=[(e^-_8,e^+_7)],\,\,\, M_2=[(e^-_5,e^+_4),(e^-_3,e^+_2),(e^-_1,e^+_0)].
 $$
 Thus we have
 \begin{equation}\label{delta2-formula}
 \delta^2(M)=-\delta^2(M)_1 +\delta^2(M)_2.
  \end{equation}
 The matching $\delta^2(M)_1$:
 $$
 \xymatrix@C=8pt@R=20pt{
 	9 \ar@{-}[d] &8\ar@{-}[d]&7\ar@{-}[d]&6\ar@{-}[d]&5\ar[dr]&4\ar[dl]&3\ar[dr]&2\ar[dl]&1\ar[dr]&0\ar[dl]\\
 	9' &8'&7'&6'&5'&4'&3'&2'&1'&0'
 }
$$
The matching $\delta^2(M)_2$:
{\small
$$
\begin{gathered}
\left \{	\xymatrix@C=8pt@R=20pt{
	9 \ar@{-}[d] &8\ar[dr]&7\ar[dl]&6\ar@{-}[d]&5\ar@{-}[d]&4\ar@{-}[d]&3\ar[dr]&2\ar[dl]&1\ar[dr]&0\ar[dl]\\
	9' &8'&7'&6'&5'&4'&3'&2'&1'&0'
}
\right \}
\\ 
-\left \{	\xymatrix@C=8pt@R=20pt{
	9 \ar@{-}[d] &8\ar[dr]&7\ar[dl]&6\ar@{-}[d]&5\ar[dr]&4\ar[dl]&3\ar@{-}[d]&2\ar@{-}[d]&1\ar[dr]&0\ar[dl]\\
	9' &8'&7'&6'&5'&4'&3'&2'&1'&0'
}
\right \}
\\
+
\left \{	\xymatrix@C=8pt@R=20pt{
	9 \ar@{-}[d] &8\ar[dr]&7\ar[dl]&6\ar@{-}[d]&5\ar[dr]&4\ar[dl]&3\ar[dr]&2\ar[dl]&1\ar@{-}[d]&0\ar@{-}[d]\\
	9' &8'&7'&6'&5'&4'&3'&2'&1'&0'
}
\right \}
\end{gathered}
$$
}
Substituting into \eqref{delta2-formula} we obtain the diagrammatic expression for $\delta^2(M)$
{\small
$$
\begin{gathered}
	\delta^2(M)=
-\left\{	\xymatrix@C=8pt@R=20pt{
		9 \ar@{-}[d] &8\ar@{-}[d]&7\ar@{-}[d]&6\ar@{-}[d]&5\ar[dr]&4\ar[dl]&3\ar[dr]&2\ar[dl]&1\ar[dr]&0\ar[dl]\\
		9' &8'&7'&6'&5'&4'&3'&2'&1'&0'
	}
\right\}
\\
+
	\left \{	\xymatrix@C=8pt@R=20pt{
		9 \ar@{-}[d] &8\ar[dr]&7\ar[dl]&6\ar@{-}[d]&5\ar@{-}[d]&4\ar@{-}[d]&3\ar[dr]&2\ar[dl]&1\ar[dr]&0\ar[dl]\\
		9' &8'&7'&6'&5'&4'&3'&2'&1'&0'
	}
	\right \}
	\\ 
	-\left \{	\xymatrix@C=8pt@R=20pt{
		9 \ar@{-}[d] &8\ar[dr]&7\ar[dl]&6\ar@{-}[d]&5\ar[dr]&4\ar[dl]&3\ar@{-}[d]&2\ar@{-}[d]&1\ar[dr]&0\ar[dl]\\
		9' &8'&7'&6'&5'&4'&3'&2'&1'&0'
	}
	\right \}
\\
	+
	\left \{	\xymatrix@C=8pt@R=20pt{
		9 \ar@{-}[d] &8\ar[dr]&7\ar[dl]&6\ar@{-}[d]&5\ar[dr]&4\ar[dl]&3\ar[dr]&2\ar[dl]&1\ar@{-}[d]&0\ar@{-}[d]\\
		9' &8'&7'&6'&5'&4'&3'&2'&1'&0'
	}
	\right \}
\end{gathered}
$$
}
 \end{example} 

The main property is that the map $\delta:=\{\delta^d\}$  is the differential of the graded module $\ZZ[{\cal M}_{\widehat{PG}_{l}}]^{+}$.
\begin{pro}\label{pro:diff}
	The graded module $\ZZ[{\cal M}_{\widehat{PG}_{l}}]^{+}$ equipped with the maps $\{\delta^d\}$ is a differential graded module, that is, one has the equation
	$$
	\delta^{d+2}\circ \delta^d =0, \forall d.
	$$
\end{pro}
\begin{pf}
	Let $M$ be a matching of degree $d$, we may assume $d\leq l-4$, otherwise there is nothing to prove. Consider the complement
	$M\setminus M(0)$ of $0$-colored edges of $M$. These are disjoint uninterrupted sequences of braided pairs of $\pm$ edges; they are ranged
	in the order encountered as we walk counterclockwise along the circle
	starting from the first $0$-colored edge in $M(0)$:
	$$
	M_1,\ldots, M_k.
	$$ 
	This gives
	$$
	\delta^d(M)=\sum^k_{i=1}(-1)^i \delta^d(M)_i.
	$$
	Applying $\delta^{d+2}$ we obtain
	$$
	\delta^{d+2}\circ\delta^d(M)=\sum^k_{i=1}(-1)^i \delta^{d+2} (\delta^d(M)_i).
	$$
	By definition $\delta^d(M)_i$ is obtained from $M$ by inserting in place of $M_i$ the alternating sum $\delta^d(M_i)$ of unbraiding operations, while all other segments $M_j$, $j\neq i$, are unchanged.
	This tells us the matchings occurring the sum $\delta^{d+2} (\delta^d(M)_i)$ are as follows
	$$
	\delta^{d+2} (\delta^d(M)_i)=\sum_{j<i} (-1)^j [\ldots \delta(M_j)\ldots \delta(M_i) \ldots] + \sum_{j>i} (-1)^{j+1} [\ldots \delta(M_i)\ldots \delta(M_j) \ldots], 
	$$
	each bracket is a matching produced from $M$ by the unbraiding operation along the segments $M_i$ and $M_j$ and the dots represent the unchanged part of $M$; we use the property 
	$$
	\delta^2(M_i)=0:
	$$
	the operation $\delta(M_i)$ is modeled on the boundary operator for simplicies. The resulting sum
	$$
	\delta^{d+2}\circ\delta^d(M)=\sum_{j<i} (-1)^{i+j} [\ldots \delta(M_j)\ldots \delta(M_i) \ldots] + \sum_{j>i} (-1)^{i+j+1} [\ldots \delta(M_i)\ldots \delta(M_j) \ldots]=0
	$$
	because of the cyclic nature of the graph $\widehat{PG}_l$ every bracket appears twice with opposite signs.
	\end{pf} 
 The above result is an indication that the perfect matching of $\widehat{PG}_l$ should be categorified. We will return to this again once we explore how the cyclic nature of $\widehat{PG}_l$ affects the representations of the quiver $PG_l$.
  
\section{Quantum-type invariants and
	the algebraic K\"ahler structures of
	  $W_{\xi}$}

We have seen that the quiver ${PG}_{l}$ comes along with the family of representations
$$
\rho_c: \{\{P^s([\xi],[\phi])\}; \alpha^{t,s}_c(\xi,\phi): P^s([\xi],[\phi]) \longrightarrow P^t([\xi],[\phi])\},
$$
attached to every point $([\xi],[\phi])\in {\mathfrak{L}}_l$; recall  $P^s([\xi],[\phi])$ are the summands of the orthogonal decomposition
$$
W_{\xi}/W^{l}_{\xi} ([\phi])=\bigoplus^{l-1}_{s=0}  P^s([\xi],[\phi]),
$$
for every nonzero $\xi$ and $\phi$ lying over $[\xi]$ and $[\phi]$ respectively, the maps  $\alpha^{t,s}_c(\xi,\phi)$ are the blocks of the map
$$
\alpha^{(2)}_{\xi,c} (\phi, \bullet): W_{\xi} \longrightarrow \HKC
$$
decorating the edge $(s)\rightarrow (t)'$ of the  quiver ${PG}_{l}$,  and $c$ is the trace parameter, that is the value
$$
tr([\xi],[\phi]) (\alpha^{(2)}_{\xi,c} (\phi, \bullet))=\sum^{l-1}_{s=0} tr(\alpha^{s,s}_c(\xi,\phi))=c.
$$

\noindent
In the previous sections we extended the graph ${PG}_{l}$ to its trivalent completion   
 $\widehat{PG}_{l}$. The issue arises
 
 \begin{equation}\label{Q:inducedreps}
 	\begin{gathered}
 	\text{\it can we extend in a meaningful way the defining representations
 		$\rho_c$ of ${PG}_{l}$}
 	\\
 	\text{\it to the representations of the
 		trivalent quiver $\widehat{PG}_{l}$?}
 	\end{gathered}
 \end{equation} 
With this question in mind we will construct interesting operators on the graded space
$$
W_{\xi}/W^{l}_{\xi} ([\phi])=\bigoplus^{l-1}_{s=0}  P^s([\xi],[\phi]).
$$
which produce decorations for the additional edges $(0) \rightarrow (l-1)'$ and $(l-1) \rightarrow (0)'$ of the trivalent quiver $\widehat{PG}_{l}$.

\vspace{0.2cm}
We will be working on the stratum ${\mathfrak L}_l$ of $\PP({\cal W}_{\Sigma^0_r})$, parametrizing points $([\xi],[\phi])$ where the length $l_{\xi}([\phi])$ of the $([\xi],[\phi])$-filtration of $W_{\xi}$  is constant equals $l$; it will be always assumed that $l$ is at least $3$. We have our graph $\widehat{PG}_{l}$ attached to such a stratum.

All edges of $\widehat{PG}_{l}$, with the exception of two edges, $e^+_{l-1}=(l-1)\to 0'$ and $e^-_0=0 \to (l-1)'$, are decorated by linear maps $\alpha^{t,s}_c(\xi,\phi)$, the summands of $\alpha^{(2)}_{\xi,c}(\phi,\bullet)$ of the representations $\rho_c(\xi\otimes\phi)$ of the graph $PG_l$; $\xi\otimes\phi$'s are nonzero vectors of the total space $\OO_{{\mathfrak L}_l}(-1)$ lying over $([\xi],[\phi])$ in ${\mathfrak L}_l$. In addition, we recall that only the diagonal components
$
\alpha^{s,s}_c (\xi,\phi)
$
depend on $c$ and this dependence is described by the formula
$$
\alpha^{s,s}_{c'}(\xi,\phi) -\alpha^{s,s}_c(\xi,\phi) =\frac{c'-c}{dim(W_{\xi}/W^l_{\xi} ([\phi]))} {\bf 1}_{P^s},
$$
for any two values $c'$ and $c$ of the trace parameter and for all $s\in [0,l-1]$.

We set
$$
P ([\xi],[\phi]):=Hom(P^0,P^{l-1}) \oplus Hom(P^{l-1},P^0),
$$
where we write $P^s$ instead of $P^s ([\xi],[\phi])$ to simplify the notation. The summands in the direct sum above will be denoted
$P^- ([\xi],[\phi])$ and $P^+ ([\xi],[\phi])$ respectively:
$$
P^- ([\xi],[\phi]):= Hom(P^0,P^{l-1}), \,\, P^+ ([\xi],[\phi]):= Hom(P^{l-1},P^0).
$$
The superscripts are chosen to match the color of edges $e^-_0=0 \to (l-1)$ and $e^+_{l-1}=(l-1)\to 0'$. The main point of the constructions presented below is to answer the question in \eqref{Q:inducedreps}: given a representation $\rho_c(\xi\otimes\phi)=\{\alpha^{t,s}_c(\xi,\phi)\}$ of the graph $PG_l$, we will assign linear maps in $P^{\mp} ([\xi],[\phi])$ to two remaining arrows $e^-_0$ and $e^+_{l-1}$ to obtain representations of the graph $\widehat{PG}_{l}$. In fact the result obtained is somewhat more general:

 we consider the category ${\mathfrak{Reps}}^{\bullet\circ}(PG_l)$ (resp.  ${\mathfrak{Reps}}^{\bullet\circ}(\widehat{PG}_l)$) of finite dimensional {\it bipartite} representations of $PG_l$ (resp. $\widehat{PG}_l$ ); a {\it bipartite representation} means that each pair $(s')$ and $(s)$ of black-white vertices of $PG_l$ is decorated with {\it the same} vector space; the defining representations $\rho_c(\xi\otimes\phi)$ attached to points of $\OO^{\times}_{{\mathfrak{L}}_l}(-1)$ are bipartite (in the sequel, the  representations considered throughout the text are bipartite, so we often drop the
 bicolored superscripts in the notation of the categories); an upshot is to construct a functor
$$
{\mathfrak{Reps}}^{\bullet\circ}(PG_l) \longrightarrow  {\mathfrak{Reps}}^{\bullet\circ}(\widehat{PG}_l)
$$
from the category ${\mathfrak{Reps}}^{\bullet\circ}(PG_l)$ to the category ${\mathfrak{Reps}}^{\bullet\circ}(\widehat{PG}_l)$ such that the composition with the obvious restriction functor
$$
{\cal R}: {\mathfrak{Reps}}^{\bullet\circ}(\widehat{PG}_l) \longrightarrow {\mathfrak{Reps}}^{\bullet\circ}(PG_l)
$$
is the identity functor:
$$
\xymatrix{
{\mathfrak{Reps}}^{\bullet\circ}(PG_l) \ar[r] \ar@/_1.5pc/[rr]_{{\cal ID}}& {\mathfrak{Reps}}^{\bullet\circ}(\widehat{PG}_l) \ar[r]^{{\cal R}} &{\mathfrak{Reps}}^{\bullet\circ}(PG_l)
}
$$
The construction presented below exhibits the whole family of such functors. They are
parametrized by pairs of linear functionals on the space of Laurent polynomials $\CC[q,q^{-1}]$, where $q$ is an auxiliary parameter. 
Our main tool will be zig-zag paths of $\widehat{PG}_{l}$ described in Lemma \ref{lem:zzpaths}. 
 
 \begin{pro}\label{pro:mu-map}
 	Let  $\rho=\{\alpha^{t,s}:P^s \longrightarrow P^t\}$ be a finite dimensional bipartite representation of the quiver $PG_l$.
 	The zig-zag paths $Z_{0,-}$ and $Z_{0,+}$ give rise, respectively, to distinguished linear maps
 	$$
 	\begin{gathered}
 		\tau^+_{\rho}: Hom(\CC[q,q^{-1}],\CC) \longrightarrow P^+ :=Hom(P^{l-1},P^0),
 		\\
 		\tau^-_{\rho}: Hom(\CC[q,q^{-1}],\CC) \longrightarrow P^-:=Hom(P^0,P^{l-1}),
 	\end{gathered}
 $$
 where $q$ is a formal parameter and $\CC[q,q^{-1}]$ the algebra of Laurent polynomials in $q$. This gives the map
 $$
 \widehat{\rho}: Hom(\CC[q,q^{-1}],\CC)\times  Hom(\CC[q,q^{-1}],\CC) \longrightarrow \{\text{Representations of $\widehat{PG}_l$ extending $\rho$}\}:
 $$
 given a pair $(F,G)$ of linear functionals on $\CC[q,q^{-1}]$ the representation  $\widehat{\rho}(F,G)$ coincides with $\rho$ on $PG_l$
 and assigns to the edges $e^{-}_0$ and $e^+_{l-1}$, respectively, the linear maps
 $$
 \text{$\tau^-_{\rho}(F):P^{0} \longrightarrow P^{l-1}$ and $\tau^+_{\rho}(G):P^{l-1} \longrightarrow P^{0}$.}
 $$
 In particular, using the natural inclusion
 $$
 \CC^{\times} \subset  Hom(\CC[q,q^{-1}],\CC)
 $$
 given by the evaluation of Laurent polynomials at a point of $\CC^{\times}$, one obtains the maps
 	$$
 	\mu_{\rho}: \CC^{\times} \times \CC^{\times} \longrightarrow P =P^- \oplus P^+  =Hom(P^0,P^{l-1}) \oplus Hom(P^{l-1},P^0).
 	$$
 	In other words, for every point $(u,v)$ of the algebraic torus $\CC^{\times} \times \CC^{\times}$ the value $\mu_{\rho}(u,v)$ is a pair of linear maps
 	$$
 	\tau^-_{\rho}(u):P^0 \longrightarrow P^{l-1} ,\,\, \tau^+_{\rho}(v):P^{l-1} \longrightarrow P^{0} ,
 	$$ 
 	decorating the arrows $e^-_0$ and $e^+_{l-1}$ and thus defining  a representation of $\widehat{PG}_{l}$ extending $\rho$; that representation will be denoted $\widehat{\rho}(u,v)$. 
 \end{pro}

The proof will be a result of a sequence of constructions presented below.

\vspace{0.2cm}
{\it Step 1.} We have a natural nondegenerate bilinear pairing between two summands
\begin{equation}\label{P+-pairing}
B:	P^- \times P^+ =Hom(P^0,P^{l-1}) \times Hom(P^{l-1},P^0) \longrightarrow Hom(P^0,P^0) \stackrel{tr}{\longrightarrow} \CC.
\end{equation}
The first arrow above is the composition of maps and the second is the trace:
$$
B(f,g):=trace(g\circ f),\,\, \forall f\in  P^- ,\,\, \forall g\in P^+.
$$
This implies the following.
\begin{lem}\label{lem:p+dualp-}
	In the direct sum decomposition
	$$
	P =P^-  \oplus P^+ 
	$$
	the two summands are naturally dual to each other under the bilinear pairing $B$ in \eqref{P+-pairing}
	$$
B^-:	P^-  \cong (P^+)^{\ast},\,\,B^+:P^+  \cong (P^-)^{\ast}:
	$$
	the isomorphisms are given respectively by the formulas
	$$
	\begin{gathered}
	P^-  \ni f \mapsto B(f,\bullet) \in (P^+ )^{\ast},
	\\
	P^+  \ni g \mapsto B(\bullet,g) \in (P^-)^{\ast}.
\end{gathered}
	$$
	In particular, we have natural isomorphisms
	$$
	P \cong P^-  \oplus \left(P^- \right)^{\ast} \cong P^+  \oplus (P^+ )^{\ast}.
	$$
\end{lem}

\begin{rem}\label{rem:sympP}
1)	The last statement can be rephrased as a natural identification
	of $P $ with the cotangent bundle $T^{\ast}_{P^{\pm}}$ of $P^{\pm}$: 
	\begin{equation}\label{two-symp-ident}
		P  \cong T^{\ast}_{P^{-}} \cong T^{\ast}_{P^{+}}.
	\end{equation}
		Thus $P$ is naturally a symplectic space with the symplectic form given by the skew-symmetric pairing
		$$
		\omega_{P }((f\oplus g),(f'\oplus g') ):=B(f,g') -B(f',g),
		$$
		where $s \oplus t$ denotes an element in $P =P^- \oplus P^+ $. 
		
		The two identifications in \eqref{two-symp-ident} correspond to switching the order of summands
		$$
		P =P^- \oplus P^+  \cong P^+ \oplus P^- :=P^{op}.
		$$
		The corresponding symplectic form on $P $ becomes
		$$
		\begin{gathered}
		\omega'((f\oplus g),(f'\oplus g') )=\omega_{P^{op}} ((g\oplus f),(g'\oplus f'))\\
		 =B(f',g) -B(f,g')=-\omega_{P}((f\oplus g),(f'\oplus g') ),
	\end{gathered}
		$$
		that is, two symplectic forms differ by a sign. Later on we will be interested in Lagrangians in $P$. From the above it follows that those are the same independently of the symplectic structure chosen.
		
		2) Though it will not be used in this writing, it should be mentioned that the representations of $PG_l$ attached to points 
		$([\xi],[\phi])$ of ${\mathfrak{L}_l}$
		$$
		\rho_c: \{\{P^s([\xi],[\phi])\}; \alpha^{t,s}_c(\xi,\phi): P^s([\xi],[\phi]) \longrightarrow P^t([\xi],[\phi])\},
		$$
		 where $P^s([\xi],[\phi])$ are the summands of the orthogonal decomposition
		$$
		W_{\xi}/W^{l}_{\xi} ([\phi])=\bigoplus^{l-1}_{s=0}  P^s([\xi],[\phi]),
		$$
	 have additional structure: this is the quaternionic structure of Nakajima, see \cite{Na}. Namely, 
		$$
		\begin{gathered}
		P ([\xi],[\phi])=P^-([\xi],[\phi]) \oplus P^+([\xi],[\phi]))\\
		=Hom(P^0([\xi],[\phi]), P^{l-1}([\xi],[\phi])) \oplus Hom(P^0([\xi],[\phi]), P^{l-1}([\xi],[\phi]))
		\end{gathered}
		$$ 
		carries a natural quaternionic structure: 
		since $P ([\xi],[\phi])$ is the direct sum of complex spaces we have the defining complex structure, the multiplication by $\imath =\sqrt{-1}$:
		$$
		I(f\oplus g)=\imath f \oplus \imath g,
		$$
		for all $f\in P^- ([\xi],[\phi]) $ and all $g \in P^+ ([\xi],[\phi])$;
		another complex structure on $P ([\xi],[\phi])$ is given by the rule
		$$
		J (f\oplus g)=(-g^{\dag}) \oplus f^{\dag}, 
		$$
		where $A^{\dag}$ denotes the homomorphism adjoint to $A$ - recall, all summands $P^{\bullet}$ are equipped with the Hermitian metric coming from the Hodge metric on $\HKC$ - the adjoint is taken with respect to this Hermitian metric;
		
		 setting $K=IJ$ gives the third complex structure and  it is straightforward to check that $(I,J,K)$ are subject to the relations of quaternions:
		 $$
		 I^2=J^2=K^2=-id_{P ([\xi],[\phi])}, \,\, IJ=-JI=K.
		 $$
\end{rem}

{\it Step 2.} We will now use the zig-zag paths $Z_{0,-}$ and $Z_{0,+}$, see Lemma \ref{lem:zzpaths}, to define linear maps
\begin{equation}\label{Cpm-step2}
C^{\pm}_{\rho}:  P^{\pm}\longrightarrow
End^{(0)} (W_{\rho}),
\end{equation}
where $W_{\rho}$ is the direct sum of the spaces $\{P^{\bullet}\}$ of the representation $\rho$:
\begin{equation}\label{W=dirsumPs}
	W_{\rho}:=\bigoplus^{l-1}_{s=0}P^s;
\end{equation}
on the right hand side of \eqref{Cpm-step2} the superscript $(0)$ means grading preserving endomorphisms of $W_{\rho}$ viewed as a graded vector space. The two maps, $C^{\pm}_{\rho}$,  are defined similarly.
We give the detailed description of the map
$$
C^{-}_{\rho} :  P^{-} \longrightarrow
End^{(0)} (W_{\rho}).
$$
Recall the zig-zag path $Z_{0,-}$:
$$
Z_{0,-}=(-e^-_0)e^0_{l-1} (-e^-_{l-1})e^0_{l-2} \cdots (-e^-_2)e^0_1 (-e^-_1)e^0_0.
$$
This is a closed path in  $\widehat{PG}_{l}$ starting and ending at the vertex $0$. We begin by reversing its orientation, that is, we consider the path
$$
-Z_{0,-}=(-e^0_0)e^-_1(-e^0_1)e^-_2 \cdots e^-_{l-2}(-e^0_{l-2})e^-_{l-1}(-e^0_{l-1})e^-_0.
$$
 Next, we label the edges $(-e^0_j)$ with the endomorphisms
$$
\alpha^{j,j} :P^j \longrightarrow P^j,
$$
for every $j\in [0,l-1]$, and the edges $e^-_j$ with the endomorphisms
$$
\alpha^{j-1,j} :P^{j} \longrightarrow P^{j-1},
$$
 for $j\in [1, l-1]$. 
  All the edges of $(-Z_{0,-})$ except $e^-_0$ are now labeled by linear maps and their composition takes us from $P^{l-1}$ to  $P^0$ as presented in the following diagram
  
 $$
 \xymatrix@C=35pt@R=35pt{
 	&&&&&P^0\\
 &&&&P^1 \ar[r]^{\alpha^{0,1}}&P^0\ar[u]_{\alpha^{0,0}}\\
 	&&\cdots \cdots&P^2 \ar[r]^{\alpha^{1,2}} &P^1 \ar[u]_{\alpha^{1,1}}&\\ 
 	&P^{l-2}\ar[r]^{\alpha^{l-3,l-2}}&P^{l-3}\ar[u]_{\alpha^{l-3,l-3}}&&&\\
 	P^{l-1} \ar[r]^{\alpha^{l-2,l-1}}&P^{l-2}\ar[u]_{\alpha^{l-2,l-2}}&&&&\\
 P^{l-1} \ar[u]_{\alpha^{l-1,l-1}}&&&&&
 }
$$

 The remaining edge $e^-_0$ of $(-Z_{0,-})$ is labeled with
 $X\in Hom(P^0, P^{l-1})$ thus completing the above diagram
 to the following
  
 $$
 \xymatrix@C=35pt@R=35pt{
 	&&&&&&P^0\\
 	&&&&&P^1 \ar[r]^{\alpha^{0,1}}&P^0\ar[u]_{\alpha^{0,0}}\\
 	&&&\cdots \cdots&\!\!\!\!\!\!\!P^2 \ar[r]^{\alpha^{1,2}} &P^1 \ar[u]_{\alpha^{1,1}}&\\ 
 	&&P^{l-2}\ar[r]^{\alpha^{l-3,l-2}}&P^{l-3}\ar[u]_{\alpha^{l-3,l-3}}&&&\\
 	&P^{l-1} \ar[r]^{\alpha^{l-2,l-1}}&P^{l-2}\ar[u]_{\alpha^{l-2,l-2}}&&&&\\
 	P^0 \ar[r]^X&P^{l-1} \ar[u]_{\alpha^{l-1,l-1}}&&&&&
 }
 $$
 
 We define $(C^{-})^0_{\rho}(X)$ by joining $P^0$ at the bottom  with $P^0$ at the top by the composite arrow

 $$
 \xymatrix@C=35pt@R=35pt{
 	&&&&&&P^0\\
 	&&&&&P^1 \ar[r]^{\alpha^{0,1}}&P^0\ar[u]_{\alpha^{0,0}}\\
 	&&&\cdots \cdots&\!\!\!\!\!\!\!\!\!\!\!\!\!P^2 \ar[r]^{\alpha^{1,2}} &P^1 \ar[u]_{\alpha^{1,1}_c}&\\ 
 	&&P^{l-2}\ar[r]^{\alpha^{l-3,l-2}}&P^{l-3}\ar[u]_{\alpha^{l-3,l-3}}&&&\\
 	&P^{l-1} \ar[r]^{\alpha^{l-2,l-1}}&P^{l-2}\ar[u]_{\alpha^{l-2,l-2}}&&&&\\
 	P^0 \ar[r]^X \ar@/^8pc/ [-5,6]&P^{l-1} \ar[u]_{\alpha^{l-1,l-1}}&&&&&
 }
 $$
Written as a formula, we have:
 $$
 (C^{-})^0_{\rho}(X)=\alpha^{0,0}\alpha^{0,1} \alpha^{1,1} \alpha^{1,2} \cdots \alpha^{l-2,l-1} \alpha^{l-1,l-1}X :P^0 \longrightarrow P^0.
 $$
This is the component of $C^{-}_{\rho}(X)$ in degree $0$.  Shifting cyclically, from left to right, we obtain
 $$
 \begin{gathered}
(C^{-})^1_{\rho}(X):=rshift((C^{-})^0_{\rho}(X))
\\
= \alpha^{1,1} \alpha^{1,2} \cdots \alpha^{l-2,l-1} \alpha^{l-1,l-1}X \alpha^{0,0}\alpha^{0,1}  :P^1 \longrightarrow P^1,
\end{gathered}
$$
the component of $C^{-}_{\rho}(X)$ in degree $1$.
Performing the right shift $i$ times gives the component of  $C^{-}_{\rho}(X)$ in degree $i$. We now define
$$
C^{-}_{\rho}(X)=\circlearrowleft((C^{-})^0_{\rho}(X)):=\sum^{l-1}_{i=0} rshift^i ((C^{-})^0_{\rho}(X)).
$$

The construction of  $C^{+}_{\rho}$ is similar. The only difference is that we use the zig-zag path $Z_{0,+}$:
$$
Z_{0,+}=e^+_{l-1}(-e^0_{l-1})e^+_{l-2}(-e^0_{l-2}) \cdots e^+_{2}(-e^0_{2})e^+_{1}(-e^0_{1})e^+_0;
$$
 we continue to label $(-e^0_j)$ with $\alpha^{j,j}$, while to $e^+_j$ we assign $\alpha^{j+1,j}$, for each $j\in [0,l-2]$; the remaining edge $e^+_{l-1}$ is labeled by an element $Y \in Hom (P^{l-1},P^0)$;
 
  we define
 the endomorphism 
$$
C^{+}_{\rho}(Y): \bigoplus^{l-1}_{i=0}P^i \longrightarrow \bigoplus^{l-1}_{i=0}P^i
$$
to be grading preserving with the component of degree $0$ 
$$
(C^{+})^0_{\rho}(Y): P^0 \longrightarrow P^0
$$
 defined by the composition formula
 $$
 (C^{+})^0_{\rho}(Y)=Y\alpha^{l-1,l-1} \alpha^{l-1,l-2} \cdots \alpha^{2,1} \alpha^{1,1} \alpha^{1,0} \alpha^{0,0}.
 $$
 Shifting cyclically, from right to left, gives the next component
 $$
 (C^{+})^1_{\rho}(Y)=lshift((C^{+})^0_{\rho}(Y)):=\alpha^{1,0} \alpha^{0,0}Y\alpha^{l-1,l-1} \alpha^{l-1,l-2} \cdots \alpha^{2,1} \alpha^{1,1}: P^1 \longrightarrow P^1.
 $$
Continuing in this fashion gives 
 $$
C^+_{\rho}(Y)= \circlearrowright((C^{+})^0_{\rho}(Y)):=\sum^{l-1}_{i=0} lshift^i ((C^{+})^0_{\rho}(Y)).
 $$
 
 To simplify the notation we let
 \begin{equation}\label{g0-Lieal}
 	{\bf \mathfrak{g}^{(0)}_{\rho}}:=End^{(0)} (W_{\rho})=\bigoplus^{l-1}_{i=0} End (P^i).
 \end{equation}

With this notation the above constructions give linear maps
\begin{equation}\label{C+-maps}
	\begin{gathered}
	C^-_{\rho}: P^- =Hom (P^0, P^{l-1}) \longrightarrow {\bf \mathfrak{g}^{(0)}_{\rho}},\\
	C^+_{\rho}: P^+ =Hom (P^{l-1}, P^{0}) \longrightarrow {\bf \mathfrak{g}^{(0)}_{\rho}}.
	\end{gathered}
\end{equation}
assigned to every bipartite representation $\rho=\{\alpha^{t,s}\}$ of the quiver $PG_l$ via a sort of `integration' along zig-zag paths $Z_{0,-}$ and $Z_{0,+}$ respectively.

\vspace{0.2cm}
{\it Step 3.} We now introduce a formal deformation parameter $q$ and consider the loop Lie algebra
$$
{\bf \mathfrak{g}^{(0)}_{\rho}}[q,q^{-1}]:={\bf \mathfrak{g}^{(0)}_{\rho}}\otimes_{\CC} \CC[q,q^{-1}]
$$
of Laurent polynomials in $q$ with coefficients in the Lie algebra ${\bf \mathfrak{g}^{(0)}_{\rho}}$. The maps in \eqref{C+-maps} extend to the linear maps
	\begin{equation}\label{C+-maps-q}
		\begin{gathered}
			C^-_{\rho}(q, \bullet): P^- =Hom (P^0, P^{l-1}) \longrightarrow {\bf \mathfrak{g}^{(0)}_{\rho}}[q,q^{-1}],\\
			C^+_{\rho} (q, \bullet): P^+ =Hom (P^{l-1}, P^{0}) \longrightarrow {\bf \mathfrak{g}^{(0)}_{\rho}}[q,q^{-1}].
		\end{gathered}
	\end{equation}
Namely, for $X \in P^-$ we put
$$
	C^-_{\rho}(q, X):=\sum^{l-1}_{j=0} (C^-)^j_{\rho}( X) q^{(-1)^j j},
	$$
	a sort of the Poincar\'e polynomial of $C^-_{\rho}( X)$, where $(C^-)^j_{\rho}( X)$ is the component of $C^-_{\rho}( X)$ in degree $j$. 
	
	Analogously, for $Y\in P^+$, we have
	$$
		C^+_{\rho}(q, Y):=\sum^{l-1}_{j=0} (C^+)^j_{\rho}( Y) q^{(-1)^j j},
	$$
	where $(C^+)^j_{\rho}( Y)$ is the component of $C^+_{\rho}( Y)$ in degree $j$.
	
	\vspace{0.5cm}
{\it Step 4.} In this step of our construction we will need the weighted trace map
	\begin{equation}\label{tr(t,.)}
		tr(t, \bullet): {\bf \mathfrak{g}^{(0)}_{\rho}} \longrightarrow \CC[ t],
	\end{equation}
where $t$ is another formal variable; it is defined by the formula
$$
tr(t,A):= \sum ^{l-1}_{m=0} tr(A^m) t^m,
$$
where for $A\in {\bf \mathfrak{g}^{(0)}_{\rho}}$ we use the decomposition of $A$ into its graded pieces according to the direct sum on the right hand side of the equation in \eqref{g0-Lieal}:  
$$
A=\sum^{l-1}_{m=0} A^m, \,\,\,\, A^m: P^m \longrightarrow P^m, \forall m\in[0,l-1].
$$
Composing the maps in \eqref{C+-maps-q} with the weighted trace map gives
\begin{equation}\label{C+-maps-q-t}
	\begin{gathered}
		C^-_{\rho}(q,t, \bullet): P^- =Hom (P^0, P^{l-1}) \longrightarrow \CC[q,q^{-1},t],\\
		C^+_{\rho} (q,t, \bullet): P^+ =Hom (P^{l-1}, P^{0}) \longrightarrow \CC[q,q^{-1},t].
	\end{gathered}
\end{equation}
Setting $t=1$, we deduce two linear maps
\begin{equation}\label{C+-maps-q-1}
	\begin{gathered}
		C^-_{\rho}(q,1, \bullet): P^- =Hom (P^0, P^{l-1}) \longrightarrow \CC[q,q^{-1}],\\
		C^+_{\rho} (q,1, \bullet): P^+ =Hom (P^{l-1}, P^{0}) \longrightarrow \CC[q,q^{-1}].
	\end{gathered}
\end{equation}
Dualize the above maps:
\begin{equation}\label{dualC+-maps-q-1}
	\begin{gathered}
		\tau^{+}_{\rho}:=	(C^-)^{\ast}_{\rho}(q,1, \bullet): Hom_{\CC} (\CC[q,q^{-1}], \CC) \longrightarrow (P^-)^{\ast} \cong P^+ ,\\
		\tau^{-}_{\rho}:=	(C^+)^{\ast}_{\rho} (q,1, \bullet): Hom_{\CC} (\CC[q,q^{-1}], \CC) \longrightarrow (P^+)^{\ast} \cong P^-.
	\end{gathered}
\end{equation}
The last isomorphism in each line above is the one from Lemma \ref{lem:p+dualp-}. The maps in \eqref{dualC+-maps-q-1} are the ones stated in Proposition \ref{pro:mu-map}. We will often omit the subscript for the representation in the definition of the above maps when no confusion is likely.

The linear space $Hom_{\CC} (\CC[q,q^{-1}], \CC)$ contains the subset of maps which are ring homomorphisms from $\CC[q,q^{-1}]$ to $\CC$, that is, the set of maximal ideals of $\CC[q,q^{-1}]$ or, equivalently, the set of closed points of the affine scheme $Spec (\CC[q,q^{-1}])$. That scheme is $\CC^{\times}$. Thus we have an obvious inclusion
$$
\CC^{\times} \subset Hom_{\CC} (\CC[q,q^{-1}], \CC).
$$
The restrictions of maps in \eqref{dualC+-maps-q-1} to this subset gives the maps
\begin{equation}\label{mu+-}
	\begin{gathered}
		\mu^{+}_{\rho}:=\tau^{+}_{\rho}|_{\CC^{\times}}: \CC^{\times} \subset Hom_{\CC} (\CC[q,q^{-1}], \CC) \longrightarrow (P^-)^{\ast} \cong P^+ ,\\
		\mu^{-}_{\rho}:=\tau^{-}_{\rho}|_{\CC^{\times}}:	\CC^{\times} \subset Hom_{\CC} (\CC[q,q^{-1}], \CC) \longrightarrow (P^+)^{\ast} \cong P^-.
	\end{gathered}
\end{equation}
Taking the Cartesian product gives the map
\begin{equation}\label{mu}
	\mu_{\rho}:=(\mu^-_{\rho}, \mu^+_{\rho}): \CC^{\times}\times \CC^{\times} \longrightarrow P^- \oplus P^+
\end{equation}
asserted in Proposition \ref{pro:mu-map}. 
${\Box}$

\vspace{0.2cm}
Recall the toric surface $X(\Delta_l)$, see Proposition \ref{pro:toricsurf}.
We can identify $\CC^{\times}\times \CC^{\times}$ with the open part of
 $X(\Delta_l)$, the complement of the toric divisors
 $$
 D_{0,+}, \, D_{0,-},  \, D_{\mp}
 $$
 corresponding respectively to the one dimensional cones
 $\RR_{+} [Z_{0,+}] $, $\RR_{+} [Z_{0,-}] $, $\RR_{+} v_{\mp} $ of $\Delta_l$, that is, we make identification
 $$
 \CC^{\times}\times \CC^{\times} \cong \stackrel{\circ}{X}(\Delta_l):=X(\Delta_l) \setminus \left(D_{0,+} \cup D_{0,-} \cup D_{\pm} \right).
 $$
 With this identification in mind we can restate Proposition \ref{pro:mu-map} as follows.
 \begin{cor}\label{cor:mu-map}
 		Let  $\rho=\{\alpha^{t,s}:P^s\longrightarrow P^t\}$, a finite dimensional bipartite representation of the quiver $PG_l$, be given. Then there is a distinguished
 		 map
 		$$
 		\mu_{\rho}: \stackrel{\circ}{X}(\Delta_l) \longrightarrow P(\rho) :=P^- (\rho) \oplus P^+ (\rho), 
 		$$
 		where $P^-(\rho) =Hom(P^0,P^{l-1})$ and $P^{+}= Hom(P^{l-1},P^0)$.
 		In other words, for every point $x$ in  $\stackrel{\circ}{X}(\Delta_l)$ the value $\mu_{\rho}(x)$ is a pair of linear maps
 		$$
 		\mu^-_{\rho}(x):P^0 \longrightarrow P^{l-1} ,\,\, \mu^+_{\rho}(x):P^{l-1} \longrightarrow P^{0} ,
 		$$ 
 		decorating the arrows $e^-_0$ and $e^+_{l-1}$ and thus defining a representation of the graph $\widehat{PG}_{l}$ extending $\rho$. 
 \end{cor}

We now specialize Proposition \ref{pro:mu-map} to the family of representations of $PG_l$ attached to the stratum $\mathfrak{L}_l$. Recall that for a value $c\in \CC$ of the trace parameter we have the map
\begin{equation}\label{rhoc-quantum}
\rho_c: \OO^{\times}_{\mathfrak{L}_l}(-1) \longrightarrow {\mathfrak{Reps}}(PG_l)
\end{equation}
which associates to a point $\xi \otimes \phi$ lying over $([\xi],[\phi])\in \mathfrak{L}_l$ the representation of $PG_l$
$$
\rho_c(\xi\otimes\phi)=\{\{P^s([\xi],[\phi])\}_s; \alpha^{t,s}_c(\xi,\phi):P^s([\xi],[\phi])\longrightarrow P^t([\xi],[\phi])\}:
$$
it labels every pair of vertices $(s)$ and $(s')$ with the vector spaces $P^s([\xi],[\phi])$ of the orthogonal decomposition
$$
W_{\xi}/W^l_{\xi}([\phi])=\bigoplus^{l-1}_{s=0} P^s([\xi],[\phi])
$$
and each edge $(s) \rightarrow (t') $ of the graph is decorated with the linear map $\alpha^{t,s}_c (\xi,\phi)$, the $(t,s)$-block of the map $\alpha^{(2)}_{\xi}(\phi,\bullet)$ in the collection $A_{[\xi]}([\phi])(\xi,\phi)$ having the trace $c$, see Lemma \ref{lem:traceparam}. 
\begin{cor}\label{cor:rhoc-ext} For every value of the trace parameter $c\in \CC$, the family of representations
	$$
	\rho_c: \OO^{\times}_{\mathfrak{L}_l}(-1) \longrightarrow {\mathfrak{Reps}}(PG_l)
	$$
	admits extensions to the representations of the graph $\widehat{PG}_l$.
	Namely, we have the map
	$$
	\widehat{\rho_c}: Hom(\CC[q,q^{-1}], \CC)\times Hom(\CC[q,q^{-1}], \CC) \longrightarrow {\mathfrak{Reps}}(\widehat{PG}_l)(\rho_c),
	$$
	where ${\mathfrak{Reps}}(\widehat{PG}_l)(\rho_c)$ denotes the representations of $\widehat{PG}_l$ extending $\rho_c$: for every pair $(F,G)$ of linear functionals on $\CC[q,q^{-1}]$, we have the family
	of representations
	$$
	\widehat{\rho_c}(F,G):\OO^{\times}_{\mathfrak{L}_l}(-1) \longrightarrow {\mathfrak{Reps}}(\widehat{PG}_l)
	$$
	extending the family $\rho_c$. More precisely, for every point $\xi\otimes\phi $ in $\OO^{\times}_{\mathfrak{L}_l}(-1)$ lying over a point
	$([\xi],[\phi]) \in {\mathfrak{L}_l}$ the representation
	$\widehat{\rho_c}(F,G)(\xi\otimes \phi)$ coincides with $\rho_c(\xi\otimes \phi)$ on $PG_l$ and labels the extra edges $e^-_0$ and $e^+_{l-1}$ by the maps
	$$
	\begin{gathered}
	\tau^{-}_{\rho_c(\xi\otimes \phi)}(F): P^0([\xi],[\phi])\longrightarrow P^{l-1}([\xi],[\phi]),
	\\
	\tau^{+}_{\rho_c(\xi\otimes \phi)}(G): P^{l-1}([\xi],[\phi])\longrightarrow P^{0}([\xi],[\phi]),
	\end{gathered}
$$
provided by Proposition \ref{pro:mu-map}. 
Specializing to the linear functionals 
$$
\CC^{\times} \subset Hom(\CC[q,q^{-1}], \CC)
$$
corresponding to the evaluation of Laurent polynomials at a point of $\CC^{\times}$ gives the maps
$$
\mu_{\rho_c(\xi\otimes \phi)} :\CC^{\times}\times \CC^{\times} \longrightarrow P([\xi],[\phi])=P^-([\xi],[\phi]) \oplus P^+([\xi],[\phi]),
$$
for all nonzero $\xi\otimes \phi$ lying over a point $([\xi],[\phi])\in\mathfrak{L}_l$; the summands $P^{\pm}([\xi],[\phi])$ are respectively
$$
\begin{gathered}
P^-([\xi],[\phi])=Hom(P^0([\xi],[\phi]), P^{l-1}([\xi],[\phi])),
\\
 P^+([\xi],[\phi])=Hom( P^{l-1}([\xi],[\phi]),P^0([\xi],[\phi])).
 \end{gathered}
$$
\end{cor}
\vspace{0.2cm}
Perhaps the reason for the adjective `quantum' in the title of this section should be explained. Recall, that we have realized our graph $\widehat{PG}_l$ on a topological torus $\mathbb{T}$. In particular, the vertices of the graph are points on the torus. Consider the pair of vertices $(l-1)$ and $(0)$ of the graph as points of $\mathbb{T}$ and cut out small disks around them. The resulting surface can be viewed as a cobordism between two boundary circles of those disks, call them $C_{l-1}$ and $C_{0}$ respectively. A representation of $\rho=\{\{P^s\}, \{\alpha^{t,s}\}\}$ of the quiver $PG_l$ assigns the vector spaces $P^{l-1}$ and $P^{0}$ to $C_{l-1}$ and to $C_{0}$ respectively. The formalism of the Topological Quantum Field Theory (TQFT), see \cite{At}, tells us that cobordisms between circles correspond to maps between linear spaces attached to circles. The constructions of the maps $\tau^{\pm}_{\rho}$ in \eqref{dualC+-maps-q-1} are reminiscent of this formalism. Hence the choice of the terminology. In addition, the map $\mu_{\rho}$ in Proposition \ref{pro:mu-map} could be viewed as a distinguished `toric' family of such maps realizing a cobordism given by a torus with two boundary components. More generally, the appearance of 
moduli for labeling  the morphisms between $P^{l-1}$ and $P^{0}$ fits with the heuristic principle of quantum mechanics - to consider all 
possible paths between two given points. In the situation at hand we have the `classical paths' between $P^{l-1}$ and $P^{0}$ given by composing in the usual way the homomorphisms of the representation $\rho=\{\alpha^{t,s}\}$ of the quiver $PG_l$:
\begin{equation}\label{classic}
\begin{gathered}
\rho (-Z'_{0,-}) : P^{l-1} \longrightarrow P^0,
\\
\rho (Z'_{0,+}) : P^{0} \longrightarrow P^{l-1},
\end{gathered}
\end{equation}
where the path $-Z'_{0,-}$ (resp. $Z'_{0,+})$ is the the path $-Z_{0,-}$ (resp. $Z_{0,+}$) {\it without} the edge $e^-_0$ (resp. $e^+_{l-1}$) and $\rho (-Z'_{0,-})$ and $\rho (Z'_{0,+})$ denote the compositions
\begin{equation}
\begin{gathered}
\rho (-Z'_{0,-})=\alpha^{0,0} \alpha^{0,1}\alpha^{1,1} \alpha^{1,2} \cdots \alpha^{l-2,l-2} \alpha^{l-2,l-1} \alpha^{l-1,l-1},
\\
\rho (Z'_{0,+})=   \alpha^{l-1,l-1}\alpha^{l-2,l-1} \alpha^{l-2,l-2} \alpha^{l-2,l-3} \cdots \alpha^{2,1} \alpha^{1,1} \alpha^{1,0} \alpha^{0,0}. 
\end{gathered}
\end{equation}
  The evaluation of the maps $C^{\pm}_{\rho}(q,1,\bullet)$ in \eqref{C+-maps-q-1} gives `classical' endomorphisms in ${\mathfrak g}^{(0)}_{\rho}[q,q^{-1}]$:
  $$
  \text{$C^-_{\rho} (q,1,\rho (Z'_{0,+}))$ and $C^+_{\rho} (q,1,\rho (-Z'_{0,-}))$.}
  $$ 
  Our construction shows that there is a `moduli' of `nonclassical' paths, that is, there is a lot of `quantum' ways to go between $P^{l-1}$ and $P^0$. That `moduli' is identified as the space of linear functionals
$Hom_{\CC} (\CC[q,q^{-1}], \CC)$ on the algebra of the Laurent polynomials
in $q$:
$$
	\begin{gathered}
		\tau^{+}_{\rho}:=	(C^-)^{\ast}_{\rho}(q,1, \bullet): Hom_{\CC} (\CC[q,q^{-1}], \CC) \longrightarrow (P^-)^{\ast} \cong P^+ =Hom_{\CC} (P^{l-1},P^0),\\
		\tau^{-}_{\rho}:=	(C^+)^{\ast}_{\rho} (q,1, \bullet): Hom_{\CC} (\CC[q,q^{-1}], \CC) \longrightarrow (P^+)^{\ast} \cong P^- =Hom_{\CC} (P^{0},P^{l-1}).
	\end{gathered}
$$
Observe that the map $\mu_{\rho}$ in Proposition \ref{pro:mu-map} makes a `preferred' choice in the space of functionals - the evaluation maps at the $\CC$-closed points of $Spec (\CC[q,q^{-1}])$. Certainly, we left out many other interesting choices. We will now discuss this aspect.

\vspace{0.2cm}
 We begin by recalling that a noncanonical isomorphism of vector spaces
\begin{equation}\label{Res}
	\CC[q,q^{-1}] \cong  Hom_{\CC} (\CC[q,q^{-1}], \CC)
\end{equation}
can be given by the symmetric bilinear pairing
$$
\CC[q,q^{-1}] \times \CC[q,q^{-1}] \longrightarrow \CC
$$
defined by the formula
$$
(f,g):=Res_0 (fg), \forall f,g \in \CC[q,q^{-1}],
$$
where $Res_0 (h)$ is the residue at $0$ of a Laurent polynomial $h$.
For $f \in \CC[q,q^{-1}]$ we denote by 
$$
R(f)=\text{the image of $f$ under the isomorphism \eqref{Res},}
$$
that is, $R(f)$ is determined by the formula
$$
R(f)(g)=(f,g)=Res_0 (fg), \,\forall g\in \CC[q,q^{-1}].
$$

	 The identification \eqref{Res} of the ring of Laurent polynomials $\CC[q,q^{-1}]$ with its dual recasts the maps $\tau^{\pm}_{\rho}$ in Proposition \ref{pro:mu-map} as the linear maps
	 \begin{equation}\label{Lpoly-Ppm}
	 	\begin{gathered}
	 	\tau^{+}_{\rho}: \CC[q,q^{-1}]\cong Hom_{\CC} (\CC[q,q^{-1}], \CC) \longrightarrow (P^-)^{\ast} \cong P^+ ,\\
	 	\tau^{-}_{\rho}: \CC[q,q^{-1}] \cong Hom_{\CC} (\CC[q,q^{-1}], \CC) \longrightarrow (P^+)^{\ast} \cong P^-,
	 \end{gathered}
 \end{equation}
 meaning that the Laurent polynomials in the variable $q$ become parameters for extending each representation $\rho$ of $PG_l$ to representations of the trivalent graph $\widehat{PG}_l$. We organize the maps above as two generating series
 {\boldmath
 \begin{equation}\label{Lpoly-Ppm-genseries}
 		\tau^{\pm}_{\rho} \mbox{\unboldmath{$(T,q)= \displaystyle \sum_{k\in \ZZ} \tau^{\pm}_{\rho}(q^k) T^k$},}
 \end{equation}}

\noindent
where $T$ is another formal variable; the above series could be viewed as a sort of formal Fourier series which encompasses all labelings of the edges $e^+_{l-1}$ and $e^-_0$ obtained from a representation $\rho$ of $PG_l$. The Fourier coefficients $\tau^{\pm}_{\rho}(q^k) $ are subject to the following.
\begin{pro}\label{pro:coefFourier}
	The homomorphisms $\tau^{\pm}_{\rho}(q^k) \in P^{\pm}$ are uniquely determined
	from the relation
	$$
	tr (\tau^{\pm}_{\rho}(q^k) \circ X)= Res_0 (q^{k} C^{\mp}_{\rho}(q,1, X)),\,\, \forall X \in P^{\mp},
	$$
	where $C^{\mp}_{\rho}(q,X)$ is ${\mathfrak g}^{(0)}$-valued Laurent polynomial
	$$
	C^{\mp}_{\rho}(q,X)= \sum^{l-1}_{j=0} (C^{\mp}_{\rho})^j (X)q^{(-1)^j j},
	$$
	see \eqref{C+-maps-q}, and $C^{\mp}_{\rho}(q,1, X)$ is the trace of $C^{\mp}_{\rho}(q,X)$:
	$$
	C^{\mp}_{\rho}(q,1, X))=\sum^{l-1}_{j=0} tr \left((C^{\mp}_{\rho})^j (X) \right) q^{(-1)^j j}.
	$$
	In particular, $\tau^{\pm}_{\rho}(q^k)=0$, unless $k=-1-2j$ or $k=2j$, for 
	$2j \in [0,l-1]$. In the latter cases we have
	$$
	tr (\tau^{\pm}_{\rho}(q^k) \circ X)=\begin{cases}
		tr((C^{\mp}_{\rho})^{2j}(X)),& \text{if $k=2j$,}\\
		tr((C^{\mp}_{\rho})^{2j+1}(X)),& if k=-1-2j.
	\end{cases}
$$
Thus the generating series in \eqref{Lpoly-Ppm-genseries} are $P^{\pm}$-valued Laurent polynomials
{\boldmath
	\begin{equation}\label{Lpoly-Ppm-genLp}
		\tau^{\pm}_{\rho} \mbox{\unboldmath{$(T,q)= \displaystyle \sum^{\left[\frac{l-1}{2} \right]}_{j=0} \tau^{\pm}_{\rho}(q^{2j}) T^{2j} +\sum^{\left[\frac{l-1}{2} \right]}_{j=0} \tau^{\pm}_{\rho}(q^{-1-2j}) T^{-1-2j}$},}
\end{equation}
}
\end{pro}
\begin{pf}
	On the left side of \eqref{Lpoly-Ppm} a Laurent polynomial $f$ is identified with the linear functional
	$$
	R(f): \CC[q,q^{-1}] \longrightarrow \CC
	$$
	defined by the formula
	$$
		R(f)(g)=Res_0 (fg), \forall g \in \CC[q,q^{-1}].
		$$
	The arrows in \eqref{Lpoly-Ppm} are dual to the maps
	$$
	C^{\mp}_{\rho}(q,1,\bullet): P^{\mp} \longrightarrow \CC[q,q^{-1}].
	$$
	Hence the linear functional
	$$
	R(f) \circ C^{\mp}_{\rho}(q,1,\bullet) :P^{\mp} \longrightarrow \CC
	$$
	defined by the formula
	\begin{equation}\label{f-functional}
	P^{\mp} \ni X \mapsto Res_0 (f C^{\mp}_{\rho}(q,1,X)) \in \CC.
\end{equation}
	The identifications
	$$
	(P^{\mp})^{\ast} \cong P^{\pm}
	$$
	on the right side of \eqref{Lpoly-Ppm} are done via the nondegenerate bilinear paring
	$$
	B: P^- \times P^+ \longrightarrow \CC
	$$
	by the formula
	$$
	B(X,Y)=tr(Y\circ X).
	$$
	Hence for the linear functional in \eqref{f-functional} there is a unique element denoted $\tau^{\pm}_{\rho}(f)$ in $P^{\pm}$ such that
	\begin{equation}\label{formulaRes=tr}
	Res_0 (f C^{\mp}_{\rho}(q,1,X))=tr (\tau^{\pm}_{\rho}(f)\circ X), \forall X \in P^{\mp}.
	\end{equation}
	The first assertion of the proposition follows by specializing $f=q^k$.
	Other statements follow easily from the definition of $C^{\mp}(q,1,X)$ and the calculation of residues.
\end{pf}

\vspace{0.2cm}
There is also a principle of `specialization' or `localization' at the points of particular interest, which often accompanies the quantum mechanical principle `summing over all histories'. We have already mentioned that given a representation $\rho=\{\{P^s\},\{\alpha^{t,s}\}\}$ of $PG_l$ the spaces $P^{+}=Hom(P^{l-1}, P^0)$ and $P^{-}=Hom(P^0,P^{l-1})$ come along with distinguished maps
$$
\rho(-Z'_{0,-}) \in P^+,\,\, \rho(Z'_{0,+}) \in P^-,
$$
coming from the `classical' paths $(-Z'_{0,-})$ and $Z'_{0,+}$ in $PG_l$ between the vertices $(0)$ and $(l-1)$, see \eqref{classic}. Hence the distinguished Laurent polynomials
\begin{equation}\label{Lp-initial}
	C^-_{\rho}(q,1,\rho(Z'_{0,+})),\,\,C^+_{\rho}(q,1,\rho(-Z'_{0,-})),
\end{equation}
attached to every representation $\rho$ of $PG_l$.
With the identification
$$
Hom(\CC[q,q^{-1}],\CC) \cong \CC[q,q^{-1}]
$$
we can create many more maps in $P^{\pm}$ as well as Laurent polynomials in $\CC[q,q^{-1}]$. Namely, with the identification above we have the diagram
\begin{equation}\label{manufacture}
	\xymatrix{
	P^- \ar[dr]^{C^-_{\rho}(q,1,\bullet)}&& &P^- \\
&\CC[q,q^{-1}] \ar@{=}[r]& Hom(\CC[q,q^{-1}],\CC) \ar[ur]^{\tau^-_{\rho}} \ar[dr]_{\tau^{+}_{\rho}}& \\
P^+ \ar[ur]_{C^{+}_{\rho}(q,1,\bullet)}&& &P^+
}
\end{equation}
which manufactures sequences of maps in $P^{\pm}$ and Laurent polynomials in $\CC[q,q^{-1}]$ starting from an initial datum. Thus with the initial data
$$
f^-_0:=C^-(q,1,\rho(Z'_{0,+})), \,\, f^+_0:=C^+(q,1,\rho(-Z'_{0,-})),
$$
which correspond to the initial maps
$$
A^-_0:=\rho(Z'_{0,+}) \in P^-,\,\, A^+_0:=\rho(-Z'_{0,-}) \in P^+
$$
we produce the `descendants' for maps
$$
\begin{gathered}
\text{$A^{+,-}_1 =\tau^+_{\rho} (f^-_0), \,\, A^{+,+}_1 =\tau^+_{\rho} (f^+_0),$ in $P^+$,}
\\
\text{$A^{-,-}_1 =\tau^-_{\rho} (f^-_0), \,\, A^{-,+}_1 =\tau^-_{\rho} (f^+_0),$ in $P^-$.}
\end{gathered}
$$
Applying the maps $C^{\pm}_{\rho}(q,1,\bullet)$ gives the Laurent polynomial `descendants':
$$
\begin{gathered}
f^{+,-}_1=C^{+}_{\rho}(q,1,A^{+,-}_1), \,\, f^{+,+}_1=C^{+}_{\rho}(q,1,A^{+,+}_1), 
\\
f^{-,-}_1=C^{-}_{\rho}(q,1,A^{-,-}_1), \,\, f^{+,+}_1=C^{-}_{\rho}(q,1,A^{-,+}_1),
\end{gathered}
$$
Continuing in this fashion we obtain that the initial data $A^-_0$ (resp. $A^+_0$) generates the tree of maps and Laurent polynomials. Below is an illustration of the tree of three `generations' with the initial datum $A^-_0$:
 
\begin{equation}\label{tree}
	 \begin{tikzpicture}
		[place/.style={circle,draw=red,thick, inner sep=0pt, minimum size=2mm},	transition/.style={circle,draw=black,fill=black, inner sep=0pt,minimum size=2mm}]
		\node (A-0) [place] at (-6,0)  [label={left:{$A^-_0$}}]{};
		\node (f-0)  [transition] at (-5,0)  [label={right:{ $f^-_0$}}]{};
		\node[fill=red] (A+-1) [place] at (-4.5,1)  [label={above:{$\scriptstyle{A^{+,-}_1}$}}]{};
		\node (A--1) [place] at (-4.5,-1)  [label={below:{ $\scriptstyle{A^{-,-}_1}$}}]{};
		\node (f+-1)  [transition] at (-3,1)  [label={above :{$\,\scriptscriptstyle{f^{+,-}_1}$}}]{};
		\node (f--1) [transition] at (-3,-1) [label={above:{$\scriptscriptstyle{f^{-,-}_1}$}}]{};
			\node[fill=red] (A+--2) [place] at (-2.5,-0.5) [label={above :{$\,\,\,\,\,\scriptscriptstyle{A^{+,-,-}_2}$}}]{};
			\node (A---2) [place] at (-2.5,-1.5) [label={below:{$\scriptscriptstyle{A^{-,-,-}_2}$}}]{};
			\node (f+--2) [transition] at (-1.5,-0.5) [label={}]{};
			\node (f---2) [transition] at (-1.5,-1.5) [label={}]{};
			\node[fill=red] (A++--3) [place] at (-0.5,-0.25) [label={}]{};
			\node (A-+--3) [place] at (-0.5,-0.75) [label={}]{};
			\node[fill=red] (A+---3) [place] at (-0.5,-1.25) [label={}]{};
			\node (A----3) [place] at (-0.5,-1.75) [label={}]{};
		\node[fill=red] (A++-2) [place] at (-2.5,1.5)  [label={above:{ $\,\,\,\scriptstyle{A^{+,+,-}_2}$}}]{};
		\node (A-+-2) [place] at (-2.5,0.5)  [label={above:{ $\,\,\,\,\,\,\,\,\,\,\scriptscriptstyle{A^{-,+,-}_2}$}}]{};
		\node[fill=red] (A+++-3) [place] at (-0.5,1.75) [label={}] {};
		\node (A-++-3)[place] at (-0.5,1.25) [label={}]{};
		\node (f++-2) [transition] at (-1.5,1.5) [label={}] {};
		\node (f-+-2)[transition] at (-1.5,0.5) [label={}]{};
		\node[fill=red] (A+-+-3) [place] at (-0.5,0.75) [label={}] {};
		\node (A--+-3)[place] at (-0.5,0.25) [label={}]{};
		\node (f+++-3) [transition] at (0.5,1.75) [label={above right:$\scriptscriptstyle{f^{+,+,+,-}_3}$}] {};
		\node (f-++-3) [transition] at (0.5,1.25) [label={ right:$\scriptscriptstyle{f^{-,+,+,-}_3}$}] {};
		\node (f+-+-3) [transition] at (0.5,0.75) [label={right:$\scriptscriptstyle{f^{+,-,+,-}_3}$}] {};
		\node (f--+-3) [transition] at (0.5,0.25) [label={right:$\scriptscriptstyle{f^{-,-,+,-}_3}$}] {};
		\node (f++--3) [transition] at (0.5,-0.25) [label={right:$\scriptscriptstyle{f^{+,+,-,-}_3}$}] {};
		\node (f-+--3) [transition] at (0.5,-0.75) [label={right:$\scriptscriptstyle{f^{-,+,-,-}_3}$}] {};
		\node (f+---3) [transition] at (0.5,-1.25) [label={right:$\scriptscriptstyle{f^{+,-,-,-}_3}$}] {};
		\node (f----3) [transition] at (0.5,-1.75) [label={right :$\scriptscriptstyle{f^{-,-,-,-}_3}$}] {};
		\begin{scope}
				[decoration={markings,
						mark =at position 0.6cm with {\arrow[black,line width=0.5mm]{stealth}}}]
		\draw[postaction={decorate}] 
		(A+++-3) to (f+++-3);
		\draw[postaction={decorate}] 
		(A-++-3) to (f-++-3);
		\draw[postaction={decorate}] 
		(A+-+-3) to (f+-+-3);
		\draw[postaction={decorate}] 
		(A--+-3) to (f--+-3);
		\draw[postaction={decorate}] 
		(A++--3) to (f++--3);
		\draw[postaction={decorate}] 
		(A-+--3) to (f-+--3);
		\draw[postaction={decorate}] 
		 (f++-2) to (A+++-3);
		\draw[postaction={decorate}] 
		 (f++-2) to (A-++-3);
		\draw[postaction={decorate}] 
		(A++-2) to (f++-2); 
		\draw[postaction={decorate}] 
		(A+---3) to (f+---3);
		\draw[postaction={decorate}] 
		(A----3) to (f----3);
		\draw[postaction={decorate}] 
		(A-+-2) to (f-+-2);
		\draw[postaction={decorate}]  
		(f-+-2) to (A+-+-3);
		\draw[postaction={decorate}] 
		(A-+-2) to (f-+-2);
		\draw[postaction={decorate}] 
		(f-+-2) to (A--+-3);
		\draw[postaction={decorate}] 
		(A-0)--(f-0);
		\draw[postaction={decorate}]  
		(f-0)--(A+-1);
		\draw[postaction={decorate}] 
		(f-0)--(A--1);
		\draw[postaction={decorate}] 
		(f-0)--(A--1);
		\draw[postaction={decorate}] 
		(A--1)--(f--1);
		\draw[postaction={decorate}] 
		(f-0)--(A+-1);
		\draw[postaction={decorate}] 
		(A+-1)--(f+-1);
		\draw[postaction={decorate}] 
		(A+--2) to (f+--2);
		\draw[postaction={decorate}] 
		(f+--2) to (A++--3);
		\draw[postaction={decorate}] 
		(A+--2) to (f+--2);
		\draw[postaction={decorate}] 
		(f+--2) to (A-+--3);
		\draw[postaction={decorate}] 
		(A---2) to (f---2);
		\draw[postaction={decorate}] 
		(f---2) to (A+---3);
		\draw[postaction={decorate}] 
		(A---2) to (f---2);
		\draw[postaction={decorate}] 
		(f---2) to (A----3);
		\end{scope}
	\begin{scope}
		[decoration={markings,
			mark =at position 0.35cm with {\arrow[black,line width=0.5mm]{stealth}}}]
	\draw[postaction={decorate}] 
	(f--1) to (A+--2);
	\draw[postaction={decorate}] 
	(f--1) to (A---2);
	\draw[postaction={decorate}] 
	(f+-1) to (A++-2);
	\draw[postaction={decorate}] 
	(f+-1) to (A-+-2);
	\end{scope}
	\end{tikzpicture}	
\end{equation}
on the drawing the red circled nodes - {\textcolor{red} 
	{$\circ$}} - are maps in $P^-$, the red dotted nodes - {\textcolor{red} 
	{$\bullet$}} - are maps in $P^+$, the black dotted nodes - {\textcolor{black} 
	{$\bullet$}} - are Laurent polynomials.
The nodes of maps (resp. polynomials) are labeled by sequences of pluses and minuses. With the convention of assigning to `plus' the value $1$ and  to `minus' the value $(-1)$ we deduce
\begin{equation}
	\begin{gathered}
	\text{\it on the $n$-th level of the tree one obtains the collection of}
	\\
	\text{\it - $2^{n-1}$ maps in $P^-$,}
	\\
	\text{\it - $2^{n-1}$ maps in $P^+$,}
	\\
	\text{\it - $2^n$ polynomials in $\CC[q,q^{-1}]$;}
	\\
	\text{\it each collection is in bijection with the vertices of the cube $[-1,1]^n$ in $\RR^n$.}
		\end{gathered}
\end{equation}

Let us recall the nonabelian Dolbeault space ${\bf H^{1,0}}(PG_l) $ associated to the graph $PG_l$: its fan is precisely the cube 
$[-1,1]^{l-1}$ in $\RR^{l-1}$; the vertices of the fan are in bijection with toric divisors of ${\bf H^{1,0}}(PG_l) $, that is, irreducible components of the Lagrangian cycle $H_0$. Thus we obtain the following. 
\begin{pro}\label{pro:trees}
	Let $Z'_{0,+}$ and $(-Z'_{0,-})$ be two paths in $PG_l$ connecting the vertices $(0)$ and $(l-1)$. Given a representation $\rho=\{\alpha^{s',s}:P^s \longrightarrow P^{s'}\}$ of the quiver $PG_l$
	we obtain trees $T(A^{\pm}_0)$ of the form depicted in \eqref{tree}:
	
	- the root vertex of $T(A^{-}_0)$ (resp. $T(A^{+}_0)$) is labeled by the map $A^{-}_0=\rho(Z'_{0,+}) \in P^-$ (resp. $A^{+}_0=\rho(-Z'_{0,-}) \in P^+$),
	
	- the labels of black vertices ($\bullet$) form collection of polynomials in $\CC[q,q^{-1}]$,
	
	- the labels of its red circled vertices (\textcolor{red}{ $\circ$}) form the collection of maps in $P^-$,
	
	- the labels of its red dotted vertices (\textcolor{red}{ $\bullet$}) form the collection of maps in $P^+$,
	
	- every pair of maps descending from a polynomial labeling a black vertex of the tree defines the extension of the representation $\rho$ to a representation of the quiver $\widehat{PG}_l$,
	
	- at each level $n$ of the tree, the maps/polynomials labeling the nodes are in bijection with the vertices of the cube $[-1,1]^n$ in $\RR^n$; in particular, at the level $(l-1)$ the collection of maps/polynomials is in bijection with the irreducible component of the Lagrangian cycle $H_0$ of the nonabelian Dolbeault space ${\bf H^{1,0}}(PG_l) $.	
\end{pro}

\begin{rem}\label{rem:evolutionDolb}
	1) Informally speaking, the quantum invariants create a cascade of maps between $P^0$ and $P^{l-1}$ starting from two `classical' paths and the maps $A^{\pm}_0$ representing them; those give extensions of a representation of $PG_l$ to the ones of $\widehat{PG}_l$. After $(l-1)$ steps, the nodes of the resulting tree are in bijection with the toric divisors of the nonabelian Dolbeault space ${\bf H^{1,0}}(PG_l) $ - to each toric divisor of ${\bf H^{1,0}}(PG_l) $ is attached a particular quantum link between $P^0$ and $P^{l-1}$ and its Laurent polynomial. Actually, as the preceding result says, it so at {\rm every} level of the trees - as though the trees keep track of {\rm all previous and future developments of the nonabelian Dolbeault spaces}.
	
	2) The labeling $\pm1$ of `plus' and `minus' can be changed to the binary $1$ and $0$:
	$$
	+ \rightarrow 1,\,\, - \rightarrow 0.
	$$
	Then at each level $n$ of the trees $T(A^{\pm}_0)$ the collection of maps/Laurent polynomials are in bijection with the vector space ${\mathbb{F}}^n_2$ over the field ${\mathbb{F}}_2$, or equivalently, with the finite field ${\mathbb{F}}_{2^n}$, the extension of ${\mathbb{F}}_2$ of degree $n$.
\end{rem}

There is another way to connect the constructions above with toric divisors of ${\bf H^{1,0}}(PG_l) $. This is based on the fact that we can link $P^0$ and $P^{l-1}$ in a representation $\rho$ of $PG_l$ by {\it broken} versions of classical paths. Namely, we can delete any number of $0$-labeled edges in $Z'_{0,+}$ (resp. $(-Z'_{0,-})$) and still obtain a valid map in $P^-$ (resp. $P^+$): for a subset $S$ of the ordered set
$[1,l-1]$ denote by $Z'_{0,+}(S)$ (resp. $(-Z'_{0,-}(S))$) the broken path
obtained from  $Z'_{0,+}$ (resp. $(-Z'_{0,-})$) by deleting the edges
$e^0_s$, for all $s$ in $S$; given the representation
 $$
\rho=\{\alpha^{t,s}:P^s \longrightarrow P^{t}\}
$$
of the quiver $PG_l$ we obtain the homomorphisms
$$
\begin{gathered}
A^-(S)_0:=\rho  (Z'_{0,+}(S)): P^0 \longrightarrow P^{l-1},
\\
A^+(S)_0:=\rho  (-Z'_{0,-}(S)): P^{l-1} \longrightarrow P^{0},
\end{gathered}
$$
 defined by the compositions
 $$
 \begin{gathered}
 	A^-(S)_0=(\alpha^{l-1,l-1} \alpha^{l-1,l-2} \alpha^{l-2,l-2} \alpha^{l-2,l-3}  \cdots \alpha^{2,2} \alpha^{2,1} \alpha^{1,1} \alpha^{1,0})_{\widehat{S}}: P^0 \longrightarrow P^{l-1},
 	\\
 	A^+(S)_0=(\alpha^{0,0}\alpha^{0,1}\alpha^{1,1} \alpha^{l,2} \alpha^{2,2}  \cdots \alpha^{l-3,l-2} \alpha^{l-2,l-2} \alpha^{l-2,l-2} \alpha^{l-2,l-1})_{\widehat{S}}: P^{l-1} \longrightarrow P^0,
 \end{gathered}
$$
where the subscript $\widehat{S}$ means that in the compositions above the maps $\alpha^{s,s}$, for all $s\in S$, are deleted. Next recall that the irreducible components $\Pi_S$ of the Lagrangian cycle $H_0$ are also labeled by the subsets $S$ of $[1,l-1] $. Applying the diagram \eqref{manufacture} to the initial data $A^{\pm}(S)_0$ gives us the trees $T(A^{\pm}(S)_0)$ with the properties of Proposition \ref{pro:trees}, for every irreducible component $\Pi_S$ of $H_0$.
\begin{pro}\label{pro:Lagrangians-trees}
	Let $\rho=\{\alpha^{t,s}: P^s \longrightarrow P^t\}$ be a representation of the quiver $PG_l$.
	Then 
	for every subset $S$ of the ordered set $[1,l-1]$, the irreducible component $\Pi_S$ of the Lagrangian cycle $H_0$ of the nonabelian Dolbeault variety comes equipped with the trees $T(A^{\pm}(S)_0)$ of the form \eqref{tree} and the properties listed in Proposition \ref{pro:trees}.
\end{pro} 

\begin{rem}\label{rem:tauSpm}
	It should be clear that the results of Proposition \ref{pro:mu-map} remain valid if in the constructions used in the proof of the proposition we replace the zig-zag paths $Z_{0,\pm}$ by their `broken' versions. More precisely, for every subset $S$ of the ordered set $[0,l-1]$ consider the paths $Z_{0,\pm}(S)$ obtained from $Z_{0,\pm}$ by deleting the edges $e^0_s$, for all $s\in S$. Then given a bipartite finite dimensional representation $\rho=\left\{\{P^s\}, \{\alpha^{t,s}:P^s\longrightarrow P^{t}\}\right\}$ of $PG_l$
	one uses the paths $Z_{0,\pm}(S)$ to create the operators
	$$
	{}^S C^{\pm}_{\rho}(q,\bullet): P^{\pm} \longrightarrow {\mathfrak{g}}^{(0)}_{\rho}[q,q^{-1}],
	$$
	where the notation $P^{\pm}$ and ${\mathfrak{g}}^{(0)}_{\rho}$ is the same as in \eqref{C+-maps-q}. Taking the trace and dualizing one arrives
	at the $S$-version of the maps $\tau^{\pm}_{\rho}$ in Proposition \ref{pro:mu-map}:
	$$
		{}^S\tau^{\pm}_{\rho}: Hom(\CC[q,q^{-1}],\CC)\longrightarrow P^{\pm}.
$$
Restricting to $\CC^{\times}\subset Hom(\CC[q,q^{-1}],\CC)$ gives the $S$-version of the map $\mu_{\rho}$ in Proposition \ref{pro:mu-map}:
$$
{}^S\mu_{\rho}: \CC^{\times}\times\CC^{\times}\longrightarrow P=P^-\oplus P^{+}.
$$
\end{rem}
\vspace{0.2cm}
 Let us also recall that there are some natural Laurent polynomials attached to the stratum ${\mathfrak{L}}_l$.
Namely, with the orthogonal decomposition
$$
W_{\xi}/W^{l}_{\xi} ([\phi])=\bigoplus^{l-1}_{s=0}  P^s([\xi],[\phi])
$$
we can associate the {\it Hilbert Laurent polynomial}
$$
h ([\xi],[\phi],q):=\sum^{l-1}_{r=0} dim(P^s([\xi],[\phi]))q^{(-1)^s s}.
$$
Recall the partition of the stratum ${\mathfrak{L}}_l$ into the strata 
${\mathfrak{L}}_l (h^l,\lambda)$, where the dimensions of the summands $\{P^s([\xi],[\phi])\}$ are constant and the ones with $s\in [0,l-1]$ arranged in the Young diagram $\lambda$. In particular, ${\mathfrak{L}}_l$ is partitioned into strata where the Laurent polynomial $h([\xi],[\phi],q)$
remains constant. This gives us {\it  apriori} a finite collection 
$$
H_{{\mathfrak{L}}_l}:=\{h_{{\mathfrak{L}}_l (h^l,\lambda)}\}
$$
of Laurent polynomials associated with ${\mathfrak{L}}_l$: the polynomials of the collection are indexed by the strata ${\mathfrak{L}}_l (h^l,\lambda)$. The discussion above gives the following.
\begin{pro}\label{pro:Hilbertreps}
	Let 
	$$
	\rho_c:\OO^{\times}_{\mathfrak{L}_l}(-1)\longrightarrow {\mathfrak Reps}(PG_l)
	$$
	the family of representations of trace $c$ attached to the stratum $\mathfrak{L}_l$, that is, for every point $([\xi],[\phi])$ in $\mathfrak{L}_l$ we have the defining representations 
	$$
	\rho_c(\xi\otimes\phi)=\{\{P^s([\xi],[\phi])\}, \alpha^{t,s}_c(\xi,\phi):P^s([\xi],[\phi]) \longrightarrow P^{t}([\xi],[\phi])\}
	$$
	for every nonzero vectors $\xi$ and $\phi$ lying over $[\xi]$ and $[\phi]$ respectively.
	 Let $H_{{\mathfrak{L}}_l}:=\{h_{{\mathfrak{L}}_l (h^l,\lambda)}\}$ be the collection of Hilbert Laurent polynomials associated with the partition of ${\mathfrak{L}}_l$ into the substrata
	${\mathfrak{L}}_l (h^l,\lambda)$. Then for each  $h_{{\mathfrak{L}}_l (h^l,\lambda)}$ of the collection, the homomorphisms
		$$
		\begin{gathered}
		\tau^{-}_{\rho_c(\xi\otimes\phi)}(h_{{\mathfrak{L}}_l (h^l,\lambda)}):P^0 ([\xi],[\phi]) \longrightarrow  P^{l-1} ([\xi],[\phi]),\\ \tau^{+}_{\rho_c(\xi\otimes\phi)}(h_{{\mathfrak{L}}_l (h^l,\lambda)}):P^{l-1} ([\xi],[\phi]) \longrightarrow  P^{0} ([\xi],[\phi]),
	\end{gathered}
		$$
		the values of the maps $\tau^{\pm}_{\rho_c(\xi\otimes\phi)}$, for $\rho= \rho_c(\xi\otimes\phi)$ in \eqref{Lpoly-Ppm}, on $h_{{\mathfrak{L}}_l (h^l,\lambda)}$, give the family of representations 
		$$
		\widehat{\rho}_{c,{\mathfrak{L}}_l (h^l,\lambda)}: \OO^{\times}_{\mathfrak{L}_l}(-1)\longrightarrow {\mathfrak Reps}(\widehat{PG}_l)
		$$
		 of the quiver $\widehat{PG}_l$ extending the family of representations $\rho_c$.
\end{pro}

From Corollary \ref{cor:rhoc-ext} it follows that whenever we have a distinguished collection of points in $\CC^{\times} \times \CC^{\times}$, the values of the map 
$$
\mu_{\rho_c(\xi\otimes\phi)}:\CC^{\times} \times \CC^{\times} \longrightarrow P([\xi],[\phi])=P^{-}([\xi],[\phi]) \oplus P^{+}([\xi],[\phi])
$$
on the points of the collection give rise to the extensions of a representation of the quiver ${PG}_l$ to the ones of $\widehat{PG}_l$.
As an example of this, we have the subgroup $\mathbb{G}_l$ of $l$-th roots of unity in $\CC^{\times}$. The subset
$$
\mathbb{G}_l \times \mathbb{G}_l \subset  \CC^{\times} \times \CC^{\times}
$$
is clearly an instance of a distinguished collection of points in
$\CC^{\times} \times \CC^{\times}$ relevant to our considerations.
Thus we deduce the following.
\begin{pro}\label{pro:lrootsreps}
		Let 
	$$
	\rho_c:\OO^{\times}_{\mathfrak{L}_l}(-1)\longrightarrow {\mathfrak Reps}(PG_l)
	$$
be the family of representations of trace $c$ of the quiver $PG_l$ attached to the stratum $\mathfrak{L}_l$.
Let $\mathbb{G}_l$ be the subgroup of $l$-th roots of unity in $\CC^{\times}$ and consider the subset
	$$
	\mathbb{G}_l \times \mathbb{G}_l \subset  \CC^{\times} \times \CC^{\times}.
	$$
	 Then for every $(\zeta,\zeta') \in \mathbb{G}_l \times \mathbb{G}_l$ the values
	$$
	\mu^{-}_{\rho_c(\xi\otimes \phi)}(\zeta):P^0 ([\xi],[\phi]) \longrightarrow  P^{l-1} ([\xi],[\phi]),\,\, \mu^{+}_{\rho_c(\xi\otimes \phi)}(\zeta'):P^{l-1} ([\xi],[\phi]) \longrightarrow  P^{0} ([\xi],[\phi]),
	$$
	of the maps $\mu^{\pm}_{\rho_c(\xi\otimes \phi)}$ in Corollary \ref{cor:rhoc-ext} on $\zeta$ and $\zeta'$, give the family of representations 
	$$
	\widehat{\rho}_{c,(\zeta,\zeta')}:\OO^{\times}_{\mathfrak{L}_l}(-1)\longrightarrow {\mathfrak Reps}(\widehat{PG}_l)
	$$
	  of the quiver $\widehat{PG}_l$ extending the representation $\rho_c$.
\end{pro}

Recall that the families of representations $\{\rho_c\}_{c\in \CC}$ have two actions on the trace parameter $c$: the natural $\CC^{\times}$-action 
$$
\CC^{\times}\times \CC \longrightarrow \CC
$$
sending $c\in \CC$ to $(vc)$, for every $v\in \CC^{\times}$, and the additive action
$$
\CC\times \CC \longrightarrow \CC
$$
given by the addition in $\CC$.  Let us see how those actions affect the extensions of the representations $\rho_c$ of the quiver $PG_l$ to the representations $\widehat{\rho_c}$ of the quiver $\widehat{PG}_l$.
\begin{pro}\label{pro:actions-on-trace}
		Let 
	$$
	\rho_c:\OO^{\times}_{\mathfrak{L}_l}(-1)\longrightarrow {\mathfrak Reps}(PG_l)
	$$
	be the family of representations of trace $c$ of the quiver $PG_l$ attached to the stratum $\mathfrak{L}_l$ and let
	$$
	\widehat{\rho_c}:Hom(\CC[q,q^{-1}],\CC)\times Hom(\CC[q,q^{-1}],\CC) \longrightarrow {\mathfrak{Reps}}(\widehat{PG}_l)(\rho_c)
	$$
	be the map extending $\rho_c$ to the representations of $\widehat{PG}_l$
	as in Corollary \ref{cor:rhoc-ext}.
	
	1) Let $v\in \CC^{\times}$ acting on $c$ by multiplication $c\mapsto (vc)$. Then the representation $\rho_{vc}$ is related to $\rho_c$ by the
	formula
	$$
	\rho_{vc}=v\rho_c,
	$$
	that is, for every $(\xi\otimes \phi)$ in $\OO^{\times}_{\mathfrak{L}_l}(-1)$ we have the representation
	 $$
	\rho_{vc}(\xi\otimes \phi)=\{\alpha^{t,s}_{vc}(\xi,\phi)\} =\{v\alpha^{t,s}_{c}(\xi,\phi)\}.
	$$
	For $F,G$ in $Hom(\CC[q,q^{-1}],\CC)$representation $\widehat{\rho}_{vc}(F,G)(\xi\otimes\phi)$ extending $\rho_{vc}(\xi\otimes \phi)$ coincides with $v\rho_{c}(\xi\otimes \phi)$ on all edges of $PG_l$ while the extra edges $e^-_0$ and $e^+_{l-1}$ are labeled by the maps
	$$
	\tau^-_{{\rho}_{vc}(\xi\otimes\phi)} (F)=v^{2l-1}\tau^-_{{\rho}_{c}(\xi\otimes\phi)} (F),\hspace{0.2cm} 	\tau^+_{{\rho}_{vc}(\xi\otimes\phi)} (G)=v^{2l-1}\tau^+_{{\rho}_{c}(\xi\otimes\phi)} (G),
	$$
	in other words the maps on extra edges of $\widehat{PG}_l$ are scaled by $(2l-1)$ power of $v$.
	
	2) Let $a\in \CC$ acting on $c$ by the translation $c\mapsto (c+a)$ and
	consider the family $\rho_c$ (resp., $\rho_{c+a}$) over the substratum $\mathfrak{L}_l(h^l,\lambda)$. Then for every $(\xi\otimes \phi)$ in $\OO^{\times}_{\mathfrak{L}_l(h^l,\lambda)}(-1)$ we have the representation
	$$
	\rho_{c+a}(\xi\otimes \phi)=\{\alpha^{t,s}_{c+a}(\xi,\phi)\}=\left\{\alpha^{t,s}_{c}(\xi,\phi)+\frac{a\delta_{t,s}}{|\lambda|}{\bf 1}_{P^s}\right\},
	$$
	where $\delta_{t,s}$ is the Kronecker symbol and $|\lambda|=g-r-h^l$ is the weight of the partition $\lambda$.
	
	For $F,G$ in $Hom(\CC[q,q^{-1}],\CC)$representation $\widehat{\rho}_{c+a}(F,G)(\xi\otimes\phi)$ coincides with $\rho_{c+a}(\xi\otimes \phi)$ on all edges of $PG_l$ while the extra edges $e^-_0$ and $e^+_{l-1}$ are labeled by the maps
	$$
	\begin{gathered}
	\tau^-_{{\rho}_{c+a}(\xi\otimes\phi)} (F)=\sum^{l}_{k=0}\left(\frac{a}{|\lambda|}\right)^k \left(A^k_{{\rho}_{c}(\xi\otimes\phi)}\right)^-(F) ,
	\\
		\tau^+_{{\rho}_{c+a}(\xi\otimes\phi)} (G)=\sum^{l}_{k=0}\left(\frac{a}{|\lambda|}\right)^k\left(A^k_{{\rho}_{c}(\xi\otimes\phi)} \right)^+ (G),
	\end{gathered}
	$$
	where the summands $\left(A^k_{{\rho}_{c}(\xi\otimes\phi)}\right)^{\pm}(\bullet) \in P^{\pm}([\xi],[\phi])$ are determined by the summation of the maps
	${}^S\tau^{\pm}_{\rho_c(\xi\otimes \phi)} (\bullet)$ over the subsets $S\subset [0,l-1]$ with $k$ elements:
	$$
	\begin{gathered}
		\left(A^k_{{\rho}_{c}(\xi\otimes\phi)}\right)^-(F)=\sum_{\substack{S\subset [0,l-1],\\ |S|=k}} {}^S\tau^{-}_{\rho_c(\xi\otimes \phi)} (F),
		\\
		\left(A^k_{{\rho}_{c}(\xi\otimes\phi)}\right)^+(G)=\sum_{\substack{S\subset [0,l-1],\\ |S|=k}} {}^S\tau^{+}_{\rho_c(\xi\otimes \phi)} (G);
	\end{gathered}
$$
for the definition of maps ${}^S\tau^{\pm}_{\rho}$ see Remark \ref{rem:tauSpm}.
\end{pro}
\begin{pf}
	For the multiplicative action the equation
	$$
	\rho_{vc}=v\rho_c
	$$
	has been seen before: it follows from the formula
	$$
	\alpha^{(2)}_{\xi,vc}(\phi,\bullet)=v\alpha^{(2)}_{\xi,c}(\phi,\bullet).
	$$
	This gives the relation for blocks 
	$$
	\alpha^{t,s}_{vc}(\xi,\phi)=v\alpha^{t,s}_c({\xi},\phi),
	$$
	that is, every arrow $\alpha^{t,s}_c({\xi},\phi)$ of the representation $\rho_c(\xi\otimes \phi)$ is scaled by $v$.
	
	The label $\tau^{-}_{\rho_{vc}(\xi\otimes\phi)}(F)$ is obtain from composing $2l$ labels of the zig-zag path $Z_{0,+}$  and $(2l-1)$ of those labels are formed by the maps
	$v\alpha^{s,s-1}_c({\xi},\phi)$, for $s\in[1,l-1]$, and $v\alpha^{s,s}_c({\xi},\phi)$, for $s\in[0,l-1]$. This implies
	$$
	\tau^{-}_{\rho_{vc}(\xi\otimes\phi)}(F)=v^{2l-1}\tau^{-}_{\rho_{c}(\xi\otimes\phi)}(F).
	$$
	The label $\tau^{+}_{\rho_{vc}(\xi\otimes\phi)}(G)$ is obtained
	 similarly from composing $2l$ labels of the zig-zag path $Z_{0,-}$  and $(2l-1)$ of those labels are formed by the maps
	 $v\alpha^{s-1,s}_c({\xi},\phi)$, for $s\in[1,l-1]$, and $v\alpha^{s,s}_c({\xi},\phi)$, for $s\in[0,l-1]$, hence the formula
	$$
	\tau^{+}_{\rho_{vc}(\xi\otimes\phi)}(G)=v^{2l-1}\tau^{+}_{\rho_{c}(\xi\otimes\phi)}(G).
	$$
	
	Next we turn to the additive action. We assume to be on the stratum ${\mathfrak L}_l (h^l, \lambda)$. We have learned that the translation of a trace parameter $c$ by
	$a$ gives the following change of maps
	$$
	\alpha^{(2)}_{\xi,c+a} (\phi,\bullet)=\alpha^{(2)}_{\xi,c}  (\phi,\bullet) +\frac{a}{|\lambda|} id_{W_{\xi}/W^l_{\xi}([\phi])}
	$$
	This implies the following equation for blocks
	$$
	\alpha^{t,s}_{c+a} (\xi,\phi)=\alpha^{t,s}_{c}  (\xi,\phi) +\frac{a \delta_{t,s}}{|\lambda|} {\bf 1}_{P^s}.
	$$
	This implies the following labeling of edges of $PG_l$:
	\begin{equation}\label{edgelabels-transl-pf}
	\begin{gathered}
	\pm e^0_s \rightarrow \alpha^{s,s}_{c}  (\xi,\phi) +\frac{a}{|\lambda|} {\bf 1}_{P^s},\,\forall s\in [0,l-1],
	\\
	e^+_s \rightarrow \alpha^{s+1,s}_{c}  (\xi,\phi), \,\forall s\in [0,l-2],
	\\
	e^-_s \rightarrow \alpha^{s-1,s}_{c}  (\xi,\phi), \,	\forall s\in [1,l-1].
	\end{gathered}
\end{equation}
	To label the edges $e^-_0$ and $e^+_{l-1}$ we have to unravel the construction of the operators
	$$
	C^{\pm}_{\rho_{c+a}(\xi\otimes \phi)}(q,\bullet): P^{\pm}\longrightarrow {\mathfrak{g}}^{(0)}[q,q^{-1}]
		$$
		corresponding to the zig-zag paths $Z_{0,\pm}$, see \eqref{C+-maps-q}. The two are completely analogous, so we provide the details for $C^-_{\rho_{c+a}(\xi\otimes \phi)}(q,\bullet)$.
		Recall that for this operator we use the zig-zag path $Z_{0,-}$ with the orientation reversed:
		$$
		-Z_{0,-}=(-e^0_0)e^-_1 (-e^0_1)e^-_2 \cdots e^-_{l-2}(-e^0_{l-1})e^-_0.
		$$
		To obtain the operator based at $P^0$ we replace the edges of $PG_l$ occurring in $-Z_{0,-}$ by the corresponding labels according to the assignment in \eqref{edgelabels-transl-pf} and the additional edge $e^-_0$ of
		$\widehat{PG}_l$ by an element $X\in P^-=Hom(P^0,P^{l-1})$:
		$$
		\begin{gathered}
		(C^-)^0_{\rho_{c+a}(\xi\otimes \phi)}(X)=
		\left(\alpha^{0,0}_c +\frac{a}{|\lambda|}{\bf 1}_{P^0}\right) \alpha^{0,1}_c \left(\alpha^{1,1}_c +\frac{a}{|\lambda|}{\bf 1}_{P^1}\right)\alpha^{1,2}_c \cdots 
		\\
		\cdots\left(\alpha^{l-2,l-2}_c +\frac{a}{|\lambda|}{\bf 1}_{P^{l-2,l-2}}\right) \alpha^{l-2,l-1}\left(\alpha^{l-1,l-1}_c +\frac{a}{|\lambda|}{\bf 1}_{P^{l-1,l-1}}\right) X
		\end{gathered}
		$$
		The expansion according  the powers of $\frac{a}{|\lambda|}$ gives
		$$
		\begin{gathered}
		(C^-)^0_{\rho_{c+a}(\xi\otimes \phi)}(X)=(C^-_{\rho_{c}(\xi\otimes \phi)})^0(X) +\frac{a}{|\lambda|} \sum_{s\in [0,l-1]}({}^sC^-_{\rho_{c}(\xi\otimes \phi)})^0(X) 
		\\
		+\left(\frac{a}{|\lambda|}\right)^2 \sum_{{s',s}\subset [0,l-1]}({}^{\{s',s\}}C^-_{\rho_{c}(\xi\otimes \phi)})^0(X) +\left(\frac{a}{|\lambda|}\right)^3 \sum_{{s'',s',s}\subset [0,l-1]}({}^{\{s'',s',s\}}C^-_{\rho_{c}(\xi\otimes \phi)})^0(X)+\cdots
	\end{gathered}
		$$
		The sum corresponding to the order $k$ term
		$$
		\left(\frac{a}{|\lambda|}\right)^k \sum_{\substack{S\subset [0,l-1]\\ |S|=k}}({}^SC^-_{\rho_{c}(\xi\otimes \phi)})^0(X)
		$$
		is taken over all subsets $S$ of $[0,l-1]$ with $k$ distinct elements; the operator ${}^SC^-_{\rho_{c}(\xi\otimes \phi)}$ is the
		one corresponding to the broken path $-Z_{0,-}(S)$ obtained from
		$-Z_{0,-}$ by omitting all arrows $(-e^0_s)$ with $s\in S$, see Remark \ref{rem:tauSpm}.
		
		 Shifting cyclically from left to right at every term of the expansion gives the terms based at $P^1$; continuing the process gives the formula
		 $$
		C^-_{\rho_{c+a}(\xi\otimes \phi)}(q,X)=\sum^{l}_{k=0}	\left(\frac{a}{|\lambda|}\right)^k \left(\sum_{\substack{S\subset [0,l-1]\\ |S|=k}}{}^SC^-_{\rho_{c}(\xi\otimes \phi)}(q,X)\right) 
		$$
		Taking the trace and dualizing we obtain
		$$
		\tau^+_{\rho_{c+a}(\xi\otimes\phi)}: Hom(\CC[q,q^{-1}],\CC) \longrightarrow P^+
			$$
			subject to the expansion
			$$
			\tau^+_{\rho_{c+a}(\xi\otimes\phi)}=\sum^{l}_{k=0}	\left(\frac{a}{|\lambda|}\right)^k \left(\sum_{\substack{S\subset [0,l-1]\\ |S|=k}}{}^S\tau^+_{\rho_{c}(\xi\otimes \phi)}\right) 
				$$
				asserted in the proposition \end{pf}

We will now use the quantum-type invariants   to construct Lagrangians in the symplectic vector space $P([\xi],[\phi])$. Before we go on with formal constructions, it is perhaps a moment to pause and discuss some guidelines of our explorations. Conceptually our refinement seems to be a part of 
higher categorical structures of IVHS invariants:

- we started with $[\xi]$ in $\Sigma^{\circ}_r$ of the Griffiths stratification, added to it the subspace $W_{\xi}$ of $\HKC$ annihilated by $\xi$; the pair $([\xi], W_{\xi})$ is the classical IVHS datum; this was further enriched by $([\xi], [\phi])$-filtrations of $W_{\xi}$ varying holomorphically with $[\phi]$ in $\PP(W_{\xi})$, the fibre over $[\xi]$ of the morphism $p_1$ of the cohomological incidence correspondence
$$
\xymatrix{
&{\bf P}\ar[dl]_(.45){p_1} \ar[dr]^{p_2}&\\
\PP(H^1(\Theta_{C}))&&\PP(\HKC)
}
$$ this gives the stratification of the incidence stratum $\PP({\cal W}_{\Sigma^{\circ}_r})=p^{-1}_1 (\Sigma^{\circ}_r)$ according to the length and the partition attached to the associated graded modules of the filtrations;
the last piece of the $([\xi], [\phi])$-filtrations was shown to have a direct bearing on the geometry of the curve $C$;
 
- we are now on a refined stratum ${\mathfrak L}(h^l,\lambda)$ of $\PP({\cal W}_{\Sigma^{\circ}_r})$ , where the length of the filtrations is $l$ and assumed to be at least $3$, the last pieces $W^l_{\xi}([\phi])$ of the filtrations of constant dimension $h^l$ and the quotient spaces
$W_{\xi}/W^l_{\xi}([\phi])$ have the induced filtration governed by the partition $\lambda$; using the Hodge metric on $\HKC$ we turned the filtration on $W_{\xi}/W^l_{\xi}([\phi])$ into the direct sum decomposition
\begin{equation}\label{dirsum-pause}
W_{\xi}/W^l_{\xi}([\phi])=\bigoplus^{l-1}_{s=0} P^s([\xi],[\phi]);
\end{equation}
this comes with the dynamics expressed by maps between the summands of the decomposition; we attach the graph $PG_l$ which captures a part of the dynamics; the graph $PG_l$ is naturally bipartite - the vertices are bicolored `black' and `white'-  and oriented; thus at this stage we arrived to the category of quiver representations; those in turn have interesting properties; 

- our quiver $PG_l$ and its representations are `open': the extremal white/black vertices $(l-1)/(l-1)'$ and $(0)/(0)'$ are not {\it directly} connected and the extremal pieces of the direct sum decomposition, the summands $P^{l-1}([\xi],[\phi])$ and $ P^0([\xi],[\phi])$ labeling those extremal vertices, have no obvious `direct' maps between them;
with an eye toward having such maps we completed $PG_l$ to the trivalent graph $\widehat{PG}_l$; our quantum-type invariants give rise to an interesting structure on the spaces of maps between $P^{l-1}([\xi],[\phi])$ and $ P^0([\xi],[\phi])$: {\it the maps naturally extending the representations of the quiver $PG_l$ to the representations of the completed quiver $\widehat{PG}_l$}; those are the maps
$$
\tau^{\pm}_{\rho_c (\xi\otimes\phi)}: Hom(\CC[q,q^{-1}],\CC) \longrightarrow P^{\pm}([\xi],[\phi]),
$$
defined for every nonzero vector $\xi\otimes \phi$ lying over a point $([\xi],[\phi])$ in ${\mathfrak L}(h^l,\lambda)$ and every value of the trace parameter $c$; recall the notation:
$$
\begin{gathered}
P^{-}([\xi],[\phi])=Hom(P^0([\xi],[\phi]),P^{l-1}([\xi],[\phi])),
\\
P^{+}([\xi],[\phi])=Hom(P^{l-1}([\xi],[\phi]),P^{0}([\xi],[\phi]),
\end{gathered}
$$
are the spaces of maps between the `ends' $P^{l-1}([\xi],[\phi])$ and $P^{0}([\xi],[\phi])$ of the decomposition \eqref{dirsum-pause}; the construction of the maps is a sort of `integration' of the representation
$\rho_c (\xi\otimes\phi)$ of
$PG_l$ along the zig-zag paths $Z_{0,\pm}$ in $\widehat{PG}_l$;
 
- once we have an interesting structure on the space of maps between the `ends' of our decompositions, it is natural to ask if there are interesting structures on the space of {\it maps between the maps of the `ends'}, that is on the spaces
$$
Hom(P^{\mp}([\xi],[\phi]), P^{\pm}([\xi],[\phi]));
$$
 and of course, the higher categories teach us, that we should not stop there; in short, we are asking for higher categorical structures of our direct sum decompositions \eqref{dirsum-pause}. 
 Guided by this line of thinking we show that the spaces
$Hom(P^{\mp}([\xi],[\phi]), P^{\pm}([\xi],[\phi]))$ have interesting structures in their turn: in the next subsection those are related to Lagrangians in the symplectic vector space $P([\xi],[\phi])=P^-([\xi],[\phi]) \oplus P^+([\xi],[\phi])$.
  
\subsection{Quantum-type invariants and Lagrangians in $P([\xi],[\phi])$}  We have already seen that the space
$$
P([\xi],[\phi])=P^-([\xi],[\phi]) \oplus P^+([\xi],[\phi])
$$
is naturally a symplectic space via the identifications
$$
P([\xi],[\phi])\cong T^{\ast}_{P^-([\xi],[\phi])} \cong T^{\ast}_{P^+([\xi],[\phi])} 
$$
with the cotangent bundles of $P^{\pm}([\xi],[\phi])$. From the symplectic point of view the objects of interest are Lagrangians in $P([\xi],[\phi])$ with respect to either symplectic structure; recall that those are independent of the symplectic structure induced on $P([\xi],[\phi])$ by two identifications above, see Remark \ref{rem:sympP}. 

A source for Lagrangians are
sections of the structure projections
$$
\pi^{\pm}: T^{\ast}_{P^{\pm}([\xi],[\phi])} \longrightarrow {P^{\pm}([\xi],[\phi])}.
$$ 
Namely, we use a well known result saying that given a regular function $f$ on a manifold $M$, the differential $df$ defines a section
$$
df: M \longrightarrow T^{\ast}_{M}
$$
of the structure projection $\pi_M:T^{\ast}_{M} \longrightarrow M$; the image of $df$ is a Lagrangian in $T^{\ast}_{M}$ with respect to the canonical symplectic structure of the cotangent bundle, see \cite{Ch-Gi}.
In the case  of $P^{\pm}([\xi],[\phi])$ we have a natural choice of functions given by the coefficients of the characteristic polynomials of endomorphisms appearing in the definition of the maps
$$
\begin{gathered}
	C^-_{\rho_c(\xi\otimes\phi)} (q,\bullet): P^-=Hom(P^0,P^{l-1}) \longrightarrow {\mathfrak g}^{(0)}[q,q^{-1}] =\bigoplus^{l-1}_{s=0} End(P^s)[q,q^{-1}], 
	\\
		C^+_{\rho_c(\xi\otimes\phi)} (q,\bullet): P^+=Hom(P^{l-1},P^{0}) \longrightarrow {\mathfrak g}^{(0)}[q,q^{-1}] =\bigoplus^{l-1}_{s=0} End(P^s)[q,q^{-1}],
\end{gathered}
$$
see {\it Step 3} of the constructions following Proposition \ref{pro:mu-map}. More explicitly, for $X \in P^{\pm}$ we have
${\mathfrak g}^{(0)}$-valued Laurent polynomials
$$
C^{\pm}_{\rho_c(\xi\otimes\phi)} (q,X)=\sum^{l-1}_{s=0}(C^{\pm})^s (X)q^{(-1)^s s}. 
$$
Taking the characteristic polynomials term-wise gives
$$
ch^{\pm}_{\rho_c(\xi\otimes\phi)} (q,X,u):=\sum^{l-1}_{s=0}det((C^{\pm})^s (X) -u{\bf 1}_{P^s})q^{(-1)^s s}.
$$
Expanding with respect to the indeterminate $u$ we obtain
$$
ch^{\pm}_{\rho_c(\xi\otimes\phi)} (q,X,u)=\sum^{l-1}_{s=0}\sum^{h^s}_{m=0}(ch^s_m)^{\pm}_{\rho_c(\xi\otimes\phi)} (X) u^mq^{(-1)^s s}
$$
where $h^s=dim(P^s)$ and $(ch^s_m)^{\pm}_{\rho_c(\xi\otimes\phi)} (X)$ is the coefficient of $u^m$ in $det((C^{\pm})^s (X) -u1_{P^s})$:
\begin{equation}\label{charfunc}
(ch^s_m )^{\pm}_{\rho_c(\xi\otimes\phi)}(X)=(-1)^m tr \left(\bigwedge^{h^s-m}(C^{\pm}(X))^s\right).
\end{equation}
We view $(ch^s_m )^{\pm}_{\rho_c(\xi\otimes\phi)}$ for various values of $s$ and $m$ as polynomial functions on 
$P^{\pm}$. The differentials of those functions, according to the result cited above, provide a collection of Lagrangians in the symplectic vector space $P([\xi],[\phi])$. We summarize this in the following statement.
\begin{pro}\label{pro:Lagrangianschar}
	The vector spaces
	$$
	P([\xi],[\phi])=P^-([\xi],[\phi])\oplus P^+([\xi],[\phi])\cong T^{\ast}_{P^{\pm}([\xi],[\phi])},
	$$
	equipped with the natural symplectic structures, come along with a distinguished collection of Lagrangians denoted $(L^s_m)^{\pm}$. Those are defined by the characteristic coefficient functions
	$$
	(ch^s_m )^{\pm}_{\rho_c(\xi\otimes\phi)}:P^{\pm}([\xi],[\phi]) \longrightarrow \CC
	$$ 
	see \eqref{charfunc}; more precisely, $(L^s_m)^{\pm}$ is the image of the differential $d(ch^s_m )^{\pm}_{\rho_c(\xi\otimes\phi)}$ viewed as the section of the cotangent bundle
	$$
	d(ch^s_m )^{\pm}_{\rho_c(\xi\otimes\phi)}:P^{\pm}([\xi],[\phi]) \longrightarrow T^{\ast}_{P^{\pm}([\xi],[\phi])}.
	$$
\end{pro}

The maps $\tau^{\pm}_{\rho_c(\xi\otimes\phi)}$ in \eqref{Lpoly-Ppm} associated to the representations $\rho_c(\xi\otimes\phi)$ also provide sections of cotangent bundles $T^{\ast}_{P^{\pm}([\xi],[\phi])}$. Indeed, recall
that those maps are obtained by dualizing the maps
$$
\begin{gathered}
	C^-_{\rho_c(\xi\otimes\phi)} (q,1,\bullet): P^-=Hom(P^0,P^{l-1}) \longrightarrow \CC[q,q^{-1}] , 
	\\
	C^+_{\rho_c(\xi\otimes\phi)} (q,1,\bullet): P^+=Hom(P^{l-1},P^{0}) \longrightarrow \CC[q,q^{-1}] ,
\end{gathered}
$$
see \eqref{C+-maps-q-1}. Putting those maps together with the ones in \eqref{Lpoly-Ppm} we obtain the periodic chain of linear maps
\begin{equation}\label{chainofmaps}
\xymatrix@C=46pt{
	P^- \ar[r]^(.4){C^- (q,1,\bullet)}&\CC[q,q^{-1}] \ar[r]^(.65){\tau^+}& 
		P^+ \ar[r]^(.4){C^+ (q,1,\bullet)}&\CC[q,q^{-1}] \ar[r]^(.65){\tau^-}& P^-,
}
\end{equation}
where the reference to $[\xi],[\phi],\rho_c(\xi\otimes\phi)$ is omitted to simplify the notation.
Composing the  first two arrows from the left gives the linear map
$$
\sigma^{+,-}_{\rho_c(\xi\otimes\phi)}:=\tau^+ \circ C^- (q,1,\bullet): P^- \longrightarrow P^+.
$$
Similarly, composing the third and the fourth arrows in \eqref{chainofmaps} gives the linear map
$$
\sigma^{-,+}_{\rho_c(\xi\otimes\phi)}:=\tau^- \circ C^+ (q,1,\bullet): P^+ \longrightarrow P^-.
$$
The above maps define the linear maps
$$
\begin{gathered}
s(\sigma^{+,-}_{\rho_c(\xi\otimes\phi)}): P^- \longrightarrow P^- \oplus P^+ \cong T^{\ast}_{P^-}
\\
s(\sigma^{-,+}_{\rho_c(\xi\otimes\phi)}): P^+ \longrightarrow P^+ \oplus P^- \cong T^{\ast}_{P^+} 
\end{gathered}
$$
given by the formulas
$$
\begin{gathered}
s(\sigma^{+,-}_{\rho_c(\xi\otimes\phi)})(x)=x + \sigma^{+,-}_{\rho_c(\xi\otimes\phi)}(x)), \,\forall x\in P^-,
\\
s(\sigma^{-,+}_{\rho_c(\xi\otimes\phi)})(y)=y+\sigma^{-,+}_{\rho_c(\xi\otimes\phi)}(y)), \,\forall y\in P^+ .
\end{gathered}
$$
Hence the maps
\begin{equation}\label{sigma+-}
	\begin{gathered}
	\sigma^{+,-}_{\rho_c(\xi\otimes\phi)}: P^- \longrightarrow P^+,
	\\
	\sigma^{-,+}_{\rho_c(\xi\otimes\phi)}: P^+ \longrightarrow P^-
	\end{gathered}
\end{equation}
can be recast as sections of cotangent bundles. The images of those sections are the graphs of the above maps
\begin{equation}\label{graphs-sigma+-}
	\begin{gathered}
		\Gamma_{\sigma^{+,-}_{\rho_c(\xi\otimes\phi)}}=\{ x + \sigma^{+,-}_{\rho_c(\xi\otimes\phi)}(x) \in P^- \oplus P^+ | \,\forall x\in P^-\},
		\\
		\Gamma_{\sigma^{-,+}_{\rho_c(\xi\otimes\phi)}}=\{ y + \sigma^{+,-}_{\rho_c(\xi\otimes\phi)}(y) \in P^+ \oplus P^- | \,\forall y\in P^+\}.
	\end{gathered}
\end{equation}
We wish to understand how to attach to those graphs Lagrangians with respect to the canonical symplectic structure on $P^- \oplus P^+$ (resp, $P^+ \oplus P^-$). The following lemma will help to clarify the situation.
\begin{lem}\label{lem:graph-Lag}
	Let $V$ be a finite dimensional complex vector space and let $V^{\ast}$ be its dual. Consider the space
	$$
	E=V\oplus V^{\ast}
	$$
	equipped with the canonical symplectic form
	$$
	\omega_E (v+t^{\ast},w+s^{\ast}):= t^{\ast}(w)-s^{\ast}(v),
	$$
	for all $v,w \in V$ and $s^{\ast},t^{\ast} \in V^{\ast}$.
	Let $A: V\longrightarrow V^{\ast}$ be a linear map and let
	$$
	\Gamma_A=\{v+A(v) \in E=V\oplus V^{\ast}| v\in V\}
	$$
	be the graph of $A$.
	Then $\Gamma_A$ is a Lagrangian subspace of $E$
	if and only if $A$ is symmetric, that is,
	$$
	A(v)(w)=A(w)(v), \forall v,w \in V.
	$$
\end{lem}
\begin{pf}
	The image of $\Gamma_A$ has the dimension $dim(V)=\HA dim(E)$, so to see that it is a Lagrangian we need to check that $im(\Gamma_A)$ is an isotropic subspace with respect to $\omega_E$:
	$$
	\omega_E(v+A(v), w+A(w))=A(v)(w)-A(w)(v).
	$$
	The last expression vanishes precisely when $A$ is symmetric.
\end{pf}

\begin{rem}\label{rem:sym}
	Assume the set up of the lemma above. The space of $\CC$-linear maps from $V$ to its dual $V^{\ast}$ is isomorphic to the space $Bil_{\CC}(V)$ of bilinear pairings
	$
	V\times V \longrightarrow \CC: 
	$
	\begin{equation}\label{V-Vdual-Bil}
	Hom_{\CC}(V,V^{\ast}) \cong Bil_{\CC}(V)
\end{equation}
	where a homomorphism $A:V\longrightarrow V^{\ast}$goes to the pairing
	$b_A: 	V\times V \longrightarrow \CC$ defined by the formula
	$$
	b_A (v,w):=A(v)(w), \forall v,w\in V;
	$$
	its inverse takes a bilinear pairing $b:V\times V \longrightarrow \CC$ to the linear map
	$$
	A_b:V\longrightarrow V^{\ast}
	$$
	defined by the formula 
	$$
	A_b(v):=b(v,\bullet), \forall v\in V.
	$$
	Clearly, the isomorphism \eqref{V-Vdual-Bil} identifies the symmetric homomorphisms with the symmetric bilinear pairings.
	
	We have the symmetrization of a bilinear pairing: given $b\in Bil(V)$, its symmetrization $b^s$ is defined by the formula
	$$
	b^s (v,w):=\HA \left(b(v,w)+b(w,v) \right), \,\,\forall v,w \in V.
	$$
	Hence to any $A \in  Hom(V,V^{\ast})$ we associate its symmetrization
	$A^s$ defined by the formula
	$$
	A^s=A_{b^s_A},
	$$
	explicitly, we have
	$$
	A^s(v)(w)=\HA(A(v)(w)+A(w)(v)), \,\,\forall v,w\in V.
	$$
\end{rem}

The above remark allows us to associate a Lagrangian subspace to {\it any}
homomorphism $A\in Hom(V,V^{\ast})$:
\begin{equation}
	L_A:=\text{ \it the Lagrangian of $A$}\stackrel{def}{=} \Gamma_{A^s}=\text{ \it the graph of the symmetrization $A^s$ of $A$}.
\end{equation}

With the above discussion in mind we can attach Lagrangians to the maps
$\sigma^{+,-}_{\rho_c(\xi\otimes\phi)}$ and $\sigma^{-,+}_{\rho_c(\xi\otimes\phi)}$ in \eqref{sigma+-}.  
\begin{pro}\label{pro:Pxiphi-Lag}
	The maps 
	$$
	\begin{gathered}
	\sigma^{+,-}_{\rho_c(\xi\otimes\phi)}: P^-([\xi],[\phi])\longrightarrow P^+([\xi],[\phi]) \cong (P^-([\xi],[\phi]))^{\ast},
	\\
\sigma^{-,+}_{\rho_c(\xi\otimes\phi)}:P^+([\xi],[\phi]) \longrightarrow P^-([\xi],[\phi])\cong (P^+([\xi],[\phi]))^{\ast}.
\end{gathered}
$$
are symmetric: the following formulas hold
$$
\begin{gathered}
	\sigma^{+,-}_{\rho_c(\xi\otimes\phi)}(X)(X')=Res_0(C^-(q,1,X) C^-(q,1,X')),\,\,\forall X,X' \in P^-,
	\\
	\sigma^{-,+}_{\rho_c(\xi\otimes\phi)}(Y)(Y')=Res_0(C^+(q,1,Y) C^+(q,1,Y')),\,\,\forall Y,Y' \in P^+.
\end{gathered}
$$
In particular, the graphs of those maps are Lagrangian subspaces of the symplectic space $P([\xi],[\phi])$:
$$
\begin{gathered}
L^{+,-}_{\rho_c(\xi\otimes\phi)}:=L_{\sigma^{+,-}_{\rho_c(\xi\otimes\phi)}}=\text{the graph of $\sigma^{+,-}_{\rho_c(\xi\otimes\phi)}$},
\\
L^{-,+}_{\rho_c(\xi\otimes\phi)}:=L_{\sigma^{-,+}_{\rho_c(\xi\otimes\phi)}}=\text{the graph of $\sigma^{-,+}_{\rho_c(\xi\otimes\phi)}$}.
\end{gathered}
$$   
\end{pro}
\begin{pf}
	The two formulas are proved in a similar way. We give the details for the first one. For $X \in P^-$ we have
	$$
	\sigma^{+,-}_{\rho_c(\xi\otimes\phi)}(X)=\tau^+(C^-(q,1,X)).
	$$
	This is a unique element of $P^+$ acting on $P^-$ by the formula
	$$
	\begin{gathered}
	\sigma^{+,-}_{\rho_c(\xi\otimes\phi)}(X)(X')=\tau^+(C^-(q,1,X))(X')
	\\
	=tr(\tau^+(C^-(q,1,X)) \circ X') =Res_0(C^-(q,1,X)C^-(q,1,X')),
	\end{gathered}
	$$
	where the last equality is \eqref{formulaRes=tr} seen in the  proof of Proposition \ref{pro:coefFourier}.
\end{pf}
\begin{rem}\label{rem:sigma-maps-explicit}
	The map $\sigma^{+,-}_{\rho_c(\xi\otimes\phi)}$ (resp. $\sigma^{+,-}_{\rho_c(\xi\otimes\phi)}$) can be written explicitly in terms of the Fourier coefficients $\tau^{\pm}(q^k)$ in \eqref{Lpoly-Ppm-genseries} as follows:
	$$
	\begin{gathered}
	\sigma^{+,-}_{\rho_c(\xi\otimes\phi)}(X)=\tau^+ (C^-(q,1,X))=\sum^{l-1}_{s=0}tr \left((C^-)^s(X) \right) \tau^+(q^{(-1)^s s}), \,\,\forall X \in P^-,
	\\
		\sigma^{-,+}_{\rho_c(\xi\otimes\phi)}(Y)=\tau^- (C^-(q,1,X))=\sum^{l-1}_{s=0}tr \left((C^+)^s(Y) \right)\tau^-(q^{(-1)^s s}), \,\,\forall Y \in P^+.
		\end{gathered}
	$$
\end{rem}

 Observe that the chain 
 $$
 \xymatrix@C=46pt{
 	P^- \ar[r]^(.4){C^- (q,1,\bullet)}&\CC[q,q^{-1}] \ar[r]^(.65){\tau^+}& 
 	P^+ \ar[r]^(.4){C^+ (q,1,\bullet)}&\CC[q,q^{-1}] \ar[r]^(.65){\tau^-}& P^-,
 }
$$
can be extended indefinitely to the right (resp. left). Then we have two full `periods'
$$
\begin{gathered}
	T^-_{\rho_c(\xi\otimes\phi)}:=\tau^- C^+ (q,1,\bullet) \tau^+ C^- (q,1,\bullet):P^- \longrightarrow P^-,
	\\
	T^+_{\rho_c(\xi\otimes\phi)}:=\tau^+ C^- (q,1,\bullet) \tau^- C^+ (q,1,\bullet):P^+ \longrightarrow P^+.
\end{gathered}
$$
With the identifications
$$
(P^{\mp})^{\ast}\cong P^{\pm}
$$
the two endomorphisms are dual to each other
\begin{equation}\label{T+-dual}
	(T^{\mp}_{\rho_c(\xi\otimes\phi)})^{\ast}=T^{\pm}_{\rho_c(\xi\otimes\phi)}
\end{equation}
The endomorphisms $T^{\mp}_{\rho_c(\xi\otimes\phi)}$ act on the spaces $Hom(P^-,P^+)$ and  $Hom(P^+,P^-)$ in the following fashion:
\begin{equation}\label{T+-act}
	\begin{gathered}
 Hom(P^-,P^+) \ni \sigma \mapsto \sigma^{(1)}:=	T^+_{\rho_c(\xi\otimes\phi)} \sigma T^-_{\rho_c(\xi\otimes\phi)} \in Hom(P^-,P^+) ,
\\
 Hom(P^+,P^-) \ni \sigma' \mapsto (\sigma')^{(1)}:=	T^-_{\rho_c(\xi\otimes\phi)} \sigma' T^+_{\rho_c(\xi\otimes\phi)} \in Hom(P^+,P^-).
 \end{gathered}
\end{equation}
\begin{pro}
	The action in \eqref{T+-act} preserves the subspace $Hom^s(P^-,P^+)$ (resp., $Hom^s(P^+,P^-)$) of the symmetric maps in $Hom(P^-,P^+)$ (resp., $Hom(P^+,P^-)$). In particular, we have the action on the corresponding Lagrangians:
	$$
	\begin{gathered}
 L^{(1)}_{\sigma}  :=L_{ {\sigma}^{(1)}}, \forall \sigma \in Hom^s(P^-,P^+) 
	\\
	L^{(1)}_{\sigma'}  :=L_{ {\sigma'}^{(1)}}, \forall \sigma' \in Hom^s(P^+,P^-).
	\end{gathered} 
	$$
\end{pro}
\begin{pf}
	Dualize $T^+_{\rho_c(\xi\otimes\phi)} \sigma T^-_{\rho_c(\xi\otimes\phi)}$ to obtain the formula
	$$
	\left(T^+_{\rho_c(\xi\otimes\phi)} \sigma T^-_{\rho_c(\xi\otimes\phi)}\right)^{\ast}=
	\left(T^-_{\rho_c(\xi\otimes\phi)}\right)^{\ast} \sigma^{\ast} 	\left(T^+_{\rho_c(\xi\otimes\phi)} \right)^{\ast}=T^+_{\rho_c(\xi\otimes\phi)} \sigma^{\ast} T^-_{\rho_c(\xi\otimes\phi)},
	$$
	where the second equality uses the duality relation \eqref{T+-dual}. This implies that  
	$T^+_{\rho_c(\xi\otimes\phi)} \sigma T^-_{\rho_c(\xi\otimes\phi)}$ is symmetric, if $\sigma$ is.
\end{pf}
	More generally, iterating the action in \eqref{T+-act} we have
$$
\begin{gathered}
\sigma^{(n)}:=	(T^+)^n_{\rho_c(\xi\otimes\phi)} \sigma (T^-)^n_{\rho_c(\xi\otimes\phi)},\,_,\forall \sigma \in Hom^s(P^-,P^+),
\\
{\sigma'}^{(n)}:=	(T^-)^n_{\rho_c(\xi\otimes\phi)} \sigma' (T^+)^n_{\rho_c(\xi\otimes\phi)}, \,\,\forall \sigma'\in Hom^s(P^+,P^-)
\end{gathered}
$$
for all $n\in \ZZ_{\geq 0}$. This gives rise to the sequence of Lagrangian subspaces
\begin{equation}
	\begin{gathered}
		L^{(n)}_{\sigma}  :=L_{ {\sigma}^{(n)}}, \forall \sigma \in Hom^s(P^-,P^+) 
		\\
		L^{(n)}_{\sigma'}  :=L_{ {\sigma'}^{(n)}}, \forall \sigma' \in Hom^s(P^+,P^-).
	\end{gathered} 
\end{equation}

The action in \eqref{T+-act} naturally extends to give a structure of $\CC[t]$-bimodule on the vector space $ Hom(P^-,P^+) $ (resp. $ Hom(P^+,P^-)$ ).
\begin{pro}\label{pro:bimodule}
	The vector spaces 
	$$ 
	\text{$Hom(P^-,P^+) $ and  $ Hom(P^+,P^-)$}
	$$
	 have a natural structure of $\CC[t]$-bimodules. This is defined by the following formulas: 
\begin{equation}\label{tntm}
	\begin{gathered}
	t^n \sigma t^m := 	
		(T^+_{\rho_c(\xi\otimes\phi)})^n	\sigma  (T^-)^m_{\rho_c(\xi\otimes\phi)} \,\, \forall \sigma \in  Hom(P^-,P^+),
		\\
		t^n \sigma' t^m := 	(T^-_{\rho_c(\xi\otimes\phi)})^n	\sigma' (T^+)^m_{\rho_c(\xi\otimes\phi)} ,\,\,\forall \sigma'\in Hom(P^+,P^-),
	\end{gathered}
\end{equation}
for every $n,m\in \ZZ_{\geq 0}$. In particular, for every polynomial
$p\in \CC[t]$, $\sigma \in  Hom^s(P^-,P^+)$ and $\sigma' \in  Hom^s(P^+,P^-)$ , we have
$$
\begin{gathered}
p \cdot\sigma:= p(t) \sigma p(t) \in Hom^s(P^-,P^+),
\\
p \cdot\sigma':= p(t) \sigma' p(t) \in Hom^s(P^+,P^-).
\end{gathered}
$$
\end{pro}

The above tells us that every polynomial $p\in \CC[t]$ gives rise to the Lagrangian subspaces
$$
\begin{gathered}
L^{p}_{\sigma}:=L_{p \cdot\sigma},\,\,\forall \sigma\in Hom^s(P^-,P^+),
\\
L^{p}_{\sigma}:=L_{p \cdot\sigma'},\,\,\forall \sigma'\in Hom^s(P^+,P^-).
\end{gathered}
$$
  Specializing $\sigma=\sigma^{+,-}_{\rho_c(\xi\otimes\phi)}$ and $\sigma'=\sigma^{-,+}_{\rho_c(\xi\otimes\phi)}$ we obtain 
 the maps from the space of  polynomials $\CC[t]$ to the variety ${\mathfrak Lag}(P([\xi],[\phi]))$ of Lagrangian subspaces of $P([\xi],[\phi]))$:
 \begin{equation}\label{p-Lag}
 	\begin{gathered}
 \CC[t]\ni p\mapsto L^p_{\sigma^{+,-}_{\rho_c(\xi\otimes\phi)}} \in {\mathfrak Lag}(P([\xi],[\phi])),
 \\
 \CC[t]\ni p\mapsto L^p_{\sigma^{-,+}_{\rho_c(\xi\otimes\phi)}} \in {\mathfrak Lag}(P([\xi],[\phi])).
\end{gathered}
\end{equation}

In a similar fashion we can define the structure of $\CC[t,t^{-1}]$-bimodule on  the spaces $Hom(P^{-},P^+)$ and $Hom(P^{+},P^-)$. For this we exponentiate the endomorphisms $T^{\pm}_{\rho_c(\xi\otimes\phi)}$
\begin{equation}\label{exp-T+-}
	{\mathfrak e}(T^{\pm}_{\rho_c(\xi\otimes\phi)}):=exp(T^{\pm}_{\rho_c(\xi\otimes\phi)})=\sum^{\infty}_{n=0} \frac{\left(T^{\pm}_{\rho_c(\xi\otimes\phi)}\right)^n}{n!}
\end{equation} 
 to obtain invertible endomorphisms of $P^{\pm}$; replacing $T^{\pm}$ by their exponentials the formulas in \eqref{tntm} can be now extended to all integers:
 \begin{equation}\label{tntm-exp}
 	\begin{gathered}
 		t^n \,{\hat{\cdot}}\,\sigma \,\hat{\cdot}\, t^m := 	
 		({\mathfrak e}(T^+_{\rho_c(\xi\otimes\phi)}))^n	\sigma  ({\mathfrak e}(T^-_{\rho_c(\xi\otimes\phi)}))^m \,\, \forall \sigma \in  Hom(P^-,P^+),
 		\\
 		t^n \,{\hat{\cdot}}\,\sigma' \,\hat{\cdot}\,  t^m := 	({\mathfrak e}(T^+_{\rho_c(\xi\otimes\phi)}))^n	\sigma' ({\mathfrak e}(T^+_{\rho_c(\xi\otimes\phi)}))^m ,\,\,\forall \sigma'\in Hom(P^+,P^-),
 	\end{gathered}
 \end{equation}
for all $m,n \in \ZZ$. 
\begin{pro}\label{pro:Lp-act}
	The formulas \eqref{tntm-exp} define a structure of $\CC[t,t^{-1}]$-bimodule on the vector spaces $Hom(P^{\mp},P^{\pm})$.
	In particular, for every Laurent polynomial $p \in \CC[t,t^{-1}]$ and every symmetric map $\sigma \in  Hom^s(P^-,P^+)$ and $\sigma' \in Hom^s(P^+,P^-)$ we have
	$$
	\begin{gathered}
		p\,\hat{\cdot} \,\sigma :=p(t) \,\hat{\cdot} \sigma \,\hat{\cdot} p(t)\in Hom^s(P^-,P^+),\\
		p\,\hat{\cdot} \,\sigma' :=p(t) \,\hat{\cdot} \sigma' \,\hat{\cdot} p(t)\in Hom^s(P^+,P^-).
	\end{gathered}
$$
\end{pro}
We now extend the correspondence with Lagrangians to the space of the Laurent polynomials:
\begin{equation}\label{Lpol-Lag}
	\begin{gathered}
		\CC[t,t^{-1}] \times Hom^s(P^-,P^+) \ni (p,\sigma) \mapsto L^{\hat{\cdot} p}_{\sigma}:= L_{p\,\hat{\cdot} \,\sigma} \in {\mathfrak Lag} (P([\xi],[\phi])),
		\\
		\CC[t,t^{-1}] \times Hom^s(P^+,P^-) \ni (p,\sigma') \mapsto L^{\hat{\cdot} p}_{\sigma'}:= L_{p\,\hat{\cdot} \,\sigma'} \in {\mathfrak Lag} (P([\xi],[\phi])).
	\end{gathered}
\end{equation}

 Specializing $\sigma=\sigma^{+,-}_{\rho_c(\xi\otimes\phi)}$ and $\sigma'=\sigma^{-,+}_{\rho_c(\xi\otimes\phi)}$ we obtain 
the maps from the space of Laurent polynomials $\CC[t,t^{-1}]$ to the space of Lagrangians ${\mathfrak Lag}(P([\xi],[\phi]))$:
\begin{equation}\label{phat-Lag}
	\begin{gathered}
		\CC[t,t^{-1}]\ni p\mapsto L^{\hat{\cdot}p}_{\sigma^{+,-}_{\rho_c(\xi\otimes\phi)}} \in {\mathfrak Lag}(P([\xi],[\phi])),
		\\
		\CC[t,t^{-1}]\ni p\mapsto L^{\hat{\cdot}p}_{\sigma^{-,+}_{\rho_c(\xi\otimes\phi)}} \in {\mathfrak Lag}(P([\xi],[\phi])).
	\end{gathered}
\end{equation} 
 
 The algebra of Laurent polynomials in one variable now have several different interpretations:
 
 $\bullet$ moduli to connect the `ends' $P^{l-1}([\xi],[\phi])$ and $P^0([\xi],[\phi])$ of the $([\xi],[\phi])$-decomposition of $W_{\xi}/W^l_{\xi}([\phi])$; those are the homomorphisms
 $$
 \tau^{\pm}_{\rho_c(\xi\otimes\phi)} :\CC[q,q^{-1}] \longrightarrow P^{\pm}([\xi],[\phi]);
 $$
 
 $\bullet$ the `orbits' of Lagrangian subspaces $\{L^{\hat{\cdot} p}_{\sigma}\}_{p \in \CC[t,t^{-1}]}$ (resp., $\{L^{\hat{\cdot} p}_{\sigma'}\}_{p \in \CC[t,t^{-1}]}$ ), for $\sigma \in Hom^s(P^-,P^+)$ (resp, $\sigma' \in Hom^s(P^+,P^-)$), see \eqref{Lpol-Lag}.
 
 \vspace{0.2cm}
 The second interpretation of course emerges from the first and could be viewed as a manifestation of the Fukaya-type categories on the spaces of
 homomorphisms between the `ends'. We have already mentioned that the homomorphisms $\tau^{\pm}_{\rho_c(\xi\otimes\phi)}$ are perhaps a manifestation of a TQFT, in a sense that they could be viewed as linear versions of cobordism between two circles on the topological torus $\mathbb{T}$. We will see shortly that the maps $\tau^{\pm}_{\rho_c(\xi\otimes\phi)}$ can be used to give {\it conformal} structure on $\mathbb{T}$, that is, the space of Laurent polynomials $\CC[q,q^{-1}]$, via the maps $\tau^{\pm}_{\rho_c(\xi\otimes\phi)}$, will become moduli for curves of genus one with $l$ marked points.

 \subsection{The zig-zag path $Z_{-,+}$}
 	Until now we were considering the zig-zag paths $Z_{0,\pm}$. Somewhat similar constructions can be performed with the path $Z_{-,+}$. This will provide additional structures on the space of maps
 	$Hom(P^{\mp}, P^{\pm})$. However, we will see momentarily that there are some differences. Throughout the subsection we consider defining representations $\rho_c(\xi\otimes \phi)$ of the quiver $PG_l$, for $\xi\otimes \phi$ a nonzero element of the tautological fibration
 	$\OO^{\times}_{\mathfrak{L}_l}(-1) \longrightarrow \mathfrak{L}_l$; the reference to $\xi$ and $\phi$ will be often omitted to simplify the notation.
 	
 	To begin with, recall, see Lemma \ref{lem:zzpaths}, that we need to distinguish between odd and even values of $l$: in the first case $Z_{-,+}$ is connected and in the second it falls into two disconnected components denoted
 	$Z^0_{-,+}$ and $Z^1_{-,+}$, where the first (resp. second) component is composed of edges $e^{\pm}_i$ with the index $i$ even (resp. odd).
 	For this reason the construction differs according to the parity of $l$.
 	
 	The other, more serious difference with the previous constructions, is
 	{\it composing} the maps $\{\alpha^{r,s}_c(\xi,\phi)\}$ labeling the edges of $Z_{-,+}$. Recall from Lemma \ref{lem:zzpaths}:
 	$$
 	Z_{-,+} =(-e^+_1)e^-_3 \cdots (-e^{+}_{l-4})e^-_{l-2} e^-_0 \cdots (-e^+_{l-5})e^-_{l-3}(-e^+_{l-3})e^-_{l-1}(-e^+_{l-1})e^-_1.
 	$$
 	The orientation of edges alternate: in the above formula, the path starts at the white vertex $(1)$ goes to the black vertex $(0)'$, then from $(0)'$ back to the white vertex $(l-1)$ and so on; that is, the path goes from `white' to `black' to `white' in an alternating pattern. However, the maps $\alpha^{s\pm 1,s}_c$ which are supposed to label the edges are always directed `white' to `black', so if we assign the labels
 	$$
 	e^-_{l-1} \rightarrow \alpha^{l-2,l-1}: P^{l-1} \longrightarrow P^{l-2}, \,\,\, -e^+_{l-3} \rightarrow \pm\alpha^{l-2,l-3}: P^{l-3} \longrightarrow P^{l-2},
 	$$
 	the concatenation $(-e^+_{l-3}) e^-_{l-1}$ must be the composition
 	$\pm\alpha^{l-2,l-3}\alpha^{l-2,l-1}$ and the latter makes no sense.
 	Thus the change of orientation $(-e^+_s)$ in the formula of $Z_{-,+}$ should also  mean that the morphism labeling it must change the direction. Since all spaces $P^s$ of the representation $\rho_c(\xi\otimes\phi)$ are Hermitian we define
 	\begin{equation}\label{adj}
 		{}^{\dag}\alpha^{s,s-1}: P^{s} \longrightarrow P^{s-1}
 	\end{equation}
 	as the adjoint of the homomorphism 
 	$$
 	\alpha^{s,s-1}: P^{s-1} \longrightarrow P^s.
 	$$
 	Explicitly, ${}^{\dag}\alpha^{s,s-1}$ is a unique $\CC$-linear map
 	from $P^s$ to $P^{s-1}$ satisfying the equation 
 	\begin{equation}\label{adj-formula}
 		(\alpha^{s,s-1}(x),y)=(x,{}^{\dag}\alpha^{s,s-1}(y) ),\,\,\forall x\in P^{s-1}, \forall y\in P^s,
 	\end{equation}
 where $(\bullet,\bullet)$ stands for the Hodge metric on $\HKC$. With this definition in mind we label the edges composing $Z_{-,+}$ as follows:
 \begin{equation}\label{Z-+-labels}
 \begin{gathered}
 e^-_s \rightarrow \alpha^{s-1,s}, \,\,\forall s\in [1,l-1],
 \\
 -e^+_t \rightarrow {}^{\dag}\alpha^{t+1,t},  \,\,\forall t\in [0,l-2].
\end{gathered}
 \end{equation}
 Given this labeling the concatenation $(-e^+_{s-2})e^-_s$ of two successive
 edges of $Z_{-,+}$ becomes the composition
 $$
 \xymatrix@C=40pt{
 	P^s \ar[r]^{\alpha^{s-1,s}}& P^{s-1} \ar[r]^{{}^{\dag}\alpha^{s-1,s-2}}& P^{s-2}
}
$$
and this is now well defined as long as $s \in [2,l-1]$ - remember we always work under assumption that the length $l\geq 3$. From now on the constructions will be similar to the ones involving $Z_{0,\pm}$.

 	\vspace{0.2cm}
 	\noindent
 	{\bf Even case.} The two unlabeled edges are $e^-_0$ and $e^+_{l-1}$ and they lie in different connected components of $Z_{-,+}$, namely, in $Z^0_{-,+}$ and $Z^1_{-,+}$ respectively. So given the representation $\rho_c=\{\alpha^{r,s}_c\}$ of $PG_l$ and labeling $e^-_0$ by an element $X$ in $P^-([\xi],[\phi])$ all arrows of $Z^0_{-,+}$
 	$$
 	Z^0_{-,+}=(-e^+_0)e^-_2(-e^+_2)\cdots(-e^+_{l-4})e^-_{l-2}(-e^+_{l-2})e^-_0
 	$$
 	 are labeled by linear maps according to the rule \eqref{Z-+-labels}.  Composing them gives the operator
 		$$
 		\begin{gathered}
 	\rho^0_{Z^0_{-,+}} (X,\rho_c):  P^0([\xi],[\phi]) \longrightarrow   P^0([\xi],[\phi]):
 	\\
 	\rho^0_{Z^0_{-,+}} (X,\rho_c)= {}^{\dag}\alpha^{1,0}\alpha^{1,2}\cdots \alpha^{l-5,l-4}\,{}^{\dag}\alpha^{l-3,l-4}\, \alpha^{l-3,l-2}\,{}^{\dag}\alpha^{l-1,l-2}X.
 \end{gathered}
 	$$
 Shifting cyclically from left to right by a single label we obtain the map
 $$
 \begin{gathered}
 	\rho^1_{Z^0_{-,+}} (X,\rho_c):  P^1([\xi],[\phi]) \longrightarrow   P^1([\xi],[\phi]):
 	\\
 	\rho^1_{Z^0_{-,+}} (X,\{\alpha^{r,s}_c\})= \alpha^{1,2}\cdots \alpha^{l-5,l-4}\,{}^{\dag}\alpha^{l-3,l-4}\, \alpha^{l-3,l-2}\,{}^{\dag}\alpha^{l-1,l-2}X{}^{\dag}\alpha^{1,0}
 \end{gathered}
 $$	
 Continuing in this fashion gives the grading preserving operator
 \begin{equation}\label{rho-Z-+}
 	\rho_{Z^0_{-,+}} (X,\rho_c):=\bigoplus \rho^s_{Z^0_{-,+}} (X,\rho_c): \bigoplus^{l-1}_{s=0} P^s \longrightarrow \bigoplus^{l-1}_{s=0} P^s
 \end{equation}
representing the path $Z^0_{-,+}$; hence we have
$$
\rho_{Z^0_{-,+}} (X,\rho_c) \in {\mathfrak g}^{(0)}=\bigoplus^{l-1}_{s=0} End(P^s).
$$ 
This gives the homomorphism
\begin{equation}\label{C-Z0-+}
C_{Z^0_{-,+},\rho_c(\xi\otimes\phi)}: P^-=P^- ([\xi],[\phi]) \longrightarrow  {\mathfrak g}^{(0)}_{[\xi],[\phi]}
\end{equation}
defined by the formula
$$
C_{Z^0_{-,+},\rho(\xi\otimes\phi)} (X)=\rho_{Z^0_{-,+}} (X,\rho(\xi\otimes\phi));
$$
to simplify the notation we often omit the reference to $([\xi],[\phi])$ and $\xi\otimes\phi$; observe also that there is no reference to the trace parameter $c$: this is because the maps $\alpha^{s\pm1,s}$ are independent of $c$.

Analogous construction can be performed with the path $Z^1_{-,+}$:
$$
Z^1_{-,+}=e^-_1 (-e^+_1)\cdots(-e^+_{l-3})e^-_{l-1}(-e^+_{l-1}).
$$
As before the change of orientation of $e^+_s$ means that we must label that edge with the adjoint operator. Hence labeling $e^+_{l-1}$ by a map $Y\in P^+$, the edge $(-e^+_{l-1})$ receives the label ${}^{\dag}Y$ and the operator representing $Z^1_{-,+}$ based at the vertex $(0)'$ has the following form
$$
		\rho^0_{Z^1_{-,+}} (Y,\rho_c)= \alpha^{0,1}({}^{\dag}\alpha^{2,1})\cdots ({}^{\dag}\alpha^{l-4,l-5})\,\alpha^{l-4,l-3}\, ({}^{\dag}\alpha^{l-2,l-3})\,\alpha^{l-2,l-1}({}^{\dag}Y).
$$
Shifting cyclically from left to right by a single label we obtain the map in degree one
$$
\begin{gathered}
\rho^1_{Z^1_{-,+}} (Y,\rho_c):  P^1([\xi],[\phi]) \longrightarrow   P^1([\xi],[\phi]):
\\
\rho^1_{Z^1_{-,+}} (Y,\rho_c)= ({}^{\dag}\alpha^{2,1})\cdots ({}^{\dag}\alpha^{l-4,l-5})\,\alpha^{l-4,l-3}\, ({}^{\dag}\alpha^{l-2,l-3})\,\alpha^{l-2,l-1}({}^{\dag}Y)\alpha^{0,1}.
\end{gathered}
$$	
Continuing in this fashion gives the grading preserving operator
\begin{equation}
\rho_{Z^1_{-,+}} (Y,\rho_c):=\bigoplus \rho^s_{Z^1_{-,+}} (Y,\rho_c): \bigoplus^{l-1}_{s=0} P^s \longrightarrow \bigoplus^{l-1}_{s=0} P^s
\end{equation}
representing the path $Z^1_{-,+}$; hence we have
$$
\rho_{Z^1_{-,+}} (Y,\rho_c) \in {\mathfrak g}^{(0)}=\bigoplus^{l-1}_{r=0} End(P^r).
$$ 
This gives the map
$$
 P^+ ({[\xi],[\phi]}) \longrightarrow  {\mathfrak g}^{(0)}_{[\xi],[\phi]}
$$
assigning to $Y\in P^+({[\xi],[\phi]})$ the grading preserving operator
$\rho_{Z^1_{-,+}} (Y,\rho_c)$ which again is independent of the trace parameter $c$. However, this is $\CC$-{\it anti-linear}, since the formulas for each component of $\rho_{Z^1_{-,+}} (Y,\rho_c)$ involve the adjoint ${}^{\dag}Y$ of $Y$. Thus the assignment
\begin{equation}
	 P^+ ({[\xi],[\phi]}) \ni Y \mapsto {}^{\dag} \left(\rho_{Z^1_{-,+}} (Y,\rho_c) \right)\in {\mathfrak g}^{(0)}_{[\xi],[\phi]}
\end{equation}
is $\CC$-linear; we denote the resulting homomorphism by $C_{Z^1_{-,+},\rho(\xi\otimes\phi)}$:

\begin{equation}\label{C-Z1-+}
C_{Z^1_{-,+},\rho(\xi\otimes\phi)}: P^+_{[\xi],[\phi]} \longrightarrow  {\mathfrak g}^{(0)}_{[\xi],[\phi]}.
\end{equation}

 	 Extend the two homomorphisms  $C_{Z^i_{-,+},\rho}$, $i=0,1$, to the loop algebra of ${\mathfrak g}^{(0)}_{[\xi],[\phi]}$ to obtain the homomorphisms
 	\begin{equation}\label{CZ0-1}
 	\begin{gathered}
 	C_{Z^0_{-,+},\rho}(q,\bullet): P^- \longrightarrow {\mathfrak{g}}^{(0)} [q,q^{-1}],
 	\\
 		C_{Z^1_{-,+},\rho}(q,\bullet): P^+ \longrightarrow {\mathfrak{g}}^{(0)} [q,q^{-1}].
 		\end{gathered}
 	\end{equation}
 	Composing with the graded trace map
 	$$
 	{\mathfrak{g}}^{(0)} [q,q^{-1}] \longrightarrow \CC[q,q^{-1}]
 	$$
 	gives the homomorphisms
 	\begin{equation}\label{CZ-+-LP}
 		\begin{gathered}
 			C^-_{Z_{-,+},\rho}(q,\bullet): P^- \longrightarrow \CC [q,q^{-1}],
 			\\
 			C^+_{Z_{-,+},\rho}(q,\bullet): P^+ \longrightarrow \CC [q,q^{-1}].
 		\end{gathered}
 	\end{equation}
 Dualizing we obtain the analogues of the maps $\tau^{\pm}_{\rho_c}$ attached to the paths $Z_{0,\pm}$:
 \begin{equation}\label{tauZ-+-LP}
 	\begin{gathered}
 		\tau^+_{Z_{-,+},\rho}:=(C^-_{Z_{-,+},\rho}(q,\bullet))^{\ast}: Hom(\CC [q,q^{-1}],\CC) \longrightarrow (P^-)^{\ast} \cong P^+,
 		\\
 		\tau^-_{Z_{-,+},\rho}:=(C^+_{Z_{-,+},\rho}(q,\bullet))^{\ast}:Hom(\CC [q,q^{-1}],\CC) \longrightarrow (P^+)^{\ast} \cong P^- .
 	\end{gathered}
 \end{equation}

The operators 	$C_{Z^0_{-,+},\rho}(q,\bullet)$ and 
$C_{Z^1_{-,+},\rho}(q,\bullet)$ in \eqref{CZ0-1} can be also `fused' together to give the bilinear pairing
$$
B^{0,1}_{\rho} (q,\bullet,\bullet):P^- \times P^+ \longrightarrow {\mathfrak{g}}^{(0)}[q,q^{-1}]
$$
defined by composing grade-wise the values of $C_{Z^0_{-,+},\rho}$ and 
$C_{Z^1_{-,+},\rho}$:
\begin{equation}\label{B01}
	B^{0,1}_{\rho} (q,X,Y):=C_{Z^1_{-,+},\rho}(q,Y)\circ C_{Z^0_{-,+},\rho}(q,X), \,\,\forall X\in P^-, \forall Y\in P^+.
\end{equation}	
Explicitly, this is given by the formula
\begin{equation}\label{B01-formula}
	\begin{gathered}
	B^{0,1}_{\rho} (q,X,Y)=C_{Z^1_{-,+},\rho}(q,Y)\circ C_{Z^0_{-,+},\rho}(q,X)
	\\
	=\left(\sum^{l-1}_{s=0} {}^{\dag}\rho^s_{Z^1_{-,+}}(Y,\rho)q^{(-1)^s s} \right)\left(\sum^{l-1}_{s=0} \rho^s_{Z^0_{-,+}}(X,\rho)q^{(-1)^s s} \right)
	\\
	=\sum^{l-1}_{s=0} \left({}^{\dag}\rho^s_{Z^1_{-,+}}(Y,\rho)\circ \rho^s_{Z^0_{-,+}}(X,\rho)\right) q^{(-1)^s s}.
	\end{gathered}
\end{equation} 
 The superscripts in the notation $	B^{0,1}_{\rho}$ indicate the order of the composition:
 from right to left, the first factors are the components of the operator attached to $Z^0_{-,+}$ and the second are the ones of the operator attached to $Z^1_{-,+}$. Inverting the order gives another bilinear paring
 \begin{equation}\label{B10}
 	\begin{gathered}
 	B^{1,0}_{\rho} (q,\bullet,\bullet):P^+ \times P^- \longrightarrow {\mathfrak{g}}^{(0)}[q,q^{-1}],
 	\\
 	B^{1,0}_{\rho} (q,Y,X):=C_{Z^0_{-,+},\rho}(q,X)\circ C_{Z^1_{-,+},\rho}(q,Y), \,\,\forall X\in P^-, \forall Y\in P^+.
 	\end{gathered}
 \end{equation}	  
Equivalently, the two bilinear pairings above give homomorphisms
\begin{equation}\label{B01-B10-tensor}
	\begin{gathered}
		B^{1,0}_{\rho} (q,\bullet):P^+ \otimes P^- \longrightarrow {\mathfrak{g}}^{(0)}[q,q^{-1}],
		\\
		B^{0,1}_{\rho} (q,\bullet):P^- \otimes P^+ \longrightarrow {\mathfrak{g}}^{(0)}[q,q^{-1}],
	\end{gathered}
\end{equation}
where by abuse of notation we denote the linear maps in the same way as the corresponding bilinear pairings.
Composing on the right with the graded trace map we obtain the linear maps
\begin{equation}\label{B01-B10-tensor-trace}
	\begin{gathered}
		tr(B^{1,0}_{\rho}):P^+ \otimes P^- \longrightarrow \CC[q,q^{-1}],
		\\
		tr(B^{0,1}_{\rho} ):P^- \otimes P^+ \longrightarrow \CC[q,q^{-1}],
	\end{gathered}
\end{equation}
Dualizing gives 
\begin{equation}\label{B01-B10-tensor-trace-dual}
	\begin{gathered}
		tr^{\ast}(B^{1,0}_{\rho}):Hom(\CC[q,q^{-1}],\CC) \longrightarrow (P^+ \otimes P^-)^{\ast} 
		\\
		\cong (P^-)^{\ast}\otimes (P^+)^{\ast}  \cong P^+\otimes (P^+)^{\ast} \cong End (P^+),
		\\
		\\
		tr^{\ast}(B^{0,1}_{\rho} ):Hom(\CC[q,q^{-1}],\CC) \longrightarrow (P^- \otimes P^+)^{\ast} 
		\\
		\cong  (P^+)^{\ast}\otimes (P^-)^{\ast}  \cong P^-\otimes (P^-)^{\ast} \cong End (P^-),
	\end{gathered}
\end{equation} 
where we use the identification $(P^-)^{\ast} \cong P^+$ (resp., $(P^+)^{\ast} \cong P^-$). We summarize the above construction in the following statement.

\begin{pro}\label{pro:Z-+leven}
	Let $\rho_c(\xi\otimes\phi)=\{\alpha^{t,s}_c(\xi,\phi): P^s \longrightarrow P^{t}\}$ be the representation of the quiver $PG_l$ associated to a nonzero element $\xi\otimes\phi$ of the total space
	of the tautological line bundle $\OO_{{\mathfrak{L}}_l}(-1)$. For $l=2k$ even, the zig-zag path
	$Z_{-,+}$ of the graph $\widehat{PG}_l$ falls into two connected components - the closed paths $Z^0_{-,+}$ and $Z^1_{-,+}$: the first (resp. second) passes through even (resp. odd) white and odd (resp. even) black vertices. These two paths give rise to two distinguished linear maps
	$$
		\begin{gathered}
		\tau^+_{Z_{-,+},\rho}:=(C^-_{Z_{-,+},\rho}(q,\bullet))^{\ast}: Hom(\CC [q,q^{-1}],\CC) \longrightarrow (P^-)^{\ast} \cong P^+,
		\\
		\tau^-_{Z_{-,+},\rho}:=(C^+_{Z_{-,+},\rho}(q,\bullet))^{\ast}:Hom(\CC [q,q^{-1}],\CC) \longrightarrow (P^+)^{\ast} \cong P^- ,
	\end{gathered}
$$
 defined in \eqref{tauZ-+-LP}. In addition, fusing the two components together gives the distinguished linear maps
 $$
 	\begin{gathered}
 	tr^{\ast}(B^{1,0}_{\rho}):Hom(\CC[q,q^{-1}],\CC) \longrightarrow   End (P^+),
 	\\
 	tr^{\ast}(B^{0,1}_{\rho} ):Hom(\CC[q,q^{-1}],\CC) \longrightarrow  End (P^-),
 \end{gathered}
$$
see \eqref{B01-B10-tensor-trace-dual}. In the above formulas $\rho=\rho_c$ and the maps above are independent of the trace parameter $c$.
\end{pro}

There is another way to `fuse' the two components $Z^0_{-,+}$ and $Z^1_{-,+}$. This is based on the fact that the white vertices of $\widehat{PG}_l$ visited by $Z^0_{-,+}$ (resp. $Z^1_{-,+}$) are even (resp. odd). This means that the components  of $\rho_{Z^0_{-,+}} (X,\rho_c)$ in \eqref{rho-Z-+} with even (resp. odd) grade are operators corresponding to the path $Z^0_{-,+}$ cyclically shifted from left to right by even number of labels.  We split $\rho_{Z^0_{-,+}} (X,\rho_c)$ according the parity of the components
$$
\rho_{Z^0_{-,+}} (X,\rho_c)=\rho^{ev}_{Z^0_{-,+}} (X,\rho_c) \oplus \rho^{odd}_{Z^0_{-,+}} (X,\rho_c),
$$
where 
$$
\begin{gathered}
\rho^{ev}_{Z^0_{-,+}} (X,\rho_c) \in {\mathfrak g}^{(0)}_{ev} =\bigoplus_{\text{$s$ even}} End(P^s),
\\
\rho^{odd}_{Z^0_{-,+}} (X,\rho_c) \in {\mathfrak g}^{(0)}_{odd} =\bigoplus_{\text{$s$ odd}} End(P^s)
\end{gathered}
$$
This gives the linear maps
$$
\begin{gathered}
	C^{ev}_{Z^0_{-,+},\rho}: P^- \longrightarrow {\mathfrak g}^{(0)}_{ev},
	\\
C^{odd}_{Z^0_{-,+},\rho}: P^- \longrightarrow {\mathfrak g}^{(0)}_{odd},
\end{gathered}
$$
defined respectively by the formulas
$$
\begin{gathered}
C^{ev}_{Z^0_{-,+},\rho}(X)=	\rho^{ev}_{Z^0_{-,+}} (X,\rho_c) \in {\mathfrak g}^{(0)}_{ev} =\bigoplus_{\text{$s$ even}} End(P^s),
	\\
C^{odd}_{Z^0_{-,+},\rho}(X)=	\rho^{odd}_{Z^0_{-,+}} (X,\rho_c) \in {\mathfrak g}^{(0)}_{odd} =\bigoplus_{\text{$s$ odd}} End(P^s).
\end{gathered}
$$
Extending those maps to the loop algebra ${\mathfrak g}^{(0)}_{\rho}[q,q^{-1}]$
gives the decomposition of  $C_{Z^0_{-,+},\rho}(q,\bullet)$ in \eqref{CZ0-1} into the regular and polar parts:
$$
C_{Z^0_{-,+},\rho}(q,\bullet)=C^{ev}_{Z^0_{-,+},\rho}(q,\bullet) \oplus C^{odd}_{Z^0_{-,+},\rho}(q,\bullet) \in {\mathfrak g}^{(0)}_{ev}[q] +{\mathfrak g}^{(0)}_{odd}[q^{-1}]q^{-1}.
$$

Similarly, $Z^1_{-,+}$ will produce the decomposition
$$
\rho_{Z^1_{-,+}} (X,\rho_c)=\rho^{ev}_{Z^1_{-,+}} (X,\rho_c) \oplus \rho^{odd}_{Z^1_{-,+}} (X,\rho_c),
$$
where 
$$
\begin{gathered}
	\rho^{ev}_{Z^1_{-,+}} (X,\rho_c) \in {\mathfrak g}^{(0)}_{ev} =\bigoplus_{\text{$r$ even}} End(P^r),
	\\
	\rho^{odd}_{Z^1_{-,+}} (X,\rho_c) \in {\mathfrak g}^{(0)}_{odd} =\bigoplus_{\text{$r$ odd}} End(P^r)
\end{gathered}
$$
This gives the linear maps
$$
\begin{gathered}
	C^{ev}_{Z^1_{-,+},\rho}: P^+\longrightarrow {\mathfrak g}^{(0)}_{ev},
	\\
	C^{odd}_{Z^1_{-,+},\rho}: P^+ \longrightarrow {\mathfrak g}^{(0)}_{odd},
\end{gathered}
$$
defined respectively by the formulas
$$
\begin{gathered}
	C^{ev}_{Z^1_{-,+},\rho}(X)=	{}^{\dag}\rho^{ev}_{Z^1_{-,+}} (X,\rho_c) \in {\mathfrak g}^{(0)}_{ev} =\bigoplus_{\text{$s$ even}} End(P^s),
	\\
	C^{odd}_{Z^1_{-,+},\rho}(X)=	{}^{\dag}\rho^{odd}_{Z^1_{-,+}} (X,\rho_c) \in {\mathfrak g}^{(0)}_{odd} =\bigoplus_{\text{$s$ odd}} End(P^s).
\end{gathered}
$$
Extending those maps to the loop algebra ${\mathfrak g}^{(0)}[q,q^{-1}]$
gives the decomposition of  {\small$C_{Z^1_{-,+},\rho}(q,\bullet)$} in \eqref{CZ0-1} into the regular and polar parts:
$$
C_{Z^1_{-,+},\rho}(q,\bullet)=C^{ev}_{Z^1_{-,+},\rho}(q,\bullet) \oplus C^{odd}_{Z^1_{-,+},\rho}(q,\bullet) \in {\mathfrak g}^{(0)}_{ev}[q] +{\mathfrak g}^{(0)}_{odd}[q^{-1}]q^{-1}.
$$

We now fuse $Z^0_{-,+}$ and $Z^1_{-,+}$ by mixing up the regular and polar parts of $C_{Z^0_{-,+},\rho}(q,\bullet)$ and $C_{Z^1_{-,+},\rho}(q,\bullet)$, that is, we define
\begin{equation}\label{mix-ev0-odd1}
	\begin{gathered}
	C^{w}_{-,+}(\rho,q):=C^{ev}_{Z^0_{-,+}}(\rho,q) \oplus C^{odd}_{Z^1_{-,+}}(\rho,q): P^- \oplus P^+ \longrightarrow {\mathfrak g}^{(0)}_{ev}[q] +{\mathfrak g}^{(0)}_{odd}[q^{-1}]q^{-1},
	\\
		C^{b}_{-,+}(\rho,q):=C^{ev}_{Z^1_{-,+}}(\rho,q) \oplus C^{odd}_{Z^0_{-,+}}(\rho,q): P^- \oplus P^+ \longrightarrow {\mathfrak g}^{(0)}_{ev}[q] +{\mathfrak g}^{(0)}_{odd}[q^{-1}]q^{-1},
		\end{gathered}
\end{equation}
where the superscripts $w/b$ stands for `white' and `black', the color of vertices of $\widehat{PG}_l$: the operators composing $C^{w}_{-,+}(\rho,q)$ correspond to the closed paths in $\widehat{PG}_l$ based at the white vertices, while the ones composing 	$C^{b}_{-,+}(\rho,q)$ are based at the black ones.

 		Composing with the graded trace map 
 		$$
 		{\mathfrak g}^{(0)}_{ev}[q] +{\mathfrak g}^{(0)}_{odd}[q^{-1}]q^{-1} \subset {\mathfrak g}[q,q^{-1}] \longrightarrow \CC[q,q^{-1}]
 		$$
 		gives the linear maps
 		$$
 		\begin{gathered}
 		C^{w,-,+}_{\rho}(q,\bullet):=tr \circ 	C^{w}_{-,+}(\rho,q): P^- ([\xi],[\phi]) \oplus P^+([\xi],[\phi]) \longrightarrow \CC[q,q^{-1}],
 		\\
 		C^{b,-,+}_{\rho}(q,\bullet):=tr \circ 	C^{b}_{-,+}(\rho,q): P^-([\xi],[\phi]) \oplus P^+([\xi],[\phi]) \longrightarrow \CC[q,q^{-1}],
 		\end{gathered}
 		$$
 		together with their duals
 		$$
 		\begin{gathered}
 		\tau^{w,-,+}_{\rho}:=\left(C^{w,-,+}_{\rho} \right)^{\ast}: Hom_{\CC}(\CC[q,q^{-1}],\CC) \longrightarrow (P^- \oplus P^+)^{\ast}
 		\\
 		=(P^-)^{\ast} \oplus (P^+)^{\ast}=P^+ \oplus P^-,
 		\\
 		\tau^{b,-,+}_{\rho}:=\left(C^{b,-,+}_{\rho}\right)^{\ast}: Hom_{\CC}(\CC[q,q^{-1}],\CC) \longrightarrow (P^- \oplus P^+)^{\ast}
 		\\
 		=(P^-)^{\ast} \oplus (P^+)^{\ast}=P^+ \oplus P^-
 	\end{gathered}
 		$$
 		Thus the second way of fusing the two connected components of $Z_{-,+}$ assigns to every linear functional on $\CC[q,q^{-1}]$ a distinguished pair of linear maps between the ends $P^0$ and $P^{l-1}$ of the orthogonal decomposition of $W_{\xi}/W^l_{\xi}([\phi])$.

 	
 	\vspace{0.2cm}
 	\noindent
 	{\bf Odd case.} The zig-zag path 
 	$$
 	Z_{-,+}=(-e^+_0)e^-_2\cdots e^-_{l-2}(-e^+_{l-3})e^-_{l-1}(-e^+_{l-1})e^-_1(-e^+_1)\cdots(-e^+_{l-4})e^-_{l-2}(-e^+_{l-2})e^-_0
 	$$ 
 	is connected and it winds two times through each pair of vertices $(i)$ and $(i)'$: it passes through a white vertex $(i)$ makes the full turn of the torus returning to the same vertical level of the graph $\widehat{PG}_l$ (meridian on the torus) but this time at the black vertex $(i)'$ and then makes another full turn to return to $(i)$.
 	For example, in the writing of $Z_{-,+}$ above we started at $(0)$, after walking through $l$ edges counterclockwise 
 	$$
 	Z'_{-,+}=e^-_1(-e^+_1)\cdots(-e^+_{l-4})e^-_{l-2}(-e^+_{l-2})e^-_0
 	$$
 	we are back at the level $0$ but at the black vertex $(0)'$, the head of the arrow $e^-_1$; the remaining portion of the path
 	$$
 	Z''_{-,+}=(-e^+_0)e^-_2\cdots e^-_{l-3}(-e^+_{l-3})e^-_{l-1}(-e^+_{l-1})
 	$$
 	 winds from $(0)'$ to $(l-1)$ and back to $(0)$. Recall that we draw the graph $\widehat{PG}_l$ on the torus $\mathbb{T}$ where the white (resp. black) vertices are placed on `equatorial' circles in the counterclockwise manner from $(l-1)$ (resp. $(l-1)'$) to $(0)$ (resp. $(0)'$). The following drawing is an illustration for $\widehat{PG}_5$.
 	 
 	 \begin{equation}\label{Z-+5}
 	 	\begin{tikzpicture}
 	 		[place/.style={circle,draw=black,thick},
 	 		transition/.style={circle,draw=black,fill=black}]
 	 		\node (white4) at (90:1.5cm) [place] [label={below:$4$}] {};
 	 		\node (black4) at (90:3cm) [transition] [label={above:$4'$}]{};
 	 		\node (white3) at (162:1.5cm) [place] [label={below right:$3$}]{};
 	 		\node (black3) at (162:3cm) [transition] [label={above left:$3'$}]{};
 	 		\node (white2) at (234:1.5cm) [place] [label={above right:$2$}]{};
 	 		\node (black2) at (234:3cm) [transition] [label={below left:$2'$}]{};
 	 		\node (white1) at (306:1.5cm) [place] [label={above:$1$}] {};
 	 		\node (black1) at (306:3cm) [transition] [label={below:$1'$}]{};
 	 		\node (white0) at (388:1.5cm) [place] [label={below left:$0$}] {};
 	 		\node (black0) at (388:3cm) [transition] [label={below right:$0'$}]{};
 	 		\draw[green, very thick]
 	 		(white0) --(black0);
 	 		\begin{scope}[gray,ultra thin]
 	 			\draw (white4) to [bend right=30] (white3);
 	 			\draw (white3) to [bend right=30](white2);
 	 			\draw (white2) to [bend right=30](white1);
 	 			\draw (white1) to [bend right=30](white0);
 	 			\draw (white0) to [bend right=30](white4);
 	 			\draw (black4) to [bend right=30](black3);
 	 			\draw (black3) to [bend right=30](black2);
 	 			\draw (black2) to [bend right=30](black1);
 	 			\draw (black1) to [bend right=30](black0);
 	 			\draw (black0) to [bend right=30](black4);
 	 		\end{scope}
 	 		\begin{scope}
 	 			[decoration={markings,
 	 				mark =at position 2cm with {\arrow[blue,line width=1mm]{stealth}}}]
 	 			\draw[blue,thick][postaction={decorate}] (white3) to [bend right=45](black2);	
 	 			\draw[blue,thick][postaction={decorate}] (white2) to [bend right=45](black1);
 	 			\draw[blue,thick][postaction={decorate}] (white1) to [bend right=45](black0);
 	 			\draw[blue,thick][postaction={decorate}](white4) to [bend right=45](black3);	
 	 			\draw[blue,thick][postaction={decorate}] (white4) to [bend right=45](black3);	
 	 			\draw[blue,thick][postaction={decorate}] (white0) to [bend right=45](black4);	
 	 		\end{scope}
 	 		\begin{scope}[decoration={markings,
 	 				mark =at position 2cm with {\arrow[red,line width=1mm]{stealth}}}]
 	 			\draw[red, dotted,thick][postaction={decorate}] (black3) to [bend right=45](white2);
 	 			\draw[red,dotted,thick][postaction={decorate}] (black2) to [bend right=45](white1);
 	 			\draw[red,dotted,thick][postaction={decorate}] (black1) to [bend right=45](white0);
 	 			\draw[red,dotted,thick][postaction={decorate}] (black4) to [bend right=45](white3);
 	 			\draw[red,dotted,thick][postaction={decorate}] (black0) to [bend right=45](white4);
 	 		\end{scope}
 	 	\end{tikzpicture}
  	\end{equation}
  	On the drawing the blue colored edges are $e^-_s$ and the red ones are $(-e^+_t)$;  the path $Z_{-,+}$ based at $(0)$ is the concatenation of those edges in the following order
  	$$
  	Z_{-,+}=(-e^+_0)e^-_2(-e^+_2)e^-_4(-e^+_4)e^-_1(-e^+_1)e^-_3(-e^+_3)e^-_0;
  	$$
  	the green segment connecting $(0)$ and $(0)'$ is the cut separating $Z_{-,+}$ into two open paths:
  	$$
  	\begin{gathered}
  		Z'_{-,+}=e^-_1(-e^+_1)e^-_3(-e^+_3)e^-_0,
  		\\
  		Z''_{-,+}=(-e^+_0)e^-_2(-e^+_2)e^-_4(-e^+_4).
  	\end{gathered}
$$  	
  	
 We return to the the general considerations. To obtain the path operator corresponding to $Z_{-,+}$ we need to take a pair 
 	$$
 	(X,Y) \in P^- \times P^+
 	$$
 	as labels for $e^-_0$ and $e^+_{l-1}$ respectively; those labels together with the representation $\rho_c(\xi\otimes\phi)=\{\alpha^{t,s}_c(\xi,\phi)\}$ give a labeling of all arrows of the graph $\widehat{PG}_{l}$. With the same rule of labeling
 	$$
 	\begin{gathered}
 e^-_s \rightarrow \alpha^{s-1,s},\,\, \forall s\in [1,l-1],
 \\
 -e^+_s \rightarrow {}^{\dag}\alpha^{s+1,s},\,\, \forall s\in [0,l-2],
\end{gathered}
 $$	
as in the even case we obtain the operators
$$
\begin{gathered}
\rho_{Z'_{-,+}}(X,\rho_c):=\alpha^{0,1}({}^{\dag}\alpha^{2,1})\cdots \alpha^{0,1}({}^{\dag}\alpha^{l-1,l-2}) X: P^0 \longrightarrow P^0,
\\
\rho_{Z''_{-,+}}(Y,\rho_c):=({}^{\dag}\alpha^{1,0})\alpha^{1,2}\cdots \alpha^{l-4,l-3}({}^{\dag}\alpha^{l-2,l-3})\alpha^{l-2,l-1} ({}^{\dag}Y):P^0 \longrightarrow P^0.
\end{gathered}
$$
It should be observed that those are {\it not} closed path operators: the labels $P^0$ of each map are placed at the vertices of different colors on the $0$-th level, the level of cutting $Z_{-,+}$ into two open paths $Z'_{-,+}$ and $Z''_{-,+}$. The path operator associated to $Z_{-,+}$ based at the vertex $(0)$  is the composition of two operators above
$$
\rho^0_{Z_{-,+}}(X,Y,\rho_c):=\rho_{Z''_{-,+}}(Y,\rho_c) \circ \rho_{Z'_{-,+}}(X,\rho_c): P^0 \longrightarrow P^0.
$$
 The procedure described above mimics the situation of the even case: we cut the path $Z_{-,+}$ into two parts so that the edges $e^-_0$ and $(-e^+_{l-1})$ lie in different parts; this decomposes the path operator as the product of two where the variables $X$ and $Y$ are {\it separated}.
 
The path operator associated to $Z_{-,+}$ based at the vertex $(1)'$ is obtained by shifting the leftmost label of $\rho^0_{Z''_{-,+}}(Y,\rho_c)$, the map ${}^{\dag}\alpha^{1,0}$, to the rightmost position in $\rho_{Z'_{-,+}}(X,\rho_c)$. We can realize this again by cutting the shifted (from left to right) path $\circlearrowleft (Z_{-,+})$
$$
\circlearrowleft (Z_{-,+})=e^-_2\cdots e^-_{l-2}(-e^+_{l-3})e^-_{l-1}(-e^+_{l-1})e^-_1(-e^+_1)\cdots(-e^+_{l-4})e^-_{l-2}(-e^+_{l-2})e^-_0 (-e^+_0)
$$
into two parts
$$
\begin{gathered}
(\circlearrowleft (Z_{-,+}))':=(-e^+_1)\cdots(-e^+_{l-4})e^-_{l-2}(-e^+_{l-2})e^-_0 (-e^+_0),
\\
(\circlearrowleft (Z_{-,+}))'':=e^-_2\cdots e^-_{l-2}(-e^+_{l-3})e^-_{l-1}(-e^+_{l-1})e^-_1.
\end{gathered}
$$
The corresponding path operators are:
$$
\begin{gathered}
	\rho_{(\circlearrowleft (Z_{-,+}))'} (X,\rho_c): P^1 \longrightarrow P^1,
	\\
	\rho_{(\circlearrowleft (Z_{-,+}))''} (Y,\rho_c): P^1 \longrightarrow P^1.
\end{gathered}
$$
The path operator associated to $Z_{-,+}$ based at the vertex $(1)'$ is the composition of the two operators above
$$
\rho^1_{Z_{-,+}}(X,Y,\rho_c)=\rho_{(\circlearrowleft (Z_{-,+}))''} (\rho_c,Y)\circ\rho_{(\circlearrowleft (Z_{-,+}))'} (\rho_c,X): P^1 \longrightarrow P^1.
$$
After  cyclically shifting $s$ times and cutting along the segment joining $(s)$ to $(s)'$ we obtain the decomposition of $\circlearrowleft^s (Z_{-,+})$ into the product of two open paths
$$
\circlearrowleft^s (Z_{-,+})=(\circlearrowleft^s (Z_{-,+}))'' (\circlearrowleft^s (Z_{-,+}))';
$$
the path operator associated to $Z_{-,+}$ based at the vertex $(s)$, if $s$ is even, and at $(s)'$, if $s$ is odd:
$$
\rho^s_{Z_{-,+}}(X,Y,\rho_c)=\rho_{(\circlearrowleft^s (Z_{-,+}))''} (Y,\rho_c)\circ\rho_{(\circlearrowleft^s (Z_{-,+}))'} (X,\rho_c): P^s \longrightarrow P^s.
$$
Taking the direct sum we obtain the grading preserving homomorphism
$$
\rho_{Z_{-,+}}(X,Y,\rho_c)=\bigoplus^{l-1}_{s=0}\rho^s_{Z_{-,+}}(X,Y,\rho_c): \bigoplus^{l-1}_{s=0} P^s \longrightarrow \bigoplus^{l-1}_{s=0} P^s
$$
 which is the (graded) composition of two maps
 \begin{equation}\label{rhoXY-+}
 \rho_{Z_{-,+}}(X,Y,\rho_c)=\rho''_{Z_{-,+}}(Y,\rho_c) \circ \rho'_{Z_{-,+}}(X,\rho_c), 
\end{equation}
 where the factors in the composition are as follows
 \begin{equation}\label{rhoXY-+factors}
 \begin{gathered}
 	\rho''_{Z_{-,+}}(Y,\rho_c):=\bigoplus^{l-1}_{s=0}\rho_{(\circlearrowleft^s (Z_{-,+}))''} (Y,\rho_c),
 	\\
 	\rho'_{Z_{-,+}}(X,\rho_c):=\bigoplus^{l-1}_{s=0}\rho_{(\circlearrowleft^s (Z_{-,+}))'} (X,\rho_c).
 \end{gathered}
\end{equation}
Observe: $\rho_{Z_{-,+}}(X,Y,\rho_c)$ as a function of $X$ and $Y$ is linear and anti-linear respectively. Replacing 
$\rho''_{Z_{-,+}}(Y,\rho_c)$ in \eqref{rhoXY-+} by its adjoint gives the $\CC$-bilinear pairing
$$
B^{+,-}_{Z_{-,+},\rho}: P^+\times P^- \longrightarrow {\mathfrak g}^{(0)}
$$
defined by the formula
\begin{equation}\label{pairing-Z-+-odd}
	B^{+,-}_{Z_{-,+},\rho}(Y,X)={}^{\dag}\rho''_{Z_{-,+}}(Y,\rho_c) \circ \rho'_{Z_{-,+}}(X,\rho_c).
\end{equation}
Passing to the loop algebra ${\mathfrak g}^{(0)}[q,q^{-1}]$ we obtain
$$
B^{+,-}_{Z_{-,+},\rho}(q,\bullet,\bullet): P^+\times P^- \longrightarrow {\mathfrak g}^{(0)}[q,q^{-1}].
$$
This factors as follows
\begin{equation}\label{pairing-Z-+-odd-q}
	B^{+,-}_{Z_{-,+},\rho}(q,Y,X)={}^{\dag}\rho''_{Z_{-,+}}(Y,\rho_c,q) \circ \rho'_{Z_{-,+}}(X,\rho_c,q),
\end{equation}
where the factors are defined by the formulas
\begin{equation}
	\begin{gathered}
			{}^{\dag}\rho''_{Z_{-,+}}(Y,\rho_c,q):=\sum^{l-1}_{s=0}{}^{\dag}\rho_{(\circlearrowleft^s (Z_{-,+}))''} (Y,\rho_c)q^{(-1)^s s},
		\\
		\rho'_{Z_{-,+}}(X,\rho_c,q):=\sum^{l-1}_{s=0}\rho_{(\circlearrowleft^s (Z_{-,+}))'} (X,\rho_c)q^{(-1)^s s}.
	\end{gathered}
\end{equation}

Exchanging the factors in \eqref{pairing-Z-+-odd-q} gives another bilinear pairing
\begin{equation}\label{B-+-pairing}
	B^{-,+}_{Z_{-,+},\rho}(q,\bullet,\bullet):P^-\times P^+ \longrightarrow {\mathfrak g}^{(0)}[q,q^{-1}]
\end{equation}
Composing with the trace map
$$
tr: {\mathfrak g}^{(0)}[q,q^{-1}] \longrightarrow \CC[q,q^{-1}]
$$
we obtain the linear maps
\begin{equation}\label{CZ-+-odd}
	\begin{gathered}
			C^+_{Z_{-,+},\rho_c}(q,\bullet):=
			tr \circ {}^{\dag}\rho''_{Z_{-+}}(\bullet,\rho_c,q):P^+ \longrightarrow \CC[q,q^{-1}],
			\\
	C^-_{Z_{-,+}\rho_c}(q,\bullet):=
	tr \circ	\rho'_{Z_{-,+}}(\bullet,\rho_c,q):  P^-\longrightarrow \CC[q,q^{-1}].
	\end{gathered}
\end{equation}
given by the formulas
$$
\begin{gathered}
C^+_{Z_{-,+}}(Y,\rho_c,q)=\sum^{l-1}_{s=0}tr \left({}^{\dag}\rho_{(\circlearrowleft^s (Z_{-,+}))''} (Y,\rho_c)\right)q^{(-1)^s s}, \,\,\forall Y\in P^+
\\
C^-_{Z_{-,+}}(X,\rho_c,q)=\sum^{l-1}_{s=0}tr \left(\rho_{(\circlearrowleft^s (Z_{-,+}))'} (X,\rho_c)\right)q^{(-1)^s s},\,\,\forall X\in P^-.
\end{gathered}
$$
Dualizing the maps in \eqref{CZ-+-odd} gives the homomorphism
\begin{equation}\label{tau-+-odd}
	\begin{gathered}
		\tau^-_{Z_{-,+},\rho}:=\left(C^+_{Z_{-,+},\rho_c}\right)^{\ast}:Hom_{\CC}(\CC[q,q^{-1}],\CC)\longrightarrow (P^+)^{\ast}\cong P^-,
		\\
		\tau^+_{Z_{-,+},\rho}:=\left(C^-_{Z_{-,+},\rho_c}\right)^{\ast}:Hom_{\CC}(\CC[q,q^{-1}],\CC)\longrightarrow (P^-)^{\ast}\cong P^+.
	\end{gathered}
\end{equation}

We also have the traced version of the pairings $B^{+,-}_{Z_{-,+},\rho}(q,\bullet,\bullet) $ and $B^{-,+}_{Z_{-,+},\rho}(q,\bullet,\bullet) $:
$$
\begin{gathered}
tr\circ B^{+,-}_{Z_{+,-},\rho}(q,\bullet,\bullet) : P^+\times P^- \longrightarrow \CC[q,q^{-1}]
\\
tr\circ B^{-,+}_{Z_{-,+},\rho}(q,\bullet,\bullet) : P^-\times P^+ \longrightarrow \CC[q,q^{-1}]
\end{gathered}
$$
defined by the formulas
\begin{equation}\label{trBZ-+odd}
	\begin{gathered}
	(tr B^{+,-}_{Z_{-,+},\rho}(q,Y,X)=tr \left({}^{\dag}\rho''_{Z_{-,+}}(Y,\rho_c,q) \circ \rho'_{Z_{-,+}}(X,\rho_c,q) \right)
	\\
	=\sum^{l-1}_{s=0} tr \left({}^{\dag}\rho_{(\circlearrowleft^s (Z_{-,+}))''} (Y,\rho_c) \circ \rho_{(\circlearrowleft^s (Z_{-,+}))'} (X,\rho_c) \right)q^{(-1)^s s}, \,\,\forall Y\in P^+, \forall X\in P^-;
	\\
		(tr B^{-,+}_{Z_{-,+},\rho}(q,X,Y)=tr \left(\rho'_{Z_{-,+}}(X,\rho_c,q) \circ {}^{\dag}\rho''_{Z_{-,+}}(Y,\rho_c,q) \right)
	\\
	=\sum^{l-1}_{s=0} tr \left(\rho_{(\circlearrowleft^s (Z_{-,+}))'} (X,\rho_c) \circ {}^{\dag}\rho_{(\circlearrowleft^s (Z_{-,+}))''} (Y,\rho_c)  \right)q^{(-1)^s s}, \,\,\forall Y\in P^+, \forall X\in P^-
\end{gathered}
\end{equation}
 \begin{example}
 	We take up the case $l=5$ and write out all the maps discussed above.
 	We have already seen the path $Z_{-,+}$ based at $(0)$ and its decomposition into two open paths $Z'_{-,+}$ and $Z''_{-,+}$:
 	$$
 	\begin{gathered}
 		Z_{-,+}=(-e^+_0)e^-_2 (-e^+_2)e^-_4 (-e^+_4)e^-_1(-e^+_1)e^-_3(-e^+_3)e^-_0,
 		\\
 		Z'_{-,+}=e^-_1(-e^+_1)e^-_3(-e^+_3)e^-_0,
 		\\
 		Z_{-,+}=(-e^+_0)e^-_2 (-e^+_2)e^-_4 (-e^+_4);
 	\end{gathered}
 $$
 the labeling 
 $$
 \begin{gathered}
 e^-_0 \rightarrow X \in P^-, \,\, e^+_4  \rightarrow Y \in P^+,
 \\
 e^-_s \rightarrow \alpha^{s-1,s}, \,\,  e^+_s \rightarrow \alpha^{s,s-1},\,\,\forall s\in [1,4],
\end{gathered}
$$
and the rule of taking the adjoint when $e^+_s$ changes the orientation gives the following path operators
$$
\begin{gathered}
\rho_{Z'_{-,+}} (X,\rho)= \alpha^{0,1}({}^{\dag}\alpha^{2,1})\alpha^{2,3}({}^{\dag}\alpha^{4,3})X,
	\\
\rho_{Z''_{-,+}} (Y,\rho)=({}^{\dag}\alpha^{1,0})\alpha^{1,2} ({}^{\dag}\alpha^{3,2})\alpha^{3,4} ({}^{\dag}Y), 
	\\
	\rho^0_{Z_{-,+}} (X,Y,\rho)=\rho_{Z''_{-,+}} (Y,\rho) \circ \rho_{Z'_{-,+}} (X,\rho)
	\\
	=\underbrace{({}^{\dag}\alpha^{1,0})\alpha^{1,2} ({}^{\dag}\alpha^{3,2})\alpha^{3,4} ({}^{\dag}Y)}_{Z''_{-,+}} \underbrace{\alpha^{0,1}({}^{\dag}\alpha^{2,1})\alpha^{2,3}({}^{\dag}\alpha^{4,3})X}_{Z'_{-,+}};
\end{gathered}
$$
shifting from left to right we obtain the component of degree one:
$$
\begin{gathered}
 \circlearrowleft (Z_{-,+})=e^-_2 (-e^+_2)e^-_4 (-e^+_4)e^-_1(-e^+_1)e^-_3(-e^+_3)e^-_0 (-e^+_0) \rightsquigarrow \rho^1_{Z_{-,+}}(X,Y,\rho)
 \\
 =\rho_{(\circlearrowleft (Z_{-,+}))''}(Y,\rho) \circ \rho_{(\circlearrowleft (Z_{-,+}))'}(X,\rho)
 =\underbrace{\alpha^{1,2} ({}^{\dag}\alpha^{3,2})\alpha^{3,4} ({}^{\dag}Y)\alpha^{0,1}}_{(\circlearrowleft (Z_{-,+}))''} \underbrace{({}^{\dag}\alpha^{2,1})\alpha^{2,3}({}^{\dag}\alpha^{4,3})X({}^{\dag}\alpha^{1,0})}_{(\circlearrowleft (Z_{-,+}))'},
 \\
 (\circlearrowleft (Z_{-,+}))'=(-e^+_1)e^-_3(-e^+_3)e^-_0 (-e^+_0) \rightsquigarrow \rho_{(\circlearrowleft (Z_{-,+}))'}(X,\rho) =({}^{\dag}\alpha^{2,1})\alpha^{2,3}({}^{\dag}\alpha^{4,3})\alpha^{1,0}X({}^{\dag}\alpha^{1,0}),
 \\
 (\circlearrowleft (Z_{-,+}))''=e^-_2 (-e^+_2)e^-_4 (-e^+_4)e^-_1\rightsquigarrow \rho_{(\circlearrowleft (Z_{-,+}))'}(Y,\rho)=\alpha^{1,2}({}^{\dag}\alpha^{3,2})\alpha^{3,4}({}^{\dag}Y\alpha^{0,1});
\end{gathered}
$$
continuing in this manner gives other graded pieces with their factorizations
$$
\begin{gathered}
\circlearrowleft^2 (Z_{-,+})= (-e^+_2)e^-_4 (-e^+_4)e^-_1(-e^+_1)e^-_3(-e^+_3)e^-_0 (-e^+_0)e^-_2  \rightsquigarrow \rho^2_{Z_{-,+}}(X,Y,\rho)
\\
=\rho_{(\circlearrowleft^2 (Z_{-,+}))''}(Y,\rho) \circ \rho_{(\circlearrowleft^2 (Z_{-,+}))'}(X,\rho)
=\underbrace{ ({}^{\dag}\alpha^{3,2})\alpha^{3,4} ({}^{\dag}Y)\alpha^{0,1}({}^{\dag}\alpha^{2,1})}_{(\circlearrowleft^2 (Z_{-,+}))''} \underbrace{\alpha^{2,3}({}^{\dag}\alpha^{4,3})X({}^{\dag}\alpha^{1,0})\alpha^{1,2}}_{(\circlearrowleft^2 (Z_{-,+}))'},
\end{gathered}
$$
$$
\begin{gathered}
	\circlearrowleft^3 (Z_{-,+})= e^-_4 (-e^+_4)e^-_1(-e^+_1)e^-_3(-e^+_3)e^-_0 (-e^+_0)e^-_2(-e^+_2)  \rightsquigarrow \rho^3_{Z_{-,+}}(X,Y,\rho)
	\\
	=\rho_{(\circlearrowleft^3 (Z_{-,+}))''}(Y,\rho) \circ \rho_{(\circlearrowleft^3 (Z_{-,+}))'}(X,\rho)
	=\underbrace{ \alpha^{3,4} ({}^{\dag}Y)\alpha^{0,1}({}^{\dag}\alpha^{2,1})\alpha^{2,3}}_{(\circlearrowleft^3 (Z_{-,+}))''} \underbrace{({}^{\dag}\alpha^{4,3})X({}^{\dag}\alpha^{1,0})\alpha^{1,2}({}^{\dag}\alpha^{3,2})}_{(\circlearrowleft^3 (Z_{-,+}))'},
\end{gathered}
$$
$$
\begin{gathered}
	\circlearrowleft^4 (Z_{-,+})=  (-e^+_4)e^-_1(-e^+_1)e^-_3(-e^+_3)e^-_0 (-e^+_0)e^-_2(-e^+_2)e^-_4  \rightsquigarrow \rho^4_{Z_{-,+}}(X,Y,\rho)
	\\
	=\rho_{(\circlearrowleft^4 (Z_{-,+}))''}(Y,\rho) \circ \rho_{(\circlearrowleft^4 (Z_{-,+}))'}(X,\rho)
	=\underbrace{  ({}^{\dag}Y)\alpha^{0,1}({}^{\dag}\alpha^{2,1})\alpha^{2,3}({}^{\dag}\alpha^{4,3})}_{(\circlearrowleft^4 (Z_{-,+}))''} \underbrace{X({}^{\dag}\alpha^{1,0})\alpha^{1,2}({}^{\dag}\alpha^{3,2})\alpha^{3,4}}_{(\circlearrowleft^4 (Z_{-,+}))'}.
\end{gathered}
$$
Collecting the factors of the formulas above we obtain the ${\mathfrak g}^{(0)}$-valued Laurent polynomials
$$
\begin{gathered}
	\rho''_{Z_{-,+}} (Y,\rho,q)=\left[({}^{\dag}\alpha^{1,0})\alpha^{1,2} ({}^{\dag}\alpha^{3,2})\alpha^{3,4} ({}^{\dag}Y)\right] +\left[\alpha^{1,2} ({}^{\dag}\alpha^{3,2})\alpha^{3,4} ({}^{\dag}Y)\alpha^{0,1} \right] q^{-1} 
	\\
	+
	\left[({}^{\dag}\alpha^{3,2})\alpha^{3,4} ({}^{\dag}Y)\alpha^{0,1}({}^{\dag}\alpha^{2,1})\right]q^2
	\\ 
	+
	\left[\alpha^{3,4} ({}^{\dag}Y)\alpha^{0,1}({}^{\dag}\alpha^{2,1})\alpha^{2,3} \right]q^{-3}
	+\left[({}^{\dag}Y)\alpha^{0,1}({}^{\dag}\alpha^{2,1})\alpha^{2,3}({}^{\dag}\alpha^{4,3})\right]q^4,
\end{gathered}
$$
$$
\begin{gathered}
	\rho'_{Z_{-,+}} (X,\rho,q)=\alpha^{0,1}({}^{\dag}\alpha^{2,1})\alpha^{2,3}({}^{\dag}\alpha^{4,3})X +\left[({}^{\dag}\alpha^{2,1})\alpha^{2,3}({}^{\dag}\alpha^{4,3})X({}^{\dag}\alpha^{1,0})\right] q^{-1}
	\\
	 +
	\left[\alpha^{2,3}({}^{\dag}\alpha^{4,3})X({}^{\dag}\alpha^{1,0})\alpha^{1,2}\right]q^2 
	\\
	+
	\left[({}^{\dag}\alpha^{4,3})X({}^{\dag}\alpha^{1,0})\alpha^{1,2}({}^{\dag}\alpha^{3,2})\right] q^{-3}
	+\left[X({}^{\dag}\alpha^{1,0})\alpha^{1,2}({}^{\dag}\alpha^{3,2})\alpha^{3,4}\right]q^4.
\end{gathered}
$$
The pairing
$$
B^{+,-}_{Z_{-,+},\rho} (q,\bullet,\bullet): P^+\times P^- \longrightarrow {\mathfrak g}^{(0)}[q,q^{-1}]
	$$
	is defined by the formula
	$$
	B^{+,-}_{Z_{-,+},\rho} (q,Y,X)={}^{\dag}\rho''_{Z_{-,+}} (Y,\rho,q)\rho'_{Z_{-,+}} (X,\rho,q),
	$$
	where the composition is performed in the graded sense. Inverting the factors gives the pairing $	B^{-,+}_{Z_{-,+},\rho} (q,\bullet,\bullet)$
	$$
	B^{-,+}_{Z_{-,+},\rho} (q,X,Y)=\rho'_{Z_{-,+}} (X,\rho,q) {}^{\dag}\rho''_{Z_{-,+}} (Y,\rho,q).
	$$
 \end{example}
 	
 	\vspace{0.2cm}
 We return to the general considerations. The $\CC$-bilinear pairing 
 	$$
 	tr \circ B^{+,-}_{Z_{-,+},\rho}: P^+ \times P^- \longrightarrow \CC[q,q^{-1}]
 	$$
 	is equivalent to the associated $\CC$-linear map denoted $	tr B^{+,-}_{Z_{-,+},\rho}$
 	\begin{equation}\label{pairZ-+odd}
 	tr B^{+,-}_{Z_{-,+},\rho}:	P^+ \otimes P^- \longrightarrow \CC[q,q^{-1}]
 \end{equation}
Dualizing gives the map
\begin{equation}\label{dpairZ-+odd}
	\begin{gathered}
	tr^{\ast} B^{+,-}_{Z_{-,+},\rho}:	Hom_{\CC}(\CC[q,q^{-1}],\CC) \longrightarrow (P^+ \otimes P^-)^{\ast} \cong (P^-)^{\ast} \otimes (P^+)^{\ast} 
	\\
 \cong P^+ \otimes (P^+)^{\ast} \cong End(P^+).
	\end{gathered}
\end{equation}	
Similarly, by exchange of factors, we obtain the map
\begin{equation}\label{-+dpairZ-+odd}
	\begin{gathered}
		tr^{\ast} B^{-,+}_{Z_{-,+},\rho}:	Hom_{\CC}(\CC[q,q^{-1}],\CC) \longrightarrow (P^- \otimes P^+)^{\ast} \cong (P^+)^{\ast} \otimes (P^-)^{\ast} 
		\\
		\cong P^- \otimes (P^-)^{\ast} \cong  End(P^-).
	\end{gathered}
\end{equation}	
 	We summarize the above considerations in the following statement.
 	\begin{pro}\label{pro:Z-+lodd}
 		Let $\rho_c(\xi\otimes\phi)=\{\alpha^{t,s}_c(\xi,\phi): P^s \longrightarrow P^{t}\}$ be the representation of the quiver $PG_l$ associated to a nonzero element $\xi\otimes\phi$ of the total space
 		of the tautological line bundle $\OO_{{\mathfrak{L}}_l}(-1)$.
 	 For $l=2k+1$ odd, the zig-zag path
 		$Z_{-,+}$ of the graph $\widehat{PG}_l$ is connected but it can be cut into two `open' paths of length $l$ each: 
 		$$
 		\circlearrowleft^s(Z_{-,+})=
 		\left(\circlearrowleft^s(Z_{-,+})\right)''\left(\circlearrowleft^s(Z_{-,+})\right)',
 		$$
 		where  $\circlearrowleft^s(Z_{-,+})$ denotes the path $Z_{-,+}$ with edges ordered so that its tail and head is at the white vertex $(s)$ of the graph $\widehat{PG}_l$. These decompositions give rise to two distinguished linear maps
 		$$
 		\begin{gathered}
 			\tau^+_{Z_{-,+},\rho}:=(C^-_{Z_{-,+,\rho}})^{\ast}: Hom(\CC [q,q^{-1}],\CC) \longrightarrow (P^-)^{\ast} \cong P^+,
 			\\
 			\tau^-_{Z_{-,+,\rho}}:=(C^+_{Z_{-,+,\rho}})^{\ast}:Hom(\CC [q,q^{-1}],\CC) \longrightarrow (P^+)^{\ast} \cong P^- ,
 		\end{gathered}
 		$$
 		defined in \eqref{tau-+-odd}. In addition, fusing the two maps together gives the distinguished linear maps
 		$$
 		\begin{gathered}
 			tr^{\ast}(B^{+,-}_{Z_{-,+,\rho}}):Hom(\CC[q,q^{-1}],\CC) \longrightarrow   End (P^+),
 			\\
 			tr^{\ast}(B^{-,+}_{Z_{-,+,\rho}} ):Hom(\CC[q,q^{-1}],\CC) \longrightarrow  End (P^-),
 		\end{gathered}
 		$$
 		see \eqref{dpairZ-+odd} and \eqref{-+dpairZ-+odd}. In the above formulas $\rho=\rho_c$ and the maps above are independent of the trace parameter $c$.
 	\end{pro}
 	\subsection{Representations of the dual quiver $Q_l$}
 	We have seen how to extend the representations of the quiver $PG_l$
 to representations of the quiver $\widehat{PG}_l$, see Proposition \ref{pro:mu-map}. The graph $Q_l$ dual to $\widehat{PG}_l$ is naturally a quiver and we studied its path algebra $\CC Q_l$. It is an interesting object: it comes with potential $W_l$ and its Jacobi algebra $A_l$ is a Calabi-Yau algebra of dimension three, see \S8. In this subsection we show how to produce representations of the path algebra $\CC Q_l$ from the representations of $PG_l$.
 	
 	Let a 
 	representation 
 	$$
 	\rho=\{\alpha^{t,s}: P^s \longrightarrow P^{t}\}
 	$$
 	 of the quiver $PG_l$ be given. We have learned how to extend this to the representations of the quiver $\widehat{PG}_l$: recall that we have constructed the linear maps
 	 $$
 	 \tau^{\pm}_{\rho}: Hom(\CC[q,q^{-1}],\CC) \longrightarrow P^{\pm},
 	 $$
 	where $P^{-}=Hom(P^0,P^{l-1})$ and $P^{+}=Hom(P^{l-1},P^0)$ the spaces of maps between the `ends' of the representation $\rho$; given two linear functionals
 	$F,G \in Hom(\CC[q,q^{-1}],\CC)$ we obtain linear maps
 	$$
 	\tau^{-}_{\rho}(F):P^0 \longrightarrow P^{l-1},\,\,\, 	\tau^{+}_{\rho}(G):P^{l-1} \longrightarrow P^{0}
 	$$
 	which are used to label respectively the edges $e^-_0$ and $e^+_{l-1}$
 	of $\widehat{PG}_l$; these labels together with the representation $\rho$ constitute a representation of the quiver $\widehat{PG}_l$
 	extending $\rho$. This representation is denoted $\widehat{\rho} (F,G)$. We will now describe how to produce from $\widehat{\rho} (F,G)$ a representation of the dual quiver $Q_l$. In fact, we will see momentarily that there is a finite collection of such representations naturally associated to $\widehat{\rho} (F,G)$.
 	
 	\begin{pro}\label{pro:repsQl}
 		Let 
 		$$
 		\rho=\{\alpha^{t,s}: P^s\longrightarrow P^{t} \}
 		$$
 		be the representation of $PG_l$.
 		For every pair of linear functionals $F,G \in Hom(\CC[q,q^{-1}],\CC)$, the representation $\widehat{\rho} (F,G)$
 		of $\widehat{PG}_l$ gives rise to a distinguished finite collection $$
 		Reps ({Q_l})(\widehat{\rho} (F,G))
 		$$
 		 of representations of the dual quiver $Q_l$:
 		$$
 		Reps(\widehat{PG}_l) \ni \widehat{\rho} (F,G) \mapsto  Reps ({Q_l})(\widehat{\rho} (F,G)) \subset Reps(Q_l).
 		$$
 		The representations of the collection $ Reps ({Q_l})(\widehat{\rho} (F,G))$ are labeled by the subsets $A$ of $[0,l-1]$ and  denoted $(\widehat{\rho} (F,G))^{\vee}_A$. 
 		In particular, to every
 		 irreducible component of the Lagrangian cycle $H_0$ of the nonabelian Dolbeault variety ${\bf H^{1,0}}(PG_l)$ is attached a representation of $Q_l$; it is given by the rule
 	$$
 H_0 \supset	\Pi_A \rightarrow	  (\widehat{\rho} (F,G))^{\vee}_A  \in Reps ({Q_l})(\widehat{\rho}_c (F,G)),
 	$$
 	where $\Pi_A$ is the irreducible component of $H_0$ corresponding to the subset $A \subset [1,l-1] \subset [0,l-1]$, see Proposition \ref{pro:NADolb} for the definition of $\Pi_A$.
 	\end{pro}
 \begin{pf}
 	We fix the representation $\widehat{\rho} (F,G)$ of $\widehat{PG}_l$ and we want to build out of it a representation of the dual quiver $Q_l$. Denote provisionally that representation $(\widehat{\rho} (F,G))^{\vee}$. We start by attaching vector spaces to the vertices of $Q_l$.
 	
 	The reader may recall that the vertices of $Q_l$ are $\{\check{F}_i\}$, for $i\in[0,l-1]$; they correspond to the faces $\{F_i\}$ of $\widehat{PG}_l$
 	and those are hexagons
 	$$
 		\begin{tikzpicture}
 			[place/.style={circle,draw=black,thick, inner sep=0pt, minimum size=2mm},
 			transition/.style={circle,draw=black,fill=black, inner sep=0pt, minimum size=2mm}]
 			\node (whitei) at (150:1.5cm) [place] [label={left:$i$}] {};
 			\node (blacki) at (210:1.5cm) [transition] [label={left:$i'$}]{};
 			\node (whitei1) at (270:1.5cm) [place] [label={below:$i-1$}]{};
 			\node (blacki2) at (330:1.5cm) [transition] [label={right:$(i-2)'$}]{};
 			\node (whitei2) at (30:1.5cm) [place] [label={ right:$i-2$}]{};
 			\node (blacki1) at (90:1.5cm) [transition] [label={above:$(i-1)'$}]{};
 			\node at (0,0) [circle, draw=red, fill=red, inner sep=0pt, minimum size=2mm] [label={below:$F_i$}]{};
 			\draw[black,ultra thick][-] (whitei) to (blacki);
 			\draw[black,ultra thick][-] (whitei1) to (blacki);
 			\draw[black,ultra thick][-] (whitei2) to (blacki2);
 			\draw[black,ultra thick][-] (whitei2) to (blacki1);
 			\draw[black,ultra thick][-] (whitei) to (blacki1);
 			\draw[black,ultra thick][-] (whitei1) to (blacki2);
 		\end{tikzpicture} 
 	$$
 	corresponding to the boundary cycles $\{B_i\}$ of $\widehat{PG}_l$. Our representation $(\widehat{\rho} (F,G))^{\vee}$ attaches to the vertex $\check{F}_i$ the vector space
 	$Hom(P^i,P^{i-2})$
 	$$
 	(\widehat{\rho}_c (F,G))^{\vee}: \check{F}_i \mapsto Hom(P^i,P^{i-2}),
 	$$
 	for each $i\in [0,l-1]$. Heuristically, we can think of this space as the ways to get from the edge $e^0_i=(i \to i')$ of the hexagon to its opposite edge $e^0_{i-2}=((i-2) \to (i-2)')$. To simplify the notation we denote the space
$ Hom(P^i,P^{i-2})$ by $P^{i,+}$.

	We now turn to labeling the edges of $Q_l$.
 	There are three edges of $Q_l$ emanating from each $\check{F}_i$
 	\begin{equation}
 		(e^+_{i-2})^{\check{}}: \check{F}_i \rightarrow \check{F}_{i-1},\,\, (e^-_{i-1})^{\check{}}: \check{F}_i \rightarrow \check{F}_{i-1},\,\,(e^0_{i})^{\check{}}: \check{F}_i \rightarrow \check{F}_{i+2}. 
 	\end{equation}
 
 In the following drawing the vertices of $Q_l$ are red dots and the edges based at $\check{F}_i$ are in blue:
 
 $$
 \begin{tikzpicture}
 	[place/.style={circle,draw=black,thick, inner sep=0pt, minimum size=2mm},
 	transition/.style={circle,draw=black,fill=black, inner sep=0pt, minimum size=2mm}]
 	\node (whitei) at (150:1.5cm) [place] [label={left:$\scriptstyle{i}$}] {};
 	\node (blacki) at (210:1.5cm) [transition] [label={left:$\scriptstyle{i'}$}]{};
 	\node (whitei1) at (270:1.5cm) [place] [label={below left:$\scriptstyle{i-1}$}]{};
 	\node (blacki2) at (330:1.5cm) [transition] [label={below:$\scriptstyle{(i-2)'}$}]{};
 	\node (whitei2) at (30:1.5cm) [place] [label={ above:$\scriptstyle{i-2}$}]{};
 	\node (blacki1) at (90:1.5cm) [transition] [label={}]{};
 	\node (Fi) at (0,0) [circle, draw=red, fill=red, inner sep=0pt, minimum size=2mm] [label={below:$F_i$}]{};
 	\draw[black,ultra thick][-] (whitei) to (blacki);
 	\draw[black,ultra thick][-] (whitei1) to (blacki);
 	\draw[black,ultra thick][-] (whitei2) to (blacki2);
 	\draw[black,ultra thick][-] (whitei2) to (blacki1);
 	\draw[black,ultra thick][-] (whitei) to (blacki1);
 	\draw[black,ultra thick][-] (whitei1) to (blacki2);
 	\begin{scope}[xshift=1.26cm, yshift=2.25cm]
 			[place/.style={circle,draw=black,thick, inner sep=0pt, minimum size=2mm},
 			transition/.style={circle,draw=black,fill=black, inner sep=0pt, minimum size=2mm}]
 			\node (whitei) at (150:1.5cm) [place] [label={above:$\scriptstyle{i-1}$}] {};
 			\node (blacki) at (210:1.5cm) [transition] [label={ below:$\scriptstyle{(i-1)'}$}]{};
 			\node (whitei1) at (270:1.5cm) [place] [label={}]{};
 			\node (blacki2) at (330:1.5cm) [transition] [label={right:$\scriptstyle{(i-3)'}$}]{};
 			\node (whitei2) at (30:1.5cm) [place] [label={ right:$\scriptstyle{i-3}$}]{};
 			\node (blacki1) at (90:1.5cm) [transition] [label={above:$\scriptstyle{(i-2)'}$}]{};
 			\node (Fi-1) at (0,0) [circle, draw=red, fill=red, inner sep=0pt, minimum size=2mm] [label={above:$F_{i-1}$}]{};
 			\draw[black,ultra thick][-] (whitei) to (blacki);
 			\draw[black,ultra thick][-] (whitei1) to (blacki);
 			\draw[black,ultra thick][-] (whitei2) to (blacki2);
 			\draw[black,ultra thick][-] (whitei2) to (blacki1);
 			\draw[black,ultra thick][-] (whitei) to (blacki1);
 			\draw[black,ultra thick][-] (whitei1) to (blacki2);
 	\end{scope}
 	\begin{scope}[xshift=1.26cm, yshift=-2.25cm]
 	[place/.style={circle,draw=black,thick, inner sep=0pt, minimum size=2mm},
 	transition/.style={circle,draw=black,fill=black, inner sep=0pt, minimum size=2mm}]
 	\node (whitei) at (150:1.5cm) [place] [label={}] {};
 	\node (blacki) at (210:1.5cm) [transition] [label={below :$\scriptstyle{(i-1)'}$}]{};
 	\node (whitei1) at (270:1.5cm) [place] [label={below:$\scriptstyle{i-2}$}]{};
 	\node (blacki2) at (330:1.5cm) [transition] [label={right:$\scriptstyle{(i-3)'}$}]{};
 	\node (whitei2) at (30:1.5cm) [place] [label={ right:$\scriptstyle{i-3}$}]{};
 	\node (blacki1) at (90:1.5cm) [transition] [label={}]{};
 	\node (Fi-1s)at (0,0) [circle, draw=red, fill=red, inner sep=0pt, minimum size=2mm] [label={below:$F_{i-1}$}]{};
 	\draw[black,ultra thick][-] (whitei) to (blacki);
 	\draw[black,ultra thick][-] (whitei1) to (blacki);
 	\draw[black,ultra thick][-] (whitei2) to (blacki2);
 	\draw[black,ultra thick][-] (whitei2) to (blacki1);
 	\draw[black,ultra thick][-] (whitei) to (blacki1);
 	\draw[black,ultra thick][-] (whitei1) to (blacki2);
 \end{scope}
\begin{scope}[xshift=2.55cm]
	[place/.style={circle,draw=black,thick, inner sep=0pt, minimum size=2mm},
	transition/.style={circle,draw=black,fill=black, inner sep=0pt, minimum size=2mm}]
	\node (whitei) at (150:1.5cm) [place] [label={}] {};
	\node (blacki) at (210:1.5cm) [transition] [label={}]{};
	\node (whitei1) at (270:1.5cm) [place] [label={}]{};
	\node (blacki2) at (330:1.5cm) [transition] [label={right:$\scriptstyle{(i-4)'}$}]{};
	\node (whitei2) at (30:1.5cm) [place] [label={ right:$\scriptstyle{i-4}$}]{};
	\node (blacki1) at (90:1.5cm) [transition] [label={}]{};
	\node at (0,0) [circle, draw=red, fill=red, inner sep=0pt, minimum size=2mm] [label={below:$F_{i-2}$}]{};
	\draw[black,ultra thick][-] (whitei) to (blacki);
	\draw[black,ultra thick][-] (whitei1) to (blacki);
	\draw[black,ultra thick][-] (whitei2) to (blacki2);
	\draw[black,ultra thick][-] (whitei2) to (blacki1);
	\draw[black,ultra thick][-] (whitei) to (blacki1);
	\draw[black,ultra thick][-] (whitei1) to (blacki2);
\end{scope}
\begin{scope}[xshift=-2.56cm]
	[place/.style={circle,draw=black,thick, inner sep=0pt, minimum size=2mm},
	transition/.style={circle,draw=black,fill=black, inner sep=0pt, minimum size=2mm}]
	\node (whitei) at (150:1.5cm) [place] [label={left:$\scriptstyle{i+2}$}] {};
	\node (blacki) at (210:1.5cm) [transition] [label={left:$\scriptstyle{(i+2)'}$}]{};
	\node (whitei1) at (270:1.5cm) [place] [label={}]{};
	\node (blacki2) at (330:1.5cm) [transition] [label={}]{};
	\node (whitei2) at (30:1.5cm) [place] [label={}]{};
	\node (blacki1) at (90:1.5cm) [transition] [label={}]{};
	\node (Fi+2) at (0,0) [circle, draw=red, fill=red, inner sep=0pt, minimum size=2mm] [label={below:$F_{i+2}$}]{};
	\draw[black,ultra thick][-] (whitei) to (blacki);
	\draw[black,ultra thick][-] (whitei1) to (blacki);
	\draw[black,ultra thick][-] (whitei2) to (blacki2);
	\draw[black,ultra thick][-] (whitei2) to (blacki1);
	\draw[black,ultra thick][-] (whitei) to (blacki1);
	\draw[black,ultra thick][-] (whitei1) to (blacki2);
\end{scope}
\begin{scope}[xshift=-1.298cm, yshift=2.25cm]
	[place/.style={circle,draw=black,thick, inner sep=0pt, minimum size=2mm},
	transition/.style={circle,draw=black,fill=black, inner sep=0pt, minimum size=2mm}]
	\node (whitei) at (150:1.5cm) [place] [label={above left:$\scriptstyle{i+1}$}] {};
	\node (blacki) at (210:1.5cm) [transition] [label={above left:$\scriptstyle{(i+1)'}$}]{};
	\node (whitei1) at (270:1.5cm) [place] [label={}]{};
	\node (blacki2) at (330:1.5cm) [transition] [label={}]{};
	\node (whitei2) at (30:1.5cm) [place] [label={}]{};
	\node (blacki1) at (90:1.5cm) [transition] [label={above:$\scriptstyle{i'}$}]{};
	\node at (0,0) [circle, draw=red, fill=red, inner sep=0pt, minimum size=2mm] [label={above:$F_{i+1}$}]{};
	\draw[black,ultra thick][-] (whitei) to (blacki);
	\draw[black,ultra thick][-] (whitei1) to (blacki);
	\draw[black,ultra thick][-] (whitei2) to (blacki2);
	\draw[black,ultra thick][-] (whitei2) to (blacki1);
	\draw[black,ultra thick][-] (whitei) to (blacki1);
	\draw[black,ultra thick][-] (whitei1) to (blacki2);	  		
\end{scope}
\begin{scope}[ultra thick]
	\draw[blue][->, >= {Stealth}] (Fi) to (Fi-1);	
		\draw[blue][->,>= {Stealth}] (Fi) to (Fi-1s);
		\draw[blue][->,>= {Stealth}] (Fi) to (Fi+2);
	\end{scope}
\begin{scope}[xshift=-1.298cm, yshift=-2.25cm]
	[place/.style={circle,draw=black,thick, inner sep=0pt, minimum size=2mm},
	transition/.style={circle,draw=black,fill=black, inner sep=0pt, minimum size=2mm}]
	\node (whitei) at (150:1.5cm) [place] [label={ left:$\scriptstyle{i+1}$}] {};
	\node (blacki) at (210:1.5cm) [transition] [label={ left:$\scriptstyle{(i+1)'}$}]{};
	\node (whitei1) at (270:1.5cm) [place] [label={below:$\scriptstyle{i}$}]{};
	\node (blacki2) at (330:1.5cm) [transition] [label={}]{};
	\node (whitei2) at (30:1.5cm) [place] [label={}]{};
	\node (blacki1) at (90:1.5cm) [transition] [label={}]{};
	\node at (0,0) [circle, draw=red, fill=red, inner sep=0pt, minimum size=2mm] [label={below:$F_{i+1}$}]{};
	\draw[black,ultra thick][-] (whitei) to (blacki);
	\draw[black,ultra thick][-] (whitei1) to (blacki);
	\draw[black,ultra thick][-] (whitei2) to (blacki2);
	\draw[black,ultra thick][-] (whitei2) to (blacki1);
	\draw[black,ultra thick][-] (whitei) to (blacki1);
	\draw[black,ultra thick][-] (whitei1) to (blacki2);	  		
\end{scope} 	
 \end{tikzpicture} 
 $$
 All edges of hexagons, the black colored edges of the drawing, are labeled by the corresponding homomorphism of the representation $\widehat{\rho}(F,G)$. We use the labels 
 $$
 \alpha^{t,s}:P^s \longrightarrow P^{t}
 $$ 
 for edges of $PG_l$ and keep in mind that we also have
 $$
 \alpha^{l-1,0}: P^0 \longrightarrow P^{l-1}, \,\,\, \alpha^{0,l-1}: P^{l-1} \longrightarrow P^{0},
 $$
 given by $\tau^-_{\rho}(F)$ and $\tau^+_{\rho}(G)$ respectively.
 
 \vspace{0.2cm}
 We begin by labeling the arrows $(e^+_{i-2})^{\check{}}$ and $(e^-_{i-1})^{\check{}}$; those should be two homomorphisms 
 $$
 \begin{gathered}
(\widehat{\rho} (F,G))^{\vee}\left((e^+_{i-2})^{\check{}}\, \right): P^{i,+}=Hom(P^i,P^{i-2}) \longrightarrow Hom(P^{i-1},P^{i-3})=P^{i-1,+},
\\
(\widehat{\rho} (F,G))^{\vee}\left((e^-_{i-1})^{\check{}}\, \right): P^{i,+}=Hom(P^i,P^{i-2}) \longrightarrow Hom(P^{i-1},P^{i-3})=P^{i-1,+}.
\end{gathered}
 $$
They are defined by the formulas
\begin{equation}\label{dualarrows-formula}
	\begin{gathered}
{\scriptstyle{(\widehat{\rho} (F,G))^{\vee}\left((e^+_{i-2})^{\check{}}\, \right)}}: Hom(P^i,P^{i-2})\ni \theta \mapsto \alpha^{i-3,i-2}\theta(\alpha^{i-1,i})^{\dag} \in Hom(P^{i-1},P^{i-3}),
\\
{\scriptstyle{(\widehat{\rho} (F,G))^{\vee}\left((e^-_{i-1})^{\check{}}\, \right)}}: Hom(P^i,P^{i-2})\ni \theta \mapsto (\alpha^{i-2,i-3})^{\dag}\theta(\alpha^{i,i-1}) \in Hom(P^{i-1},P^{i-3}),
\end{gathered}
\end{equation}
where $(\cdot)^{\dag}$ stands for the adjoint. 

\vspace{0.2cm}
Similarly, we can label the edges $(e^+_{i-2})^{\check{}}$ and $(e^-_{i-1})^{\check{}}$ with orientation {\it reversed} by the maps:
$$
 \begin{gathered}
	(\widehat{\rho} (F,G))^{\vee}\left(-(e^+_{i-2})^{\check{}}\, \right): P^{i-1,+}=Hom(P^{i-1},P^{i-3}) \longrightarrow Hom(P^i,P^{i-2})=P^{i,+},
	\\
	(\widehat{\rho} (F,G))^{\vee}\left(-(e^-_{i-1})^{\check{}}\, \right): P^{i-1,+}= Hom(P^{i-1},P^{i-3})\longrightarrow Hom(P^i,P^{i-2})=P^{i,+}.
\end{gathered}
$$
defined respectively by the formulas
\begin{equation}\label{dualarrows-reversed-formula}
	\begin{gathered}
		{\scriptstyle{(\widehat{\rho} (F,G))^{\vee}\left(-(e^+_{i-2})^{\check{}}\, \right)}}: Hom(P^{i-1},P^{i-3})\ni \eta \mapsto (\alpha^{i-3,i-2})^{\dag} \eta \alpha^{i-1,i} \in Hom(P^{i},P^{i-2}),
		\\
		{\scriptstyle{(\widehat{\rho} (F,G))^{\vee}\left(-(e^-_{i-1})^{\check{}}\, \right)}}: Hom(P^{i-1},P^{i-3})\ni \eta \mapsto \alpha^{i-2,i-3}\eta(\alpha^{i,i-1})^{\dag} \in Hom(P^{i},P^{i-2}),
	\end{gathered}
\end{equation}

\vspace{0.2cm}
We now have labels for all $\pm$-colored edges of $Q_l$ in 
{\it both} directions. To label the
$0$-colored edges we use the fact that they enter the boundary of
adjacent faces of $Q_l$: the edge $(e^0_i)^{\check{}}$ is the edge of adjacency of the faces $\check{i}$ and  $\check{i'}$; this is depicted in the next drawing below:
  
$$
\begin{tikzpicture}
	[place/.style={circle,draw=black,thick, inner sep=0pt, minimum size=2mm},
	transition/.style={circle,draw=black,fill=black, inner sep=0pt, minimum size=2mm}]
	\node (whitei) at (150:1.5cm) [place] [label={left:$\scriptstyle{i}$}] {};
	\node (blacki) at (210:1.5cm) [transition] [label={left:$\scriptstyle{i'}$}]{};
	\node (whitei1) at (270:1.5cm) [place] [label={below left:$\scriptstyle{i-1}$}]{};
	\node (blacki2) at (330:1.5cm) [transition] [label={below:$\scriptstyle{(i-2)'}$}]{};
	\node (whitei2) at (30:1.5cm) [place] [label={ above:$\scriptstyle{i-2}$}]{};
	\node (blacki1) at (90:1.5cm) [transition] [label={}]{};
	\node (Fi) at (0,0) [circle, draw=red, fill=red, inner sep=0pt, minimum size=2mm] [label={below:$F_i$}]{};
	\draw[black,ultra thick][-] (whitei) to (blacki);
	\draw[black,ultra thick][-] (whitei1) to (blacki);
	\draw[black,ultra thick][-] (whitei2) to (blacki2);
	\draw[black,ultra thick][-] (whitei2) to (blacki1);
	\draw[black,ultra thick][-] (whitei) to (blacki1);
	\draw[black,ultra thick][-] (whitei1) to (blacki2);
	\begin{scope}[xshift=1.26cm, yshift=2.25cm]
		[place/.style={circle,draw=black,thick, inner sep=0pt, minimum size=2mm},
		transition/.style={circle,draw=black,fill=black, inner sep=0pt, minimum size=2mm}]
		\node (whitei) at (150:1.5cm) [place] [label={above:$\scriptstyle{i-1}$}] {};
		\node (blacki) at (210:1.5cm) [transition] [label={ below:$\scriptstyle{(i-1)'}$}]{};
		\node (whitei1) at (270:1.5cm) [place] [label={}]{};
		\node (blacki2) at (330:1.5cm) [transition] [label={right:$\scriptstyle{(i-3)'}$}]{};
		\node (whitei2) at (30:1.5cm) [place] [label={ right:$\scriptstyle{i-3}$}]{};
		\node (blacki1) at (90:1.5cm) [transition] [label={above:$\scriptstyle{(i-2)'}$}]{};
		\node (Fi-1) at (0,0) [circle, draw=red, fill=red, inner sep=0pt, minimum size=2mm] [label={above:$F_{i-1}$}]{};
		\draw[black,ultra thick][-] (whitei) to (blacki);
		\draw[black,ultra thick][-] (whitei1) to (blacki);
		\draw[black,ultra thick][-] (whitei2) to (blacki2);
		\draw[black,ultra thick][-] (whitei2) to (blacki1);
		\draw[black,ultra thick][-] (whitei) to (blacki1);
		\draw[black,ultra thick][-] (whitei1) to (blacki2);
	\end{scope}
	\begin{scope}[xshift=1.26cm, yshift=-2.25cm]
		[place/.style={circle,draw=black,thick, inner sep=0pt, minimum size=2mm},
		transition/.style={circle,draw=black,fill=black, inner sep=0pt, minimum size=2mm}]
		\node (whitei) at (150:1.5cm) [place] [label={}] {};
		\node (blacki) at (210:1.5cm) [transition] [label={below :$\scriptstyle{(i-1)'}$}]{};
		\node (whitei1) at (270:1.5cm) [place] [label={below:$\scriptstyle{i-2}$}]{};
		\node (blacki2) at (330:1.5cm) [transition] [label={right:$\scriptstyle{(i-3)'}$}]{};
		\node (whitei2) at (30:1.5cm) [place] [label={ right:$\scriptstyle{i-3}$}]{};
		\node (blacki1) at (90:1.5cm) [transition] [label={}]{};
		\node (Fi-1s)at (0,0) [circle, draw=red, fill=red, inner sep=0pt, minimum size=2mm] [label={below:$F_{i-1}$}]{};
		\draw[black,ultra thick][-] (whitei) to (blacki);
		\draw[black,ultra thick][-] (whitei1) to (blacki);
		\draw[black,ultra thick][-] (whitei2) to (blacki2);
		\draw[black,ultra thick][-] (whitei2) to (blacki1);
		\draw[black,ultra thick][-] (whitei) to (blacki1);
		\draw[black,ultra thick][-] (whitei1) to (blacki2);
	\end{scope}
	\begin{scope}[xshift=2.55cm]
		[place/.style={circle,draw=black,thick, inner sep=0pt, minimum size=2mm},
		transition/.style={circle,draw=black,fill=black, inner sep=0pt, minimum size=2mm}]
		\node (whitei) at (150:1.5cm) [place] [label={}] {};
		\node (blacki) at (210:1.5cm) [transition] [label={}]{};
		\node (whitei1) at (270:1.5cm) [place] [label={}]{};
		\node (blacki2) at (330:1.5cm) [transition] [label={right:$\scriptstyle{(i-4)'}$}]{};
		\node (whitei2) at (30:1.5cm) [place] [label={ right:$\scriptstyle{i-4}$}]{};
		\node (blacki1) at (90:1.5cm) [transition] [label={}]{};
		\node at (0,0) [circle, draw=red, fill=red, inner sep=0pt, minimum size=2mm] [label={below:$F_{i-2}$}]{};
		\draw[black,ultra thick][-] (whitei) to (blacki);
		\draw[black,ultra thick][-] (whitei1) to (blacki);
		\draw[black,ultra thick][-] (whitei2) to (blacki2);
		\draw[black,ultra thick][-] (whitei2) to (blacki1);
		\draw[black,ultra thick][-] (whitei) to (blacki1);
		\draw[black,ultra thick][-] (whitei1) to (blacki2);
	\end{scope}
	\begin{scope}[xshift=-2.56cm]
		[place/.style={circle,draw=black,thick, inner sep=0pt, minimum size=2mm},
		transition/.style={circle,draw=black,fill=black, inner sep=0pt, minimum size=2mm}]
		\node (whitei) at (150:1.5cm) [place] [label={left:$\scriptstyle{i+2}$}] {};
		\node (blacki) at (210:1.5cm) [transition] [label={left:$\scriptstyle{(i+2)'}$}]{};
		\node (whitei1) at (270:1.5cm) [place] [label={}]{};
		\node (blacki2) at (330:1.5cm) [transition] [label={}]{};
		\node (whitei2) at (30:1.5cm) [place] [label={}]{};
		\node (blacki1) at (90:1.5cm) [transition] [label={}]{};
		\node (Fi+2) at (0,0) [circle, draw=red, fill=red, inner sep=0pt, minimum size=2mm] [label={below:$F_{i+2}$}]{};
		\draw[black,ultra thick][-] (whitei) to (blacki);
		\draw[black,ultra thick][-] (whitei1) to (blacki);
		\draw[black,ultra thick][-] (whitei2) to (blacki2);
		\draw[black,ultra thick][-] (whitei2) to (blacki1);
		\draw[black,ultra thick][-] (whitei) to (blacki1);
		\draw[black,ultra thick][-] (whitei1) to (blacki2);
	\end{scope}
	\begin{scope}[xshift=-1.298cm, yshift=2.25cm]
		[place/.style={circle,draw=black,thick, inner sep=0pt, minimum size=2mm},
		transition/.style={circle,draw=black,fill=black, inner sep=0pt, minimum size=2mm}]
		\node (whitei) at (150:1.5cm) [place] [label={above left:$\scriptstyle{i+1}$}] {};
		\node (blacki) at (210:1.5cm) [transition] [label={above left:$\scriptstyle{(i+1)'}$}]{};
		\node (whitei1) at (270:1.5cm) [place] [label={}]{};
		\node (blacki2) at (330:1.5cm) [transition] [label={}]{};
		\node (whitei2) at (30:1.5cm) [place] [label={}]{};
		\node (blacki1) at (90:1.5cm) [transition] [label={above:$\scriptstyle{i'}$}]{};
		\node (Fi+1) at (0,0) [circle, draw=red, fill=red, inner sep=0pt, minimum size=2mm] [label={above:$F_{i+1}$}]{};
		\draw[black,ultra thick][-] (whitei) to (blacki);
		\draw[black,ultra thick][-] (whitei1) to (blacki);
		\draw[black,ultra thick][-] (whitei2) to (blacki2);
		\draw[black,ultra thick][-] (whitei2) to (blacki1);
		\draw[black,ultra thick][-] (whitei) to (blacki1);
		\draw[black,ultra thick][-] (whitei1) to (blacki2);	  		
	\end{scope}
	\begin{scope}[xshift=-1.298cm, yshift=-2.25cm]
		[place/.style={circle,draw=black,thick, inner sep=0pt, minimum size=2mm},
		transition/.style={circle,draw=black,fill=black, inner sep=0pt, minimum size=2mm}]
		\node (whitei) at (150:1.5cm) [place] [label={ left:$\scriptstyle{i+1}$}] {};
		\node (blacki) at (210:1.5cm) [transition] [label={ left:$\scriptstyle{(i+1)'}$}]{};
		\node (whitei1) at (270:1.5cm) [place] [label={below:$\scriptstyle{i}$}]{};
		\node (blacki2) at (330:1.5cm) [transition] [label={}]{};
		\node (whitei2) at (30:1.5cm) [place] [label={}]{};
		\node (blacki1) at (90:1.5cm) [transition] [label={}]{};
		\node (Fi+1s) at (0,0) [circle, draw=red, fill=red, inner sep=0pt, minimum size=2mm] [label={below:$F_{i+1}$}]{};
		\draw[black,ultra thick][-] (whitei) to (blacki);
		\draw[black,ultra thick][-] (whitei1) to (blacki);
		\draw[black,ultra thick][-] (whitei2) to (blacki2);
		\draw[black,ultra thick][-] (whitei2) to (blacki1);
		\draw[black,ultra thick][-] (whitei) to (blacki1);
		\draw[black,ultra thick][-] (whitei1) to (blacki2);	  		
	\end{scope}
\begin{scope}[ultra thick]
	\draw[blue][->, >= {Stealth}] (Fi+2) to (Fi+1);	
	\draw[blue][->,>= {Stealth}] (Fi+1) to (Fi);
	\draw[blue][->,>= {Stealth}] (Fi) to (Fi+2);
		\draw[blue][->,>= {Stealth}] (Fi+2) to (Fi+1s);	
		\draw[blue][->,>= {Stealth}] (Fi+1s) to (Fi);	
\end{scope} 	
\end{tikzpicture} 
$$
Thus we can label $(e^0_i)^{\check{}}$ with the compositions of maps labeling $-(e^-_i)^{\check{}}$ and $-(e^+_i)^{\check{}}$ (resp. $-(e^+_{i-1})^{\check{}}$ and $-(e^-_{i+1})^{\check{}}$). More precisely, we have the bilinear map
$$
\xymatrix{
Hom(P^{i,+}, P^{(i+1),+}) \times Hom(P^{(i+1),+},P^{(i+2),+} )\ar[r]& 
Hom(P^{i,+},P^{(i+2),+})
}
$$
given by composition of maps. In particular, we assign to $(e^0_i)^{\check{}}$ the compositions 
$$
\begin{gathered}
{}^{\circ}(\widehat{\rho})^{\vee}((e^0_i)^{\check{}}\,) =(\widehat{\rho})^{\vee}(-(e^+_i)^{\check{}}\,) \circ (\widehat{\rho})^{\vee}(-(e^-_i)^{\check{}}):Hom(P^i,P^{i-2})\longrightarrow Hom(P^{i+2},P^{i}),
\\
{}^{\bullet}(\widehat{\rho})^{\vee}((e^0_i)^{\check{}}\,) =(\widehat{\rho})^{\vee}(-(e^-_{i+1})^{\check{}}\,) \circ (\widehat{\rho})^{\vee}(-(e^+_{i-1})^{\check{}}):Hom(P^i,P^{i-2})\longrightarrow Hom(P^{i+2},P^{i}).
\end{gathered}
$$
The first one corresponds to circling around the white vertex $i$ and the second around the black vertex $i'$.
Explicitly, we have the formulas
\begin{equation}
	\begin{gathered}
	{\scriptstyle{{}^{\circ}(\widehat{\rho})^{\vee}((e^0_i)^{\check{}}\,)}}: {\scriptstyle{Hom(P^i,P^{i-2})}} \ni \theta \mapsto (\alpha^{i-1,i})^{\dag} \alpha^{i-1,i-2} \theta (\alpha^{i+1,i})^{\dag} \alpha^{i+1,i+2} \in {\scriptstyle{Hom(P^{i+2},P^{i})}},
	\\
	{\scriptstyle{{}^{\bullet}(\widehat{\rho})^{\vee}((e^0_i)^{\check{}}\,)}}: {\scriptstyle Hom(P^i,P^{i-2})} \ni \theta \mapsto (\alpha^{i,i-1}) (\alpha^{i-2,i-1})^{\dag} \theta (\alpha^{i,i+1}) (\alpha^{i+2,i+1})^{\dag} \in {\scriptstyle Hom(P^{i+2},P^{i})}.
	\end{gathered} 
\end{equation}
We provide the paths labeled by the homomorphisms in the above compositions:

- for ${}^{\circ}(\widehat{\rho})^{\vee}((e^0_i)^{\check{}}\,)$: start at the white vertex $(i+2)$, follow the two top edges of the hexagon $F_{i+2}$ to arrive to $i$, from there we jump to the white vertex $(i-2)$ on the opposite side of $F_i$ and we return to $i$ via the two top edges of
that hexagon;

- for ${}^{\bullet}(\widehat{\rho})^{\vee}((e^0_i)^{\check{}}\,)$: start at the black vertex $(i+2)'$, follow the two bottom edges of the hexagon $F_{i+2}$ to arrive to $i'$, from there we jump to the black vertex $(i-2)'$ on the opposite side of $F_i$ and we return to $i$ via the two bottom edges of
that hexagon.

The following drawing shows the paths described above:

$$
\begin{tikzpicture}
	[place/.style={circle,draw=black,thick, inner sep=0pt, minimum size=2mm},
	transition/.style={circle,draw=black,fill=black, inner sep=0pt, minimum size=2mm}]
	\node (whitei) at (150:1.5cm) [place] [label={left:$\scriptstyle{i}$}] {};
	\node (blacki) at (210:1.5cm) [transition] [label={left:$\scriptstyle{i'}$}]{};
	\node (whitei1) at (270:1.5cm) [place] [label={below left:$\scriptstyle{i-1}$}]{};
	\node (blacki2) at (330:1.5cm) [transition] [label={below:$\scriptstyle{(i-2)'}$}]{};
	\node (whitei2) at (30:1.5cm) [place] [label={ above:$\scriptstyle{i-2}$}]{};
	\node (blacki1) at (90:1.5cm) [transition] [label={}]{};
	\node (Fi) at (0,0) [circle, draw=red, fill=red, inner sep=0pt, minimum size=2mm] [label={below:$F_i$}]{};
	\draw[black,ultra thick][-] (whitei) to (blacki);
	\draw[orange,ultra thick][->,>={Stealth}] (whitei1) to (blacki);
	\draw[orange,ultra thick][<-,>={Stealth}] (blacki2) to (blacki);
	\draw[black,ultra thick][-] (whitei2) to (blacki2);
	\draw[green,ultra thick][->,>={Stealth}] (whitei2) to (blacki1);
	\draw[green,ultra thick][<-,>={Stealth}] (whitei) to (blacki1);
	\draw[orange,ultra thick][<-,>={Stealth}] (whitei1) to (blacki2);
	\draw[green,ultra thick][->, >= {Stealth}] (whitei) to (whitei2);
	\begin{scope}[xshift=1.26cm, yshift=2.25cm]
		[place/.style={circle,draw=black,thick, inner sep=0pt, minimum size=2mm},
		transition/.style={circle,draw=black,fill=black, inner sep=0pt, minimum size=2mm}]
		\node (whitei) at (150:1.5cm) [place] [label={above:$\scriptstyle{i-1}$}] {};
		\node (blacki) at (210:1.5cm) [transition] [label={ below:$\scriptstyle{(i-1)'}$}]{};
		\node (whitei1) at (270:1.5cm) [place] [label={}]{};
		\node (blacki2) at (330:1.5cm) [transition] [label={right:$\scriptstyle{(i-3)'}$}]{};
		\node (whitei2) at (30:1.5cm) [place] [label={ right:$\scriptstyle{i-3}$}]{};
		\node (blacki1) at (90:1.5cm) [transition] [label={above:$\scriptstyle{(i-2)'}$}]{};
		\node (Fi-1) at (0,0) [circle, draw=red, fill=red, inner sep=0pt, minimum size=2mm] [label={above:$F_{i-1}$}]{};
		\draw[black,ultra thick][-] (whitei) to (blacki);
		\draw[green,ultra thick][->,>={Stealth}] (whitei1) to (blacki);
		\draw[black,ultra thick][-] (whitei2) to (blacki2);
		\draw[black,ultra thick][-] (whitei2) to (blacki1);
		\draw[black,ultra thick][-] (whitei) to (blacki1);
		\draw[black,ultra thick][-] (whitei1) to (blacki2);
	\end{scope}
	\begin{scope}[xshift=1.26cm, yshift=-2.25cm]
		[place/.style={circle,draw=black,thick, inner sep=0pt, minimum size=2mm},
		transition/.style={circle,draw=black,fill=black, inner sep=0pt, minimum size=2mm}]
		\node (whitei) at (150:1.5cm) [place] [label={}] {};
		\node (blacki) at (210:1.5cm) [transition] [label={below :$\scriptstyle{(i-1)'}$}]{};
		\node (whitei1) at (270:1.5cm) [place] [label={below:$\scriptstyle{i-2}$}]{};
		\node (blacki2) at (330:1.5cm) [transition] [label={right:$\scriptstyle{(i-3)'}$}]{};
		\node (whitei2) at (30:1.5cm) [place] [label={ right:$\scriptstyle{i-3}$}]{};
		\node (blacki1) at (90:1.5cm) [transition] [label={}]{};
		\node (Fi-1s)at (0,0) [circle, draw=red, fill=red, inner sep=0pt, minimum size=2mm] [label={below:$F_{i-1}$}]{};
		\draw[black,ultra thick][-] (whitei) to (blacki);
		\draw[black,ultra thick][-] (whitei1) to (blacki);
		\draw[black,ultra thick][-] (whitei2) to (blacki2);
		\draw[black,ultra thick][-] (whitei2) to (blacki1);
		\draw[orange,ultra thick][<-, >= {Stealth}] (whitei) to (blacki1);
		\draw[black,ultra thick][-] (whitei1) to (blacki2);
	\end{scope}
	\begin{scope}[xshift=2.55cm]
		[place/.style={circle,draw=black,thick, inner sep=0pt, minimum size=2mm},
		transition/.style={circle,draw=black,fill=black, inner sep=0pt, minimum size=2mm}]
		\node (whitei) at (150:1.5cm) [place] [label={}] {};
		\node (blacki) at (210:1.5cm) [transition] [label={}]{};
		\node (whitei1) at (270:1.5cm) [place] [label={}]{};
		\node (blacki2) at (330:1.5cm) [transition] [label={right:$\scriptstyle{(i-4)'}$}]{};
		\node (whitei2) at (30:1.5cm) [place] [label={ right:$\scriptstyle{i-4}$}]{};
		\node (blacki1) at (90:1.5cm) [transition] [label={}]{};
		\node at (0,0) [circle, draw=red, fill=red, inner sep=0pt, minimum size=2mm] [label={below:$F_{i-2}$}]{};
		\draw[black,ultra thick][-] (whitei) to (blacki);
		\draw[black,ultra thick][-] (whitei1) to (blacki);
		\draw[black,ultra thick][-] (whitei2) to (blacki2);
		\draw[black,ultra thick][-] (whitei2) to (blacki1);
		\draw[black,ultra thick][-] (whitei) to (blacki1);
		\draw[black,ultra thick][-] (whitei1) to (blacki2);
	\end{scope}
	\begin{scope}[xshift=-2.56cm]
		[place/.style={circle,draw=black,thick, inner sep=0pt, minimum size=2mm},
		transition/.style={circle,draw=black,fill=black, inner sep=0pt, minimum size=2mm}]
		\node (whitei) at (150:1.5cm) [place] [label={left:$\scriptstyle{i+2}$}] {};
		\node (blacki) at (210:1.5cm) [transition] [label={left:$\scriptstyle{(i+2)'}$}]{};
		\node (whitei1) at (270:1.5cm) [place] [label={}]{};
		\node (blacki2) at (330:1.5cm) [transition] [label={}]{};
		\node (whitei2) at (30:1.5cm) [place] [label={}]{};
		\node (blacki1) at (90:1.5cm) [transition] [label={}]{};
		\node (Fi+2) at (0,0) [circle, draw=red, fill=red, inner sep=0pt, minimum size=2mm] [label={below:$F_{i+2}$}]{};
		\draw[black,ultra thick][-] (whitei) to (blacki);
		\draw[orange,ultra thick][<-, >= {Stealth}] (whitei1) to (blacki);
		\draw[black,ultra thick][-] (whitei2) to (blacki2);
		\draw[green,ultra thick][<-,>={Stealth}] (whitei2) to (blacki1);
		\draw[green,ultra thick][->,>={Stealth}] (whitei) to (blacki1);
		\draw[orange,ultra thick][->, >= {Stealth}] (whitei1) to (blacki2);
	\end{scope}
	\begin{scope}[xshift=-1.298cm, yshift=2.25cm]
		[place/.style={circle,draw=black,thick, inner sep=0pt, minimum size=2mm},
		transition/.style={circle,draw=black,fill=black, inner sep=0pt, minimum size=2mm}]
		\node (whitei) at (150:1.5cm) [place] [label={above left:$\scriptstyle{i+1}$}] {};
		\node (blacki) at (210:1.5cm) [transition] [label={above left:$\scriptstyle{(i+1)'}$}]{};
		\node (whitei1) at (270:1.5cm) [place] [label={}]{};
		\node (blacki2) at (330:1.5cm) [transition] [label={}]{};
		\node (whitei2) at (30:1.5cm) [place] [label={}]{};
		\node (blacki1) at (90:1.5cm) [transition] [label={above:$\scriptstyle{i'}$}]{};
		\node (Fi+1) at (0,0) [circle, draw=red, fill=red, inner sep=0pt, minimum size=2mm] [label={above:$F_{i+1}$}]{};
		\draw[black,ultra thick][-] (whitei) to (blacki);
		\draw[green,ultra thick][<-,>={Stealth}] (whitei1) to (blacki);
		\draw[black,ultra thick][-] (whitei2) to (blacki2);
		\draw[black,ultra thick][-] (whitei2) to (blacki1);
		\draw[black,ultra thick][-] (whitei) to (blacki1);
		\draw[green,ultra thick][<-,>={Stealth}] (whitei1) to (blacki2);	  		
	\end{scope}
	\begin{scope}[xshift=-1.298cm, yshift=-2.25cm]
		[place/.style={circle,draw=black,thick, inner sep=0pt, minimum size=2mm},
		transition/.style={circle,draw=black,fill=black, inner sep=0pt, minimum size=2mm}]
		\node (whitei) at (150:1.5cm) [place] [label={ left:$\scriptstyle{i+1}$}] {};
		\node (blacki) at (210:1.5cm) [transition] [label={ left:$\scriptstyle{(i+1)'}$}]{};
		\node (whitei1) at (270:1.5cm) [place] [label={below:$\scriptstyle{i}$}]{};
		\node (blacki2) at (330:1.5cm) [transition] [label={}]{};
		\node (whitei2) at (30:1.5cm) [place] [label={}]{};
		\node (blacki1) at (90:1.5cm) [transition] [label={}]{};
		\node (Fi+1s) at (0,0) [circle, draw=red, fill=red, inner sep=0pt, minimum size=2mm] [label={below:$F_{i+1}$}]{};
		\draw[black,ultra thick][-] (whitei) to (blacki);
		\draw[black,ultra thick][-] (whitei1) to (blacki);
		\draw[black,ultra thick][-] (whitei2) to (blacki2);
		\draw[orange,ultra thick][->, >= {Stealth}] (whitei2) to (blacki1);
		\draw[orange,ultra thick][->, >= {Stealth}] (whitei) to (blacki1);
		\draw[black,ultra thick][-] (whitei1) to (blacki2);	  		
	\end{scope}
	\begin{scope}[ultra thick]
		\draw[blue][->, >= {Stealth}] (Fi+2) to (Fi+1);	
		\draw[blue][->,>= {Stealth}] (Fi+1) to (Fi);
		\draw[blue][->,>= {Stealth}] (Fi) to (Fi+2);
		\draw[blue][->,>= {Stealth}] (Fi+2) to (Fi+1s);	
		\draw[blue][->,>= {Stealth}] (Fi+1s) to (Fi);	
	\end{scope} 	
\end{tikzpicture} 
$$
on the drawing the green (resp. orange) path represents the composition
${}^{\circ}(\widehat{\rho})^{\vee}((e^0_i)^{\check{}}\,)$ (resp. ${}^{\bullet}(\widehat{\rho})^{\vee}((e^0_i)^{\check{}}\,)$); the `jumps' are the horizontal green and orange arrows and they are labeled by $\theta$ in the formulas.

The labeling of edges of $Q_l$ is now completed. The procedure described above tells us that the representations of $Q_l$ obtained from $\widehat{\rho}$ depend only on choices of labels for the $0$-colored edges: for each edge $(e^0_i)^{\check{}}$ we have two possibilities
$$
\xymatrix@C=10pt@R=10pt{
	&{}^{\circ}(\widehat{\rho})^{\vee}((e^0_i)^{\check{}}\,)\\
	(e^0_i)^{\check{}} \ar[ru]\ar[rd]&\\
	&{}^{\bullet}(\widehat{\rho})^{\vee}((e^0_i)^{\check{}}\,)
}
$$
For a subset $A$ of $[0,l-1]$ denote by $(\widehat{\rho})^{\vee}_A$ the representation of $Q_l$ described above which puts labels
${}^{\circ}(\widehat{\rho})^{\vee}((e^0_i)^{\check{}}\,)$, for $i\in A$, and the label ${}^{\bullet}(\widehat{\rho})^{\vee}((e^0_j)^{\check{}}\,)$, for $j \in A^c$, in the complement $A^c$ of $A$ in $[0,l-1]$:
\begin{equation}\label{repQl-A}
	(\widehat{\rho})^{\vee}_A((e^0_i)^{\check{}}\,)=\begin{cases}
		{}^{\circ}(\widehat{\rho})^{\vee}((e^0_i)^{\check{}}\,),&i\in A,\\
		{}^{\bullet}(\widehat{\rho})^{\vee}((e^0_i)^{\check{}}\,),&i\in A^c.
	\end{cases}
\end{equation}
This gives the collection of representations of $Q_l$
$$
Reps(Q_l)(\widehat{\rho}):=\left\{(\widehat{\rho})^{\vee}_A \Big|\,\, A \subset [0,l-1] \right\}
$$
asserted in the proposition.
\end{pf}
 
 Let us recall that to every stratum $\mathfrak{L}_l$ are attached families
 of representations
 $$
 \rho_c: \OO^{\times}_{\mathfrak{L}_l}(-1)\longrightarrow \mathfrak{Reps}(PG_l)
 $$
 where $c$ is the trace parameter. The reader may recall that $\rho_c$ assigns to a vector $\xi\otimes\phi$ in the fibre of $\OO^{\times}_{\mathfrak{L}_l}(-1)$ lying over a point $([\xi],[\phi])\in {\mathfrak{L}_l}$ the representation
 $$
 \rho_c(\xi\otimes\phi)=\left\{\{P^s([\xi],[\phi])\}, \{\alpha^{t,s}_c(\xi,\phi):P^s([\xi],[\phi])\longrightarrow P^{t}([\xi],[\phi])\}\right\}.
 $$
 In Corollary \ref{cor:rhoc-ext} we learned  how to extend those to the families of representations
 $$
 \rho_c(F,G): \OO^{\times}_{\mathfrak{L}_l}(-1)\longrightarrow \mathfrak{Reps}(\widehat{PG}_l)
 $$
 of the quiver $\widehat{PG}_l$, for every pair of linear functionals $F$ and $G$ on $\CC[q,q^{-1}]$. Specializing Proposition \ref{pro:repsQl} to the family
$\rho_c(F,G)$ we obtain the families of representations of the dual quiver
$Q_l$.
\begin{cor}
	Let $c$ be a value of the trace parameter and let
	 $$
	\rho_c: \OO^{\times}_{\mathfrak{L}_l}(-1)\longrightarrow \mathfrak{Reps}(PG_l)
	$$
	be the corresponding family of the representations of $PG_l$. A pair of linear functionals $(F,G)$ on the space of Laurent polynomials
	$\CC[q,q^{-1}]$ determines the family of representations of $\widehat{PG}_l$
	$$
	\widehat{\rho_c}(F,G):\OO^{\times}_{\mathfrak{L}_l}(-1)\longrightarrow \mathfrak{Reps}(\widehat{PG}_l),
	$$
	see Corollary \ref{cor:rhoc-ext}. This in turn gives rise to the finite collection $$
	\{(\widehat{\rho_c}(F,G))^{\vee}_A \big| A \subset [0,l-1] \}
	$$
	of families of representations of the dual quiver $Q_l$: the families
	are labeled by subsets $A$ of $[0,l-1]$, the set of pairs of black-white vertices of
	$\widehat{PG}_l$, and each family is the map
	$$
	(\widehat{\rho_c}(F,G))^{\vee}_A: \OO^{\times}_{\mathfrak{L}_l}(-1)\longrightarrow \mathfrak{Reps}(Q_l),
	$$
	which assigns to every nonzero vector $\xi\otimes\phi$ in $\OO_{\mathfrak{L}_l}(-1)$ the representation
	$$
		(\widehat{\rho_c}(F,G)(\xi\otimes\phi))^{\vee}_A \in {\mathfrak{Reps}(Q_l)}
		$$
		associated to the representation $\widehat{\rho_c}(F,G)(\xi\otimes\phi)$ of $\widehat{PG}_l$ and the subset $A\subset [0,l-1]$ according to the procedure of the proof of Proposition \ref{pro:repsQl}.
\end{cor} 
   
 	The ideas, in part related to TQFT, led us to uncover an interesting feature of the  representations of the quiver $PG_l$: they can be extended to
 	the representations of the completed graph $\widehat{PG}_l$. This is done through linear maps
 	from the space of linear functionals on $\CC[q,q^{-1}]$ to the space of homomorphisms between $P^0$ and $P^{l-1}$, the ends of the representations of $PG_l$. 
 	
 	The representations of the dual quiver described in Proposition \ref{pro:repsQl} invokes cutting the topological torus into pieces - here the pieces correspond to `tiles' of the cell decomposition of $\mathbb{T}$. More generally, the issue seems to concern additional structures of the ring of endomorphisms
 	$$
 	End_{\CC} (W_{\rho})
 	$$
 	where  $\rho=\{\{P^s\},\{\alpha^{t,s}:P^s \longrightarrow P^t\}\}$ is a bipartite representation of the quiver $PG_l$ and $W_{\rho}$ is the graded space of the representation
 	$$
 	W_{\rho}=\bigoplus^{l-1}_{s=0} P^s.
 	$$
 Extending a representation $\rho$ to $\widehat{\rho}$ of the completed quiver $\widehat{PG}_l$ allows to view the space $W_{\rho}$ as a vector space graded by the cyclic group $\ZZ/l\ZZ$. In particular, the two pieces $P^0$ and $P^{l-1}$ of the decomposition become {\it adjacent} and 
 our constructions produce the maps of degree $(\pm1)$ between them. In the next subsection we show how to shift the construction to produce the operators of degree $(\pm1)$ on the whole graded space
  	$$
  W_{\rho}=\bigoplus_{s\in \ZZ/l\ZZ} P^s.
  $$
  In addition, we will see that there are `higher' order multiplications related to the perfect matchings of $\widehat{PG}_l$; here the perfect matchings provide pattern to cut-and-paste the torus into pieces.
  
 	\subsection{Shifted quantum-type invariants and perfect matchings}
{\bf 11.4.1. Shifted quantum-type invariants.} The constructions of the previous subsections can be summarized by the following diagram
 $$
 \begin{tikzpicture}
 	[place/.style={circle,draw=black,thick},
 	transition/.style={circle,draw=black,fill=black}]
 	\begin{scope}[scale=0.7, shift={(-2,-0.5)}]
 	\node (white3) at (-4,2) [place] [label={above:$l-1$}] {};
 	\node (black3) at (-4,0) [transition] [label={below:$(l-1)'$}]{};
 	\node (white0) at (0,2) [place] [label={above:$0$}] {};
 	\node (black0) at (0,0) [transition] [label={below:$0'$}]{};
 	\node (PG) at (-2,1.5) [label={below:$PG_l$}]{};
 	\draw[black, dotted]
 	(black3) to (black0);
 	\draw[black, dotted]
 	(white3) to (white0);
 	\draw[black, very thick]
 	(white3) to (black3);
 	\draw[black, thick]
 	(white0) to (black0);
 		\draw[black, dotted]
 	(black3) to (black0);
 	\draw[black, dotted]
 	(white3) to (white0);
 	\draw[black, thick]
 	(white3) to (black3);
 	\draw[black, very thick]
 	(white0) to (black0);
 	\end{scope}
 	\begin{scope} [scale=0.7, shift={(6,0)}]
 		\node (white3) at (3,1) [place] [label={below right:$l-1$}] {};
 		\node (black3) at (2,2) [transition] [label={above:$(l-1)'$}]{};
 		\node (white0) at (5,1) [place] [label={below right:$0$}]{};
 		\node (black0) at (6,2) [transition] [label={above right:$0'$}]{};
 		\node (hatPG)  at (4,-2.7) [label={$\widehat{PG}_l$}]{};
 		\draw[black, thick]
 		(white0) to (black0);
 	\draw[black, thick]
 		(white3) to (black3);
 		\draw[black,dotted]
 		(4,0) circle [radius=1.4];
 		\draw[black,dotted]
 		(4,0) circle [radius=2.8];
 		\draw [red, densely dotted]
 		(white3) to [bend left=45] (black0);
 		\draw [blue, thick]
 		(black3) to [bend left=45] (white0);
 		\end{scope}
 	\draw[->,decorate,decoration={snake, amplitude=.4mm, segment length=2mm, post length=1mm}, >={Stealth}]
 	(0.5,0.5)--(3.5,0.5)
 	node[above,align=center,midway]
 	{\textcolor{red} {Quantum invariants}};
 \end{tikzpicture}
$$
On the left side of the diagram is the `open' graph $PG_l$ thought as a quiver together with the category of its bipartite representations: we remind the reader that this means that each pair of vertically aligned black-white vertices of $PG_l$ is labeled by the same vector space; all representations considered here are bipartite and finite dimensional, so we will omit these specifications. 

The arrow `Quantum invariants' turns $PG_l$ into the cyclic graph $\widehat{PG}_l$ by putting the vertices on two concentric circles, `equatorial' circles of the torus ${\mathbb{T}}$, and joining the end white vertex $(l-1)$ (resp. $0$) to the black one $(0)'$, (resp. $(l-1)'$), the dotted red (resp. solid blue) edge of the diagram. The quantum invariants give a functorial procedure of extending a representation of $PG_l$ to the one of $\widehat{PG}_l$. Going in the other direction amounts to cutting the added edges by a radial segment and erasing the half edges at each end; on the level of representations it comes down to the restricting a representation of 
$\widehat{PG}_l$ to the one of $PG_l$ by forgetting the maps labeling the erased edges. But in $\widehat{PG}_l$ all neighboring pairs of `black-white' vertices are on equal footing, so the restriction/forgetting can be performed at any level:
 $$
\begin{tikzpicture}
	[place/.style={circle,draw=black,thick},
	transition/.style={circle,draw=black,fill=black}]
	\begin{scope}[scale=0.7, shift={(-2,0.5)}]
		\node (white3) at (-4,2) [place] [label={above:$s$}] {};
		\node (black3) at (-4,0) [transition] [label={below:$(s)'$}]{};
		\node (white0) at (0,2) [place] [label={above:$s-1$}] {};
		\node (black0) at (0,0) [transition] [label={below:$(s-1)'$}]{};
		\node (PG) at (-2,1.5) [label={below:$PG_l(s)$}]{};
		\draw[black, dotted]
		(black3) to (black0);
		\draw[black, dotted]
		(white3) to (white0);
		\draw[black, very thick]
		(white3) to (black3);
		\draw[black, thick]
		(white0) to (black0);
		\draw[black, dotted]
		(black3) to (black0);
		\draw[black, dotted]
		(white3) to (white0);
		\draw[black, thick]
		(white3) to (black3);
		\draw[black, very thick]
		(white0) to (black0);
	\end{scope}
	\begin{scope} [scale=0.7, shift={(6,0)}, rotate=30]
		\node (white3) at (3,1) [place] [label={below right:$s-1$}] {};
		\node (black3) at (2,2) [transition] [label={above left:$(s-1)'$}]{};
		\node (white0) at (5,1) [place] [label={below right:$s$}]{};
		\node (black0) at (6,2) [transition] [label={above right:$(s)'$}]{};
		\node (hatPG)  at (4,-2.7) [label={$\widehat{PG}_l$}]{};
		\draw[black, thick]
		(white0) to (black0);
		\draw[black, thick]
		(white3) to (black3);
		\draw[green, thick, dashed]
		(4,1)--(4,3.5);
		\draw[black,dotted]
		(4,0) circle [radius=1.4];
		\draw[black,dotted]
		(4,0) circle [radius=2.8];
		\draw [red, densely dotted]
		(white3) to [bend left=45] (black0);
		\draw [blue, thick]
		(black3) to [bend left=45] (white0);
	\end{scope}
\begin{scope}[yshift=10]
	\draw[->,decorate,decoration={snake, amplitude=.4mm, segment length=2mm, post length=1mm}, >={Stealth}]
	(4,0.5)--(-0.5,0.5)
	node[above,align=center,midway]
	{\textcolor{red} {Cutting and Forgetting}};
	\end{scope}
\end{tikzpicture}
$$
the diagram above represents the functor of cutting between the levels $s$ and $(s-1)$ of the graph $\widehat{PG}_l$ - the green dashed line of the diagram. The resulting open graph denoted $PG_l(s)$ on the left side of the diagram, endowed with the same representation $\rho$ will now have the ends labeled by $P^s$ and $P^{s-1}$. The `Quantum Invariants' functor will restore the ties between the $s$ and $(s-1)$ levels but will produce {\it  apriori} {\it different} labels on the restored edges. More precisely, we have the following.
\begin{pro}\label{pro:tauZ0+-s-maps}
	Let $PG_l(s)$ be the graph
	$$
	\begin{tikzpicture}
		[place/.style={circle,draw=black,thick},
		transition/.style={circle,draw=black,fill=black}]
		\begin{scope}[scale=0.7, shift={(-2,0.5)}]
			\node (white3) at (-4,2) [place] [label={above:$s$}] {};
			\node (black3) at (-4,0) [transition] [label={below:$(s)'$}]{};
			\node (white0) at (0,2) [place] [label={above:$s-1$}] {};
			\node (black0) at (0,0) [transition] [label={below:$(s-1)'$}]{};
			\node (PG) at (-2,1.5) [label={below:$PG_l(s)$}]{};
			\draw[black, dotted]
			(black3) to (black0);
			\draw[black, dotted]
			(white3) to (white0);
			\draw[black, very thick]
			(white3) to (black3);
			\draw[black, thick]
			(white0) to (black0);
			\draw[black, dotted]
			(black3) to (black0);
			\draw[black, dotted]
			(white3) to (white0);
			\draw[black, thick]
			(white3) to (black3);
			\draw[black, very thick]
			(white0) to (black0);
		\end{scope}
	\end{tikzpicture}
$$
obtained from the cyclic graph $\widehat{PG}_l$ by cutting between the levels $s$ and $(s-1)$. Assume we are given a representation
$$
\rho(s)=\{\alpha^{i,j}: P^j \longrightarrow P^i\}
$$
of the quiver $PG_l(s)$. Then the zig-zag paths $Z_{0,\pm}$ give rise to two
distinguished linear maps
$$
	\tau^{\pm}_{\rho(s)}: Hom(\CC[q_s,q^{-1}_s],\CC) \longrightarrow P^{\pm}(s),
$$
where $P^{+}(s)=Hom(P^{s-1}, P^s)$ and  $P^{-}(s)=Hom(P^{s}, P^{s-1})$. For $s=0$, the spaces $P^{\pm}(0)=P^{\pm}$ and the maps $\tau^{\pm}_{\rho(0)}$ coincide with $\tau^{\pm}_{\rho}$ in \eqref{dualC+-maps-q-1}.
\end{pro}

Putting the maps $\tau^{\pm}_{\rho(s)}$ together gives the following.
\begin{cor}\label{cor:defrepPGl}
	A representation $\rho=\{\alpha^{t,s}:P^s \longrightarrow P^t\}$ of the quiver $PG_l$ admits the extension to the cyclic quiver $\widehat{PG}_l$ via the maps
	$$
		\tau^{\pm}_{\rho(0)}: Hom(\CC[q_0,q^{-1}_0],\CC) \longrightarrow P^{\pm}(0)=P^{\pm}.
		$$
		This in turn gives rise to the linear map
		$$
	\tau^{\pm}_{\rho(\bullet)}:	\bigotimes^{l-1}_{s=0}Hom(\CC[q_s,q^{-1}_s],\CC) \longrightarrow {\mathfrak g}^{(\pm1)}_{\rho}
		$$
		where ${\mathfrak g}_{\rho}=End_{\CC}(W_{\rho})$ is endowed with the grading coming from the direct sum decomposition
		$$
		W_{\rho}=\bigoplus_{s\in \ZZ/l\ZZ} P^s,
		$$
		that is, ${\mathfrak g}^{(i)}_{\rho}$ is the subspace of ${\mathfrak g}_{\rho}$ of endomorphisms of degree $i$ with respect to the $\ZZ/l\ZZ$ grading of
		$W_{\rho}$. More precisely, the map $\tau^{-}_{\rho(\bullet)}$ sends a vector 
		$$
		f_{l-1} \otimes f_{l-2} \otimes \cdots \otimes f_{1} \otimes f_0,\,\, f_i \in Hom(\CC[q_i,q^{-1}_i], \CC),
		$$
	into the cyclic chain of maps
		$$
		\xymatrix@C=37pt{
		P^{l-1} \ar[r]^{\tau^-_{\rho(l-1)}(f_{l-1})}&P^{l-2} \ar[r]^(.55){\tau^-_{\rho(l-2)}(f_{l-2})}&\cdots \ar[r]&P^i \ar[r]^{\tau^-_{\rho(i)}(f_i)}& P^{i-1} \ar[r]& \cdots \ar[r]&P^1 \ar[r]^{\tau^-_{\rho(1)}(f_1)}& P^0; \ar@/^2.5pc/[0,-7]^{\tau^-_{\rho(0)}(f_0)}
	}
$$
and the map $\tau^+_{\rho(\bullet)}$ is similar to the above with the directions of arrows reversed:
$$
\xymatrix@C=36pt{
	P^{l-1} \ar@/_2.5pc/[0,7]^{\tau^+(f_{l-1})}&P^{l-2} \ar[l]_{\tau^+(f_{l-2})}&\cdots \ar[l]&P^i \ar[l]_{\tau^+(f_i)}& P^{i-1} \ar[l]& \cdots \ar[l]&P^1 \ar[l]_{\tau^+(f_1)}& P^0 \ar[l]_{\tau^+(f_0)};
}
$$
in the diagram above the reference to the representation $\rho(\bullet)$ is omitted to simplify the notation.
\end{cor}

Informally, the above statement says that starting with a representation of a quiver $PG_l$ we have a natural way to deform it with the deformation parameters being the space
$$
\bigotimes^{l-1}_{s=0}Hom(\CC[q_s,q^{-1}_s],\CC).
$$
For every factor above we have the inclusion
$$
ev_s: \CC^{\times} \hookrightarrow Hom(\CC[q_s,q^{-1}_s],\CC)
$$
given by the evaluation at a closed point of $\CC^{\times}$. This gives a geometric realization of deformations for a representation of $PG_l$.
\begin{cor}\label{cor:defrepPGl-ev}
	A representation $\rho=\{\alpha^{t,s}:P^s \longrightarrow P^t\}$ of the quiver $PG_l$ admits the extension to the cyclic quiver $\widehat{PG}_l$ via the maps
	$$
	\mu^{\pm}_{\rho(0)}: \CC^{\times} \longrightarrow P^{\pm}(0)=P^{\pm}.
	$$
	This in turn gives rise to the algebraic maps
	$$
	\mu^{\pm}_{\rho(\bullet)}:=\prod^{l-1}_{i=0} \mu^{\pm}_{\rho(i)}:	(\CC^{\times} )^l \longrightarrow {\mathfrak g}^{(\pm1)}_{\rho}
	$$
	where ${\mathfrak g}_{\rho}$ and ${\mathfrak g}^{(\pm1)}_{\rho}$ are as in Corollary \ref{cor:defrepPGl}, the factors of  $(\CC^{\times} )^l$ are labeled by the white vertices of $PG_l$ and $\mu^{\pm}_{\rho(i)}$ are defined by the formula
	$$
	\mu^{\pm}_{\rho(i)}=\tau^{\pm}_{\rho(i)} \circ ev_{i}: \CC^{\times} \longrightarrow P^{\pm}(i).
	$$
	 More precisely, the map $\mu^{-}_{\rho(\bullet)}$ sends a point 
	$$
	{\bf z}=(z_{l-1}, z_{l-2}, \ldots, z_{1}, z_0) \in 	(\CC^{\times} )^l,
	$$
	into the cyclic chain of maps
	$$
	\xymatrix@C=37pt{
		P^{l-1} \ar[r]^(.55){\mu^-_{\rho(l-1)}(z_{l-1})}&P^{l-2} \ar[r]^(.55){\mu^-_{\rho(l-2)}(z_{l-2})}&\cdots \ar[r]&P^i \ar[r]^{\mu^-_{\rho(i)}(z_i)}& P^{i-1} \ar[r]& \cdots \ar[r]&P^1 \ar[r]^{\mu^-_{\rho(1)}(z_1)}& P^0; \ar@/^2.5pc/[0,-7]^{\mu^-_{\rho(0)}(z_0)}
	}
	$$
	and the map $\mu^+_{\rho(\bullet)}$ is similar to the above with the directions of arrows reversed:
	$$
	\xymatrix@C=37pt{
		P^{l-1} \ar@/_2.5pc/[0,7]^{\mu^+(z_{l-1})}&P^{l-2} \ar[l]_{\mu^+(z_{l-2})}&\cdots \ar[l]&P^i \ar[l]_{\mu^+(z_i)}& P^{i-1} \ar[l]& \cdots \ar[l]&P^1 \ar[l]_{\mu^+(z_1)}& P^0 \ar[l]_{\mu^+(z_0)};
	}
	$$
	in the diagram above the subscripts $\rho(i)$ in the notation $\mu^+_{\rho(i)}(z_i)$ are omitted for simplicity.
\end{cor}

\vspace{0.4cm}
\noindent{\bf 11.4.2. Perfect matchings.} The graph $\widehat{PG}_l$ is bipartite and as such it comes with the set of perfect matchings. The reader may recall that we already discussed them in \S10. The main emphasis there was on the connections between the perfect matchings and closed paths in $\widehat{PG}_l$. The quantum-type construction gives a different perspective
and suggests to view the perfect matchings as a devise to cut - and - paste
the graph $\widehat{PG}_l$. We recall that the set of perfect matchings of
$\widehat{PG}_l$ is denoted by ${\cal M}_{\widehat{PG}_l}$; for a perfect matching $M$, its edges colored by $0$ form the subset $M(0)$ of $M$; if that set is nonempty we write
$$
M(0)=\{e^0_{i_1}, \ldots,e^0_{i_m}, \}
$$
for the $0$-colored edges composing $M(0)$; recall $e^0_s$ is the edge
$(s) \rightarrow (s)'$, while $e^{\pm}_s$ denote the edges $(s) \rightarrow (s\pm1)'$; the number $m$ must be of the same parity as $l$; also recall our convention: the vertices $(s)$ and $(s)'$ are thought as points on two equatorial circles of the torus $\mathbb{T}$ placed in the counterclockwise order from $(l-1)$ to $0$; the edge $(s) \rightarrow (s)'$ is an arc of the meridian circle of the torus. The following drawing illustrates the situation.
 
$$
\begin{tikzpicture}
		[place/.style={circle,draw=black,thick},
	transition/.style={circle,draw=black,fill=black}]
	\begin{scope} 
			\draw[black,dotted]
		(4,0) circle [radius=1.4];
		\draw[black,dotted]
		(4,0) circle [radius=2.8];
	\node (white3) at (3,1) [place] [label={below right:$l-1$}] {};
	\node (black3) at (2,2) [transition] [label={above:$(l-1)'$}]{};
		\draw[black, thick]
	(white3) to (black3);
\end{scope}	
\begin{scope} [rotate around={-40 : (4,0)}]
	\draw[black,dotted]
	(4,0) circle [radius=1.4];
	\draw[black,dotted]
	(4,0) circle [radius=2.8];
	\node (white3) at (3,1) [place] [label={below right:$0$}] {};
	\node (black3) at (2,2) [transition] [label={above:$(0)'$}]{};
	\draw[black, thick]
	(white3) to (black3);
\end{scope}	
\begin{scope} [rotate around={-80 : (4,0)}]
	\draw[black,dotted]
	(4,0) circle [radius=1.4];
	\draw[black,dotted]
	(4,0) circle [radius=2.8];
	\node (white3) at (3,1) [place] [label={below:$1$}] {};
	\node (black3) at (2,2) [transition] [label={above:$(1)'$}]{};
	\draw[black, thick]
	(white3) to (black3);
\end{scope}	
	\begin{scope} [rotate around={40 : (4,0)}]
	\draw[black,dotted]
	(4,0) circle [radius=1.4];
	\draw[black,dotted]
	(4,0) circle [radius=2.8];
	\node (white3) at (3,1) [place] [label={below right:$l-2$}] {};
	\node (black3) at (2,2) [transition] [label={above left:$(l-2)'$}]{};
	\draw[black, thick]
	(white3) to (black3);
\end{scope}	
\begin{scope} [rotate around={120:(4,0)}]
	\draw[black,dotted]
	(4,0) circle [radius=1.4];
	\draw[black,dotted]
	(4,0) circle [radius=2.8];
	\node (white3) at (3,1) [place] [label={above:$s$}] {};
	\node (black3) at (2,2) [transition] [label={below:$(s)'$}]{};
	\draw[black, thick]
	(white3) to (black3);
\end{scope}	
\begin{scope} [rotate around={160:(4,0)}]
	\draw[black,dotted]
	(4,0) circle [radius=1.4];
	\draw[black,dotted]
	(4,0) circle [radius=2.8];
	\node (white3) at (3,1) [place] [label={above:$s-1$}] {};
	\node (black3) at (2,2) [transition] [label={below:$(s-1)'$}]{};
	\draw[black, thick]
	(white3) to (black3);
\end{scope}	
\draw[very thick,loosely dotted] (1.7,0) arc [start angle=180, end angle= 240, radius=2.3];
\begin{scope} [rotate around={150:(4,0)}]
	\draw[very thick,loosely dotted] (1.7,0) arc [start angle=180, end angle= 240, radius=2.3];
\end{scope}
\end{tikzpicture}
$$

Let $M$ be a perfect matching with $M(0)=\{e^0_{i_1}, \ldots,e^0_{i_m} \}$ nonempty, proper subset of $M$ and ordered counterclockwise in the decreasing order
\begin{equation}\label{M0seq}
l-1 \geq i_1 > i_2 > \ldots > i_m \geq 0.
\end{equation}
We think of $M(0)$ as a pattern to cut the torus $\mathbb{T}$ into cylinders. Namely, the complement $M\setminus M(0)$ is nonempty and consists of disjoint union of uninterrupted string of pairs of $\pm$-colored edges. Call those $\pm$-strings
$
M_1, \ldots, M_s
$
and arrange them in order they are encountered when we walk around the circle counterclockwise starting from $i_1$:
$$
M_1 >\cdots >M_s.
$$  
We think of each $M_i$ as a cylinder by inserting the cuts immediately before (resp. after) the first (resp. last) vertex of $M_i$.
\begin{example}
	Take $l=10$. Since the value is even a matching $M$ must contain an even number $d$ of $0$-colored edges; consider the case $d=6$.
	
	1)  Take $M(0)$ as follows 
	$$
M(0) =\{	e^0_9, e^0_8, e^0_3, e^0_2,e^0_1, e^0_0 \},
$$
that is, the corresponding sequence is
$$
9>8>3>2>1>0;
$$
the complement of $M(0)$ is a single uninterrupted string of pairs of $\pm$-edges
$$
M_1= \{(e^-_7,e^+_6), (e^-_5,e^+_4)\}
$$
we cut out the cylinder containing the arc from the vertex $7$ to $4$ included:
$$
\begin{tikzpicture}
	[place/.style={circle,draw=black,thick},
	transition/.style={circle,draw=black,fill=black}]
	\begin{scope} 
		\draw[black,dotted]
		(4,0) circle [radius=1.4];
		\draw[black,dotted]
		(4,0) circle [radius=2.8];
		\node (white3) at (3,1) [place] [label={below right:$9$}] {};
		\node (black3) at (2,2) [transition] [label={above:$(9)'$}]{};
		\draw[black, thick]
		(white3) to (black3);
	\end{scope}	
	\begin{scope} [rotate around={-40 : (4,0)}]
		\draw[black,dotted]
		(4,0) circle [radius=1.4];
		\draw[black,dotted]
		(4,0) circle [radius=2.8];
		\node (white3) at (3,1) [place] [label={below:$0$}] {};
		\node (black3) at (2,2) [transition] [label={above:$(0)'$}]{};
		\draw[black, thick]
		(white3) to (black3);
	\end{scope}	
	\begin{scope} [rotate around={-80 : (4,0)}]
		\draw[black,dotted]
		(4,0) circle [radius=1.4];
		\draw[black,dotted]
		(4,0) circle [radius=2.8];
		\node (white3) at (3,1) [place] [label={below left:$1$}] {};
		\node (black3) at (2,2) [transition] [label={above right:$(1)'$}]{};
		\draw[black, thick]
		(white3) to (black3);
	\end{scope}	
	\begin{scope} [rotate around={40 : (4,0)}]
		\draw[black,dotted]
		(4,0) circle [radius=1.4];
		\draw[black,dotted]
		(4,0) circle [radius=2.8];
		\node (white3) at (3,1) [place] [label={right:$8$}] {};
		\node (black3) at (2,2) [transition] [label={left:$(8)'$}]{};
		\draw[black, thick]
		(white3) to (black3);
	\end{scope}	
\begin{scope} [rotate around={80 : (4,0)}]
	\draw[black,dotted]
	(4,0) circle [radius=1.4];
	\draw[black,dotted]
	(4,0) circle [radius=2.8];
	\node (white3) at (3,1) [place] [label={above right:$7$}] {};
	\node (black3) at (2,2) [transition] [label={below left:$(7)'$}]{};
	\draw[black, thick]
	(white3) to (black3);
\end{scope}	
	\begin{scope} [rotate around={-120:(4,0)}]
		\draw[black,dotted]
		(4,0) circle [radius=1.4];
		\draw[black,dotted]
		(4,0) circle [radius=2.8];
		\node (white3) at (3,1) [place] [label={left:$2$}] {};
		\node (black3) at (2,2) [transition] [label={right:$(2)'$}]{};
		\draw[black, thick]
		(white3) to (black3);
	\end{scope}	
	\begin{scope} [rotate around={-160:(4,0)}]
	\draw[black,dotted]
	(4,0) circle [radius=1.4];
	\draw[black,dotted]
	(4,0) circle [radius=2.8];
	\node (white3) at (3,1) [place] [label={left:$3$}] {};
	\node (black3) at (2,2) [transition] [label={right:$(3)'$}]{};
	\draw[black, thick]
	(white3) to (black3);
\end{scope}	
\begin{scope} [rotate around={-200:(4,0)}]
	\draw[black,dotted]
	(4,0) circle [radius=1.4];
	\draw[black,dotted]
	(4,0) circle [radius=2.8];
	\node (white3) at (3,1) [place] [label={above left:$4$}] {};
	\node (black3) at (2,2) [transition] [label={below right:$(4)'$}]{};
	\draw[black, thick]
	(white3) to (black3);
\end{scope}	
\begin{scope} [rotate around={60:(4,0)}]
\draw[green, very thick, dashed]
(3,1) -- (2,2);
\end{scope}
\begin{scope} [rotate around={180:(4,0)}]
	\draw[green, very thick, dashed]
	(3,1) -- (2,2);
\end{scope}
	\begin{scope} [rotate around={45:(4,0)}]
		\draw[very thick,loosely dotted] (1.7,0) arc [start angle=180, end angle= 240, radius=2.3];
\end{scope}
\end{tikzpicture}
$$
the green lines on the picture show the cuts.

2) Take $M(0)$ given by the sequence
$$
8>7>4>3>2>1
$$
there are two uninterrupted $\pm$-strings
$$
M_1=\{e^-_6,e^+_5\} > M_2=\{e^-_0, e^+_9\};
$$
there are four cuts: two for the cylinder containing $M_1$ and two more for the cylinder containing $M_2$: 
$$
\begin{tikzpicture}
	[place/.style={circle,draw=black,thick},
	transition/.style={circle,draw=black,fill=black}]
	\begin{scope} 
		\draw[black,dotted]
		(4,0) circle [radius=1.4];
		\draw[black,dotted]
		(4,0) circle [radius=2.8];
		\node (white3) at (3,1) [place] [label={below right:$9$}] {};
		\node (black3) at (2,2) [transition] [label={above:$(9)'$}]{};
		\draw[black, thick]
		(white3) to (black3);
	\end{scope}	
	\begin{scope} [rotate around={-40 : (4,0)}]
		\draw[black,dotted]
		(4,0) circle [radius=1.4];
		\draw[black,dotted]
		(4,0) circle [radius=2.8];
		\node (white3) at (3,1) [place] [label={below:$0$}] {};
		\node (black3) at (2,2) [transition] [label={above:$(0)'$}]{};
		\draw[black, thick]
		(white3) to (black3);
	\end{scope}	
	\begin{scope} [rotate around={-80 : (4,0)}]
		\draw[black,dotted]
		(4,0) circle [radius=1.4];
		\draw[black,dotted]
		(4,0) circle [radius=2.8];
		\node (white3) at (3,1) [place] [label={below left:$1$}] {};
		\node (black3) at (2,2) [transition] [label={above right:$(1)'$}]{};
		\draw[black, thick]
		(white3) to (black3);
	\end{scope}	
	\begin{scope} [rotate around={40 : (4,0)}]
		\draw[black,dotted]
		(4,0) circle [radius=1.4];
		\draw[black,dotted]
		(4,0) circle [radius=2.8];
		\node (white3) at (3,1) [place] [label={right:$8$}] {};
		\node (black3) at (2,2) [transition] [label={left:$(8)'$}]{};
		\draw[black, thick]
		(white3) to (black3);
	\end{scope}	
	\begin{scope} [rotate around={80 : (4,0)}]
		\draw[black,dotted]
		(4,0) circle [radius=1.4];
		\draw[black,dotted]
		(4,0) circle [radius=2.8];
		\node (white3) at (3,1) [place] [label={above right:$7$}] {};
		\node (black3) at (2,2) [transition] [label={below left:$(7)'$}]{};
		\draw[black, thick]
		(white3) to (black3);
	\end{scope}	
	\begin{scope} [rotate around={-120:(4,0)}]
		\draw[black,dotted]
		(4,0) circle [radius=1.4];
		\draw[black,dotted]
		(4,0) circle [radius=2.8];
		\node (white3) at (3,1) [place] [label={left:$2$}] {};
		\node (black3) at (2,2) [transition] [label={right:$(2)'$}]{};
		\draw[black, thick]
		(white3) to (black3);
	\end{scope}	
	\begin{scope} [rotate around={-160:(4,0)}]
		\draw[black,dotted]
		(4,0) circle [radius=1.4];
		\draw[black,dotted]
		(4,0) circle [radius=2.8];
		\node (white3) at (3,1) [place] [label={left:$3$}] {};
		\node (black3) at (2,2) [transition] [label={right:$(3)'$}]{};
		\draw[black, thick]
		(white3) to (black3);
	\end{scope}	
	\begin{scope} [rotate around={-200:(4,0)}]
		\draw[black,dotted]
		(4,0) circle [radius=1.4];
		\draw[black,dotted]
		(4,0) circle [radius=2.8];
		\node (white3) at (3,1) [place] [label={above left:$4$}] {};
		\node (black3) at (2,2) [transition] [label={below right:$(4)'$}]{};
		\draw[black, thick]
		(white3) to (black3);
	\end{scope}	
	\begin{scope} [rotate around={15:(4,0)}]
	\draw[green, very thick, dashed]
	(3,1) -- (2,2);
\end{scope}
\begin{scope} [rotate around={-55:(4,0)}]
	\draw[green, very thick, dashed]
	(3,1) -- (2,2);
\end{scope}
	\begin{scope} [rotate around={100:(4,0)}]
		\draw[green, very thick, dashed]
		(3,1) -- (2,2);
	\end{scope}
	\begin{scope} [rotate around={140:(4,0)}]
		\draw[green, very thick, dashed]
		(3,1) -- (2,2);
	\end{scope}
	\begin{scope} [rotate around={65:(4,0)}]
		\draw[very thick,loosely dotted] (1.7,0) arc [start angle=180, end angle= 200, radius=2.3];
	\end{scope}
\end{tikzpicture}
$$
\end{example}

Once we cut out the the cylinders corresponding to $M_i$'s, what remains are cylinders corresponding to the uninterrupted strings of $0$-colored edges.
Thus to every matching $M$ with $M(0)$ a proper subset of $M$ we associate
a slicing of the torus $\mathbb{T}$ into cylinders of two types: the cylinders containing the edges of $M(0)$ and the ones containing $\pm$-colored edges of $M$. We denote $C^0_J$ the cylinder of type $0$ stretching over the edges $e^0_j$ with $j\in J$. Similarly a cylinder of type $\pm$ is denoted $C^{\pm}_I$ - it is composed of  pairs of edges $e^-_i$, $e^{+}_{i-1}$, for $i\in I$. The edges in each cylinder are subject to counterclockwise ordering convention. Observe that $C^0_J$ (resp. $C^{\pm}_I$)
is modeled on the graph $PG_{|J|}$ (resp. $PG_{|I|}$) with vertices labeled by the set $J$ (resp. $I$); to get the edges of $C^0_J$ (resp. $C^{\pm}_I$) one forgets $\pm$-colored (resp. $0$-colored) edges of $PG_{|J|}$ (resp. $PG_{|I|}$).
 
 By construction ever cylinder of type $C^{\pm}_I$ is sandwiched between two cylinders of type $0$ which are denoted $C^0_{I^+}$ and $C^0_{I^-}$: the first one is {\it preceding} $C^{\pm}_I$ and the second one succeeding it in the counterclockwise order; observe that $C^0_{I^+}$ and $C^0_{I^-}$ coincide in case the matching $M$ cuts the torus into precisely two cylinders.
  
  We assume that a representation $\rho=\{\alpha^{t,s}:P^s \longrightarrow P^t\}$ of the quiver $PG_l$ is given and it has been extended, according to Proposition \ref{pro:mu-map}, to a representation $\widehat{\rho}$ of $\widehat{PG}_l$.
  Using quantum-type invariant construction we show how the cylinder
  $C^{\pm}_I$ provides maps between the `end' spaces $P^{I^+}$ and $P^{I^-}$, where for an ordered (counterclockwise) subset $A$ in $[0,l-1]$
  the space $P^A$ stands for the ordered tensor product of $P^a$, $a\in A$:
  $$
  P^A=\bigotimes_{a\in A} P^a.
  $$
  
  The counterclockwise orientation of the circles on which the vertices of
  $\widehat{PG}_l$ are placed provides the order of the labels
  $I^{\pm}$:
  $$
  I^+=\{ i^+_0, i^+_0-1, \ldots, i^+_0-|I^+|+1\}, \,\, I^-=\{ i^-_0, i^-_0-1, \ldots, i^-_0-|I^-|+1\}.
  $$
  For the counterclockwise ordered uninterrupted subset $A$ of $[0,l-1]$ we set $t(A)$ (resp. $h(A)$), the tale (resp. the head) of $A$ to denote the first (resp. the last) label in $A$; for example, in the sets $I^{\pm} $ above we have
  $$
  t(I^{\pm})=i^{\pm}_0,\,\, h(I^{\pm})=i^{\pm}_0-|I^{\pm}|+1.
  $$
   Thus the three consecutive cylinders $C^0_{I^-}$, $C^{\pm}_I$, $C^0_{I^+}$ correspond to the following sequence of graphs.
   $$
   \begin{tikzpicture}
   	[place/.style={circle,draw=black,thick},
   	transition/.style={circle,draw=black,fill=black}]
   	\begin{scope}[scale=0.7, shift={(-2,0.5)}]
   		\node (white3) at (-4,2) [place] [label={}] {};
   		\node (black3) at (-4,0) [transition] [label={}]{};
   		\node (white0) at (0,2) [place] [label={}] {};
   		\node (black0) at (0,0) [transition] [label={}]{};
   		\node (PG) at (-2,1.5) [label={below:$PG_{|I|}(t(I))$}]{};
   		\draw[black, dotted]
   		(black3) to (black0);
   		\draw[black, dotted]
   		(white3) to (white0);
   		\draw[black, very thick]
   		(white3) to (black3);
   		\draw[black, thick]
   		(white0) to (black0);
   		\draw[black, dotted]
   		(black3) to (black0);
   		\draw[black, dotted]
   		(white3) to (white0);
   		\draw[black, thick]
   		(white3) to (black3);
   		\draw[black, very thick]
   		(white0) to (black0);
   	\end{scope}
   \begin{scope}[scale=0.7, shift={(-7,0.5)}]
   	\node (white3) at (-4,2) [place] [label={above:$i^+_0$}] {};
   	\node (black3) at (-4,0) [transition] [label={below:$(i^+_0)'$}]{};
   	\node (white0) at (0,2) [place] [label={above:$h(I^+)$}] {};
   	\node (black0) at (0,0) [transition] [label={below:$(h(I^+))'$}]{};
   	\node (PG) at (-2,1.5) [label={below:$PG_{|I^+|}(i^+_0)$}]{};
   	\draw[black, dotted]
   	(black3) to (black0);
   	\draw[black, dotted]
   	(white3) to (white0);
   	\draw[black, very thick]
   	(white3) to (black3);
   	\draw[black, thick]
   	(white0) to (black0);
   	\draw[black, dotted]
   	(black3) to (black0);
   	\draw[black, dotted]
   	(white3) to (white0);
   	\draw[black, thick]
   	(white3) to (black3);
   	\draw[black, very thick]
   	(white0) to (black0);
   \end{scope}
 \begin{scope}[scale=0.7, shift={(3,0.5)}]
	\node (white3) at (-4,2) [place] [label={above:$i^-_0$}] {};
	\node (black3) at (-4,0) [transition] [label={below:$(i^-_0)'$}]{};
	\node (white0) at (0,2) [place] [label={above:$h(I^-)$}] {};
	\node (black0) at (0,0) [transition] [label={below:$(h(I^-))'$}]{};
	\node (PG) at (-2,1.5) [label={below:$PG_{|I^-|}(i^-_0)$}]{};
	\draw[black, dotted]
	(black3) to (black0);
	\draw[black, dotted]
	(white3) to (white0);
	\draw[black, very thick]
	(white3) to (black3);
	\draw[black, thick]
	(white0) to (black0);
	\draw[black, dotted]
	(black3) to (black0);
	\draw[black, dotted]
	(white3) to (white0);
	\draw[black, thick]
	(white3) to (black3);
	\draw[black, very thick]
	(white0) to (black0);
\end{scope}
   \end{tikzpicture}
   $$
 Let $s^-$ be a label in $I^-$ and $t^+$ a label in $I^+$ always assumed to be different from $s^-$ (the last assumption is superfluous as long as the intervals $I^-$ and $I^+$ are distinct; if $I^+=I^-$, we assume that  walking from $t^+$  counterclockwise  we arrive to $s^-$ before completing the full turn, that is, the arc $(t^+ \rightarrow s^-) $ is less than $2\pi$). The representation $\rho$ attaches the vector spaces $P^{s^-}$ and $P^{t^+}$ to those labels. By `fusing' the graph $PG_{|I|}(t(I))$, on the left
 
 - with the part
 of the graph $PG_{|I^+|} (i^+_0)$ stretching from $t^+$ to $h(I^+) $,
 
 and on the right
 
 - with the part of the graph $PG_{|I^-|}(i^-_0)$ from $i^-_0$ to $s^-$,
 
  we obtain the graph
 
 $$
  \begin{tikzpicture}
 	[place/.style={circle,draw=black,thick},
 	transition/.style={circle,draw=black,fill=black}]
 	\begin{scope}[scale=0.7, shift={(-2,0.5)}]
 		\node (white3) at (-4,2) [place] [label={}] {};
 		\node (black3) at (-4,0) [transition] [label={}]{};
 		\node (white0) at (0,2) [place] [label={}] {};
 		\node (black0) at (0,0) [transition] [label={}]{};
 		\node (PG) at (-2,1.5) [label={below:$PG_{|I|}(t(I))$}]{};
 		\draw[black, dotted]
 		(black3) to (black0);
 		\draw[black, dotted]
 		(white3) to (white0);
 		\draw[black, very thick]
 		(white3) to (black3);
 		\draw[black, thick]
 		(white0) to (black0);
 		\draw[black, dotted]
 		(black3) to (black0);
 		\draw[black, dotted]
 		(white3) to (white0);
 		\draw[black, thick]
 		(white3) to (black3);
 		\draw[black, very thick]
 		(white0) to (black0);
 	\end{scope}
 	\begin{scope}[scale=0.7, shift={(-3.5,0.5)}]
 	\node (white3) at (-4,2) [place] [label=above:{$h(I^+)$}] {};
 	\node (black3) at (-4,0) [transition] [label=below:{$(h(I^+))'$}]{};
 	\draw[black, very thick]
 	(white3) to (black3);
 	\draw[blue, very thick]
 	(-4,2) --(-2.7,0);
 	\draw[red, dashed]
 	(-4,0) --(-2.7,2);
 	\end{scope}
 	\begin{scope}[scale=0.7, shift={(-6,0.5)}]
 	\node (white3) at (-4,2) [place] [label=above:{$t^+$}] {};
 	\node (black3) at (-4,0) [transition] [label=below:{$(t^+)'$}]{};
 	\draw[black, very thick]
 	(white3) to (black3);
 	\draw[black, dotted]
 	(-4,0)--(-1.8,0);
 	\draw[black, dotted]
 	(-4,2)--(-1.8,2);
 	\end{scope}
 	\begin{scope}[scale=0.7, shift={(-0.5,0.5)}]
 	\node (white0) at (0,2) [place] [label=above:{$i^-_0$}] {};
 	\node (black0) at (0,0) [transition] [label=below:{$(i^-_0)'$}]{};
 	\draw[black, very thick]
 	(white0) to (black0); 	
 	\draw[black, dotted]
 	(0.2,2)--(2.2,2);
 	\draw[black, dotted]
 	(0.2,0)--(2.2,0);
 	\draw[red,dashed]
 	(0,2)--(-1.3,0);
 	\draw[blue,very thick]
 	(0,0)--(-1.3,2);
 	\end{scope}
 
 \begin{scope}[scale=0.7, shift={(2,0.5)}]
 	\node (white0) at (0,2) [place] [label=above:{$s^-$}] {};
 	\node (black0) at (0,0) [transition] [label=below:{$(s^-)'$}]{};
 	\draw[black, very thick]
 	(white0) to (black0); 	
 	\end{scope}
\end{tikzpicture}
$$
The quantum-type invariants for this graph give the maps
$$
\begin{gathered}
\tau^{-}_{\rho}(s^-,t^+): Hom_{\CC} (\CC[q_{t^+ s^-}, q^{-1}_{t^+ s^-}], \CC) \longrightarrow Hom (P^{s^-},P^{t^+}),
\\
\tau^{+}_{\rho} (s^-,t^+): Hom_{\CC} (\CC[q_{t^+ s^-}, q^{-1}_{t^+ s^-}], \CC) \longrightarrow Hom (P^{t^+},P^{s^-}),
\end{gathered}
$$
where $q_{t^+ s^-}$ is a formal variable. All those maps can be arranged in the matrix form as follows
\begin{equation}\label{I+-matrix}
	\begin{gathered}
	\tau^-_{\rho}(I^-,I^+)= \begin{pmatrix}
		\tau^{-}_{\rho}(t(I^-),t(I^+)) &\cdots&\tau^{-}_{\rho}(h(I^-),t(I^+))\\
		\vdots &\vdots&\vdots\\
			\tau^{-}_{\rho}(t(I^-),h(I^+)) &\cdots&\tau^{-}_{\rho}(h(I^-),h(I^+))
	\end{pmatrix}
\\
\\
	\tau^+_{\rho}(I^-,I^+)= \begin{pmatrix}
	\tau^{+}_{\rho}(t(I^-),t(I^+)) &\cdots&\tau^{+}_{\rho}(t(I^-),h(I^+))\\
	\vdots &\vdots&\vdots\\
	\tau^{+}_{\rho}(h(I^-),t(I^+)) &\cdots&\tau^{+}_{\rho}(h(I^-),h(I^+))
\end{pmatrix}
\end{gathered}
\end{equation} 
In the first matrix the rows (resp. columns) are labeled by the ordered set $I^+$  (resp. $I^-$) from the top (resp. left) to bottom (resp. right); in the second matrix the roles of rows and columns are inverted. 

We can use the above maps to constructs various `higher order multiplications' in the graded tensor algebra 
$$
\text{$T(W_{\rho}):=\bigoplus^{\infty}_{n=0} \left(\bigoplus^{l-1}_{s=0} P^s \right)^{\otimes n},$ where $W_{\rho}=\bigoplus^{l-1}_{s=0} P^s$.}
$$
Namely, for every perfect matching $M$ with $M(0) =\{e^0_{i_1}, e^0_{i_2}, \ldots, e^0_{i_m} \}$ a proper subset of $M$,
we divide the the torus $\mathbb{T}$ into the the triples of disjoint cylinders $\{C^0_{I^-}, C^{\pm}_I, C^0_{I^+}\}$ as described above, where $I$ ranges through the `connected' components (uninterrupted strings of integers) of the set
$$
[0,l-1]\setminus \{i_1>i_2> \cdots>i_m\};
$$
for every weakly order preserving map
$$
\sigma^{I^-}_{I^+} : I^+ \longrightarrow I^-
$$
we have linear maps
$$
m^-_{M,\sigma^{I^-}_{I^+}} \left(\bigotimes_{t\in I^+} f_t \right): P^{im (\sigma^{I^-}_{I^+})}=\bigotimes_{t\in I^+} P^{\sigma^{I^-}_{I^+}(t)} \longrightarrow \bigotimes_{t\in I^+} P^t= P^{I^+},
$$
where $f_t \in \CC[q_{\sigma^{I^-}_{I^+}(t)t }, q^{-1}_{\sigma^{I^-}_{I^+}(t)t }]$, for $t\in I^+$, and  the map is defined by the equation
$$
m^-_{M,\sigma^{I^-}_{I^+}} \left(\bigotimes_{t\in I^+} f_t \right)=\bigotimes_{t\in I^+} \tau^-_{\rho}(\sigma^{I^-}_{I^+}(t),t)( f_t).
$$
Those are the maps from the cylinder $C^0_{I^-}$ to the cylinder $C^0_{I^+}$, labeled with the representation $\rho$. Similarly we have maps going in the other direction: 
 $$
 m^+_{M,\sigma^{I^-}_{I^+}} \left(\bigotimes_{t\in I^+} f_t \right): P^{I^+}=\bigotimes_{t\in I^+} P^t \longrightarrow \bigotimes_{t\in I^+} P^{\sigma^{I^-}_{I^+}(t)} =P^{im(\sigma^{I^-}_{I^+})},
 $$
 
It should be clear that the organizing principle behind the above discussion is the one of categories:
\begin{equation}\label{match-cat}
	\text{\it the set of perfect matchings ${\cal M}_{\widehat{PG}_l}$ of the graph $\widehat{PG}_l$ should be categorified.}
\end{equation}
 Actually, the constructions above suggest that a perfect matching $M$ with
$M(0)$ a proper subset of $M$ could itself be considered as a category: its objects are cylinders $C^0_{I^{\pm}}$, where $I$ runs through the set of connected components of the complement in $[0,l-1]$ of the set of vertices covered by the edges in $M(0)$; its morphisms are cylinders $C^{\pm}_I$ and their concatenations. The particularity of the category is the following: given a representation $\rho=\{a^{t,s}:V^s \longrightarrow V^t\}$ of $\widehat{PG}_l$ in the category of finite dimensional vector spaces over a field $k$, subject to the bipartite condition that the same space $V^i$ is attached to the black-white pair $(i)-(i)'$, for every $i \in [0,l-1]$, the tensor algebra
$$
T(V)=\bigoplus^{\infty}_{n=0} \left(\bigoplus^{l-1}_{i=0} V^i\right)^{\otimes n}
$$
is equipped with maps between the graded pieces which are linked to geometry and combinatorics associated to $\widehat{PG}_l$; the construction of those maps as we already mentioned is reminiscent to TQFT.

  Let us return to this analogy and see how it looks from the point of view of cutting the torus $\mathbb{T}$ into cylinders: recall that the graph $\widehat{PG}_l$ is embedded into the torus $\mathbb{T}$, and we suggested to view the end vertices $(0)$ and $(l-1)$ of $PG_l$ as points on the torus, cut out two small open disks $D_0$ and $D_{l-1}$ around those points and then view our quantum-type invariants as maps associated to the cobordism   between the boundary circles $\partial D_0$ and $\partial D_{l-1}$ provided by the surface $\mathbb{T} \setminus (D_0 \cup D_{l-1})$.
In the interpretation of matchings as a pattern to cut the torus into cylinders, we think of two edges $e^0_0 $ and $e^0_{l-1}$ to be  on the boundary of the cylinder $C^0_{I^-} =C^0_{I^+}$ and the cylinder $C^{\pm}_I$ being the complementary one and where the set $I$ is the complement of $\{0,l-1\}$:
$$
I=[0,l-1]\setminus\{0,l-1\},
$$
that is we insert two cuts as indicated by dashed green segments in the following drawing 
$$
\begin{tikzpicture}
	[place/.style={circle,draw=black,thick},
	transition/.style={circle,draw=black,fill=black}]
	\begin{scope} 
		\draw[black,dotted]
		(4,0) circle [radius=1.4];
		\draw[black,dotted]
		(4,0) circle [radius=2.8];
		\node (white3) at (3,1) [place] [label={below right:$l-1$}] {};
		\node (black3) at (2,2) [transition] [label={above:$(l-1)'$}]{};
		\draw[black, thick]
		(white3) to (black3);
	\end{scope}	
	\begin{scope} [rotate around={-40 : (4,0)}]
		\draw[black,dotted]
		(4,0) circle [radius=1.4];
		\draw[black,dotted]
		(4,0) circle [radius=2.8];
		\node (white3) at (3,1) [place] [label={below right:$0$}] {};
		\node (black3) at (2,2) [transition] [label={above:$(0)'$}]{};
		\draw[black, thick]
		(white3) to (black3);
	\end{scope}	
	\begin{scope} [rotate around={-80 : (4,0)}]
		\draw[black,dotted]
		(4,0) circle [radius=1.4];
		\draw[black,dotted]
		(4,0) circle [radius=2.8];
		\node (white3) at (3,1) [place] [label={below:$1$}] {};
		\node (black3) at (2,2) [transition] [label={above:$(1)'$}]{};
		\draw[black, thick]
		(white3) to (black3);
	\end{scope}	
	\begin{scope} [rotate around={40 : (4,0)}]
		\draw[black,dotted]
		(4,0) circle [radius=1.4];
		\draw[black,dotted]
		(4,0) circle [radius=2.8];
		\node (white3) at (3,1) [place] [label={below right:$l-2$}] {};
		\node (black3) at (2,2) [transition] [label={above left:$(l-2)'$}]{};
		\draw[black, thick]
		(white3) to (black3);
	\end{scope}	
	\begin{scope} [rotate around={120:(4,0)}]
		\draw[black,dotted]
		(4,0) circle [radius=1.4];
		\draw[black,dotted]
		(4,0) circle [radius=2.8];
		\node (white3) at (3,1) [place] [label={above:$s$}] {};
		\node (black3) at (2,2) [transition] [label={below:$(s)'$}]{};
		\draw[black, thick]
		(white3) to (black3);
	\end{scope}	
	\begin{scope} [rotate around={160:(4,0)}]
		\draw[black,dotted]
		(4,0) circle [radius=1.4];
		\draw[black,dotted]
		(4,0) circle [radius=2.8];
		\node (white3) at (3,1) [place] [label={above:$s-1$}] {};
		\node (black3) at (2,2) [transition] [label={below:$(s-1)'$}]{};
		\draw[black, thick]
		(white3) to (black3);
	\end{scope}	
\begin{scope} [rotate around={-60 : (4,0)}]
	\node (white3) at (3,1)  [label={}] {};
	\node (black3) at (2,2)  [label={}]{};
	\draw[very thick, green, dashed]
	(white3) to (black3);
\end{scope}	
\begin{scope} [rotate around={20 : (4,0)}]
	\node (white3) at (3,1)  [label={}] {};
	\node (black3) at (2,2)  [label={}]{};
	\draw[very thick, green, dashed]
	(white3) to (black3);
\end{scope}	
	\draw[very thick,loosely dotted] (1.7,0) arc [start angle=180, end angle= 240, radius=2.3];
	\begin{scope} [rotate around={150:(4,0)}]
		\draw[very thick,loosely dotted] (1.7,0) arc [start angle=180, end angle= 240, radius=2.3];
	\end{scope}
\end{tikzpicture}
$$
The quantum type invariant is a result of fusing the edges $e^0_{l-1}$ and $e^0_0$ on the left and right respectively to the graph $PG_I$, that is, getting back the graph $PG_l$. On the level of cylinders, the operation can be 
envisaged as taking the cylinder $C^{\pm}_I$ and pasting on its left end an annulus with the inside boundary labeled $(l-1)$ and on its right end an annulus with the inside boundary labeled $(0)$:
$$
\begin{tikzpicture}
	\draw (0,0) ellipse [x radius =15pt, y radius =20pt];
	\draw (-3,0.6) arc [start angle =90, end angle =270, x radius=15pt, y radius =20pt];
	\draw [dotted]
	(-3,0.6) arc [start angle =90, end angle =-90, x radius=15pt, y radius =20pt];
	\draw
	(-3,0.6) --(-0.15,0.68);
	\draw
	(0,-0.7) --(-3,-0.78);
	\node
	(PG) at (-1.5, -0.5) [label={$C^{\pm}_I$}]{};
	\begin{scope}[xshift=2.5cm]
		\node
		(0) at (0,0.1) [label={$\scriptscriptstyle{(0)}$}]{};
		\draw (0,0) circle [ radius =20pt];
		\draw[red] (0,0) circle [ radius =10pt];
	\end{scope}
\begin{scope}[xshift=-5.5cm]
	\node (L) at (0,0.1) [label={$\scriptscriptstyle{(l-1)}$}]{};
	\draw (0,0) circle [ radius =20pt];
	\draw[red] (0,0) circle [radius =10pt];
\end{scope}
	\draw[->,decorate,decoration={snake, amplitude=.4mm, segment length=2mm, post length=1mm}, >={Stealth}]
(2,0)--(0.25,0)
node[above,align=center,midway]
{\textcolor{red} {Glue}};
	\draw[->,decorate,decoration={snake, amplitude=.4mm, segment length=2mm, post length=1mm}, >={Stealth}]
(-5,0)--(-3.25,0)
node[above,align=center,midway]
{\textcolor{red} {Glue}};
\end{tikzpicture}
$$
the quantum invariants themselves tell us how to join those inside boundary circles, the red circles on the drawing, along an `outside' tube which one needs to glue in to get back the torus.
This can be generalized to morphisms between $C^0_{I^+}$ and $C^0_{I^-}$ given by $C^{\pm}_I$ for any interval $I$ associated to a matching $M$ with $M(0)$ a proper subset of $M$. Namely,
 put on each side of the cylinder $C^{\pm}_I$ the cylinders $C^0_{I^+}$ and $C^0_{I^-}$ viewed as annuli with sequences of expanding circles marked by the sets $I^+$ and $I^-$ respectively:

$$
\begin{tikzpicture}
	\draw (0,0) ellipse [x radius =15pt, y radius =20pt];
	\draw (-3,0.6) arc [start angle =90, end angle =270, x radius=15pt, y radius =20pt];
	\draw [dotted]
	(-3,0.6) arc [start angle =90, end angle =-90, x radius=15pt, y radius =20pt];
	\draw
	(-3,0.6) --(-0.15,0.68);
	\draw
	(0,-0.7) --(-3,-0.78);
	\node
	(PG) at (-1.5, -0.5) [label={$C^{\pm}_I$}]{};
	\begin{scope} [xshift=3cm]
		\draw (0,0) circle [ radius =20pt];
		\draw [red] (0,0) circle [ radius =15pt];
		\draw [red](0,0) circle [ radius =5pt];
		\draw [very thick, red, loosely dotted] 
		(0,0.25) -- (0,0.5);
		\draw [rotate=90, very thick, red, loosely dotted] 
		(0,0.25) -- (0,0.5);
		\draw [rotate=-90, very thick, red, loosely dotted]
		(0,0.25) -- (0,0.5); 
		\draw [rotate=180, very thick, red, loosely dotted] 
		(0,0.25) -- (0,0.5);  
		\node (I-) at (1.5,-0.5) [label={$C^0_{I^-}$}]{};
	\end{scope}
	\begin{scope}[xshift=-6cm]
			\draw (0,0) circle [ radius =20pt];
		\draw [red] (0,0) circle [ radius =15pt];
		\draw [red](0,0) circle [ radius =5pt];
		\draw [very thick, red, loosely dotted] 
		(0,0.25) -- (0,0.5);
		\draw [rotate=90, very thick, red, loosely dotted] 
		(0,0.25) -- (0,0.5);
		\draw [rotate=-90, very thick, red, loosely dotted]
		(0,0.25) -- (0,0.5); 
		\draw [rotate=180, very thick, red, loosely dotted] 
		(0,0.25) -- (0,0.5);  
		\node (I+) at (-1.5,-0.5) [label={$C^0_{I^+}$}]{};
	\end{scope}
	\draw[->,decorate,decoration={snake, amplitude=.4mm, segment length=2mm, post length=1mm}, >={Stealth}]
	(2.3,0)--(0.25,0)
	node[above,align=center,midway]
	{\textcolor{red} {Glue}};
	\draw[->,decorate,decoration={snake, amplitude=.4mm, segment length=2mm, post length=1mm}, >={Stealth}]
	(-5.25,0)--(-3.25,0)
	node[above,align=center,midway]
	{\textcolor{red} {Glue}};
\end{tikzpicture}
$$
the marked circles appear in red in the drawing and they expand from
the head $h(I^+)$ to the tale $t(I^+)$ on the left and from $t(I^-)$ to $h(I^-)$ on the right; the maps $\tau^{\pm}_{\rho} (s^-,t^+)$ tell us how to join the circles labeled $s^-$ and $t^+$ on the right and left respectively by traveling on the `outside' tubes. 

 As an application of the above ideas we construct families of elliptic curves attached to the representations 
 $$
 \rho_c: \OO^{\times}_{\mathfrak{L}_l}(-1) \longrightarrow \mathfrak{Reps}(PG_l)
 $$ 
 of the quiver $PG_l$ associated to the stratum ${\mathfrak L}_l$.
 
 \begin{pro}\label{pro:match-ell}
 	Let ${\mathfrak L}_l$ be the stratum of $\PP({\cal W}_{\Sigma^0_r})$ corresponding to $l\geq 5$. Let $\xi \otimes\phi \in \OO^{\times}_{\mathfrak{L}_l}(-1)$ and let $\rho_c(\xi\otimes\phi)=\{\alpha^{t,s}(\xi,\phi):P^s \longrightarrow P^t\}$ be the corresponding representation of the quiver $PG_l$ corresponding, where the spaces $\{P^s([\xi],[\phi])\}$ are the summands of the orthogonal decomposition
 	$$
 	W_{\xi}/W^l_{\xi}([\phi])=\bigoplus^{l-1}_{s=0} P^s([\xi],[\phi])
 	$$
 	at the point $([\xi],[\phi])$ of ${\mathfrak L}_l$ underlying $\xi \otimes\phi$. Extend $\rho_c(\xi\otimes\phi)$ to the representation $\widehat{\rho}_c(\xi\otimes\phi)$ of the quiver $\widehat{PG}_l$ using the maps
 	$$
 	\tau^{\pm}_{\rho_c(\xi\otimes\phi)}: Hom_{\CC}(\CC[q,q^{-1}],\CC) \longrightarrow P^{\pm}([\xi],[\phi]),
 	$$
 	according to Corollary \ref{cor:rhoc-ext}. Let $M$ be a perfect matching of $\widehat{PG}_l$ such that its decomposition into cylinders contains a triple
 	$C^0_{I^+}, C^{\pm}_I, C^0_{I^-}$ with the interval $I^+$ or $I^-$ having length at least three. Then for every linear functional 
 	$F \in Hom_{\CC}(\CC[q,q^{-1}],\CC) $ there is a distinguished holomorphic map
 	$$
 	\gamma_F: \PP^1 \setminus \{\infty,0, z_F\} \longrightarrow {\mathfrak{M}}^4_1,
 	$$
 	where ${\mathfrak{M}}^4_1$ is the moduli space of curves of genus $1$ with four marked points, and where the point $z_F \in \PP^1\setminus \{\infty,0\} \cong \CC^{\times}$ is specified in the proof.
 \end{pro}
\begin{pf}
	Assume for the sake of being explicit that it is the interval $I^{+}$ that has length at least three. Then we can choose a triple of successive edges 
	$$
	e^0_t, e^0_{t-1}, e^0_{t-2}
	$$
	in $M(0)$ with the indices $t,t-1,t-2$ in $I^+$. This means that the graph
	$PG_3$
	\begin{equation}\label{PG3qin} 
	\begin{tikzpicture}
		[place/.style={circle,draw=black,thick},
		transition/.style={circle,draw=black,fill=black}]
		\node (white3) at (0,2) [place] [label={above:$t$}] {};
		\node (black3) at (0,0) [transition] [label={below:$(t)'$}]{};
		\node (white2) at (2,2) [place] [label={above:$t-1$}]{};
		\node (black2) at (2,0) [transition] [label={below:$(t-1)'$}]{};
		\node (white1) at (4,2) [place] [label={above:$t-2$}]{};
		\node (black1) at (4,0) [transition] [label={below:$(t-2)'$}]{};
		\draw[red,ultra thick] (white1) to (black1);
		\draw[red,ultra thick] (white2) to (black2);
		\draw[red,ultra thick] (white3) to (black3);
		\begin{scope}[thick]
			\draw[blue][-] (white3) to (black2);	
			\draw[blue][-] (white2) to (black1);
		\end{scope}
		\draw[red, dotted,thick][-] (white2) to (black3);
		\draw[red,dotted,thick][-] (white1) to (black2);
	\end{tikzpicture}
\end{equation}
	is a subgraph of the cylinder $C^0_{I^+}$. At this point we need to recall that the nonabelian Dolbeault space ${\bf H^{1,0}}(PG_3)$ is the complete intersection of two quadrics
	$$
	X_t Y_t =T^2,\,\,\, X_{t-1} Y_{t-1}=T^2
	$$
	in the projective space $\PP(V)=\PP^4$, where $V$ is the vector space spanned by the edges
	$$
	e^-_t, e^+_{t-1}, e^-_{t-1}, e^+_{t-2}, e^0_{t-2}
	$$
	while $X_t,Y_t,X_{t-1},Y_{t-1},T$ linear forms on $\PP(V)$ dual to the basis above. The hyperplane sections of ${\bf H^{1,0}}(PG_3)$ are space curves of degree $4$ and arithmetic genus one. Let
	$$
	L(a,b,c,d)=aX_t +bY_t +cX_{t-1}+dY_{t-1}
	$$
	be a linear form on $V$ and consider the corresponding hyperplane section of ${\bf H^{1,0}}(PG_3)$
	$$
	C(a,b,c,d)={\bf H^{1,0}}(PG_3) \bigcap (L(a,b,c,d)=0).
	$$
	Observe the following:
	\begin{equation}\label{abcd}
		\begin{gathered}
		\text{\it $C(a,b,c,d)$ is smooth if the matrix $\begin{pmatrix}
				a&c\\
				d&b
			\end{pmatrix}$ is invertible}
		\\
		\text{\it and all its entries are nonzero.}
		\end{gathered}
	\end{equation}
For this consider the Jacobian matrix of the complete intersection $C(a,b,c,d)$
$$
J(a,b,c,d) = \begin{pmatrix}
	Y_t &X_t&0&0&2T\\
	0&0&Y_{t-1}&X_{t-1}&2T\\
	a&b&c&d&0
\end{pmatrix}.
$$
It is easy to verify that the conditions \eqref{abcd} ensure that the rank of $J(a,b,c,d)$ is three at every point of $C(a,b,c,d)$. 
 
 Thus the hyperplane sections $C(a,b,c,d)$ with $(a,b,c,d)$ subject to \eqref{abcd} provide smooth curves of genus one. The four distinguished points on $C(a,b,c,d)$ are obtained as the intersection with the Lagrangian cycle
 $H_0={\bf H^{1,0}}(PG_3) \bigcap (T=0)$, that is, the divisor
 $$
 Z(a,b,c,d) =C(a,b,c,d) \bigcap (T=0).
 $$
 
	To complete the argument we need to show how quantum invariants manufacture the quadruples $(a,b,c,d)$ subject to \eqref{abcd}. For this fix $s\in I^-$. The quantum invariants give us maps
	$$
	\begin{gathered}
	\tau^+_{\rho_c(\xi\otimes\phi)}(s,t-j):Hom(\CC[q_{s,t-j},q^{-1}_{s,t-j}],\CC) \longrightarrow Hom(P^{t-j}, P^s),
	\\
	\tau^-_{\rho_c(\xi\otimes\phi)}(s,t-j):Hom(\CC[q_{s,t-j},q^{-1}_{s,t-j}],\CC) \longrightarrow Hom(P^s,P^{t-j}),
	\end{gathered}
$$
for $j=0,1,2$; in the sequel we omit the reference to $\rho_c(\xi\otimes\phi)$ to simplify the notation. We collect the maps in the following diagram
$$
\xymatrix{
P^t \ar@/_1pc/[rd]  & & \\
&P^s \ar@/_1pc/[lu] \ar@/_1pc/[ld] \ar@/_1pc/[rr]&&P^{t-2} \ar@/_1pc/[ll]\\
P^{t-1}\ar@/_1pc/[ru]&&
}
$$
 Compose the arrows of the diagram to obtain endomorphisms of the spaces at the nodes and then assign them to the four edges 
$e^-_t, e^+_{t-1}, e^-_{t-1}, e^+_{t-2}$ of the graph $PG_3$ in \eqref{PG3qin}:
$$
\begin{gathered}
e^-_t \rightarrow A^-_t:=\tau^-(s,t)(\bullet) \tau^{+}(s,t-1)(\bullet)\tau^-(s,t-1)( \bullet) \tau^+(s,t)( \bullet):P^t \longrightarrow P^t,
\\
e^+_{t-1}\rightarrow A^+_{t-1}:=\tau^-(s,t-1)( \bullet)\tau^+(s,t) (\bullet)\tau^-(s,t)(\bullet) \tau^+(s,t-1)(\bullet): P^{t-1}\longrightarrow P^{t-1},
\\
e^-_{t-1} \rightarrow A^-_{t-1}:=\tau^-(s,t-1)(\bullet) \tau^{+}(s,t-2)(\bullet)\tau^-(s,t-2)( \bullet) \tau^+(s,t-1)( \bullet):P^{t-1} \longrightarrow P^{t-1},
\\
e^+_{t-2} \rightarrow A^+_{t-2}:= \tau^-(s,t-2)(\bullet) \tau^{+}(s,t-1)(\bullet)\tau^-(s,t-1)(\bullet) \tau^+(s,t-2)( \bullet): P^{t-2} \longrightarrow P^{t-2};
\end{gathered}
$$
the logic of the above compositions is to go through the shortest loop of the diagram based at the tale of each edge and passing through its head. 
Thus with every triple of linear functionals on the spaces of Laurent polynomials
$$
(P,Q,R)\in \prod^2_{j=0} Hom(\CC[q_{s,t-j},q^{-1}_{s,t-j}],\CC)
$$
we obtain endomorphisms
$$
A^-_t(P,Q), A^+_{t-1}(P,Q),  A^-_{t-1}(Q,R), A^+_{t-2}(Q,R),
$$
labeling the four slanted edges of the graph $PG_3$ in \eqref{PG3qin}. 
We are now ready to have the quadruples. Define the map
$$
E(s,t,\bullet):\prod^2_{j=0} Hom(\CC[q_{s,t-j},q^{-1}_{s,t-j}],\CC) \longrightarrow M_2(\CC)
$$
which assigns to each triple $(P,Q,R)$ in the domain the matrix
$$
\begin{pmatrix}exp(tr(A^-_t(P,Q)))& exp(tr(A^-_{t-1}(Q,R)))\\
 exp(tr(A^+_{t-2}(Q,R)))&	exp(tr(A^+_{t-1}(P,Q))) 
\end{pmatrix}.
$$
Multiplying the first entry by the additional parameter $z \in \CC^{\times}$ gives the matrices
$$
E_z(s,t, P,Q,R)=\begin{pmatrix}zexp(tr(A^-_t(P,Q)))& exp(tr(A^-_{t-1}(Q,R)))\\
	exp(tr(A^+_{t-2}(Q,R)))&	exp(tr(A^+_{t-1}(P,Q))) 
\end{pmatrix}
$$
with the determinant
$$
det (E_z(s,t, P,Q,R))=zexp(tr(A^-_t(P,Q)+A^+_{t-1}(P,Q))) -exp(tr(A^-_{t-1}(Q,R)+A^+_{t-2}(Q,R))).
$$
This is nonzero unless
$$
\begin{gathered}
z=\frac{exp(tr(A^-_{t-1}(Q,R)+A^+_{t-1}(Q,R)))}{exp(tr(A^-_{t-1}(Q,R)+A^+_{t-2}(Q,R)))}
\\
=exp(tr(A^-_{t-1}(Q,R)+A^+_{t-1}(Q,R)-A^-_{t-1}(Q,R)-A^+_{t-2}(Q,R))).
\end{gathered}
$$
We denote that value $z_{P,Q,R}$. Thus the matrices $E_z (s,t,P,Q,R)$ for all nonzero values $z\neq z_{P,Q,R}$ are subject to the smoothness criterion \eqref{abcd}. Hence the linear forms
$$
\begin{gathered}
L_z(P,Q,R):=
zexp(tr(A^-_t(P,Q))) X_t +exp(tr(A^+_{t-1}(P,Q)))Y_t
\\
 +exp(tr(A^-_{t-1}(Q,R))) X_{t-1} + exp(tr(A^+_{t-2}(Q,R)))Y_{t-1}
\end{gathered}
$$
define smooth hyperplane sections
$$
C_z(P,Q,R):={\bf H^{1,0}}(PG_3) \bigcap (L_z(P,Q,R)=0), \,\forall z\neq z_{P,Q,R}, 0. 
$$
We now define the map
$$
\gamma_{P,Q,R}:\PP^1 \setminus \{\infty,0,z_{P,Q,R}\}=\CC^{\times} \setminus \{z_{P,Q,R}\} \longrightarrow {\mathfrak M}^4_1
$$
by sending $z\in \CC^{\times} \setminus \{z_{P,Q,R}\}$ to the point of ${\mathfrak M}^4_1$ corresponding to the pair
$$
(C_z(P,Q,R), C_z(P,Q,R) \cdot H_0).
$$
The map in the proposition is obtained by polarizing the $q$-parameters and the functionals, that is, we set $q_{s,t-j}=q$, for all values of $j$, and
$ P=Q=R=F \in Hom(\CC[q,q^{-1}],\CC)$.
\end{pf}

Heuristically, the above result tells us how to connect the ends of an `open' quiver of type $PG_m$ via the families of elliptic curves. Of course, the statement can be generalized to perfect matchings having triples
$(C^0_{I^+},C^{\pm}_I,C^0_{I^+})$ with the intervals $I^{\pm}$ having any length $k\geq 3$ and then creating families of hyperplane sections of the nonabelian Dolbeault space ${\bf H^{1,0}}(PG_k)$. Those are singular Calabi-Yau varieties of dimension $(k-2) $. For example, for $k=4$ we will find families of $K3$-surfaces in Example \ref{exDolbl=4}. All this clearly connects the refined IVHS invariants of a canonical curve of genus $g$ with moduli spaces of Calabi-Yau varieties and their periods.

 \subsection{Number theoretic connections}\label{rem:numth}
 	 This subsection is of speculative nature: it suggests a possible link of structures emerging from the refinement of IVHS with number theory. At this stage it is very sketchy and more elaborate discussion will appear elsewhere. The basic point is that number fields can be used to build the extensions of the representations of the quiver $PG_l$ to the ones of $\widehat{PG}_l$.
 	 
 	 Let $K$ be a number field, a finite field extension of the field of rational numbers $\mathbb{Q}$. Let $\mathfrak{O}_K$ be its ring of integers; assume that $K$ is not purely real and consider two complex embeddings of $K$
 	$$
 	\sigma, \sigma': K\hookrightarrow \CC;
 	$$
 	in such a situation one could fix $a\in \mathbb{N}$ and consider the subset 
 	$$
 	S_{\sigma,\sigma'} (a):=\{ (\sigma(s),\sigma'(s) ) \in \CC^{\times}\times \CC^{\times} | s\in \mathfrak{O}_K,\, |N_{K} (s)|=a \},
 	$$
 	where $N_K$ is the norm on $K$ over ${\mathbb{Q}}$. It is well known that the set of {\rm nonassociate} elements of  $\mathfrak{O}_K$ of a given norm is at most finite, see \cite{B-S}. With the understanding that $s$ in the definition of $S_{\sigma,\sigma'} (a)$ runs through the set of distinct nonassociate solutions of the equation
 	$$
 	|N_{K} (s)|=a,
 	$$
 	we attach to the data $(K,\sigma,\sigma',a)$ a distinguished finite subset  $S_{\sigma,\sigma'} (a)$ of $\CC^{\times}\times \CC^{\times}$.
 	Now we connect this with the family of representations of $PG_l$
 	$$
 	\rho_c: \OO^{\times}_{{\mathfrak{L}}_l}(-1) \longrightarrow \mathfrak{Reps}(PG_l)
 	$$
 	over the stratum ${\mathfrak{L}}_l$. In Corollary \ref{cor:rhoc-ext} we have learned to extend such families to the representations of the quiver $\widehat{PG}_l$. In particular, the algebraic torus $\CC^{\times}\times \CC^{\times}$ constitutes a set of natural parameters for such extensions. Namely, for every point $\xi\otimes \phi$ in $\OO^{\times}_{{\mathfrak{L}}_l}(-1)$ lying over a point $([\xi],[\phi]) \in {\mathfrak{L}}_l$ we have the map
 	$$
 	\mu_{\rho_c(\xi\otimes\phi)}:\CC^{\times}\times \CC^{\times} \longrightarrow P([\xi],[\phi])=P^-([\xi],[\phi])\oplus P^+([\xi],[\phi])
 	$$
 	 which assigns labels to additional edges $e^-_0$ and $e^+_{l-1}$ of $\widehat{PG}_l$; recall the notation
 	 $$
 	 \begin{gathered}
 	 P^-([\xi],[\phi])=Hom (P^0([\xi],[\phi]),P^{l-1}([\xi],[\phi])),
 	 \\
 	   P^+([\xi],[\phi])=Hom (P^{l-1}([\xi],[\phi]),P^{0}([\xi],[\phi])).
 	 \end{gathered}
  $$
 	  The values 
 	  $$
 	  \mu_{\rho_c(\xi\otimes\phi)}(z,z')=\mu^-_{\rho_c(\xi\otimes\phi)}(z)+\mu^+_{\rho_c(\xi\otimes\phi)}(z'), \,\,\forall (z,z')\in \CC^{\times}\times \CC^{\times}, 
 	  $$
 	 determine the representations $\widehat{\rho}_c(z,z') (\xi\otimes \phi)$ of $\widehat{PG}_l$ extending $\rho_c (\xi\otimes \phi)$ and receiving the labels $\mu^-_{\rho_c(\xi\otimes\phi)}(z)$ and $\mu^+_{\rho_c(\xi\otimes\phi)}(z')$ respectively on edges $e^-_0$ and $e^+_{l-1}$. In particular, the subset $S_{\sigma,\sigma'} (a)$ produces the family of representations
 	 $$
 	 \widehat{\rho}_c(S_{\sigma,\sigma'}(a)): \OO^{\times}_{{\mathfrak{L}}_l}(-1) \longrightarrow \mathfrak{Reps}(\widehat{PG}_l)
 	 $$
 	 of $\widehat{PG}_l$ extending the family $\rho_c$ and parametrized by $S_{\sigma,\sigma'}(a)$. More precisely, for every $\xi\otimes \phi$ in
 	 $\OO^{\times}_{{\mathfrak{L}}_l}(-1)$ we have the collection
 	 $$
 	 \widehat{\rho_c}(S_{\sigma,\sigma'}(a))(\xi\otimes \phi):=\{\widehat{\rho}_c (\sigma(s), \sigma'(s))(\xi\otimes \phi)| s\in S_{\sigma,\sigma'}(a)\}
 	 $$
 	  of representations of the quiver $\widehat{PG}_l$ extending $\rho_c(\xi\otimes \phi)$ and  
 	 parametrized by the algebraic nonassociate integers $s\in {\mathfrak{O}_K}$ subject to $|N_K(s)|=a$.
 	  
 	 The situation in Proposition \ref{pro:lrootsreps} is the special case when $K=\mathbb{Q}(\zeta)$, where $\zeta$ is a primitive $l$-th root of unity, a number field naturally appearing in our considerations.
 	 Of course, if our curve $C$ is defined over some number field $K$ with complex embeddings, that number field provides a distinguished choice for the above constructions as well.  
 	In the next section we will see other instances of number fields appearing in our constructions.

 	\vspace{0.2cm}
 	It is well known that the number fields have cousins - the function fields of curves defined over finite fields, see \cite{W}. This gives another way to build representations of $\widehat{PG}_l$ extending $\rho_c$.  
 	To be specific, let us briefly recall some well known facts, see \cite{KS} and the literature therein for details:
 	
 	 let $p$ be a prime number and let $\Gamma$ be a smooth, proper, geometrically connected curve of genus $\gamma \geq 1$ over the finite field $\mathbb{F}_q$ of characteristic $p$; then one has the Artin zeta function attached to $\Gamma$
 	$$
 	Z(\Gamma/\mathbb{F}_q,T):=exp \left(\sum^{\infty}_{n=1}\frac{N_nT^n}{n}\right),
 	$$
 	where $N_n=\text{\rm the number of points of $\Gamma$ over the field $\mathbb{F}_{q^n}$}$; it is known that this is a rational function of $T$:
 	$$
 	Z(\Gamma/\mathbb{F}_q,T)=\frac{P_{\Gamma}(T)}{(1-T)(1-qT)},
 	$$
 	where $P_{\Gamma} (q)$ is a polynomial of degree $2\gamma$ with integer coefficients; furthermore, writing the factorization of that polynomial
 	$$
 	P_{\Gamma}(T)=\prod^{2\gamma}_{j=1} (1-q_j T),
 	$$
 	where the reciprocals $q^{-1}_j$ are the roots of $P_{\Gamma}(T)$ in $\CC$, one knows that they are situated on the circle centered at $0$ of radius ${q}^{-\HA}$, that is,
 	$|q_j|=q^{\HA}$, for all $j$. Thus the correspondence
 	$$
 	\Gamma/\mathbb{F}_q \longrightarrow \{P_{\Gamma}(T), \{q_j\}\}.
 	$$
 	Returning to our considerations, we see that the polynomial $P_{\Gamma} (T)$ or its Laurent polynomial version $\displaystyle{\frac{P_{\Gamma} (T)}{T^{\gamma}}}$  (resp., the set  $Q_{\Gamma}=\{q_j\}$) can be fed into the maps $\tau^{\pm}_{\rho_c}$ in \eqref{Lpoly-Ppm} (resp. $\mu_{\rho_c}$ in Corollary \ref{cor:rhoc-ext}) to obtain the representations of $\widehat{PG}_l$ extending the representations $\rho_c$ and labeled by the data $\{P_{\Gamma}(T), \{q_j\}\}$.
  
\subsection{Zeta functions of zig-zag paths}
To the stratum ${\mathfrak{L}}_l$ parametrizing the algebraic K\"{a}hler structures of fixed length $l \geq 3$ we have attached the graph $\widehat{PG}_l$. This comes equipped with zig-zag paths
$$
Z_{0,-},\,Z_{0,+},\, Z_{-,+}.
$$
We have seen how to use these paths to decorate {\it all} edges of $\widehat{PG}_l$ with linear maps between the summands of the direct sum decompositions
$$
W_{\rho} =\bigoplus^{l-1}_{j=0} P^j ,
$$
 for every representation $\rho$ of the quiver $PG_l$.  This is provided in particular by the map
 $$
 \mu_{\rho}=(\mu^-_{\rho}, \mu^+_{\rho}): \CC^{\times}\times \CC^{\times} \longrightarrow P^-\oplus P^+  
 $$
 in Proposition \ref{pro:mu-map}. Inspired by the finite field case discussed in the previous subsection, we put the above ingredients together to define a sort of Artin zeta functions attached to points of
 $ \OO^{\times}_{{\mathfrak{L}}_l} (-1) \times (\CC^{\times})^3$. This also will be a complement
 to the zeta function $\mathfrak{z}(\beta(\phi,\bullet),t)$ attached to the last term $W^l_{\xi}([\phi])$ of $([\xi],[\phi])$-filtration, see the discussion and the formula in the end of \S5, page \pageref{Zetabeta}.
 
 \vspace{0.2cm}
 The construction goes as follows. We start with the family of representations
 $$
 \rho_c: \OO^{\times}_{{\mathfrak{L}}_l} (-1) \longrightarrow {\mathfrak{Reps}}(PG_l)
 $$
 of $PG_l$ which attaches the representation
 $$
 \rho_c(\xi\otimes\phi)=\{\alpha^{t,s}_c(\xi,\phi):P^s([\xi],[\phi])\longrightarrow P^t([\xi],[\phi])\}
 $$
 to every point $\xi\otimes\phi$ in $\OO^{\times}_{{\mathfrak{L}}_l} (-1)$ lying over a point $([\xi],[\phi]) \in {\mathfrak{L}}_l$; recall: the representation space
 $$
 W_{\rho_c(\xi\otimes\phi)} =W_{\xi}/W^l_{\xi}([\phi])=\bigoplus^{l-1}_{j=0} P^j ([\xi],[\phi])
 $$
 and the maps $\alpha^{t,s}_c(\xi,\phi)$ are the $(t,s)$-blocks of the map
 $\alpha^{(2)}_{\xi,c}(\phi,\bullet)$ with respect to the orthogonal decomposition of $W_{\xi}/W^l_{\xi}([\phi])$ above.
   
 We go back to the maps $C^{\pm}_{\rho_c(\xi\otimes\phi)} (q,\bullet)$ in \eqref{C+-maps-q}:
 $$
  C^{\pm}_{\rho_c(\xi\otimes\phi)} (q,\bullet): P^{\pm}([\xi],[\phi]) \longrightarrow {\mathfrak{g}}^{(0)}[q,q^{-1}],
  $$
  where ${\mathfrak{g}}^{(0)}$ denotes the Lie subalgebra of $End(W_{\rho_c(\xi\otimes\phi)})$ of grading preserving endomorphisms of
  the graded vector space $W_{\rho_c(\xi\otimes\phi)}$. 
  In Corollary \ref{cor:rhoc-ext} we have the map
  $$
  \mu_{\rho_c(\xi\otimes\phi)}:\CC^{\times}\times \CC^{\times} \longrightarrow P^{-}([\xi],[\phi])\oplus P^{+}([\xi],[\phi])
  $$
  attaching the labels to two additional edges of $\widehat{PG}_l$ and thus defining the representation the family of representations
  $$
  \widehat{\rho}_c(\bullet,\bullet)(\xi\otimes\phi): \CC^{\times}\times \CC^{\times} \longrightarrow {\mathfrak{Reps}}(\widehat{PG}_l)
  $$
 extending $\rho_c(\xi\otimes\phi)$; the projections of $\mu_{\rho_c(\xi\otimes\phi)}$ onto the summands $P^{\pm}([\xi],[\phi])$ will be denoted $\mu^{\pm}_{\rho_c(\xi\otimes\phi)}$.
 Here we will be concerned with the compositions
   $$
  \begin{gathered}
  	C^{-}_{\rho_c(\xi\otimes\phi)} (q,\bullet) \circ \mu^-_{\rho_c(\xi\otimes\phi)}:\CC^{\times} \longrightarrow {\bf \mathfrak{g}^{(0)}}[q,q^{-1}],
  	\\
  	C^{+}_{\rho_c(\xi\otimes\phi)} (q,\bullet) \circ \mu^+_{\rho_c(\xi\otimes\phi)} :\CC^{\times} \longrightarrow {\bf \mathfrak{g}^{(0)}}[q,q^{-1}]. 
  \end{gathered}
  $$
  Evaluate those at a point $(u,v) \in \CC^{\times}\times \CC^{\times}$ we obtain 
 \begin{equation}\label{(u,v)-formula}
 \CC^{\times} \ni u \mapsto  C^{-}_{\rho_c(\xi\otimes\phi)} (q, \mu^-_ {\rho_c(\xi\otimes\phi)} (u))=\sum^{l-1}_{j=0} A^-_j (\rho_c(\xi\otimes\phi);u)q^{(-1)^j j} \in {\bf \mathfrak{g}^{(0)}}[q,q^{-1}],
\end{equation}
 where  $A^-_j (\rho_c(\xi\otimes\phi);u):= (C^{-})^j_{\rho_c(\xi\otimes\phi)} (q, \mu^-_{\rho_c(\xi\otimes\phi)}  (u))$ is the component of degree $j$. Analogously, for the second
 \begin{equation}\label{(u,v)+formula}
 	\CC^{\times} \ni v \mapsto  C^{+}_{\rho_c(\xi\otimes\phi)} (q, \mu^+_{\rho_c(\xi\otimes\phi)}(v))=\sum^{l-1}_{j=0} A^+_j (\rho_c(\xi\otimes\phi);v)q^{(-1)^j j} \in {\bf \mathfrak{g}^{(0)}}[q,q^{-1}],
 \end{equation}
 where  $A^+_j (\rho_c(\xi\otimes\phi);v):= (C^{+})^j_{\rho_c(\xi\otimes\phi)} (q, \mu^+_{\rho_c(\xi\otimes\phi)} (v))$ is the component of degree $j$.
 
 \vspace{0.2cm}
Next we recall: given an endomorphism $f: V \longrightarrow V$ of a finite dimensional vector space $V$, following the standard construction one defines the zeta function of $f$
 $$
 {\bf \mathfrak{z}} (f, t):=exp \left(\sum^{\infty}_{n=1} \frac{tr(f^n) t^n}{n}\right),
 $$
 where $t$ is a formal variable. It is again a standard fact to express the above in a closed form
 $$
 {\bf \mathfrak{z}} (f, t) = (det({\bf 1}-tf))^{-1},
 $$
 where ${\bf 1}$ stands for the identity isomorphism of $V$. In case the vector space $V$ is graded
 $$
 V=\bigoplus^N_{m=0} V_m
 $$
 and $f$ is grading preserving endomorphism, that is, the restriction $f_m$ of $f$ to $V_m$ takes values in $V_m$, one attaches to $(f,V)$ the product of zeta functions
 $$
 {\bf \mathfrak{z}} (f, t):=\prod^N_{m=0} ({\bf \mathfrak{z}} (f_m, t))^{(-1)^m}.
 $$
 Using the closed form of the zeta function we deduce
 $$
{\bf \mathfrak{z}} (f, t)=\prod^N_{m=0} ({\bf \mathfrak{z}} (f_m, t))^{(-1)^m}=\frac{\displaystyle{\prod_{2i+1 \in [1, N]}} det(1-tf_{2i+1})}{\displaystyle{\prod_{2i \in [0, N]}} det(1-tf_{2i})}
$$

Applying this to the case of the graded vector spaces
$$
W_{\xi}/W^l_{\xi}([\phi])=\bigoplus^{l-1}_{j=0} P^j ([\xi],[\phi])
$$
and grading preserving endomorphisms
$$
\begin{gathered}
C^{-}_{\rho_c(\xi\otimes\phi)} (q, \mu^-_{\rho_c(\xi\otimes\phi)} (u))=\sum^{l-1}_{j=0} A^-_j (\rho_c(\xi\otimes\phi);u)q^{(-1)^j j}, \\
C^{+}_{\rho_c(\xi\otimes\phi)} (q, \mu^+_{\rho_c(\xi\otimes\phi)} (v))=\sum^{l-1}_{j=0} A^+_j (\rho_c(\xi\otimes\phi);v)q^{(-1)^j j},
\end{gathered}
$$
we obtain two zeta functions
\begin{equation}\label{zeta-formulas}
	\begin{gathered}
	{\bf \mathfrak{z^-}}_{\rho_c(\xi\otimes\phi)}( u,q,t):=	{\bf \mathfrak{z}}(C^{-}_{\rho_c(\xi\otimes\phi)} (q, \mu^-_{\rho_c(\xi\otimes\phi)} (u)), t)
	\\
	 =
	\frac{\displaystyle{\prod_{2i+1 \in [1, l-1]}} det(1-tq^{-2i-1}A^-_{2i+1} (\rho_c(\xi\otimes\phi);u))}{\displaystyle{\prod_{2i \in [0, l-1]}} det(1-tq^{2i}A^-_{2i}(\rho_c(\xi\otimes\phi);u))},
	\\
	\\
	{\bf \mathfrak{z^+}}_{\rho_c(\xi\otimes\phi)}(v,q,t):=	{\bf \mathfrak{z}}(C^{+}_{\rho_c(\xi\otimes\phi)} (q, \mu^+_{\rho_c(\xi\otimes\phi)} (v)), t)
	\\
	 =
	\frac{\displaystyle{\prod_{2i+1 \in [1, l-1]}} det(1-tq^{-2i-1}A^+_{2i+1}(\rho_c(\xi\otimes\phi);v))}{\displaystyle{\prod_{2i \in [0, l-1]}} det(1-tq^{2i}tA^+_{2i}(\rho_c(\xi\otimes\phi);v))}.
\end{gathered}
\end{equation} 
These are viewed as functions of $(\xi\otimes \phi; u,v,q) \in \OO^{\times}_{{\mathfrak{L}}_l}(-1) \times (\CC^{\times})^3$ with values in the field $\CC(t)$ of rational functions of $t$. They encode the spectral properties of natural families of representations of the quiver $\widehat{PG}_l$ on the IVHS side.

We will now use the ingredients of these zeta functions to construct similar objects but on the {\it geometric side} of the two part diagram
\begin{equation}\label{two-partdiag-zeta} 
\xymatrix@C=18pt{
	&{\bf P}\ar_{p_1}[dl]\ar^{p_2}[dr]&&{\cal Z}_C\ar_{q_1}[dl]\ar^{q_2}[dr]&&\\
	\PP(H^1 (\Theta_{C}))&&\PP(\HKC)&&C \ar@{^{(}->}[r]& \PP(\HKC^{\ast})
}
\end{equation}
which has guided our thinking. Namely, under some `genericity' assumptions we will construct conjugacy classes of the fundamental group of the complement of the dual variety $C^{\vee}$ in $\PP(\HKC)$. Those deliver conjugacy classes of the automorphism group of the
fibres of $q_1$ via the monodromy representation; these classes are used to attach Artin-type zeta functions to finite dimensional representations of the fundamental group of the complement of $C^{\vee}$. 

We begin by discussing the genericity assumption. Recall the orthogonal decomposition
$$
\HKC=W_{\xi} \oplus (W_{\xi})^{\perp}=\left(\bigoplus^{l}_{s=0}P^s([\xi],[\phi])\right)\oplus (W_{\xi})^{\perp}
$$
for $([\xi],[\phi])\in {\mathfrak{L}}_l$ and varying in $C^{\infty}$ manner with $([\xi],[\phi])$. We impose the condition
\begin{equation}\label{gencond-zeta}
	\begin{gathered}
\text{\it there is a Zariski dense open subset ${\mathfrak{L}}^{\circ}_l$ of ${\mathfrak{L}}$}
\\
\text{\it such that $\PP(P^s([\xi],[\phi])) \nsubseteq C^{\vee},\,\,\forall s\in[0,l-1]$ and $([\xi],[\phi])\in {\mathfrak{L}}^{\circ}_l$.}
\end{gathered}
\end{equation}
\begin{rem}\label{rem:gencond-zeta}
	The condition is obviously satisfied if the subspaces $\PP(P^s([\xi],[\phi]))$'s, viewed as linear subsystems of the canonical system $|K_C|$, are base point free. So the failure of \eqref{gencond-zeta} means that the summands of $([\xi],[\phi])$-decomposition are geometric: the geometry is delivered by the base loci of $\PP(P^s([\xi],[\phi]))$. 
\end{rem}

The quiver representation $\rho_c (\xi\otimes \phi)=\{(\alpha^{t,s}_c(\xi,\phi): P^s([\xi],[\phi]) \longrightarrow P^t([\xi],[\phi]))\}$ of $PG_l$ gives the maps
$$
\alpha^{s-1,s}(\xi,\phi): P^s([\xi],[\phi]) \longrightarrow P^{s-1}([\xi],[\phi]), \,\, \forall s\in [1,l-1];
$$
recall those are independent of the trace parameter $c$ and are injective, see Proposition \ref{pro:xi-phi-filt},2). In particular, the dimensions of the summands $P^s$'s form the monotone sequence
$$
h^{s}=dim P^s([\xi],[\phi]) \leq h^{s-1}=dim P^{s-1}([\xi],[\phi]),\,\,\forall s\in[1,l-1].
$$
Our next assumption is that we are {\it not} in the case of the maximal ladder, in other words we assume
\begin{equation}\label{hrbiger2}
	\text{$h^s \geq 2$ for some $s \in [0,l-1]$.}
\end{equation} 
(Recall: the case of maximal ladders was studied in \S4, see Proposition \ref{pro:maxladder}, Corollary \ref{cor:Ml-complement}). Let $s_0$ be the first (in the decreasing order) summand $P^{s_0}$ of the orthogonal decomposition having dimension bigger than one. By the monotone property of dimensions $h^s$'s we have
$$
h^s \geq 2, \forall s\in [0,s_0].
$$
We now turn to the endomorphisms $A^{\pm}_s (\rho_c(\xi\otimes\phi);\bullet)$ appearing in the formulas \eqref{zeta-formulas}. Here our genericity assumptions will be the semi-simplicity of endomorphisms in the algebra generated
by $A^{\pm}_s (\rho_c(\xi\otimes\phi);\bullet)$, for $ s\in [0,s_0]$ and $\bullet=w\in \CC^{\times}$. Precisely, for every point $([\xi],[\phi])\in {\mathfrak{L}}_l$ define
$$
{\bf A}^s([\xi],[\phi])=\txt{ \it the subalgebra of $End(P^s([\xi],[\phi]))$\\ \it generated by $A^{\pm}_s (\rho_c(\xi\otimes\phi);w), \forall w\in \CC^{\times}$},
$$ 
and assume
\begin{equation}\label{gencond2-zeta}
	\txt{\it a general $a\in {\bf A}^s([\xi],[\phi]))$ is {\rm regular, semisimple},} 
\end{equation}
that is, a general endomorphism of ${\bf A}^s([\xi],[\phi])$ is diagonalizable with distinct eigen values. 
Our final requirement will be putting the conditions \eqref{gencond-zeta} and \eqref{gencond2-zeta} together:
\begin{equation}\label{gencond3-zeta}
	\txt{\it there is a Zariski dense open subset ${\mathfrak{L}}^{\circ\circ}_l$ of ${\mathfrak{L}}_l$ such that
		\\
		\it for all
		$([\xi],[\phi])\in {\mathfrak{L}}^{\circ\circ}_l$ the following hold:
	\\
(1)	 $\PP(P^s([\xi],[\phi])) \nsubseteq C^{\vee},\,\,\forall s\in[0,l-1]$, 
 \\
 \\
 (2) \it for all $s\in [0,l-1]$, a general $a\in {\bf A}^s([\xi],[\phi])$ is regular, semisimple
 \\
 \it and, in addition,
   the eigen spaces of $a$ are points in $\PP(P^s([\xi],[\phi]))$ not lying in $C^{\vee}$,
  \\
  \\
  (3) \it there is $s_0 \in [0,l-1]$ such that $h^{s_0}([\xi],[\phi])\geq 2$.
}
\end{equation}
With the above conditions in mind we construct closed loops in the complement of the dual variety $C^{\vee}$.
\begin{lem}\label{lem:loops}
	Assume the condition \eqref{gencond3-zeta} holds. Fix a point $([\xi],[\phi])$ in ${\mathfrak{L}}^{\circ\circ}_l$ and, for $s\geq s_0$, a regular semi-simple endomorphism
	$a\in {\bf A}^s([\xi],[\phi])$ subject to the condition $(2)$ of \eqref{gencond3-zeta}. Fix an order of the eigen spaces of $a$. Then one obtains a loop $\gamma_a$ in $\PP(\HKC)\setminus C^{\vee}$ passing through the points of $\PP(\HKC)$ corresponding to the eigen spaces of $a$. 
\end{lem}
\begin{pf}
	Let $([\xi],[\phi])$, $s$ and $a$ be as in the statement of the lemma. Let $\{v_1,\ldots,v_{h^s}\}$ be a basis of $P^s([\xi],[\phi])$ formed by eigen vectors of $a$ and ranged according to the fixed order. This gives the ordered collection of points
	$$
	[v_1],\ldots,[v_{h^s}]
	$$
	in $\PP(P^s([\xi],[\phi]))$ and outside the dual variety $C^{\vee}$. Draw the lines $L_i:=\langle [v_{i}],[v_{i+1}]\rangle$, for $i=1,\ldots,h^s-1$. The union of those lines form reducible curve whose Dynkin diagram is $A_{h^s-1}$. Complete it to a cycle by the line $L_{h^s}=\langle [v_{h^s}],[v_{1}]\rangle$ joining the last point $[v_{h^r}]$ to $[v_1]$.  The union 
	$$
	L(a)=\bigcup^{h^s}_{i=1} L_i
	$$
	is a cycle of lines in  $\PP(P^s([\xi],[\phi]))$ with Dynkin diagram $\widetilde{A}_{h^s-1}$. We now draw a loop in $L(a)$ as asserted in the lemma.
	For this in each line $L_i$ choose a simple path $\gamma_i$ joining $[v_i]$ to $[v_{i+1}]$ and lying outside $C^{\vee}$: start with the distance minimizing geodesic $\gamma'_i$ connecting $[v_i]$ to $[v_{i+1}]$;
	$L_i$  intersects the dual variety at finite number of points, if none of those points lie on the geodesic, then we take $\gamma_i=\gamma'_i$;
	if $\gamma'_i$ passes through some of those points, modify $\gamma'_i$ by an arc of a small circle in the complement of $C^{\vee}$ and going around each point of $C^{\vee}$ encountered by $\gamma'_i$; set $\gamma_i$ to be this modified path; the following drawing is an illustration:
	$$
	\begin{tikzpicture}
		\draw (2,0)
		arc [start angle=0, end angle =180, x radius=2, y radius=1];
		\filldraw [black] (2,0) circle [radius=2pt] node[anchor= west]{$v_{i+1}$};
		\filldraw [black] (-2,0) circle [radius=2pt] node[anchor=east]{$v_i$};
		\filldraw [red] (-1,0.85) circle [radius=2pt];
		\draw  [blue] (-0.8,0.9) arc [start angle=0, end angle=203, radius=6pt];
		\filldraw [red] (0,1) circle [radius=2pt];
		\draw  [blue] (0.2,0.98) arc [start angle=0, end angle=180, radius=6pt];
		\filldraw [red] (1,0.85) circle [radius=2pt];
			\draw  [blue] (1.2,0.8) arc [start angle=-10, end angle=155, radius=6pt];
	\end{tikzpicture} 
$$
the path $\gamma'_i$ joining $v_i$ to $v_{i+1}$ is in black, the red dots are the points of $C^{\vee}$ lying on the path, the blue arcs around each of those points are the modifications needed to obtain the path $\gamma_i$.

Define
	the loop
	$$
	\gamma_a:[0,1] \longrightarrow L(a)\setminus C^{\vee} \subset \PP(\HKC)\setminus C^{\vee}
	$$
	by concatenation of the paths $\gamma_1, \ldots,\gamma_{h^s}$. This gives a loop based at $[v_1]$ and passing through all other points, the eigen spaces of $a$.
\end{pf}

The morphism $q_1:{\cal Z}_C \longrightarrow \PP(\HKC)$ in \eqref{two-partdiag-zeta} is a finite morphism
of degree $(2g-2)$ ramified over the dual variety $C^{\vee}$: the fibre of 
$q_1$ over a point $[\phi]\in \PP(\HKC)$ is identified via the projection $q_2$ with the divisor $Z_{\phi} =(\phi=0)$ on $C$. Thus over the complement
of $C^{\vee}$ we have the unramified covering
$$
q_1: {\cal Z}_C \setminus q^{-1}(C^{\vee}) \longrightarrow \PP(\HKC)\setminus C^{\vee}=\left(C^{\vee}\right)^c.
$$
So if we fix a point $[v]\in \left(C^{\vee}\right)^c$, the fundamental group
$\pi_1(\left(C^{\vee}\right)^c,[v])$ acts on $Z_v$ by permutations, that is we have the monodromy representation
\begin{equation}\label{monodromy}
	\varrho_{{}_C}: \pi_1(\left(C^{\vee}\right)^c,[v]) \longrightarrow Aut(Z_{[v]})\cong S_{2g-2},
\end{equation}
where $S_{m}$ stands for the group of permutations on $m$ letters. 

Assume the condition \eqref{gencond3-zeta}. Over the Zariski open subset
${\mathfrak{L}}^{\circ\circ}_l$ we have the sheaf of algebras 
$$
{\bf A}^s_{{\mathfrak{L}}^{\circ\circ}_l}:= \bigcup_{([\xi],[\phi])\in {\mathfrak{L}}^{\circ\circ}_l} {\bf A}^s ([\xi],[\phi]) \longrightarrow {\mathfrak{L}}^{\circ\circ}_l.
$$
Denote by $\left({\bf A}^s_{{\mathfrak{L}}^{\circ\circ}_l}\right)^{rss}$ the subspace of ${\bf A}^s_{{\mathfrak{L}}^{\circ\circ}_l}$ formed by regular semi-simple endomorphisms. From Lemma \ref{lem:loops} we deduce the following.
\begin{pro}\label{pro:conj}
	Assume the condition \eqref{gencond3-zeta} and denote by 
	${\mathfrak{L}}^{\circ\circ}_l (s_0)$ the stratum of ${\mathfrak{L}}^{\circ\circ}_l$ of closed points $([\xi],[\phi])$  subject to $h^s([\xi],[\phi]) \geq 2$, for all $0\leq s\leq s_0$. Fix a point
	$[v_0]\in (C^{\vee})^c$ and let $\pi_1 ((C^{\vee})^c,[v_0])$ be the fundamental group of $(C^{\vee})^c$ based at $[v_0]$. Then for every $s \in [0,s_0]$ there is a correspondence between regular semisimple elements of ${\bf A}^s_{{\mathfrak{L}}^{\circ\circ}_l}$ and finite subsets
	of conjugacy classes of the fundamental group $\pi_1 ((C^{\vee})^c,[v_0])$:
	$$
	\left({\bf A}^s_{{\mathfrak{L}}^{\circ\circ}_l(r_0)}\right)^{rss} \ni a\mapsto C(a) \subset Conj(\pi_1 ((C^{\vee})^c,[v_0])),
	$$
	where $ Conj(\pi_1 ((C^{\vee})^c,[v_0]))$ denotes the set of conjugacy classes of $\pi_1 ((C^{\vee})^c,[v_0])$.
\end{pro}
\begin{pf}
	Let $s\in [0,s_0]$. For every point $([\xi],[\phi])$ in ${\mathfrak{L}}^{\circ\circ}_l (s_0)$ take a regular semisimple endomorphism $a$ in ${\bf A}^s ([\xi],[\phi])$. For an order $o$ of its eigen spaces fixed, Lemma \ref{lem:loops} produces a loop $\gamma^o_a$ in $(C^{\vee})^c$ passing through the points 
	\begin{equation}\label{v-o}
	[v_1],\ldots, [v_{h^s}]
\end{equation}
	of $(C^{\vee})^c$ corresponding to the eigen spaces of $a$ and ranged according the order $o$. Let $[v]$ be one of those points and view $\gamma^o_a$ as a path based at $[v]$. Let $\gamma([v_0],[v]) $ be a path in $(C^{\vee})^c$ connecting the base point $[v_0]$ to $[v]$. Denote by $[\gamma^o_a ([v_0])]$ the homotopy class
	of the path $\gamma^{-1}([v_0],[v])\gamma^o_a \gamma([v_0],[v])$. The conjugacy class $conj([\gamma^o_a ([v_0])])$ is independent of the choices made. A change of order of eigen spaces corresponds to a permutation of the points $[v_i]$'s in \eqref{v-o}. Thus the collection of the conjugacy classes $conj([\gamma^o_a ([v_0])])$ is labeled by the permutations $S_{h^s}$ and we obtain the correspondence
	$$
		\left({\bf A}^s_{{\mathfrak{L}}^{\circ\circ}_l(s_0)}\right)^{rss} \ni a\mapsto C(a):=\{conj([\gamma^o_a ([v_0])])\}_{o\in S_{h^s}} \subset Conj(\pi_1 ((C^{\vee})^c,[v_0]))
	$$
\end{pf}

The above proposition combined with the monodromy representation \eqref{monodromy} implies the following.
\begin{cor}\label{cor:conj}
	With the assumptions and notation of Proposition \ref{pro:conj}, the correspondence
	$$
	\left({\bf A}^s_{{\mathfrak{L}}^{\circ\circ}_l(s_0)}\right)^{rss} \ni a \mapsto C(a)\subset Conj(\pi_1 ((C^{\vee})^c,[v_0]))
	$$
	composed with the map
	$$
	Conj(\pi_1 ((C^{\vee})^c,[v_0])) \stackrel{\varrho_{{}_C}}{\longrightarrow} Conj(Aut(Z_{v_0}))
	$$
induced by the monodromy representation	\eqref{monodromy} determines a finite collection  of conjugacy classes in
	$ Aut(Z_{v_0})$, for every $s\in [0,s_0]$. That collection will be denoted $Conj^s_{{\mathfrak{L}}^{\circ\circ}_l(s_0)}$. 
\end{cor}
\begin{pf}
	Follows immediately from Proposition \ref{pro:conj} and the finiteness of the group $Aut(Z_{v_0})$. 
\end{pf}

Consider the constant sheaf $\underline{\CC}$ on ${\cal Z}_C$. The direct image $(q_{1})_{\ast} (\underline{\CC})$ gives a local system on $\left(C^{\vee}\right)^c$. It corresponds to the representation
$$
\varrho^{\CC}_{{}_C}:\pi_1(\left(C^{\vee}\right)^c,[v_0]) \longrightarrow {\bf GL}(H^0(\OO_{Z_{v_0}}))
$$ 
determined by the monodromy representation $\varrho_{{}_C}$ via the regular action of the group $Aut(Z_{v_0})$ on the space of functions on $Z_{v_0}$: 
$$
\varrho^{\CC}_{{}_C} (g) (f)=f\circ \varrho_{{}_C}(g^{-1}),\,\, \forall f\in H^0(\OO_{Z_{v_0}}),
$$
that is the representation $\varrho^{\CC}_{{}_C}$ factors through the monodromy representation as follows:
$$
\xymatrix{
\pi_1(\left(C^{\vee}\right)^c,[v_0]) \ar[r]^(.6){\varrho_{{}_C}}\ar@/_2pc/[rr]_{\varrho^{\CC}_{{}_C}}& Aut(Z_{v_0}) \ar[r]^(.4)R& {\bf GL}(H^0(\OO_{Z_{v_0}})),
}
$$ 
where $R$ denotes the regular representation of  $Aut(Z_{v_0})$ on $H^0(\OO_{Z_{v_0}})$.

Using the collection of conjugacy classes $Conj^s_{{\mathfrak{L}}^{\circ\circ}_l(s_0)}$ provided by Corollary \ref{cor:conj} define the following $s$-factor of the forthcoming Artin-type function:
\begin{equation}\label{r-factors}
	L^s ({\mathfrak{L}}_l,\varrho^{\CC}_{{}_C};t):=\prod_{\sigma\in Conj^s_{{\mathfrak{L}}^{\circ\circ}_l(s_0)}} det({\bf 1}_{2g-2}-t\varrho^{\CC}_{{}_C}(\sigma))^{-1},
\end{equation}
for every $s\in [0,s_0]$. The final output is the product of those $s$-factors 
\begin{equation}\label{L-geomside}
L({\mathfrak{L}}_l,\varrho^{\CC}_{{}_C};t)= \prod_{s \in [0,s_0]} L^s ({\mathfrak{L}}_l,\varrho^{\CC}_{{}_C};t).
\end{equation}
Thus we arrive to the following correspondence between zeta functions
${\mathfrak z}^{\pm}$ defined on the IVHS side of the diagram \eqref{two-partdiag-zeta} and the Artin-type functions $L({\mathfrak{L}}_l,\varrho^{\CC}_{{}_C};t)$ on the geometric side of that diagram.
\begin{pro}\label{pro:Lcorresp}
	Assume ${\mathfrak{L}}_l$ satisfies the conditions \eqref{gencond3-zeta}. Then the zeta functions ${\mathfrak z}^{\pm}$ defined in \eqref{zeta-formulas} in terms of the refined IVHS data give rise to the function $L({\mathfrak{L}}_l,\varrho^{\CC}_{{}_C};t)$ in \eqref{L-geomside} attached to the representation $\varrho^{\CC}_{{}_C}$ of the fundamental group of the complement of the dual variety $C^{\vee}$ in $\PP(\HKC)$. 
\end{pro}

\begin{rem}\label{rem:Lcorresp}
1)	The correspondence of Proposition \ref{pro:Lcorresp} is evocative of the Langlands philosophy: the zeta functions ${\mathfrak z}^{\pm}$ are determined in terms of the representations of certain reductive groups on the cohomological side of the diagram \eqref{two-partdiag-zeta} and we arrive to the function $L({\mathfrak{L}}_l,\varrho^{\CC}_{{}_C};t)$ on the side of the morphism $q_1:{\cal Z}_C \setminus q^{-1}_1 (C^{\vee}) \longrightarrow \left(C^{\vee}\right)^c$ which could be viewed as the Galois side of the correspondence. 

2) The `$s$-factors' $L^s({\mathfrak{L}}_l,\varrho^{\CC}_{{}_C};t)$ in \eqref{r-factors} could be viewed as local contributions in the following sense: in the proof of Lemma \ref{lem:loops} we really construct loops in the subspaces
$$
\PP(P^s) \setminus C^{\vee} \bigcap \PP(P^s),
$$
 for every $s\in [0,s_0]$; so the conjugacy classes in Proposition \ref{pro:conj} belong to the fundamental group $\pi_1(\PP(P^s) \setminus C^{\vee} \bigcap \PP(P^s) )$ and its representation associated to the restriction of the local system $(q_{1})_{\ast} (\underline{\CC})$ to the space $\PP(P^s) \setminus C^{\vee} \bigcap \PP(P^s)$. Of course, we have the homomorphisms
 $$
 \pi_1 (\PP(P^s) \setminus C^{\vee} \bigcap \PP(P^s)) \longrightarrow \pi_1 (\PP(\HKC)\setminus C^{\vee})
 $$
 induced by the inclusions
 $$
  \PP(P^s) \setminus C^{\vee} \bigcap \PP(P^s) \hookrightarrow \PP(\HKC)\setminus C^{\vee},
  $$
  and the Lefschetz type theorem for the complement of a hypersurface in projective space tells us that the above homomorphism of the fundamental groups is an isomorphism if $h^s\geq 3$ and $\PP(P^s)$ is {\rm general},
  see \cite{Ha-Le}. However, there is no reason for summands $\PP(P^s)$'s to
  be general. So for each $s$, {\rm apriori}, we have conjugacy classes in different groups. It is an interesting issue to compare them. In fact, the linear maps
  $$
  \alpha^{s-1,s}: P^s \longrightarrow P^{s-1} 
  $$
  which are part of the quiver representations on the IVHS side may provide a natural comparison of the conjugacy classes for different values of $s$:
  those maps can be used to construct `{\rm transversal}' loops in $\PP(\HKC)$ interpolating between $\PP(P^s)$ for different $s$. We hope to address those issues elsewhere. 
  
  3) In the construction of $L({\mathfrak{L}}_l,\varrho^{\CC}_{{}_C};t)$ one can replace $\varrho^{\CC}_{{}_C}$ by {\rm any} finite dimensional representation of the monodromy group. Thus we can have a host of different functions on the geometric side associated to the zeta functions
  ${\mathfrak z}^{\pm}$ on the cohomological side.
\end{rem}

We now turn to another aspect of quantum-type constructions which will 
link the refined IVHS invariants with the moduli of elliptic curves with marked points.

\section{Metrized ribbon graph $\widehat{PG}_l$ and the moduli space ${\mathfrak M}^l_{1}$}
 
 We have seen that the graph $\widehat{PG}_l$ with its natural structure of a ribbon graph can be realized as a graph on the topological torus
 ${\mathbb{T}}$ in a such a way that the complement
 $$
 {\mathbb{T}} \setminus \widehat{PG}_l
 $$
 is a disjoint union of $l$ two dimensional cells. In this section we add a metric to the ribbon graph $\widehat{PG}_l$. This equips ${\mathbb{T}}$ with
 
 - a complex structure, 
 
 - $l$ marked points, the barycenters of two cells,
 
 - assigns the length for each boundary cycle. 

\noindent 
 Thus the {\it metrized} ribbon graph $\widehat{PG}_l$ connects to the space 
 ${\mathfrak M}^l_1 \times \RR^l_{+}$, where ${\mathfrak M}^l_1$ is the moduli space of curves of genus $1$ with $l$ marked points. All of this is a particular case of the general theory relating metrized ribbon graphs and the moduli spaces of curves, the theory which goes back to the works of J.Harer, D.Mumford, R.Penner, W.Thurston, see \cite{Har} and \cite{MuP} for overviews.
  
  \vspace{0.2cm}
  Recall: a metric on a ribbon graph $G$ is a positive valued function on the set of edges $E_G$ of $G$:
  $$
  m :E_G \longrightarrow \RR_{+}=(0, +\infty).
  $$
  Thus for every edge $e\in E_G$ a metric $m$ associates a positive real number $m(e)$ which is thought of as the length of $e$.
   A metrized ribbon graph $G$ is a pair $(G,m)$: a ribbon graph $G$ together with a metric $m$ on it.
   
  The main point is that a metric defines a complex structure on the compact oriented surface $S_G$ associated to the ribbon graph. The resulting smooth complex projective curve will be denoted $\Gamma_{G,m}$.
  In addition, it comes with with marked points determined by the boundary cycles of $G$:
   the topological surface $S_G$ comes with the collection $B_G$ of boundary cycles which are ordered once and for all
  $$
  B_G=\{B_1, \ldots, B_n\}.
  $$
  We recall that the surface $S_G$ is obtained from the ribbon graph by gluing in the closed $2$-dimensional disks along those boundary cycles.
  This gives the partition of $S_G$ into two dimensional cells. The barycenters of these cells give an ordered collection of points
  $$
  \{x_1, \ldots, x_n\}
  $$  
on the curve $\Gamma_{G,m}$. In addition, to every point $x_i$ is attached the length $m_i$ of the boundary cycle circling around $x_i$:
$$
m_i:= \sum_{\text{$e$ an edge of $B_i$}} m (e),
	$$
where the sum is taken over the edges of $B_i$. Thus we have the correspondence
$$
(G,B_1,\ldots, B_n,m) \mapsto (\Gamma_{G,m}, x_1,\ldots,x_n, m_1,\ldots,m_n)
$$
going from ribbon graphs $G$ with metric $m$ and an ordering of its boundary cycles $\{B_1,\ldots, B_n\}$ to complex projective curves $\Gamma_{G,m}$ with marked points $x_1,\ldots,x_n$, and the vector
$(m_1,\ldots,m_n)$ in $\RR^n_{+}$.  
  
  Fix the genus $s$ of the compact oriented surface $S=S_G$. It is related to the numerical data of $G$ via the topological Euler characteristic of $S$
  $$
  2-2s=|V_G|-|E_G|+|B_G|,
  $$
  where $V_G$ (resp., $E_G$ and $B_G$) stands for the set of vertices of $G$
  (resp., the set of edges and boundary cycles).
  
   The space of metrics on $G$ is $\RR_{+}^{E_G}$. The above correspondence gives rise to the continuous map
   \begin{equation}\label{met-moduli}
   \RR_{+}^{E_G} \longrightarrow {\mathfrak M}^n_s \times \RR^n_{+},
\end{equation}
  where $n=|B_G|$ is the number of the boundary cycles of $G$. The group of automorphism $Aut(G)$ of the graph $G$ acts on the space of metrics. It has the subgroup $Aut_B(G)$ preserving the fixed order
  $$
  B_1, \ldots,B_n
  $$
  of the boundary cycles of $G$.  The map \eqref{met-moduli} factors through the quotient space $\RR_{+}^{E_G} /Aut_B(G)$:
  $$
  \RR_{+}^{E_G} \longrightarrow \RR_{+}^{E_G} / Aut_B(G) \hookrightarrow {\mathfrak M}^n_s \times \RR^n_{+}.
  $$
  The image is an orbifold cell of  ${\mathfrak M}^n_s \times \RR^n_{+}$ of dimension $|E_G|$. 
   The fact that those orbifold cells cover ${\mathfrak M}^n_s \times \RR^n_{+}$ comes from the theory of Strebel differentials. This tells us how to go from a smooth projective curve with weighted marked points to ribbon graphs with boundary cycles circling around marked points and having the lengths equal to weights, see \cite{S}, \cite{MuP} for an overview. For convenience of the reader we recall the main result.
   
   Let $S$ be smooth projective curve of genus $s$ and let $p_1,\ldots,p_n$ be an ordered collection of $n$ distinct points on $S$; the pair $(s,n)$ is subject to the condition
   $$
   2s-2 +n >0;
   $$
   furthermore, to each point $p_i$ of the collection is assigned a positive real number $a_i$. 
   Then the theory of Strebel differentials tells us that there is a unique global section $q$ of $\OO_S (2K_S +2p_1 +\cdots+2p_n)$ subject to the following conditions 
   
   $\bullet$ $q$ is holomorphic on $S\setminus\{p_1,\ldots,p_n\}$,
   
   $\bullet$ $q$ has a pole of order two at every point $p_i$ of the collection,
   
   $\bullet$ for every $i\in [1,n]$ one has
   $$
   a_i =\int_{\gamma} \sqrt{q}
   $$
   for every compact horizontal leaf $\gamma$ of $q$ circling once around the point $p_i$; a horizontal leaf of $q$ means a parametrized curve
   $\gamma: (0,1) \longrightarrow S\setminus\{p_1,\ldots,p_n\}$ subject to
   \begin{equation}\label{qdiff}
   \gamma^{\ast}(q)>0;
\end{equation}
  more precisely, locally on $S\setminus\{p_1,\ldots,p_n\}$ we write $q$ as a local holomorphic section of $\OO_S(2K_S)=\Omega^2_S$:
  $$
  q=f(z)(dz)^2,
  $$
  where $z$ is a local parameter on $S$ and $f(z)$ is holomorphic; the condition \eqref{qdiff} means that
  $$
  f(\gamma(t)) \left(\frac{d\gamma(t)}{dt}\right)^2 >0, \forall t\in(0,1);
  $$
  a parametrized curve $\gamma: (0,1) \longrightarrow S\setminus\{p_1,\ldots,p_n\}$ is called a vertical leaf of $q$ if we have
   $$
  f(\gamma(t)) \left(\frac{d\gamma(t)}{dt}\right)^2 <0, \forall t\in(0,1).
  $$
   The horizontal (resp. vertical) leaves of the quadratic differential $q$ define the foliation of $S$ outside the poles and zeros of $q$. The metrized ribbon graph $G_q$ associated to $q$ is obtained by taking the zeros of $q$ as vertices, the edges as {\it closures} of horizontal leaves of $q$ connecting a pair of zeros of $q$; the orientation of $S$ determines a cyclic order of edges at every vertex. This makes $G_q$ into a ribbon graph. The metric on $G_q$ is defined by assigning the length
   $$
   m(e)=\int_{\gamma_e} \sqrt{q}
   $$
   to an edge $e$ connecting two zeros $p$ and $p'$, where $\gamma_e$ is a horizontal leaf of $q$ connecting  $p$ and $p'$, that is,
   $$
   \gamma_e:(0,1) \longrightarrow S\setminus\{p_1,\ldots,p_n\}
   $$
   is a horizontal leaf of $q$ subject to
    $$
    p=\lim_{t\to0^+} \gamma_e (t),\,\,\, p'=\lim_{t\to1^-} \gamma_e (t).
    $$
    
    \vspace{0.2cm}   
    We now return to our situation of the graph $\widehat{PG}_l$ with its natural ribbon structure. The graph has $2l$ vertices, $3l$ edges and $l$ boundary cycles. Furthermore, by definition the graph comes with a natural order of vertices, edges and boundary cycles:
    $$
    \begin{gathered}
    V_{\widehat{PG}_l}=\{(0), (0)',(1),(1)',\ldots, (l-1),(l-1)'\},
    \\
    E_{\widehat{PG}_l}=\{e^{-}_i,e^0_i, e^{+}_i | i=0,1, \ldots, l-1\},
 \\
 B_{\widehat{PG}_l} =\{B_0,B_1, \ldots, B_{l-1}\}, 
  \end{gathered}
$$
where $B_i$ is the unique boundary cycle of $\widehat{PG}_l$ containing the oriented edge $(i)\rightarrow (i)'$. With the order of edges fixed as above the space of metrics on $\widehat{PG}_l$ is identified with
$\RR^{3l}_{+}$. By the discussion above we have the continuous map
$$
\RR^{E_{\widehat{PG}_l}}_{+}\cong \RR^{3l}_{+} \longrightarrow {\mathfrak M}^l_1 \times \RR^{l}_{+}.
$$
This produces an orbifold cell of dimension $3l$ in ${\mathfrak M}^l_1 \times \RR^{l}_{+}$. On the other hand
$$
dim_{\RR} ({\mathfrak M}^l_1 \times \RR^{l}_{+})=dim_{\RR} ({\mathfrak M}^l_1) +l=2l+l=3l.
$$
 We summarize the above discussion in the following. 
\begin{pro}\label{pro:topcell}
	The space of metrics on the graph $\widehat{PG}_l$ with its natural ribbon structure determines an orbifold cell of dimension $3l$ in
	${\mathfrak M}^l_1 \times \RR^{l}_{+}$. This is a top dimensional cell of ${\mathfrak M}^l_1 \times \RR^{l}_{+}$.
\end{pro}
Actually, in our situation we can be more explicit: a metric $m$ on $\widehat{PG}_l$ gives the complex structure on $\mathbb{T}$ which is the curve $\Gamma_{\widehat{PG}_l,m}$; denote by $\omega_{m}$ its unique, up to a nonzero multiple, holomorphic differential; we also have the zig-zag paths $Z_{0,+}$ and $Z_{0,-}$ providing a $\ZZ$-basis of the homology group $H_1(\Gamma_{\widehat{PG}_l,m},\ZZ)$:
$$
\gamma_{-}:=Z_{0,-}, \,\, \gamma_{+}:=\frac{1}{l}(Z_{0,+}-Z_{0,-});
$$
this gives two periods
$$
\Omega^-_{m}:=\int_{\gamma_{-}} \omega_{m},\,\,\,\Omega^+_{m}:=\int_{\gamma_{+}} \omega_{m}
$$
and the lattice 
$$
\Lambda_m:=\ZZ\{\Omega^-_{m}, \Omega^+_{m}\} \subset \CC.
$$
The curve $\Gamma_{\widehat{PG}_l,m}$ is the elliptic curve
$\CC/\Lambda_m$. Thus over every point $m \in \RR^{E_{\widehat{PG}_l}}_{+}$
we have the elliptic curve $\CC/\Lambda_m$ whose lattice is equipped with a preferred basis $\{\Omega^-_{m}, \Omega^+_{m}\} $. We can do even better by ordering the basis so that the quotient of the basis belongs to the upper half plane
${\mathfrak{H}}$. Namely, we have the continuous map
$$
\RR^{E_{\widehat{PG}_l}}_{+}=\RR^{3l}_{+} \longrightarrow \RR^{\times}
$$
defined by the rule
$$
\RR^{3l}_{+} =\RR^{E_{\widehat{PG}_l}} \ni m \mapsto \Im \left(\frac{\Omega^+_m}{\Omega^-_m}\right) \in \RR^{\times}.
	$$
Thus the image of the map lies in the {\it same} connected component of $\RR^{\times}$ and we 
order the periods ${\Omega^{\pm}_m}$ once and for all so that their quotient has positive imaginary part; this gives the normalized period of $\Gamma_{\widehat{PG}_l,m}$;
it will be denoted $\tau_{m}$. We obtain the continuous map
$$
\RR^{E_{\widehat{PG}_l}}_{+}=\RR^{3l}_{+} \longrightarrow {\mathfrak H} = \text{the upper-half plane,}
$$
which sends a metric $m$ to the normalized period $\tau_{m}$; this gives the period lattice
$$
\Lambda_{\tau_{m}} := \ZZ \oplus \ZZ\tau_{m}
$$
 and the identification
$$
 \Gamma_{\widehat{PG}_l,m} \cong \CC/ \Lambda_{\tau_{m}};
 $$
under the identification the marked points $\{x_0(m),\ldots,x_{l-1}(m)\}$ go to the points
$$
 \{(x_i(m))-(x_0 (m)) | \,i\in [0,l-1]\} \subset \CC/\Lambda_{\tau_{m}}.
 $$
 
 In particular, we can write down the Weierstrass function and its derivative
 $$
 \begin{gathered}
 \wp(z)=\frac{1}{z^2}+\sum_{\omega \in \Lambda'_{\tau_{m}}} \frac{1}{(z-\omega)^2} -\frac{1}{\omega^2},
 \\
 \wp'(z)=-2\sum_{\omega \in \Lambda_{\tau_{m}}} \frac{1}{(z-\omega)^3},
 \end{gathered}
 $$
where $\Lambda'_{\tau_{m}}=\Lambda_{\tau_{m}}\setminus \{0\}$, and have the Weierstrass equation for $\Gamma_{\widehat{PG}_l,m}$:
$$
Y^2= 4X^3 -g_2(\tau_{m}) X -g_3(\tau_{m}),
$$
where the coefficients are the standard multiples of values of automorphic forms $E_4$ and $E_6$ at $\tau_{m}$:
$$
\begin{gathered}
g_2(\tau_{m}) =60 E_4(\tau_{m})=\sum_{(k,n)\in \ZZ^2 \setminus (0,0)} \frac{1}{(k+n\tau_{m})^4},
\\
 g_3(\tau_{m}) =140 E_6(\tau_{m})=\sum_{(k,n)\in \ZZ^2 \setminus (0,0)} \frac{1}{(k+n\tau_{m})^6},
 \end{gathered}
$$
see \cite{Sh}. From this we obtain the value of the modular function
\begin{equation}
	J(\tau_{m})=12^3 \frac{g^3_2 (\tau_{m})}{\Delta(\tau_{m})},
\end{equation}
where $\Delta(\tau_{m})$ is the discriminant
$$
\Delta(\tau_{m})=g^3_2 (\tau_{m}) -27g^2_3 (\tau_{m})
$$ 
of the cubic polynomial on the right side of the Weierstrass equation. This gives the modular function 
$$
\begin{gathered}
J_{\widehat{PG}_l} : \RR^{3l}_{+} \longrightarrow \CC \cong SL_2(\ZZ) \diagdown {\mathfrak{H}},
\\
m \mapsto J(\tau_{m}),
\end{gathered}
$$
attached to the graph $\widehat{PG}_{l}$. Thus in the situation of the graph
$\widehat{PG}_l$ we have the elliptic fibration over $\RR^{3l}_{+}$:
\begin{equation}\label{univell}
	\xymatrix{
	{\mathfrak G}_{\widehat{PG}_l}
	\ar[d]_{\pi_l}&\\
	\RR^{E_{\widehat{PG}_l}}_{+} \ar[r]^(.55){\psi_l}&{\mathfrak M}^l_1,
}
\end{equation}
where ${\mathfrak G}_{\widehat{PG}_l}$ is the universal elliptic curve
$$
{\mathfrak G}_{\widehat{PG}_l}=\left(\RR^{E_{\widehat{PG}_l}}_{+} \times \CC\right)/\sim;
$$
the equivalence `$\sim$' is defined as follows:
$$
(m,z) \sim (m',z') \Leftrightarrow \begin{cases}
	m=m',&\\
	(z-z')\in \Lambda_{\tau_{m}}.&
\end{cases}
$$
In addition, the marked points on the fibres of $\pi_l$ give rise to $l$ disjoint sections of $\pi_l$ 
$$
\varkappa_i: \RR^{E_{\widehat{PG}_l}}_{+} \longrightarrow {{\mathfrak G}_{\widehat{PG}_l}},
$$
for $i\in [0,l-1]$, defined by the rule
$$
\RR^{E_{\widehat{PG}_l}}_{+} \ni m \mapsto \varkappa_i(m)=(x_i(m)) - (x_0(m)) \in \CC/\Lambda_{\tau_{m}}=\pi^{-1}_l(m).
$$
With this additional structure we can complete \eqref{univell} to the commutative diagram
\begin{equation}\label{univell1}
	\xymatrix{
		{\mathfrak G}_{\widehat{PG}_l}
		\ar[d]_{\pi_l} \ar[r]&{\mathfrak M}^l_1 \times \RR^l_{+} \ar[d]\\
		\RR^{E_{\widehat{PG}_l}}_{+} \ar[r]&{\mathfrak M}^l_1,
	}
\end{equation}
where the top horizontal map sends the fibre $\pi^{-1}_l(m)$ over a point $m \in \RR^{E_{\widehat{PG}_l}}_{+}$ to the isomorphism class
of $(\CC/\Lambda_{\tau_{m}},\{\varkappa_i(m)\}_{i\in [0,l-1]})$
together with the vector $(m_0,\ldots,m_{l-1})$, where
$$
m_i=\sum_{e\in B_i} m(e)=m(e^0_i) +m(e^+_{i-1})+m(e^-_{i-1})+m(e^0_{i-2})+m(e^+_{i-2})+m(e^-_i).
$$

We will now connect the metrized graph $\widehat{PG}_{l}$ and the corresponding complex structures on the topological torus $\mathbb{T}$ with
the refinement of IVHS of the curve $C$. Namely, we consider the stratum ${\mathfrak{L}}_l$ of $\PP({\cal W}_{\Sigma^0_r})$ parametrizing points $([\xi],[\phi])$, where the $([\xi], [\phi])$-filtration of $W_{\xi}$
has length $l$. Over ${\mathfrak{L}}_l$ we consider the fibration 
$$
\pi_{{\mathfrak{L}}_l}:\OO^{\times}_{{\mathfrak{L}}_l} (-1) \longrightarrow {\mathfrak{L}}_l
$$
where $\OO^{\times}_{{\mathfrak{L}}_l} (-1)$ is the total space of the tautological line bundle $\OO_{{\mathfrak{L}}_l} (-1)$ with the zero section removed. We have seen that the IVHS refinement defines the family of representations 
$$
\rho_c :\OO^{\times}_{{\mathfrak{L}}_l} (-1) \longrightarrow {\mathfrak{Reps}}(PG_l).
$$
of the quiver $PG_l$ for every value $c\in \CC$, the value of the trace parameter: at every point $([\xi], [\phi])$ of ${\mathfrak{L}}_l$ we have the orthogonal decomposition
$$
W_{\xi}/W^l_{\xi}([\phi])=\bigoplus^{l-1}_{s=0} P^s([\xi], [\phi])
$$
and the summands are used to label the pairs of black-white vertices of $PG_l$; we fix a value $c$ of the trace parameter, then for every point $\xi\otimes \phi$ of $\OO^{\times}_{{\mathfrak{L}}_l} (-1)$ lying over $([\xi], [\phi])$ we have  linear maps
$$
\alpha^{t,s}_c(\xi,\phi): P^s ([\xi], [\phi])\longrightarrow P^t([\xi], [\phi]),
$$
for $t=s+j$, $j=0, \pm1$, and defined as long as $s$ and $t$ are in the interval $[0,l-1]$; thus the representation $\rho_c(\xi\otimes \phi)$ of $PG_l$ attached to a point $\xi\otimes \phi$
$$
\rho_c(\xi\otimes \phi)=\{\alpha^{t,s}_c(\xi,\phi): P^s ([\xi], [\phi])\longrightarrow P^t([\xi], [\phi])\}.
$$
We have learned that those representations admit natural extensions
to the representations $\widehat{\rho}_c(\xi\otimes\phi)$ of the quiver $\widehat{PG}_l$, that is, we know how, starting from $\rho_c(\xi\otimes \phi)$, to label the additional edges
$e^-_0$ and $e^+_{l-l}$ with linear maps
$$
\text{$P^0 ([\xi], [\phi])\longrightarrow P^{l-1}([\xi], [\phi])$ and $P^{l-1} ([\xi], [\phi])\longrightarrow P^{0}([\xi], [\phi])$.}
$$
Those labelings are provided by by our quantum-type invariants
$$
\begin{gathered}
	\tau^+_{\rho_c(\xi\otimes\phi)}: Hom_{\CC} (\CC[q,q^{-1}],\CC) \longrightarrow P^+([\xi], [\phi])=Hom(P^{l-1} ([\xi], [\phi]),P^{0} ([\xi], [\phi])),\\
	\tau^-_{\rho_c(\xi\otimes\phi)}: Hom_{\CC} (\CC[q,q^{-1}],\CC) \longrightarrow P^-([\xi], [\phi])=Hom(P^{0} ([\xi], [\phi]),P^{l-1} ([\xi], [\phi])),
\end{gathered}
$$
see Corollary \ref{cor:rhoc-ext}, also \eqref{Lpoly-Ppm} and the discussion leading to this construction for details. Thus we have the map
$$
 \widehat{\rho}_c : Hom_{\CC} (\CC[q,q^{-1}],\CC)\times  Hom_{\CC} (\CC[q,q^{-1}],\CC) \longrightarrow Maps\left(\OO^{\times}_{{\mathfrak{L}}_l} (-1), {\mathfrak{Reps}}(\widehat{PG}_l)\right)
 $$
 which assigns to a pair $(F,G)$ of linear functionals on $\CC[q,q^{-1}]$
 the family of representations
 $$
 \widehat{\rho}_c(F,G):\OO^{\times}_{{\mathfrak{L}}_l} (-1)\longrightarrow {\mathfrak{Reps}}(\widehat{PG}_l)
 $$
 of the quiver $\widehat{PG}_l$ extending the representations $\rho_c$; the value of $\widehat{\rho}_c(F,G)$ at a point $\xi\otimes\phi$ of $\OO^{\times}_{{\mathfrak{L}}_l} (-1)$ is denoted $\widehat{\rho_c}(F,G)(\xi\otimes\phi)$. This is the representation of $\widehat{PG}_l$ which coincides with ${\rho_c}(\xi\otimes\phi)$ on $PG_l$
 and receives the labels
 $$
 e^-_0 \rightarrow \tau^-_{\rho_c(\xi\otimes\phi)}(F),\hspace{0.2cm}  e^+_{l-1} \rightarrow \tau^+_{\rho_c(\xi\otimes\phi)}(G),
$$
on the additional edges $e^-_0$ and $e^+_{l-1}$ of $\widehat{PG}_l$. 

 Heuristically, we think of the values of the maps $\tau^{\pm}_{\rho_c(\xi\otimes\phi)}$ as `quantum paths' between the end spaces $P^{0}([\xi], [\phi])$ and $P^{l-1}([\xi], [\phi])$ of the representation $\rho_c(\xi\otimes\phi)$ of the quiver $PG_l$. Metrizing the graph $\widehat{PG}_l$ will put the `complex structures' on those `quantum paths'.
More precisely, every pair $F$ and $G$ of elements in $Hom_{\CC} (\CC[q,q^{-1}],\CC)$ give the representation $\widehat{\rho}_c (F,G)(\xi\otimes\phi)$ of the quiver $\widehat{PG}_l$ extending the representation $\rho_c(\xi\otimes\phi)$.  We attach the metric to $\widehat{\rho}_c (F,G)(\xi\otimes\phi)$
$$
m_{c, (F,G),\xi\otimes \phi} : E_{\widehat{PG}_l} \longrightarrow \RR_{+}
$$
by assigning to each edge of $\widehat{PG}_l$ the exponential of the (square of the) norm of the linear map
of $\widehat{\rho}_c (F,G)(\xi\otimes\phi)$ labeling that edge:
$$
\begin{gathered}
	e^0_s \mapsto exp(\parallel \alpha^{s,s}_c(\xi,\phi) \parallel^2), \,\forall s\in [0,l-1],
	\\
	e^-_s \mapsto exp(\parallel \alpha^{s-1,s}(\xi,\phi)  \parallel^2), \,\forall s\in [1,l-1],
	\\
	e^-_0 \mapsto exp(\parallel \tau^-_{\rho_c(\xi\otimes\phi)}(F)\parallel^2), 
	\\
	e^+_s \mapsto exp(\parallel \alpha^{s+1,s}(\xi,\phi)  \parallel^2), \,\forall s\in [0,l-2],
	\\
	e^+_{l-1} \mapsto exp(\parallel \tau^+_{\rho_c(\xi\otimes\phi)}(G) \parallel^2).
\end{gathered}
$$
The norm of maps is defined using the Hermitian metric on the spaces $\{P^s([\xi],[\phi])\}$ coming from the Hodge metric on $\HKC$, that is, for a linear map $A:P^s \longrightarrow P^t$ we let
$$
\parallel A\parallel^2 :=trace(AA^{\dag}),
$$
where $A^{\dag}$ is the adjoint of $A$ with respect to the Hermitian metrics on 
$P^s$ and $P^t$; the reference to $([\xi],[\phi])$ has been omitted to simplify the notation.

Thus we obtain the $C^{\infty}$ map
\begin{equation}\label{metrl}
metr_{{\mathfrak{L}}_l}:	\OO^{\times}_{{\mathfrak{L}}_l}(-1) \times \CC \times \left(Hom_{\CC} (\CC[q,q^{-1}],\CC)\right)^2 \longrightarrow \RR^{E_{\widehat{PG}_l}}_{+},
\end{equation}
sending a point $(\xi\otimes\phi, c, (F,G))$ in the domain to the metric
$m_{c,(F,G),\xi\otimes\phi}$ defined above. Composing with the map $\psi_l$ in \eqref{univell} relates ${\mathfrak{L}}_l $ to the moduli space ${\mathfrak{M}}^l_1$.
\begin{pro}\label{pro:plmap}
	Let ${\mathfrak{L}}_l$ be the stratum of $\PP({\cal W}_{\Sigma^0_r})$ parametrizing points $([\xi],[\phi])$, where the $([\xi], [\phi])$-filtration of $W_{\xi}$
	has length $l$. Then we have a continuous map
	$$
	\wp_l:\OO^{\times}_{{\mathfrak{L}}_l}(-1) \times \CC \times \left(Hom_{\CC} (\CC[q,q^{-1}],\CC)\right)^2 \longrightarrow {\mathfrak{M}}^l_1
	$$
	which sends  a point $(\xi\otimes\phi, c, (F,G))$ in the domain of the map
	to the isomorphism class of the curve $\Gamma_{\widehat{PG}_l, m_{c, (F,G),\xi\otimes\phi} }$ of genus $1$ with $l$ marked points.
	The map $\wp_l$ is the composition of the maps in \eqref{metrl} and \eqref{univell}
	$$
	\xymatrix@C=46pt{
	\OO^{\times}_{{\mathfrak{L}}_l}(-1)\times \CC \times \left(Hom_{\CC}(\CC[q,q^{-1}],\CC)\right)^2 \ar@/_2.5pc/[rr]^(.55){\wp_l} \ar[r]^(.75){metr_{{\mathfrak{L}}_l}} & \RR^{E_{\widehat{PG}_l}}_{+} \ar^(.55){\psi_l}[r]& {\mathfrak{M}}^l_1,
}
$$
where the map $metr_{{\mathfrak{L}}_l}$ is differentiable and $\psi_l$ is the composition of the natural maps
$$
\RR^{E_{\widehat{PG}_l}}_{+} \longrightarrow \RR^{E_{\widehat{PG}_l}}_{+}/ Aut_{\partial}(\widehat{PG}_l) \hookrightarrow {\mathfrak{M}}^l_1 \times \RR^l_{+} \longrightarrow {\mathfrak{M}}^l_1;
$$
the notation $Aut_{\partial}(\widehat{PG}_l)$ means the subgroup of automorphisms of the graph $\widehat{PG}_l$ preserving the order of its boundary cycles $B_{\widehat{PG}_l}=\{B_0,B_1,\ldots,B_{l-1}\}$.
\end{pro}

\begin{rem}\label{rem:lifttoH}
	As we remarked before there is a lifting of the map $\psi_l$ above to the upper-half plane. This gives the map
$$
\widetilde{\wp}_l :\OO^{\times}_{{\mathfrak{L}}_l}(-1) \times \CC \times \left(Hom_{\CC} (\CC[q,q^{-1}],\CC)\right)^2 \longrightarrow {\mathfrak{H}}
$$
which sends the metric $m_{c,(F,G),\xi\otimes\phi}$ to its normalized period $\tau_{m_{c,(F,G),\xi\otimes\phi}}$; that period will be denoted $\tau_{c,(F,G),\xi\otimes\phi}$; the corresponding elliptic curve
is the complex torus
$$
\Gamma_{\widehat{PG}_l,m_{c, (F,G),\xi\otimes\phi}}= \CC/\Lambda_{\tau_{c, (F,G),\xi\otimes\phi}},
$$
where the lattice $\Lambda_{\tau_{c,(F,G),\xi\otimes\phi}}=\ZZ\oplus \ZZ\tau_{c,(F,G),\xi\otimes\phi}$. We will denote this elliptic curve $\Gamma_{c,(F,G),\xi\otimes\phi}$:
$$
\Gamma_{c,(F,G),\xi\otimes\phi}= \CC/\Lambda_{\tau_{c,(F,G),\xi\otimes\phi}}.
$$
	
	Composing with the modular function $J$,
we obtain the modular function attached to the variety 
$\OO^{\times}_{{\mathfrak{L}}_l}(-1) \times \CC \times \left( \CC[q,q^{-1}]^{\vee}\right)^2$, where $\CC[q,q^{-1}]^{\vee}$ denotes the space $Hom_{\CC} (\CC[q,q^{-1}],\CC)$ of linear functionals on $\CC[q,q^{-1}]$:
$$
J_{{\mathfrak{L}}_l} : \OO^{\times}_{{\mathfrak{L}}_l}(-1) \times \CC \times \left( \CC[q,q^{-1}]^{\vee}\right)^2 \longrightarrow \CC
$$
which sends a point $(\xi\otimes\phi,c,(F,G))$ in the domain to the value
of the $J$-invariant $J({\tau_{c,(F,G),\xi\otimes\phi}})$.
 
Fixing a value of the trace parameter $c$ and a pair of linear functionals $(F,G)$ on $\CC[q,q^{-1}]$, that is, fixing a point $(c,(F,G)) \in \CC \times \left(\CC[q,q^{-1}]^{\vee}\right)^2 $ gives the map
$$
\widetilde{\wp}_{l } (\bullet,c,(F,G)): \OO^{\times}_{{\mathfrak{L}}_l}(-1) \longrightarrow {\mathfrak H},
$$
the restriction of $\widetilde{\wp}_l$ to the slice $\OO^{\times}_{{\mathfrak{L}}_l}(-1) \times \{(c,(F,G))\}$. The fibres of the map divide the fibre space $\OO^{\times}_{{\mathfrak{L}}_l}(-1)$ over the stratum ${\mathfrak{L}}_l$ into slices - the fibres of $\widetilde{\wp}_{l } (\bullet,c,(F,G))$ - to which is attached a specific elliptic curve: 
for each value
$\tau$ of $\widetilde{\wp}_{l } (\bullet,c,(F,G))$, we denote by 
$\widetilde{{\mathfrak{L}}}_l (\tau,c,(F,G))$ the inverse image of $\tau$ under the map $\widetilde{\wp}_{l } (\bullet,c,(F,G))$:
$$
\widetilde{{\mathfrak{L}}}_l (\tau,c,(F,G)):=\big(\widetilde{\wp}_{l } (\bullet,c,(F,G))
\big)^{-1}(\tau);
$$
this is the subset of $\OO^{\times}_{{\mathfrak{L}}_l}(-1)$ to which is attached the elliptic curve
$$
\Gamma_{\tau}:=\CC/\Lambda_{\tau}
$$
together with $l$ marked points. Those vary on $\Gamma_{\tau}$ as we move the points in $\widetilde{{\mathfrak{L}}}_l (\tau,c,(F,G))$. This gives the continuous map
\begin{equation}\label{delta-tau}
\delta_{\tau}: \widetilde{{\mathfrak{L}}}_l (\tau,c,(F,G)) \longrightarrow \Gamma^{l}_{\tau}=\text{the $l$-th power of $\Gamma_{\tau}$.}
\end{equation}
It should be also stressed that at every point $(\xi\otimes\phi) \in \widetilde{{\mathfrak{L}}}_l (\tau,c,(F,G)) $ the elliptic curve $\Gamma_{\tau}$ comes with additional structure: each marked point in
$\delta_{\tau}(\xi\otimes\phi)$ has a positive value attached to it, the length of the boundary cycle of $\widehat{PG}_l$ corresponding to the marked point. This is given by the sum of the values
of the metric $m_{c,(F,G),\xi\otimes\phi}$ on the edges of the boundary cycle:
$$
\delta^i_{\tau}(\xi\otimes\phi) \rightarrow m_i=\sum_{e\in B_i} m_{c, (F,G),\xi\otimes\phi} (e),
$$
where $\delta^i_{\tau}(\xi\otimes\phi)$ is the $i$-th coordinate of $\delta_{\tau}(\xi\otimes\phi)$ and $i\in [0,l-1]$. As we recalled, this gives a unique Strebel differential, call it $\omega_{c, (F,G),\xi\otimes\phi}$. In the case at hand this can be identified with an elliptic function on $\Gamma_{\tau}$ having poles of order $2$ precisely at the marked points
$\{\delta^i_{\tau}(\xi\otimes\phi)\}$.  More precisely, take a global section $s(\xi\otimes\phi)$ of $\OO_{\Gamma_{\tau}} (D_{\delta_{\tau}(\xi\otimes\phi)})$ vanishing on the divisor 
$$
D_{\delta_{\tau}(\xi\otimes\phi)}=\sum^{l-1}_{i=0}\delta^i_{\tau}(\xi\otimes\phi);
$$
 there is
a global section $\widehat{\omega}_{c, (F,G),\xi\otimes\phi}$ of the line bundle $\OO_{\Gamma_{\tau}} (2D_{\delta_{\tau}(\xi\otimes\phi)})$, such that the elliptic function $\omega_{c,(F,G),\xi\otimes\phi}$ has the form
$$
\omega_{c,(F,G),\xi\otimes\phi}=\frac{\widehat{\omega}_{c, (F,G),\xi\otimes\phi }}{s^2(\xi\otimes\phi)}.
$$
The function $\omega_{c,(F,G),\xi\otimes\phi}$ has $2l$ zeros, the zeros of the section $\widehat{\omega}_{c, (F,G),\xi\otimes\phi}$. Those are vertices of the graph $\widehat{PG}_l$; since the graph is trivalent all zeros are simple.  The above can be expressed as a continuous map
$$
S(\tau,c,(F,G)): \widetilde{\mathfrak{L}}_l (\tau,c,(F,G)) \longrightarrow \PP(H^0(\OO_{\Gamma_{\tau}} (2D_{\delta_{\tau}(\xi\otimes\phi)}))).
$$
In addition, outside the set of poles and zeros of $\omega_{c,(F,G),\xi\otimes\phi}$, the curve $\Gamma_{\tau}$ is foliated by horizontal and vertical leaves of $\omega_{c, (F,G),\xi\otimes\phi}$; the edges of $\widehat{PG}_l$ are the closures of the horizontal leaves connecting the zeros of $\omega_{c,(F,G),\xi\otimes\phi}$ into $l$ pairs.

\vspace{0.2cm} 

Another way to think about the $\widetilde{\wp}_{l }$ is as a continuous family
of elliptic curves
$$
\xymatrix{
{\mathfrak{G}}_{\widetilde{\wp}_{l }}
\ar[d]_{\pi_{\widetilde{\wp}_{l}}}&\\
\OO^{\times}_{{\mathfrak{L}}_l}(-1) \times \CC \times \left(\CC[q,q^{-1}]^{\vee}\right)^2 \ar[r]^(.8){\wp_l} & {\mathfrak M}^l_1,
}
$$
where the fibre of $\pi_{\widetilde{\wp}_{l}}$ over a point 
$(\xi\otimes\phi,c,(F,G))$ is the elliptic curve 
$$
\Gamma_{c,(F,G),\xi\otimes\phi}= \CC/\Lambda_{\tau_{c, (F,G),\xi\otimes\phi}} =\pi^{-1}_{\widetilde{\wp}_{l}}(\xi\otimes\phi,c,(F,G)).
$$
The projection $\pi_{\widetilde{\wp}_{l}}$ comes with $l$ disjoint sections.
They will be denoted $\varkappa^{\widetilde{\wp}_l}_i$:
$$
\varkappa^{\widetilde{\wp}_l}_i: \OO^{\times}_{{\mathfrak{L}}_l}(-1) \times \CC \times \left(\CC[q,q^{-1}]^{\vee}\right)^2 \longrightarrow {\mathfrak G}_{\widetilde{\wp}_l},\,\,i\in[0,l-1]
$$

One can equally change a perspective and consider the map $\widetilde{\wp}_{l }$ as the map from $\OO^{\times}_{{\mathfrak{L}}_l}(-1)$ to the space
of continuous maps from $\CC \times \left(\CC[q,q^{-1}]^{\vee}\right)^2$ to the upper-half plane
$$
\widetilde{\wp}_{l }: \OO^{\times}_{{\mathfrak{L}}_l}(-1) \longrightarrow {\mathfrak H}^{\CC \times \left(\CC[q,q^{-1}]^{\vee}\right)^2}.
$$ 
In this version each point $([\xi], [\phi])$ of ${\mathfrak{L}}_l$ acquires
moduli, the map
$$
\widetilde{\wp}_{l }([\xi], [\phi]): \CC^{\times}\times\CC \times \left(\CC[q,q^{-1}]^{\vee}\right)^2 \longrightarrow  {\mathfrak H},
$$
where the fibre of $\OO^{\times}_{{\mathfrak{L}}_l}(-1)$ over $([\xi], [\phi])$ is identified with $\CC^{\times}$. 

Recall that the factor $ \left(\CC[q,q^{-1}]^{\vee}\right)^2 $ contains the `geometric' part, the linear functionals of evaluation on points of
$\CC^{\times}$. This gives the maps
$$
\widetilde{\wp}_{l }([\xi], [\phi]): \CC^{\times} \times \CC \times (\CC^{\times})^2  \longrightarrow  {\mathfrak H}.
$$
Thus for every value $c \in \CC$ of the trace parameter we obtain
$$
\widetilde{\wp}_{l }(([\xi], [\phi]),c): \CC^{\times} \times (\CC^{\times})^2  \longrightarrow  {\mathfrak H},
$$
or, equivalently, a continuous family of elliptic curves over $\CC^{\times} \times (\CC^{\times})^2$. The reader may recall that one of the geometric objects attached to the graph $\widehat{PG}_l$ is the toric surface $X(\Delta_l)$, see \S8. So the algebraic torus $(\CC^{\times})^2$ can be viewed as the open orbit $X^{\circ}(\Delta_l)$ of $(\CC^{\times})^2$-action on $X(\Delta_l)$. From this point of view, the maps $\widetilde{\wp}_{l }(([\xi], [\phi]),c)$ can be viewed as continuous elliptic fibrations over that open part of $X(\Delta_l)$ parametrized by $\CC^{\times}$. The resulting map will be written as follows:
$$
{\widetilde{\wp}}^{X(\Delta_l)}_{l }(([\xi], [\phi]),c): \CC^{\times} \times X^{\circ}(\Delta_l)  \longrightarrow  {\mathfrak H}.
$$
The corresponding map to the moduli space $\mathfrak{M}^l_1$ is denoted
\begin{equation}
	{\wp}^{X(\Delta_l)}_{l }(([\xi], [\phi]),c): \CC^{\times} \times X^{\circ}(\Delta_l)  \longrightarrow\mathfrak{M}^l_1. 
\end{equation}
Furthermore, recall that on $X(\Delta_l)$ we have the Kasteleyn curves
$$
{\mathfrak C}_l(x)=\overline{(K_l(x,z,w)=0)},
$$
where $x$ is a vector of weights attached to the edges of the graph $\widehat{PG}_l$, $(z,w)$ parameters on the open part $(\CC^{\times})^2=X^{\circ}(\Delta_l)$ and $K_l(x,z,w)$ is the Kasteleyn determinant, see Proposition \ref{pro:detKastelyan} and Remark \ref{rem:Kastelyanl}. Thus the open part ${\mathfrak C}^{\circ}_l(x)=(K_l(x,z,w)=0)$ of each Kasteleyn curve is equipped with the
{\rm modular} maps
$$
{\wp}^{X(\Delta_l),x}_{l}(([\xi], [\phi]),c): \CC^{\times}\times{\mathfrak C}^{\circ}_l(x)   \longrightarrow\mathfrak{M}^l_1,
$$
where ${\wp}^{X(\Delta_l),x}_{l}(([\xi], [\phi]),c)$ denotes the restriction of ${\wp}^{X(\Delta_l)}_{l}(([\xi], [\phi]),c)$ to $\CC^{\times}\times {\mathfrak C}^{\circ}_l(x) $.
\end{rem}

\vspace{0.2cm}
\noindent
{\bf Connection with ${\bf H^{1,0}}(PG_l)$.} 
We will now connect the maps $\wp_l$ with the nonabelian Dolbeault variety ${\bf H^{1,0}}(PG_l)$. Recall that it comes with a distinguished anticanonical divisor $H_0$. Consider its complement
$$
{{\bf \stackrel{\circ}{H}^{1,0}}}(PG_l):={\bf H^{1,0}}(PG_l)\setminus H_0 \cong (\CC^{\times})^{l-1}.
$$
Over it we have the total space ${\mathfrak K}(PG_l)$ of the canonical bundle
$$
\xymatrix{
{\mathfrak K}(PG_l) \ar[r] \ar@{=}[d]& {{\bf \stackrel{\circ}{H}^{1,0}}}(PG_l) \ar@{=}[d]\\
\CC \times (\CC^{\times})^{l-1} \ar[r]& (\CC^{\times})^{l-1}.
}.
$$
One has a sort of `{\it action}' of  ${\mathfrak K}(PG_l)$ on the map $\wp_l$. More precisely, we can extend the map $\wp_l$ to the Cartesian product ${\mathfrak K}(PG_l) \times \OO^{\times}_{{\mathfrak L}_l}(-1)$. 
\begin{pro}\label{pro:cantwpl}
	The map $\wp_l$ admits the extension
	$$
	t(\wp_l): {\mathfrak K}(PG_l) \times \OO^{\times}_{{\mathfrak L}_l}(-1) \longrightarrow \left({\mathfrak M}^l_1\right)^{\left(\CC[q,q^{-1}]^{\vee}\right)^2}
	$$
	In particular, for every pair $(F,G) \in \left(\CC[q,q^{-1}]^{\vee}\right)^2$ we have the continuous map
$$
t(\wp_l)(F,G): {\mathfrak K}(PG_l) \times \OO^{\times}_{{\mathfrak L}_l}(-1) \longrightarrow {\mathfrak M}^l_1;
$$
The map $t(\wp_l)$ will be called {\rm the canonical twisting of $\wp_l$}. 
\end{pro}
\begin{pf}
	Use the identification 
	$$
	{ \mathfrak K}(PG_l)\cong \CC \times (\CC^{\times})^{l-1}.
	$$
	For a point $(c, {\bf z})=(c,z_1,\ldots,z_{l-1})$ in ${ \mathfrak K}(PG_l)$  and a point $\xi\otimes\phi$ in $\OO^{\times}_{{\mathfrak L}_l}(-1)$ lying over a point $([\xi],[\phi]) \in {\mathfrak L}_l$ we consider the quiver representation 
	$$
	\rho_c(\xi\otimes\phi)=\{\alpha^{t,s}_c(\xi,\phi):P^s([\xi],[\phi]) \longrightarrow P^t([\xi],[\phi])\}
	$$
	with the trace parameter equals $c$. Define the {\it twisting} of $\rho_c$ by ${\bf z}$
	$$
	  {\bf z} \cdot \rho_c
	  $$
	  as the following modification of the representation $\rho_c$:
	  
	  - the maps $\alpha^{s,s}_c:P^s([\xi],[\phi]) \longrightarrow P^s([\xi],[\phi])$ are unchanged,
	  
	  - the map $\alpha^{s,s-1}:P^s([\xi],[\phi]) \longrightarrow P^{s-1}([\xi],[\phi])$ is scaled by $z_s$, for each $s\in [1,l-1]$,
	  
	  - the map $\alpha^{s,s+1}:P^s([\xi],[\phi]) \longrightarrow P^{s+1}([\xi],[\phi])$ is scaled by $z^{-1}_{s+1}$, for each $s\in [0,l-2]$; 
	here and further on the reference to $\xi\otimes\phi$ is omitted to simplify the notation.
	
Take a pair $(F,G) \in \left(\CC[q,q^{-1}]^{\vee}\right)^2$ and define the extension $\widehat{{\bf z}\cdot \rho_c}(F,G)$ of the representation  ${\bf z} \cdot \rho_c(\xi\otimes\phi)$ to the one of the quiver $\widehat{PG}_l$. The representation $\widehat{{\bf z}\cdot \rho_c}(F,G)$ gives rise to the metric $m_{\widehat{{\bf z}\cdot \rho_c}(F,G)}$ on $\widehat{PG}_l$. The map $\psi_l:\RR^{E_{\widehat{PG}_l}}_{+} \longrightarrow {\mathfrak{M}}^l_1$ evaluated at $m_{\widehat{{\bf z}\cdot \rho_c}(F,G)}$ is defined to be the value of $t(\wp_l)$ at $((F,G),(c,{\bf z}),\xi\otimes\phi)$, that is, the map $t(\wp_l) $ is defined by the formula
	$$
	t(\wp_l)(F,G)((c,{\bf z}),\xi\otimes\phi) =\psi_l (m_{\widehat{{\bf z}\cdot \rho_c}(F,G)}), 
	$$ 
	for every point $((c,{\bf z}), \xi\otimes\phi) \in { \mathfrak K}(PG_l) \times \OO^{\times}_{{\mathfrak L}_l}(-1)$ and every pair $(F,G) \in \left(\CC[q,q^{-1}]^{\vee}\right)^2$. 
\end{pf}

The above result tells us that the stratum ${\mathfrak L}_l$ and pairs of linear functionals on $\CC[q,q^{-1}]$ provide moduli for continuous maps from ${ \mathfrak K}(PG_l)$ to ${\mathfrak M}^l_1$:
at every point $([\xi],[\phi])\in {\mathfrak L}_l$ denote $\CC^{\times}_{([\xi],[\phi]) }$ the fibre of $\OO^{\times}_{{\mathfrak L}_l }(-1)$ over $([\xi],[\phi])$, then we have the maps
$$
t(\wp_l)_{([\xi],[\phi])}(F,G):{ \mathfrak K}(PG_l)\times \CC^{\times}_{([\xi],[\phi]) } \longrightarrow {\mathfrak M}^l_1
$$
labeled by pairs $(F,G) \in \left(\CC[q,q^{-1}]^{\vee}\right)^2$.
Furthermore, a choice of $\xi\otimes\phi$ in the fibre $\CC^{\times}_{([\xi],[\phi]) } $ gives the maps
$$
t(\wp_l)_{\xi\otimes\phi}(F,G):{ \mathfrak K}(PG_l) \longrightarrow {\mathfrak M}^l_1
$$
We will now discuss how to obtain maps on ${{\bf \stackrel{\circ}{H}^{1,0}}}(PG_l)$ with values in ${\mathfrak M}^l_1$.

Given a section $\omega: {{\bf \stackrel{\circ}{H}^{1,0}}}(PG_l) \longrightarrow { \mathfrak K}(PG_l)$ of the canonical bundle of 
${{\bf \stackrel{\circ}{H}^{1,0}}}(PG_l)$ we obtain the map
$$
 t(\wp_l)_{\xi\otimes\phi}(F,G)\circ \omega: {{\bf \stackrel{\circ}{H}^{1,0}}}(PG_l)\longrightarrow {\mathfrak M}^l_1.
$$
There is a distinguished `constant' section given by the $(\CC^{\times})^{l-1}$-invariant form
$$
\omega_0=\frac{dz_1}{z_1}  \wedge \cdots \wedge\frac{dz_{l-1}}{z_{l-1}},
$$
 then any other global algebraic $(l-1)$-form on ${{\bf \stackrel{\circ}{H}^{1,0}}}(PG_l)$ can be written
 $$
 \omega=f \omega_0,
 $$
 where $f$ is a regular function on ${{\bf \stackrel{\circ}{H}^{1,0}}}(PG_l)=(\CC^{\times})^{l-1}$, that is, a $f$ is a  polynomial in $\CC[z_1,z^{-1}_1,\ldots, z_{l-1},z^{-1}_{l-1}]$; denoting the corresponding form by $\omega_f$, gives the map
 $$
 \CC[(\CC^{\times})^{l-1}]\cong \CC[z_1,z^{-1}_1,\ldots, z_{l-1},z^{-1}_{l-1}] \longrightarrow \left({\mathfrak M}^l_1\right)^{{{\bf \stackrel{\circ}{H}^{1,0}}}(PG_l)}
$$ 
sending the function $f \in \CC[z_1,z^{-1}_1,\ldots, z_{l-1},z^{-1}_{l-1}]$
to the map
\begin{equation}\label{f-tpl}
 t(\wp_l)_{\xi\otimes\phi}(F,G)\circ\omega_f: {{\bf \stackrel{\circ}{H}^{1,0}}}(PG_l)\longrightarrow {\mathfrak M}^l_1.
\end{equation} 

Let us recall that  ${{\bf {H}^{1,0}}}(PG_l)$ is a toric variety with the fan $\Delta({{\bf {H}^{1,0}}}(PG_l))$ generated by the vertices of the cube
$[-1,1]^{l-1}$. This gives a natural choice for Laurent polynomials spanned by the monomials corresponding to the vertices of the cube, that is, for a vertex $\epsilon=(\epsilon_1,\ldots,\epsilon_{l-1})$ of the cube $[-1,1]^{l-1}$ set
$$
z^{\epsilon}:=\prod^{l-1}_{i=1} z^{\epsilon_i}_i
$$
the monomial corresponding to a vertex $\epsilon$; then the natural choice alluded to above is the nonzero polynomials of the form
$$
f=\sum_{\epsilon} a_{\epsilon}z^{\epsilon},
$$
where the sum is taken over the vertices $\epsilon$ of the cube and the coefficients  $a_{\epsilon}$ are in a subring $R$ of $\CC$, 

For $f$ fixed we also have {\it natural} choices for $(F,G)$ in \eqref{f-tpl} determined by 
the polynomial coefficients of the Picard-Fuchs operator $PF_f$ of $f$. We briefly recall the definition, see \cite{CCGGK} for details: for a Laurent polynomial $f\in \CC[z_1,z^{-1}_1,\ldots, z_{l-1},z^{-1}_{l-1}]$ one defines the {\rm classical period of $f$} 
$$
\pi_f(q) =\left(\frac{1}{2\pi i}\right)^{l-1} \int_{|z_1|=\cdots=|z_{l-1}|=1} \frac{\omega_0}{1-qf};
	$$
	it is known that there is a differential operator 
	$$
	L=\sum p_j (q) D^j
	$$
	which is a polynomial in $D=q\frac{d}{dq}$, the coefficients $p_j(q)$'s are polynomials in $q$, and it annihilates the classical period
	$$
	L \cdot \pi_f =0;
	$$
	one obtains the Picard-Fuchs operator $PF_f$ by choosing $L$ above of the smallest possible degree with respect to $D$, say $k$,
	$$
	PF_f=\sum^k_{j=0} p_j (q) D^j
	$$
	and having the coefficient $p_k(q)$ of the smallest possible degree; such an operator is clearly unique up to a nonzero scalar multiple.
	
	Returning to the maps in \eqref{f-tpl}, what we referred to as a natural choice for $F,G$, are the polynomials $\{p_j(q)\}$, the coefficients of the Picard-Fuchs operator $PF_f$ of $f$; we are using here the identification
	$$
	Hom_{\CC} (\CC[q,q^{-1}],\CC) \cong \CC[q,q^{-1}].
	$$
	We summarize the above discussion in the following.
	\begin{pro}\label{pro:Dolb-PF-moduli}
		For every polynomial $f$ in $\CC[z_1,z^{-1}_1,\ldots, z_{l-1},z^{-1}_{l-1}]$ and every pair $(F,G)$ in $\left(\CC[q,q^{-1}]^{\vee}\right)^2$, there is a distinguished continuous map
		$$
	t(\wp_l)(f;F,G):=	 t(\wp_l)(F,G)\circ\omega_f : \OO^{\times}_{{\mathfrak L}_l}(-1) \times  {{\bf \stackrel{\circ}{H}^{1,0}}}(PG_l) \longrightarrow {\mathfrak M}^l_1;
	$$
	defined by the rule
	$$
\OO^{\times}_{{\mathfrak L}_l}(-1) \times  {{\bf \stackrel{\circ}{H}^{1,0}}}(PG_l) \ni	(\xi\otimes\phi, {\bf z}) \mapsto t(\wp_l)_{\xi\otimes\phi}(F,G) (\omega_f ({\bf z})) \in {\mathfrak M}^l_1.
$$
 Furthermore, let 
	$$
	PF_f=\sum^k_{j=0} p_j (q) D^j
	$$
	be the Picard-Fuchs operator of $f$, then the coefficients $\{p_j(q)\}$ provide natural labels for the maps from ${{\bf \stackrel{\circ}{H}^{1,0}}}(PG_l)$ to ${\mathfrak M}^l_1$; namely, we have 
	$$
	t(\wp_l)(f;p_j,p_{j'}): \OO^{\times}_{{\mathfrak L}_l}(-1)  \times {{\bf \stackrel{\circ}{H}^{1,0}}}(PG_l) \longrightarrow {\mathfrak M}^l_1,
	$$
	for every pair $p_j$ and $p_{j'}$ of coefficients of the Picard-Fuchs operator of $f$.
	In addition, specializing $f$ to a nonzero polynomial
	$$
	f_{\bf a}=\sum_{\epsilon} a_{\epsilon}z^{\epsilon},
	$$
	where the sum runs through the set of vertices $V^0_{l-1}$ of the cube $[-1,1]^{l-1}$, the fan {\small$\Delta({{\bf {H}^{1,0}}}(PG_l))$} of ${{\bf {H}^{1,0}}}(PG_l)$, and ${\bf a}=(\epsilon \mapsto a_{\epsilon})_{\epsilon \in V^0_{l-1}}$, a nonzero function on the set of vertices $V^0_{l-1}$,
gives the maps
$$
	t(\wp_l)({\bf a};p^{\bf a}_j,p^{\bf a}_{j'}):=	t(\wp_l)(f_{\bf a};p^{\bf a}_j,p^{\bf a}_{j'})	: \OO^{\times}_{{\mathfrak L}_l}(-1)  \times{{\bf \stackrel{\circ}{H}^{1,0}}}(PG_l) \longrightarrow {\mathfrak M}^l_1,
$$
where $\{p^{\bf a}_j\}$ the coefficients of the Picard-Fuchs operator $PF_{f_{\bf a}}$ of $f_{\bf a}$.
	\end{pro}
\begin{example}
	For a vertex $\epsilon$ of the cube $[-1,1]^{l-1}$ take ${\bf a}=1_{\epsilon}$, the function on the set of vertices taking the value $1$ at $\epsilon$ and zero otherwise. The corresponding polynomial
	$$
	f_{1_{\epsilon}}=z^{\epsilon}.
	$$
	For $\epsilon =(-1,-1,\ldots,-1)$, we have
	$$
	f_{1_{\epsilon}}=z^{-1}_1 z^{-1}_2 \cdots z^{-1}_{l-1}; 
	$$
	the classical period $\pi_{z^{\epsilon}}(q)=1$. The Picard-Fuchs operator $PF_1=D$; the corresponding collection of polynomials is
	$\{p_0=0,p_1=1\}$ and we obtain the maps
	$$
	t(\wp_l)(1_{\epsilon};0,0), \,\,	t(\wp_l)(1_{\epsilon};0,1), \,\,	t(\wp_l)(1_{\epsilon};1,0),\,\,	t(\wp_l)(1_{\epsilon};1,1).
	$$
	defined on $\OO^{\times}_{{\mathfrak L}_l}(-1)  \times  {{\bf \stackrel{\circ}{H}^{1,0}}}(PG_l)$ with values in ${\mathfrak M}^l_1$.
\end{example}

	Another natural choice for $F$ and $G$ in the maps of the above proposition comes from the stratum ${\mathfrak L}_l$ itself. Recall that it is the disjoint union of the substrata
	${\mathfrak L}_l(h^l,\lambda)$ where the dimensions of the graded pieces of $([\xi],[\phi])$-filtrations remain constant and determined by
	$h^l$ and the partition $\lambda$; to such a stratum we attached the Laurent polynomial
	$$
	h_{{\mathfrak L}_l(h^l,\lambda)}=\sum^{l-1}_{s=0}\lambda'_s q^{(-1)^s s},
	$$
	where $\lambda'$ is the partition conjugate to $\lambda$ and $\lambda'_s$ its $s$-th part.
	
	Thus we obtain the following. 
	\begin{pro}\label{pro:substrat-moduli}
		The collection of the substrata  $\{{\mathfrak L}_l(h^l, \lambda)\}$ of ${\mathfrak L}_l$ gives rise to the collection of the Laurent polynomials
		$$
		\{h_{{\mathfrak L}_l(h^l,\lambda)}\}.
		$$
		For every pair of the partitions $(\lambda,\mu)$ occurring in the collection we obtain the continuous maps 
		$$
		t(\wp_l)(f;\lambda,\mu):=	 t(\wp_l)(f;h_{{\mathfrak L}_l(h^l,\lambda)},h_{{\mathfrak L}_l(h^l,\mu)}): \OO^{\times}_{{\mathfrak L}_l}(-1)  \times {{\bf \stackrel{\circ}{H}^{1,0}}}(PG_l) \longrightarrow {\mathfrak M}^l_1,
		$$
		for every polynomial $f$ in $\CC[z_1,z^{-1}_1,\ldots, z_{l-1},z^{-1}_{l-1}]$. 
	\end{pro}

\vspace{0.2cm}
\noindent
{\bf Connection to the maps $\lambda_{\phi}$.} Let us recall \S6.4 where we constructed the map
$$
\lambda_{\phi}: \CC\times \OO^{\times}_{{\mathfrak L}_l (\phi)}(-1) \longrightarrow \CC(C)[q,q^{-1}],
$$
where $\phi$ is an admissible global section of the canonical bundle $\OO_C(K_C)$, see Definition \ref{def:phi-admiss}, that is, for points $([\xi],[\psi]) \in {\mathfrak L}_l $ we consider the orthogonal decomposition
$$
\HKC=W_{\xi}\oplus W^{\perp}_{\xi}=\left(\bigoplus^l_{s=0} P^s([\xi],[\psi])
\right) \oplus W^{\perp}_{\xi}
$$
and write $\phi$ as the sum of its orthogonal components
\begin{equation}\label{phi-dec-Lphi}
\phi=\sum^{l-1}_{s=0}\phi^s([\xi],[\psi]) + \phi^{\perp}([\xi])
\end{equation}
with respect to the direct sum above; the subset of points $([\xi],[\psi])$ in ${\mathfrak L}_l $ where all components $\phi^s([\xi],[\psi])$ are nonzero is denoted  by ${\mathfrak L}_l (\phi)$ and called the domain of $\phi$; we call $\phi$ admissible if its domain ${\mathfrak L}_l (\phi)$ is nonempty. The map $\lambda_{\phi}$ assigns
to a point $(c,\xi\otimes\psi)$ in $\CC\times \OO_{{\mathfrak L}_l (\phi)}(-1)$ a nonzero Laurent polynomial
\begin{equation}\label{lambdaphi-remind}
\lambda_{\phi}(c) (\xi\otimes\psi)=\sum^{l-1}_{s=0} F_s(c,\xi\otimes\psi)q^{(-1)^s s}
\end{equation}
where the coefficients $F_s(c,\xi\otimes\psi)$ are rational functions on $C$, determined by the decomposition \eqref{phi-dec-Lphi} and the representation
$\rho_{c} (\xi\otimes\psi)$ of the quiver $PG_l$, see \eqref{lambdaphi-formula} for the precise formula. 
Putting this together with the map $\wp_l$ in Proposition \ref{pro:plmap} we obtain the following.
\begin{pro}\label{pro:pl-C-moduli}
	Let $\phi \in \HKC$ be admissible and let ${\mathfrak L}_l (\phi)$ be its domain. Then to every point $([\xi],[\psi])$ in ${\mathfrak L}_l (\phi)$ is attached: 
	 
	 - the `polar' divisor $D(\phi,[\xi],[\psi]) \in |lK_C|$,
	 
	 - the continuous map
	$$
	\wp^{\phi}_l (\xi\otimes\psi): \CC \times \left(C\setminus D(\phi,[\xi],[\psi])\right) \longrightarrow {\mathfrak M}^l_1
	$$
	for every nonzero $\xi\otimes \psi $ in  the fibre of the total space $\OO_{{\mathfrak L}_l (\phi)}(-1)$ lying over the point $([\xi],[\psi])$ in ${\mathfrak L}_l (\phi)$. 
\end{pro} 
\begin{pf}
	Let $([\xi],[\psi])$ in ${\mathfrak L}_l (\phi)$. Evaluate the map
	$$
	\lambda_{\phi}: \CC \times \OO^{\times}_{{\mathfrak L}_l (\phi)}(-1)  \longrightarrow \CC(C)[q,q^{-1}]
	$$
	at a point $(c,\xi\otimes\psi)$ of $\CC \times\OO^{\times}_{{\mathfrak L}_l (\phi)}(-1) $ to obtain the Laurent polynomial
	\begin{equation}\label{lambdaphi-psi}
	\lambda_{\phi}(c)(\xi\otimes\psi)=\sum^{l-1}_{s=0}F_s(c,\xi\otimes\psi)q^{(-1)^ss},
\end{equation}
	where the coefficients  $F_r(c,\xi\otimes\psi)$ are rational functions on $C$. More precisely, those functions have the form
	$$
	F_s(c,\xi\otimes\psi) =\frac{\omega^s(c,\xi\otimes\psi)}{\phi^s([\xi],[\psi])},
	$$
	where on the right side is the quotient of global sections of the canonical bundle $\OO_C(K_C)$; the one in the denominator is as in the decomposition \eqref{phi-dec-Lphi} and for the numerator we refer to the discussion preceding \eqref{lambdaphi-formula}; we do not need to know its expression, it is important to know that it is nonzero for all $s\geq 1$. 
	
	The `polar' divisor in the statement of the proposition is the sum of the zero divisors of the denominators of the functions $F_s(c,\xi\otimes\psi)$:
	$$
	D(\phi,[\xi],[\psi]):= \sum^{l-1}_{s=0} Z_{\phi^s([\xi],[\psi])}
	$$
	where $Z_{\phi^s([\xi],[\psi])}=(\phi^s([\xi],[\psi])=0)$ is a divisor in $|K_C|$. Hence $D(\phi,[\xi],[\psi])$ is a divisor in $|lK_C|$ as asserted and depends only on the point $([\xi],[\psi])$ in ${\mathfrak L}_l (\phi)$.
	
	From the definition of $D(\phi,[\xi],[\psi])$ it follows that all functions $F_s(c,\xi\otimes\psi)$ are regular on the complement of that divisor in $C$. Hence $\lambda_{\phi}(c)(\xi\otimes\psi)$ can be viewed
	as a function
	$$
	 C\setminus D(\phi,[\xi],[\psi]) \longrightarrow \CC[q,q^{-1}];
	 $$
	 defined by the rule
	 $$
	 C\setminus D(\phi,[\xi],[\psi]) \ni x \mapsto \lambda_{\phi}(c)(\xi\otimes\psi)(x)=\sum^{l-1}_{s=0}F_s(c,\xi\otimes\psi)(x)q^{(-1)^s s} \in \CC[q,q^{-1}].
	 $$
	 We use the same notation $\lambda_{\phi}(c)(\xi\otimes\psi)$ for this function. 
	 
	 The Laurent polynomials $\lambda_{\phi}(c)(\xi\otimes\psi)(x)$ thus obtained are used to evaluate the maps 
	 $$
	 \tau^{\pm}_{\rho_c(\xi\otimes\psi)}: \CC[q,q^{-1}] \longrightarrow P^{\pm}([\xi],[\psi]),
	 $$
	 see Corollary \ref{cor:rhoc-ext}. This gives the representations denoted $\widehat{\rho}_{c}(\xi\otimes\psi) (x)$ of the quiver
	 $\widehat{PG}_l$ extending the representation $\rho_{c}(\xi\otimes\psi)$ of $PG_l$. Assigning the metric of the graph $\widehat{PG}_l$ to each representation $\widehat{\rho}_{c}(\xi\otimes\psi) (x)$ gives the map
	 $$
	 C\setminus D(\phi,[\xi],[\psi]) \longrightarrow \RR^{E_{\widehat{PG}_l}}_{+}
	 $$
	 defined by the formula
	 $$
	C\setminus D(\phi,[\xi],[\psi])\ni x\mapsto m_{\widehat{\rho}_{c}(\xi\otimes\psi) (x)}.
	$$
	Composing with the map $\psi_l: \RR^{E_{\widehat{PG}_l}}_{+} \longrightarrow {\mathfrak M}^l_1$ we obtain the map
	$$
	C\setminus D(\phi,[\xi],[\psi]) \longrightarrow {\mathfrak M}^l_1
	$$
	which we denote $\wp^{\phi}_l(c,\xi\otimes\psi)$. Varying the trace parameter $c \in \CC$ produces the map
$$
\wp^{\phi}_l(\xi\otimes\psi):	\CC\times C\setminus D(\phi,[\xi],[\psi]) \longrightarrow {\mathfrak M}^l_1
	$$
	defined by the formula
	$$
	\wp^{\phi}_l(\xi\otimes\psi)(c,x)=\wp^{\phi}_l(c,\xi\otimes\psi)(x), \forall (c,x)\in \CC \times C\setminus D(\phi,[\xi],[\psi]).  
	$$
This is the map asserted in the statement of the proposition\end{pf}

Let us fix the trace parameter $c$. The result of the proposition above tells us that we have continuous maps
\begin{equation}\label{wp-xipsi}
\wp^{\phi}_l(\xi\otimes\psi)(c,\bullet): C\setminus D(\phi,[\xi],[\psi]) \longrightarrow {\mathfrak M}^l_1,
\end{equation}
where the parameters for the maps  are nonzero points $\xi\otimes\psi$ in the fibre of the tautological line bundle $\OO_{{\mathfrak L}_l (\phi)}(-1)$ over points $([\xi],[\psi])$ in ${\mathfrak L}_l (\phi)$.
This is reminiscent of a nonlinear $\sigma$-model in string theory. Under the assumption that the polar divisor 
$D(\phi,[\xi],[\psi])$ is reduced, we will attach to the triple
$(C,D(\phi,[\xi],[\psi]), \wp_l(\xi\otimes\psi)(c,\bullet) ) $ the following objects:

\vspace{0.2cm}
$\bullet$ a Strebel differential on the curve $C$ with the double poles precisely along the divisor $D(\phi,[\xi],[\psi])$,

$\bullet$ the trivalent graph on $C\setminus D(\phi,[\xi],[\psi])$,

$\bullet$ elliptic curves with $l$ marked points attached to vertices of the graph in the item above.

\vspace{0.2cm}
We begin the discussion by attaching Strebel differentials and trivalent graphs to each summand of the polar divisor $D(\phi,[\xi],[\psi])$. Recall,
 for a point $([\xi],[\psi]) \in {\mathfrak L}_l (\phi)$, the polar divisor $D(\phi,[\xi],[\psi])$ is the sum of the divisors 
$$
Z_{\phi^s([\xi],[\psi])}=(\phi^s([\xi],[\psi])=0)
$$
 and where
$\phi^s([\xi],[\psi]) $, for $s\in [0,l-1]$, are the components of $\phi$
with respect to the direct sum decomposition
$$
W_{\xi}/W^l([\xi],[\psi]) =\bigoplus^{l-1}_{s=0}P^s([\xi],[\psi])
$$
at the point $([\xi],[\psi])$. Thus we have the continuous map
\begin{equation}\label{l-phi-map}
{\mathfrak d}_l:{\mathfrak L}_l (\phi) \longrightarrow \left(\PP(\HKC)\right)^l
\end{equation}
which is defined by the rule
$$
{\mathfrak L}_l (\phi)\ni ([\xi],[\psi]) \mapsto ([\phi^0([\xi],[\psi])], [\phi^1([\xi],[\psi])], \ldots,[\phi^{l-1}([\xi],[\psi])]) \in \left(\PP(\HKC)\right)^l.
$$
 We wish to have some of the zero divisors $Z_{\phi^s([\xi],[\psi])}$ reduced. Precisely, the nonzero global sections of $\OO_C(K_C)$ whose zero divisor is not reduced are parametrized by the dual variety $C^{\vee}$ of $C \subset \PP(\HKC^{\ast})$, the variety of hyperplanes in $\PP(\HKC^{\ast})$ which are tangent to $C$ at some point of $C$. Set
$$
p_s: \left(\PP(\HKC)\right)^l \longrightarrow \PP(\HKC)
$$
 the projection on the $s$-th factor, for $s\in [0,l-1]$. The pull back $p^{\ast}_s (C^{\vee})$
 gives the divisor parametrizing $l$-tuples in $\left(\PP(\HKC)\right)^l$ with the hyperplane on the $s$-th coordinate being tangent to $C$. Set
 $$
 {\mathfrak L}^{(s)}_l (\phi):={\mathfrak d}^{-1}_l (p^{\ast}_s (C^{\vee})),
 $$
 the inverse image of that divisor under the map ${\mathfrak d}_l$. We wish the subset $ {\mathfrak L}^{(s)}_l (\phi)$ to be smaller than ${\mathfrak L}_l (\phi)$ for some $s \in [0,l-1]$. More generally, for every nonempty subset $I$ of
 $[0,l-1]$ define
 $$
 {\mathfrak L}^{(I)}_l (\phi)=\bigcup_{s\in I} {\mathfrak L}^{(s)}_l (\phi).
 $$
 and let
 \begin{equation}\label{LphiI}
 {\mathfrak L}_l (\phi,I):= {\mathfrak L}_l (\phi) \setminus {\mathfrak L}^{(I)}_l (\phi),
\end{equation}
 the complement of ${\mathfrak L}^{(I)}_l (\phi)$ in
  ${\mathfrak L}_l (\phi)$. Thus, if ${\mathfrak L}_l (\phi,I)$ is nonempty, it parametrizes the points $([\xi],[\psi]) \in {\mathfrak L}_l (\phi)$ with the property that in the decomposition
  $$
  \phi =\phi^{\perp} ([\xi],[\psi]) +\sum^{l}_{s=0} \phi^s ([\xi],[\psi]) 
  $$
  according to the orthogonal sum
  $$
  \HKC =W^{\perp}_{\xi}\oplus \left(\bigoplus^l_{s=0} P^s ([\xi],[\psi])\right)
  $$
  the terms $\phi^s ([\xi],[\psi]) \in P^s ([\xi],[\psi])$ are subject to
  \begin{equation}\label{phir-cond}
  	[\phi^s ([\xi],[\psi])] \in \PP(\HKC)\setminus C^{\vee},\,\, \forall s\in I. 
  \end{equation}
In other words the zero divisor
$$
Z_{\phi^s ([\xi],[\psi])}=(\phi^s ([\xi],[\psi])=0)\,\,\, \text{\it consists of $(2g-2)$ distinct points on $C$, for every $s \in I$.}
$$
We will now attach weights to the points of  those $Z_{\phi^s ([\xi],[\psi])}$. For this we use the fact that $\xi$ annihilates $\phi^s ([\xi],[\psi])$, for all $s\in [0,l]$, that is the cup product with $\xi$ vanishes
$$
\xi \phi^s ([\xi],[\psi])=0,\,\,\text{ in $H^1(\OO_C)$ for all $s\in [0,l]$.}
$$
This allows to identify $\xi$ with a function on $Z_{\phi^s ([\xi],[\psi])}$: we are using here the map \eqref{Fs-map}. The construction is formalized as follows.

\begin{pro}\label{pro:LphiI}
	Assume the locus ${\mathfrak L}_l (\phi)$ is nonempty. Then we have the $C^{\infty}$ map
	$$
	{\mathfrak d}_l:{\mathfrak L}_l (\phi) \longrightarrow \PP(\HKC)^l
	$$
defined in \eqref{l-phi-map}. Let $\widetilde{{\mathfrak L}}_l (\phi)$ be the restriction to ${\mathfrak L}_l (\phi)$ of the pull back
$p^{\ast}_1 (\OO_{\PP(H^1(\Theta_{C})}(-1))$ of the tautological bundle of
$\PP(H^1(\Theta_{C}))$. Then
the diagram
$$
\xymatrix{
\widetilde{{\mathfrak L}}_l (\phi) \ar[d]& \left(q_{1\ast}\OO_{{\cal Z}_C}\right)^l \ar[d]\\
{\mathfrak L}_l (\phi) \ar[r]^(.35){{\mathfrak d}_l}& \PP(\HKC)^l
}
$$
where the vertical arrows are the structure projections, can be completed to the commutative diagram
 \begin{equation}\label{tildLlphi-functions}
\xymatrix{
	\widetilde{{\mathfrak L}}_l (\phi) \ar[d]\ar[r]^{F_l}&\left(q_{1\ast}\OO_{{\cal Z}_C}\right)^l \ar[d]\\
	{\mathfrak L}_l (\phi) \ar[r]^(.35){{\mathfrak d}_l}& \PP(\HKC)^l
}
\end{equation}
where the map $F_l$ sends the line $	\widetilde{{\mathfrak L}}_l (\phi)_{([\xi],[\psi])}=\CC\xi$ over a point $([\xi],[\psi])$ in ${\mathfrak L}_l (\phi)$ into the Cartesian product 
$$
\prod^{l-1}_{s=0}H^0(\OO_{Z_{\phi^s([\xi],[\psi])}}),
$$
the fibre of $\left(q_{1\ast}\OO_{{\cal Z}_C}\right)^l$ over the point 
$$
{\mathfrak d}_l(([\xi],[\psi]))=([\phi^0([\xi],[\psi])],[\phi^1([\xi],[\psi])],\ldots, [\phi^{l-1}([\xi],[\psi])]);
$$
that map is defined by the formula
$$
(F_l)_{([\xi],[\psi])} (\xi')=\left(F(\phi^s([\xi],[\psi]),\xi')\right)_{\{s=0,1,\ldots,l-1\}},
$$
where $\xi' \in \CC \xi$ and the map $F(\cdot,\cdot)$ is as in \eqref{Fs-map}.
\end{pro}

Let $I$ be a nonempty subset of $[0,l-1]$ and assume the locus
${\mathfrak L}_l (\phi,I)$ is nonempty; recall, this means that for
every point $([\xi],[\psi]) \in {\mathfrak L}_l (\phi,I)$ the divisor
$$
Z_{\phi^s([\xi],[\psi])}=(\phi^s([\xi],[\psi])=0)
$$
is reduced for each $s\in I$. 
Consider the diagram \eqref{tildLlphi-functions} over the locus ${\mathfrak L}_l (\phi,I)$
$$
\xymatrix{
	\widetilde{{\mathfrak L}}_l (\phi,I) \ar[d]\ar[r]^{F_l}&\left(q_{1\ast}\OO_{{\cal Z}_C}\right)^l \ar[d]\\
	{\mathfrak L}_l (\phi,I) \ar[r]^(.35){{\mathfrak d}_l}& \PP(\HKC)^l
}
$$
Over a point $([\xi],[\psi]) \in {\mathfrak L}_l (\phi,I) $ lies the line
$\widetilde{{\mathfrak L}}_l (\phi,I)_{([\xi],[\psi])}=\CC \xi$; the map
$F_l$ in the diagram sends that line into the product
$$
\prod^{l-1}_{s=0}H^0(\OO_{Z_{\phi^s([\xi],[\psi])}}) 
$$
according to the rule
$$
F_l(u\xi)=u(F(\phi^s([\xi],[\psi]),\xi))_{s=0,\ldots,l-1} \in \prod^{l-1}_{s=0}H^0(\OO_{Z_{\phi^s([\xi],[\psi])}}), \,\forall u\in \CC; 
$$
for every $s\in I$ the zero divisor $Z_{\phi^s([\xi],[\psi])}$ is reduced
and we use the function $F(\phi^s([\xi],[\psi]),\xi)$ to assign weights to
every point of the divisor:
\begin{equation}\label{weight-Zphir}
w_{{}_{Z_{\phi^s([\xi],[\psi])}}}:	Z_{\phi^s([\xi],[\psi])}\ni p \mapsto exp(|F(\phi^s([\xi],[\psi]),\xi)(p)|) \in \RR_{+}.
\end{equation}
In other words the weight function $w_{{}_{Z_{\phi^s([\xi],[\psi])}}}$ assigned to $Z_{\phi^s([\xi],[\psi])}$ is the exponential of the absolute value $|F(\phi^s([\xi],[\psi]),\xi)|$:
\begin{equation}\label{weightfZphi-r}
	w_{{}_{Z_{\phi^s([\xi],[\psi])}}}=e^{|F(\phi^s([\xi],[\psi]),\xi)|}.
\end{equation}

 Once an ordering of points of $Z_{\phi^s([\xi],[\psi])}$ is chosen we obtain a unique Strebel differential on $C$ with double poles along $Z_{\phi^s([\xi],[\psi])}$ and specified weights.
 	\begin{pro}\label{pro:Strebel xipsir}
 		Assume ${\mathfrak L}_l (\phi,I) $ is nonempty and let $([\xi],[\psi])$ be a point in ${\mathfrak L}_l (\phi,I) $ and
 		consider the value at $([\xi],[\psi])$ of the map ${\mathfrak{d}}_l$ defined in \eqref{l-phi-map}:
 		$$
 		{\mathfrak{d}}_l(([\xi],[\psi]))=([\phi^0([\xi],[\psi])],[\phi^1([\xi],[\psi])],\ldots, [\phi^{l-1}([\xi],[\psi])]) \in (\PP(\HKC))^l.
 		$$
 		1) For every $\xi$ overlying $[\xi]$ and every $s\in I$, the zero divisor $	Z_{\phi^s([\xi],[\psi])}$ comes equipped with the weight function
 		$$
 		w_{{}_{Z_{\phi^s([\xi],[\psi])}}}:	Z_{\phi^s([\xi],[\psi])}\ni p \mapsto exp(|F(\phi^s([\xi],[\psi]),\xi)(p)|) \in \RR_{+}
 		$$
 		as defined in \eqref{weight-Zphir}.
 		
 		2) For $s\in I$, fix a marking of the points of $Z_{\phi^s([\xi],[\psi])}$, that is, a bijective map 
 	$$
 	o_{{}_{Z_{\phi^s([\xi],[\psi])}}}: (1,2,\ldots,2g-2) \longrightarrow Z_{\phi^s([\xi],[\psi])}.
 $$
 Then there exists a unique Strebel differential denoted $\Psi^s_{([\xi],[\psi])}(\xi,o_{{}_{Z_{\phi^s([\xi],[\psi])}}} )$ assigned to the triple
 $(C,o_{{}_{Z_{\phi^s([\xi],[\psi])}}}, w_{{}_{Z_{\phi^s([\xi],[\psi])}}} \circ o_{{}_{Z_{\phi^s([\xi],[\psi])}}}) 
 $:
 $$
 (C,o_{{}_{Z_{\phi^s([\xi],[\psi])}}}, w_{{}_{Z_{\phi^s([\xi],[\psi])}}} \circ o_{{}_{Z_{\phi^s([\xi],[\psi])}}}) \mapsto \Psi^s_{([\xi],[\psi])}(\xi,o_{{}_{Z_{\phi^s([\xi],[\psi])}}} ),
 $$
 that is, $\Psi^s_{([\xi],[\psi])}(\xi,o_{{}_{Z_{\phi^s([\xi],[\psi])}}} )$
 is a unique quadratic differential on $C$ subject to the following conditions:
 
 - $\Psi^s_{([\xi],[\psi])}(\xi,o_{{}_{Z_{\phi^s([\xi],[\psi])}}} )$ is holomorphic on $C\setminus Z_{\phi^s([\xi],[\psi])}$,
 
 - $\Psi^s_{([\xi],[\psi])}(\xi,o_{{}_{Z_{\phi^s([\xi],[\psi])}}} )$ has a pole of order two at every point of $Z_{\phi^s([\xi],[\psi])}$,
 
 - every compact horizontal leaf $\gamma$ is a simple loop circling around
 one of the points of $Z_{\phi^s([\xi],[\psi])}$ and the length of the loop
 with respect to $\Psi^s_{([\xi],[\psi])}(\xi,o_{{}_{Z_{\phi^s([\xi],[\psi])}}} )$ is the the value of the weight function at that point; the last condition means 
 that for every $p\in Z_{\phi^s([\xi],[\psi])}$ we have
 $$
 \int_{\gamma} \sqrt{\Psi^s_{([\xi],[\psi])}(\xi,o_{{}_{Z_{\phi^s([\xi],[\psi])}}} )}
 =w_{{}_{Z_{\phi^s([\xi],[\psi])}}}(p).
 $$
 
 3) There is a unique global section $T^s_{([\xi],[\psi])}(\xi,o_{{}_{Z_{\phi^r([\xi],[\psi])}}} )$ of $\OO_C(4K_C)$ such that 
 $$
 \Psi^s_{([\xi],[\psi])}(\xi,o_{{}_{Z_{\phi^s([\xi],[\psi])}}} )=\frac{T^s_{([\xi],[\psi])}(\xi,o_{{}_{Z_{\phi^s([\xi],[\psi])}}} )}{(\phi^s([\xi],[\psi]))^2}.
 $$
 	\end{pro}
 \begin{pf}
 	All statements come from the theory of Strebel differentials, see \cite{S} and \cite{MuP} for an overview.
 \end{pf}

We also formulate topological and diffeo-geometric consequences of the Strebel theory.
\begin{cor}\label{cor:C-celldecomp}
Let	${\mathfrak L}_l (\phi,I) $ and $([\xi],[\psi])$ be as in Proposition \ref{pro:Strebel xipsir}. For every $s\in I$ choose a marking 
$$
o_{{}_{Z_{\phi^s([\xi],[\psi])}}}: (1,2,\ldots,2g-2) \longrightarrow Z_{\phi^s([\xi],[\psi])}
$$
of the zero divisor $Z_{\phi^s([\xi],[\psi])}$ to obtain the Strebel differential  $\Psi^s_{([\xi],[\psi])}(\xi,o_{{}_{Z_{\phi^s([\xi],[\psi])}}} )$ associated to the triple
$(C,o_{{}_{Z_{\phi^s([\xi],[\psi])}}}, w_{{}_{Z_{\phi^s([\xi],[\psi])}}} \circ o_{{}_{Z_{\phi^s([\xi],[\psi])}}}) 
$. That differential in turn gives rise to the ribbon graph 
$\Gamma_{([\xi],[\psi]),\xi,o_{{}_{Z_{\phi^s([\xi],[\psi])}}}}$ embedded into $C$ so that its complement in $C$ is disjoint union of the open 2-cells 
$$
U_1,U_2,\ldots,U_{2g-2};
$$
 each $U_i$ is an open neighborhood of the point $p_i:=o_{{}_{Z_{\phi^s([\xi],[\psi])}}}(i)$ of $Z_{\phi^s([\xi],[\psi])}$ and $U_i \setminus \{p_i\}$ is swept out by the compact horizontal leaves of $\Psi^s_{([\xi],[\psi])}(\xi,o_{{}_{Z_{\phi^s([\xi],[\psi])}}} )$ circling around $p_i$.
 
 The vertices of the graph $\Gamma_{([\xi],[\psi]),\xi,o_{{}_{Z_{\phi^s([\xi],[\psi])}}})}$ are zeros of the Strebel differential or, equivalently, zeros of the global section
 $T^s_{([\xi],[\psi])}(\xi,o_{{}_{Z_{\phi^s([\xi],[\psi])}}} )$ in Proposition \ref{pro:Strebel xipsir}. In particular, the number of vertices
 of $\Gamma_{([\xi],[\psi]),\xi,o_{{}_{Z_{\phi^s([\xi],[\psi])}}})}$ is
 $deg(4K_C)=8(g-1).$ All those points are distributed among the boundaries
 $\partial \overline{U_i}$'s, where $\overline{U_i}$ is the closure of $U_i$ in the $C^{\infty}$ topology of $C$.  
 
 For each $i$ the boundary $\partial \overline{U_i}$ is a polygon formed by
 the vertices of the graph lying on that boundary and edges being the closures of non compact horizontal leaves: if $\gamma:[0,1] \longrightarrow \partial \overline{U_i}$ is the edge joining a vertex
 $p=\gamma(0)$ to a vertex $p'=\gamma(1)$, then $\gamma:(0,1) \longrightarrow \partial \overline{U_i}$ is a horizontal leaf of the Strebel differential and
 $p$ and $p'$ are the limits of $\gamma(t)$, when $t\to 0$ and $1$ respectively:
 $$
 \lim_{t\to0^+} \gamma (t)=p,\,\,\,\lim_{t\to1^-} \gamma (t)=p'.
 $$
 The number of edges $e$ of the graph $\Gamma_{([\xi],[\psi]),\xi,o_{{}_{Z_{\phi^s([\xi],[\psi])}}})}$ is 
 $e=12(g-1)$. In particular, the graph is $3$-valent or, equivalently, all zeros of the Strebel differential $\Psi^s_{([\xi],[\psi])}(\xi,o_{{}_{Z_{\phi^s([\xi],[\psi])}}} )$ (resp., the global section $T^s_{([\xi],[\psi])}(\xi,o_{{}_{Z_{\phi^s([\xi],[\psi])}}} )$ ) are simple.
\end{cor}
\begin{pf}
	All statements, perhaps with exception of the last one, are standard facts of the theory of Strebel differentials. The last one follows from
	the combinatorial properties of a ribbon graph: for $\Gamma$ a ribbon graph on a curve of genus $g$ with $n$ marked points one has the relation
	$$
	2-2g=v(\Gamma)-e(\Gamma) +n,
	$$
	where $v(\Gamma)$ and $e(\Gamma)$ are respectively the number of vertices and edges of the graph; in addition, since every vertex must have degree at least three one has the upper bound
	$$
	3v(\Gamma) \leq 2e(\Gamma);
	$$
	solving for $v(\Gamma)$ and substituting into the relation one obtains
	$$
	e(\Gamma) \leq 3(2g-2 +n),
	$$
	with equality if and only if $\Gamma$ is trivalent; in our case
	$n=2g-2$ and $e=12(g-1)$, hence the equality holds. Thus the graph
	$\Gamma_{([\xi],[\psi]),\xi,o_{{}_{Z_{\phi^s([\xi],[\psi])}}})}$ is trivalent.
	
	The statement about zeros of the Strebel differential being simple, comes from the well known fact that the degree $d_v$ of the vertex $v$ of the graph associated to the Strebel differential is related to the multiplicity $mult_v$ of zero of the differential at $v$ (remember, the vertices of the graph are the zeros of the Strebel differential) by the formula
	$$
	 mult_v=d_v-2.
	 $$
	 Hence trivalent vertices are precisely the simple zeros of the Strebel differential.
\end{pf}

The ribbon graph $\Gamma=\Gamma_{([\xi],[\psi]),\xi,o_{{}_{Z_{\phi^s([\xi],[\psi])}}})}$ obtained
in the proposition above produces the orbifold cell of the moduli space
${\mathfrak{M}}^{2g-2}_g$ of curves of genus $g$ with $(2g-2)$ marked points. This is given by metrizing the graph: setting $E(\Gamma)$, the set of edges of $\Gamma$ we have
a continuous map
$$
\RR^{E(\Gamma)}_{+} \cong \RR^{12(g-1)}_{+} \longrightarrow {\mathfrak{M}}^{2g-2}_g \times \RR^{2g-2}_{+}.
$$
which maps onto a top dimensional orbifold cell of ${\mathfrak{M}}^{2g-2}_g \times \RR^{2g-2}_{+}$ containing the point $(C,o_{{}_{Z_{\phi^s([\xi],[\psi])}}}, w_{{}_{Z_{\phi^s([\xi],[\psi])}}} \circ o_{{}_{Z_{\phi^s([\xi],[\psi])}}})$ as in Corollary \ref{cor:C-celldecomp}. The following statement summarizes the discussion.
\begin{cor}\label{cor:orbicell-C}
	With the assumptions and notation of Corollary \ref{cor:C-celldecomp},
	for every $s \in I$ the ribbon graph  $\Gamma=\Gamma_{([\xi],[\psi]),\xi,o_{{}_{Z_{\phi^s([\xi],[\psi])}}})}$
	on $C$ determines the continuous map
	$$
	\RR^{E(\Gamma)}_{+} \cong \RR^{12(g-1)}_{+} \longrightarrow {\mathfrak{M}}^{2g-2}_g \times \RR^{2g-2}_{+}.
	$$
	whose image is a top dimensional orbifold cell in ${\mathfrak{M}}^{2g-2}_g \times \RR^{2g-2}_{+}$. The cell contains in its interior the  point corresponding  to $(C,o_{{}_{Z_{\phi^s([\xi],[\psi])}}}, w_{{}_{Z_{\phi^s([\xi],[\psi])}}} \circ o_{{}_{Z_{\phi^s([\xi],[\psi])}}})$.
\end{cor}

\begin{rem}\label{rem:I-collection}
	Assume that the set $I$ in ${\mathfrak{L}}_l (\phi,I)$ has at least two elements, $|I|\geq 2$. Then the above constructions produce the collection of trivalent graphs
	$\Gamma_s:=\Gamma_{([\xi],[\psi]),\xi,o_{{}_{Z_{\phi^s([\xi],[\psi])}}})}
	$, for $s\in I$. For $s$ and $t$, two distinct elements of $I$, one can ask 
	
	- if the poles  $Z_{\phi^s([\xi],[\psi])}$ and  $Z_{\phi^{t}([\xi],[\psi])}$ intersect,
	
	- if the poles $Z_{\phi^s([\xi],[\psi])}$ (resp. $Z_{\phi^{t}([\xi],[\psi])}$) interact with the graph $\Gamma_{t}$ (resp. $\Gamma_s$), that is, if  $Z_{\phi^s([\xi],[\psi])}$ (resp. $Z_{\phi^{t}([\xi],[\psi])}$) intersects the graph $\Gamma_{t}$ (resp. $\Gamma_{s}$).
	
	\vspace{0.2cm}
	This might be an interesting source of combinatorial and geometric relations. One could also envisage such interactions as {\rm morphisms}
	of a new category whose objects are the pairs $(Z_{\phi^s([\xi],[\psi])}, \Gamma_s)$, as $s$ runs through the set $I$.
	
	It should also be observed that when $Z_{\phi^s([\xi],[\psi])}$ and $Z_{\phi^{t}([\xi],[\psi])}$ do not `interact', that is, have no points in common, one can perform the constructions involved in Proposition \ref{pro:Strebel xipsir} and subsequent corollaries
	with the (disjoint) union of $Z_{\phi^s([\xi],[\psi])}$ and $Z_{\phi^{t}([\xi],[\psi])}$. The outcome is a new trivalent graph
	$\Gamma_{s,t}$ with $2(2g-2)$ boundary cycles each winding around a corresponding point of the disjoint union $Z_{\phi^s([\xi],[\psi])}\coprod Z_{\phi^{t}([\xi],[\psi])}$. The graph
	$\Gamma_{s,t}$ could be thought of as a sort of fusion of the graphs 
	$\Gamma_{s}$ and $\Gamma_{t}$.
\end{rem}

We will now combine the maps 
$$
\wp^{\phi}_l(\xi\otimes\psi)(c,\bullet): C\setminus D(\phi,[\xi],[\psi]) \longrightarrow {\mathfrak M}^l_1
$$
obtained in Proposition \ref{pro:pl-C-moduli} with the construction of the Strebel differentials and the corresponding ribbon graphs on $C$. The outcome will be a map of that ribbon graph into the moduli space ${\mathfrak M}^l_1$. 

 The construction will assume that the polar divisor $D(\phi,[\xi],[\psi])$ in 
Proposition \ref{pro:pl-C-moduli} is {\it reduced}. This means that in the
decomposition
\begin{equation}\label{polar-disjoint}
	D(\phi,[\xi],[\psi])=\sum^{l-1}_{s=0} Z_{\phi^s([\xi],[\psi])}
\end{equation}
each zero divisor $Z_{\phi^s([\xi],[\psi])}=(\phi^s([\xi],[\psi])=0)$ is reduced and they are mutually disjoint; in other words, we are in the case when there is no interaction between the poles $Z_{\phi^s([\xi],[\psi])}$'s in a sens of Remark \ref{rem:I-collection}. Thus we can `fuse' all the graphs $\Gamma_s$'s to obtain one single graph: the following construction will serve as an illustration of the (multiple) fusion operation with trivalent graphs alluded to in the remark above.

 To implement the construction, we put together the weight functions
$$
w_{{}_{Z_{\phi^s([\xi],[\psi])}}}=e^{|F(\phi^s([\xi],[\psi]),\xi)|}
$$
in \eqref{weightfZphi-r} to define the weight function
\begin{equation}\label{weight-polar}
w_{D(\phi,[\xi],[\psi])}: D(\phi,[\xi],[\psi]) \longrightarrow \RR_{+}
\end{equation}
on the polar divisor: this assigns to every point $p \in D(\phi,[\xi],[\psi])$ the value
$$
w_{D(\phi,[\xi],[\psi])}(p):=w_{{}_{Z_{\phi^{s_p}([\xi],[\psi])}}} (p),
$$
where $Z_{\phi^{s_p}([\xi],[\psi])}$ is a unique summand in the decomposition \eqref{polar-disjoint} containing $p$. To obtain the Strebel differential, it remains to fix an ordering of points in $D(\phi,[\xi],[\psi])$: the components in the decomposition of $D(\phi,[\xi],[\psi])$ come with the prescribed order
$$
Z_{\phi^{l-1}([\xi],[\psi])}, Z_{\phi^{l-2}([\xi],[\psi])}, \ldots, Z_{\phi^{1}([\xi],[\psi])}, Z_{\phi^{0}([\xi],[\psi])}
$$
given by the (white/black) vertices of $\widehat{PG}_l$, so the ordering of points is determined by a choice of bijection
$$
o_{{}_{Z_{\phi^{s}([\xi],[\psi])}}}:\{1,2,\ldots, 2g-2\} \longrightarrow Z_{\phi^{s}([\xi],[\psi])},
$$
for every $s\in [0,l-1]$. With those fixed we obtain the order on the polar divisor $D(\phi,[\xi],[\psi])$:
\begin{equation}\label{o-polar}
o_{{}_{D(\phi,[\xi],[\psi])}}:\coprod^{l-1}_{s=0}\{1+s,2+s,\ldots, 2g-2 +s \} \longrightarrow D(\phi,[\xi],[\psi])=\coprod^{l-1}_{s=0} Z_{\phi^{s}([\xi],[\psi])}
\end{equation}
defined by the formula
$$
o_{{}_{D(\phi,[\xi],[\psi])}}(j+s)=o_{{}_{Z_{\phi^{s}([\xi],[\psi])}}}(j),
$$
for every $j\in[1,2g-2]$ and every $s\in[0,l-1]$.
 \begin{pro}\label{pro:St-Graph-polar}
 	Assume the set of points in ${\mathfrak{L}}_l(\phi)$ having the polar divisor reduced is nonempty; denote that set by ${\mathfrak{L}}^{red}_l(\phi)$. Let $\widetilde{{\mathfrak{L}}}^{red}_l(\phi)$ denote the restriction of $p^{\ast}_1 \left(\OO_{\PP(H^1(\Theta_{C}))}(-1)\right)$ to ${\mathfrak{L}}^{red}_l(\phi)$.
 	Then the map
 	$$
 	F_l:\widetilde{{\mathfrak{L}}}^{red}_l(\phi) \longrightarrow (q_{1\ast} \OO_{{\cal Z}_C} )^l
 	$$
 	in Proposition \ref{pro:LphiI} provides the weight function
 	$$
 	w_{D(\phi,[\xi],[\psi])}(\xi): D(\phi,[\xi],[\psi]) \longrightarrow \RR_{+}
 	$$
 	for every $([\xi],[\psi])$ in ${\mathfrak{L}}^r_l(\phi)$ and every nonzero $\xi \in \widetilde{{\mathfrak{L}}}^{red}_l(\phi) $ lying over $([\xi],[\psi])$, see \eqref{weight-polar} and the paragraph that follows. With an ordering
 	$o_{{}_{D(\phi,[\xi],[\psi])}}$ of points in $D(\phi,[\xi],[\psi])$ chosen as in \eqref{o-polar}, one obtains a unique Strebel differential
 	$\Psi^{D(\phi,[\xi],[\psi])}(\xi,o_{{}_{D(\phi,[\xi],[\psi])}})$ associated to the triple
 	$(C,o_{{}_{D(\phi,[\xi],[\psi])}}, w_{D(\phi,[\xi],[\psi])} (\xi)\circ o_{{}_{D(\phi,[\xi],[\psi])}} )$. This has the form
 	\begin{equation}\label{Strebel-polar}
 	\Psi^{D(\phi,[\xi],[\psi])}(\xi,o_{{}_{D(\phi,[\xi],[\psi])}})=T^{D(\phi,[\xi],[\psi])}(\xi,o_{{}_{D(\phi,[\xi],[\psi])}})\left(\prod_{s=0}^{l-1}\phi^s([\xi],[\psi])\right)^{-2},
 \end{equation}
 	where $T^{D(\phi,[\xi],[\psi])}(\xi,o_{{}_{D(\phi,[\xi],[\psi])}})$ is a uniquely determined global section of $\OO_C (2(l+1)K_C) $.
 	
 	The Strebel differential in \eqref{Strebel-polar} determines the ribbon
 	graph $\Gamma_{D(\phi,[\xi],[\psi]),\xi,o_{{}_{D(\phi,[\xi],[\psi])}}}$ subject to the following properties:
 	
 	- the complement of the graph in $C$ is the disjoint union of open disks
 	$$
 	U_1,U_2, \ldots, U_{l(2g-2)},
 	$$
 	where for each $j\in [1,l(2g-2)]$ the disk $U_j$ is centered at the point $p_j:=o_{{}_{D(\phi,[\xi],[\psi])}}(j)$ of $D(\phi,[\xi],[\psi])$,
 	
 	- the vertices of the graph are zeros of the Strebel differential or, equivalently, of the global section $T^{D(\phi,[\xi],[\psi])}(\xi,o_{{}_{D(\phi,[\xi],[\psi])}})$; in particular, the graph has $v_{\Gamma} =4(l+1)(g-1)$ vertices,
 	
 	- the edges of the graph form the boundaries of the disks in the first item; their number $e_{\Gamma}=6(l+1)(g-1)$.
 	
 	In particular, the graph $\Gamma_{D(\phi,[\xi],[\psi]),\xi,o_{{}_{D(\phi,[\xi],[\psi])}}}$ is trivalent and it determines a top dimensional orbifold cell in
 	${\mathfrak{M}}^{l(2g-2)}_g \times \RR^{l(2g-2)}_{+}$.
 \end{pro}

With the additional structure of the ribbon graph $\Gamma:=\Gamma_{D(\phi,[\xi],[\psi]),\xi,o_{{}_{D(\phi,[\xi],[\psi])}}}$ on $C$, we return to the map
$$
\wp^{\phi}_l(\xi\otimes\psi)(c,\bullet): C\setminus D(\phi,[\xi],[\psi]) \longrightarrow {\mathfrak M}^l_1
$$
obtained in Proposition \ref{pro:pl-C-moduli}. Since the polar divisor is in the complement of the graph, the restriction of $\wp^{\phi}_l(\xi\otimes\psi)(c,\bullet)$ to the graph gives a well defined function
$$
\wp^{\phi}_l(\xi\otimes\psi)(c,\bullet)\big|_{\Gamma}: \Gamma \longrightarrow {\mathfrak M}^l_1
$$
From this we deduce the following.
\begin{cor}\label{cor:ribbon-M1l}
	Let $([\xi],[\psi])$ be a point of ${\mathfrak{L}}_l(\phi)$ with the polar divisor $D(\phi,[\xi],[\psi])$ reduced. Let 
	$
	\Gamma:=\Gamma_{D(\phi,[\xi],[\psi]),\xi,o_{{}_{D(\phi,[\xi],[\psi])}}}
	$
	be the trivalent ribbon graph in Proposition \ref{pro:St-Graph-polar} and let
	$$
	\wp^{\phi}_l(\xi\otimes\psi)(c,\bullet)\big|_{\Gamma}: \Gamma \longrightarrow {\mathfrak M}^l_1
	$$
	be the restriction to $\Gamma$ of the map 
	$$
	\wp^{\phi}_l(\xi\otimes\psi)(c,\bullet): C\setminus D(\phi,[\xi],[\psi]) \longrightarrow {\mathfrak M}^l_1
	$$
	obtained in Proposition \ref{pro:pl-C-moduli}. Then the following hold:
	
	1) to every vertex $v$ of $\Gamma$ is attached 
	 the elliptic curve $E_{\wp^{\phi}_l(\xi\otimes\psi)(c,v)}$ with $l$ marked points, the image $\wp^{\phi}_l(\xi\otimes\psi)(c,v)$ of $\wp^{\phi}_l(\xi\otimes\psi)(c,\bullet)$ at $v$; 
	 in particular, to the triple 
	 $$
	 (C,\wp^{\phi}_l(\xi\otimes\psi)(c,\bullet), \Gamma_{D(\phi,[\xi],[\psi]),\xi,o_{{}_{D(\phi,[\xi],[\psi])}}})
	 $$
	 is attached the abelian variety
	 $$
	 A^{\phi}(\xi\otimes\psi,\Gamma):= \prod_{v\in V_{\Gamma}}E_{\wp^{\phi}_l(\xi\otimes\psi)(c,v)},
	 $$
	 where the product is taken over the set of vertices $V_{\Gamma}$ of the graph $\Gamma$;
	  
	2) every boundary cycle $\partial U_i$ of $\Gamma$ give rise to the loop
	 $$
	\sigma_i:= \wp^{\phi}_l(\xi\otimes\psi)(c,\bullet)\big|_{\Gamma} (\partial U_i)
	$$
	in ${\mathfrak M}^l_1$ passing through the images of the vertices of $\Gamma$ lying on the boundary cycle $\partial U_i$; here $U_i$'s are the disks in Proposition \ref{pro:St-Graph-polar}.	   
\end{cor}

\begin{rem}\label{rem:singpl-C-polar}
	There is another interesting issue related to the maps
	$$
	\wp^{\phi}_l(\xi\otimes\psi)(c,\bullet): C\setminus D(\phi,[\xi],[\psi]) \longrightarrow {\mathfrak M}^l_1:
	$$
	what happens when we approach to the points of the polar divisor? In other words, the study of singularities of the map  	$\wp_l([\xi],[\psi])(c,\bullet)$.
	
	 Consider the Deligne-Mumford compactification
	$\overline{{\mathfrak M}^l_1}$ of ${\mathfrak M}^l_1$ and view the above
	map as the map to $\overline{{\mathfrak M}^l_1}$:
		$$
	\wp^{\phi}_l(\xi\otimes\psi)(c,\bullet): C\setminus D(\phi,[\xi],[\psi]) \longrightarrow {\mathfrak M}^l_1 \subset \overline{{\mathfrak M}^l_1}.
	$$
	At a point $p$ in the support of $D(\phi,[\xi],[\psi])$, the map may have the limit and then we can extend $\wp^{\phi}_l(\xi\otimes\psi)(c,\bullet)$ continuously through $p$ and thus remove that point from the list of singularities of $\wp^{\phi}_l(\xi\otimes\psi)(c,\bullet)$. Otherwise, a point $p$ is a {\rm non-removable} singularity of the map and one may wish to
	understand the set of limits in $\overline{{\mathfrak M}^l_1}$ of $\wp^{\phi}_l(\xi\otimes\psi)(c,\bullet)$ as one approaches such a non-removable
	singular point. In other words, in the Cartesian product
	$C\times \overline{{\mathfrak M}^l_1}$ consider the graph of the map
	$\wp^{\phi}_l(\xi\otimes\psi)(c,\bullet)$:
	$$
	I_{\wp^{\phi}_l(\xi\otimes\psi)(c,\bullet)}=\{(x,\wp^{\phi}_l(\xi\otimes\psi)(c,x)\in C\times \overline{{\mathfrak M}^l_1}\,\, \big | \,\,\forall x\in C\setminus D(\phi,[\xi],[\psi])  )\};
	$$
	take its closure $\overline{I}_{\wp^{\phi}_l(\xi\otimes\psi)(c,\bullet)}$ in the usual topology of $C\times \overline{{\mathfrak M}^l_1}$; the restriction of projections on each factor give the diagram
	$$
	\xymatrix{
	&\overline{I}_{\wp^{\phi}_l(\xi\otimes\psi)(c,\bullet)}\ar[rd]^(.6){p_{\overline{{\mathfrak M}}^l_1}}\ar[ld]_(.55){p_C}&\\
C&&\overline{{\mathfrak M}}^l_1
}
$$
and we are asking for the properties of the fibres of the map $p_C$ over the points $p$ in the polar divisor $D(\phi,[\xi],[\psi])$; the image of those under the map $p_{\overline{{\mathfrak M}}^l_1}$ give the sets of limiting values of $\wp^{\phi}_l(\xi\otimes\psi)(c,\bullet)$ as we approach the points of the polar divisor. Those limiting sets and their invariants might be interesting invariants of the the pair $(C,D(\phi,[\xi],[\psi]))$.  
\end{rem}

\vspace{0.2cm}
\noindent
{\bf The cyclic automorphism of $\widehat{PG}_l$.}
The graph $\widehat{PG}_l$ has an additional structure: it is equipped with the automorphism
$$
\sigma_l: \widehat{PG}_l \longrightarrow \widehat{PG}_l
$$
acting on the vertices by cyclic permutation
$$
\sigma_l (i)=i+1, \, \sigma_l (i')=(i+1)',
$$
where the indices $i$ are taken modulo $l$ with values in the set 
$\{0,1,\ldots, l-1\}$. This determines the permutation on the edges
$$
\sigma_l(e^0_i) =e^0_{i+1}, \, \sigma_l(e^{\pm}_i) =e^{\pm}_{i+1}
$$
as well as on the boundary cycles
$$
\sigma_l(B^i)=B^{i+1}.
$$
The metrics on $\widehat{PG}_l$ invariant with respect to the action of $\sigma_l$ are the functions
$$
m: E_{\widehat{PG}_l} \longrightarrow \RR_{+}
$$
which are constant on the edges of the same color $c\in \{-,0,+\}$:
$$
m (e^c_i)=m^c, \,\forall i\in [0,l-1].
$$
Thus we obtain
\begin{lem}\label{lem:metricinv}
	Under the identification
	$$
	\RR^{E_{\widehat{PG}_l}} \cong \RR^{3l}_{+}
	$$
	the subset of metrics on $\widehat{PG}_l$ invariant with respect to the action of $\sigma_l$ is identified with the subset of vectors
	$$
	\left(\RR^{E_{\widehat{PG}_l}}\right)^{\sigma_l}= \{(\underbrace{\underbrace{x,y,z},\underbrace{x,y,z},\ldots,\underbrace{x,y,z}}_{\text{$l$ times}}) \big| \, x,y, z \in \RR_{+}\} \cong \RR^3_{+}.
	$$  
\end{lem}

Let $m$ be a   $\sigma_l$-invariant metric on $\widehat{PG}_l$, then the complex projective curve $\Gamma_{\widehat{PG}_l, m}$ with marked points $(x_0, x_1, \ldots, x_{l-1})$, associated to $m$, comes with an automorphism
$$
f_{\sigma_l, m} : \Gamma_{\widehat{PG}_l, m} \longrightarrow \Gamma_{\widehat{PG}_l, m}
$$
which induces the cyclic permutation
$$
f_{\sigma_l, m}(x_i)=x_{i+1},
$$
with the understanding that $i$ is taken modulo $l$ with values in $[0,l-1]$.

\begin{pro}\label{pro:torsion}
	Let $m$ be a $\sigma_l$-invariant metric on $\widehat{PG}_l$ and let  $f_{\sigma_l, m} : \Gamma_{\widehat{PG}_l, m} \longrightarrow \Gamma_{\widehat{PG}_l, m}$ be the corresponding automorphism of $\Gamma_{\widehat{PG}_l, m}$. Then $f_{\sigma_l, m}$ has order $l$, that is, $l$ is the smallest positive integer
	subject to
	$$
	f^l_{\sigma_l, m}=id_{\Gamma_{\widehat{PG}_l, m}}.
	$$
	Let $\{x_0,\ldots,x_{l-1}\}$ be the marked points of the genus one curve $\Gamma_{\widehat{PG}_l, m}$. Fix the point $x_0$ thus identifying $\Gamma_{\widehat{PG}_l, m}$ with its Jacobian
	$$
	\Gamma_{\widehat{PG}_l, m}\cong \CC/\Lambda_{\tau_m},
	$$
	where $\Lambda_{\tau_m}=\ZZ \oplus \ZZ \tau_{m}$ is the period lattice with
	the normalized periods $\{1, \tau_{m}\}$ with $\tau_{m}$ lying in the upper half plane.
	Then from the assumption $l\geq 3$, it follows:
	
	- either $f_{\sigma_l, m}$ is a translation by a torsion point 
	$p$ of order $l$ in $\CC/\Lambda_{\tau_m}$ and that point corresponds to the point $((x_1)-(x_0))$ of the Jacobian of $\Gamma_{\widehat{PG}_l, m}$; in particular, the marked points
	$\{x_0, x_1, \ldots, x_{l-1}\}$ correspond to the torsion points
	$$
	\{j p | j=0,1,\ldots, (l-1)p\} \subset \CC/\Lambda_{\tau_m};
	$$
	
	- or $l=3$ and $f_{\sigma_l, m}$, up to a translation in $\CC/\Lambda_{\tau_m}$, is the complex multiplication by the primitive third root of unity $\omega$. In addition, there is a $2\times 2$ integer matrix  $\begin{pmatrix}
		a&b\\
		c&d
	\end{pmatrix}$ of determinant $\pm1$, such that the period $\tau_{m} =\frac{\omega -d}{c} \in {\mathbb{Q}[\omega]}$; in addition, the integer
$d$ is relatively prime to $c$ and $b$, while
 $(d^2 +d+1)$ is divisible by $b$ and $c$; 
 the following facts about this curve are well known:
 
$\bullet$ the field 
 $$
 End_{\mathbb Q} (\CC/\Lambda_{\tau_m})={\mathbb Q}(\omega),
 $$

$\bullet$  the ring of endomorphism 
 $End (\CC/\Lambda_{\tau_m})=\ZZ[\omega]$, the ring of integers of the field ${\mathbb Q}(\omega)$,
 
 $\bullet$  the group of automorphisms 
 $$
 Aut (\CC/\Lambda_{\tau_m})=\left(\ZZ[\omega]\right)^{\times},
 $$
 the group of units of $\ZZ[\omega]$, that group consists of six elements
 $$
 Aut (\CC/\Lambda_{\tau_m})=\left(\ZZ[\omega]\right)^{\times}=\{\pm1,\pm\omega,\pm\omega^2\}, 
 $$
 
 $\bullet$ the $j$-invariant of $\CC/\Lambda_{\tau_m}$ is zero; in particular, the curve $\CC/\Lambda_{\tau_m}$ is defined over $\mathbb{Q}$.
\end{pro}
\begin{pf}
	The automorphism $\sigma_l: \widehat{PG}_l \longrightarrow \widehat{PG}_l$ of the ribbon graph underlying the automorphism $f_{\sigma_l, m}$ has degree $l$. Hence the same holds for $f_{\sigma_l, m}$. 
	To identify $f_{\sigma_l, m}$ precisely we make the identification
	$$
	\Gamma_{\widehat{PG}_l, m}\cong \CC/\Lambda_{\tau_m}
	$$
	and write $f_{\sigma_l, m}$ as a composition of an invertible 
	endomorphism of $\CC/\Lambda_{\tau_m}$, call it $h$, together with a translation by a point of $\CC/\Lambda_{\lambda}$, call that point $p$,
	$$
	f_{\sigma_l, m}(x) =h(x)+p, \forall x\in \CC/\Lambda_{\tau_m}.
	$$
Taking the $l$-th power we obtain
$$
x=h^l(x) +h^{l-1}(p)+\cdots+h(p)+p,
$$
for all $x$ in $\CC/\Lambda_{\lambda}$. Evaluating at $x=0$, implies
\begin{equation}\label{c-a}
	h^{l-1}(p)+\cdots+h(p)+p=0,
\end{equation}
 and the complex multiplication $h$ is of order at most $l$, that is, $h$ is defined by the multiplication by an $l$-th root of unity. Since the multiplication
by $h$ sends the lattice $\Lambda_{\tau_m}=\ZZ\oplus \ZZ \tau_{m}$ to itself, we have
\begin{equation}\label{h-tau}  
	h\tau_{m}=a\tau_m+b,\,\,h=c\tau_{m}+d,
\end{equation}
where the matrix
$$
\begin{pmatrix}
	a&b\\
	c&d
\end{pmatrix}
$$
is integral and unimodular. This implies that the field $\mathbb{Q}(\tau_{m})$ is an imaginary quadratic extension of $\mathbb{Q}$ and $h$ lies in that field. This tells us that $h$ is a root of unity of degree at most $3$, that is, $h=\pm1$ or $h$ is a primitive third root of unity. 

The case $h=1$ and the equation \eqref{c-a} tell us that $f_{\sigma_l, m}$ is the translation by the point $p$ and that point is of order $l$ on the elliptic curve $\CC/\Lambda_{\tau_m}$.  If $h=-1$, then we obtain that $f_{\sigma_l, m}$ has order $l=2$. But this contradicts the assumption $l\geq 3$.
Thus the remaining case $h=\omega$, a primitive third root of unity. In this case the order of $f_{\sigma_l, m}$ is $3$ and from the equations in \eqref{h-tau} we have
$$
\tau_{m}=\frac{\omega-d}{c}, \,\, a=-1-d,\,\, bc=-(1+d+d^2).
$$

Once we know the multiplication by $\omega$ is an automorphism of $\CC/\Lambda_{\tau_m}$, it follows that
$$
Aut(\CC/\Lambda_{\tau_m}) \supset \{\pm1,\pm\omega,\pm\omega^2\} 
$$
and the theory of curves with complex multiplication implies the equality as well as all the facts listed in the proposition, see \cite{Sh}.
\end{pf}
\begin{rem}\label{rem:constmetric}
	A ribbon graph can always be thought as {\rm metrized} by assigning length one
	to all of its edges. If all vertices of such a graph have degree at least $3$, it is known that the corresponding curve is a Belyi curve, that is, has a morphism onto $\PP^1$ branched precisely over three points. By a remarkable theorem of Belyi, such a curve is defined over $\overline{\mathbb{Q}}$, the field of algebraic numbers over $\mathbb{Q}$, see \cite{MuP}. In the case of the graph $\widehat{PG}_l$,
assigning to it a constant metric, that is, a metric
$$
m_a: E_{\widehat{PG}_l} \longrightarrow \RR_{+}
$$
for which $m_a(e)=a >0$, for all $e\in E_{\widehat{PG}_l} $, gives rise to
an extra automorphism, say $\upsilon$, of the metrized graph $(\widehat{PG}_l, m_a) $: it acts on vertices by the rule 
$$
\upsilon: i\rightarrow (i-1)', \,\, i' \rightarrow (i-1), \,\,\forall i\in[0,l-1],
$$
on edges
$$
\upsilon: e^0_i \rightarrow -e^0_{i-1},\,\, e^-_{i} \rightarrow -e^+_{i-2},\,\,e^+_{i}\rightarrow -e^+_{i},\,\, \forall i\in[0,l-1],
$$
and it sends a boundary cycle $B_i$ to the boundary cycle $B_{i-1}$ rotated by the angle $\frac{2\pi}{6} =\frac{\pi}{3}$; 

\noindent
this implies that after factoring out by the cyclic group generated by the automorphism $f_{\sigma_l,m_a}$, we obtain an elliptic curve  having an automorphism given by the multiplication by $\rho=e^{\frac{\pi \sqrt{-1}}{3}}$. Hence, that curve is isomorphic to the curve in the last part of Proposition \ref{pro:torsion}; thus a constant metric $m_a$ on $\widehat{PG}_l$ define  the elliptic curve $\Gamma_{\widehat{PG}_l,m_a}$
which is either isomorphic to the curve $E_0$ in the last part of Proposition \ref{pro:torsion} or $f_{\sigma_l,m_a}$ is a translation by a torsion point $p_a$ of order $l$ on $\Gamma_{\widehat{PG}_l,m_a}$ such that the quotient
$\Gamma_{\widehat{PG}_l,m_a}/\langle p_a \rangle$ is isomorphic to $E_0$; the notation $\langle p_a \rangle$ stands for the group generated by $p_a$.  
\end{rem}

\begin{cor}\label{cor:torsion}
	Assume $l\geq 4$. Then every $\sigma_l$-invariant metric $m$ on $\widehat{PG}_l$ defines the elliptic curve $\Gamma_{\widehat{PG}_l,m}$ together with a distinguished point $p_{m}$ of order $l$. This gives a continuous map
	$$
	\left(\RR^{E_{\widehat{PG}_l}}\right)^{\sigma_l}=	\RR^3_{+} \longrightarrow Y_1 (l)=\Gamma_1 (l) \diagdown {\mathfrak H} 
	$$
	where ${\mathfrak H}$ is the upper-half plane and $\Gamma_1 (l)$ is the congruence subgroup
	$$
	\Gamma_1 (l)=\left\{ \begin{pmatrix}
		a&b\\
		c&d
	\end{pmatrix} \in SL_2 (\ZZ) |\, c\equiv 0(mod \,l), a\equiv d \equiv 1(mod \,l) \right\}.
$$
\end{cor}

It should be clear that it is interesting to have $\sigma_l$-invariant metrics on $\widehat{PG}_l$. To construct those we can start with any metric
$$
m: E_{\widehat{PG}_l} \longrightarrow \RR_{+}
$$
Remember the edges come in three colors, `$0$', `$+$' and `$-$'. Let $c$ be one of these colors and let $E^c_{\widehat{PG}_l}$ be the subset of edges of 
$E_{\widehat{PG}_l}$ having color $c$:
$$
 E^c_{\widehat{PG}_l} =\{e^c_i \,| \, i\in [0,l-1]\}.
 $$
 The restriction of the metric $m$ to $E^c_{\widehat{PG}_l}$ will be denoted $m^c$. Symmetrizing those functions by taking for example the sum of powers
 $$
 S^{a_c}(m^c):= \sum^{l-1}_{i=0} (m (e^c_i))^{a_c},
 $$
 for $a_c\in \NN$, gives us the $\sigma_l$-invariant metric
 $$
 S^{\bf a}(m):E_{\widehat{PG}_l} \longrightarrow \RR_{+}
 $$
 for every ${\bf a}=(a_0,a_{+},a_{-}) \in \NN^3$, defined by the equation
 $$
 S^{\bf a}(m)(e^c_i) =S^{a_c}(m^c), \,\forall i \in [0,l-1], \forall c\in \{0,+,-\}.
 $$
 More generally, consider the set of {\it positive} symmetric polynomials $s$ in $l$ variables
 {\small$(t_0,\cdots,t_{l-1})$}, that is, $s$ is symmetric in variables $t_0,\ldots,t_{l-1}$ and takes positive values on $\RR^l_{+}$: 
 \begin{equation}
 	s(x_0,\ldots,x_{l-1}) >0, \forall  x_0,\ldots,x_{l-1} \in \RR_{+}. 
 \end{equation}
Then we have the following.
\begin{lem}\label{lem:symm}
	Choose any three positive symmetric polynomials and denote them $s_0, s_{+}, s_{-}$ according to the colors of the edges of $\widehat{PG}_l$.
	Then any metric $m:E_{\widehat{PG}_l} \longrightarrow \RR_{+}$
	admits $(s_0, s_{+}, s_{-})$-symmetrization, which assigns to the edges
	$e^c_i$ of color $c$ the length
	$$
	s_c(m):=s_{c}(m(e^c_0), \ldots, m(e^c_{l-1})).
	$$
	This gives a $\sigma_l$-invariant metric on $\widehat{PG}_l$ which takes the value
	$(s_0(m), s_{+}(m), s_{-}(m)) \in \RR^3_{+}$. The $(s_0, s_{+}, s_{-})$-symmetrization of a metric $m$ will be denoted $m_{s_0, s_{+}, s_{-}}$; in case $s_0= s_{+}= s_{-}=s$, the notation is simplified to $m_s$.
\end{lem}

Recall that standard bases for symmetric functions, such as monomial, elementary, power functions, are given by positive symmetric functions and labeled by partitions: for example, given a partition
$$
\eta=(\eta_1,\eta_2, \ldots,\eta_n) 
$$
with $n$ parts one associates the symmetric polynomial
$$
m_{\eta} =\sum x^{k_1}_1 x^{k_2}_2 \cdots x^{k_n}_n,
$$
where the sum is taken over $k=(k_1,k_2, \ldots,k_n)$, all distinct permutations of $\eta$; this is the monomial symmetric polynomial associated to the partition $\eta$; 

the elementary symmetric polynomial associated to the partition $\eta$:
$$
e_{\eta}=e_{\eta_1}e_{\eta_2} \cdots e_{\eta_n},
$$
where the factors $e_{k}$ are
$$
e_k= \sum_{i_1 <\cdots<i_k} x_{i_1} \cdots x_{i_k};
$$

the power symmetric polynomial associated to the partition $\eta$:
$$
p_{\eta}=p_{\eta_1}p_{\eta_2} \cdots p_{\eta_n},
$$
where the factors $p_{k}$ are
$$
p_k= \sum_{i} x^k_{i};
$$
see \cite{Mac} for in depth treatment of symmetric functions. 

With the above reminder, we can use partitions to label symmetrizations of metrics of $\widehat{PG}_l$.
\begin{cor}\label{cor:metric-partitions}
	Choose three partitions of length $l$ and label them
	$\eta^0, \eta^{+}, \eta^{-}$ according to the colors of the edges of $\widehat{PG}_l$.
	Then any metric $m:E_{\widehat{PG}_l} \longrightarrow \RR_{+}$
	admits $(s_{\eta^0}, s_{\eta^{+}}, s_{\eta^{-}})$-symmetrization as in Lemma \ref{lem:symm}, where $s_{\eta^c}$ is a positive symmetric function in variables $(t_0,t_1,\ldots, t_{l-1})$ associated to the partition $\eta^c$. 
	This gives a $\sigma_l$-invariant metric on $\widehat{PG}_l$ associated to $\eta^0, \eta^{+}, \eta^{-}$ and it will be denoted 
	$m_{s_{\eta^0}, s_{\eta^{+}}, s_{\eta^{-}}}$; it takes the value
	$$
	(s_{\eta^0}(m), s_{\eta^{+}}(m), s_{\eta^{-}}(m)) \in \RR^3_{+}.
	$$ 
	In case the symmetric functions associated to partitions are specified
	the notation is simplified to $m_{{\eta^0},{\eta^{+}},{\eta^{-}}}$; if three partitions are all equal $\eta^0= \eta^{+}= \eta^{-}=\eta$ the associated symmetrized metric is denoted by $m_{s_{\eta}}$ or $m_{\eta}$, if there is no ambiguity on the type of symmetric functions.   
\end{cor}

From now on we agree that whenever a metric on $\widehat{PG}_l$ is symmetrized it means that a basis of {\it positive symmetric functions} labeled by partitions is chosen
$$
Sym:=\{s_{\eta}\}.
$$
Such a choice will be called a {\it symmetrization pattern for $\widehat{PG}_l$}.

Next recall that the stratum ${{\mathfrak{L}}_l }$ admits a stratification in its own turn: the substrata are ${\mathfrak{L}}_l(h^l,\lambda)$, where the dimensions of the summands $\{P^r([\xi],[\phi])\}$ of the orthogonal decomposition
$$
W_{\xi} =\bigoplus^l_{r=0} P^r([\xi],[\phi])
$$
are constant and the ones with $r\in [0,l-1]$ are arranged in the Young diagram $\lambda$. Recall that to such a stratum we associate the Hilbert-Laurent polynomial
\begin{equation}\label{LH-lambda}
h_{{\mathfrak{L}}_l(h^l,\lambda)}(q)=\sum^{l-1}_{s=0}dim(P^s([\xi],[\phi]))q^{(-1)^s s}.
\end{equation}
Also recall the identification
$$
Hom_{\CC}(\CC[q,q^{-1}], \CC) \cong \CC[q,q^{-1}]
	$$
	via the pairing defined by the residue at zero
	$$
	(f,g)=Res_{0}(fg).
	$$
Thus we think of the Laurent polynomials $h_{{\mathfrak{L}}_l(h^l,\lambda)}$  as linear functionals on $\CC[q,q^{-1}]$ and obtain the associated maps 
$$
\widetilde{\wp}_{l} (\bullet,c,h_{{\mathfrak{L}}_l(h^l,\lambda)},h_{{\mathfrak{L}}_l('h^l,\mu)}): \OO^{\times}_{\mathfrak{L}_l} (-1) \longrightarrow \mathfrak{H},
$$
for any pair $h_{{\mathfrak{L}}_l(h^l,\lambda)},h_{{\mathfrak{L}}_l('h^l,\mu)}$ of Laurent polynomials occurring in the stratification of ${\mathfrak{L}}_l$. To simplify the notation, we denote those maps by $\widetilde{\wp}_{l} (\bullet,c,\lambda,\mu)$; when $(h^l,\lambda)=('h^l,\mu)$, this is simplified further to $\widetilde{\wp}_{l} (\bullet,c,\lambda)$.

In addition, with a symmetrization pattern for $\widehat{PG}_l$ {\it chosen}, the Young diagrams $\lambda$ provide the symmetrizations for the metrics. Let us recall that the Young diagram $\lambda$ has $l$ columns. This means that the Young diagram (and the corresponding partition) $\lambda'$ {\it conjugate} to $\lambda$ has $l$ parts, with the parts
$$
\lambda'_i =dim(W^i_{\xi}([\phi])/W^{i+1}_{\xi}([\phi])),
$$
for $i\in[0,l-1]$. Choosing all three partitions in Corollary \ref{cor:metric-partitions} to be equal to $\lambda'$ we obtain the following.
\begin{pro}\label{pro:plmap-ltorsion}
	Assume $l\geq 4$ and let ${\mathfrak{L}}_l(h^l,\lambda)$ be the stratum of ${\mathfrak{L}}_l$ parametrizing points $([\xi],[\phi])$, where the $([\xi], [\phi])$-filtrations of $W_{\xi}$
	have length $l$ and the dimensions of the graded pieces are constant: $h^l$ for $W^l_{\xi}([\phi])$ and the conjugate partition $\lambda'$ for all other pieces. Choose a symmetrization pattern for $\widehat{PG}_l$ which associates the positive symmetric function $s_{\lambda'}$ to $\lambda'$. Then we have a continuous map
	$$
	\wp_{l,\lambda,c}:\OO^{\times}_{{\mathfrak{L}}_l}(-1) \longrightarrow Y_1(l)
	$$
	which sends  a point $\xi\otimes\phi$ in the domain of the map
	to the isomorphism class of the elliptic curve $\Gamma_{{m(c, h_{\mathfrak L(h^l,\lambda)},\xi\otimes\phi)}_{s_{\lambda'}}} $  with a distinguished torsion point of order $l$.
\end{pro}

The appearance of the modular curve $Y_1(l)$ is very much in line with the ties of generating functions of Gromov-Witten invariants with modular forms. The essential ingredient in our construction was a symmetrization of
metrics of the graph $\widehat{PG}_l$ and this consists of a triple of partitions of length $l$: given a triple of partitions $(\eta^0,\eta^+, \eta^-)$ with $l$ parts allows us to go from the space of metrics
$\RR^{E_{\widehat{PG}_l}}_+$ to $\RR^3_+$. To be more precise we think of the axes of $\RR^3_+$ as corresponding to the three colors $\{+,0,-\}$ of edges
of $\widehat{PG}_l$, that is the coordinate axes of $\RR^3_+$  are labeled
$(x_+,x_0,x_{-})$, and to remember how to symmetrize, we attach the partition of the same color to each axis:
$$
x_{\pm} \rightarrow \eta^{\pm},\,\,x_0\rightarrow \eta^{0}.
$$  
The data of the positive octant $\RR^3_+$ with a partition attached to each of its axes might be familiar: the data give rise to a three-dimensional partitions which in the direction of each axis are asymptotic to the partition attached to that axis. And this purely combinatorial entity produces the topological vertex! In the next subsection we discuss the relation of the refinement constructions with the topological vertex.

\subsection{Relation to the topological vertex}
The topological vertex originates in the work of physicists, \cite{AKMV}, and it was invented for counting holomorphic curves in toric Calabi-Yau varieties, that is, Gromov-Witten invariants. The work of Okounkov, Reshetikhin and Vafa, \cite{ORV}, put it on mathematical basis by giving it a spectacular reinterpretation as `crystal melting'. Since then the topological vertex became an important part in the theory of Gromov-Witten and Donaldson-Thomas invariants, representation theory, combinatorics. We give a minimum of notions from the subject to lead the reader to the results below; we borrow relevant definitions from
\cite{BKY}. 

In our setting we work in the closure of the octant $\RR^3_+$ with the coordinate axes labeled $(x_{+}, x_0, x_{-})$ according to the color of edges of the graph $\widehat{PG}_l$. Let $\eta=(\eta^{+},\eta^{0},\eta^{-})$
be a triple of partitions; each partition is attached to the axis of the octant having the same color. We think of partitions as their Young diagram
drawn in the plane with rows running horizontally and descending vertically.
For purposes of describing a three-dimensional partition ($3d$-partition for short) we use the coordinates $(i,j)$ to describe an ordinary $2d$-partition: the corner of the Young diagram is $(0,0)$; $j$ is the row number and $i$ is the position in the $j$-th row.

\begin{example}\label{ex:2d-partition}
The partition	$(5,3,2,1)$ is the subset of points
$$
\begin{matrix}
	(0,0)&(1,0)&(2,0)&(3,0)&(4,0)\\
	(0,1)&(1,1)&(2,1)&&\\
	(0,2)&(1,2)&&&\\
	(0,3)
\end{matrix}
$$
Drawing boxes around points above clearly gives the usual picture of the Young diagram of the partition.
\end{example}

With the above conventions we define a $3d$-partition $\pi$ {\it asymptotic} to
$(\eta^{+},\eta^{0},\eta^{-})$ as the subset of $\ZZ^3_{\geq 0}$ subject to the following conditions:

1) if one of the points $(i+1,j,k)$, or $(i,j+1,k)$, or $(i,j,k+1)$ lies in $\pi$, then $(i,j,k) \in \pi$,

2) (i) for all $i>>0$ a point $(i,j,k)\in \pi$ if and only if $(j,k)\in \eta^{+}$;

(ii)  for all $j>>0$ a point $(i,j,k)\in \pi$ if and only if $(k,i)\in \eta^{0}$;

(iii) for all $k>>0$ a point $(i,j,k)\in \pi$ if and only if $(i,j)\in \eta^{-}$.

\vspace{0.2cm}
To visualize, we think of $\pi$ as a collection of unit cubes stacked in the corner of the positive octant $\RR^3_{\geq 0}$; the rule 1) says that the pile of cubes is `stable' and pulled in the direction $(-1,-1,-1)$; the rule 2) says that the collection extends infinitely along each axis; furthermore, if one
goes far enough along $x_{+}$-axis and cut with the plane orthogonal to that axis and passing through $(i,0,0)$, then one obtains the $2d$-partition $\eta^{+}$ - the condition 2, (i); (ii) and (iii) correspond to cutting along the remaining two axes. Observe that the $2d$-partitions in the cross-section   
are placed according the cyclic rule:

- $\eta^+$ has its rows along $x_0$-direction,

- $\eta^0$ has its rows in $x_{-}$-direction,

- $\eta^-$ has its rows in $x_{+}$-direction.

With this rule in mind, one defines the {\it legs} of $\pi$:

$\bullet$ the subset
$$
\{(i,j,k)| (j,k)\in \eta^+\} \subset \pi =\text{ \it the leg of $\pi$ at $i$},
$$ 

$\bullet$ the subset
$$
\{(i,j,k)| (k,i)\in \eta^0\} \subset \pi =\text{ \it the leg of $\pi$ at $j$},
$$ 

$\bullet$
the subset
$$
\{(i,j,k)| (i,j)\in \eta^-\} \subset \pi =\text{ \it the leg of $\pi$ at $k$}.
$$ 

To a $3d$-partition $\pi$ is attached the weight function 
$$
w_{\pi}: \ZZ^3_{\geq0} \longrightarrow \ZZ
$$
defined by the rule
$$
w_{\pi}(i,j,k)=1-\text{\it the number of legs of $\pi$ passing through $(i,j,k)$.}
$$
This is used to define $|\pi|$, {\it the normalized volume of $\pi$}:
$$
|\pi|:=\sum_{(i,j,k) \in \pi} w_{\pi}(i,j,k),
$$
that is, the total weight of points in $\pi$.

\vspace{0.2cm} 
The topological vertex $V_{\eta^{+}\eta^{0}\eta^{-}}$ of $(\eta^{+},\eta^{0},\eta^{-})$ is a certain generating series counting all $3d$-partitions asymptotic to $(\eta^{+},\eta^{0},\eta^{-})$. Precisely, let ${\cal P}(\eta^{+},\eta^{0},\eta^{-})$ be the set of $3d$-partitions asymptotic to $(\eta^{+},\eta^{0},\eta^{-})$ and let $t$ be an indeterminate; one defines the generating function
$$
V_{\eta^{+}\eta^{0}\eta^{-}}= \sum_{\pi \in {\cal P}(\eta^{+},\eta^{0},\eta^{-})} t^{|\pi|}= \text{ \it the topological vertex of $(\eta^{+},\eta^{0},\eta^{-})$} .
	$$

\begin{example}\label{McMah}
	Take $\eta^{+}=\eta^{0}=\eta^{-}=\emptyset$, the empty partition. Then
	${\cal P}(\emptyset,\emptyset,\emptyset)$ is the set of all finite collections of cubes stacked in the corner of $\RR^3_{\geq 0}$ according to the rule 1) above; the normalized volume of $\pi$ is just the number of cubes in $\pi$; the topological vertex in this case is the famous McMahon function
	$$
	V_{\emptyset \emptyset \emptyset} = \sum_{\pi \in {\cal P}(\emptyset,\emptyset,\emptyset)} t^{|\pi|}=\prod^{\infty}_{n=1} (1-t^n)^{-n}.
	$$
\end{example}

From the refinement construction we have the stratification of ${\mathfrak{L}}_l$ by the strata ${\mathfrak{L}}_l(h^l,\lambda)$. They are labeled by the partitions $\lambda$: the conjugate partition $\lambda'$ of $\lambda$ counts the dimensions of the graded pieces of our filtrations.
Thus to ${\mathfrak{L}}_l$ is naturally attached the collection ${\cal P}_{{\mathfrak{L}}_l}$ of partitions $\lambda$, the labels of the strata
${\mathfrak{L}}_l(h^l,\lambda)$ occurring in ${\mathfrak{L}}_l$. Our discussion implies
\begin{equation}\label{l-TV}
	\begin{gathered}
	\text{\it to every triple $(\lambda,\mu,\nu)$ of partitions in ${\cal P}_{{\mathfrak{L}}_l}$ is attached the topological vertex}
	\\ 
		V_{\lambda\mu\nu}(t)=\sum_{\pi \in {\cal P}(\lambda,\mu,\nu)} t^{|\pi|}
		\end{gathered}
\end{equation}

The topological vertex in the context of Gromov-Witten theory is supposed to `count' curves, or in Donaldson-Thomas theory objects of some suitable categories. What do we count here? The answer is not clear at the time of writing and will require further understanding. Let us suggest one direction. The reader may recall, that the partitions $\lambda$ labeling the strata ${\mathfrak{L}}_l(h^l,\lambda)$ have a precise meaning:

-geometrically, they count cones over rational normal curves passing through the hyperplane sections of the canonical embedding of our curve,

- algebraically, irreducible ${\mathfrak{sl}}(2)$ modules occurring in $W_{\xi}/W^l_{\xi}([\phi])$, for points $([\xi],[\phi])$ in the stratum ${\mathfrak{L}}_l(h^l,\lambda)$.

So the labels for the topological vertices `count' {\it visible} objects related to the geometry of $C$. Logically, the topological vertex is supposed to count `quantum' objects. The way we arrived to the positive octant $\RR^3_{+}$, where live the $3d$-partitions of the topological vertex, was
via the quantum invariants: the points of $\RR^3_{+}$ in our setting are symmetrized metrics of the graph $\widehat{PG}_l$ and those are elliptic curves with a distinguished torsion point, see Proposition \ref{pro:plmap-ltorsion}. So one answer to `what do we count in \eqref{l-TV}' - we count elliptic curves corresponding to the metrics given by the centers of the unit cubes of the $3d$-partitions. Pursuing this direction one could refine our count by taking into consideration only cubes which are relevant to the strata ${\mathfrak{L}}_l$ and the substrata ${\mathfrak{L}}_l(h^l,\lambda)$. Precisely, fix a triple $(\lambda,\mu,\nu)$
in the set of partitions ${\cal P}_{{\mathfrak{L}}_l}$; every stratum
labeled by a partition in that set contributes the Hilbert-Laurent polynomial of the stratum, see \eqref{LH-lambda}, thus our triple $(\lambda,\mu,\nu)$ provides the Laurent polynomials
$$
h_{{\mathfrak{L}}_l(h^l,\epsilon)}, \text{for $\epsilon=\lambda,\mu,\nu$.}
$$
To simplify the notation we write $h_{\epsilon}$ instead of $h_{{\mathfrak{L}}_l(h^l,\epsilon)}$. Take a pair of such polynomials, $h_{\epsilon}, h_{\epsilon'}$, for $\epsilon,\epsilon' \in \{\lambda,\mu,\nu\}$, and use them to obtain metrics
 $$
 m_{\widehat{\rho}_c(\xi\otimes\phi) (h_{\epsilon},h_{\epsilon'})}:E_{\widehat{PG}_l} \longrightarrow \RR_{+},
$$
the value of the map
$$
metr_{\mathfrak{L}_l}: \OO^{\times}_{\mathfrak{L}_l} (-1) \times \CC \times \left(\CC[q,q^{-1}]\right)^2 \longrightarrow \RR^{E_{\widehat{PG}_l}}_{+}
	$$
in \eqref{metrl}	at the point $(\xi\otimes\phi,c,h_{\epsilon},h_{\epsilon'}) $; thus each pair $(h_{\epsilon}, h_{\epsilon'})$ produces the map
\begin{equation}\label{metrlh-h}
metr_{\mathfrak{L}_l,h_{\epsilon}, h_{\epsilon'} }: \OO^{\times}_{\mathfrak{L}_l} (-1) \times \CC \longrightarrow \RR^{E_{\widehat{PG}_l}}_{+}.
\end{equation}
Using the triple of partitions $(\lambda,\mu,\nu)$ to symmetrize the metrics
we arrive to $\RR^3_{+}$:
$$
s_{\lambda\mu\nu}:\RR^{E_{\widehat{PG}_l}}_{+} \longrightarrow \RR^3_{+};
$$
that is, we specify positive symmetric functions $s_{\lambda'}$, $s_{\mu'}$,
$s_{\nu'}$ to {\it conjugate} partitions and define the value of $s_{\lambda\mu\nu}$ for a metric $m\in \RR^{E_{\widehat{PG}_l}}_{+}$ by the rule
$$
s_{\lambda\mu\nu} (m)=(s_{\lambda'}(m^+), s_{\mu'}(m^0), s_{\nu'}(m^-)),
$$
recall: for a color $c\in\{+,0,-\}$ we denote $m^c$ the vector of values of
the metric $m$ on the edges colored by $c$:
$$
m^c=(m(e^c_0),m(e^c_1),\ldots, m(e^c_{l-1})) \in \RR^l_{+}.
$$ 
Composing the map \eqref{metrlh-h} with $s_{\lambda\mu\nu} $ gives
\begin{equation}\label{metrl-3part}
	metr^{\lambda\mu\nu}_{\mathfrak{L}_l,h_{\epsilon}, h_{\epsilon'} }: \OO^{\times}_{\mathfrak{L}_l} (-1) \times \CC \longrightarrow \RR^{3}_{+}
\end{equation}
 
 We want to modify the topological vertex $V_{\lambda\mu\nu}$ in \eqref{l-TV} to take into account the maps $metr^{\lambda\mu\nu}_{\mathfrak{L}_l,h_{\epsilon}, h_{\epsilon'} }$, that is, we will sum only over the $3d$-partitions in 
${\cal P}(\lambda,\mu,\nu)$ which `see' the images of the maps \eqref{metrl-3part} for various choices of pairs $h_{\epsilon}, h_{\epsilon'} $. The formal definition is as follows.

For a nonempty subset $A$ in $\RR^3_{\geq 0}$ and a $3d$-partition $\pi$ in
${\cal P}(\lambda,\mu,\nu)$ define the set
$$
\pi^A:=\{(i,j,k)\in \pi | \Box_{ijk} \bigcap A \neq \emptyset\},
$$
where $\Box_{ijk}$ stands for a unit box in $\RR^3_{\geq 0}$ having $(i,j,k)$ as a vertex. Consider the subset ${\cal P}^A(\lambda,\mu,\nu)$ of  those $3d$-partitions $\pi \in {\cal P}(\lambda,\mu,\nu)$ for which $\pi^A$ is nonempty:
$$
 {\cal P}^A(\lambda,\mu,\nu):=\{\pi \in {\cal P}(\lambda,\mu,\nu) | \pi^A \neq \emptyset\};
 $$
 
 define the $A$-topological vertex
\begin{equation}\label{A-TV}
	V^A_{\lambda\mu\nu}(t):=\sum_{\pi\in {\cal P}^A(\lambda,\mu,\nu) } t^{|\pi|}
\end{equation}

To a triple of partitions $\lambda,\mu,\nu$ in ${\cal P}_{\mathfrak{L}_l}$ we assign their {\it characteristic set} $A_{\lambda\mu\nu}$ as the closure
of the union of images of the maps $metr^{\lambda\mu\nu}_{\mathfrak{L}_l,h_{\epsilon}, h_{\epsilon'} }$
as $\epsilon$ and $\epsilon'$ range through $\{\lambda,\mu,\nu\}$:
\begin{equation}\label{char-3partitions}
	A_{\lambda\mu\nu}:=\overline{\bigcup_{\epsilon,\epsilon'\in \{\lambda,\mu,\nu\}} im(metr^{\lambda\mu\nu}_{\mathfrak{L}_l,h_{\epsilon}, h_{\epsilon'} })};
\end{equation}
the closure is taken with respect to the metric topology of $\RR^3$. We now attach to $\lambda,\mu,\nu$ in ${\cal P}_{\mathfrak{L}_l}$ the $A_{\lambda\mu\nu}$-topological vertex
\begin{equation}\label{A-TV-3part}
	V^{A_{\lambda\mu\nu}}_{\lambda\mu\nu}(t) =\sum_{\pi \in {\cal P}^{A_{\lambda\mu\nu}}(\lambda,\mu,\nu)}t^{|\pi|}.
\end{equation}
The following statement summarizes the discussion above.
\begin{pro}\label{pro:l-TV}
	Let ${\cal P}_{\mathfrak{L}_l}$ be the set of partitions labeling the substrata ${\mathfrak{L}}_l(h^l,\lambda)$ of $\mathfrak{L}_l$. For every triple of partitions $(\lambda,\mu,\nu)\in ({\cal P}_{\mathfrak{L}_l})^3$ is attached the characteristic set
	$$
	A_{\lambda\mu\nu} \subset \RR^3_{\geq 0}
	$$
	and $A_{\lambda\mu\nu}$-topological vertex 
	$$
	V^{A_{\lambda\mu\nu}}_{\lambda\mu\nu}(t) =\sum_{\pi \in {\cal P}^{A_{\lambda\mu\nu}}(\lambda,\mu,\nu)}t^{|\pi|}.
	$$ 
\end{pro}

\begin{rem}\label{rem:TV}
	1) The topological vertices $V^{A_{\lambda\mu\nu}}_{\lambda\mu\nu}(t)$ could be considered as a sort of characters of the maps
	$$
	\wp_{l,\lambda,c}: \OO^{\times}_{\mathfrak{L}_l} (-1) \longrightarrow Y_1(l)
	$$
	in Proposition \ref{pro:plmap-ltorsion}.
	
	2) There should be `gluing' rules between topological vertices corresponding to different triples. At this stage the author has no understanding what those could be.
	
	3) There is a circle of fascinating interrelations between $3d$-partitions, hexagonal tilings of the plane corresponding to bipartite dimers and their perfect matchings. They come from different
	sources: combinatorics, statistical mechanics, toric geometry - see \cite{ORV}, \cite{KeOkSh}, \cite{Ken},\cite{JWY}. In our case these relations are part of the constructions: the $3d$-partitions are drawn in the octant $\RR^3_{+}$, the range of symmetrized metrics of the graph $\widehat{PG}_l$, this graph is a dimer model on a torus and its boundary cycles are hexagons, thus giving hexagonal tilings of the plane, the universal covering of the torus; and the perfect matchings of $\widehat{PG}_l$ play important role in several constructions we have encountered along the way. 
\end{rem}

In the remark above we alluded to the `gluing rules' between different triples of partitions involved in the strata  ${\mathfrak{L}}_l(h^l,\lambda)$ occurring in ${\mathfrak{L}}_l$. In the next subsection we discuss one aspect of the interaction
between those strata.  This seems to be unrelated to the discussion of the topological vertex, but hopefully it will illustrate the interest to explore those interactions.

\subsection{The incidence strata $\widetilde{\mathfrak{L}}_l(\lambda,\mu)$ and Weil pairing} 
We will be working under assumption $l\geq 4$. For every stratum ${\mathfrak{L}}_l(h^l,\lambda)$ occurring in ${\mathfrak{L}}_l$ we have the lifting
$$
\widetilde{\wp}_{l,\lambda,c}:\OO^{\times}_{\mathfrak{L}_l} (-1) \longrightarrow \mathfrak{H}
$$
assigning to a point $\xi\otimes\phi$ in $\OO^{\times}_{\mathfrak{L}_l} (-1)$ the period
 of the elliptic curve $\Gamma_{{m(c, h_{\mathfrak L(h^l,\lambda)},\xi\otimes\phi)}_{s_{\lambda'}} }$ in Proposition \ref{pro:plmap-ltorsion}. Equivalently, the map can be recast as a continuous family of elliptic curves
 $$
 \xymatrix{
 {\mathfrak G}_{l,\lambda,c} \ar[d]_{\pi_{l,\lambda,c}}&\\
 \OO^{\times}_{\mathfrak{L}_l} (-1) \ar[r]^(.6){\widetilde{\wp}_{l,\lambda,c}}&{\mathfrak H}
}
$$
Furthermore, the projection $\pi_{l,\lambda,c}$ comes with a distinguished section 
$$
\varkappa^{\lambda}_1: \OO^{\times}_{\mathfrak{L}_l} (-1) \longrightarrow {\mathfrak G}_{l,\lambda,c}
$$
which maps a point $\xi\otimes\phi$ in $\OO^{\times}_{\mathfrak{L}_l} (-1)$ to the torsion point $\varkappa^{\lambda}_1(\xi\otimes\phi)\in \Gamma_{{ m(c, h_{\mathfrak L(h^l,\lambda)},\xi\otimes\phi)}_{s_{\lambda'}} }$ as asserted in Proposition \ref{pro:plmap-ltorsion}. We will also denote ${\mathfrak G}_{l,\lambda,c} (l)$ the relative subgroup of torsion points of order $l$. In particular, the section $\varkappa^{\lambda}_1$ takes its values in ${\mathfrak G}_{l,\lambda,c} (l)$:
$$
\varkappa^{\lambda}_1:  \OO^{\times}_{\mathfrak{L}_l} (-1)\longrightarrow {\mathfrak G}_{l,\lambda,c} (l) \subset {\mathfrak G}_{l,\lambda,c}.
$$

Let ${\mathfrak{L}}_l(h^l,\lambda)$ and ${\mathfrak{L}}_l('h^l,\mu)$ be two, not necessarily distinct, substrata of ${\mathfrak{L}}_l$. Then we form the fibre product
$$
\xymatrix{
\OO^{\times}_{\mathfrak{L}_l} (-1) \times_{{\mathfrak{H}}} \OO^{\times}_{\mathfrak{L}_l} (-1) \ar[r]\ar[d]&\OO^{\times}_{\mathfrak{L}_l} (-1) \ar[d]^{\widetilde{\wp}_{l,\lambda,c}}\\
\OO^{\times}_{\mathfrak{L}_l} (-1) \ar[r]^{\widetilde{\wp}_{l,\mu,c}}& {\mathfrak{H}} 
}
$$
This fibre product will be denoted ${\widetilde{\mathfrak{L}}_l}(\lambda,\mu)$:
\begin{equation}\label{lambda-mu-stratum}
\widetilde{\mathfrak{L}}_l(\lambda,\mu)=\left\{(\xi\otimes\phi, \xi'\otimes\phi')\in \OO^{\times}_{\mathfrak{L}_l} (-1) \times \OO^{\times}_{\mathfrak{L}_l} (-1) \big|\,\, \widetilde{\wp}_{l,\mu,c}(\xi\otimes\phi)=\widetilde{\wp}_{l,\lambda,c}(\xi'\otimes\phi')\right\}.
\end{equation}
Over the fibre product $\widetilde{\mathfrak{L}}_l(\lambda,\mu)$ we have the fibration  
$$
\xymatrix{
	{\mathfrak G}_{l,\lambda,c} \times_{\mathfrak{H}} {\mathfrak G}_{l,\mu,c} \ar[d]_{\pi_{l,\lambda,\mu,c}}&\\
\widetilde{\mathfrak{L}}_l(\lambda,\mu) \ar[r]^(.6){\widetilde{\wp}_{l,\lambda,\mu,c}}&{\mathfrak{H}}
}
$$
in abelian surfaces. Precisely, over a point $\tau \in {\mathfrak{H}}$ in the image of $\widetilde{\wp}_{l,\lambda,\mu,c}$ lies the Cartesian product
$$
\widetilde{\mathfrak{L}}^{\lambda}_l(\tau) \times \widetilde{\mathfrak{L}}^{\mu}_l(\tau)={\widetilde{\wp}}^{-1}_{l,\lambda,\mu,c}(\tau),
$$
where $\widetilde{\mathfrak{L}}^{\lambda}_l(\tau)$ and $\widetilde{\mathfrak{L}}^{\mu}_l(\tau)$ are respectively the fibre of ${\widetilde{\wp}}_{l,\lambda,c}$ and ${\widetilde{\wp}}_{l,\mu,c}$ over $\tau$;
over the Cartesian product $\widetilde{\mathfrak{L}}^{\lambda}_l(\tau) \times \widetilde{\mathfrak{L}}^{\mu}_l(\tau)$ the family $\pi_{l,\lambda,\mu,c}$ is isotrivial with the fibre $\Gamma_{\tau} \times \Gamma_{\tau}$, where $\Gamma_{\tau}=\CC/(\ZZ+\ZZ\tau)$. 
 
 In addition, the abelian surface $\Gamma_{\tau} \times \Gamma_{\tau}$ comes with distinguished points of order $l$:
$$
\varkappa^{\lambda,\mu}(\tau):=\varkappa^{\lambda}_1|_{\widetilde{\mathfrak{L}}^{\lambda}_l(\tau) } \times \varkappa^{\mu}_1|_{\widetilde{\mathfrak{L}}^{\mu}_l(\tau) } :\widetilde{\mathfrak{L}}^{\lambda}_l(\tau) \times \widetilde{\mathfrak{L}}^{\mu}_l(\tau) \longrightarrow \Gamma_{\tau}(l)\times \Gamma_{\tau}(l) \subset \Gamma_{\tau} \times \Gamma_{\tau}
$$
provided by sections $\varkappa^{\lambda}_1$ and  $\varkappa^{\mu}_1$. Since the map is continuous, it is constant on every connected component of 
$\widetilde{\mathfrak{L}}^{\lambda}_l(\tau) \times \widetilde{\mathfrak{L}}^{\mu}_l(\tau)$. Thus we obtain the map
\begin{equation}\label{lambda-mu-torsion}
\varkappa^{\lambda,\mu}(\tau)	:\pi_0 \big(\widetilde{\mathfrak{L}}^{\lambda}_l(\tau) \times \widetilde{\mathfrak{L}}^{\mu}_l(\tau) \big)=\pi_0(\widetilde{\mathfrak{L}}^{\lambda}_l(\tau)) \times \pi_0(\widetilde{\mathfrak{L}}^{\mu}_l(\tau) ) \longrightarrow \Gamma_{\tau}(l)\times \Gamma_{\tau}(l).
\end{equation}
Recall that the group of $l$-torsion points $\Gamma_{\tau}(l)$ is equipped
with the Weil pairing
\begin{equation}\label{Weil-pair}
{\mathfrak{e}}_l:\Gamma_{\tau}(l) \times \Gamma_{\tau}(l) \longrightarrow {\mathbb{G}_l},
\end{equation}
where ${\mathbb{G}_l}$ is the group of $l$-th roots of unity in the algebraic closure of the field of definition of the elliptic curve and which assigns to each pair $(p,p')$ of $l$-torsion points of the curve an $l$-th root of unity ${\mathfrak{e}}_l(p,p')$, see \cite{Sh}. Thus we obtain the following.
\begin{pro}\label{pro:Weilpair}
	Let ${\mathfrak{L}}_l(h^l,\lambda)$ and ${\mathfrak{L}}_l('h^l,\mu)$ be two substrata of ${\mathfrak{L}}_l$. Then we have the incidence correspondence
	$$
	\widetilde{\mathfrak{L}}_l(\lambda,\mu) \subset \OO^{\times}_{{\mathfrak{L}}_l}(-1)\times \OO^{\times}_{{\mathfrak{L}}_l}(-1)
	$$
	defined in \eqref{lambda-mu-stratum} together with the continuous period map
	$$
\widetilde{\wp}_{l,\lambda,\mu,c}:	\widetilde{\mathfrak{L}}_l(\lambda,\mu) \longrightarrow {\mathfrak{H}}.
$$
For every value $\tau$ in the range of $\widetilde{\wp}_{l,\lambda,\mu,c}$
denote $\widetilde{\mathfrak{L}}^{\lambda}_l(\tau)$ and $\widetilde{\mathfrak{L}}^{\mu}_l(\tau)$
the fibre of $\widetilde{\wp}_{l,\lambda,c}$ and $\widetilde{\wp}_{l,\mu,c}$
respectively.
Then we obtain the map
$$
\varkappa^{\lambda,\mu}(\tau):\pi_0 \big(\widetilde{\mathfrak{L}}^{\lambda}_l(\tau) \times \widetilde{\mathfrak{L}}^{\mu}_l(\tau) \big)=\pi_0(\widetilde{\mathfrak{L}}^{\lambda}_l(\tau)) \times \pi_0(\widetilde{\mathfrak{L}}^{\mu}_l(\tau) ) \longrightarrow \Gamma_{\tau}(l)\times \Gamma_{\tau}(l)
$$
which assigns to every connected component of $\widetilde{\mathfrak{L}}^{\lambda}_l(\tau) \times \widetilde{\mathfrak{L}}^{\mu}_l(\tau)$ a pair of $l$-torsion points of $\Gamma_{\tau}$. The composition of this map with the Weil pairing \eqref{Weil-pair} gives the map
$$
{\mathfrak{e}}^{\lambda,\mu}_l (\tau):={\mathfrak{e}}_l \circ \varkappa^{\lambda,\mu}(\tau): \pi_0 \big(\widetilde{\mathfrak{L}}^{\lambda}_l(\tau) \times \widetilde{\mathfrak{L}}^{\mu}_l(\tau) \big)=\pi_0(\widetilde{\mathfrak{L}}^{\lambda}_l(\tau)) \times \pi_0(\widetilde{\mathfrak{L}}^{\mu}_l(\tau) ) \longrightarrow {\mathbb G}_l
$$
The values of this map are either $1$ or a primitive $l$-th root of unity:
let $(p,p')$ be in the image of $\varkappa^{\lambda,\mu}(\tau)$, then

- the first value occurs when the pair $(p,p')$  generates the same subgroup in $\Gamma_{\tau}(l)$, that is,
$\langle p\rangle =\langle p'\rangle$, where $\langle t\rangle$ denotes the subgroup of $\Gamma_{\tau}$ generated by a point $t$; 

- the second value occurs when $p$ and $p'$ generate $\Gamma_{\tau}(l)$:
$$
\langle p\rangle \times\langle p'\rangle =\Gamma_{\tau}(l).
$$
\end{pro}
\begin{pf}
	Only the last assertion about values of the Weil pairing needs a proof. Recall that upon a choice of two generators, say $x$ and $y$, we have an isomorphism
	$$
	\Gamma_{\tau}(l) \cong \ZZ/l\ZZ \times \ZZ/l\ZZ. 
	$$
	and the {\it noncanonical} Weil pairing can be computed by sending a pair $(p,p')\in \Gamma_{\tau}(l)\times \Gamma_{\tau}(l)$ to the determinant
	$$
	det(p,p'):=\begin{vmatrix}
		a&a'\\
		b&b'
	\end{vmatrix} \in \ZZ/l\ZZ,
$$
where $p=ax+by$ and $p'=a'x+b'y$ are expressed in terms of the generators $x$ and $y$. The value ${\mathfrak e}_l(p,p') $ is $\zeta^{det(p,p')}$, where $\zeta$ is a primitive $l$-th root of unity.

We now take $(p,p')$ in the image of $\varkappa^{\lambda,\mu}(\tau)$. We know that both $p$ and $p'$ are points of order precisely $l$ in $\Gamma_{\tau}(l)$. Hence they generate the subgroups $\langle p \rangle$
and $\langle p' \rangle$ of order $l$ in  $\Gamma_{\tau}(l)$ and we have the commutative diagram

$$
\xymatrix{
	&&0\ar[d]&&\\
	&&\langle p' \rangle \ar[d]\ar[dr]&&\\
0\ar[r]&\langle p \rangle \ar[r]\ar[dr]&\Gamma_{\tau}(l) \ar[r]\ar[d]&\Gamma_{\tau}(l)/\langle p \rangle\ar[r] &0\\
&& \Gamma_{\tau}(l)/\langle p' \rangle\ar[d]&&\\
&&0&&
}
$$
The two possibilities arise

$\bullet$ the slated arrow is injective and hence an isomorphism,

$\bullet$ the slanted arrow has a nontrivial kernel which is the intersection of subgroups
$\langle p \rangle $ and $\langle p' \rangle$.

The first item means that $p$ and $p'$ are generators of  $\Gamma_{\tau}(l)$.
Hence $det(p,p')$ is invertible in $\ZZ/l\ZZ$ and $\zeta^{det(p,p')}$ is a primitive $l$-th root of unity.

The second item and the fact that $p$ (resp., $p'$) is a generator of $\langle p \rangle$ (resp., $\langle p' \rangle$) imply the equality
$$
 \langle p \rangle =\langle p' \rangle.
 $$
 Hence $p'=sp$, for some (invertible) $s\in \ZZ/l\ZZ$. The determinant
 $$
 det(p,p')=0
 $$
 with the conclusion ${\mathfrak e}_l(p,p')=1 $.
\end{pf}

The outcome of the above is that the period maps taking into account the strata 
${\mathfrak{L}}_l(h^l,\lambda)$ not only see the elliptic curves but also detect subtler invariants such as torsion points of order $l$ on those curves: for every value $\tau$ of the period map, 
  the set of connected components
$$
\pi_0 \big(\widetilde{\mathfrak{L}}^{\lambda}_l(\tau) \times \widetilde{\mathfrak{L}}^{\mu}_l(\tau)\big)=\pi_0 \big(\widetilde{\mathfrak{L}}^{\lambda}_l(\tau)\big) \times \pi_0\big(\widetilde{\mathfrak{L}}^{\mu}_l(\tau)\big)
$$
comes equipped with the map
$$
{\mathfrak{e}}^{\lambda,\mu}_l (\tau):\pi_0 \big(\widetilde{\mathfrak{L}}^{\lambda}_l(\tau)\big) \times \pi_0\big(\widetilde{\mathfrak{L}}^{\mu}_l(\tau)\big) \longrightarrow {\mathbb{G}}_l.
$$
Either the map is constant of value $1$ and then all the torsion points
labeling the connected components of $\widetilde{\mathfrak{L}}^{\lambda}_l(\tau)$ and $\widetilde{\mathfrak{L}}^{\mu}_l(\tau)$ define the {\it same} subgroup of order $l$ in $\Gamma_{\tau}(l)$, or, the map is non-constant and every pair $([x],[y])$ of connected components for which ${\mathfrak{e}}^{\lambda,\mu}_l (\tau)$ is not $1$ give

- ${\mathfrak{e}}^{\lambda,\mu}_l (\tau)([x],[y])$ a primitive $l$-th root of unity,

- a full $l$-level structure on $\Gamma_{\tau}$, that is, a basis
of the group $\Gamma_{\tau}(l)$ given by the pair of torsion points
$$
\varkappa^{\lambda,\mu}([x],[y]) \in \Gamma_{\tau}(l) \times \Gamma_{\tau}(l);
$$
in other words $(\Gamma_{\tau},\varkappa^{\lambda,\mu}([x],[y]) )$ gives a point of the modular curve $Y(l)=\Gamma(l)\backslash {\mathfrak H}$, where $\Gamma(l)$ is the full congruence subgroup of level $l$:
$$
\Gamma(l)=ker\big(SL_2(\ZZ) \longrightarrow SL_2(\ZZ/l\ZZ)\big).
$$
We summarize the above in the following statement.
\begin{pro}\label{pro:Llambdamu-Yl}
	Let ${\mathfrak{L}}_l(h^l,\lambda)$ and ${\mathfrak{L}}_l('h^l,\mu)$ be two substrata of ${\mathfrak{L}}_l$. Then we have the incidence
	$$
	\widetilde{\mathfrak{L}}_l(\lambda,\mu) \subset \OO^{\times}_{{\mathfrak{L}}_l}\times \OO^{\times}_{{\mathfrak{L}}_l}
	$$
	together with the continuous period map
	$$
	\widetilde{\wp}_{l,\lambda,\mu,c}:	\widetilde{\mathfrak{L}}_l(\lambda,\mu) \longrightarrow {\mathfrak{H}}.
	$$
	Furthermore, $\widetilde{\mathfrak{L}}_l(\lambda,\mu)$ admits the map
	$$
{\mathfrak{e}}^{\lambda,\mu}_l:	\widetilde{\mathfrak{L}}_l(\lambda,\mu) \longrightarrow {\mathbb{G}}_l
	$$
	into the group ${\mathbb{G}}_l$ of $l$-th roots of unity. Set
	$$
	\widetilde{\mathfrak{L}}^{1}_l(\lambda,\mu) :=\left({\mathfrak{e}}^{\lambda,\mu}_l\right)^{-1}(1)
	$$
	the preimage of $1$ and consider the complement
	$$
			\widetilde{\mathfrak{L}}_l(\lambda,\mu)\setminus\widetilde{\mathfrak{L}}^{1}_l(\lambda,\mu).
			$$
			If the latter is nonempty it is the disjoint union
			$$
			\widetilde{\mathfrak{L}}_l(\lambda,\mu)\setminus\widetilde{\mathfrak{L}}^{1}_l(\lambda,\mu)=\coprod_{\varepsilon \in im'({\mathfrak{e}}^{\lambda,\mu}_l)}
			\widetilde{\mathfrak{L}}^{\varepsilon}_l(\lambda,\mu),
			$$
			where $im'({\mathfrak{e}}^{\lambda,\mu}_l)=im({\mathfrak{e}}^{\lambda,\mu}_l)\setminus \{1\}$ consists of primitive $l$-th roots of unity and where $\widetilde{\mathfrak{L}}^{\varepsilon}_l(\lambda,\mu)$ is
			the preimage $\left({\mathfrak{e}}^{\lambda,\mu}_l\right)^{-1}(\varepsilon)$. Furthermore, for every $\varepsilon$ occurring in the above union we have the continuous map
			$$
				{\wp}^{\varepsilon}_{l,\lambda,\mu,c}:	\widetilde{\mathfrak{L}}^{\varepsilon}_l(\lambda,\mu) \longrightarrow Y(l)=\Gamma(l)\backslash{\mathfrak{H}}
			$$
			fitting into the following commutative diagram
			$$
			\xymatrix{
				&{\mathfrak{H}}\ar[dd]\\
				\widetilde{\mathfrak{L}}^{\varepsilon}_l(\lambda,\mu)\ar[ur]^{\widetilde{\wp}_{l,\lambda,\mu,c}}\ar[dr]_{{\wp}^{\varepsilon}_{l,\lambda,\mu,c}}&\\
				&Y(l)
			}
		$$
		where the vertical map is the natural projection.
\end{pro}

\subsection{Relating $C$ to the modular curve $Y_1(l)$}
The appearance of modular curves in our constructions reinforces the connections with number theory, see a brief indication about such a connection in \S11.5. Recall the continuous maps 
$$
\wp^{\phi}_{l }(\xi\otimes\psi)(c,\bullet): C\setminus D(\phi,[\xi],[\psi]) \longrightarrow {\mathfrak{M}}^l_1
 $$
 attached to the points $\xi\otimes\psi$ in $\OO^{\times}_{{\mathfrak{L}}_l }(-1)$ lying over the stratum ${\mathfrak{L}}_l (\phi)$, see Proposition \ref{pro:pl-C-moduli} and the subsequent discussion. By construction the map factors through the space of metrics $\RR^{E_{\widehat{PG}_l}}_{+}$ on $\widehat{PG}_l$
 $$
 C\setminus D(\phi,[\xi],[\psi]) \longrightarrow \RR^{E_{\widehat{PG}_l}}_{+} \longrightarrow {\mathfrak{M}}^l_1.
 $$
 Assume $l\geq 4$, then we have seen that symmetrizing metrics allows us to obtain maps  to the modular curve $Y_1(l)$ 
 $$
  \OO^{\times}_{{\mathfrak{L}}_l}(-1) \longrightarrow Y_1 (l).
 $$
 Applying this to metrics defining the values of $\wp^{\phi}_{l }(\xi\otimes\psi)(c,\bullet)$ will give its symmetrized version
 $$
 {}^S\wp^{\phi}_{l,\lambda }(\xi\otimes\psi)(c,\bullet): C\setminus D(\phi,[\xi],[\psi]) \longrightarrow Y_1(l),
 $$
 for every $\lambda$, the partition occurring in a substrata ${\mathfrak{L}}_l(h^l,\lambda)$ of ${\mathfrak{L}}_l$; the superscript $S$ on the left stands for `symmetrization'.
This depends on a symmetrization pattern chosen, that is, we agree once and for all to fix a basis
 $$
 Sym=\{s_{\eta}\}
 $$
 of positive symmetric functions indexed by partitions. 
 \begin{pro}\label{pro:pl-C-modcurve}
 	Let $\phi \in \HKC$ be admissible and let ${\mathfrak L}_l (\phi)$ be its domain. Then to every point $([\xi],[\psi])$ in ${\mathfrak L}_l (\phi)$ is attached: 
 	
 	- the `polar' divisor $D(\phi,[\xi],[\psi]) \in |lK_C|$,
 	
 	- the continuous map
 	$$
 	\wp^{\phi}_l (\xi\otimes\psi): \CC \times \left(C\setminus D(\phi,[\xi],[\psi])\right) \longrightarrow {\mathfrak M}^l_1,
 	$$
 	for every nonzero $\xi\otimes\psi$ in the fibre of the total space
 	$\OO_{{\mathfrak{L}}_l}(-1)$ lying over the point $([\xi],[\psi])$ in ${\mathfrak L}_l (\phi)$.
 	In addition, with a symmetrization pattern chosen, each substratum 
 	${\mathfrak L}_l (h^l,\lambda)$ occurring in ${\mathfrak L}_l$ determines the partition $\lambda$. Its conjugate partition $\lambda'$ gives the positive symmetric function $s_{\lambda'}$. This gives rise to
 	the continuous map
 	$$
 	 {}^S\wp^{\phi}_{l,\lambda}(\xi\otimes\psi)(c,\bullet): C\setminus D(\phi,[\xi],[\psi]) \longrightarrow Y_1(l),
 	 $$
 	 for every nonzero $\xi\otimes\psi$ in $\OO_{{\mathfrak{L}}_l (\phi)}(-1)$ lying over a point $([\xi],[\psi]$ in ${\mathfrak{L}}_l (\phi)$.
 \end{pro} 
 \begin{pf}
 	The first part of the statement is Proposition \ref{pro:pl-C-moduli}.
 	To obtain the second part, the map ${}^S\wp^{\phi}_{l,\lambda}(\xi\otimes\psi)(c,\bullet)$, we symmetrize the metrics defined by the map
 	$$
 	 \wp^{\phi}_l (\xi\otimes\psi)(c,\bullet):C\setminus D(\phi,[\xi],[\psi]) \longrightarrow {\mathfrak M}^l_1
 	 $$
 	 using the symmetric function $s_{\lambda'}$ associated to the partition
 	 $\lambda'$ conjugate to $\lambda$ occurring as a label in the substratum ${\mathfrak{L}}_l(h^l,\lambda)$. This gives the map
 	 $$
 	{}^S\wp^{\phi}_{l,\lambda} (\xi\otimes\psi)(c,\bullet):C\setminus D(\phi,[\xi],[\psi]) \longrightarrow Y_1(l). 
 	$$
 	as stated in the proposition.
 \end{pf}

 In \S11.5 we discussed how the data of a number field can be used to construct representations of the quiver $\widehat{PG}_l$. We have learned that each such representation defines a metric on the ribbon graph $\widehat{PG}_l$ and hence points on the moduli space ${\mathfrak{M}}^l_1$ or on the modular curve $Y_1 (l)$, after appropriately symmetrizing the metric. Thus we have the values  in ${\mathfrak{M}}^l_1$ (resp. $Y_1 (l)$) associated to number fields via our refinement constructions. We formulate this as the following principle:
 \begin{equation}\label{SK-princile}
 	\begin{gathered}
 		\text{\it associate with a number field $K$ a finite collection of points
 			$S_K \subset \CC^{\times}\times \CC^{\times}$,}
 		\\
 		 \text{\it then obtain the continuous map}
 		\\
 		\wp^{S_K}_{l,c}: \OO^{\times}_{{\mathfrak{L}}_l} (-1)\longrightarrow ({\mathfrak{M}}^l_1)^{S_K}
 		\\
 		\text{\it as well as its symmetrized versions}
 		\\
 		{}^S\wp^{S_K}_{l,c,\lambda}: \OO^{\times}_{{\mathfrak{L}}_l} (-1)\longrightarrow (Y_1(l))^{S_K}.
 	\end{gathered}
 \end{equation}
The map $\wp^{S_K}_{l,c}$ is obtained from the map
$$
\wp_{l}: \CC \times\OO^{\times}_{{\mathfrak{L}}_l}(-1)\times(\CC[q,q^{-1}]^{\vee})^2 \longrightarrow {\mathfrak{M}}^l_1
$$
in Proposition \ref{pro:plmap} by fixing a value $c\in \CC$ of the trace parameter and where the functionals on $\CC[q,q^{-1}]$ are specified by the inclusions
$$
 S_K\subset (\CC^{\times} )^2 \subset (\CC[q,q^{-1}]^{\vee})^2.
 $$
In other words, to every point $\xi\otimes\phi$ in $\OO^{\times}_{{\mathfrak{L}}_l}(-1)$ and every $(u,v)\in S_K$  we associate the metric
$m(c,(u,v),\xi\otimes\phi)$ on the ribbon graph $\widehat{PG}_l$ and hence  the curve
$\Gamma_{m(c,(u,v),\xi\otimes\phi)}$; the collection of curves
$$
\{\Gamma_{m(c,(u,v),\xi\otimes\phi)}|\,(u,v)\in S_K\}
$$
determines the value
  $\wp^{S_K}_{l,c}(\xi\otimes\phi)$. 
  
  The symmetrized version of the above is obtained by fixing a symmetrization pattern, that is, a basis $S=\{s_{\eta}\}$ of positive symmetric functions labeled by partitions $\eta$. Then for every partition $\eta$ with $l$ parts the corresponding symmetric function $s_{\eta}$ is used to symmetrize the metric $m(c,(u,v),\xi\otimes\phi)$,
  see Lemma \ref{lem:symm}. The corresponding metric denoted
  $m(c,(u,v),\xi\otimes\phi)_{s_{\eta}}$ produces the elliptic curve
  $\Gamma_{m(c,(u,v),\xi\otimes\phi)_{s_{\eta}}}$ with a distinguished torsion point $p_{m(c,(u,v),\xi\otimes\phi)_{s_{\eta}}}$ of order $l$, see Proposition \ref{pro:torsion}. The collection 
  $$
  \{(\Gamma_{m(c,(u,v),\xi\otimes\phi)_{s_{\eta}}},p_{m(c,(u,v),\xi\otimes\phi)_{s_{\eta}}}) | \,(u,v)\in S_K\}
  $$ 
  determines the value ${}^S \wp_{l,c,\eta}(\xi\otimes\phi)$.
  
  The points of ${\mathfrak{M}}^l_1$ (resp., $Y_1(l)$) associated to $S_K$ via the refined IVHS constructions could be considered as {\it IVHS $K$-points of 
${\mathfrak{M}}^l_1$ (resp., $Y_1 (l)$)}:
\begin{equation}\label{IVHS-Kpoints}
	\begin{gathered}
		\OO^{\times}_{{\mathfrak{L}}_l}(-1) \ni \xi\otimes\phi \mapsto \{\wp^{S_K}_{l,c}(\xi\otimes\phi )\}= \text{\it IVHS $K$-points of ${\mathfrak{M}}^l_1$,}	
		\\
	\OO^{\times}_{{\mathfrak{L}}_l}(-1) \ni\xi\otimes\phi \mapsto \{{}^S\wp^{S_K}_{l,c,\lambda}(\xi\otimes\phi )\}= \text{\it IVHS $K$-points of $Y_1 (l)$.}
	\end{gathered}
\end{equation}
Furthermore, using the maps ${}^S\wp^{\phi}_{l,\lambda}(\xi\otimes\phi)(c,\bullet)$ in Proposition \ref{pro:pl-C-modcurve} we can pull those IVHS $K$-points back to $C$. Though our curve $C$ may have no relation to a number field, we obtain the loci related to those number fields through refined IVHS constructions.
 \begin{pro}\label{pro:IVHS-K-points}
 	Let $K$ be a number field and let $S_K$ be a finite subset of $\CC^{\times}\times \CC^{\times}$ associated to $K$ via some natural arithmetic constructions, see for example the discussion in \S11.5.
 	Then the map ${}^S\wp^{S_K}_{l,c,\lambda}$ in \eqref{SK-princile} assigns
 	to every point $(\xi'\otimes\phi' )$ of $\OO^{\times}_{{\mathfrak{L}}_l}(-1)$ the set
 	$\{{}^S\wp^{S_K}_{l,c,\lambda}(\xi'\otimes\phi')\}$ of IVHS $K$-points of $Y_1(l)$. This set will be denoted
 	$$
 	S^{l,c,\lambda}_K (\xi'\otimes\phi'):=\{{}^S\wp^{S_K}_{l,c,\lambda}(\xi'\otimes\phi')\}.
 	$$
 	 Furthermore, for a point $([\xi],[\psi])$ in ${\mathfrak{L}}_l(\phi)$
 	 and $\xi\otimes\psi$ in $\OO^{\times}_{{\mathfrak{L}}_l(\phi)}(-1)$ lying over it we have the continuous map
 	$$
 	{}^S\wp^{\phi}_{l,\lambda}(\xi\otimes\psi)(c,\bullet): C\setminus D(\phi,[\xi],[\psi]) \longrightarrow Y_1(l)
 	$$
 	in Proposition \ref{pro:pl-C-modcurve}; the closure of the preimage
 	of  $S^{l,\lambda,c}_K (\xi'\otimes\phi')$ under the map ${}^S\wp^{\phi}_{l,\lambda,c} (\xi\otimes\psi)(c,\bullet)$ is a compact subset of $C$.
 	This subset will be denoted
 	$$
 	(\phi,\xi\otimes\psi)S^{l,\lambda,c}_K (\xi'\otimes\phi'):=\overline{\left({}^S\wp^{\phi}_{l,\lambda}(\xi\otimes\psi)(c,\bullet)\right)^{-1}(S^{l,c,\lambda}_K (\xi'\otimes\phi'))}
 		$$
 		and we call it IVHS $K$-locus of $C$ associated to the map ${}^S\wp^{\phi}_{l,\lambda}(\xi\otimes\psi)(c,\bullet)$ and lying over $S^{l,\lambda,c}_K (\xi'\otimes\phi')$.
 	 \end{pro}
 
 Thus we now have `arithmetic' loci of $C$ attached to a number field $K$ and the data coming from IVHS - the IVHS-induced $K$-points of $C$.
 Of course $Y_1 (l)$ is defined over $\mathbb{Q}$, so given a number field $K$, we can consider $Y_1 (l)(K)$,
 the $K$-valued points of $Y_1(l)$, that is, the isomorphism classes of pairs $(E(K), p(K)) $ of elliptic curves $E$ with an $l$-torsion point $p$ on it, both defined over $K$. Taking the preimage of those points with respect to ${}^S\wp^{\phi}_{l,\lambda}(\xi\otimes\psi)(c,\bullet)$ gives apriori another collection of loci on $C$ attached to the number field $K$. It seems interesting to understand each type of loci and to compare them.
  
\section{Representations $\widehat{\rho}_c$ and Local systems}
We have seen how representations  $\widehat{\rho}_c$ of the quiver $\widehat{PG}_l$ give rise to metrics on the ribbon graph $\widehat{PG}_l$ thus turning the topological torus $\mathbb{T}$ into elliptic curves with $l$ marked points. In this section we use quantum-type invariants to define
local systems on those elliptic curves.

We assume to be on the stratum ${\mathfrak L}_l$, for $l\geq 3$. Let
$([\xi],[\phi])$ be a point in  ${\mathfrak L}_l$ and let $\xi\otimes\phi$
be a point of $\OO^{\times}_{{\mathfrak{L}}_l}(-1)$ lying over it.
Then we have the representation 
$$
\rho_c(\xi\otimes\phi)=\{\alpha^{t,s}_c(\xi,\phi):P^s ([\xi],[\phi])\longrightarrow P^t([\xi],[\phi])\}
$$ 
 of the quiver $PG_l$ associated to the $([\xi],[\phi])$-filtration of $W_{\xi}$. We have constructed the maps
$$
\tau^{\pm}_{\rho_c(\xi\otimes\phi)}: Hom_{\CC}(\CC[q,q^{-1}], \CC) \longrightarrow P^{\pm}([\xi],[\phi]),
$$
where $P^{\pm}([\xi],[\phi])$ are the spaces of linear maps between the end-spaces of the representation $\rho_c(\xi\otimes\phi)$:
$$
P^-([\xi],[\phi])=Hom(P^0([\xi],[\phi]),P^{l-1}([\xi],[\phi])), \,\, P^+([\xi],[\phi])=Hom(P^{l-1}([\xi],[\phi]),P^0([\xi],[\phi])),
$$
see Proposition \ref{pro:mu-map} and the details of the construction leading to \eqref{dualC+-maps-q-1}. Given two linear functionals $F$ and $G$ on the space of Laurent polynomials
$\CC[q,q^{-1}]$ we obtain the extension $\widehat{\rho}_c (F,G)(\xi\otimes\phi)$ of the representation $\rho_c(\xi\otimes\phi)$ of $PG_l$ to the quiver $\widehat{PG}_l$, where the additional
edges $e^-_0$ and $e^+_{l-1}$ are labeled by linear maps
$\tau^-_{\rho_c(\xi\otimes\phi)}(F)$ and $\tau^+_{\rho_c(\xi\otimes\phi)}(G)$ respectively. We learned in the previous sections that this gives us the metric
$$
m_{c, (F,G),\xi\otimes\phi}: {E_{\widehat{PG}_l}}\longrightarrow \RR_{+}
$$
and the associated elliptic curve $\Gamma_{c, (F,G),\xi\otimes\phi }$; to simplify the notation the curve will be denoted $\Gamma$. 
We now describe a procedure which will turn the vector spaces $\{P^s([\xi],[\phi])\}$ of the representation
$\widehat{\rho}_c (F,G)(\xi\otimes\phi)$ into local systems on the curve $\Gamma$; in the constructions below a point
$([\xi],[\phi])$ will be fixed, so the reference to $\xi$ and $\phi$ will be omitted, if no ambiguity is likely. 

Recall: to have a local system of rank $n$ on $\Gamma$ is the same as to have a representation
$$
\mu: \pi_1 (\Gamma, p) \longrightarrow GL(V),
$$
where $V$ is a complex vector space of dimension $n$ and $ \pi_1 (\Gamma,p)$ is the fundamental group of the curve  $\Gamma$ based at a point $p$ of the curve. Since our curve is of genus one, its fundamental group is abelian and hence independent of a base point; more precisely, we have
$$
\pi_1 (\Gamma)\cong \ZZ^2.
$$
In particular, to a give a representation of $\pi(\Gamma))$ on a vector space $V$, it is enough to assign to two generators of the fundamental group two commuting invertible
endomorphisms of $V$. As we have noticed before the zig-zag paths $Z_{0,\pm}$ provide a distinguished basis of  $\pi_1 (\Gamma)$:
$$
\gamma_{-}:=[Z_{0,-}],\,\, \gamma_{+}:=\frac{1}{l}[Z_{0,+}-Z_{0,-}].
$$
Next recall that given a representation $\rho_c(\xi\otimes\phi)=\{\alpha^{t,s}:P^s \longrightarrow P^t\}$ of the quiver $PG_l$ , the zig-zag paths $Z_{0,\pm}$
define the linear maps
$$
C^{\pm}_{\rho_c(\xi\otimes\phi)}:P^{\pm} \longrightarrow {\mathfrak g}^{(0)}=\bigoplus ^{l-1}_{s=0} End(P^s),
$$
see \eqref{C+-maps}. Thus given $X^{+}$ in $P^+$ and $X^-$ in $P^-$, the zig-zag paths give us
endomorphisms
$$
(C^+)^s(X^+):P^s \longrightarrow P^s, \,\,\,  (C^-)^s(X^-):P^s \longrightarrow P^s,
$$
for every $s\in [0,l-1]$. In particular, for the curve  $\Gamma$ we have natural choice for $X^{\pm}$: $$
\text{$X^+=\tau^+(G)$ and $X^{-}=\tau^-(F)$.}
$$
 With this choice made, we proceed with the construction  of the representations of the fundamental group of $\Gamma$ for each vector space $P^s$. Namely, to the path $\gamma_{-}$ we assign the exponential of the endomorphism
 $(C^-)^s(\tau^-(F))$:
 $$
 \gamma_{-}\mapsto A^s_{-}:= exp((C^-)^s(\tau^-(F))):P^s \longrightarrow P^s;
 $$
to the path $\gamma_{+}$ we assign the exponential of a Laurent polynomial
$Q$ of $A^s_{-}$:
$$
 \gamma_{+}\mapsto A^s_{+}:= exp(Q(A^s_{-})):P^s \longrightarrow P^s.
$$
Thus we obtain the following.
\begin{pro}\label{pro:LS}
	Let $([\xi],[\phi]) \in {\mathfrak{L}}_l$. Given a representation
	$$
	\rho_c(\xi\otimes\phi)=\{\alpha^{t,s}_c(\xi,\phi):P^s \longrightarrow P^t\}
	$$
	 of the quiver $PG_l$
	and its extension $\widehat{\rho}_c(F,G)(\xi\otimes\phi)$ gives an elliptic curve
	$$
	\Gamma:=\Gamma_{c, (F,G),\xi\otimes\phi }
	$$ 
	with $l$ marked points $\{x_i:=x_i (c, (F,G),\xi\otimes\phi )\}_{i\in[0,l-1]}$. This curve
	comes along with families of local systems: for every vector space $P^s$
	of the quiver representation $\widehat{\rho}_c(F,G)(\xi\otimes\phi)$ we have
	$$
{\mathfrak{ls}}^{P^s}_{\Gamma}:	\CC[t,t^{-1}] \longrightarrow {\mathfrak{Loc}}(P^s)
	$$
	which assigns to a Laurent polynomial $Q$ in $\CC[t,t^{-1}]$ the local system ${\mathfrak{ls}}^{P^s}_{\Gamma} (Q)$ determined by the representation of
	the fundamental group of ${\Gamma}$:
	$$
	Rep^{P^s}(\widehat{\rho}_c (F,G)(\xi\otimes\phi)): \pi_1(\Gamma)=\ZZ\{\gamma_{-},\gamma_{+}\} \longrightarrow GL(P^s)
	$$
	which sends the generators $\gamma_{\pm}$ to the following linear automorphisms of $P^s$:
	$$
	\gamma_{-}\mapsto A^s_{-}:= exp((C^-)^s(\tau^-(F))):P^s \longrightarrow P^s,
	$$
	$$
	\gamma_{+}\mapsto A^s_{+}:= exp(Q(A^s_{-})):P^s \longrightarrow P^s.
	$$
	In addition, the marked points $\{x_i \}_{i\in[0,l-1]}$
	are vertices of the dual graph $(\widehat{PG}_l)^{\vee}$ whose dual edges
	give rise to `morphisms' between those points, that is, every edge
	$e^{\vee}$ of the dual graph $(\widehat{PG}_l)^{\vee}$ gives rise to the
	isomorphism
	$$
	A(e^{\vee}):  {\mathfrak{ls}}^{P^s}_{\Gamma}(Q)_{t(e^{\vee})} \longrightarrow {\mathfrak{ls}}^{P^s}_{\Gamma}(Q)_{h(e^{\vee})},
	$$
	where $t(e^{\vee})$ (resp. $h(e^{\vee})$) is the tale (resp. head) of the edge $e^{\vee}$, and ${\mathfrak{ls}}^{P^s}_{\Gamma_{ m_{\widehat{\rho}_c (F,G)} }}(Q)_x$ denotes the fibre of the local system ${\mathfrak{ls}}^{P^s}_{\Gamma_{ m_{\widehat{\rho}_c (F,G)} }}(Q)$ at a point $x \in  \Gamma$.
\end{pro}
\section{Sheaf version of $([\xi],[\phi])$-filtrations}
This section is of a technical nature: we construct the sheaf version of $([\xi],[\phi])$-filtrations.
 We work on the open stratum $\Sigma^{\circ}_r$ which will be fixed once and for all, so we simplify the notation by writing $\Sigma$ for that stratum. The reader will have to  recall  the constructions of \S3:
 over the Cartesian product $Y_C=C\times \PP(H^1(\Theta_{C}))$ we have the universal extension on $Y_C$
 \begin{equation}\label{ext-univ1}
 	\xymatrix{
 		0\ar[r]& \pi^{\ast}_2 \OO_{\PP(H^1(\Theta_{C}))}(1) \ar[r]&{\EE_{\xi_{univ}}} \ar[r]& \pi^{\ast}_1 \OO_C (K_C) \ar[r]&0,  
 	}
 \end{equation}
where $\pi_i$, $i=1,2$, are the projections of $Y_C$ on the $i$-th factor.
 Taking the direct image under $\pi_2$ of the universal extension in \eqref{ext-univ1} we obtain
 $$
 \xymatrix@C=12pt{
 	0\ar[r]&  \OO_{\PP(H^1(\Theta_{C}))}(1) \ar[r]&\pi_{2 \ast} {\EE_{\xi_{univ}}} \ar[r]& H^0 (\OO_C (K_C)) \otimes \OO_{\PP(H^1(\Theta_{C}))} \ar[r]&H^1(\OO_C) \otimes \OO_{\PP(H^1(\Theta_{C}))} (1).  
 }
 $$
 Restricting to the stratum $\Sigma=\Sigma^0_r$ we deduce the exact sequence
 \begin{equation}\label{WSig-seq1}
 	\xymatrix{
 		0\ar[r]&  \OO_{\Sigma}(1) \ar[r]&\pi_{2 \ast} {\EE_{\xi_{univ}}} \otimes \OO_{\Sigma} \ar[r]& {\cal W}_{\Sigma} \ar[r]&0,  
 	}
 \end{equation}
see \eqref{WSig-seq} for details. This should be thought of as a sheaf version of a splitting maps used in the construction of $([\xi],[\phi])$-filtrations; the sheaf ${\cal W}_{\Sigma}$, by construction, is the kernel of the morphism
$$
\xymatrix{
H^0 (\OO_C (K_C)) \otimes \OO_{\Sigma} \ar[r]&H^1(\OO_C) \otimes \OO_{\Sigma} (1);  
}
$$
hence it is a subsheaf, actually a subbundle, of $H^0 (\OO_C (K_C)) \otimes \OO_{\Sigma}$.  
 Then we went further and constructed the morphism
\begin{equation}\label{sheaf-cxi}
	\xymatrix{
	\bigwedge^2 {\cal W}_{\Sigma} \ar[r]& \big(\HKC\otimes\OO_{\Sigma}/{\cal W}_{\Sigma}\Big)\otimes \OO_{\Sigma}(1), }
\end{equation}
the sheaf version of the map $c_{\xi}$ in Lemma \ref{lem:cxi-map}; this morphism will be used for the initial step of our construction. Before we proceed let us agree to simplify the notation as follows:

$$
V:=\HKC, \,\, V_{\Sigma}:=V\otimes\OO_{\Sigma},\,\,{\cal W}:={\cal W}_{\Sigma}, \,\,E_{\Sigma}:=\pi_{2 \ast} {\EE}_{\xi_{univ}}\otimes \OO_{\Sigma}.
$$
With these simplifications \eqref{WSig-seq1} and \eqref{sheaf-cxi} take the following form 
\begin{equation}\label{simp-eq1}
	\xymatrix@C=15pt{
 0\ar[r]& \OO_{\Sigma}(1)\ar[r]&  E_{\Sigma} \ar[r]&{\cal W}\ar[r]&0,
}
\end{equation}
\begin{equation}\label{simp-eq2}
\xymatrix@C=12pt{	
	\bigwedge^2 {\cal W} \ar[r]& \Big(V_{\Sigma}/{\cal W}\Big)\otimes \OO_{\Sigma}(1).& }
\end{equation}
We also remind the reader that the second morphism is related to the second exterior power of the exact sequence \eqref{simp-eq1} via the following  commutative diagram
\begin{equation}\label{wedge2-diag-Sig}
	\xymatrix@C=12pt{
		0\ar[d]&0\ar[d]\\
	{\cal W}\otimes \OO_{\Sigma}(1)\ar@{=}[r]\ar[d]& {\cal W}\otimes \OO_{\Sigma}(1)\ar[d]\\
	\bigwedge^2 E_{\Sigma} \ar[r]\ar[d]& V_{\Sigma}\otimes \OO_{\Sigma}(1)\ar[d]\\
	\bigwedge^2 {\cal W} \ar[r]\ar[d]& \Big(V_{\Sigma}/{\cal W}\Big)\otimes \OO_{\Sigma}(1)\ar[d]\\
	0&0
}
\end{equation}
As we said before our construction is taking place on the incidence correspondence
$$
\xymatrix{
&\PP({\cal W})\ar[ld]_{p_1} \ar[rd]^{p_2}&\\
\Sigma&&\PP(V)
}
$$
We choose the tautological line bundle $\OO_{\PP({\cal W})} (1)$ subject to
$$
p_{1\ast} \Big(\OO_{\PP({\cal W})} (1)\Big)={\cal W}^{\ast}.
$$
This gives the tautological inclusion on $\PP({\cal W})$
\begin{equation}\label{tautincl}
\OO_{\PP({\cal W})} (-1) \hookrightarrow p^{\ast}_1 ({\cal W}).
\end{equation}
Consider the pull back of the morphism \eqref{simp-eq2} by $p_1$:
$$
\mbox{$p^{\ast}_1 (\bigwedge^2 {\cal W} ) \longrightarrow p^{\ast}_1 \left(V_{\Sigma}/{\cal W} \otimes \OO_{\Sigma} (1) \right)$.}
$$
This together with the inclusion \eqref{tautincl} give the commutative diagram
$$
\xymatrix{
	\OO_{\PP({\cal W})} (-1) \otimes p^{\ast}_1 ({\cal W}) \ar[rd] \ar@{^{(}->}[r]&p^{\ast}_1 ({\cal W}) \otimes p^{\ast}_1 ({\cal W})  \ar[d]&\\
&p^{\ast}_1 (\bigwedge^2{\cal W}) \ar[r]&   p^{\ast}_1 \Big(V_{\Sigma}/{\cal W}\otimes \OO_{\Sigma}(1)\Big)
}
$$
from which we deduce the morphism
\begin{equation}\label{c0-sheaf}
	\xymatrix{
	{c}_{\Sigma} :	\OO_{\PP({\cal W})} (-1) \otimes p^{\ast}_1 ({\cal W}) \ar[r] &p^{\ast}_1 \Big(V_{\Sigma}/{\cal W}\otimes \OO_{\Sigma}(1)\Big).
	}
\end{equation}
The value of this morphism at a closed point $([\xi],[\phi])$ of $\PP({\cal W})$ is the map $c_{\xi}(\phi, \bullet)$ in the proof of Lemma \ref{lem:phi-filt}. The above gives us the initial step of our filtration
$$
{\bf W^0}:=\OO_{\PP({\cal W})} (-1) \otimes p^{\ast}_1 ({\cal W}),\,\, 
{\bf W^1}:=ker({c}_{\Sigma}).
$$
The inclusion \eqref{tautincl} and the skew symmetry of the morphism in \eqref{simp-eq2} imply the inclusion
\begin{equation}\label{taut-in-W1}
\OO_{\PP({\cal W})} (-2)=\OO_{\PP({\cal W})} (-1)\otimes \OO_{\PP({\cal W})} (-1) \hookrightarrow {\bf W^1}.
\end{equation}
To go further we need a lift of $c_{\Sigma}$ to $p^{\ast}_1 (V_{\Sigma}\otimes \OO_{\Sigma}(1))$. For this we combine the tautological inclusion \eqref{tautincl} with the pull back by $p_1$ of the exact sequence \eqref{simp-eq1} to obtain a commutative diagram 
$$
 \xymatrix@C=12pt{
 	&&0\ar[d]&0\ar[d]&\\
 	0\ar[r]& p^{\ast}_1 (\OO_{\Sigma}(1))\ar[r] \ar@{=}[d]& \widetilde{E} \ar[r] \ar[d]& \OO_{\PP({\cal W})} (-1)\ar[r] \ar[d]&0\\
 	0\ar[r]&p^{\ast}_1 ( \OO_{\Sigma}(1))\ar[r]&  p^{\ast}_1 (E_{\Sigma} ) \ar[r]&p^{\ast}_1 ({\cal W})\ar[r]&0,
 }
$$
where $\widetilde{E}$ is the inverse image of $\OO_{\PP({\cal W})} (-1)$ under the epimorphism of the bottom exact sequence. In particular, that sheaf is locally free, so tensoring the bottom exact sequence with $\widetilde{E}$ gives the exact sequence
$$
\xymatrix@C=12pt{
	0\ar[r]& \widetilde{E}\otimes p^{\ast}_1 (\OO_{\Sigma}(1) )\ar[r] & \widetilde{E} \otimes p^{\ast}_1 ( E_{\Sigma})\ar[r] & \widetilde{E}\otimes p^{\ast}_1 ( {\cal W})\ar[r] &0.
}
$$
We now bring in the morphism
\begin{equation}\label{wedge2-Sig-1}
\mbox{$\bigwedge^2 E_{\Sigma} \longrightarrow V_{\Sigma}\otimes \OO_{\Sigma}(1)$,}
\end{equation}
see the middle row of the diagram \eqref{wedge2-diag-Sig}, to obtain the the vertical arrow
\begin{equation}\label{vertarrow}
\xymatrix@C=12pt{
	0\ar[r]& \widetilde{E}\otimes p^{\ast}_1 (\OO_{\Sigma}(1))\ar[r] & \widetilde{E} \otimes p^{\ast}_1 ( E_{\Sigma})\ar[r] \ar[d]& \widetilde{E}\otimes p^{\ast}_1 ( {\cal W})\ar[r] &0\\
0\ar[r]	&p^{\ast}_1 ({\cal W}\otimes \OO_{\Sigma}(1))\ar[r]&p^{\ast}_1 (V_{\Sigma}\otimes \OO_{\Sigma}(1))\ar[r]&p^{\ast}_1 (V_{\Sigma}/{\cal W} \otimes \OO_{\Sigma}(1))\ar[r]&0
}
\end{equation}
Let us examine its restriction to $\widetilde{E}\otimes p^{\ast}_1 (\OO_{\Sigma}(1))$:
$$
\xymatrix@C=12pt{
	0\ar[r]& p^{\ast}_1 (\OO_{\Sigma}(2))\ar[r] & \widetilde{E}\otimes p^{\ast}_1 (\OO_{\Sigma}(1)) \ar[r]\ar[d] & \OO_{\PP({\cal W})} (-1)\otimes p^{\ast}_1 (\OO_{\Sigma}(1))\ar[r] &0\\
	&&p^{\ast}_1 (V_{\Sigma}\otimes \OO_{\Sigma}(1))&&
}
$$
where the top exact sequence is the defining exact sequence for $\widetilde{E}$ tensored with $p^{\ast}_1 (\OO_{\Sigma}(1))$. By the skew-symmetry of \eqref{wedge2-Sig-1}, the first term from the left in the exact sequence goes to zero, while the quotient sheaf goes to the tautological subsheaf
$$
 \OO_{\PP({\cal W})} (-1)\otimes p^{\ast}_1 (\OO_{\Sigma}(1)) \stackrel{\cong}{\longrightarrow} \OO_{\PP({\cal W})} (-1)\otimes p^{\ast}_1 (\OO_{\Sigma}(1)) \hookrightarrow p^{\ast}_1 ({\cal W})\otimes  p^{\ast}_1 (\OO_{\Sigma}(1)))=p^{\ast}_1 ({\cal W}\otimes \OO_{\Sigma}(1)).
 $$
 Furthermore, according to \eqref{taut-in-W1} the tautological subsheaf sits in ${\bf W^1} \otimes \OO_{\PP({\cal W})} (1)$
 $$
 \OO_{\PP({\cal W})} (-1) \hookrightarrow {\bf W^1} \otimes \OO_{\PP({\cal W})} (1) \subset p^{\ast}_1 ({\cal W})={\bf W^0}\otimes \OO_{\PP({\cal W})} (1).
$$
Thus the restriction of the vertical arrow in \eqref{vertarrow} to $\widetilde{E}\otimes p^{\ast}_1 (\OO_{\Sigma}(1))$ is as follows
$$
\widetilde{E}\otimes p^{\ast}_1 (\OO_{\Sigma}(1))\longrightarrow \OO_{\PP({\cal W})} (-1)\otimes p^{\ast}_1 (\OO_{\Sigma}(1)) \hookrightarrow {\bf W^1} \otimes \OO_{\PP({\cal W})} (1) \otimes p^{\ast}_1 (\OO_{\Sigma}(1)) \hookrightarrow p^{\ast}_1 (V_{\Sigma}\otimes \OO_{\Sigma}(1)).
$$
This implies that the diagram \eqref{vertarrow} can be completed to the commutative diagram
$$
\xymatrix@C=12pt{
	0\ar[r]& \widetilde{E}\otimes p^{\ast}_1 (\OO_{\Sigma}(1))\ar[r]\ar[d] & \widetilde{E} \otimes p^{\ast}_1 ( E_{\Sigma})\ar[r] \ar[d]& \widetilde{E}\otimes p^{\ast}_1 ( {\cal W})\ar[r] \ar[d]&0\\
	0\ar[r]	&p^{\ast}_1 ({\cal W}\otimes \OO_{\Sigma}(1))\ar[r]&p^{\ast}_1 (V_{\Sigma}\otimes \OO_{\Sigma}(1))\ar[r]&p^{\ast}_1 (V_{\Sigma}/{\cal W} \otimes \OO_{\Sigma}(1))\ar[r]&0
}
$$
In addition, the vertical arrow on the left factors through the tautological
line bundle
$$
\xymatrix@C=12pt{
	0\ar[r]& \widetilde{E}\otimes p^{\ast}_1 (\OO_{\Sigma}(1))\ar[r]\ar[dd] \ar[ld]& \widetilde{E} \otimes p^{\ast}_1 ( E_{\Sigma})\ar[r] \ar[dd]& \widetilde{E}\otimes p^{\ast}_1 ( {\cal W})\ar[r] \ar[dd]&0\\
	\OO_{\PP({\cal W})} (-1)\otimes p^{\ast}_1 (\OO_{\Sigma}(1))\ar[rd]&&&&\\
	0\ar[r]	&p^{\ast}_1 ({\cal W}\otimes \OO_{\Sigma}(1))\ar[r]&p^{\ast}_1 (V_{\Sigma}\otimes \OO_{\Sigma}(1))\ar[r]&p^{\ast}_1 (V_{\Sigma}/{\cal W} \otimes \OO_{\Sigma}(1))\ar[r]&0
}
$$
The vertical arrow on the right in the above diagram factors through the morphism $c_{\Sigma}$ as follows:
$$
\xymatrix@C=12pt{
	\widetilde{E}\otimes p^{\ast}_1 ( {\cal W})\ar[rd] \ar[dd]&\\
	&\OO_{\PP({\cal W})} (-1)\otimes p^{\ast}_1 ({\cal W} )\ar^{c_{\Sigma}}[ld]\\
	p^{\ast}_1 (V_{\Sigma}/{\cal W} \otimes \OO_{\Sigma}(1))&
}
	$$
	Combining with the previous diagram gives the commutative diagram
	{\small
	\begin{equation}\label{sheaf-lift}
		\xymatrix@C=9pt@R=10pt{
			0\ar[r]& \widetilde{E}\otimes p^{\ast}_1 (\OO_{\Sigma}(1))\ar[r]\ar[dd] \ar[ld]& \widetilde{E} \otimes p^{\ast}_1 ( E_{\Sigma})\ar[r] \ar[dd]& \widetilde{E}\otimes p^{\ast}_1 ( {\cal W})\ar[r] \ar[dd]\ar[rd]&0\\
			{\scriptscriptstyle{\OO_{\PP({\cal W})} (-1)\otimes p^{\ast}_1 (\OO_{\Sigma}(1))}}\ar[rd]&&&&{\scriptscriptstyle{\OO_{\PP({\cal W})} (-1)\otimes p^{\ast}_1 ({\cal W})\ar[ld]^{c_{\Sigma}}}}\\
			0\ar[r]	&p^{\ast}_1 ({\cal W}\otimes \OO_{\Sigma}(1))\ar[r]&p^{\ast}_1 (V_{\Sigma}\otimes \OO_{\Sigma}(1))\ar[r]&p^{\ast}_1 (V_{\Sigma}/{\cal W} \otimes \OO_{\Sigma}(1))\ar[r]&0
		}
	\end{equation}
}
\noindent
The vertical arrow in the middle of the diagram will play the role of a lifting of $c_{\Sigma}$. In particular, it will enable us to define further steps of the sheaf version of $([\xi],[\phi])$-filtration. Namely, we take the subsheaf 
$$
{\bf W^1}\otimes \OO_{\PP({\cal W} }(1)=ker(c_{\Sigma}) \otimes \OO_{\PP({\cal W})} (1) \subset {\bf W^0} \otimes \OO_{\PP({\cal W})} (1)= p^{\ast}_1 ({\cal W})
$$
and look at its preimage $\widetilde{{\bf W^1}}$ in $ p^{\ast}_1 ( E_{\Sigma}) $ to obtain the diagram
\begin{equation}\label{W1-diag}
\xymatrix@C=12pt{
	&&0\ar[d]&0\ar[d]&\\
	0\ar[r]&p^{\ast}_1 (\OO_{\Sigma} (1)) \ar[r]\ar@{=}[d]&\widetilde{{\bf W^1}}\ar[r]\ar[d]&{\bf W^1}\otimes \OO_{\PP({\cal W} }(1)\ar[r]\ar[d]&0\\
	0\ar[r]&p^{\ast}_1 (\OO_{\Sigma} (1)) \ar[r]&p^{\ast}_1 (E_{\Sigma} ) \ar[r]&p^{\ast}_1 ({\cal W})\ar[r]&0
}
\end{equation}
	Tensoring with $\widetilde{E}$ gives
	$$
	\xymatrix@C=12pt{
		&&0\ar[d]&0\ar[d]&\\
		0\ar[r]&\widetilde{E}\otimes p^{\ast}_1 (\OO_{\Sigma} (1)) \ar[r]\ar@{=}[d]&\widetilde{E}\otimes \widetilde{{\bf W^1}}\ar[r]\ar[d]&\widetilde{E}\otimes{\bf W^1}\otimes \OO_{\PP({\cal W} }(1)\ar[r]\ar[d]&0\\
		0\ar[r]&\widetilde{E}\otimes p^{\ast}_1 (\OO_{\Sigma} (1)) \ar[r]&\widetilde{E}\otimes p^{\ast}_1 (E_{\Sigma} ) \ar[r]&\widetilde{E}\otimes p^{\ast}_1 ({\cal W})\ar[r]&0
	}
	$$
 Applying the vertical arrows of the diagram \eqref{sheaf-lift} to the top row and using the fact that $c_{\Sigma}$ vanishes on ${\bf W^{1}}$ we deduce
 {\small
 $$
 \xymatrix@C=12pt{
 	0\ar[r]&\widetilde{E}\otimes p^{\ast}_1 (\OO_{\Sigma} (1)) \ar[r]\ar[d]&\widetilde{E}\otimes \widetilde{{\bf W^1}}\ar[r]\ar[d]&\widetilde{E}\otimes{\bf W^1}\otimes \OO_{\PP({\cal W} }(1)\ar[r]\ar[d]&0\\
 	0\ar[r]&\OO_{\PP({\cal W})}(-1)\otimes p^{\ast}_1 (\OO_{\Sigma} (1)) \ar[r]& p^{\ast}_1 ({\cal W} \otimes \OO_{\Sigma} (1) ) \ar[r]& p^{\ast}_1 ({\cal W}) / \Big({\bf W^1}\otimes \OO_{\PP({\cal W} }(1)\Big) \otimes p^{\ast}_1 (\OO_{\Sigma} (1) )&
 }
 $$
}
 The right vertical arrow factors to give the commutative triangle
 $$
\xymatrix@C=12pt{
&\widetilde{E}\otimes{\bf W^1}\otimes \OO_{\PP({\cal W} }(1)\ar[rd]\ar[dd]&\\
&&{\bf W^1} \ar[ld]\\
\Big({\bf W^0}/{\bf W^1}\Big)\otimes \OO_{\PP({\cal W} }(1)\otimes p^{\ast}_1 (\OO_{\Sigma} (1)\ar@{=}[r] &p^{\ast}_1 ({\cal W}) / \Big({\bf W^1}\otimes \OO_{\PP({\cal W} }(1)\Big) \otimes p^{\ast}_1 (\OO_{\Sigma} (1) )
}
$$
where we use the equality ${\bf W^0}\otimes\OO_{\PP({\cal W} }(1) =p^{\ast}_1 ({\cal W})$. The resulting arrow
$$
c^1_{\Sigma}: {\bf W^1} \longrightarrow \Big({\bf W^0}/{\bf W^1}\Big)\otimes \OO_{\PP({\cal W} }(1)\otimes p^{\ast}_1 (\OO_{\Sigma} (1))
$$
is the sheaf version of the map $c^{(1)}_{\xi} (\phi,\bullet)$ in Proposition \ref{pro:xi-phi-filt} and the next step is defined by setting
\begin{equation}\label{filtW2-Sigm}
	{\bf W^2}:=ker(c^1_{\Sigma}).
\end{equation}
Replacing ${\bf W^1}$ by ${\bf W^2}$ in the construction of the diagram \eqref{W1-diag} and repeating the subsequent steps give us the next term of filtration.	 This gives an inductive procedure to obtain the filtration.
We record this in the following statement.
\begin{pro}
	There is a distinguished morphism
	$$
	c_{\Sigma}: \OO_{\PP({\cal W})} (-1)\otimes p^{\ast}_1 ({\cal W}) \longrightarrow  p^{\ast}_1 \Big(V_{\Sigma}/{\cal W} \otimes \OO_{\Sigma} (1)\Big)
	$$
defined in \eqref{c0-sheaf}. Starting with 
	$$
	{\bf W^0}:=\OO_{\PP({\cal W})} (-1)\otimes p^{\ast}_1 ({\cal W}), \,\,
	{\bf W^1}:=ker(c_{\Sigma}),
	$$
	we have an inductive procedure to define for every $i\geq 1$ a morphism
	$$
	c^i_{\Sigma}: {\bf W^i} \longrightarrow \Big({\bf W^{i-1}}/{\bf W^i}\Big)\otimes \OO_{\PP({\cal W} }(1)\otimes p^{\ast}_1 (\OO_{\Sigma} (1).
	$$
	 Setting
	 $$
	 {\bf W^{i+1}}:=ker(c^i_{\Sigma})
	 $$
	 gives the decreasing filtration
	 $$
	 {\bf W^0} \supset {\bf W^1} \supset \cdots \supset  {\bf W^{i}} \supset {\bf W^{i+1}} \supset \cdots
	 $$
	 
	 At every closed point $([\xi],[\phi]) \in \PP({\cal W})$ the fibres
	 ${\bf W^{i}}([\xi],[\phi])$ of ${\bf W^{i}}$ at $([\xi],[\phi])$ give the $([\xi],[\phi])$-filtration in Proposition \ref{pro:xi-phi-filt}.
\end{pro}

For each open stratum $\Sigma^{0}_r$ of the Griffiths stratification \eqref{rkstrat}, for $r\in [2,g-3]$, the incidence correspondence 
$$
\PP({\cal W}_{\Sigma^{0}_r})=\{ ([\xi],[\phi]) \in \Sigma^{0}_r \times \PP(\HKC) | \, \text{$\xi \phi =0$ in $H^1(\OO_C)$}\}
	$$
acquires an additional stratification into locally closed strata parametrizing points {\small$([\xi],[\phi])$}, where the partition
${\bf h}([\xi],[\phi])$ associated to the {\small$([\xi],[\phi])$}-filtration   remains constant. Thus we have the decomposition of $\PP({\cal W}_{\Sigma^{0}_r})$ into disjoint strata labeled by certain partitions 
$$
{\bf h}=(h^0 , h^1 ,\ldots)
$$
of $(g-r)=rk ({\cal W}_{\Sigma^{0}_r})$, where the parts $h^i$ are the dimensions  of the quotients {\small$W^{i}_{\xi} ([\phi])/ W^{i+1}_{\xi}([\phi])$}.

Observe that all sheaves in the filtration $\{{\bf W^{i}}\}$ are torsion free. So their ranks are well defined and we have the stratum, call it
$U_{{\Sigma^{0}_r}}$, where all sheaves of the filtration are locally free. This is a Zariski dense open subset of  $\PP({\cal W}_{\Sigma^{0}_r})$ labeled by the partition denoted
$$
{\bf h}_{\Sigma^{0}_r}=(h^0_{\Sigma^{0}_r} , h^1_{\Sigma^{0}_r} ,\ldots)
$$
where the parts are 
$$
h^i_{\Sigma^{0}_r}=rk({\bf W^{i}})-rk({\bf W^{i+1}}).
$$
Thus the open stratum comes along with the morphism
\begin{equation}\label{pSig}
p_{\Sigma^{0}_r} : U_{{\Sigma^{0}_r}} \longrightarrow {\bf F}({\bf h}_{\Sigma^{0}_r};{\cal W}_{\Sigma^{0}_r} )
\end{equation}
where ${\bf F}({\bf h}_{\Sigma^{0}_r};{\cal W}_{\Sigma^{0}_r} )$ is the relative variety of flags of ${\cal W}_{\Sigma^{0}_r} $ of type determined by the partition ${\bf h}_{\Sigma^{0}_r}$, that is, ${\bf F}({\bf h}_{\Sigma^{0}_r};{\cal W}_{\Sigma^{0}_r})$ is the fibre bundle
$$
\pi: {\bf F}({\bf h}_{\Sigma^{0}_r};{\cal W}_{\Sigma^{0}_r}) \longrightarrow \Sigma^0_r
$$
over $\Sigma^0_r$ whose fibre over $[\xi] \in \Sigma^0_r$ is the variety of flags  ${\bf F}({\bf h}_{\Sigma^{0}_r};W_{\xi})$ of type ${\bf h}_{\Sigma^{0}_r}$ in $W_{\xi}$. The morphism $p_{\Sigma^{0}_r}$ fits into the diagram
$$
\xymatrix{
 U_{{\Sigma^{0}_r}} \ar[r]^(.4){p_{\Sigma^{0}_r}} \ar[dr]_{p_1}& {\bf F}({\bf h}_{\Sigma^{0}_r};{\cal W}_{\Sigma^{0}_r} ) \ar[d]^{\pi}\\
 &\Sigma^{0}_r
}
$$
where $p_1$ is the restriction to $U_{{\Sigma^{0}_r}}$ of the structure projection 
$$
p_1 : \PP({\cal W}_{\Sigma^{0}_r}) \longrightarrow \Sigma^{0}_r.
$$
The subset $U_{{\Sigma^{0}_r}}$ constitutes the top stratum of our refinement stratification and the map $p_{\Sigma^{0}_r}$ could be viewed as a sort of {\it secondary} infinitesimal variation of Hodge structure attached to it.

\vspace{1cm}
\begin{flushright}         
	Universit\'e d'Angers\\
	D\'epartement de Math\'ematiques
	\\
	2, boulevard Lavoisier\\
	49045 ANGERS Cedex 01 \\
	FRANCE\\
	{\em{E-mail addres:}} reider@univ-angers.fr
\end{flushright} 
\end{document}